\documentclass [12pt,twoside]{book}%
\usepackage[T1]{fontenc}		
\usepackage[utf8]{inputenc}
\usepackage{amsmath,amsfonts}
\usepackage{amsthm, amssymb, mathrsfs, setspace, fancyhdr}
\usepackage{pifont}
\usepackage{yfonts}
\usepackage{newpxtext,newpxmath}
\usepackage{enumerate}              
\usepackage[all,2cell]{xy}\UseAllTwocells 
\usepackage{mathtools}		
\usepackage{thmtools} 		
\usepackage{bm}					  
\usepackage{parskip}				
\usepackage{graphicx}
\usepackage{xcolor}
\usepackage{array}          
\usepackage{hhline}
\usepackage{multicol}
\usepackage{booktabs,makecell}  			
\usepackage{verbatim}     
\usepackage{mathdots}			
\usepackage{float}
\usepackage{anyfontsize}
\DeclareMathAlphabet{\EuRoman}{U}{eur}{m}{n}     
\SetMathAlphabet{\EuRoman}{bold}{U}{eur}{m}{n}   %
\usepackage{euscript}   
\usepackage{titlesec}
\usepackage{framed}				
\usepackage{index}
\makeindex%
\newindex{not}{ndx}{nnd}{Notation Index}  
\newindex{acr}{cdx}{cnd}{Acronyms} 
\usepackage{appendix}

\def\cmykBlack{0 0 0 1}

\def\pdfsetcolor#1{\pdfliteral{#1 k}}

\def\maincolor{\cmykBlack}
\pdfsetcolor{\maincolor}

\def\makefootline{
    \baselineskip24pt
    \line{\pdfsetcolor{\maincolor}\the\footline}}

\def\makeheadline{%
    \edef\M{\topmark}
    \ifx\M\empty\let\M=\maincolor\fi
    \vbox to 0pt{\vskip-22.5pt
        \line{\vbox to8.5pt{}%
        \pdfsetcolor{\maincolor}\the\headline\pdfsetcolor{\M}}\vss}%
    \nointerlineskip}
\usepackage
[colorlinks=true,
  hyperindex=true,  
  pagebackref=true,   
  anchorcolor=black, 
  linkcolor=blue, 
  citecolor=magenta, 
  plainpages=true, 
  ]{hyperref}
\hypersetup{pdfauthor=George Peschke and Tim Van der Linden}
\hypersetup{pdftitle=A Homological View of Categorical Algebra}
\definecolor{theoremcaption}{rgb}{0,0,0}		
\definecolor{proofcaption}{rgb}{0,0,0}				
\newtheoremstyle{captionstyle}
{2ex}{}
{}{0pt}
{\bfseries}     
{}                  
{\newline}
{\color{theoremcaption} \thmnumber{#2} \thmname{#1}\hfill\thmnote{#3} } 

\makeatletter		
\define@key{thmdef}{tag}[{}]{%
	\thmt@trytwice{}{%
		\addtotheorempreheadhook[\thmt@envname]{#1 AA}%
		\addtotheorempostfoothook[\thmt@envname]{BB}%
	}%
}
\define@key{thmt@tag}{PSA}{\textbf{{\color{Maroon}(PSA)}}}
\makeatother

\declaretheoremstyle[
	spaceabove=2ex,
	spacebelow=2ex,
	headfont=\normalfont\bfseries,
	notefont=\mdseries\bfseries,
	bodyfont=\normalfont,
	notebraces={\hfill}{},
	postheadspace={\newline}
]{tagstyle}

\theoremstyle{captionstyle}  
\newtheorem{theorem}[subsection]{Theorem}
\newtheorem{definition}[theorem]{Definition}
\newtheorem{proposition}[theorem]{Proposition}
\newtheorem{corollary}[theorem]{Corollary}
\newtheorem{lemma}[theorem]{Lemma}
\newtheorem{example}[theorem]{Example}
\newtheorem{remark}[theorem]{Remark}
\newtheorem{exercise}[theorem]{Exercise}

\newtheorem{notation}[theorem]{Notation}

\newtheorem{terminology}[theorem]{Terminology}

\titleformat{\subsection}
{\normalfont\bfseries}{\thesubsection}{1em}{}
\makeatletter
\renewenvironment{proof}[1][\proofname]{\vspace{-2ex}\par       
	\pushQED{\qed}%
	\normalfont \topsep6\p@\@plus6\p@\relax
	\trivlist
	\item[\hskip\labelsep
	            \color{proofcaption}\bfseries                
	            #1\@addpunct{\quad}]\ignorespaces
}{%
	\popQED\endtrivlist\@endpefalse
}
\makeatother

\newcommand{\NoProof}{{\unskip\nobreak\hfil\penalty 50\hskip 2em\hbox{}
			\nobreak\hfil$\lozenge$\parfillskip=0pt\finalhyphendemerits=0\par}}
\newcommand{\NoProofDiag}{\tag*{$\lozenge$}}
\definecolor{light-gray}{gray}{0.95}
\newenvironment{colbox}{%
	\MakeFramed{\advance\hsize-\width \FrameRestore}}
{\endMakeFramed}

\newenvironment{ulist}{			
	\begin{itemize}}{
	\end{itemize}
}

\newenvironment{thmlist}{		
	\begin{enumerate}[(i)]}{
	\end{enumerate}
}
\newenvironment{tfae}{		
	\begin{enumerate}[(I)]}{
	\end{enumerate}
}
\newenvironment{subordinate}[2][Notes]{%
	\textbf{#1\qquad #2}
	\vskip 0.1ex
	\footnotesize{%
	}
}

\newenvironment{subsubordinate}[1]{%
	\footnotesize{%
		\textbf{#1}\quad      }
}

\newenvironment{exercises}{%
	\MakeFramed{\hsize=0.95\linewidth\advance\hsize-\width\FrameRestore%
		\bigskip
		\textbf{Exercises}\vspace{-2ex}\footnotesize{
		}}
}
{\endMakeFramed}
\newdir{>>}{{}*!/3.5pt/:(1,-.2)@^{>}*!/3.5pt/:(1,+.2)@_{>}*!/7pt/:(1,-.2)@^{>}*!/7pt/:(1,+.2)@_{>}}
\newdir{ >>}{{}*!/8pt/@{|}*!/3.5pt/:(1,-.2)@^{>}*!/3.5pt/:(1,+.2)@_{>}}
\newdir{ |>}{{}*!/-3.5pt/@{|}*!/-8pt/:(1,-.2)@^{>}*!/-8pt/:(1,+.2)@_{>}}
\newdir{ >}{{}*!/-8pt/@{>}}
\newdir{>}{{}*:(1,-.2)@^{>}*:(1,+.2)@_{>}}
\newdir{<}{{}*:(1,+.2)@^{<}*:(1,-.2)@_{<}}

\newcommand{\KernelArr}[4]{\xymatrix@R=5ex@C=#3em{ #1 \ar@{{ |>}->}[r]^-{#4} & #2 } }
\newcommand{\CoKernelArr}[4]{\xymatrix@R=5ex@C=#3em{ #1 \ar@{-{ >>}}[r]^-{#4} & #2 } }
\newcommand{\PullLU}[1]{\ar@{}[#1]|-{%
\begin{picture}(10,10)%
\put(1,1){\line(1,0){8}}%
\put(9,1){\line(0,1){8}}%
\put(1,1){\circle*{2}}%
\put(9,1){\circle*{2}}%
\put(9,9){\circle*{2}}%
\put(1,9){\circle{2}}
\end{picture} } }
\newcommand{\PullRU}[1]{\ar@{}[#1]|-{%
\begin{picture}(10,10)%
\put(1,1){\line(1,0){8}}%
\put(1,1){\line(0,1){8}}%
\put(1,1){\circle*{2}}%
\put(9,1){\circle*{2}}%
\put(9,9){\circle{2}}%
\put(1,9){\circle*{2}}
\end{picture} } }
\newcommand{\PullLD}[1]{\ar@{}[#1]|-{%
\begin{picture}(10,10)%
\put(1,9){\line(1,0){8}}%
\put(9,1){\line(0,1){8}}%
\put(1,1){\circle{2}}%
\put(9,1){\circle*{2}}%
\put(9,9){\circle*{2}}%
\put(1,9){\circle*{2}}
\end{picture} } }
\newcommand{\PullRD}[1]{\ar@{}[#1]|-{%
\begin{picture}(10,10)%
\put(1,9){\line(1,0){8}}%
\put(1,1){\line(0,1){8}}%
\put(1,1){\circle*{2}}%
\put(9,1){\circle{2}}%
\put(9,9){\circle*{2}}%
\put(1,9){\circle*{2}}
\end{picture} } }
\newcommand{\PushLU}[1]{\ar@{}[#1]|-{%
\begin{picture}(10,10)%
\put(1,1){\line(1,0){8}}%
\put(9,1){\line(0,1){8}}%
\put(1,1){\circle*{2}}%
\put(9,1){\circle*{2}}%
\put(9,9){\circle*{2}}%
\put(1,9){\circle{2}}
\end{picture} } }
\newcommand{\PushRU}[1]{\ar@{}[#1]|-{%
\begin{picture}(10,10)%
\put(1,1){\line(1,0){8}}%
\put(1,1){\line(0,1){8}}%
\put(1,1){\circle*{2}}%
\put(9,1){\circle*{2}}%
\put(9,9){\circle{2}}%
\put(1,9){\circle*{2}}
\end{picture} } }
\newcommand{\PushLD}[1]{\ar@{}[#1]|-{%
\begin{picture}(10,10)%
\put(1,9){\line(1,0){8}}%
\put(9,1){\line(0,1){8}}%
\put(1,1){\circle{2}}%
\put(9,1){\circle*{2}}%
\put(9,9){\circle*{2}}%
\put(1,9){\circle*{2}}
\end{picture} } }
\newcommand{\PushRD}[1]{\ar@{}[#1]|-{%
\begin{picture}(10,10)%
\put(1,9){\line(1,0){8}}%
\put(1,1){\line(0,1){8}}%
\put(1,1){\circle*{2}}%
\put(9,1){\circle{2}}%
\put(9,9){\circle*{2}}%
\put(1,9){\circle*{2}}
\end{picture} } }
\newcommand{\BiCart}[1]{\ar@{}[#1]|-{%
\begin{picture}(10,10)%
\put(1,1){\line(1,0){8}}%
\put(1,9){\line(1,0){8}}%
\put(1,1){\line(0,1){8}}%
\put(9,1){\line(0,1){8}}%
\put(1,1){\circle*{2}}%
\put(9,1){\circle*{2}}%
\put(9,9){\circle*{2}}%
\put(1,9){\circle*{2}}
\end{picture} } }
\newcommand{\IndSep}{\qquad}

\def\rtop{\text{\tiny\rotatebox[origin=c]{-45}{$\bot$}}}
\def\ltop{\text{\tiny\rotatebox[origin=c]{45}{$\bot$}}}
\newcommand{\del}{\partial}
\newcommand{\To}{\Rightarrow}
\renewcommand{\implies}{\To}
\newcommand{\normal}{\ensuremath{\lhd}}
\newcommand{\TT}{\ensuremath{\mathbb{T}}}
\newcommand{\Defn}[1]{\emph{#1}}
\newcommand{\DefEq}{\coloneq} 		
\newcommand{\EqDef}{\eqcolon} 		
\newcommand{\hy}{\text{-}}													
\newcommand{\Poincare}{Poincaré}
\newcommand{\BCPtdCat}{pointed category}  
\newcommand{\NEM}{normal epi/mono\ }

\newcommand{\NENMComp}[1]{\mathit{EM}(#1)}		
\newcommand{\XRA}[1]{\xrightarrow{\ #1\ }}
\newcommand{\XLA}[1]{\xleftarrow{\ #1\ }}

\newcommand{\from}{\colon}				
\newcommand{\Comp}{\raisebox{0.1ex}{\ensuremath{\,\scriptstyle{\circ}}\,}}
\newcommand{\IdMap}{1}												
\newcommand{\IdMapOn}[1]{1_{#1}}	
\newcommand{\InclsnOf}[1]{\textit{i}_{#1}}		
\newcommand{\Mono}{\rightarrowtail}			
\newcommand{\Epi}{\twoheadrightarrow}			
\newcommand{\NEpi}{-\hspace{-0.8ex}-\hspace{-0.7ex}\triangleright}	
\newcommand{\NMono}{\ \triangleright\hspace{-0.9ex}\rightarrow}			
\newcommand{\SctndEpi}[2]{{(#1\downarrow #2)}}				
\newcommand{\RtrctdMono}[2]{{(#1\uparrow #2)}}				

\newcommand{\Dgnl}{\Delta}						
\newcommand{\DgnlOn}[1]{\Delta_{#1}}	
\newcommand{\FoldOn}[1]{\nabla_{#1}}				
\newcommand{\PrjctnOnto}[1]{\textit{pr}_{#1}} 	
\newcommand{\Set}[1]{\left\{#1\right\}}		
\newcommand{\SetSlct}[2]{\left\{#1\mid #2 \right\}}		
\newcommand{\union}{\cup}									
\newcommand{\FamUnion}[2]{\bigcup_{#1}\, #2}
\newcommand{\PwrSt}[1]{\EuScript{P}(#1)}			
\newcommand{\NNr}{\mathbb{N}}		
\newcommand{\ZNr}{\mathbb{Z}}		
\newcommand{\ZMod}[1]{\mathbb{Z}\hskip -.2em/\hskip -.1em #1}	
\newcommand{\QNr}{\mathbb{Q}}		
\newcommand{\RNr}{\mathbb{R}}		
\newcommand{\CyclcGrp}[1]{C_{#1}}					
\newcommand{\EndRng}[1]{\textit{End}(#1)}         

\newcommand{\Hom}[2]{\textit{Hom}\left(#1,#2\right)}
\newcommand{\HomIn}[3]{\textit{Hom}_{#1}\left(#2,#3\right)}
\newcommand{\HomBsd}[2]{\textit{Hom}_{\ast}\left(#1,#2\right)}      
\newcommand{\Ext}{\ensuremath{\mathit{Ext}}}
\newcommand{\AdjUnit}{\eta}										     
\newcommand{\AdjUnitOn}[1]{\eta_{#1}}					     

\newcommand{\SmallCats}{\EuRoman{C{\kern-0.17ex}a{\kern-0.12ex}t}}			
\newcommand{\LocSmallCats}{\EuRoman{C{\kern-0.08ex}A{\kern-0.25ex}T}}		
\newcommand{\Sets}{\EuRoman{S{\kern-0.12ex}e{\hskip-0.14ex}t}}			                       
\newcommand{\SetsBsd}{\EuRoman{S{\kern-0.12ex}e{\hskip-0.14ex}t}_{\ast}\,}	
\newcommand{\SetsBsdOp}{\EuRoman{S{\kern-0.12ex}e{\hskip-0.14ex}t}_{\ast}^{\op}\,}	
\newcommand{\Magmas}{\EuRoman{M{\kern -.12ex}a{\kern -.15ex}g}}			
\newcommand{\Monoids}{\EuRoman{M{\kern -.12ex}n{\kern -.15ex}d}}			
\newcommand{\CMon}{\ensuremath{\EuRoman{CMon}}}
\newcommand{\TopMonoids}{\EuRoman{Mnd}(\EuRoman{Top})}			
\newcommand{\XMod}{\ensuremath{\EuRoman{XMod}}}
\newcommand{\Grps}{\EuRoman{G{\kern -.15ex}rp}}				
\newcommand{\GrpsTF}{\EuRoman{G{\kern -.15ex}rp}_{\mathit{tf}}}				
\newcommand{\Mlc}{\ensuremath{\EuRoman{Mlc}}}
\newcommand{\KProj}{\text{$\mathbb{K}$-$\EuRoman{Proj}$}}
\newcommand{\CHopfK}{\EuRoman{CHopf}_\mathbb{K}}
\newcommand{\CCHopfK}{\EuRoman{CCHopf}_\mathbb{K}}
\newcommand{\Loop}{\EuRoman{Loop}}
\newcommand{\Tops}{\EuRoman{Top}}											
\newcommand{\HComp}{\ensuremath{\EuRoman{HComp}}}		
\newcommand{\TopGrps}{\EuRoman{Grp}(\EuRoman{Top})}							
\newcommand{\CHTops}{\EuRoman{CHTop}}						 
\newcommand{\TopsBsd}{\EuRoman{Top}_{\ast}}	
\newcommand{\Alg}{\ensuremath{\EuRoman{Alg}}}

\newcommand{\Cat}{\ensuremath{\EuRoman{Cat}}}
\newcommand{\cat}{\ensuremath{\EuRoman{cat}}}
\newcommand{\PSACats}{\EuRoman{PS{\kern-0.15ex}A{\kern-0.15ex}C{\kern-0.15ex}a{\kern-0.15ex}t}}				
\newcommand{\SACats}{\EuRoman{S{\kern-0.15ex}A{\kern-0.15ex}C{\kern-0.15ex}a{\kern-0.15ex}t}}				
\newcommand{\ACats}{\EuRoman{AC{\kern-0.15ex}a{\kern-0.15ex}t}}	                         			
\newcommand{\AbGrps}{\EuRoman{A{\kern-0.15ex}b}}	
\newcommand{\AbGrpsTF}{\EuRoman{A{\kern-0.15ex}b}_\mathit{tf}}	
\newcommand{\Rngs}{\EuRoman{ Rn{\kern-.1ex}g}}                        
\newcommand{\RngsUnit}{\EuRoman{ Rn{\kern-.1ex}g}^{1}}          
\newcommand{\URngs}{\EuRoman{ U{\kern-.1ex}Rn{\kern-.1ex}g}}          
\newcommand{\CRngs}{\EuRoman{ Rn{\kern-.1ex}g}_{c}} 	          
\newcommand{\CRngsUnit}{\EuRoman{ Rn{\kern-.1ex}g}^{1}_{c}}     
\newcommand{\LModules}[1]{{#1}\text{-}\EuRoman{Mod}}    
\newcommand{\RModules}[1]{\EuRoman{Mod}\text{-}#1}				
\newcommand{\ModulesOver}[1]{\EuRoman{Mod}(#1)}						
\newcommand{\OrdGrp}{\ensuremath{\EuRoman{OrdGrp}}}
\newcommand{\SetPBij}{\ensuremath{\EuRoman{SetPBij}}} 
\newcommand{\Lie}{\ensuremath{\EuRoman{Lie}}}
\newcommand{\Vect}{\ensuremath{\EuRoman{Vect}}}
\newcommand{\op}{\ensuremath{\mathit{op}}}
\newcommand{\Ord}[1]{[#1]}										          
\newcommand{\RtrctCat}{\textswab{r}}			
\newcommand{\ZeroObject}{0}                           
\newcommand{\ZeroMap}{0}                                
\newcommand{\DiagObj}{\square}
\newcommand{\Ctgry}[1]{\EuScript{#1}}					
\newcommand{\SmallCtgry}[1]{#1}											
\newcommand{\SACtgry}[1]{\EuScript{#1}}			
\newcommand{\ZExact}{z-exact}									
\newcommand{\Objcts}[1]{\mathit{Ob}(#1)}	
\newcommand{\SubObjcts}[1]{\mathit{Sub}(#1)}	
\newcommand{\NSubObjcts}[1]{\mathit{NSub}(#1)}	
\newcommand{\SubObjctCat}[1]{\EuRoman{S{\kern-0.12ex}O{\kern-0.12ex}b}({\EuScript{#1}})}	
\newcommand{\SubObjctsPullFunc}{\textit{Sub}^{\ast}}									
\newcommand{\NQuoObjcts}[1]{\textit{NQuo}(#1)}	
\newcommand{\ArrowCat}[1]{\EuRoman{A{\kern-0.15ex}r{\hskip-0.1ex}r}(#1)}         
\newcommand{\NMonoCat}[1]{\EuRoman{NM}(\Ctgry{#1})}				

\newcommand{\NEpiCat}[1]{\EuRoman{NE}(\Ctgry{#1})}				

\newcommand{\SEpisIn}[1]{\EuRoman{SEpi}(\Ctgry{#1})}								
\newcommand{\SEpisInOver}[2]{\EuRoman{SEpi}_{#2}(\Ctgry{#1})}								
\newcommand{\ANPCat}[1]{\EuRoman{ANP}(\EuScript{#1})}							

\newcommand{\SESCat}[1]{\EuRoman{SES}(\Ctgry{#1})}				
\newcommand{\DExCat}[1]{\EuRoman{DEx}(\EuScript{#1})}							
\newcommand{\HExCat}[2]{\EuRoman{Ex}^{#1}(\EuScript{#2})} 
\newcommand{\SSESCat}[2][{}]{\EuRoman{SSES}_{#1}(\Ctgry{#2})}				
\newcommand{\SESCatn}[1]{\EuRoman{SES}^n(\Ctgry{#1})}				
\newcommand{\GrothendieckCContravOf}[1]{\mathfrak{Gr}^{\circ}(#1)}						
\newcommand{\FamPrdct}[2]{\prod_{#1}#2}		
\newcommand{\prdct}{\times} 					
\newcommand{\Prdct}[2]{#1 \times #2}	 	
\newcommand{\PPrdct}[3]{#1 \times #2 \times #3}	 	
\newcommand{\PrdctMapInto}[1]{( #1)}			
\newcommand{\FamCoPrdct}[2]{\coprod_{#1}#2}	
\newcommand{\CoPrdct}[2]{#1 + #2}					
\newcommand{\BiPrdct}[2]{#1\oplus #2}
\newcommand{\SumMapOutOf}[1]{\langle #1\rangle}     
\newcommand{\MapOutOf}[1]{\langle #1\rangle}     
\newcommand{\SumProdComp}[2]{\gamma_{#1,#2}}			
\newcommand{\meet}{\ensuremath{\wedge}}
\newcommand{\join}{\ensuremath{\vee}}
\newcommand{\bigmeet}{\ensuremath{\bigwedge}}
\newcommand{\bigjoin}{\ensuremath{\bigvee}}

\newcommand{\CoLim}{\textit{colim}}       								
\newcommand{\CoLimOf}[1]{\textit{colim}\left(#1\right)}		
\newcommand{\CoLimOfOver}[2]{\textit{colim}^{#2}\left(#1\right)}		
\newcommand{\Lim}{\textit{lim}}             										
\newcommand{\LimOf}[1]{\textit{lim}\left(#1\right)}        			
\newcommand{\LimOfOver}[2]{\textit{lim}_{#2}\left(#1\right)}    
\newcommand{\CoEq}[2]{\textit{coeq}\left(#1,#2\right)}		
\newcommand{\KerFunc}{\textit{ker}}		               	
\newcommand{\Ker}[1]{\textit{K}(#1)}		     	
\newcommand{\KerMap}[1]{\textit{ker}(#1)}		     	
\newcommand{\KKerMap}[1]{\textit{ker}\left(#1\right)}		     	
\newcommand{\CoKerFunc}{\textit{coker}}								
\newcommand{\CoKer}[1]{\textit{Q}(#1)}               
\newcommand{\CoKerMap}[1]{\textit{coker}(#1)}						        
\newcommand{\CCoKerMap}[1]{\textit{coker}\left(#1\right)}        
\newcommand{\Img}[1]{\textit{I}(#1)}	               
\newcommand{\ImgMap}[1]{\textit{im}(#1)}	      
\newcommand{\CoImg}[1]{\textit{coI}(#1)}	                
\newcommand{\CoImgMap}[1]{\textit{coim}(#1)}	      
\newcommand{\KrnlPr}[1]{\textit{KP}\left(#1\right)}			
\newcommand{\Hmlgy}[2]{H_{#1}(#2)}			         	
\newcommand{\HmlgyKer}[2]{H^{k}_{#1}#2}           
\newcommand{\HmlgyCoKer}[2]{H^{c}_{#1}#2}         
\newcommand{\ExtMini}[2]{\textit{ext}\left( #1,#2 \right)}
\newcommand{\Ab}{\mathit{ab}}
\newcommand{\AbCoreOf}[1]{\Ab(\Ctgry{#1})}	
\newcommand{\OneMapOn}[1]{1_{#1}}		
\newcommand{\ANKTag}{ - $\EuRoman{A{\kern-0.2ex}N{\kern-0.15ex}K}$}								
\newcommand{\ZExactTag}{ - {\color{Cerulean} $\EuRoman{z\hy Ex}$}}
\newcommand{\HSDTag}{ - {\color{Cerulean} $\EuRoman{HSD}$}}													
\newcommand{\DExTag}{ - {\color{Cerulean} $\EuRoman{DEx}$}}			
\newcommand{\DPNTag}{ - {\color{Cerulean} $\EuRoman{DPN}$}} 
\newcommand{\NDPNTag}{ - {\color{Brown} $\EuRoman{N}$} + {\color{Cerulean} $\EuRoman{DPN}$}} 
\newcommand{\NTag}{ - {\color{Brown} $\EuRoman{N}$}}																				
\newcommand{\HTag}{ - {\color{Brown} $\EuRoman{H}$}}																					
\newcommand{\SATag}{ - {\color{MidnightBlue} $\EuRoman{S{\kern-0.15ex}A}$}}			
\newcommand{\ACTag}{ - {\color{OliveGreen} $\EuRoman{A{\kern-0.15ex}C}$}}			
\newcommand{\AddTag}{ - {\color{OliveGreen} $\EuRoman{A{\kern-0.15ex}dd}$}}			
\newcommand{\AbTag}{ - {\color{OliveGreen} $\EuRoman{A{\kern-0.15ex}b}$}}			
\newcommand{\HSDInline}{(HSD)}																
\newcommand{\DPNInline}{(DPN)}																
\newcommand{\ANNInline}{(ANN)}																
\newcommand{\PNEInline}{(PNE)}																
\newcommand{\AENInline}{(AEN)}																
\newcommand{\KSGInline}{(KSG)}																
\begin{document}
\pagestyle{fancy}   
\renewcommand{\sectionmark}[1]{\markright{ #1}}
\renewcommand{\chaptermark}[1]{\markboth{ #1}{}}
\fancyhead{}
\fancyfoot{}
\fancyhead[LO]{\bfseries\footnotesize \thechapter\ \leftmark}
\fancyhead[LE]{\bfseries\footnotesize \thechapter\ \leftmark}
\fancyhead[C]{}
\fancyhead[RO]{\bfseries\footnotesize \thesection\ \rightmark}
\fancyhead[RO]{\bfseries\footnotesize \thesection\ \rightmark}
\fancyfoot[C]{\bfseries\footnotesize $\blacktriangleleft$\qquad \thepage\qquad $\blacktriangleright$}

\definecolor{Maroon}{cmyk}{0,0.87,0.68,0.32}
\definecolor{BrickRed}{cmyk}{0,0.89,0.94,0.28}
\definecolor{Red}{cmyk}{0,1,1,0}
\definecolor{OrangeRed}{cmyk}{0,1,0.50,0}
\definecolor{Salmon}{cmyk}{0,0.53,0.38,0}
\definecolor{CarnationPink}{cmyk}{0,0.63,0,0}
\definecolor{Magenta}{cmyk}{0,1,0,0}
\definecolor{VioletRed}{cmyk}{0,0.81,0,0}
\definecolor{Rhodamine}{cmyk}{0,0.82,0,0}
\definecolor{Mulberry}{cmyk}{0.34,0.90,0,0.02}
\definecolor{RedViolet}{cmyk}{0.07,0.90,0,0.34}
\definecolor{Fuchsia}{cmyk}{0.47,0.91,0,0.08}
\definecolor{Lavender}{cmyk}{0,0.48,0,0}
\definecolor{Thistle}{cmyk}{0.12,0.59,0,0}
\definecolor{Orchid}{cmyk}{0.32,0.64,0,0}
\definecolor{DarkOrchid}{cmyk}{0.40,0.80,0.20,0}
\definecolor{Purple}{cmyk}{0.45,0.86,0,0}
\definecolor{Plum}{cmyk}{0.50,1,0,0}
\definecolor{Violet}{cmyk}{0.79,0.88,0,0}
\definecolor{RoyalPurple}{cmyk}{0.75,0.90,0,0}
\definecolor{BlueViolet}{cmyk}{0.86,0.91,0,0.04}
\definecolor{Periwinkle}{cmyk}{0.57,0.55,0,0}
\definecolor{CadetBlue}{cmyk}{0.62,0.57,0.23,0}
\definecolor{CornflowerBlue}{cmyk}{0.65,0.13,0,0}
\definecolor{MidnightBlue}{cmyk}{0.98,0.13,0,0.43}
\definecolor{NavyBlue}{cmyk}{0.94,0.54,0,0}
\definecolor{RoyalBlue}{cmyk}{1,0.50,0,0}
\definecolor{Blue}{cmyk}{1,1,0,0}
\definecolor{Cerulean}{cmyk}{0.94,0.11,0,0}
\definecolor{Cyan}{cmyk}{1,0,0,0}
\definecolor{ProcessBlue}{cmyk}{0.96,0,0,0}
\definecolor{SkyBlue}{cmyk}{0.62,0,0.12,0}
\definecolor{Turquoise}{cmyk}{0.85,0,0.20,0}
\definecolor{TealBlue}{cmyk}{0.86,0,0.34,0.02}
\definecolor{Aquamarine}{cmyk}{0.82,0,0.30,0}
\definecolor{BlueGreen}{cmyk}{0.85,0,0.33,0}
\definecolor{Emerald}{cmyk}{1,0,0.50,0}
\definecolor{JungleGreen}{cmyk}{0.99,0,0.52,0}
\definecolor{SeaGreen}{cmyk}{0.69,0,0.50,0}
\definecolor{Green}{cmyk}{1,0,1,0}
\definecolor{ForestGreen}{cmyk}{0.91,0,0.88,0.12}
\definecolor{OliveGreen}{cmyk}{0.64,0,0.95,0.40}
\definecolor{RawSienna}{cmyk}{0,0.72,1,0.45}
\definecolor{Sepia}{cmyk}{0,0.83,1,0.70}
\definecolor{Brown}{cmyk}{0,0.81,1,0.60}
\definecolor{Tan}{cmyk}{0.14,0.42,0.56,0}
\definecolor{Gray}{cmyk}{0,0,0,0.50}
\definecolor{Black}{cmyk}{0,0,0,1}
\title{\bfseries\Huge A Homological View\\ of Categorical Algebra}
\author{\bfseries\Large George Peschke and Tim Van der Linden}
\date{\textbf{\today}}
\maketitle

\pagenumbering{roman}
\thispagestyle{empty}
\newpage

\pagestyle{fancy}   
%
\lhead{\bfseries\footnotesize
  Preface}
\rhead{\bfseries\footnotesize
}


{\bf Abstract}\quad We provide a foundation for working with homological and homotopical methods in categorical algebra. This involves two mutually complementary components, namely (a) the strategic selection of suitable axiomatic frameworks, some well known and some new, and (b) the development of categorical tools for effective reasoning and computing within those frameworks.

The selection of axiomatic frameworks begins `from the ground up' with z-exact categories. These are pointed categories in which every morphism has a kernel and cokernel. Then we progress all the way to abelian categories, en route meeting contexts such as Borceux--Bourn homological categories and Janelidze--M\'arki--Tholen semiabelian categories. We clarify the relationship between these axiomatic frameworks by direct comparison, but also by explaining how concrete examples fit into the selection. The outcome is a fine-grained set of criteria by which one can map varieties of algebras (in the sense of Universal Algebra) and topological models of algebraic theories into the various frameworks.

The categorical tools for effective computation deal mostly with situations involving universal factorizations of morphisms, with exact sequences, and with the homology of chain complexes.

Further, the categorical tools for computation include the `basic diagram lemmas' of homological algebra, that is the (Short) $5$-Lemma, the $(3\times 3)$-Lemma, the Snake Lemma, and applications thereof. We find that these tools are even available in some surprisingly weak categorical environments, such as the category of pointed sets. In discovering such features, we make systematic use of what we call the self-dual axis of a category.
\chapter*{Preface}
\addcontentsline{toc}{part}{Preface}

From the ground up, we develop those aspects of Categorical Algebra which allow a categorical approach to homological methods in potentially non-abelian categories. In doing so, we build on several earlier developments, and we add a number of new insights.

To recap, by the 1950's Universal Algebra had evolved into a well established subject. Then arrived the highly influential notion of \emph{abelian category} \cite{Buchsbaum:ExactCats,Tohoku}. It quickly motivated efforts to create a non-commutative analogue relevant to corresponding types of universal algebras and beyond. The conflicting demands for broadest possible applicability vs.\ strong categorical tools for computation resulted a number of approaches such as the works of Gerstenhaber, Huq, and Orzech amongst others \cite{Huq, Gerstenhaber, Orzech}.

It took till the turn of the millennium until these efforts came to a kind of confluence when Janelidze, Márki, and Tholen introduced \emph{semiabelian} categories \cite{Janelidze-Marki-Tholen}. Their work took into account then recent novelties in Category Theory such as the concepts of Barr exactness~\cite{Barr-Grillet-vanOsdol} and Bourn protomodularity~\cite{DBourn1991}. Soon thereafter, Borceux and Bourn~\cite{FBorceuxDBourn2004} proposed the broader and, as a consequence, computationally weaker environment of \emph{homological} categories. Its scope reaches beyond semiabelian categories and includes categories of topological-algebraic objects such as topological groups.

Initially, our goal was to create a systematic foundation for the fruitful use of these new categorical methods within the framework of Quillen's homotopical algebra; relying specifically on his remarkable Theorem 4 in Section \cite[II.4]{DGQuillen1967}. Then we realized that a systematic homological view of categorical algebra opens up new opportunities to the effect that key tools for computations with homological invariants are actually available in settings where, historically, they were not expected to be valid. For example `diagram lemmas' such as the (Short) $5$-Lemma, the $(\Prdct{3}{3})$-Lemma, the Snake Lemma hold true in the category $\SetsBsd$ of pointed sets.

Accordingly, we shifted our approach, and developed the from-the-ground-up view of homological methods presented here. As a `side effect' the reader will also find an ample display of categorical methods for investigation and computation; methods which unify processes which had previously been based on element manipulation. - For further details about these developments, we refer the reader to the introductions to parts, chapters, and sections.

An environmental scan reveals a rich web of axiom systems which are related to our objectives in one way or another. We dedicate an entire chapter to the purpose of making this web more transparent to the non-specialist.

\begin{subordinate}[The Role of Exercises]{}
  Most sections end with a set of exercises of varying purpose and varying level of difficulty. The reader will find ample opportunity to practice the use of the categorical methods presented here, often in contrast with arguments based on the manipulation of elements in underlying sets. Then there is the development of explicit examples in familiar categories, such as monoids, groups, etc. We also turned certain questions which we currently are unable to answer into exercises; these are labelled `ANK' for `answer not known'.
\end{subordinate}

\begin{subordinate}[Acknowledgements]{}
  We are grateful to the \emph{Pacific Institute for the Mathematical Sciences} for inviting us to visit their beautiful environment to work on this manuscript, and to the \emph{Institut de Recherche en Mathématique et Physique} of the Université catholique de Louvain for inviting George Peschke to Louvain-la-Neuve on several occasions. Many thanks to those who have read a very preliminary version of this text: Florent Afsa, Maxime Culot and Bo Shan Deval. Many thanks also to Marino Gran and Joost Vercruysse for fruitful discussions.

  Tim Van der Linden is a Senior Research Associate of the Fonds de la Recherche Scientifique-FNRS.
\end{subordinate}

\begin{subordinate}{}
  \begin{subsubordinate}{On the arXiv editions}
    For the time being, this is a living document. More topics will be added\footnote{Chapters on varieties of algebras and internal actions are in preparation, as is a more extensive biblio\-graphy.}, and the current document will evolve accordingly to form a coherent exposition. We welcome any comments and suggestions you may have.
  \end{subsubordinate}
\end{subordinate}

\bigskip\bigskip\bigskip
{\footnotesize
  \begin{multicols}{2}\raggedright\makeatletter
    George Peschke\newline
    Department of Math.\ and Stat.\ Sciences\newline
    University of Alberta\newline
    Edmonton\newline
    Canada T6G 2G1\newline
    e-mail: \href{mailto:george.peschke@ualberta.ca}{george.peschke@ualberta.ca}\newline
    {\ }\newline
    Tim Van der Linden\newline
    Institut de Recherche en Math{\'e}matique et Physique\newline
    Universit{\'e} catholique de Louvain\newline
    chemin du cyclotron 2 bte L7.01.02\newline
    B-1348 Louvain-la-Neuve\newline
    Belgium\newline
    e-mail: \href{mailto:tim.vanderlinden@uclouvain.be}{tim.vanderlinden@uclouvain.be}
  \end{multicols}   }

\pagestyle{plain} 

\tableofcontents

\newpage
\lhead{\bfseries\footnotesize
  Introduction to Part I}
\rhead{\bfseries\footnotesize
}
\pagenumbering{arabic}
\setcounter{page}{1}
\part[Homology in Pointed Categories]{Homology in Pointed Categories}
\label{part:Homology}%
\chapter*{Introduction to Part I}
\label{chap:IntroPart-I}%


\subsubsection*{A peek at the history of homology}
To provide context for what the reader will find here, let us review some highlights in the history of homology and homological methods\footnote{For more information on the history of homology, see for instance \cite{JDieudonne1989} and \cite{CAWeibel1994-HomAlg}.}.

The origins of homology are purely geometric. In Analysis Situs \cite[\S 5]{HPoincare1895}, \Poincare\ used embedded cobordisms in a given manifold to define when a formal linear combination of compact boundaryless $k$-dimensional submanifolds was homologous to another such linear combination. The evolution of combinatorial analogues followed via the study of polyhedra: Around the 1920's a chain complex of finitely generated free abelian groups was associated to a polyhedron to extract Betti numbers and torsion coefficients of its underlying topological space. Combined, these numerical invariants characterize the homology groups of the chain complex up to isomorphism by providing the rank of its free summand plus the torsion summand. %

Emmy Noether realized \cite{FHirzebruch1999-Noether},\cite[p.\ 174]{ADick1981} that these groups themselves should be the objects of primary attention. In fact, according to P.\ Alexandrov, E.\ Noether responded with these ideas to lectures given by H.\ Hopf in 1926 and 1927; see \cite[p.\ 174]{ADick1981}. An element of such a homology group in dimension $k$ could then be viewed as a combinatorial analogue of  \Poincare's homology class of formal linear combinations of compact boundaryless $k$-dimensional submanifolds of a given manifold; see \cite{HPoincare1895}.

Increasingly, the use of homology of chain complexes broadened beyond its topological-geometric origins. Eilenberg and Mac Lane introduced the theory of categories to provide a formal foundation for `natural transformations'. As part of the confluence of these developments, in 1956 H.~Cartan and S.~Eilenberg presented the first text on homological algebra ever \cite{Cartan-Eilenberg}. %

Derived functors on categories of modules, algebras, groups, etc.\ were its dominant theme. Commonalities of homological constructions in those environments were soon abstracted into `abelian categories' by  A. Grothendieck \cite{Tohoku};
see \cite{MacLane:Duality} and \cite{Buchsbaum:ExactCats} for earlier steps in that direction. %

Via (simplicial) model categories \cite{DGQuillen1967}, Quillen created a foundational framework for \emph{Homotopical Algebra}. Its vast scope encompasses the essentials of homological algebra, as well as homotopical invariants from algebraic topology and the homotopy theory of simplicial sets. Further examples include varieties of algebras/categories of models of Lawvere theories, next to what are today called \emph{Mal'tsev categories} (with enough projectives and coequalizers of kernel pairs), as all of these are simplicial Quillen model categories \cite[II.4]{DGQuillen1967}. %

\subsubsection*{Our objectives}
While wanting to develop a direct and effective interface between modern categorical algebra and Quillen's simplicial model category structure on Mal'tsev categories/varieties of algebras, we observed more and more advantages in building homology of chain complexes up systematically from categories which need only meet minimal structural prerequisites, namely the presence of a zero object along with the existence of kernels and cokernels. We call such a category \ZExact. The environment of a \ZExact\ category is sufficient for a discussion of universal factorizations of a given morphism and, consequently, of (short) exact sequences, of chain complexes, and homology as a measure for their failure to be exact. Remarkably, a version of the (Short) $5$-Lemma for normal maps holds in this context.---The basic properties of \ZExact\ categories are developed in Chapter~\ref{chap:PointedCats}.

To be able to compute effectively with homological invariants, one frequently relies on tools such as the Snake Lemma and the border cases of the $(\Prdct{3}{3})$-Lemma. Perhaps surprisingly, only one additional structural hypothesis\footnote{A form of this condition occurs in the definition of a \emph{w-exact category} due to Burgin~\cite{MSBurgin1970-Involution} and reappears as one of the `old axioms' in~\cite{Janelidze-Marki-Tholen}.} suffices to ensure that these tools are available: We assume that every composite
\begin{equation*}
  f=\DiagObj \overset{\mu}{\NMono} \DiagObj \overset{\varepsilon}{\NEpi} \DiagObj
\end{equation*}
in which $\mu$ is a kernel and $\varepsilon$ a cokernel admits a factorization as $f=\DiagObj \overset{e}{\NEpi} \DiagObj \overset{m}{\NMono}\DiagObj$, where $e$ is a cokernel and $m$ is a kernel. Calling a composite of a normal epimorphism with a normal monomorphism such as $m\Comp e$ a \emph{normal} map, and a composite of a normal monomorphism with a normal epimorphism such as $\varepsilon\Comp \mu$ an \emph{antinormal} map, this amounts to asking that \emph{every antinormal composite is a normal map}\footnote{In comparison: an abelian category is a \ZExact\ category with binary product and coproducts, in which \emph{every map is normal}.}. We refer to this as the \ANNInline-property. A~\ZExact\ category which also satisfies the \ANNInline-property is called \emph{di-exact}. ---The basic properties of di-exact categories are developed in Chapter \ref{chap:Di-ExactCats}. Moreover, we identify several structural axioms which are weaker than the \ANNInline-condition, and which still support homologically meaningful results.  We also explain connections with concepts such as p-exact categories.

\subsubsection*{Homological and semiabelian categories}
To explain the relationship with existing literature, the authors of \cite{FBorceuxDBourn2004} presented \emph{homological categories} as an environment in which the classical diagram lemmas of homological algebra hold: (Short) $5$-Lemma, $(\Prdct{3}{3})$-Lemma, Noether isomorphism theorems, Snake Lemma, which is the key to establishing the long exact homology sequence from a short exact sequence of normal chain complexes. We view a homological category as a category with zero object which is finitely bicomplete\footnote{This means that our definition of `homological category' differs slightly from earlier definitions such as the one in~\cite{FBorceuxDBourn2004}, where the only colimits required to exist are coequalizers of kernel pairs.} and satisfies two more structural axioms:
\begin{ulist}
  \item pullbacks preserve normal epimorphisms, and
  \item for every morphism $f\from X\to Y$ with section $s\from Y\to X$, the maps $s$ and $\KerMap{f}$ generate $X$; see Definition \ref{def:HomologicalCategory}.
\end{ulist}
Homological categories enjoy some features which are beyond the scope of di-exact categories. These are developed in Chapter \ref{chap:HomologicalCats}.

On the other hand, di-exact categories have qualities, such as Theorem~\ref{thm:(Co)Ker(ProperMapLESs)}, which are beyond the scope homological categories. This is thanks to the surprising strength of the \ANNInline-property. The synthesis of di-exact categories and homological categories is given by \emph{semiabelian categories}\footnote{The inventors of semiabelian categories chose this terminology because of the non-self-dual nature of this axiom system: It is half of abelian, in the sense that a category $\Ctgry{X}$ is abelian if and only if  both $\Ctgry{X}$ and its opposite $\Ctgry{X}^{\op}$ are semiabelian.}$^{,\thinspace}$\footnote{The term `semiabelian category' is also in use with two separate different meanings: (1) Palamodov calls a preabelian category \emph{semiabelian} if for every morphism  $f\from X\to Y$, the natural comparison map ${\CoKer{\Ker{f}} \to \Ker{\CoKer{f}}}$ is simultaneously a monomorphism and an epimorphism~\cite{VPPalamodov1968,VPPalamodov1971,Gruson,Rump0}. These are closely related to, but different from~\cite{Rump}, (2) \emph{Raïkov semiabelian} categories~\cite{Raikov:Semi-Abelian}, which coincide with the \emph{almost abelian} categories of Rump~\cite{Rump0} and the \emph{quasi-abelian categories} of~\cite{Yoneda-Exact-Sequences,Rump0,MR1779315}. As explained by G.~Janelidze~\cite{JRosickyWTholen2007}, the latter may be characterized as those categories which are at the same time homological and co-homological.}: A category $\Ctgry{X}$ is semiabelian exactly when it is homological and satisfies the \ANNInline-property.

With a view toward varieties of algebras, the semiabelian framework is particularly attractive for the following reason. In every homological category regular epimorphisms and normal epimorphisms coincide. In a variety of algebras  an epimorphism is regular if and only if its underlying set theoretic function is surjective, which makes it `easy' to recognize. ---Remarkably, every homological variety of algebras is automatically semiabelian. ---Foundational material of semiabelian categories is presented in Chapter \ref{chap:SACats}.

Still, \ZExact\ categories which are lower structured than `semiabelian' have their place: Thanks to Lawvere's semantic view of universal algebra, a Lawvere theory has models in many categories. For example, the category $\TopGrps$ of topological groups is homological but fails to be semiabelian.

In Part I of this monograph we work towards these convenient settings from the bottom up, starting with the context of pointed categories, in a way which is meant to keep the focus on the most basic homological-algebraic concepts---such as kernels and cokernels, (short) exact sequences, chain complexes and their homology, diagram lemmas---while assuming as little as possible in terms of conditions on the surrounding category. At the same time, this minimizes the amount of background knowledge we ask from the reader. The purpose of this methodology is to arrive at minimal requirements for our results to hold, requirements which moreover arise out of a clear homological-algebraic need. As a side-effect, the scope of the results in the first chapters encompasses part of the environment (which is not directly comparable to the semiabelian setting) considered by Grandis in~\cite{Grandis-HA2} for slightly different purposes.

\subsubsection*{Self-duality}
Along the way, we learned to pay attention to what we call the self-dual axis of a category. Let us sketch the underlying idea\footnote{This material is of independent interest. A more detailed treatment is in preparation.}: In the language of category theory, consider the collection of all self-dual statements about a category. Then the self-dual axis of a category $\Ctgry{X}$ is the collection of all those self-dual statements about $\Ctgry{X}$ which are true. Thus, we arrive at the tautological \dots

\emph{Metatheorem}\quad A category and its opposite always have the same self-dual axis.

This is so because a self-dual statement about $\Ctgry{X}$ is true if and only if it is true about $\Ctgry{X}^{\op}$.

For example, from the classical diagram lemmas we obtain self-dual statements about a \ZExact\ category. So, a \ZExact\ category $\Ctgry{X}$ contains those classical diagram statements in its self-dual axis if and only if they hold in $\Ctgry{X}$.

Immediately, this leads us to an application of the Metatheorem: It is known that $\SetsBsdOp$, the opposite of the category $\SetsBsd$ of pointed sets is a semiabelian category. So, the classical diagram lemmas belong to the self-dual axis of $\SetsBsdOp$ and, hence, to the self-dual axis of the category $\SetsBsd$ itself\ \dots\ even though $\SetsBsd$ is far from being semiabelian.

We organized the material in this part of the work so that we can relate ongoing developments more explicitly to the self-dual axis of a \ZExact\ category. This point of view naturally makes us arrive at the axioms of semiabelian categories from an unexpected new angle. Since the basic constructions of homological algebra are self-dual, we emphasize working with self-dual structural axioms---just as in the abelian framework.

\subsubsection*{Preliminaries}

We assume that the reader is familiar with basic category theoretical notions. In order to fix terminology and notation, we collect selected materials in an Appendix~\ref{chap:Categorical-Preliminaries}. For further background on category theory, we recommend \cite{SMacLane1998,Borceux:Cats,AHS:Cats}.

\newpage
\rhead{\bfseries\footnotesize
  Overview of exactness conditions}
\subsection*{An overview of exactness conditions and categorical structures}%
\label{sec:Examples-Overview}
\label{sec:SACatRecognize}

\begin{figure}[h!]
  \begin{equation*}
    \vcenter{
    \xymatrix@!@R=-6ex@C=3em{
    & \text{\color{OliveGreen} abelian} \ar[dl] \ar@[OliveGreen][d] \ar[dr] \\
    \text{\color{blue} semiabelian} \ar@[blue][d] \ar[dr] &
    \text{\color{OliveGreen} p-exact}  \ar@[OliveGreen][d] &
    \text{\color{red} co-semiabelian} \ar@[red][d] \ar[dl] \\
    \text{\color{blue} homological}  \ar@[blue][d] \ar[dr] &
    \text{\color{OliveGreen} di-exact}  \ar@[OliveGreen][d] &
    \text{\color{red} co-homological} \ar[dl] \ar@[red][d] \\
    \text{\color{blue} normal} \ar[dr] &
    \txt{\color{OliveGreen}(DPN)}  \ar@[OliveGreen][d] &
    \text{\color{red} co-normal} \ar[dl] \\
    & \txt{\color{OliveGreen}homologically\\\color{OliveGreen} self-dual}  \ar@[OliveGreen][d] \\
    & \text{\color{OliveGreen} z-exact}}
    }
  \end{equation*}
  \caption{Relations between categorical structures and exactness conditions. Categories on the right satisfy the dual set of structural axioms of those on the left and vice versa. Self dual structures appear in the center.}\label{fig:relations}
\end{figure}
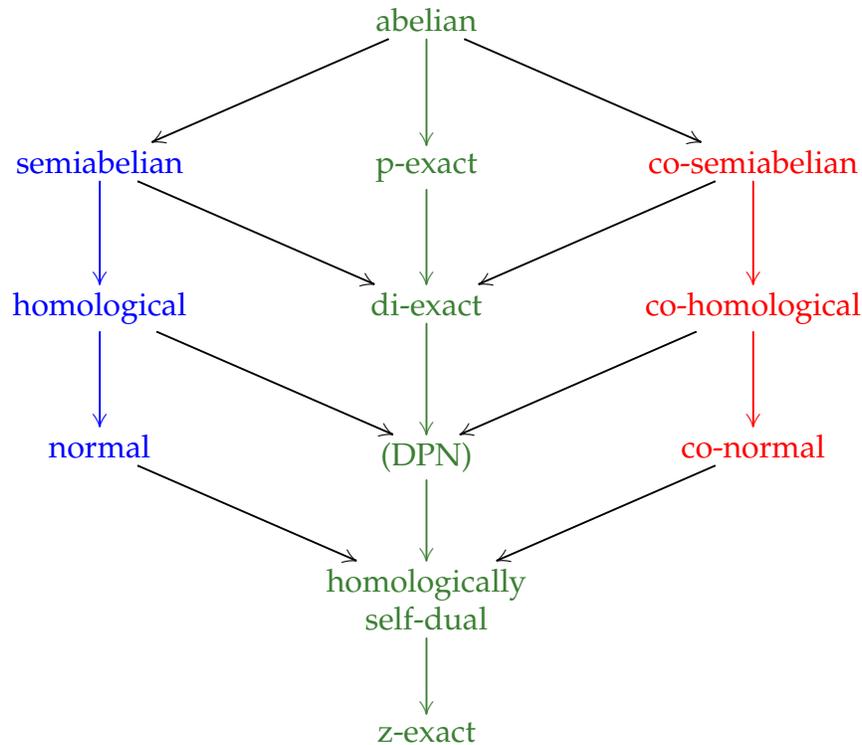

We will be working with various sets of structural axioms, each with certain underlying exactness conditions; see Figure~\ref{fig:relations} for a diagrammatical presentation of the relationships between those structures. Then we show how certain well-studied categories of algebraic structures fit into this diagram; see Tables~\ref{tab:Examples} and \ref{tab:Examples2} below.

The categorical structures in the left hand column are familiar from the literature, as are their duals on the right. Homologically relevant properties of so structured categories are systematically developed in Chapters~\ref{chap:NormalCategories}, \ref{chap:HomologicalCats} and \ref{chap:SACats}.

\newpage

We use Table~\ref{tab:Examples} to show how categories of familiar algebraic objects fit into the left hand column of Figure \ref{fig:relations}.

\begin{table}[bth!]
  \resizebox{\textwidth}{!}
  {\begin{tabular}{lccc}
      \toprule
      \textbf{semiabelian}                                        & \txt{varieties of $\Omega$-groups                                         \\
      Orzech categories of interest}                              & \txt{ $\Alg_{\TT}(\Sets)$, $\Alg_{\TT}(\HComp)$                           \\
      $\Grps$, $\Lie_R$, $\CCHopfK$ }                             & \txt{$\Loop$                                                              \\
      $\Sets^{\op}_{*}$, $\Cat^n(\Grps)$}                                                                                                     \\
      \midrule
      \textbf{homological}                                        & \txt{$\Grps(\Tops)$                                                       \\
      almost abelian categories                                 } & $\Alg_{\TT}(\Tops)$                             & $\AbGrpsTF$             \\
      \midrule
      \textbf{normal}                                             &                                                 &             & $\OrdGrp$ \\
      \bottomrule
    \end{tabular}}
  \caption{Rows separate algebraic structures: Referring to the left hand column of Figure~\ref{fig:relations}, entries in higher rows are also entries in lower rows, while entries of lower rows are in general \emph{not} entries of higher rows. }\label{tab:Examples}
\end{table}

The middle column in Figure~\ref{fig:relations} is reserved for self dual categorical structures. Amongst them abelian categories and p-exact categories are classically known; we elaborate on those in Chapter~\ref{chap:AbelianCategories}. We systematically introduce the remaining structures in Chapter~\ref{chap:Di-ExactCats}, and we explain their significance. Strongest among these conditions is di-exactness. It plays a special role for two reasons: (a) when satisfied, homology behaves really well; see e.g.\ Theorems \ref{thm:LES-Homology} and \ref{thm:(Co)Ker(ProperMapLESs)}; and (b) it separates homological categories from semiabelian ones.  A category $\Ctgry{X}$ is di-exact if and only if the following axioms hold.
\begin{ulist}
  \item  $\Ctgry{X}$ has a zero object; see (\ref{def:0-Object}).
  \item  Any morphism in $\Ctgry{X}$ has a (functorially chosen) kernel and a cokernel; see (\ref{sec:Kernel/CoKernel}).
  \item \emph{\ANNInline\ antinormal composites are normal:} whenever a morphism $f$ can be written as a composite $f=\epsilon\mu$, with a normal monomorphism $\mu$ followed by a normal epimorphism $\epsilon$, then $f$ may be factored as $f=me$, where $m$ is a normal monomorphism and $e$ is a normal epimorphism.
  \begin{equation*}
    \xymatrix@R=5ex@C=4em{
    \DiagObj \ar@{{ |>}->}[r]^-{\mu}   \ar@{.{ >>}}[d]_-{e}  &
    \DiagObj \ar@{-{ >>}}[d]^-{\epsilon} \\
    \DiagObj  \ar@{{ |>}.>}[r]_{m} &
    \DiagObj}
  \end{equation*}
\end{ulist}

The lowest structured categories we will be working with are the \ZExact\ ones. We only assume the existence of a zero object and that every map has a kernel and a cokernel. Thus \ZExact\ categories are self dual.

A \Defn{Puppe-exact} or \Defn{p-exact} category is a pointed category in which every morphism $f$ is normal. Thus $f$ admits an essentially unique decomposition $f=me$ with $e$ a normal epimorphism and $m$ a normal monomorphism. This implies the existence of kernels and cokernels, so that a p-exact category is always z-exact. For information on p-exact categories see~\cite{Grandis-HA1,Grandis-HA2,Borceux-Grandis,Alligators}, in particular the citations there to the foundational work of Puppe, Mitchell and others. A p-exact category which admits binary products or binary coproducts is abelian; see (\ref{thm:AdditivePlusp-exactIsAbelian}).

A di-exact additive category is abelian (\ref{thm:AdditivePlusDiexactIsAbelian}) because in an a pointed additive category, a morphism $f\from X\to Y$ may be factored as the normal monomorphism $\PrdctMapInto{\IdMapOn{X},0}\colon X\to X\oplus Y$ followed by the normal epimorphism $\MapOutOf{f, \IdMapOn{Y}}\colon X\oplus Y\to Y$. By the \ANNInline-condition, every morphism is normal. Hence the given category is both p-exact and additive, which makes it abelian.

We use Table~\ref{tab:Examples2} to show how categories of familiar algebraic objects fit into the middle column of Figure \ref{fig:relations}.

\begin{table}[h!]
  \resizebox{\textwidth}{!}
  {\begin{tabular}{lccc}
      \toprule
      \textbf{abelian}                 & $\LModules{R}$, $\RModules{R}$   & sheaves in $\LModules{R}$, $\RModules{R}$                  & $\AbGrps_{\mathit{fg}}$ \\
      \midrule
      \textbf{p-exact}                 &                                  & $\KProj$, $\Mlc$                                           & $\SetPBij$              \\
      \midrule
      \textbf{di-exact}                & \txt{semiabelian categories                                                                                             \\ w-exact categories}             &             & $\SetsBsd$, $\SetsBsdOp$ \\
      \midrule
      \textbf{homologically self-dual} & normal categories                & {$\SESCat{\Ctgry{X}}$ for $\Ctgry{X}$ di-exact}, $\CHopfK$ & $\CMon$,    $\TopsBsd$  \\
      \midrule
      \textbf{z-exact}                 & \txt{varieties of algebras whose                                                                                        \\ theory admits a single constant}                             &                                                 & $\SESCat{\CMon}$        \\
      \bottomrule
    \end{tabular}}
  \caption{Referring to the middle column of Figure~\ref{fig:relations}, entries in higher rows are also entries in lower rows, while entries of lower rows are in general \emph{not} entries of higher rows.}\label{tab:Examples2}
\end{table}

Here is how some well known and some less familiar categories of algebraic structures relate to the tables:
\begin{enumerate}
  \item Group objects $\Grps(\Tops)$ in the category $\Tops$ of topological spaces form a homological category which is \emph{not} semiabelian, as do models of a Lawvere theory $\TT$ in $\Tops$, provided its models in $\Sets$ form a semiabelian category.
  \item Group objects $\Grps(\Sets)$ in the category $\Sets$ of sets and functions form a semiabelian category, as do varieties whose objects have an underlying group structure and a unique constant: so-called \Defn{varieties of $\Omega$-groups}. Among the varieties of $\Omega$-groups we find:
        \begin{ulist}
          \item The category $\XMod$ of crossed modules.
          \item The category $\Rngs$ of (non-unital) rings and, more generally, associative and non-associative algebras over any given ring, such as the category $\Lie_R$ of Lie algebras over any ring $R$.
          \item Algebras (in $\Vect_\mathbb{K}$) over any algebraic operad.
        \end{ulist}
  \item All \emph{Orzech categories of interest}~\cite{Orzech} are semiabelian.
  \item Cocommutative Hopf algebras over a field form a semiabelian category~\cite{GSV,GKV}, while commutative Hopf algebras are at least homologically self-dual (see Example~\ref{exa:CocommHopf}; they are co-protomodular by~\cite{GM-VdL1}, but regularity is currently not clear).
  \item Torsion-free groups and torsion-free abelian groups form homological categories, denoted $\GrpsTF$ and $\AbGrpsTF$ respectively. Neither category is semiabelian. The variety $\Loop$ of loops is semiabelian but not a variety of $\Omega$-groups.
  \item A category is \emph{almost abelian}, \cite{Rump0}, exactly when it is at the same time homological and co-homological~\cite{JRosickyWTholen2007}. Among the almost abelian categories we find the categories of real or complex Banach spaces, and the category of locally compact abelian groups.

  \item The category $\SetsBsd$ of pointed sets is not homological because it fails to satisfy the \KSGInline-condition. However, it is co-semiabelian, see Exercise~\ref{exe:Set_*CoSemiabelian}. Consequently, $\SetsBsdOp$ is semiabelian, and both categories are di-exact. Neither category is p-exact because not every epimorphism is normal.
  \item For any $n\geq 1$, the category $\Cat^n(\Grps)$ of $n$-fold iterated internal categories in the category of groups (equivalent to Loday's $\cat^n$-groups~\cite{Loday}) is also semiabelian.
  \item The category $\AbGrps_{\mathit{fg}}$ of finitely generated \emph{abelian} groups is abelian.
  \item The category $\OrdGrp$ of preordered groups is known to be normal but not homological~\cite{MMClementinoNMartinsFerreiraAMontoli2019-Preordered}: a preordered group being a group equipped with a preorder (a reflexive and transitive relation) for which the group operation is monotone; arrows are monotone group homomorphisms.
  \item The category $\SetPBij$ of sets and partial bijections (bijections from a subset of the domain to a subset of the codomain) is known to be p-exact but not abelian~\cite{Grandis-HA2}.
  \item The category $\Mlc$ of modular lattices and modular connections is p-exact; see~\cite{Borceux-Grandis} and the references there.
  \item The category $\KProj$ of projective objects over a field $\mathbb{K}$, useful in $K$-theory, is one of the classical examples of a non-abelian p-exact category.
  \item \emph{Weakly p-exact} or \emph{w-exact} categories were introduced by Burgin in~\cite{MSBurgin1970-Involution}. They form a context encompassing both p-exact and semiabelian categories. For a definition and quick introduction in English, see~\cite{Borceux-Grandis}; note that \ANNInline\ is amongst the axioms defining w-exact categories.
  \item The category $\CMon$ of commutative monoids is homologically self-dual as shown in (\ref{thm:CMonIsHSD}). In (\ref{exa:HSD-not-DPN}) it is shown that this category is not \DPNInline, so that it cannot be di-exact.
  \item (\ref{exe:Top_*IsHSD}) explains that the category $\TopsBsd$ of pointed topological spaces is \HSDInline. Of course, the same holds for the dual of these last two categories.
  \item The category $\SESCat{\Ctgry{X}}$ of short exact sequences in any di-exact category $\Ctgry{X}$ is homologically self-dual (\ref{thm:ANN->HSDOf-NM(X),SES(X),NE(X)}). On the other hand~\ref{exa:NotHSD}, if we take short exact sequences in the category $\CMon$ of commutative monoids, then we find a category which is z-exact but not \HSDInline.
\end{enumerate}

\newpage

\begin{center}
  \textbf{Leitfaden for Part \ref{part:Homology}}
\end{center}

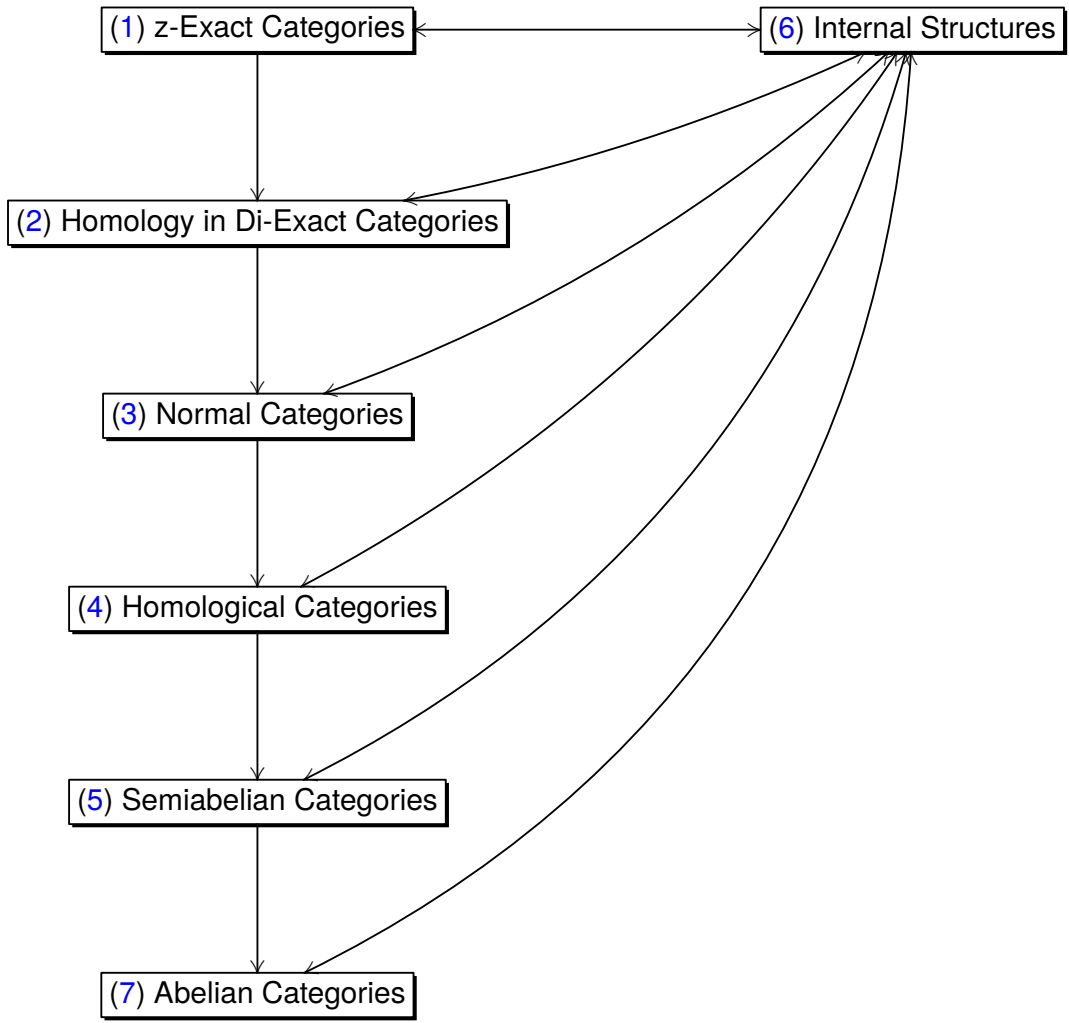
\bigskip

\begin{equation*}
  \xymatrix@R=10ex@C=8em{
  *+[F-,]{\txt{\sffamily (\ref{chap:ZExactCats}) z-Exact Categories}} \ar@{<->}[r] \ar[d] &
  *+[F-,]{\txt{\sffamily (\ref{chap:InternalStructures}) Internal Structures}} \ar@{<->}@/^2ex/[dl] \ar@{<->}@/^3ex/[ddl] \ar@{<->}@/^4ex/[dddl] \ar@{<->}@/^8ex/[ddddl] \ar@{<->}@/^12ex/[dddddl] \\
  *+[F-,]{\txt{\sffamily (\ref{chap:Di-ExactCats}) Homology in Di-Exact Categories}} \ar[d] \\
  *+[F-,]{\txt{\sffamily (\ref{chap:NormalCategories}) Normal Categories}} \ar[d] \\
  *+[F-,]{\txt{\sffamily (\ref{chap:HomologicalCats}) Homological Categories}} \ar[d] \\
  *+[F-,]{\txt{\sffamily (\ref{chap:SACats}) Semiabelian Categories}} \ar[d] \\
  *+[F-,]{\txt{\sffamily (\ref{chap:AbelianCategories}) Abelian Categories}}
  }
\end{equation*}
\fancyhead[LO]{\bfseries\footnotesize \thechapter\ \leftmark}
\fancyhead[LE]{\bfseries\footnotesize \thechapter\ \leftmark}
\fancyhead[C]{}
\fancyhead[RO]{\bfseries\footnotesize \thesection\ \rightmark}
\fancyhead[RO]{\bfseries\footnotesize \thesection\ \rightmark}

\chapter[z-Exact Categories]{z-Exact Categories}
\label{chap:ZExactCats}%
\label{chap:PointedCats}

A category is \ZExact\footnote{This terminology is specialized from Grandis in~\cite{Grandis-HA2}. He works with maps which assume the role of kernels and cokernels in a relativized setting: A class $\EuScript{Z}$ of objects which may contain non-zero objects is designated to act as zero-objects. Then `zero maps' are defined relative to $\EuScript{Z}$.} if it has a zero object, and if every morphism $f$ has a kernel and a cokernel which depends functorially on $f$. If a morphism $m$ represents the kernel of some map, then we call it a \Defn{normal monomorphism}. Dually, if a morphism $e$ represents the cokernel of some map, then we call it a \Defn{normal epimorphism}. %
\index{normal monomorphism}\index{normal epimorphism}%

In a \ZExact\ category every morphism $f$ has two universal factorizations:
\begin{ulist}
  \item The class of all compositions $f=mu$, with $m$ a normal monomorphism, is not empty and contains an initial element. We call it the \Defn{normal mono factorization} of $f$.
  \item The class of all compositions $f=ve$, with $e$ a normal epimorphism, is not empty and contains a terminal element. We call the \Defn{normal epi factorization } of $f$.
\end{ulist}
The factorizations come from the following categories associated to a \ZExact\ category~$\Ctgry{X}$:
\begin{enumerate}
  \item The category $\ArrowCat{\Ctgry{X}}$ of morphisms (arrows) in $\Ctgry{X}$.
  \item The category $\NMonoCat{X}$ of normal monomorphisms in $\Ctgry{X}$, which is contained in $\ArrowCat{\Ctgry{X}}$ as a full reflective subcategory.
  \item The category $\NEpiCat{X}$ of normal epimorphisms in $\Ctgry{X}$, which is contained in $\ArrowCat{\Ctgry{X}}$ as a full coreflective subcategory.
\end{enumerate}
On this foundation, we introduce the concept of exact sequence, as well as its dual: a coexact sequence.

We end this chapter by discussing internal magmas, monoids, and groups in a \ZExact\ category; see Sections \ref{sec:CatSESs-I} and \ref{sec:BiProducts}. We also explain how \ZExact\ categories are related to normal categories; see Section \ref{sec:NormalCats}.

\begin{center}
  \textbf{Leitfaden for Chapter \ref{chap:PointedCats}}
\end{center}

\bigskip

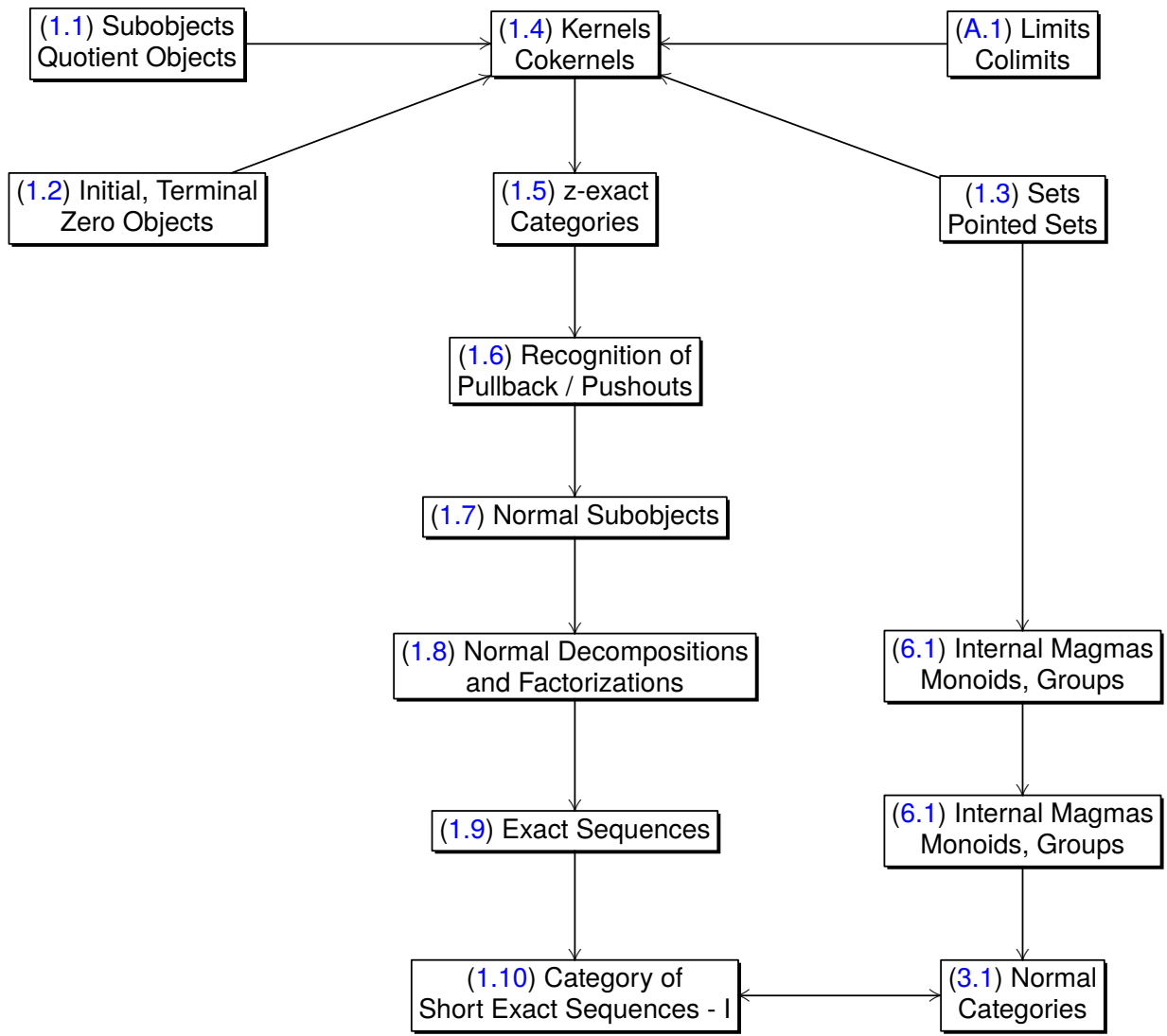
\begin{equation*}
  \resizebox{.98\textwidth}{!}{$
    \xymatrix@R=7ex@C=4.5em{
    *+[F-,]{\txt{\sffamily (\ref{sec:SubObjects-QuotientObjects}) Subobjects\\ \sffamily Quotient Objects}} \ar[r] &
    *+[F-,]{\txt{\sffamily (\ref{sec:Kernel/CoKernel}) Kernels \\ \sffamily Cokernels}} \ar[d] &
    *+[F-,]{\txt{\sffamily (\ref{sec:Limits-CoLimits}) Limits \\ \sffamily Colimits}} \ar[l] \\
    *+[F-,]{\txt{\sffamily (\ref{sec:Initial/Terminal/ZeroObjects}) Initial, Terminal \\ \sffamily Zero Objects}} \ar[ur] &
    *+[F-,]{\txt{\sffamily (\ref{sec:z-ExactCats}) \ZExact \\ \sffamily Categories}} \ar[d] &
    *+[F-,]{\txt{\sffamily (\ref{sec:PointedSets}) Sets \\ \sffamily Pointed Sets}} \ar[lu] \ar[ddd] \\
    & *+[F-,]{\txt{\sffamily (\ref{sec:PullbackPushout-Recognition}) Recognition of \\ \sffamily Pullback / Pushouts}} \ar[d] \\
    & *+[F-,]{\txt{\sffamily (\ref{sec:NormalSubobjects}) Normal Subobjects}} \ar[d] \\
    & *+[F-,]{\txt{\sffamily (\ref{sec:NEM-ImageFac}) Normal Decompositions \\ \sffamily and Factorizations}} \ar[d] &
    *+[F-,]{\txt{\sffamily (\ref{sec:InternalMagmas-Monoids-Groups}) Internal Magmas \\ \sffamily Monoids, Groups}} \ar[d] \\
    & *+[F-,]{\txt{\sffamily (\ref{sec:ExactSeqs}) Exact Sequences}} \ar[d] &
    *+[F-,]{\txt{\sffamily (\ref{sec:InternalMagmas-Monoids-Groups}) Internal Magmas \\ \sffamily Monoids, Groups}} \ar[d] \\
    & *+[F-,]{\txt{\sffamily (\ref{sec:CatSESs-I}) Category of  \\ \sffamily Short Exact Sequences - I}} \ar@{<->}[r]  &
    *+[F-,]{\txt{\sffamily (\ref{sec:NormalCats}) Normal \\ \sffamily Categories}} \\
    }
  $}
\end{equation*}
\section{Subobjects and Quotient Objects}
\label{sec:SubObjects-QuotientObjects}

We introduce the concepts of monomorphism, epimorphism, as well as the associated concepts of subobject and quotient object. These notions are valid in any category.

\subsection*{Monomorphisms}

\begin{definition}[Monomorphism]
  \label{def:Monomorphism}
  In any category $\Ctgry{X}$, a morphism $m\from M\to X$ is a \Defn{monomorphism}, a \Defn{mono} or  \Defn{is monic} if, whenever $u$, $v\from {U\to M}$ satisfy $mu=mv$, then $u=v$; notation: a tailed arrow $M\Mono X$. %
  \index{monic morphism}\index{monomorphism}\index[not]{m!$X\Mono Y$\IndSep monomorphism}
\end{definition}

\begin{definition}[Jointly monomorphic family]
  \label{def:JointlyMonomorphicFamily}
  A family of morphisms $(f_i\from M\to X_i)_{i\in I}$ in $\SACtgry{X}$ is \Defn{jointly monomorphic} if, whenever $u$, $v\from {U\to M}$ satisfy $f_iu=f_iv$ for all $i\in I$, then $u=v$. %
  \index{jointly!monomorphic family of maps}%
\end{definition}

For example, if $\lambda\from \LimOf{F}\Rightarrow F$ is a limit cone of a functor $F$ from a small category into $\SACtgry{X}$, then the family of morphisms $\lambda_j$ of the cone $\lambda$ is jointly monomorphic.

\subsection*{Subobjects}
In a category $\SACtgry{X}$ the class of monomorphism with fixed codomain $X$ carries a canonical partial ordering: We say that $m\from M\Mono X$ is \Defn{less than or equal to} $n\from N\Mono X$, in symbols $m\leq n$, if $m$ factors through $n$; i.e.\ there exists $p\from {M\to N}$, with $m=n\Comp p$. In this situation $p$ is unique, and it is monic by (\ref{thm:Monomorphism-Properties}.ii). We say that $m$ is \Defn{equivalent} to $n$ if $m\leq n$ and $n\leq m$. It follows that $p$ with $m=n\Comp p$ is an isomorphism. So, the class of monomorphisms with codomain $X$ carries a canonical equivalence relation.

\begin{definition}[Subobject]
  \label{def:Subobject}
  A \Defn{subobject} of an object $X$ in a category $\Ctgry{X}$ is an equivalence class of monomorphisms with codomain $X$. %
  \index{subobject}%
\end{definition}

It is common practice to refer to a morphism $m\from M\to X$ as the subobject it represents. If $m$ is determined by context, we may even refer to $M$ itself as a subobject of $X$, and write $M\leq X$. We write $\SubObjcts{X}$ for the (possibly large) set of all subobjects of $X$. It is equipped with the partial ordering induced by the $\leq$-relation on monomorphisms with codomain~$X$. %
\index[not]{s!$M\leq X$\IndSep $M$ is subobject of $X$}%
\index[not]{s!$\SubObjcts{X}$\IndSep class of subobjects of $X$}%

\begin{example}[Monomorphisms and subobjects of sets]
  \label{exa:Subobjects-Set}
  A morphism in the category $\Sets$ of sets is monic if and only if it is an injective function. Using the image of a function $f\from A\to X$, we see that every subobject of $X$ is uniquely represented by a subset of $X$. Consequently,  $\SubObjcts{X}$ is naturally equivalent to the power set $\PwrSt{X}$, equipped with the partial ordering by inclusion. %
  \index{subobject!in $\Sets$}\index[not]{p!$\PwrSt{X}$\IndSep power set of set $X$}
\end{example}

\begin{example}[Monomorphisms and subobjects of topological spaces]
  \label{exa:Subobjects-Top}
  In the category $\Tops$ of topological spaces, a monomorphism $m\from M\to X$ is given by a continuous function which is injective on underlying sets. Consequently, a subobject of $X$ may be constructed in two steps. In step 1, select a subset $S$ of the set underlying $X$, and equip $S$ with the subspace topology. In step 2 choose any refinement of the topology of $S$. For example, refining the topology of $X$ itself yields a space $X'$ for which the identity function $X'\to X$ of underlying sets is a continuous bijection, not necessarily a homeomorphism.  Thus, in $\Tops$, not all monomorphisms are subspace inclusions. For instance, the inclusion of the set of positive reals, equipped with the discrete topology, into $\RNr$ is such a monomorphism. %
  \index{subobject!in $\Tops$}\index[not]{t!$\Tops$\IndSep category of topological spaces}
\end{example}

\begin{example}[Monomorphisms and subobjects of compact Hausdorff spaces]
  \label{exa:Subobjects-CHTop}%
  In the category $\CHTops$ of compact Hausdorff spaces, a monomorphism $m\from M\to X$ is given by an embedding of a compact Hausdorff space $M$ as a closed subspace of $X$. Consequently, a subobject of $X$ may be uniquely represented by the inclusion mapping of a closed subspace of $X$.

  Indeed, any proper refinement of the topology of a compact Hausdorff space yields a Hausdorff space which is no longer compact. This is so because in a compact Hausdorff space, a subset is closed if and only if it is compact. %
  \index{subobject!in $\Tops$}\index[not]{c!$\CHTops$\IndSep category of compact Hausdorff spaces}
\end{example}

\begin{example}[Monomorphisms and subobjects in a variety of algebras]
  In a variety of algebras, a monomorphism $m\from M\to X$ is given by a morphism of algebras which is injective on underlying sets. Thus the subobjects of an algebra correspond to its subalgebras. %
\end{example}

\subsection*{Intersection and union of subobjects}
\label{subsec:Intersection/Union-Subobject}

In a category $\Ctgry{X}$, we introduce the notions of `meet' and `join' of subobjects. These are categorical versions of `intersection' and `union' of subsets.

\begin{definition}[Meet and join of subobjects]
  \label{def:Meet/Join-Subobjects}%
  Given an object $X$ in a category $\Ctgry{X}$, consider two subobjects, respectively represented by monomorphisms $m\from {M\to X}$ and $n\from {N\to X}$.
  \begin{thmlist}
    \item The \Defn{meet} or \Defn{intersection} of $m$ and $n$ is their greatest lower bound in the partial ordering of $\SubObjcts{X}$, provided it exists; notation: $M\meet N$ or, more precisely, $m\meet n$. %
    \index{intersection!of subobjects}\index{meet!of subobjects: definition}\index{subobjects!meet}\index{subobjects!intersections}%
    \index[not]{i!$M\meet N$\IndSep intersection/meet of $M$ and $N$}%
    \index[not]{m!$M\meet N$\IndSep intersection/meet of $M$ and $N$}%
    \item The \Defn{join} or \Defn{union} of $m$ and $n$ of $X$ is their least upper bound in $\mathit{Sub}(X)$, provided it exists; notation $M\join N$. %
    \index{union!of subobjects}\index{join!of subobjects: definition}\index{subobjects!union}\index{subobjects!join}%
    \index[not]{j!$M\join N$\IndSep join/union of subobjects $M$ and $N$}%
    \index[not]{u!$M\join N$\IndSep join/union of subobjects $M$ and $N$}%
  \end{thmlist}
\end{definition}

The intersection of $M$ and $N$ always exists whenever $\Ctgry{X}$ admits pullbacks; see (\ref{exe:IntersectionConstruct-Via-Pullback}). The join $M\join N$ does exist whenever $\Ctgry{X}$ admits binary coproducts and image factorizations; see (\ref{thm:JoinMorphisms-Construction}). If $\SACtgry{X}$ admits arbitrary limits, we have:

\begin{proposition}[Arbitrary meets and joins in a complete category]
  \label{thm:Meet/Join-InCompleteCat}
  For a family $\mathcal{S}=(S_{\lambda} | \lambda\in \Lambda)$ of subobjects of $X$ in a complete category $\Ctgry{X}$ the following hold:
  \begin{thmlist}
    \item The meet of the objects in $\mathcal{S}$ exists and may be constructed as %
    \index{meet!of subobjects: construction}
    \begin{equation*}
      \bigmeet_{\lambda\in \Lambda} S_{\lambda} = \LimOfOver{S_{\lambda}\Mono X}{\Lambda}
    \end{equation*}
    \item If $\SACtgry{X}$ is well powered, then the join of $\mathcal{S}$ exists and may be constructed as %
    \index{join!of subobjects: construction}
    \begin{equation*}
      \bigjoin_{\lambda\in \Lambda} S_{\lambda} = \bigmeet_{i\in I} T_i
    \end{equation*}
    with $(T_i | i\in I)$ the family of all those subobjects of $X$ which contain every $S_{\lambda}$. \NoProof
  \end{thmlist}
\end{proposition}

\subsection*{Epimorphisms}

Dual to the concepts of monomorphism and subobject are the concepts of epimorphism and quotient object:

\begin{definition}[Epimorphism]
  \label{def:Epimorphism}%
  In any category, a morphism $e\from Y\to Z$ is an \Defn{epimorphism}, an \Defn{epi} or \Defn{is epic} if, whenever $u$, $v\from Z\to W$ satisfy $ue=ve$, then $u=v$; notation: a double-headed arrow $Y\Epi Z$.%
  \index{epic morphism}\index{epimorphism}\index[not]{e!$Y\Epi Z$\IndSep epimorphism}%
\end{definition}

We call a family $(g_i\from A_i\to Y)_{i\in I}$ \Defn{jointly epimorphic} if, whenever $u$, $v\from Y\to Z$ satisfy $ug_i=vg_i$ for all $i\in I$, then $u=v$. For example, if $\gamma\from F\Rightarrow \CoLimOf{F}$ is a colimit cocone of a functor $F$ from a small category into $\SACtgry{X}$, then the family of morphisms $\gamma_j$ of the cocone $\gamma$ is jointly epimorphic.%
\index{jointly!epimorphic family of maps}

\subsection*{Quotient objects}

In a category $\SACtgry{X}$, an epimorphism $e\from X\to Q$ is said to be \Defn{greater than or equal to} an epimorphism $\varepsilon\from X\to R$ if there exists $q\from Q\to R$ with $qe=\varepsilon$. Whenever such $q$ exists, it is unique and, by (\ref{exe:Epimorphisms-Composite}.ii) epic. If $e$ is greater than or equal to $\varepsilon$ and $\varepsilon$ is greater than or equal to $e$, we say that $e$ and $\varepsilon$ are \Defn{equivalent}. %
\index{equivalent!monomorphisms}

A quotient object of $X$ is an equivalence class of epimorphisms with domain $X$. It is common practice to refer to the epimorphism $e\from X\to Q$ or the object $Q$ as the quotient object it represents. %
\index{quotient!object}

Unexpected quotient objects can occur. For example, in the category of unital rings the field of rationals numbers $\QNr$ appears as a quotient object of $\ZNr\to \QNr$; see (\ref{exa:Z>->Q}). We conclude that the concept of epimorphism does not characterize surjective algebra morphisms. For this reason, quotients determined by mere epimorphisms may exhibit pathological properties. Such pathologies do not occur in quotients determined by stronger types of epimorphisms, such as regular epimorphisms in the context of a variety of algebras, or normal epimorphisms as in Section~\ref{sec:ImageFactorizations}.%

\begin{exercises}

\begin{exercise}[Composite of monomorphisms]
  \label{exe:Monomorphisms-Composite}
  \label{thm:Monomorphism-Properties}
  \label{thm:Monomorphism-Composites}
  In an arbitrary category $\SACtgry{X}$, consider composable morphisms $X\XRA{f} Y \XRA{g} Z$, and show the following:%
  \index{monomorphism!properties}%
  \begin{thmlist}
    \item If both $f$ and $g$ are monomorphisms, then $gf$ is a monomorphism.
    \item If $gf$ is a monomorphism, then $f$ is a monomorphism.
    \item If $gf$ is a monomorphisms, then $g$ need not be a monomorphism.
  \end{thmlist}
\end{exercise}

\begin{exercise}[Monomorphisms in $\CHTops$]
  \label{exe:Monomorphis-CHTop}
  Show that in the category $\CHTops$ of compact Hausdorff spaces every monomorphism $m\from M\to X$ is uniquely represented by the inclusion of a closed subspace $S$ of $X$ into $X$.
\end{exercise}

\begin{exercise}[Composite of epimorphisms]
  \label{exe:Epimorphisms-Composite}
  In a category $\SACtgry{X}$, consider composable morphisms $X\XRA{f} Y \XRA{g} Z$, and show the following:%
  \index{epimorphism!properties}%
  \begin{thmlist}
    \item If both $f$ and $g$ are epimorphisms, then $gf$ is an epimorphism.
    \item If $gf$ is an epimorphism, then $g$ is an epimorphism.
    \item If $gf$ is an epimorphisms, then $f$ need not be an epimorphism.
  \end{thmlist}
\end{exercise}

\begin{exercise}[Epimorphisms in $\Grps$]
  \label{exe:Epimorphisms-In-Grps}
  \cite[p.~21]{SMacLane1998}\quad Show that in the category $\Grps$ of groups every epimorphism is surjective.
\end{exercise}

\begin{exercise}[Limit construction of join]
  \label{exe:Join-LimitConstruction}
  Let $\Ctgry{C}$ be a complete category in which for each object $X$ the collection of subobjects is a set; i.e.\ $\Ctgry{C}$ is \Defn{well powered}. If $M$ and $N$ are subobjects of $X$, show that
  \begin{equation*}
    M\join N = \LimOfOver{\Phi}{J}
  \end{equation*}
  where $J$ is the category spanned by all those representatives of subobjects of $X$ which contain both $M$ and $N$. This makes $J$ a subcategory of $\mathcal{C}$, and $\Phi\from J\to \mathcal{C}$ is the inclusion functor.
\end{exercise}
\end{exercises}
\section[Initial, Terminal, Zero Objects]{Initial, Terminal, Zero Objects}
\label{sec:Initial/Terminal/ZeroObjects}
\label{sec:PointedCategories}
\label{sec:Kernels/Cokernels} 

We introduce the concepts of initial, terminal, and zero objects, and identify such objects in a sampling of categories.

\begin{definition}[Initial, terminal, zero object]
  \label{def:0-Object}
  An object $i$ in a category $\Ctgry{X}$ is called \Defn{initial} if $\HomIn{\Ctgry{X}}{i}{X}$ is a singleton set for all objects $X$ in~$\Ctgry{X}$. Dually, an object $t$ in $\Ctgry{X}$ is \Defn{terminal} if $\HomIn{\Ctgry{X}}{X}{t}$ is a singleton set for all  $X$. A  \Defn{zero object} in $\Ctgry{X}$ is an object $\ZeroObject$ which is both initial and terminal. %
  \index{zero!object}\index{terminal!object}\index{initial!object}\index{object!initial}\index{object!terminal}\index{object!zero}%
\end{definition}

\begin{definition}[Pointed category]
  \label{def:PointedCat}%
  We say that a category $\Ctgry{X}$ is \Defn{pointed} if it has a zero object. %
  \index{pointed!category}\index{category!pointed}%
\end{definition}

\begin{example}[Terminal / initial object in $\Sets$]
  \label{exa:Set-Terminal/Initial-Object}%
  The category $\Sets$ of sets has a unique initial object, namely the empty set $\emptyset$. It also has terminal objects, namely any $1$-element set. As $1$-element sets are non-empty, $\Sets$ does not have a zero object.
\end{example}

In a category $\Ctgry{X}$ with zero object every hom-set $\HomIn{\Ctgry{X}}{X}{\ZeroObject}$ contains a unique element, denoted $\ZeroMap_{\ZeroObject X}$, and $\HomIn{\Ctgry{X}}{\ZeroObject}{Y}$ contains a unique element, denoted $\ZeroMap_{Y\ZeroObject}$. Consequently, $\HomIn{\Ctgry{X}}{ X }{ Y }$ contains the composite $\ZeroMap_{YX}\DefEq \ZeroMap_{Y\ZeroObject}\Comp \ZeroMap_{\ZeroObject X}$, called the \Defn{zero map} from $X$ to $Y$. It is the unique map from $X$ to $Y$ which factors through $\ZeroObject$. We write $\ZeroMap$ in place of $\ZeroMap_{YX}$ if there is no risk of confusion. An object in $\Ctgry{X}$ is a zero object if and only if $\IdMapOn{X}=\ZeroMap_{XX}$; see (\ref{exe:0-Object-Recognition}).
\index[not]{z!$\ZeroMap_{YX}$\IndSep zero morphism from $X$ into $Y$}\index[not]{z!$\ZeroObject$\IndSep zero object}%
\index{zero!morphism}\index{morphism!zero}%

Composition of morphisms preserves the zero map in the following sense:

\begin{proposition}[Zero maps are absorbant under composition]
  \label{thm:ZeroMapPreserved-Composition}
  In category with zero object the following hold for any morphisms $u\from A\to X$ and ${v\from Y\to Z}$:
  \begin{equation*}
    \ZeroMap_{0X}\Comp u = \ZeroMap_{0A} \qquad\text{and}\qquad v\Comp \ZeroMap_{Y0} = \ZeroMap_{Z0}  \NoProofDiag
  \end{equation*}
\end{proposition}

Here are examples of categories with a zero object.

\begin{example}[The category of sets with base point is pointed]
  \label{exa:Set_*PointedCat}%
  Associated to $\Sets$ is the category $\SetsBsd$ of sets with a chosen base point, called `pointed sets' for short: An object in $\SetsBsd$ is a pair $(X,x)$ consisting of a set $X$ and an element $x\in X$, called the base point of $(X,x)$. A morphism $f\from (X,x)\to (Y,y)$ is a function $f\from X\to Y$ such that $f(x)=y$. It follows that any pair $(\Set{x},x)$ is a zero object in $\SetsBsd$. %
  \index[not]{s!$\SetsBsd$\IndSep category of sets with base point}%
  \NoProof
\end{example}

\begin{example}[$\Grps$ has is pointed]
  \label{exa:Groups-PointedCat}%
  The category $\Grps$ of groups has a zero object: any one-element group is one. The zero maps are those which send every element of the domain to the neutral element of the codomain. %
  \index[not]{g!$\Grps$\IndSep category of groups}%
  \NoProof
\end{example}

\begin{example}[$\Rngs$ is pointed]
  \label{exa:UnitaryRings-Pointedness}
  \label{exa:Rings-Pointed}%
  In the category $\Rngs$ of rings with or without a unit we allow that the multiplicative unit $1$ equals the additive $0$. A morphism in $\Rngs$ preserves the sum and the multiplication, but need not preserve the unit (if a unit exists). The zero ring is a zero object. %
  \index[not]{r!$\Rngs$\IndSep category of rings with or without unit} \NoProof
\end{example}

\begin{example}[$\URngs$ has distinct initial and terminal objects]
  \label{exa:UnitaryRings-NotPointed}
  In contrast, let $\URngs$ be the category whose objects are rings with a multiplicative $1$, called the unit; again, we do \emph{not exclude} the possibility that $1$ equals the additive $0$. Morphisms are set theoretic functions which preserve the unit and commute with addition and multiplication. Then $\URngs$ has the integers $\ZNr$ as an initial object, and it has the zero ring, consisting of a single element `$0$', as a terminal object. As $\ZNr\neq 0$, $\URngs$ does not have a zero object. %
  \index[not]{r!$\URngs$\IndSep category of unital rings}%
  \index[not]{z!$\ZNr$\IndSep number system of integers}
\end{example}

\begin{subordinate}{Special role of the category $\SetsBsd$}
  Whenever we work with pointed categories, the category $\SetsBsd$ of pointed sets is implicitly involved. This is so because every hom-set is canonically pointed by the zero map, and zero maps are absorbant under composition (\ref{thm:ZeroMapPreserved-Composition}). So, pointed categories are enriched in $\SetsBsd$.

  We note that the concept of a pointed category is self-dual. %
\end{subordinate}

\begin{exercises}
\begin{exercise}[Initial object as a colimit]
  \label{exe:InitialObject=CoLim(empty)}
  In any category $\EuScript{X}$, show that an object $A$ is an initial object in $\EuScript{X}$ if and only if it is a colimit of the empty functor from the empty category into $\EuScript{X}$.
\end{exercise}

\begin{exercise}[Terminal object as a limit]
  \label{exe:TerminalObject=Lim(empty)}
  In any category $\EuScript{X}$, show that an object $Z$ is a terminal object in $\EuScript{X}$ if and only if it is a limit of the empty functor from the empty category into $\EuScript{X}$.
\end{exercise}

\begin{exercise}[Recognizing a zero object\ZExactTag]
  \label{exe:0-Object-Recognition}
  In a category with a zero object $0$, show that the following conditions are equivalent:
  \begin{tfae}
    \item $X$ is a zero object;
    \item there exists a monomorphism $X\to 0$;
    \item there exists an epimorphism $0\to X$;
    \item $\IdMapOn{X}=\ZeroMap_{XX}$.
  \end{tfae}
\end{exercise}

\begin{exercise}[Sum projections yield pushout]
  \label{exe:SumProjectionsYieldPushout}
  In a pointed category with sums use arbitrary objects $X$ and $Y$ to form this commutative square.
  \begin{equation*}
    \xymatrix@R=5ex@C=4em{
    X+Y \ar[r]^-{\SumMapOutOf{0,\IdMapOn{Y}}} \ar[d]_{\SumMapOutOf{\IdMapOn{X},0}} &
    Y \ar[d] \\
    X \ar[r] &
    0
    }
  \end{equation*}
  Show that this square is a pushout.
\end{exercise}

\begin{exercise}[Product via pullback]
  \label{exe:ProductViaPullback}
  In a pointed category with pullbacks show that the universal object in the pullback diagram below is a product of $X$ and $Y$. %
  \index{product!as pullback}\index{pullback!construction of product}
  \begin{equation*}
    \xymatrix@R=5ex@C=4em{
    P \ar[r]^-{x} \ar[d]_{y} \PullLU{rd} &
    X \ar[d]^{\ZeroMap} \\
    Y \ar[r]_-{\ZeroMap} &
    \ZeroObject
    }
  \end{equation*}
\end{exercise}

\begin{exercise}[Alternative construction of the category of pointed sets]
  The category of pointed sets $\SetsBsd$ may be constructed as follows. If $\ast\in \Sets$ is any terminal object, i.e.\ any $1$-element set, then the category $\ast/\Sets$ of \Defn{object under $\ast$}---also called \Defn{slice category} or \Defn{comma category}, an object of $\ast/\Sets$ is a morphism $x\from \ast\to X$, and a morphism from $x\from \ast\to X$ to $y\from \ast\to Y$ is a function $f\from X\to Y$ such that $f\Comp x=y$---is isomorphic to $\SetsBsd$.
\end{exercise}
\end{exercises}
\section[Sets and Pointed Sets]{Sets and Pointed Sets}
\label{sec:PointedSets}

Here we present some basic properties of the categories $\Sets$ of sets and $\SetsBsd$ of pointed sets. We explain in detail how functorial limits and colimits can be constructed in $\Sets$. Then we use the forgetful functor $U\from \SetsBsd\to \Sets$ and its left adjoint, given by adding a disjoint point `$+$' to a given set $X$\*, to obtain functorial limits and colimits in $\SetsBsd$; see (\ref{thm:Sets_*/Sets-Adjunction}),   (\ref{thm:Set_*Complete}), and (\ref{thm:Set_*CoComplete}). %
\index[not]{s!$\Sets$\IndSep category of sets}\index[not]{s!$\SetsBsd$\IndSep category of sets with base point}

Any product over the empty category is the set $\Set{\emptyset}$. If $J$ is a discrete $1$-element category, then the limit of $F\from J\to \Sets$ is $F(j)$. The product of two sets $X$ and $Y$ is the set of all ordered pairs $(x,y)\DefEq \Set{x,\Set{x,y}}$. With binary products, the concept of function $f\from X\to Y$ is defined as a subset of $\Prdct{X}{Y}$. If $F\from J\to \Sets$ is a functor, where $J$ is a discrete category with $3$ or more objects, then
\begin{equation*}
  \LimOf{F} = \FamPrdct{j\in J}{F(j)} \DefEq \Set{f\from X\to Y\ |\ f\ \text{is function}}
\end{equation*}
For $n\geq 0$, we construct the $n$-th power $P_n\from \Sets\to \Sets$ by letting $P_n(X)$ be the product  determined by the constant functor $p\from \Set{1,\dots ,n}\to \Sets$ with value $X$.

If $J$ is an arbitrary small category, then the limit $F\from J\to \Sets$ is obtained as
\begin{equation*}
  \LimOf{F} \DefEq \bigl\{x \in \FamPrdct{j\in J}{F(j)}\mid \text{$\forall\ u\from i\to j$ in $J$, $F_{u}(x_i) = x_j$}\bigr\}
\end{equation*}
Given a pointed category $\Ctgry{X}$, we observe that the Hom-functor on $\Ctgry{X}$ factors through the category $\SetsBsd$:
\begin{equation*}
  \xymatrix@R=6ex@C=3em{
  \Prdct{\Ctgry{X}^{\op}}{\Ctgry{X}} \ar[dr]_{\HomBsd{-}{-}} \ar[rr]^-{\HomIn{\Ctgry{X}}{-}{-}} &&
  \Sets \\
  & \SetsBsd \ar[ru]_{U}
  }
\end{equation*}
We take the existence of such a factorization as a light way of saying that $\Ctgry{X}$ is enriched in the category of pointed sets.

Turning to colimits, any coproduct over the empty category is the empty set. If $J$ is a discrete $1$-element category, then the colimit of $F\from J\to \Sets$ is $F(j)$. If $J$ is a discrete category with $2$ or more elements, then
\begin{equation*}
  \FamCoPrdct{j\in J}{F(j)} \DefEq \Bigl\{ (j,x)\in \Objcts{J}\prdct \bigl( \FamUnion{j\in J}{F(j)}\bigr)\mid x\in F(j)\Bigr\}
\end{equation*}
If $J$ is an arbitrary small category, then the colimit $F\from J\to \Sets$ is defined as the quotient set
\begin{equation*}
  \CoLimOf{F} \DefEq \left. \FamCoPrdct{j\in J}{F(j)} \right/ \sim
\end{equation*}
where `$\sim$' is the equivalence relation generated by $(i,x)\sim (j,F_{u}x)$, for all morphisms $u\from i\to j$ in $J$.

\begin{theorem}[$\Sets$ is closed symmetric monoidal]
  \label{thm:Set-ClosedSymmetricMonoidal}%
  The category $\Sets$  is symmetric monoidal via `product'. The commutative diagram shows that $\FamCoPrdct{a\in -}{B} \cong \Prdct{-}{B}$  is left adjoint to  $\Hom{B}{-}$.
  \begin{equation*}
    \xymatrix@R=5ex@C=3em{
    \Hom{\Prdct{A}{B}}{C} \ar@{<->}[r] \ar@{<->}[d]  &
    \Hom{A}{\Hom{B}{C}} \ar@{<->}[d] \\
    \Hom{\FamCoPrdct{a\in A}{B}}{C} \ar@{<->}[r] &
    \FamPrdct{a\in A}{ \Hom{B}{C} }
    }
  \end{equation*}
\end{theorem}

Let us now turn to the relationship between the categories of sets and pointed sets.

\begin{proposition}[$\SetsBsd \leftrightarrows \Sets$ adjunction]
  \label{thm:Sets_*/Sets-Adjunction}
  The base point forget functor $U\from \SetsBsd\to \Sets$ has a left adjoint $X\mapsto X_+$ which sends a set $X$ to the pointed set $X_+\DefEq(X+ \Set{X},X)$, which is the set $X$ with the base-point $X$ adjoined via coproduct (= disjoint union) in $\Sets$. \NoProof
\end{proposition}

\begin{theorem}[$\SetsBsd$ is complete]
  \label{thm:Set_*Complete}
  The forgetful functor $U\from \SetsBsd\to \Sets$ preserves and reflects limits. --- If $(X_{i},x_{i})_{i\in I}$ is a family of pointed sets, then
  \begin{equation*}
    \FamPrdct{i\in I}{(X_i,x_i)} = \Bigl( \FamPrdct{i\in I}{X_i}, (x_i)_{i\in I}\Bigr)
  \end{equation*}
  Subsequently, the limit of a diagram of pointed sets is canonically a pointed subset of the product of pointed sets at the nodes of the diagram. In particular, $\SetsBsd$ is functorially complete. \NoProof
\end{theorem}

\begin{definition}[Wedge operation in $\SetsBsd$]
  \label{def:Wedge-Set_*}
  In $\SetsBsd$ the \Defn{wedge} or \Defn{$1$-point union} of a family $(X_{k},x_{k})_{k\in K}$ of pointed sets is %
  \index{wedge in $\SetsBsd$}\index{$1$-point union of pointed sets}%
  \index[not]{w!$\bigvee_{k\in K} (X_k,x_k)$\IndSep wedge/coproduct in $\SetsBsd$}%
  \index[not]{c!$\bigvee_{k\in K} (X_k,x_k)$\IndSep wedge/coproduct in $\SetsBsd$}%
  \begin{equation*}
    \bigvee_{k\in K} (X_k,x_k) \DefEq \Bigl( \left. \FamCoPrdct{k\in K}{X_k}\right/ (x_{k_1}\sim x_{k_2})\; , \ast \Bigr)
  \end{equation*}
  with $\ast$ the equivalence class of any $x_k$.
\end{definition}

Let $q\from \FamCoPrdct{k\in K}{X_k}\to \bigvee_{k\in K}(X_k,x_k)$ be the quotient map. If  $\InclsnOf{k}\from X_k\to \FamCoPrdct{k\in K}{X_k}$ are the structure maps for the coproduct in $\Sets$, then $q\InclsnOf{k}$ are structure maps for the coproduct in $\SetsBsd$. As the coproduct in $\Sets$ is associative and commutative, so is the wedge operation in $\SetsBsd$.

\begin{theorem}[$\SetsBsd$ is cocomplete]
  \label{thm:Set_*CoComplete}
  The category $\SetsBsd$ is functorially cocomplete. The colimit of a diagram in $\SetsBsd$ is given by the appropriate quotient of the wedge formed from the pointed sets at the nodes of the diagram. %
  \index{colimit!in $\SetsBsd$}%
  \NoProof
\end{theorem}

\begin{definition}[Smash product in $\SetsBsd$]
  \label{def:SmashInSet_*}
  In $\SetsBsd$, the \Defn{smash (product)} of $(X,x_0)$ by $(Y,y_0)$ is %
  \index{smash product!in $\SetsBsd$}
  \begin{equation*}
    (X,x_0)\wedge (Y,y_0) \DefEq \left( (\Prdct{X}{Y})/((\Prdct{X}{\Set{y_0}})\union (\Prdct{\Set{x_0}}{Y}))\, ,\, \Set{(\Prdct{X}{\Set{y_0}})\union (\Prdct{\Set{x_0}}{Y})} \right)
  \end{equation*}
  We write $x\wedge y$ for the image of $(x,y)\in \Prdct{X}{Y}$ in $(X,x_0)\wedge (Y,y_0)$. %
  \index[not]{s!$(X,x_0)\wedge (Y,y_0)$\IndSep smash of pointed sets}
\end{definition}

Thus the smash Of $(X,x)$ by $(Y,y)$ is constructed by collapsing the subset $(\Prdct{X}{\Set{y}})\union (\Prdct{\Set{x}}{Y})$ of $\Prdct{X}{Y}$ to a point, which becomes the base point of $(X,x)\wedge (Y,y)$.

\begin{theorem}[$\SetsBsd$ with $-\wedge-$ is closed symmetric monoidal]
  \label{thm:Sets_*ClosedSymmetricMonoidal}
  The category $\SetsBsd$  is closed symmetric monoidal via `smash product'. A natural bijection of based sets %
  \index{closed symmetric monoidal!structure in $\SetsBsd$}
  \begin{equation*}
    \HomBsd{(X,x_0)\wedge (Y,y_0)}{(Z,z_0)} \longleftrightarrow \HomBsd{(X,x_0)}{\HomBsd{(Y,y_0)}{(Z,z_0)}}
  \end{equation*}
  is given by
  \begin{thmlist}
    \item $f\from (X,x_0)\wedge (Y,y_0)\to (Z,z_0)$ gets sent to
    \begin{equation*}
      \bar{f}\from (X,x_0)\to \HomBsd{(Y,y_0)}{(Z,z_0)},\qquad [\bar{f}(x)](y) \DefEq f(x\wedge y)
    \end{equation*}
    \item $g\from (X,x_0)\to \HomBsd{(Y,y_0)}{(Z,z_0)}$ gets sent to
    \begin{equation*}
      \underline{g}\from (X,x_0)\wedge (Y,y_0)\to (Z,z_0),\qquad \underline{g}(x\wedge y)\DefEq [g(x)](y)
    \end{equation*}
  \end{thmlist}
\end{theorem}

Now consider a pointed category $\Ctgry{X}$. Every Hom-set $\HomIn{\Ctgry{X}}{X}{Y}$ in $\Ctgry{X}$ is canonically equipped with a base point, namely the zero-map $\ZeroMap_{YX}$. Let us write $\HomBsd{X}{Y}$ for this pointed Hom-set.

\begin{theorem}[Pointed Hom-sets of a pointed category]
  \label{thm:PointedCat-EnrichedInSet_*-Light}
  In a pointed category $\Ctgry{C}$ the following hold.
  \begin{enumerate}[(i)]
    \item $\HomIn{\Ctgry{X}}{X}{Y}=U\Comp \HomBsd{X}{Y}$, where $U\from \SetsBsd\to \Sets$ is the forgetful functor.
    \item For $u\from X\to Y$ and $v\from Y\to Z$, $\ZeroMap_{ZY}u = \ZeroMap_{ZX} = v\ZeroMap_{YX}$.
  \end{enumerate}
\end{theorem}

Thus, in particular $\SetsBsd$ is enriched in itself. In the following lemma we confirm that a pointed category can only be enriched in $\SetsBsd$ by taking the zero map $\ZeroMap_{XY}\from X\to Y$ as the base point of $\Hom{X}{Y}$.

\begin{lemma}[In pointed category: base point in $\Hom{X}{Y}$ unique]
  \label{thm:PointedCat-UniqueBasePointIn-Hom(X,Y)}
  If a category enriched over pointed sets has a zero-object, then in any hom-set $\Hom{X}{Y}$ the base point is the zero map.
\end{lemma}
\begin{proof}
  Let $\ast_{YX}\from X\to Y$ denotes the base point of $\Hom{X}{Y}$. Then, for arbitrary ${g\from Y\to Z}$ and $f\from U\to X$, we must have $g\Comp \ast_{YX}=\ast_{ZX}$ and $\ast_{YX}\Comp f=\ast_{YU}$. Therefore,
  \begin{equation*}
    \ast_{YX} = \ast_{Y0}\Comp \ast_{0X}=\ZeroMap_{Y0}\Comp 0_{0X} = \ZeroMap_{YX}
  \end{equation*}
  This is what we wanted to show.
\end{proof}

In a category with zero object, by (\ref{thm:PointedCat-UniqueBasePointIn-Hom(X,Y)}), we use the notation $0_{XY}$ for both the map which factors through a zero object and the base-point in a pointed hom-set $\HomBsd{X}{Y}$.

\begin{exercises}

\begin{exercise}[Split monomorphisms in $\SetsBsd$]
  \label{exe:SplitMonoSet*}
  Show that in the category $\SetsBsd$ of pointed sets and base-point preserving functions, a morphism is a split monomorphism if and only if it is an injective base-point preserving function. %
  \index{split!monomorphism in $\SetsBsd$}%
\end{exercise}

\begin{exercise}[One-point union versus cartesian product in $\SetsBsd$]\label{exe:OnePointUnion}
  In $\SetsBsd$, the coproduct is called the \Defn{one-point union}: given two pointed sets $(X,x)$ and $(Y,y)$, it is obtained as the disjoint union of $X$ and $Y$, modulo identification of $x$ and $y$. Check this. Then use \eqref{exe:SplitMonoSet*} to prove that the canonical comparison map $(X,x)+(Y,y)\to (X,x)\times(Y,y)$ is a split monomorphism.
\end{exercise}

The following two exercises provide examples of non-pointed categories which are canonically enriched in $\SetsBsd$.

\begin{exercise}[Non-zero rings are enriched in $\SetsBsd$]
  \label{exe:NonZeroRingsEnriched-Set_*}
  Show that the category of non-zero rings ($\Rngs$ from which the zero object is removed) is non-pointed category which is enriched over pointed sets.

  Formulate a generalization of the above which starts with an arbitrary pointed category.
\end{exercise}

\begin{exercise}[Epi / monomorphism in $\Sets$ and in $\SetsBsd$]
  \label{exe:Epi/Mono-Set/Set_*}
  About a function $f\from X\to Y$ in $\Sets$ or in $\SetsBsd$, show that the following hold:
  \begin{enumerate}[(i)]
    \item $f$ is a monomorphism if and only if it is injective.
    \item $f$ is an epimorphism if and only if it is surjective.
  \end{enumerate}
\end{exercise}

\begin{exercise}[Surjective functions are pullback stable]
  \label{exe:Surjective->PullbackStable}
  Show that the pullback of a surjective function $g\from Z\to Y$ along any function $f\from X\to Y$ is again surjective.
\end{exercise}

\begin{exercise}[Image factorization in $\Sets$ and in $\SetsBsd$]
  \label{exe:ImageFactorization-Set/Set_*}
  About a function $f\from X\to Y$ in $\Sets$ or in $\SetsBsd$, show that the following hold:
  \begin{enumerate}[(i)]
    \item $E\DefEq \KrnlPr{f}$ is an equivalence relation on $X$.
    \item $q\from X\to X/E$ is a coequalizer for the projections $\PrjctnOnto{1},\PrjctnOnto{2}\from \KrnlPr{f}\to X$.
    \item $f$ factors through $q$ as $f=i\Comp q$, where $i$ is injective.
  \end{enumerate}
\end{exercise}

Below, we collect a guided suite of properties of the category $\Tops$ of topological spaces and the category $\TopsBsd$ of based topological spaces which are relevant to our purposes, in particular to working in topological varieties.

\begin{exercise}[Forget $\TopsBsd\to \SetsBsd$ preserves limits and colimits]
  \label{exe:Set_*->Top_*}
  Show that the following: %
  \index{limit!in $\Tops$}\index{limit!in $\TopsBsd$}\index{colimit!in $\Tops$}\index{colimit!in $\TopsBsd$}%
  \begin{thmlist}
    \item The forget functors $\Tops\to \Sets$ and $\TopsBsd\to \SetsBsd$ preserve and reflect limits.
    \item The forget functors $\Tops\to \Sets$ and $\TopsBsd\to \SetsBsd$ preserve colimits.
  \end{thmlist}
  Hint: Each of these forget functors has a left and a right adjoint.
\end{exercise}

\begin{exercise}[Subobjects in $\Tops$ and $\TopsBsd$]
  \label{exe:Top/Top_*-Subobjects}
  Show that a subobject of a (based) topological space $X$ is given by the inclusion of a subset of $X$ equipped with the subspace topology or a refinement of the subspace topology. %
  \index{subobject!in $\Tops$}\index{subobject!in $\TopsBsd$}%
\end{exercise}

\begin{exercise}[Compact Hausdorff spaces can not be refined]
  \label{exe:CHD-NotRefinable}
  In the category $\Tops$, if $X$ is a compact Hausdorff space, and $f\from A\to X$ is a continuous bijection, then exactly one of the following two applies:
  \begin{thmlist}
    \item $f$ is a homeomorphism.
    \item $A$ is not compact.
  \end{thmlist}
\end{exercise}

\begin{exercise}[Compact Hausdorff spaces]
  \label{exe:CompactHausdorff}
  In the category $\CHTops$ of compact Hausdorff spaces show the following:
  \vspace{-2ex}
  \begin{thmlist}
    \item A subobject of $X$ is given by a closed subspace of $X$ with the subspace topology.
    \item A proper refinement of the topology of $X$ yields a space which is no longer compact.
  \end{thmlist}
\end{exercise}
\end{exercises}
\section[Kernels and Cokernels]{Kernels and Cokernels}
\label{sec:Kernel/CoKernel}

In this section we introduce the categorical concepts of kernel and cokernel, and we develop their basic properties from the ground up. For example, we show that every kernel has an equivalent characterization as a certain kind of pullback. Dually, every cokernel may equivalently be expressed as a certain kind of pushout; see (\ref{thm:Ker/CoKerExist}). Further, by (\ref{thm:Ker(Mono)=0}) the kernel of a monomorphism is $0$, as is the cokernel of an epimorphism (\ref{thm:CoKer(Epi)=0}).  Additional properties can be found in (\ref{thm:NormalMono-Props}) and in (\ref{thm:NormalEpi-Props}). - Let us now turn to details.

\begin{definition}[Kernel]
  \label{def:Kernel}
  In category with a zero object a \Defn{kernel} of a morphism $f\from X\to Y$ is given by a map $\kappa\from K\to X$ such that these two conditions hold: %
  \index{kernel}\index[not]{k!$\KerMap{f}$\IndSep kernel of $f$}\index[not]{k!$\Ker{f}$\IndSep kernel object of $f$}%
  \begin{thmlist}
    \item $f\kappa=\ZeroMap$, and
    \item any map $\kappa'\from K'\to X$ with $f\kappa'=0$ factors uniquely through $\kappa$.
  \end{thmlist}
  In this situation we write $\KerMap{f}$ for $\kappa$ and $\Ker{f}$ for $K$.
\end{definition}

Dualizing, we obtain:

\begin{definition}[Cokernel]
  \label{def:CoKernel}
  A \Defn{cokernel} of a morphism $f\from X\to Y$ is given by a map $\pi\from Y\to Q$ such that: %
  \index{cokernel}\index[not]{c!$\CoKerMap{f}$\IndSep cokernel of $f$}\index[not]{c!$\CoKer{f}$\IndSep cokernel object of $f$}%
  \begin{thmlist}
    \item $\pi f=\ZeroMap$, and
    \item any map $\pi'\from Y\to Q'$ with $\pi' f=0$ factors uniquely through $\pi$.
  \end{thmlist}
  In this situation we write $\CoKerMap{f}$ for $\pi$ and $\CoKer{f}$ for $Q$.
\end{definition}

In diagram language, a map $\kappa\from K\to X$ represents the kernel of $f\from X\to Y$ if and only if, on the left below, every diagram of solid arrows admits a unique filler $\lambda$ as shown.
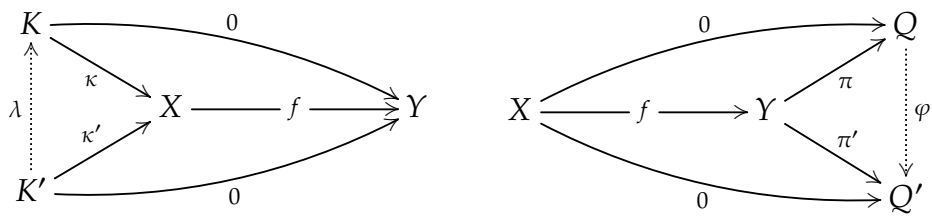
\begin{equation*}
  \xymatrix@R=3ex@C=3em{
  K \ar@/^2ex/[rrrd]^-{\ZeroMap} \ar[rd]_-{\kappa} &&& \\
  & X \ar[rr]|-{\ f\ } &&
  Y \\
  K' \ar[ru]^-{\kappa'} \ar@{.>}[uu]^{\lambda} \ar@/_2ex/[rrru]_-{\ZeroMap}
  }\qquad
  \xymatrix@R=3ex@C=3em{
  &&& Q \ar@{.>}[dd]^-{\varphi} \\
  X \ar[rr]|-{\ f\ } \ar@/^2ex/[rrru]^-{\ZeroMap} \ar@/_2ex/[rrrd]_-{\ZeroMap} &&
  Y \ar[ru]_{\pi} \ar[rd]^{\pi'} \\
  &&& Q'
  }
\end{equation*}
Dually, a map $\pi\from Y\to Q$ represents the cokernel of $f\from X\to Y$ if and only if, on the right above, every diagram of solid arrows admits a unique filler $\varphi$ as shown.  Consequently (\ref{exe:Kernel->Monomorphism}), if $\kappa\from K\to X$ is a kernel of $f$, then $\kappa$ is a monomorphism. Moreover, any two kernels of $f$ are related by a unique isomorphism.

This means that, whenever a kernel of $f$ exists, it represents a special and unique kind of subobject of $X$, which is what we call \emph{the} kernel of $f$. It is common practice to accept the inaccuracy of writing $\kappa=\KerMap{f}$ or, even shorter, $K=\Ker{f}$. We call a morphism $k$ \Defn{a kernel} when it represents the kernel of some map.

\begin{definition}[Normal monomorphism]
  \label{def:NormalMono}%
  A morphism $\kappa$ in a category with a zero object is called a \Defn{normal monomorphism} if it is the kernel of some arrow. We write $\kappa\from K\NMono X$ to indicate that $\kappa$ is a normal monomorphism, and say that $K$ is a normal subobject \emph{of $X$}, denoted $K\normal X$. %
  \index{normal!monomorphism}%
  \index{monomorphism!normal}%
  \index{normal!subobject}\index[not]{n!$K\normal X$\IndSep $K$ is normal subobject of $X$}%
  \index[not]{n!$\kappa\from K\NMono X$\IndSep $\kappa$ is normal monomorphism}%
\end{definition}

Dually, in (\ref{def:CoKernel}) it follows that the map $\pi$ is an epimorphism and, further, that any two cokernel maps are related by a unique isomorphism. This means that, whenever a cokernel of $f$ exists, it represents a special and unique kind of quotient object of $Y$, called \emph{the} cokernel of $f$. We write $\pi=\CoKerMap{f}$ or $Q=\CoKer{f}$ to say that $\pi$, respectively $Q$, is a cokernel of $f$. %
If $K\normal X$ represented by a normal monomorphism $k\from K\to X$, then we write $X/K$ for the quotient of $X$ by $K$, the cokernel $Q$ of $k$.

\begin{definition}[Normal epimorphism]
  \label{def:NormalEpi}%
  In a category with a zero object, a morphism $\pi\from Y\to Q$ is called a \Defn{normal epimorphism} if it is the cokernel of some arrow. We write $\pi\from Y\NEpi Q$ to indicate that $q$ is a normal epimorphism. %
  \index{normal!epimorphism}%
  \index{epimorphism!normal}%
  \index[not]{n!$\pi\from Y\NEpi Q$\IndSep $\pi$ is normal epimorphism}%
\end{definition}

Let us now turn to the question of whether kernels, respectively cokernels, exist:

\begin{proposition}[(Co)Kernels and pushout/pullback diagrams]
  \label{thm:Ker/CoKerExist}
  \label{thm:Ker/CoKerVsPullhPush}%
  For a morphism $f\from X\to Y$ in a category with a zero object the following hold.
  \begin{thmlist}
    \item $f$ has a kernel if and only if the diagram $X\XRA{f} Y \longleftarrow 0$ admits a pullback.
    \item $f$ has a cokernel if and only if the diagram $0 \longleftarrow X\XRA{f} Y$ admits a pushout.
  \end{thmlist}
\end{proposition}
\begin{proof}
  The commutative square on the left below has the pullback property exactly when $\kappa$ has the universal property of the kernel of $f$.
  \[
    \xymatrix@R=5ex@C=3em{
    K \ar@{{ |>}->}[d]_{\kappa} \ar[r] \PullLU{rd} &
    0 \ar[d] \\
    X \ar[r]_-{f} &
    Y
    }\qquad\qquad\qquad
    \xymatrix@R=5ex@C=3em{
    X \ar[r]^-{f} \ar[d] \PushRD{rd} &
    Y \ar@{-{ >>}}[d]^{\pi} \\
    0 \ar[r] &
    C
    }
  \]
  Similarly, the commutative square on the right above has the pushout property exactly when $\pi$ has the universal property of the cokernel of $f$.
\end{proof}

\begin{example}[Various kernels and cokernels]
  \label{exa:Kernels/CoKernels-Examples}%
  Familiar kernels and cokernels include:
  \begin{thmlist}
    \item In $\Grps$, kernels correspond to inclusions of normal subgroups; in $\Rngs$, as in $\Lie_R$, a kernel is the inclusion of an ideal.
    \item In $\SetsBsd$ the kernel of $f\from {(X,x)\to (Y,y)}$ is the pointed set $(\Set{z\in X\mid f(z)=y},x)$. The cokernel of $f$ is the pointed set $(Y/f(X),\Set{f(X)})$ where $Y/f(X)$ is the set $Y$ modulo the equivalence relation $R$ defined by $aRb$ iff $a=b$ or $a$, $b\in f(X)$.
    \item In $\Monoids$, the addition morphism $+\from \Prdct{\NNr}{\NNr}\to \NNr$ has kernel $\Set{(0,0)}$.
    \item Any pointed category admits a kernel and a cokernel for any isomorphism, namely~$0$. Conversely, $\IdMapOn{A}$ is a kernel of $A\to 0$ and a cokernel of $0\to A$.
  \end{thmlist}
\end{example}

Kernels and cokernels determine each other in the following sense.

\begin{proposition}[Kernel of cokernel of kernel]
  \label{thm:Ker(CoKer)-CoKer(Ker)}
  In a category with a zero object, the following hold:
  \begin{thmlist}
    \item if $k$ is a normal monomorphism which admits a cokernel, then $k=\KerMap{\CoKerMap{k}}$;
    \item if $f$ is a normal epimorphism which admits a kernel, then $f=\CoKerMap{\KerMap{f}}$.
  \end{thmlist}
  Hence, in a pointed category, any normal monomorphism is a kernel of its cokernel, while any normal epimorphism is a cokernel of its kernel.
\end{proposition}
\begin{proof}
  We prove (i), as (ii) is its dual. Suppose $k\from K\to X$ is a kernel of $f\from X\to Y$ and consider a cokernel $\varepsilon\from X\to Q$ of $k$. Then there is a unique $g\from Q\to Y$ such that $g\Comp \varepsilon=f$. We prove that $k$ is a kernel of $\varepsilon$ by verifying the universal property. We let $h\from Z\to X$ satisfy $\varepsilon\Comp h=0$. Then also $f\Comp h=g\Comp \varepsilon\Comp h$ is trivial, so that a unique $l\from Z\to K$ exists for which $k\Comp l=h$. This proves the claim.
\end{proof}

It can be helpful to interpret (\ref{thm:Ker(CoKer)-CoKer(Ker)}) via commutative diagrams:
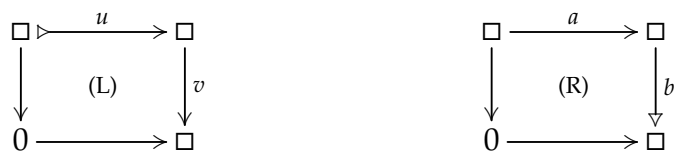
\begin{equation*}
  \xymatrix@R=5ex@C=4em{
  \DiagObj \ar@{{ |>}->}[r]^{u} \ar[d] \ar@{}[dr]|-{\text{(L)}} &
  \DiagObj \ar[d]^-{v} &&
  \DiagObj \ar[r]^-{a} \ar[d] \ar@{}[dr]|-{\text{(R)}} &
  \DiagObj \ar@{-{ >>}}[d]^{b} \\
  0 \ar[r] &
  \DiagObj &&
  0 \ar[r] &
  \DiagObj
  }
\end{equation*}
If the square (L) is a pushout with $u$ a normal monomorphism, then $v=\CoKer{u}$  and (L) is also a pullback. If the square (R) is a pullback with $b$ a normal epimorphism, then $a=\Ker{b}$ and (R) is also a pushout.

\begin{subordinate}{}
  \begin{subsubordinate}{On the concept of normal monomorphism}
    Some authors define the term normal monomorphism differently; see e.g.\  \cite{DBourn2000}, \cite[Lemma~5.9]{Bourn-Gran-CategoricalFoundations}, \cite{FBorceuxDBourn2004}. These works conceptualize a different aspect of `being normal'---namely, \emph{being a class of an internal equivalence relation}. In a category which is not Barr exact, this latter notion does not agree with the concept of a kernel.

    One example where this happens is the category of finitely generated free abelian groups. This category is additive and admits kernels. However, there are monomorphisms which are not kernels. Yet all monomorphisms are normal in the sense of Bourn, see \cite[Lemma~5.9]{Bourn-Gran-CategoricalFoundations}. One can show that in a semiabelian category a monomorphism is a kernel, if and only if it normal in the sense of \cite{DBourn2000}, \cite[Lemma~5.9]{Bourn-Gran-CategoricalFoundations}, \cite{FBorceuxDBourn2004}. Also, a morphism is monic if and only if its kernel vanishes; see (\ref{thm:Mono<->0-Kernel}).
  \end{subsubordinate}
\end{subordinate}

\begin{exercises}

\begin{exercise}[Isomorphism and normal epi / mono]
  \label{exe:Iso->NormalEpi/Mono}
  In any pointed category, show that an isomorphism is both, a normal epimorphism and a normal monomorphism.
\end{exercise}

\begin{exercise}[Equivalent condition for kernel\ZExactTag]
  \label{exe:Kernel-EquivalentCondn}
  Given a map $f\from X\to Y$, show that the following are equivalent for a map $m\from M\to X$:
  \begin{tfae}
    \item $m$ represents $\Ker{f}$.
    \item $m$ represents the greatest subobject of $X$ for which $f\Comp m=\ZeroMap$.
  \end{tfae}
\end{exercise}

\begin{exercise}[Equivalent condition for cokernel\ZExactTag]
  \label{exe:CoKernel-EquivalentCondn}
  Given a map $f\from X\to Y$, show that the following are equivalent for a map $q\from Y\to Q$:
  \begin{tfae}
    \item $q$ represents $\CoKer{f}$.
    \item $q$ represents the greatest subobject of $y$ for which $q\Comp f=\ZeroMap$.
  \end{tfae}
\end{exercise}

\begin{exercise}[Normal monomorphism is monic\ZExactTag]
  \label{exe:Kernel->Monomorphism}
  Show that a kernel of $f\from X\to Y$ is a monomorphism. %
  \index{kernel!is monomorphism}
\end{exercise}

\begin{exercise}[$\Ker{\Ker{f}}=\ZeroObject$]
  \label{exe:Ker(Ker(f))=0}
  Given  $ f\from X\to Y$ with its kernel $ \mu\from \Ker{f}\to X$. Show that  $ \Ker{\mu}=0$.
\end{exercise}

\begin{exercise}[Composition of normal monomorphisms\ZExactTag]
  \label{exe:Kernels,Composition}%
  About composable maps  $ K \XRA{\mu} L \XRA{\nu} M$  prove the following: %
  \index{kernel!composition}%
  \begin{thmlist}
    \item If $\mu$ and $\nu$ are normal monomorphisms, show that  $ \nu\mu$  need not be a normal monomorphism.
    \item If  $ \nu\mu$  is a a normal monomorphism and $\nu$ is monic, show that $\mu$ is a normal monomorphism.
    \item Show that (ii) is false when $\nu$ is not monic.
  \end{thmlist}
\end{exercise}

\begin{exercise}[Some special cokernels]
  \label{exe:CoKer(0),CoKer(Epi)}%
  In a category with a zero object $\Ctgry{X}$ prove the following:
  \begin{thmlist}
    \item For $X$ in $\Ctgry{X}$ arbitrary, $\IdMapOn{X}\from X\to X$  is a cokernel of  $ z\from 0\to X$. %
    \index{cokernel!of $\ZeroObject\to X$}%
    \item If  $ f\from X\to Y$ is an isomorphism in $\Ctgry{X}$, then $f$ is a cokernel of  $ z\from 0\to X$. %
    \index{cokernel!of $\ZeroObject\to X$}%
    \item If  $ f\from X\to Y$  is a epimorphism, then  $\zeta\from Y\to 0$  is a cokernel of $f$. %
    \index{cokernel!of an epimorphism}%
  \end{thmlist}
\end{exercise}

\begin{exercise}
  \label{exe:Ker(0)}%
  In a category with a zero object show that, for arbitrary $X$, the identity $\IdMapOn{X}\from X\to X$ is a kernel of  $q\from X\to 0$ and the cokernel of $m\from \ZeroObject\to X$.
\end{exercise}

\begin{exercise}[$\CoKer{\CoKer{f}}=\ZeroMap$]
  \label{exe:CoKer(CoKer(f))=0}%
  Given  $ f\from X\to Y$, let  $ \varepsilon\from Y \to \CoKer{f}$  be its cokernel. Show that  $ \CoKer{\varepsilon}=0$.
\end{exercise}

\begin{exercise}[Subobject / quotient object of $\ZeroObject$\ZExactTag]
  \label{exe:Sub/QuotientOf-0}
  About the zero object in a category which has a zero object, show the following:
  \begin{thmlist}
    \item $\ZeroObject$ is the one and only subobject of $\ZeroObject$, and it is normal.
    \item $\ZeroObject$ is the one and only quotient object of $\ZeroObject$, and it is normal.
  \end{thmlist}
\end{exercise}

\begin{exercise}[$\Ker{f}=0$ does not imply $f$ is monomorphism]
  \label{exe:Ker(f)=0,But-f-not-mono}
  Give an example of a pointed category containing a morphism $f$ such that $\Ker{f}=\ZeroMap$, yet $f$ is not a monomorphism.
\end{exercise}

\begin{exercise}[Cokernel vanishes yet not epimorphic]
  \label{exe:Coker=0/NotEpi}
  In the category $\Grps$ find a non-epimorphic morphism whose cokernel vanishes.
\end{exercise}

\begin{exercise}[Normal epimorphism is epic]
  \label{exe:NormalEpi->Epi}
  In any pointed category, a normal epimorphism  $f\from Y\to Q$  is epic.
\end{exercise}

\begin{exercise}[Normal epi cancellation]
  \label{exe:NormalEpiCancellation}
  In any pointed category, if $gf$ is a normal epimorphism and $f$ is epic, then $g$ is a normal epimorphism.
\end{exercise}

\begin{exercise}[Composition of normal epimorphisms]
  \label{exe:NormalEpis-Composition}
  Find an example of a pointed category with normal epimorphisms $f$ and $g$ such that the composite $gf$ is \emph{not} normal epic.
\end{exercise}

\begin{exercise}[Non-cancellation for normal monomorphisms]
  \label{exe:KernelChangingComposites}
  Find a pointed category and a normal monomorphism composite of $X\XRA{f} Y \XRA{g} Z$ in which $f$ is not normal; compare (\ref{thm:MonomorphismCancellationInKernel}) in Proposition~\ref{thm:NormalMono-Props}.
\end{exercise}

\begin{exercise}[Normal epimorphisms / monomorphisms in $\SetsBsd$]
  \label{thm:NormalEpi/Mono-Set_*}%
  In the category $\SetsBsd$ of pointed sets, show a function $f\from (X,x)\to (Y,y)$ has the following properties:
  \begin{thmlist}
    \item $f$ is a normal epimorphism if and only if, for every $y\neq t\in Y$, its fiber $f^{-1}(t)$ under $f$ is a singleton. \index{normal!epimorphism in $\SetsBsd$}
    \item $f$ is a normal monomorphism if and only if it is an injective function. %
    \index{normal!monomorphism in $\SetsBsd$}
  \end{thmlist}
\end{exercise}
\end{exercises}
\section[z-Exact Categories]{z-Exact Categories}
\label{sec:z-ExactCats}

To unburden the subsequent exposition from existence considerations, we then specialize to categories with a zero object in which every map has a kernel and a cokernel. We say that such a category is \ZExact. Perhaps surprising is that, in \ZExact\ categories a primordial version of, what is commonly known as, the Third Isomorphism Theorem holds; see (\ref{thm:ThirdIsoThm - I}).

\begin{definition}[\ZExact\ category]
  \label{def:PointedCatWith(Co)Kernels}%
  A category with a zero object in which every morphism has a (functorially chosen) kernel and cokernel is called an \ZExact\ category. %
  \index{\ZExact\ category}%
\end{definition}

The label {\color{Cerulean} $\EuRoman{z\hy Ex}$} in the caption of a definition or of a result indicates that, what is captioned, applies to \ZExact\ categories. For example: %
\index[acr]{z!{\color{Cerulean} $\EuRoman{z\hy Ex}$}\IndSep \ZExact\ category}

\begin{corollary}[Normal subobject / quotient object inversion\ZExactTag]
  \label{thm:NormalSubObject/QuotientObjectInversion}
  For every object $X$ in a \ZExact\ category, the class $\NSubObjcts{X}$ of normal subobjects of $X$ is in bijective correspondence with the class $\NQuoObjcts{X}$ via the operations %
  \index[not]{n!$\NSubObjcts{X}$\IndSep category of normal subobjects of object $X$}%
  \index[not]{n!$\NQuoObjcts{X}$\IndSep category of normal quotient objects of object $X$}
  \begin{equation*}
    \xymatrix@R=5ex@C=4em{
    \NSubObjcts{X} \ar@<0.5ex>[r]^-{\CoKerFunc} &
    \NQuoObjcts{X} \ar@<0.5ex>[l]^-{\KerFunc}
    }
  \end{equation*}
\end{corollary}

Motivated by the effect of the correspondence in (\ref{thm:NormalSubObject/QuotientObjectInversion}) on the orders of normal subgroups / quotient groups of a finite group, we introduce the following terminology.

\begin{terminology}[Normal subobject / quotient object inversion]
  \label{term:NormalSubobject/QuotientObjectInversion}%
  We refer to the bijection in (\ref{thm:NormalSubObject/QuotientObjectInversion}) as the \Defn{normal subobject / quotient object inversion over $X$}. %
  \index{inversion!normal subobjects / quotient objects}\index{normal subobject / quotient object inversion}
\end{terminology}

Existence of kernels and cokernels suffices for the availability of certain pullbacks and pushouts.

\begin{proposition}[Existence of special pullbacks / pushouts\ZExactTag]
  \label{thm:Pullback/Pushout-Existence}%
  \label{thm:PushoutPreservesCokernel} 
  \label{thm:PushoutPreservesNormalEpis}%
  \label{thm:PullbackPreservesNormalMonos}%
  Let $\Ctgry{X}$ be a category with zero object in which every map has a kernel and a cokernel. Then the following hold:
  \begin{thmlist}
    \item The pullback of a normal monomorphism along an arbitrary morphism exists; such pullbacks preserve normal monomorphisms.
    \item The pushout of a normal epimorphism along an arbitrary morphism exists; such pushouts preserve normal epimorphisms.
  \end{thmlist}
\end{proposition}
\begin{proof}
  In the diagram below, we start from a normal monomorphism $\mu$ and an arbitrary morphism $f$.
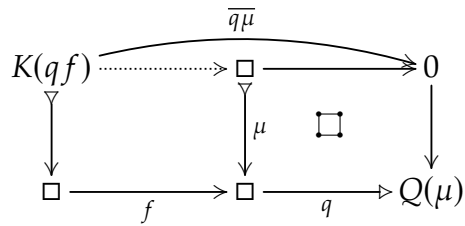
  \begin{equation*}
    \xymatrix@R=5ex@C=4em{
    \Ker{q f} \ar@{.>}[r] \ar@{{ |>}->}[d] \ar@/^2ex/[rr]^{\overline{q\mu}} &
    \DiagObj \ar[r] \ar@{{ |>}->}[d]^{\mu} \BiCart{rd} &
    \ZeroObject \ar[d] \\
    \DiagObj \ar[r]_-{f} &
    \DiagObj \ar@{-{ >>}}[r]_-{q} &
    \CoKer{\mu}
    }
  \end{equation*}
  Then $\mu$, together with its cokernel yields the bicartesian square on the right. The left hand square is constructed by first taking the kernel of the composite $qf$, then adding the uniquely induced dotted arrow. Then, pullback cancellation (\ref{thm:Pullbacks,ConcatenatedSquares}) shows that the square on the left is a pullback. Thus the pullback of the normal monomorphism $\mu$ along the arbitrary map $f$ exists and is again a normal monomorphism.

  The proof of (ii) is dual.
\end{proof}

In (\ref{thm:CoKer(Epi)=0}) and its dual (\ref{thm:Ker(Mono)=0}), we present ecessary conditions for being a (normal) epimorphism, respectively a (normal) monomorphism.

\begin{proposition}[$\CoKer{\textit{epi}}=0$]
  \label{thm:CoKer(Epi)=0}%
  In a pointed category, an epimorphism $q\from X\to Y$ has $\CoKerMap{q} = (Y\to \ZeroObject)$. \NoProof
\end{proposition}

\begin{proposition}[$\Ker{\mathit{mono}}=0$]
  \label{thm:Ker(Mono)=0}%
  In a pointed category, a monomorphism $f\from X\to Y$ has $\KerMap{f} = (\ZeroObject\to X)$.%
  \index{kernel!of monomorphism: $\ZeroObject$.} \NoProof
\end{proposition}

The converse of (\ref{thm:Ker(Mono)=0}) or (\ref{thm:CoKer(Epi)=0}) does not always hold, unless the category satisfies additional conditions; see (\ref{thm:Mono<->0-Kernel}). For example, in the category $\SetsBsd$ of pointed sets a function whose kernel is $\ZeroObject$ need not be a monomorphism; see exercise (\ref{exe:Ker(f)=0,But-f-not-mono}). Similarly, in $\Grps$ there are non-epimorphic maps whose cokernel vanishes; see Exercise~\ref{exe:Coker=0/NotEpi}.

Very useful is the following recognition criterion for isomorphisms.

\begin{proposition}[Isomorphism recognition]
  \label{thm:IsomorphismRecognition}%
  In a category with a zero object, the following are equivalent for a morphism $f\from X\to Y$: %
  \index{isomorphism!recognition}
  \begin{tfae}
    \item $f$ is an isomorphism;
    \item $f$ is a normal epimorphism and a monomorphism;
    \item $f$ is an epimorphism and a normal monomorphism.
  \end{tfae}
\end{proposition}
\begin{proof}
  The implication (i) implies (ii) is clear. To see that (ii) implies (i) suppose $f$ is a normal epimorphism and a monomorphism. Then the kernel of $f$ is $0$ by (\ref{thm:Ker(Mono)=0}). By (\ref{thm:Ker(CoKer)-CoKer(Ker)}), the normal epimorphism $f$ is the cokernel of $0\to X$, so that it is an isomorphism by item (iv) in Example~\ref{exa:Kernels/CoKernels-Examples}.

  Now (iii) implies (i) by duality. To see that (ii) implies (iii) we see with (\ref{thm:CoKer(Epi)=0}) that $\CoKerMap{f}=Y\to 0$. With (\ref{thm:Ker(CoKer)-CoKer(Ker)}) we have that $f=\KerMap{Y\to 0}$, and so $f$ is an isomorphism.
\end{proof}

\subsection*{Computations Involving (Co-)Kernels}

\begin{proposition}[Properties of normal monomorphisms\ZExactTag]
  \label{thm:NormalMono-Props}%
  \label{thm:KernelFunctor-Props}
  The following hold: %
  \index{normal!monomorphism: properties}%
  \vspace{-2ex}
  \begin{thmlist}
    \item \label{thm:Kernel(gf)}%
    Given two composable maps $ X \XRA{f} Y \XRA{g} Z$ in $\Ctgry{X}$, the kernel of $ gf$ is the pullback along $f$ of the kernel of $g$. %
    \index{kernel!of a composite of maps}%
    \item \label{thm:Ker(MonoComp f)}%
    If $g$ is a monomorphism, then $\KerMap{gf}=\KerMap{f}$.
    \item \label{thm:Ker(f)=mono-Comp-u}%
    In a composite $K \XRA{\varphi} U \XRA{m} X \XRA{f} Y$, if $\KerMap{f\from X\to Y} = m\varphi$ and $m$ is monic, then $\varphi = \KerMap{fm}$.
    \item \label{thm:MonomorphismCancellationInKernel}%
    If a composite $vu$ is a normal monomorphism and $v$ is a monomorphism, then $u$ is a normal monomorphism.
  \end{thmlist}
\end{proposition}
\begin{proof}
  (\ref{thm:Kernel(gf)})\quad See the proof of (\ref{thm:Pullback/Pushout-Existence}).

  (\ref{thm:Ker(MonoComp f)})\quad We know from (\ref{thm:Ker(Mono)=0}) that $\KerMap{g}=\ZeroMap$. By (\ref{thm:KernelFunctor-Props}), $\KerMap{gf}$ is the pullback of $\KerMap{g}=\ZeroMap$ along $f$. But that is exactly $\KerMap{f}$.

  (\ref{thm:Ker(f)=mono-Comp-u})\quad We verify the kernel property of $\varphi$ directly. Consider the situation depicted in the diagram below.
  
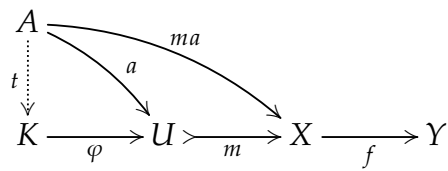
\begin{equation*}
    \xymatrix@R=5ex@C=3em{
    A \ar@{.>}[d]_{t} \ar@/^1ex/[rd]^(.6){a} \ar@/^2ex/[rrd]^{ma} \\
    K \ar[r]_-{\varphi} &
    U \ar@{{ >}->}[r]_-{m} &
    X \ar[r]_-{f} &
    Y
    }
  \end{equation*}
  If $(f m)a=0$, then $ f(ma)=0$, and so there exists $ t\from A\to K$, unique with the property $m\varphi t=ma$. As $m$ is monic, $\varphi t=a$. If $ t'\from A\to K$ is another map with $ \varphi t'=a$, then $ m\varphi t'=ma$, and the universal property of the kernel $m\varphi$ yields $t=t'$. This proves the claim.

  (\ref{thm:MonomorphismCancellationInKernel})\quad This follows from (\ref{thm:Ker(f)=mono-Comp-u}) by setting $f$ the cokernel of $vu$.
\end{proof}

In (\ref{thm:NormalMono-Props}.\ref{thm:MonomorphismCancellationInKernel}) the assumption that $m$ is monic is essential; see exercise (\ref{exe:KernelChangingComposites}). Dualizing Proposition~\ref{thm:NormalMono-Props}.\ref{thm:Ker(f)=mono-Comp-u} results in a valid statement about cokernels.

\begin{proposition}[Properties of normal epimorphisms\ZExactTag]
  \label{thm:NormalEpi-Props}%
  The following hold:
  \vspace{-2ex}
  \begin{thmlist}
    \item \label{thm:CoKer(gf)}%
    Given two composable maps $X \XRA{f} Y \XRA{g} Z$, the cokernel of $gf$ is the pushout  along $g$ of the cokernel of $f$. %
    \index{cokernel!of a composite of maps}%
    \item \label{thm:CoKer(gComp Epi)}%
    If $f$ is an epimorphism, then $\CoKerMap{gf}=\CoKerMap{g}$.
    \item \label{thm:NormalEpiComposite}%
    If a composite $vu$ is a normal epimorphism,  and $u$ is epic then $v$ is a normal epimorphism.\NoProof
  \end{thmlist}
\end{proposition}

\begin{corollary}[Factoring the cokernel of a map\ZExactTag]
  \label{thm:FactoringCoKer}
  Given composable maps $X \XRA{f} Y \XRA{g} Z$, let $\gamma\from \CoKerMap{f}\to \CoKerMap{gf}$ be the induced map. Then
  \begin{equation*}
    \CoKerMap{g} = \CoKerMap{\gamma}\Comp \CoKerMap{g\Comp f}
  \end{equation*}
\end{corollary}
\begin{proof}
  We consider the following diagram (cf.\ (\ref{thm:NormalEpi-Props}.\ref{thm:CoKer(gf)})).
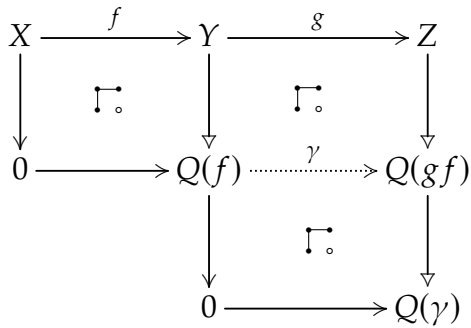
  \begin{equation*}
    \xymatrix@R=6ex@C=4em{
    X \ar[r]^-{f} \ar[d] \PushRD{rd} &
    Y \ar[r]^-{g} \ar@{-{ >>}}[d] \PushRD{rd} &
    Z \ar@{-{ >>}}[d] \\
    \ZeroObject \ar[r] &
    \CoKer{f} \ar@{.>}[r]^-{\gamma} \ar[d] \PushRD{rd}&
    \CoKer{gf} \ar@{-{ >>}}[d] \\
    & 0 \ar[r] &
    \CoKer{\gamma}
    }
  \end{equation*}
  The concatenation of the two pushout squares on the right is again a pushout by (\ref{thm:PushOuts,ConcatenatedSquares}). So, the composite of the vertical maps on the right represents the cokernel of $g$.
\end{proof}

If $X$ is a subobject of $Y$, write $Y/X$ for $\CoKer{X\Mono Y}$. With this notation, a very special case of (\ref{thm:FactoringCoKer}) is:

\begin{corollary}[Primordial Third Isomorphism Theorem\ZExactTag]
  \label{thm:ThirdIsoThm - I}
  Given subobjects $X< Y< Z$, then %
  \index{Isomorphism Theorem!Third}%
  \begin{equation*}
    Z/Y \cong \CoKer{(Y/X) \to (Z/X)}
  \end{equation*}
\end{corollary}

We will frequently encounter situations in which we commute one limit with another, or one colimit with another. Such reasoning is supported by `Fubini's Theorem',  \cite[Section~IX.2]{SMacLane1998}. Here is a particularly useful sampling of such results. As limits commute, kernels commute with products:

\begin{proposition}[Kernels commute with products\ZExactTag]
  \label{thm:KernelsCommuteProducts}%
  If in a \ZExact\ category the product of morphisms morphisms $f\from A\to B$ and $g\from X\to Y$ exists, then the identity $\Ker{\Prdct{f}{g}} = \Prdct{\Ker{f}}{\Ker{g}}$ holds. \NoProof%
  \index{kernel!of product of maps}%
\end{proposition}

Dualizing, we obtain:

\begin{proposition}[Cokernels commute with sums\ZExactTag]
  \label{thm:CoKernelsCommuteSums}%
  If in a \ZExact\ category the sum of morphisms morphisms $f\from A\to B$ and $g\from X\to Y$ exists, then the identity $\CoKer{f+g} = \CoKer{f}+\CoKer{g}$ holds. \NoProof%
  \index{cokernel!of sum of maps}%
\end{proposition}

\begin{proposition}[Cokernel of the inclusion into a sum\ZExactTag]
  \label{thm:CoKer(InclusionInSum)}
  For any two objects $X$ and $Y$ which admit a coproduct, the cokernel of the structure map ${\InclsnOf{X}\from X\to  X+Y}$ is the universal map $\SumMapOutOf{\ZeroMap,\OneMapOn{Y}}\from X+Y\to Y$. %
  \index{normal epimorphism!given by $\SumMapOutOf{\IdMap,\ZeroMap}\from X+Y\NEpi X$}%
\end{proposition}
\begin{proof}
  The concatenation of pushouts below is again a pushout:
  
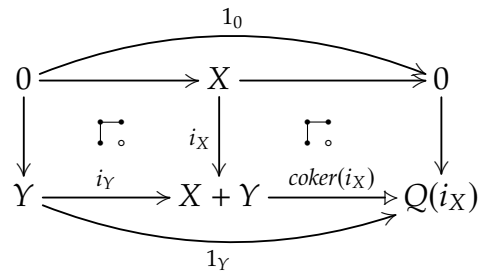
\begin{equation*}
    \xymatrix@R=5ex@C=4em{
    0 \ar[r] \ar[d] \ar@/^3ex/[rr]^-{\OneMapOn{0}} \PushRD{rd} &
    X \ar[d]_{\InclsnOf{X}} \ar[r] \PushRD{rd} &
    0 \ar[d] \\
    Y \ar[r]^-{\InclsnOf{Y}} \ar@/_3ex/[rr]_-{\OneMapOn{Y}}&
    X+Y \ar@{-{ >>}}[r]^-{\CoKerMap{\InclsnOf{X}} } &
    \CoKer{\InclsnOf{X}}
    }
  \end{equation*}
  As the composite along the top is the identity, so is the composite along the bottom. The diagram then explains why $\CoKerMap{\InclsnOf{X}} = \SumMapOutOf{0,\IdMapOn{Y}}$.
\end{proof}

Dually:

\begin{proposition}[Kernel of the projection from a product\ZExactTag]
  \label{thm:Ker(ProjectionFromProduct)}
  For any two objects $X$ and $Y$ which admit a product, the kernel of the projection map $\PrjctnOnto{X}\from \Prdct{X}{Y}\to X$ is the universal map $\PrdctMapInto{\ZeroMap,\IdMapOn{Y}}\from Y\to \Prdct{X}{Y}$. %
  \index{normal monomorphism!given by $\PrdctMapInto{\IdMap,\ZeroMap}\from X\to \Prdct{X}{Y}$}%
  \NoProof
\end{proposition}

\begin{subordinate}{}
  \begin{subsubordinate}{The Primordial Third Isomorphism Theorem (\ref{thm:ThirdIsoThm - I})}\quad is valid in all \ZExact\ categories. In the literature the Third Isomorphism Theorem involves normal subobjects $X\normal Y\normal Z$, with $X\normal Z$, and concludes with a short exact sequence
    \begin{equation*}
      (Y/X) \NMono (Z/X) \NEpi (Z/Y)
    \end{equation*}
    Such a conclusion may not be available in an arbitrary \ZExact\ category.  In Section \ref{sec:HomologicalSelfDuality} we identify a necessary and sufficient condition under which this conclusion holds.
  \end{subsubordinate}

  \begin{subsubordinate}{On the concept of \ZExact\ categories}
    In (\ref{def:PointedCatWith(Co)Kernels}) we defined a \ZExact\ category as a category with a zero-object, and where every morphism has a kernel and a cokernel. As such it is a special and motivating special case of the concept of $Z$-exact category in the sense of Grandis~\cite{Grandis-HA2}, where $Z$ is a chosen class of objects which serve as zeroes. Ours is the special case where $Z$ consists of the object $0$.
  \end{subsubordinate}
\end{subordinate}

\begin{exercises}
\begin{exercise}[Cokernels commute with sums / kernels commmute with products]
  \label{exe:CoKernelsCommuteSums/KernelsCommuteProducts}
  Prove propositions (\ref{thm:KernelsCommuteProducts}) and (\ref{thm:CoKernelsCommuteSums}).
\end{exercise}

\begin{exercise}[Kernel from kernel pair]
  \label{thm:Kernel-KernelPair}%
  Let $f\from X\to Y$ be a morphism and $(\KrnlPr{q},\PrjctnOnto{1},\PrjctnOnto{2})$ a kernel pair, a pullback of $f$ along itself. Then a kernel of $f$ is represented by the composite
  \begin{equation*}
    \PrjctnOnto{1}\Comp \KerMap{\PrjctnOnto{2}}\from \Ker{\PrjctnOnto{2}}\longrightarrow X.
  \end{equation*}
\end{exercise}
\end{exercises}
\section[Pullback and Pushout Recognition]{Recognizing Pullback Squares and Pushout Squares}
\label{sec:PullbackPushout-Recognition}

We develop frequently used facts about pullbacks and pushouts in their relation to kernels and cokernels.

\begin{proposition}[Pullback yields isomophic kernels\ZExactTag]
  \label{thm:Pullback->IsoOfKernels}%
  In the commutative diagram
  \begin{equation*}
    \xymatrix@R=5ex@C=3em{
    \Ker{\bar{f}} \ar@{{ |>}->}[r] \ar[d]_{\tilde{g}} &
    P \ar[r]^-{\bar{f}} \ar[d]_{\bar{g}} \ar@{}[rd]|-{\text{(R)}} &
    U \ar[d]^{g} \\
    \Ker{f} \ar@{{ |>}->}[r] &
    V \ar[r]_-{f} &
    W,
    }
  \end{equation*}
  if the square (R) is a pullback, then $\tilde{g}$ is an isomorphism.
\end{proposition}
\begin{proof}
  Consider the commutative diagram below.
  $$
    
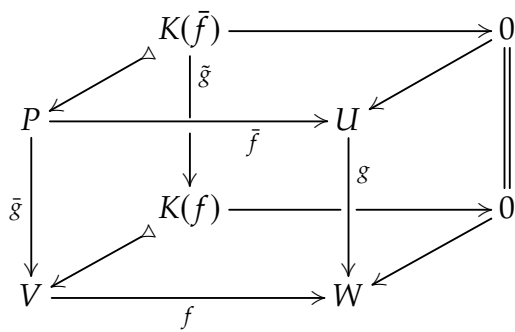
\xymatrix@R=3ex@!C=3em{
    & \Ker{\bar{f}} \ar@{{ |>}->}[ld] \ar[rr] \ar'[d]^{\tilde{g}}[dd] &&
    0 \ar@{=}[dd] \ar[ld] \\
    P \ar[rr]_(.7){\bar{f}} \ar[dd]_{\bar{g}} &&
    U \ar[dd]^(.3){g} \\
    & \Ker{f} \ar'[r][rr] \ar@{{ |>}->}[ld] &&
    0 \ar[ld] \\
    V \ar[rr]_-{f} &&
    W
    }
  $$
  The front, top, and bottom faces of this cube are pullbacks by design. Therefore by (\ref{thm:Pullbacks,ConcatenatedSquares}.i), the rectangle from the top back edge to the front bottom edge is a pullback as well. But then pullback cancellation (\ref{thm:Pullbacks,ConcatenatedSquares}.ii) says that also the back face is a pullback. The identity map $ 0\to 0$ pulls back to an isomorphism, namely $\tilde{g}$. --- This was to be shown.
\end{proof}

\begin{proposition}[Pullback recognition: kernel side I\ZExactTag]
  \label{thm:PullbackRecognition-KernelSide-1}%
  In the commutative diagram below, let $\alpha=\Ker{q}$, and $\alpha'$ a monomorphism such that $q'\alpha'=0$.%
  \index{pullback!recognition, kernel side I}
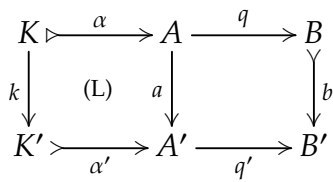
  $$
    \xymatrix@R=5ex@C=3em{
    K \ar@{{ |>}->}[r]^-{\alpha} \ar[d]_{k} \ar@{}[rd]|-{\text{(L)}} &
    A \ar[r]^-{q} \ar[d]_{a} &
    B \ar@{{ >}->}[d]^{b} \\
    K' \ar@{{ >}->}[r]_-{\alpha'} &
    A' \ar[r]_-{q'} &
    B'
    }
  $$
  If $b$ is a monomorphism then the square (L) is a pullback.
\end{proposition}
\begin{proof}
  Given arrows $u$ and $v$ as shown, we prove that the square on the left has the required universal property.
  $$
    
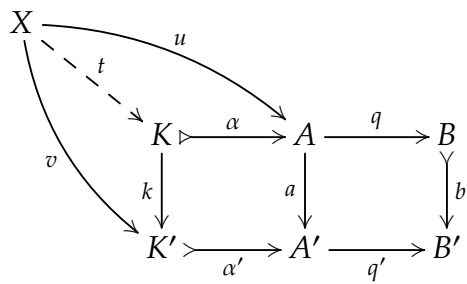
\xymatrix@R=5ex@C=3em{
    X \ar@/^4mm/[rrd]^{u} \ar@/_4mm/[rdd]_{v} \ar@{-->}[rd]^{t} & \\
    & K \ar@{{ |>}->}[r]^-{\alpha} \ar[d]_-{k} &
    A \ar[r]^-{q} \ar[d]_-{a} &
    B \ar@{ >->}[d]^-{b} \\
    & K' \ar@{{ >}->}[r]_-{\alpha'} &
    A' \ar[r]_-{q'} &
    B'
    }
  $$
  We find that $ b q u=q'au=q'\alpha' v=0$. As $b$ is a monomorphism, $q u=0$. Via the universal property of the kernel $\alpha$, there is a map $ t\from X\to K$, unique with the property $ \alpha t=u$. It remains to show that $ kt=v$. This follows from the defining property of the monomorphism $\alpha'$, upon observing that
  $$
    \alpha' kt = a \alpha t = au = \alpha' v.
  $$
  So $ kt=v$, as required.
\end{proof}

We can say more in the context of a homological category: see (\ref{thm:PullbackRecognition-KernelSide}).

The following propositions on pushouts and epimorphisms are dual to their respective counterparts on pullbacks and monomorphisms.

\begin{proposition}[Pushout yields isomorphism of cokernels\ZExactTag]%
  \label{thm:Pushout->IsoOfCoKers}%
  In the commutative diagram below, suppose $q=\CoKerMap{\alpha}$ and $q'=\CoKerMap{\alpha'}$. %
  \index{pushout!isomophic cokernels}%
  $$
    \xymatrix@R=5ex@C=3em{
    R \ar[r]^-{\alpha} \ar[d]_{r} &
    X \ar@{-{ >>}}[r]^-{q} \ar[d]_{x} &
    Y \ar[d]^{y}\\
    R' \ar[r]_-{\alpha'} &
    X' \ar@{-{ >>}}[r]_-{q'} &
    Y'
    }
  $$
  If the left hand square is a pushout then $y$ is an isomorphism. \NoProof
\end{proposition}

\begin{proposition}[Pushout recognition\ZExactTag]
  \label{thm:PushoutRecognize-Categorical}%
  In the commutative diagram below, suppose $q'=\CoKerMap{\alpha'}$, and $q$ is an epimorphism such that $q{\alpha}=0$. %
  \index{pushout!recognition - cokernel side}
  \begin{equation*}
    \xymatrix@R=5ex@C=3em{
    R \ar[r]^-{\alpha} \ar@{-{>>}}[d]_{r} &
    X \ar@{-{>>}}[r]^-{q} \ar[d]_{x} &
    Y \ar[d]^{y}\\
    R' \ar[r]_-{\alpha'} &
    X' \ar@{-{ >>}}[r]_-{q'} &
    Y'
    }
  \end{equation*}
  If $r$ is an epimorphism then the right hand square is a pushout. \NoProof
\end{proposition}
\newpage
\section[Normal Subobjects]{Normal Subobjects and Normal Closure}
\label{sec:NormalSubobjects}
\label{sec:NormalSubobjects/NormalClosure}

By Definition~\ref{def:NormalMono}, a subobject of an object $X$ is normal if it is the kernel of some map with domain $X$. In a \ZExact\ category we know that the intersection of two normal subobjects always exists by (\ref{thm:Pullback/Pushout-Existence}); here we prove that it is again normal, see (\ref{thm:NormalSubobjects-Intersection}). Then we show that every every subobject of $X$ has a unique normal closure.

\begin{proposition}[Intersection of normal subobjects\ZExactTag]
  \label{thm:NormalSubobjects-Intersection}%
  \label{thm:KernelsIntersection}
  \index{intersection!of normal subobjects}
  \index{meet!of normal subobjects}%
  The intersection of two normal subobjects is again a normal subobject. %
  \index[not]{m!$M\meet N$\IndSep intersection/meet of $M$ and $N$}%
  \index[not]{j!$M\meet N$\IndSep intersection/meet of $M$ and $N$}%
\end{proposition}
\begin{proof}
  Let $m\from M\to X$ be the kernel of $f$, and $n\from N\to X$ the kernel of $g$. From (\ref{subsec:Intersection/Union-Subobject}), we know that the pullback of $m$ and $n$ constructs their intersection $M\meet N$:
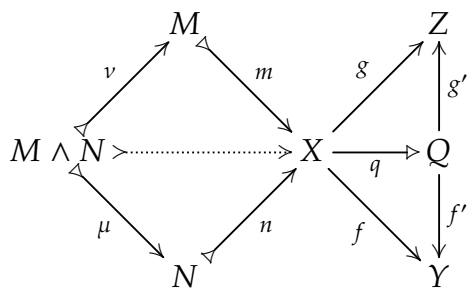
  \begin{equation*}
    \xymatrix@!0@=4em{
    & M \ar@{{ |>}->}[dr]^-{m} && Z\\
    M\meet N \ar@{{ |>}->}[ru]^-{\nu} \ar@{{ |>}->}[dr]_-{\mu} \ar@{ >.>}[rr] &&
    X \ar@{-{ >>}}[r]_{q} \ar[rd]_-f \ar[ru]^-g &
    Q \ar[u]_-{g'} \ar[d]^-{f'}
    \\
    & N \ar@{{ |>}->}[ru]_-{n} && Y
    }
  \end{equation*}
  We take the cokernel $q$ of $m\nu=n\mu$ and notice that $g$ factors through $q$ as a morphism $g'$, while $f$ factors through $q$ as a morphism $f'$.

  Consider a morphism $u\from {U\to X}$ such that $q\Comp u=0$. This implies that $f\Comp u=0$ and $g\Comp u=0$. So, there exist morphisms $x\from U\to M$ and $y\from {U\to N}$ for which $m\Comp x=u=n\Comp y$. Hence there exists a unique morphism $\PrdctMapInto{x,y}\from U\to M\meet N$ such that $\nu\Comp \PrdctMapInto{x,y}=x$ and $\mu\Comp \PrdctMapInto{x,y}=y$. We see that $(m\Comp \nu)\Comp \PrdctMapInto{x,y}=m\Comp x=u$, which proves that the intersection of $m$ and $n$ is the kernel of $q\from X \to Q$.
\end{proof}

Examples in the category $\Grps$ of groups show that a composite of two normal monomorphisms need not represent a normal subobject in general. On the other hand, as in $\Grps$, every subobject has a unique normal closure.

\begin{definition}[Normal closure of a subobject]
  \label{def:NormalClosure}%
  In a pointed category, let $m\from M\Mono X$ represent a subobject of $X$. Provided it exists, the \Defn{normal closure} $\bar{m}\from \bar{M}\NMono X$ of $m$ is the least normal subobject of $X$ containing $m$.  %
  \index{normal!closure of a subobject}\index{subobject!normal closure}%
\end{definition}

In a \ZExact\ category every subobject has a unique normal closure:

\begin{lemma}[Existence of normal closure\ZExactTag]
  \label{thm:NormalClosure-Exists}
  Let $m\from M\Mono X$ represent a subobject of $X$. Then $k\from \Ker{\CoKer{u}}\NMono X$ is the least normal subobject of $X$ containing $m$. %
  \index{normal!closure - existence}
\end{lemma}
\begin{proof}
  To see this, consider this commutative diagram:
  
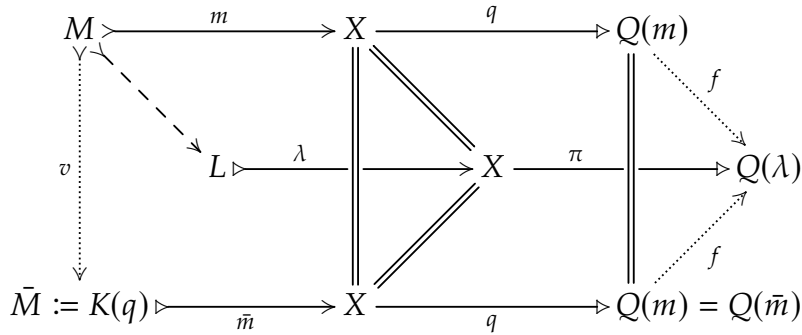
\begin{equation*}
    \xymatrix@!=5ex{
    M \ar@{{ >}->}[rr]^-{m} \ar@{{ >}.>}[dd]_{v} \ar@{{ >}-->}[rd] &&
    X \ar@{-{ >>}}[rr]^-{q} \ar@{=}[dd] \ar@{=}[rd]&&
    **[r]\CoKer{m}\ar@{=}[dd] \ar@{.>}[rd]^{f} \\
    & L \ar@{{ |>}->}[rr]|\hole^(0.3){\lambda} &&
    X \ar@{-{ >>}}[rr]|\hole^(0.3){\pi} &&
    \CoKer{\lambda} \\
    \bar{M}\DefEq \Ker{q} \ar@{{ |>}->}[rr]_-{\bar{m}} &&
    X \ar@{-{ >>}}[rr]_-{q} \ar@{=}[ru] &&
    **[r] \CoKer{m} = \CoKer{\bar{m}} \ar@{.>}[ru]_{f}
    }
  \end{equation*}
  Via the front facing part, we see that $m<\bar{m}$. To see that $\bar{m}$ represents the least normal subobject of $X$ containing $m$, consider an arbitrary normal subobject $\lambda\from L\NMono X$ containing $m$. It follows that $\pi$ factors uniquely through $q$ as shown. Consequently, $\bar{m}$ factors uniquely through $\lambda$; i.e., $\bar{m}<\lambda$. --- This completes the proof.
\end{proof}

Dually:

\begin{proposition}[Pushout of normal epimorphisms\ZExactTag]
  \label{thm:PushoutOfNormalEpis}%
  In a pushout diagram of normal epimorphisms %
  \index{pushout!of normal epimorphisms}%
  \begin{equation*}
    \xymatrix@R=5ex@C=4em{
    \DiagObj \ar@{-{ >>}}[r]^-{e} \ar@{-{ >>}}[d]_{\varepsilon} &
    \DiagObj \ar@{-{ >>}}[d]^{\underline{\varepsilon}} \\
    \DiagObj \ar@{-{ >>}}[r]_-{\underline{e}} &
    \DiagObj
    }
  \end{equation*}
  the composite $\underline{\varepsilon} e = \underline{e}\varepsilon$ is a normal epimorphism. \NoProof
\end{proposition}

\begin{exercises}
\begin{exercise}[Intersection when binary products exist] 
  \label{exe:NormalSubobjects-Intersection}%
  Give an alternate proof of (\ref{thm:NormalSubobjects-Intersection}): In a pointed category with binary products, the intersection of the kernels of maps $f\from {X\to Y}$ and $g\from {X\to Z}$ is the kernel of the map $\PrdctMapInto{f,g}\from X \to \Prdct{Y}{Z}$.
\end{exercise}
\end{exercises}
\section[Normal Decompositions and Factorizations]{Normal Decompositions and Factorizations}
\label{sec:NEM-ImageFac}
\label{sec:NormalDecompositionsFactorizations}

In a \ZExact\ category $\Ctgry{X}$, we present canonical ways of expressing a given morphism as a composite involving a strategically placed normal monomorphism, respectively normal epimorphism. This development forms a good foundation for a discussion of exact sequences, chain complexes, and homology. To explain the underlying idea, we use the interplay between three categories which are associated to $\Ctgry{X}$:
\begin{ulist}
  \item The \Defn{arrow category} $\ArrowCat{\Ctgry{X}}$ of $\Ctgry{X}$. Its objects are morphisms $\DiagObj\XRA{u} \DiagObj$, where $u$ is a morphism in $\Ctgry{X}$. A morphism $(a,b)\from u\to v$ is given by a commutative square in $\Ctgry{X}$: %
  \index{category!of arrows}\index{category!of morphisms}\index{arrow!category}%
  \index[not]{a!$\ArrowCat{\Ctgry{X}}$\IndSep category of morphisms in $\Ctgry{X}$}%
  \begin{equation*}
    \xymatrix@R=4ex@C=3em{
    \DiagObj \ar[d]_-{u} \ar[r]^-{a} &
    \DiagObj \ar[d]^{v} \\
    \DiagObj \ar[r]_-{b} &
    \DiagObj
    }
  \end{equation*}
  \item The category $\NMonoCat{X}$ is the full subcategory of $\Ctgry{X}$ whose objects are normal monomorphisms in $\Ctgry{X}$. %
  \index[not]{n!$\NMonoCat{X}$\IndSep category of normal monomorphisms in $\Ctgry{X}$}%
  \item The category $\NEpiCat{X}$ is the full subcategory of $\Ctgry{X}$ whose objects are normal epimorphisms in $\Ctgry{X}$. %
  \index[not]{n!$\NEpiCat{X}$\IndSep category of normal epimorphisms in $\Ctgry{X}$}%
\end{ulist}

An equivalence between the categories $\NMonoCat{X}$ and $\NEpiCat{X}$ is given by the functors
\begin{equation*}
  \CoKerFunc\from \NMonoCat{X} \longrightarrow \NEpiCat{X}  \qquad \text{and}\qquad \KerFunc\from \NEpiCat{X}\longrightarrow \NMonoCat{X}
\end{equation*}
Further, $\ArrowCat{\Ctgry{X}}$ contains $\NMonoCat{X}$ as a reflective subcategory, and it contains $\NEpiCat{X}$ as a coreflective subcategory; see Propositions \ref{thm:NMono(X)ReflectiveInArr(X)} and \ref{thm:NEpi(X)CoReflectiveInArr(X)}.

\begin{proposition}[$\NMonoCat{X}$ reflective in $\ArrowCat{\Ctgry{X}}$\ZExactTag]
  \label{thm:NMono(X)ReflectiveInArr(X)}%
  The category $\NMonoCat{X}$ is reflective in $\Ctgry{X}$, and the adjunction unit on an object $f$ is given by the commutative square on the left.
  \begin{equation*}
    \xymatrix@R=5ex@C=4em{
    \DiagObj \ar[r]^-{f} \ar@{.>}[d]_-{\eta^{0}_{f}} &
    \DiagObj \ar@{-{ >>}}[r]^-{q\DefEq \CoKerMap{f}} \ar@{=}[d] &
    \CoKer{f} \ar@{=}[d] \\
    \Ker{q} \ar@{{ |>}->}[r]_-{\KerMap{q}} &
    \DiagObj \ar@{-{ >>}}[r]_-{q} &
    \CoKer{f}
    }
  \end{equation*}
\end{proposition}
\begin{proof}
  This follows directly from the definitions.
\end{proof}

Dually:

\begin{proposition}[$\NEpiCat{X}$ coreflective in $\ArrowCat{\Ctgry{X}}$\ZExactTag]
  \label{thm:NEpi(X)CoReflectiveInArr(X)}%
  The category $\NEpiCat{X}$ is coreflective in $\Ctgry{X}$, and the adjunction counit on an object $u$ is given by the commutative square on the right.
  \begin{equation*}
    \xymatrix@R=5ex@C=4em{
    \Ker{f} \ar@{{ |>}->}[r]^-{k} \ar@{=}[d] &
    \DiagObj \ar@{=}[d] \ar@{-{ >>}}[r]^-{\CoKerMap{k}} &
    \CoKer{k} \ar@{.>}[d]^{\varepsilon_{f}^{1}}\\
    \Ker{f} \ar@{{ |>}->}[r]_-{k\DefEq \KerMap{f}} &
    \DiagObj \ar[r]_-{f} &
    \DiagObj
    }
  \end{equation*}
\end{proposition}

We are now ready to analyze decompositions $f=uv$ where $u$ is a normal monomorphism or $v$ is a normal epimorphism.

\begin{definition}[Normal decompositions\ZExactTag]
  \label{def:NEM-Image-Fac}
  \label{def:NormalDecompositions}%
  Consider a morphism $f\from X\to Y$ in a \ZExact\ category.
  \begin{thmlist}
    \item A composite $f=m u$ in which $m$ is a normal monomorphism is called a \Defn{normal mono decomposition of $f$}. %
    \index{normal!mono decomposition}%
    \item A composite $f=ve$ in which $e$ is a normal epimorphism is called a \Defn{normal epi decomposition of $f$}. %
    \index{normal!epi decomposition}%
  \end{thmlist}
\end{definition}

Thus $f=\IdMapOn{Y}\Comp f = f\Comp \IdMapOn{X}$ is a normal mono, respectively a normal epi, decomposition of $f$. Whenever we have simultaneously a normal epi and a normal mono decomposition of a morphism, then the two are uniquely related as follows.

\begin{proposition}[Relationship between normal epi / mono decompositions\ZExactTag]
  \label{thm:NormalEpi/MonoDecompositions-Relationship}%
  Given a morphism $f$ with a given normal mono decomposition $f=mu$ and a normal epi decomposition $f=ve$, there is a unique map $\varphi$ which renders the diagram below commutative.
  \begin{equation*}
    \xymatrix@R=5ex@C=4em{
    \DiagObj \ar@{-{ >>}}[r]_-{e} \ar@/^2ex/[rr]^-{f} &
    \DiagObj \ar[r]_-{v} \ar@{.>}[d]_{\varphi} &
    \DiagObj \\
    \DiagObj \ar[r]^-{u} \ar@{=}[u] \ar@/_2ex/[rr]_-{f}&
    \DiagObj \ar@{{ |>}->}[r]^-{m}  &
    \DiagObj \ar@{=}[u]
    }
  \end{equation*}
\end{proposition}
\begin{proof}
  Let $k\DefEq\KerMap{e}$. Then $e=\CoKerMap{k}$. Noting that $muk=\ZeroMap$, the monic property of $m$ yields $uk=\ZeroMap$. So, there exists $\alpha$, unique with $u=\alpha e$.

  Let $q\DefEq \CoKerMap{m}$. Then $m=\KerMap{q}$. Noting that $qve=\ZeroMap$, the epic property of $e$ yields $qv=\ZeroMap$. So, there exists $\beta$, unique with $v=m\beta$.

  We infer $\alpha=\beta$ from the computation $mae = mu = f = ve=m\beta e$. So, $\varphi\DefEq \alpha=\beta$ is the claimed map.
\end{proof}

Among all normal epi decompositions of a given map $f$, there is one class which is distinguished by a being terminal. Dually, among all normal mono decompositions of $f$, there is an initial class:

\begin{lemma}[Initial / terminal normal decompositions\ZExactTag]
  \label{thm:NormalDecompositions-Initial/Terminal}
  Every morphism $f\from X\to Y$ has an initial normal mono decomposition given by the construction on the left.%
  \index{initial!normal mono decomposition}\index{normal!mono decomposition, initial}%
  \index{terminal!normal epi decomposition}\index{normal!epi decomposition, terminal}%
  \begin{equation*}
    \xymatrix@R=5ex@C=4em{
    X \ar[r]^-{f} \ar@{=}[d] &
    Y \ar@{-{ >>}}[r]^-{q} &
    \CoKer{f}&
    \Ker{f} \ar@{{ |>}->}[r]^-{k} &
    X \ar[r]^-{f} \ar@{-{ >>}}[d]_-{\CoKerMap{k}\EqDef e}&
    Y \ar@{=}[d] \\
    X \ar@{.>}[r]_-{u\DefEq \eta_{f}^{0}} &
    K \ar@{{ |>}->}[u]_{m\DefEq \KerMap{q}} &&&
    Q \ar@{.>}[r]_-{v\DefEq \varepsilon_{f}^{1}} &
    Y
    }
  \end{equation*}
  Dually, $f$ has a terminal normal epi decomposition given by the construction on the right. %
\end{lemma}
\begin{proof}
  Apply the universal property of the unit $\eta$ of the reflection $\ArrowCat{\Ctgry{X}}\to \NMonoCat{X}$, respectively the counit $\varepsilon$ of the coreflection $\ArrowCat{\Ctgry{X}}\to \NEpiCat{X}$.
\end{proof}

Via its universal property, the initial normal mono decomposition of a map $f$ is unique up to isomorphism. We refer to it as the \Defn{normal mono factorization of $f$}. Similarly, we refer to the terminal normal epi decomposition of $f$ as the \Defn{normal epi factorization of $f$}. %
\index{normal!mono factorization}\index{normal!epi factorization}%

\begin{notation}[Normal epi / normal mono factorization comparison map]
  \label{not:NormalEpi/NormalMonoFactorizationComparison}
  Given a map $f$ in a \ZExact\ category, Proposition \ref{thm:NormalEpi/MonoDecompositions-Relationship} yields a unique map $\NENMComp{f}$ relating its normal epi factorization to its normal mono factorization. We call it the normal epi / normal mono factorization comparison map. %
  \index[not]{e!$\NENMComp{f}$\IndSep normal epi / normal mono factorization comparison map}%
\end{notation}

\begin{lemma}[Recognizing special initial/terminal decompositions\ZExactTag]
  \label{thm:Initial/TerminalDecomposition-Recognize}
  For a composite $f=vu$ the following hold:%
  \index{normal!factorization, recognize}%
  \begin{thmlist}
    \item If $u$ is an epimorphism and $v$ is a normal monomorphism, then $vu$ is the normal mono factorization of $f$.
    \item If $u$ is a normal epimorphism and $v$ is a monomorphism, then $vu$ is the normal epi factorization of $f$.
  \end{thmlist}
\end{lemma}
\begin{proof}
  (i)\quad The universal property of the normal mono factorization of $f$ yields this commutative diagram:
  \begin{equation*}
    \xymatrix@R=5ex@C=4em{
    X \ar[r]^-{u'} \ar@{=}[d] &
    K \ar@{{ |>}->}[r]^-{m} \ar@{.>}[d]^{\alpha} &
    Y \ar@{=}[d] \\
    X \ar@{-{>>}}[r]_-{u} &
    L \ar@{{ |>}->}[r]_-{v} &
    Y
    }
  \end{equation*}
  With (\ref{thm:IsomorphismRecognition}), we see that $\alpha$ is an isomorphism because (a) commutativity of the square on the right shows (\ref{thm:NormalMono-Props}) that $\alpha$ is a normal monomorphism, while (b) commutativity of the square on the left shows (\ref{exe:Epimorphisms-Composite}) that $\alpha$ is an epimorphism.
\end{proof}

We single out three special cases involving the normal epi factorization and the normal mono factorization of a given morphism.

\begin{definition}[Image factorization\ZExactTag]
  \label{def:ImageFactorization}%
  A morphism $f$ admits an \Defn{image factorization} if it can be written as a composite $f=v e$, with $e$ a normal epimorphism and $v$ a monomorphism.
  \begin{equation*}
    \xymatrix@R=5ex@C=4em{
    \DiagObj \ar@/^2ex/[rr]^-{f} \ar@{-{ >>}}[r]_-{e} &
    \Img{f} \ar@{{ >}->}[r]_-{v} &
    \DiagObj
    }
  \end{equation*}
  In this situation, the domain of $v$ is called the \Defn{image} of $f$. %
  \index{image!in a normal epi factorization}%
  \index[not]{i!$\ImgMap{f}$\IndSep image map of $f$}\index[not]{i!$\Img{f}$\IndSep domain of image map of $f$}%
\end{definition}

In (\ref{def:ImageFactorization}), the term image factorization generalizes what we know about functions of sets: Every function has an essentially unique factorization into `surjection' followed by `injection'. The situation here is analogous: Whenever $f$ admits an image factorization as $v e$, then (\ref{thm:Initial/TerminalDecomposition-Recognize}) tells us that it is a special kind of normal epi factorization. It is available if and only if the comparison map $\NENMComp{f}$ is monic.

Less commonly used is the dual `coimage factorization' of `image factorization': An epimorphism followed by a normal monomorphism. In this case, the normal mono factorization $f=mu$ the map $u$ happens to be an epimorphism. This occurs if and only if the comparison map $\NENMComp{f}$ is epic. Thus $\CoImgMap{f}\DefEq u$ is called the \Defn{coimage} of $f$, and $\CoImg{f}$ is the codomain of $\CoImgMap{f}$. %
\index{coimage!in a normal mono factorization}%
\index[not]{c!$\CoImgMap{f}$\IndSep coimage map of $f$}\index[not]{c!$\CoImg{f}$\IndSep codomain of coimage map of $f$}%

Given a map $f$, we are particularly interested in the special case where a normal mono factorization and a normal epi factorization coincide or, more precisely, where the normal epi / normal mono factorization comparison map $\NENMComp{f}$ is an isomorphism:

\begin{definition}[Normal morphism\ZExactTag]%
  \label{def:NormalMap}%
  A morphism $f\from X\to Y$ is called \Defn{normal}, also \Defn{proper}, if it admits a decomposition $f=me$ in which $m$ is a normal  monomorphism and $e$ is a normal epimorphism.%
  \index{proper!morphism}\index{normal!map}\index{morphism!normal}%
  \begin{equation*}
    \xymatrix@R=5ex@C=3em{
    X \ar[rr]^-{f} \ar@{-{ >>}}[rd]_{e} &&
    Y \\
    & I \ar@{{ |>}->}[ru]_{m}
    }
  \end{equation*}
\end{definition}

With (\ref{thm:Initial/TerminalDecomposition-Recognize}) we see that, for a normal morphism, its normal mono factorization coincides with its normal epi factorization. In particular, the decomposition $f=me$, with $m=\KerMap{\CoKerMap{f}}$ and $e=\CoKerMap{\KerMap{f}}$, is unique up to unique isomorphism. Thus we may refer to the decomposition $f=me$ as the \Defn{normal factorization} of $f$. --- We develop basic properties of normal maps. %
\index{normal!factorization}%

\begin{proposition}[Recognizing the normal factorization\ZExactTag]
  \label{thm:NormalFactorization-Recognize}%
  Suppose a normal map $f$ admits a decomposition $f=np$ with $n$ monic and $p$ epic. Then $n$ is a normal monomorphism, and $p$ is a normal epimorphism. %
  \index{normal!factorization - recognize}%
\end{proposition}
\begin{proof}
  In the diagram below, the bottom row is the normal factorization of $f$.
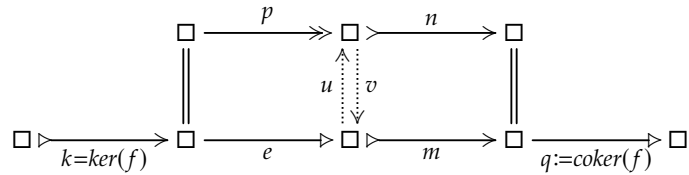
  \begin{equation*}
    \xymatrix@R=5ex@C=4em{
    & \DiagObj \ar@{-{>>}}[r]^-{p} \ar@{=}[d] &
    \DiagObj \ar@{{ >}->}[r]^-{n} \ar@{.>}@<0.5ex>[d]^{v} &
    \DiagObj \ar@{=}[d] \\
    \DiagObj \ar@{{ |>}->}[r]_-{k=\KerMap{f}} &
    \DiagObj \ar@{-{ >>}}[r]_-{e} &
    \DiagObj \ar@{{ |>}->}[r]_-{m} \ar@{.>}@<0.5ex>[u]^{u}&
    \DiagObj \ar@{-{ >>}}[r]_-{q\DefEq\CoKerMap{f}} &
    \DiagObj
    }
  \end{equation*}
  We know that $e=\CoKerMap{k}$, and that $npk = fk=0$. With the monic property of $n$ we conclude that $pk=\ZeroMap$. There exists $u$, unique with $p=ue$ because $e=\CoKerMap{k}$. Dually, $qnp=qf=\ZeroMap$. With the epic property of $p$ we conclude that $qn=\ZeroMap$. There exists $v$, unique with $n=mv$ because $m=\KerMap{q}$.

  The claim will follow once we show that $u$ and $v$ are inverses of each other. Indeed, using the monic and epic properties of $m$ and $e$,
  \begin{equation*}
    mvue = mvp = np = me \quad \text{and so}\quad e=vue\quad \text{and so}\quad vu=\IdMap
  \end{equation*}
  Similarly,
  \begin{equation*}
    nuvp = mvuvp = mvp = np
  \end{equation*}
  Using the monic and epic properties of $n$ and $p$, $uv=\IdMap$, and the proof is complete.
\end{proof}

\begin{proposition}[Factoring  normal map through a subobject of its codomain\ZExactTag]
  \label{thm:FactorNormalThroughSubObject->Normal}%
  If a normal map $f\from X\to Y$ is a composite $f=\mu f'$ with $\mu$ a monomorphism, then $f'$ is a normal map. %
  \index{normal!map factored through monic map}%
\end{proposition}
\begin{proof}
  We reason using this diagram:
  
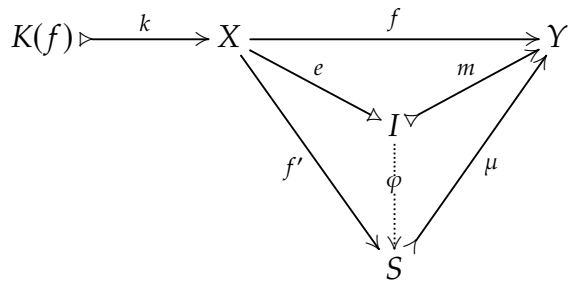
\begin{equation*}
    \xymatrix@R=3ex@C=4em{
    \Ker{f} \ar@{{ |>}->}[r]^-{k} &
    X \ar[rr]^-{f} \ar@{-{ >>}}[rd]^{e} \ar[rddd]_{f'}&&
    Y \\
    && I \ar@{{ |>}->}[ru]^{m} \ar@{.>}[dd]|(0.4){\varphi} & \\ \\
    && S \ar@{{ >}->}[ruuu]_{\mu}
    }
  \end{equation*}
  First note that $k=\Ker{f'}$ because $\mu$ is monic; see (\ref{thm:NormalMono-Props}.ii).  Then we know from (\ref{thm:NormalDecompositions-Initial/Terminal}) that $e=\CoKer{k}$. So there is $\varphi\from I\to S$, unique with $f'=\varphi e$. The identity $m=\mu\varphi$ follows from the epic property of $e$ and the computation
  \begin{equation*}
    \mu\varphi e = \mu f' = f = me
  \end{equation*}
  With (\ref{thm:NormalMono-Props}.iv) we infer that $\varphi$ is a normal monomorphism. Thus the composite $\varphi e$ is a normal decomposition of $f'$; i.e., $f'$ is a normal map.
\end{proof}

In other words, (\ref{thm:FactorNormalThroughSubObject->Normal})  says that, whenever in a normal composite  $X\XRA{f} Y \XRA{g} Z$ the map $g$ is a monomorphism, then $f$ is a normal map. - Dually:

\begin{proposition}[Factoring a normal map through quotient of domain\ZExactTag]
  \label{thm:FactorNormalQuotientObject->Normal}
  If a normal map $f\from X\to Y$ is a composite $f= \bar{f}\varepsilon$ with $\varepsilon$ an epimorphism, then $\bar{f}$ is a normal map.%
  \index{normal!map factored through epic map}%
  \NoProof
\end{proposition}

By using the universal properties of kernels and cokernels, we see that normal factorizations are functorial:

\begin{proposition}[Normal factorizations are functorial\ZExactTag]
  \label{thm:NormalFactorization-Functorial}%
  If in a commutative square of morphisms, as on the left, the horizontal arrows are normal, then there is a unique map which renders their normal factorizations commutative, as on the right below.
  
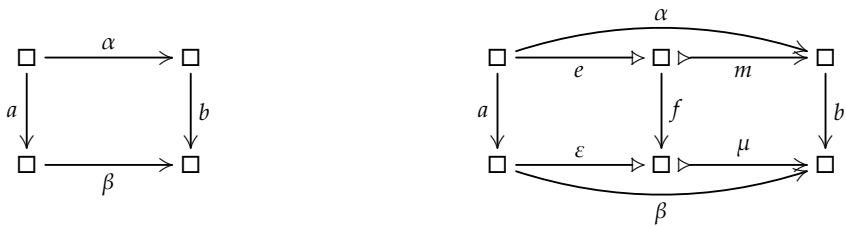
\begin{equation*}
    \xymatrix@R=5ex@C=4em{
    \DiagObj \ar[r]^-{\alpha} \ar[d]_{a}&
    \DiagObj \ar[d]^{b} &&
    \DiagObj \ar@{-{ >>}}[r]_-{e} \ar[d]_{a} \ar@/^2ex/[rr]^-{\alpha} &
    \DiagObj \ar@{{ |>}->}[r]_-{m} \ar[d]^{f} &
    \DiagObj \ar[d]^{b} \\
    \DiagObj \ar[r]_-{\beta} &
    \DiagObj &&
    \DiagObj \ar@{-{ >>}}[r]^-{\varepsilon} \ar@/_2ex/[rr]_-{\beta} &
    \DiagObj \ar@{{ |>}->}[r]^{\mu} &
    \DiagObj \\
    }
  \end{equation*}
\end{proposition}

We now consider decompositions of order $2$. Initially perhaps surprising, such decompositions occur `naturally'.

\begin{definition}[Subnormal decompositions\ZExactTag]
  \label{def:SubnormalDecomposition}%
  A \Defn{subnormal mono decomposition} of a morphism $f$ is given by a sequence $f=m_{1}m_{2}u$ in which $m_{1}$ and $m_{2}$ are normal monomorphisms. Dually, a \Defn{subnormal epi decomposition} of $f$ is given by a sequence $f=ve_{2}e_{1}$ in which $e_{1}$ and $e_{2}$ are normal epimorphisms. %
  \index{subnormal!mono decomposition}\index{subnormal!epi decomposition}%
\end{definition}

The subnormal mono decompositions of $f$ contain an initial object which is obtained in two steps, namely
\begin{enumerate}
  \item Into the normal mono factorization $f=m_{1}u_{1}$ of $f$, insert
  \item  the normal mono factorization $f=m_{1}(m_{2}u_{2})$ of $u_{1}$.
\end{enumerate}
Dually, the subnormal epi decompositions of $f$ contain a terminal object $f=(v_{2}e_{2})e_{1}$  constructed from the normal epi factorization $f=v_{1}e_{1}$ of $f$.

We call $f$ \Defn{subnormal} if it admits a subnormal mono decomposition $f=m_{1}m_{2}e$ in which $e$ is a normal epimorphism. We call $f$ \Defn{cosubnormal} if it admits a subnormal epi decomposition $f=me_{2}e_{1}$ in which $m$ is a normal monomorphism.%
\index{subnormal!morphism}\index{cosubnormal!morphism}\index{morphism!subnormal}\index{morphism!cosubnormal}%

For instance, in the category $\Grps$ of groups, any composite of two normal subgroup inclusions is a subnormal map, while any cosubnormal map is automatically normal. In general, we have:

\begin{proposition}[Subnormal and cosubnormal implies normal\ZExactTag]
  \label{thm:Subnormal+Cosubnormal->Normal}%
  A morphism which is both subnormal and cosubnormal is normal.
\end{proposition}
\begin{proof}
  Assume that $me_{2}e_{1}=f=m_{1}m_{2}e$ factors as $me_{2}e_{1}$, where $e,e_{1},e_{2}$ are normal epimorphisms and $m,m_{1},m_{2}$ are normal monomorphisms. Thus we have this commutative diagram:
  
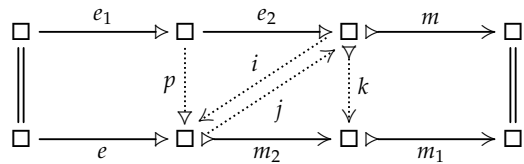
\begin{equation*}
    \xymatrix@R=5ex@C=4em{
    \DiagObj \ar@{=}[d] \ar@{-{ >>}}[r]^-{e_1} &
    \DiagObj \ar@{.{ >>}}[d]_-{p} \ar@{-{ >>}}[r]^-{e_2} &
    \DiagObj \ar@{{ |>}.>}[d]^{k} \ar@{.>}@<-0.5ex>[ld]_{i}  \ar@{{ |>}->}[r]^-{m} &
    \DiagObj \ar@{=}[d] \\
    \DiagObj  \ar@{-{ >>}}[r]_-{e} &
    \DiagObj \ar@{{ |>}->}[r]_-{m_{2}} \ar@{.>}@<-0.5ex>[ru]_{j} &
    \DiagObj \ar@{{ |>}->}[r]_-{m_{1}} &
    \DiagObj
    }
  \end{equation*}
  The composite $m(e_{2}e_{1})$ is the normal mono factorization of $f$. So, there exists $k$ unique with $m=m_{2}k$ and $ke_{2}e_{1}=m_{2}e$. Then $k$ is a normal monomorphism because the right hand square commutes. Also, $(m_{1}m_{2})e$ is the normal epi factorization of $f$. So, there exists $p$ unique with $e=pe_{1}$ and $m_{1}m_{2}p=me_{2}$. Then $p$ is a normal epimorphism because the left hand square commutes. Via the monic property of $m_{1}$ or the epic property of $e_{1}$, we see that the middle square commutes as well.

  Now recall that normal epimorphisms are strong. So, the center square has two unique fillers $i$ and $j$ which render the resulting triangles commutative. It follows that $j=i^{-1}$, implying that the composite $e_{2}e_{1}$ is a normal epimorphism and that $m_{1}m_{2}$ is a normal monomorphism. Thus $f$ is a normal map.
\end{proof}

We close this section by clarifying how to compute (co)limits in the categories $\NMonoCat{X}$ and $\NEpiCat{X}$. Here, we rely on the general principles governing limits and colimits in reflective, respectively coreflective, subcategories; see Section \ref{sec:ReflectiveSubCats},

\begin{proposition}[(Co)limits in $\NMonoCat{X}$ and in $\NEpiCat{X}$\ZExactTag]
  \label{thm:(Co)Limits-NMono(X)/NEpi(X)}%
  In a \ZExact\ category $\Ctgry{X}$ the following hold:
  \begin{thmlist}
    \item If $\Ctgry{X}$ admits limits over a small category $D$, then so does $\NMonoCat{X}$, and these are computed object-wise. %
    \index{limit!in $\NMonoCat{X}$}%
    \item If $\Ctgry{X}$ admits colimits over a small category $D$, then so does $\NMonoCat{X}$. To compute the colimit of a diagram $\Phi\from D\to \NMonoCat{X}$ of normal monomorphisms, compute the colimit of $\Phi$ in $\ArrowCat{\Ctgry{X}}$, then apply the reflector $\ArrowCat{\Ctgry{X}}\to \NMonoCat{X}$. %
    \index{colimit!in $\NMonoCat{X}$}
    \item If $\Ctgry{X}$ admits colimits over a small category $D$, then so does $\NEpiCat{X}$, and these are computed object-wise. %
    \index{colimit!in $\NEpiCat{X}$}%
    \item If $\Ctgry{X}$ admits limits over a small category $D$, then so does $\NEpiCat{X}$. To compute the limit of a diagram $\Psi\from D\to \NEpiCat{X}$ of normal monomorphisms, compute the limit of $\Psi$ in $\ArrowCat{\Ctgry{X}}$, then apply the coreflector $\ArrowCat{\Ctgry{X}}\to \NEpiCat{X}$. %
    \index{limit!in $\NEpiCat{X}$}
  \end{thmlist}
\end{proposition}
\begin{proof}
  (i) is true because kernels commute with limits. Dually, (iii) is true because cokernels commute with colimits.

  (ii) is an application of (\ref{thm:ReflectiveSubCat-Colimits}), and (iv) is its dual.
\end{proof}

\begin{proposition}[Normal mono in $\NMonoCat{X}$ / normal epi in $\NEpiCat{X}$ \ZExactTag]
  \label{thm:KernelsInNMono(X)-CoKernelsInNEpi(X)}
  In a \ZExact\ category $\Ctgry{X}$ the following hold:
  \begin{thmlist}
    \item A morphism in $\NEpiCat{X}$ is a normal epimorphism if and only if its underlying square in $\Ctgry{X}$ is a pushout of normal epimorphisms. %
    \index{normal!epimorphism in $\NEpiCat{X}$}%
    \item A morphism in $\NMonoCat{X}$ is a normal monomorphism if and only if its underlying square in $\Ctgry{X}$ is a pullback of normal monomorphisms. %
    \index{normal!monomorphism in $\NMonoCat{X}$}%
  \end{thmlist}
\end{proposition}
\begin{proof}
  We prove (ii). In the commutative diagram below, assume that the left hand square represents the kernel of the right hand square in $\NMonoCat{X}$.
  \begin{equation*}
    \xymatrix@R=5ex@C=3em{
    K \PullLU{rd} \ar[r]^-{\alpha} \ar@{{ |>}->}[d]_{k}  &
    A \ar[r]^-{q} \ar@{{ |>}->}[d]_{a} &
    B \ar@{{ |>}->}[d]^{b} \\
    K' \ar[r]_-{\alpha'} &
    A' \ar[r]_-{q'} &
    B'
    }
  \end{equation*}
  Then $\alpha=\KerMap{q}$ in $\Ctgry{X}$, and $\alpha'=\KerMap{q'}$. By (\ref{thm:PullbackRecognition-KernelSide-1}), the left hand square is a pullback.

  Conversely, suppose the square on the left is a pullback of normal monomorphisms. Construct the right hand square as the cokernel of $(\alpha,\alpha')$ in $\NMonoCat{X}$. Then $q'=\CoKerMap{\alpha'}$ and, hence, $\alpha'=\KerMap{q'}$. We show directly that $(\alpha,\alpha')=\KerMap{q,q'}$. Indeed, if $(x,x')\from m\to a$ is such that $(q,q')(x,x')=(qx,q'x')=(\ZeroMap,\ZeroMap)$, then $x'$ factors uniquely through $\alpha'$ via $t'$. The pullback property of the left hand square yields a unique factorization of $x$ through $\alpha$ via $t$. It follows that $(t,t')$ is the required unique factorization of $(x,x')$ through $(\alpha,\alpha')$ in $\NMonoCat{X}$.
\end{proof}

\begin{subordinate}{On Puppe-exact categories}
  A \Defn{Puppe-exact} or \Defn{p-exact} category is a \ZExact\ category in which all morphisms are normal; see~\cite{Grandis-HA1,Grandis-HA2} and the citations there to the foundational work of Puppe, Mitchell and others. Whenever a p-exact category admits binary products or binary coproducts, then it is abelian; see~\cite{Alligators} for an explicit proof of this fact.
\end{subordinate}

\begin{exercises}

\bigskip
\begin{definition}[Subnormal decompositions of order $r$\ZExactTag]
  \label{def:SubnormalDecomposition-Order-r}%
  A \Defn{subnormal mono decomposition of order $r$} of $f\from X\to Y$ is a diagram %
  \index{normal!mono decomposition of order $r$}%
  \begin{equation*}
    X \XRA{u} K_{r} \overset{m_{r}}{\NMono} \cdots \overset{m_{1}}{\NMono} Y
  \end{equation*}
  in which $m_{r},\dots , m_{1}$ are normal monomorphisms. Dually, a \Defn{subnormal epi decomposition of order $r$} of $f$ is a diagram %
  \index{normal!epi decomposition of order $r$}%
  \begin{equation*}
    X \overset{e_{1}}{\NEpi} \cdots \overset{e_{r}}{\NEpi} Q_{r} \XRA{v} Y
  \end{equation*}
  in which $e_{1},\dots ,e_{r}$ are normal epimorphisms.
\end{definition}

\begin{exercise}[Initial / terminal subnormal decompositions\ZExactTag]
  \label{exe:SubnormalDecomposition-Initial/Terminal}
  Given a map $f$ in a pointed category, show that iterating the constructions in Lemma \ref{thm:NormalDecompositions-Initial/Terminal} yields:
  \begin{thmlist}
    \item A normal mono decomposition of order $r$ of $f$, which is initial among all such decompositions, and
    \item a normal epi decomposition of order $r$ of $f$, which is terminal among all such decompositions.
  \end{thmlist}
\end{exercise}

\begin{definition}[Subnormal / cosubnormal map of order $r$\ZExactTag]
  \label{def:SubNormal/CoSubNormal-Order-r}
  We say that $f$ is \Defn{subnormal of order $r$} if it admits a subnormal mono decomposition of order $r$ in which the map $u$ is a normal epimorphism. Dually, $f$ is \Defn{cosubnormal of order $r$} if it admits a subnormal epi decomposition of order $r$ in which the map $v$ is a normal monomorphism.
\end{definition}

\begin{exercise}[Subnormal and cosubnormal of order $r$ implies normal]
  \label{exe:SubNormal/CoSubNormal-r-=>Normal}%
  A map $f$ which is both subnormal and cosubnormal of order $r$ is normal.
\end{exercise}

\begin{exercise}[Solvable group via subnormal factorization]
  \label{exe:SolvableGroup-SubNormalFactorization}
  In the category $\Grps$ of groups, show that a group $G$ is solvable if and only if there is $r\geq 1$ for which the zero map $\ZeroObject\to G$ admits a subnormal decomposition of order $r$ for which the quotient groups $Q_{i}\DefEq N_{i}/N_{i+1}$ in the diagram below are abelian.
  \begin{equation*}
    \xymatrix@R=5ex@C=4em{
    0 \ar@{{ |>}->}[r] &
    N_{r} \ar@{{ |>}->}[r]^-{m_{r}} \ar@{-{ >>}}[d] &
    N_{r-1} \ar@{{ |>}->}[r]^-{m_{r-1}} \ar@{-{ >>}}[d] &
    \cdots \ar@{{ |>}->}[r]^-{m_{2}} &
    N_{1} \ar@{{ |>}->}[r]^-{m_{1}} \ar@{-{ >>}}[d]  &
    G \ar@{-{ >>}}[d]  \\
    & Q_{r} &
    Q_{r-1} &&
    Q_{1} &
    Q_{0}
    }
  \end{equation*}
\end{exercise}

\begin{exercise}[Nilpotent group via cosubnormal factorization]
  \label{exe:NilpotentGroup-CoSubNormalFactorization}%
  In the category $\Grps$ of groups, show that a group $G$ is nilpotent of order less than or equal to $r$ if and only if the zero map $G\to \ZeroObject$ admits a cosubnormal decomposition of order $r$ for which the kernel groups $K_{i}\DefEq \Ker{Q_{i}\to Q_{i+1}}$ in the diagram below are central in $Q_{i}$: %
  \index{nilpotent!group via cosubnormal factorization}%
  \begin{equation*}
    \xymatrix@R=5ex@C=4em{
    C_{0} \ar@{{ |>}->}[d] &
    C_{1} \ar@{{ |>}->}[d] &&
    C_{r-1} \ar@{{ |>}->}[d] &
    C_{r} \ar@{=}[d] \\
    G \ar@{-{ >>}}[r]_-{e_{1}} &
    Q_{1} \ar@{-{ >>}}[r]_-{e_{2}} &
    \cdots \ar@{-{ >>}}[r]_-{e_{r-2}} &
    Q_{r-1} \ar@{-{ >>}}[r]_-{e_{r-1}}&
    Q_{r} \ar@{-{ >>}}[r]_-{e_{r}} &
    0
    }
  \end{equation*}
\end{exercise}

\begin{exercise}[Colimits in $\NMonoCat{X}$, alternate computation\ZExactTag]
  \label{exe:CoLimsInNMono(X)-Alternat}%
  In a pointed category $\Ctgry{X}$, show that the colimit of a diagram $\Phi\from D\to \NMonoCat{X}$ may be computed as:
  \begin{equation*}
    \CoLimOfOver{\Phi}{D} \cong \KerFunc \left(\CoLimOfOver{\CoKerFunc\Comp \Phi)}{D}\right)
  \end{equation*}
  In detail:
  \begin{enumerate}
    \item Let $\Psi\DefEq\CoKerFunc\Comp \Phi\from D\to \NEpiCat{X}$ be the diagram in $\NEpiCat{X}$ obtained by applying normal subobject / quotient inversion to the diagram $\Phi$ in $\NMonoCat{X}$.
    \item Compute $\CoLimOf{\Psi}$ in $\NEpiCat{X}$ by taking object-wise colimits.
    \item Apply normal subobject / quotient inversion to the result of (ii).
  \end{enumerate}
  Dualized the above to a computation of the limit of a diagram $\Psi\from D\to \NEpiCat{X}$.
\end{exercise}
\end{exercises}
\section[Exact Sequences]{Exact Sequences}
\label{sec:ExactSeqs}

In this section we introduce exact sequences and we develop their basic properties. So, what say here applies even in low structured settings such as the category $\SetsBsd$ of pointed sets or the category $\Magmas$ of unital magmas.

\begin{definition}[Short exact sequence/extension]
  \label{def:ShortExactSequence-Basic}%
  A \Defn{short exact sequence} in a pointed category is a diagram%
  \index{sequence!short exact}%
  \index{short exact sequence}%
  \index{extension!of $K$ by $Q$}\index{cover!of $Q$ by $K$}%
  \begin{equation*}
    \xymatrix@R=5ex@C=2em{
    E &&
    K \ar@{{ |>}->}[rr]^-{k} &&
    X \ar@{-{ >>}}[rr]^-{q} &&
    Q
    }
  \end{equation*}
  where $k$ is a kernel of $q$, and $q$ is a cokernel of $k$. We also say that $E$ is an \Defn{extension of $K$ over~$Q$}, or a \Defn{cover of $Q$ by $K$}.
\end{definition}

In a short exact sequence $E$, as in (\ref{def:ShortExactSequence-Basic}), we sometimes write $Q$ as the quotient $X/K$. Such a sequence provides a tool for the structural analysis of objects: In a z-exact category, where kernels and cokernels are always available, the middle object $X$, together with either of the end objects determines the other end object up to isomorphism. --- When dealing with short exact sequences, we generally assume that kernels and cokernels exist. With (\ref{exe:Sub/QuotientOf-0}) we see:

\begin{example}[Extension with middle object $\ZeroObject$ \ZExactTag]
  \label{exa:Extending-0-By-Epi}
  A sequence $K\to 0\to Q$ is short exact if and only if $K=0=Q$. \NoProof
\end{example}

\begin{example}[Extension by $0$'s\ZExactTag]
  \label{exa:Extension-By-0's}%
  A sequence $0\to X\to 0$ is short exact if and only if  $X=0$. \NoProof
\end{example}

\begin{example}[Constructing a short exact sequence from a normal mono/epi\ZExactTag]
  \label{exa:ConstructingSESFromNormalEpi/Mono}
  If $f\from X\NMono Y$ is a normal monomorphism, then the sequence on the left below is short exact. If $g\from Y\NEpi Z$ is a normal epimorphism, then the sequence on the right is short exact.
  \begin{equation*}
    \xymatrix@R=5ex@C=3em{
    X	 \ar@{{ |>}->}[r]^-{f} &
    Y \ar@{-{ >>}}[r] &
    \CoKer{f} &
    \Ker{g} \ar@{{ |>}->}[r] &
    Y \ar@{-{ >>}}[r]^-{g} &
    Z
    }
  \end{equation*}
\end{example}
\begin{proof}
  This is so because a normal monomorphism is the kernel of its cokernel, respectively, a normal epimorphism is the cokernel of its kernel; see (\ref{thm:Ker(CoKer)-CoKer(Ker)}).
\end{proof}

More generally:

\begin{definition}[Exactness of composable maps\ZExactTag]
  \label{def:ExactnessComposableMaps}
  A sequence of composable maps
  \begin{equation*}
    A \XRA{f} B \XRA{g} C
  \end{equation*}
  is \Defn{exact (in position $B$)} if $f$ admits a normal factorization $f=me$ such that $m=\KerMap{g}$. %
  It is \Defn{coexact (in position $B$)} if $g$ admits a normal factorization $g=kp$ such that $p=\CoKerMap{f}$. %
  \index{exact!sequence}
\end{definition}

\begin{lemma}[Recognizing exactness\ZExactTag]
  \label{thm:Exactness-Recognize}
  For composable \emph{normal} morphisms $g$ following $f$ the conditions below are equivalent. %
  \index{exact!node recognize}%
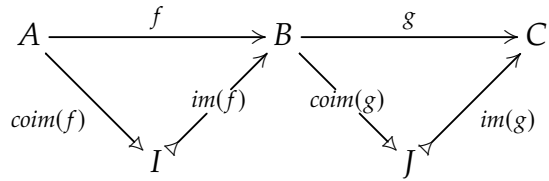
  \begin{equation*}
    \xymatrix@!0@=4em{
    A \ar[rr]^-{f} \ar@{-{ >>}}[rd]_-{\CoImgMap{f}} &&
    B \ar[rr]^-{g} \ar@{-{ >>}}[rd]|-{\CoImgMap{g}} &&
    C \\
    & I \ar@{{ |>}->}[ru]|-{\ImgMap{f}}
    && J \ar@{{ |>}->}[ru]_-{\ImgMap{g}}
    }
  \end{equation*}
  \begin{tfae}
    \item The sequence is exact in position $B$; i.e., $\ImgMap{f}=\KerMap{g}$.
    \item The sequence is coexact in position $B$; i.e., $\CoKerMap{f}= \CoImgMap{g}$
    \item The sequence $I \NMono B \NEpi J$ is short exact. \NoProof
  \end{tfae}
\end{lemma}

\begin{corollary}[Recognizing a short exact sequence\ZExactTag]
  \label{thm:ShortExactSequence-Recognize}\label{Def:ShortExactSequence}%
  A sequence of normal morphisms of the form %
  \index{exact!short exact}%
  \index{short exact sequence}%
  \begin{equation*}
    \text{(S)}\qquad
    \xymatrix@R=5ex@C=3em{
    0 \ar[r] &
    A \ar@{->}[r]^-{f} &
    B \ar@{->}[r]^-{g} &
    C \ar[r] &
    0
    }
  \end{equation*}
  is exact if and only if the sequence $A\XRA{f}B \XRA{g}$ is short exact. \NoProof
\end{corollary}

Extending Definition \ref{def:ExactnessComposableMaps}, we say that a sequence of finite length
\begin{equation*}
  X_0\XRA{f_0} X_1 \XRA{f_1} \cdots \to X_k\XRA{f_k} \cdots \to X_{n}\XRA{f_n} X_{n+1}
\end{equation*}
is \Defn{exact} if, for every $1\leq k\leq n$, the sequence $X_{k-1}\to X_k\to X_{k+1}$ is exact. This implies that the maps $f_0,$ \dots , $f_{n-1}$ are normal maps.

Dually, the sequence is \Defn{coexact} if, for every $1\leq k\leq n$, the sequence $X_{k-1}\to X_k\to X_{k+1}$ is coexact.   This implies that the maps $f_1$, \dots , $f_{n}$ are normal maps. Extending to infinite sequences:

\begin{definition}[Long exact sequence\ZExactTag]
  \label{def:LES}
  A sequence of morphisms $\cdots \to X_{n+1}\to X_n\XRA{f_n} X_{n-1}\to \cdots $ is called \Defn{long exact} if it can be spliced together from short exact sequences $Z_n \overset{k_n}{\NMono} X_n \overset{q_n}{\NEpi} Z_{n-1}$ as in Figure~\ref{Fig:ExactViaSplice}. %
  \index{long exact sequence}
  \begin{figure}[H]
    \[
      \xymatrix@!0@=3.5em{
      &&&& Z_n \ar@{{ |>}->}@/^1.5ex/[rd]^-{k_n} &&&&
      Z_{n-2} \ar@{{ |>}->}@/^1.5ex/[rd] \\
      \cdots \ar[r] &
      X_{n+2} \ar@{-{ >>}}@/_1.5ex/[rd] \ar[rr]^-{f_{n+2}} &&
      X_{n+1} \ar@{-{ >>}}@/^1.5ex/[ru]^-{q_{n+1}} \ar[rr]_-{f_{n+1}} &&
      X_n \ar@{-{ >>}}@/_1.5ex/[rd]^-{q_{n}} \ar[rr]^-{f_{n}} &&
      X_{n-1} \ar[rr]_-{f_{n-1}} \ar@{-{ >>}}@/^1.5ex/[ru]^-{q_{n-1}} &&
      X_{n-2} \ar[r] &
      \cdots \\
      && Z_{n+1} \ar@{{ |>}->}@/_1.5ex/[ru]^-{k_{n+1}} &&&&
      Z_{n-1} \ar@{{ |>}->}@/_1.5ex/[ru]^-{k_{n-1}}}
    \]
    \caption{An exact sequence admits splicing by short exact sequences}\label{Fig:ExactViaSplice}
  \end{figure}
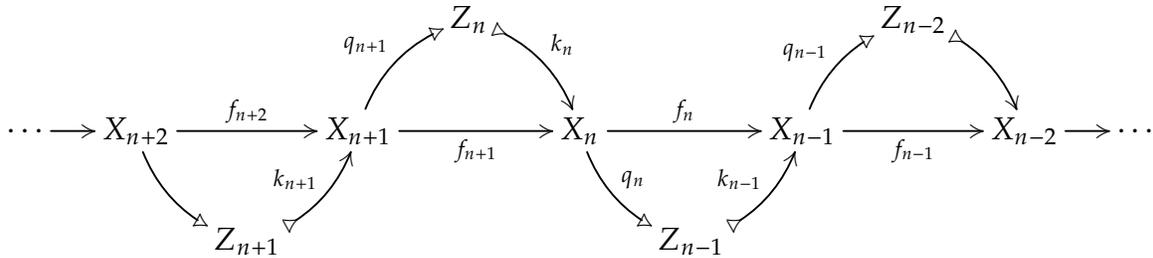
\end{definition}

This means that every morphism in a long exact sequence admits a normal factorization, and factorizations of adjacent morphisms form short exact sequences as shown. Note that the sequence of morphisms in (\ref{def:LES}) is exact if and only if it is coexact (i.e., exact in the opposite category). For historic reasons we favor using the term `exact'. With (\ref{thm:NormalMono-Props}.ii) and (\ref{thm:NormalEpi-Props}.ii) we see that
\begin{equation*}
  k_n=\KerMap{q_n}= \KerMap{f_n}  \qquad \text{and}\qquad q_n=\CoKerMap{k_n} = \CoKerMap{f_{n+1}}
\end{equation*}
Whenever a sequence of morphisms admits a decomposition as in Figure \ref{Fig:ExactViaSplice}, then this decomposition is unique up isomorphism, in agreement with Section  \ref{sec:NormalDecompositionsFactorizations}.

\begin{subordinate}{}

  \begin{subsubordinate}{Use of the term `extension'}
    The use of the term `extension' is conflicted in the literature. Some authors refer to a short exact sequence $K\overset{k}{\NMono} X \overset{q}{\NEpi} Q$ as an extension of $Q$.
  \end{subsubordinate}

  \begin{subsubordinate}{First use of the term `exact'}
    Curiously, it seems that exactness of a sequence of abelian groups was only used more than a decade after the introduction of homology of a chain complex; see a talk by Hurewicz \cite{WHurewicz1941-DualityAbstract}. This is according to the Kelley-Pitcher \cite[p.~682]{JLKelleyEPitcher1947} who appear to be the first to use the term `exact sequence' \cite[Def. 3.1]{JLKelleyEPitcher1947}, but see also what Hilton \cite[p.~291]{PHilton1988-HomologyHistory} explains regarding J.\ H.\ C.\ Whitehead's work.
  \end{subsubordinate}

\end{subordinate}

\begin{exercises}

\begin{exercise}[Factoring a morphism of zero maps\ZExactTag]
  \label{exe:FactorMorphisZeroMaps}
  \cite[p.~272ff]{FBorceuxDBourn2004}\quad In the diagram below, assume that the front face commutes, and that the bottom right square is a pullback. %
  \index{morphism!of zero maps - factorization}
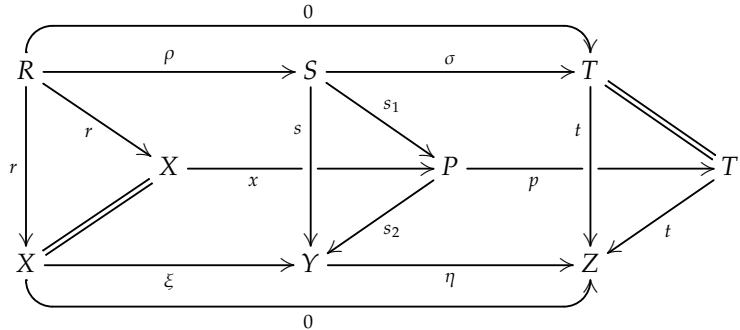
  \begin{equation*}
    \xymatrix@R=5ex@C=4em{
    R \ar[rr]^-{\rho} \ar `u[rrrr]`[rrrr]^{\ZeroMap}[rrrr] \ar[dd]_{r} \ar[rd]_{r} &&
    S \ar[rr]^-{\sigma} \ar[dd]_(0.3){s} \ar[rd]^{s_1} &&
    T \ar[dd]_(0.3){t} \ar@{=}[rd] \\
    & X \ar@{=}[ld] \ar[rr]|\hole_(0.3){x} &&
    P \ar[ld]^{s_2} \ar[rr]|\hole _(0.3){p} &&
    T \ar[ld]^{t} \\
    X \ar[rr]_-{\xi} \ar@{-<}`d[rrrr]`[rrrr]_{\ZeroMap}[rrrr] &&
    Y \ar[rr]_-{\eta} &&
    Z
    }
  \end{equation*}
  Then the following hold:
  \begin{thmlist}
    \item There exist unique maps $s_1$ and $x$ rendering the entire diagram commutative.
    \item If $\xi=\KerMap{\eta}$, then $x=\KerMap{p}$.
    \item If $\xi=\KerMap{\eta}$ and $\rho=\KerMap{\sigma}$, then the upward facing left hand square is a pullback.
    \item If $\xi=\KerMap{\eta}$ and the upper sequence is short exact, then $X\XRA{x} P \XRA{p} T$ is short exact.
  \end{thmlist}
\end{exercise}

\begin{exercise}[Short exact sequence in $\SetsBsd$]
  \label{exe:SES-In-Set_*}
  In the category $\SetsBsd$ of pointed sets show that, in a short exact sequence $N\overset{i}{\NMono} X \overset{q}{\NEpi} Q$ the following hold: %
  \index{short exact sequence!in $\SetsBsd$}
  \begin{thmlist}
    \item $q$ has a unique section $s\from Q\to X$, and
    \item $i$ and $s$ form a cocone for $X$ as the coproduct of $N$ and $Q$.
  \end{thmlist}
\end{exercise}

\begin{exercise}[Factoring the cokernel of a map\ZExactTag]
  \label{exe:FactoringCoKer}
  Given composable maps  $f\from X\to Y$ and $g\from Y\to Z$, show that the diagram of below has exact rows and columns and commutes.
  
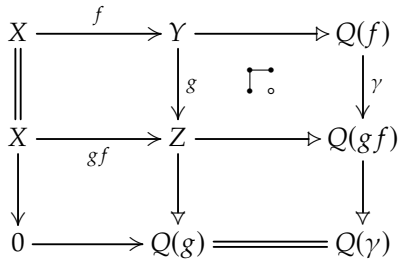
\begin{equation*}
    \xymatrix@R=5ex@C=4em{
    X \ar[r]^-{f} \ar@{=}[d] &
    Y \ar@{-{ >>}}[r] \ar[d]^{g} \PushRD{rd} &
    \CoKer{f} \ar[d]^{\gamma}\\
    X \ar[r]_-{gf} \ar[d] &
    Z \ar@{-{ >>}}[r] \ar@{-{ >>}}[d] &
    \CoKer{gf} \ar@{-{ >>}}[d] \\
    0 \ar[r] &
    \CoKer{g} \ar@{=}[r] &
    \CoKer{\gamma}
    }
  \end{equation*}
  Moreover, the square on the upper right is a pushout. Hint: Corollary \ref{thm:FactoringCoKer}.
\end{exercise}
\end{exercises}
\section[The Category of Short Exact Sequences - I]{The Category of Short Exact Sequences - I}
\label{sec:CatSESs-I}

Associated to a \ZExact\ category $\Ctgry{X}$ is the category $\SESCat{X}$ of short exact sequences in $\Ctgry{X}$. It has: %
\index[not]{s!$\SESCat{X}$\IndSep category of short exact sequences in $\Ctgry{X}$}%
\begin{ulist}
  \item \emph{objects}\quad short exact sequences in $\Ctgry{X}$, which we denote \quad ($\varepsilon$)\quad $M\NMono X \NEpi Q$;
  \item \emph{morphisms}\quad commutative diagrams in $\Ctgry{X}$ whose rows are short exact sequences:
  \begin{equation*}
    \xymatrix@R=5ex@C=4em{
    \varepsilon \ar[d]_{(\kappa,\xi,\rho)} &
    K \ar@{{ |>}->}[r]^-{k} \ar[d]_{\kappa} &
    X \ar@{-{ >>}}[r]^-{q} \ar[d]_{\xi} &
    Q \ar[d]^{\rho} \\
    \varphi &
    L \ar@{{ |>}->}[r]_-{l} &
    Y \ar@{-{ >>}}[r]_-{r} &
    R
    }
  \end{equation*}
\end{ulist}
We show (\ref{thm:SES(C)IsP-Category}) that $\SESCat{X}$ is again a \ZExact\ category. A zero object is given by the short exact sequence $0\to 0\to 0$. If (co)limits over some small category $\SmallCtgry{D}$ exists in $\Ctgry{X}$, then the constructions in (\ref{thm:SESCat(X)(Finite)BiComplete}) provide (co)limits over $\SmallCtgry{D}$ in $\SESCat{X}$. Such computations are facilitated by using that the category $\SESCat{X}$ is equivalent to the category $\NMonoCat{X}$ of normal monomorphisms in $\EuScript{X}$ as well as to the category $\NEpiCat{X}$ of normal epimorphisms in $\EuScript{X}$; see Section \ref{sec:NormalDecompositionsFactorizations}. %
\index[not]{s!$\SESCat{X}$\IndSep category of short exact sequences in $\Ctgry{X}$}

Indeed, forgetting the normal monomorphism of a short exact sequence yields a functor $C\from \SESCat{X}\to \NEpiCat{X}$. Similarly, forgetting the normal epimorphism of a short exact sequence yields a functor $K\from \SESCat{X}\to \NMonoCat{X}$. Using functorial kernels, respectively cokernels, yields equivalence inverses %
\index[not]{n!$\NMonoCat{X}$\IndSep category of normal monomorphisms in $\Ctgry{X}$}%
\begin{equation*}
  \NEpiCat{X}\longrightarrow \SESCat{X}  \qquad \text{and}\qquad \NMonoCat{X} \longrightarrow \SESCat{X}
\end{equation*}
\index[not]{n!$\NEpiCat{X}$\IndSep category of normal epimorphisms in $\Ctgry{X}$}%
Here is a diagrammatic view of these equivalences.
\begin{equation*}
  \xymatrix@!0@R=6em@C=6em{
  \NMonoCat{X} \ar@<1ex>[rr]^-{\CoKerFunc} \ar@<1ex>[rd]^-{} \ar@{}[rr]|-{\bot} \ar@{}[rd]|-{\rtop} &&
  \NEpiCat{X} \ar@<1ex>[ld]^-{} \ar@<1ex>[ll]^-{\KerFunc} &
  K \overset{k}{\NMono} X \ar@{|->}@<0.5ex>[rr]^-{\CoKerMap{k}} \ar@{|->}@<0.5ex>[rd] &&
  X \overset{q}{\NEpi} Q \ar@{|->}@<0.5ex>[ll]^-{\KerMap{q}}\ar@{|->}@<0.5ex>[ld] \\
  & \SESCat{X} \ar@<1ex>[lu]^{K} \ar@<1ex>[ru]^{C} \ar@{}[ru]|-{\ltop}&&&
  K \overset{k}{\NMono} X \overset{q}{\NEpi} Q \ar@{|->}@<0.5ex>[ru] \ar@{|->}@<0.5ex>[lu]
  }
\end{equation*}

\begin{theorem}[$\SESCat{X}$ is a \ZExact\ category\ZExactTag]
  \label{thm:SES(C)IsP-Category}
  If $\Ctgry{X}$ has (co)limits over a small category $\SmallCtgry{D}$, then so does $\SESCat{X}$.
\end{theorem}
\begin{proof}
  If $\ZeroObject$ is a zero object in $\Ctgry{X}$, then the short exact sequence $0\to 0\to 0$ is a zero object in $\SESCat{X}$. To see what (co)limits exist in $\SESCat{X}$ we may use the equivalence between $\SESCat{X}$ and either of the categories $\NMonoCat{X}$ and $\NEpiCat{X}$, combined with Proposition \ref{thm:(Co)Limits-NMono(X)/NEpi(X)}.
\end{proof}

Let us explain explicitly how (co)limits in $\SESCat{X}$ are computed: A diagram of short exact sequences modeled on a small category $\SmallCtgry{D}$ is given by a functor $F\from \SmallCtgry{D}\to \SESCat{X}$ . Then every object $d$ in $\SmallCtgry{D}$ yields a short exact sequence
\begin{equation*}
  \xymatrix@R=5ex@C=4em{
  Fd & K_{Fd} \ar@{{ |>}->}[r]^-{\kappa_{Fd}} &
  X_{Fd} \ar@{-{ >>}}[r]^-{q_{Fd}} &
  Q_{Fd}
  }
\end{equation*}
With this notational convention we have:

\begin{corollary}[Computation of (co)limits in $\SESCat{X}$\ZExactTag]%
  \label{thm:SESCat(X)(Finite)BiComplete}
  \label{thm:SESCat(X)-(Co)Limits}%
  For a pointed category $\SACtgry{X}$ in which (co)limits modeled on a small category $\SmallCtgry{D}$ exist:
  \begin{thmlist}
    \item $D$-colimits in $\SESCat{X}$ exist. If $F\from D\to \SESCat{X}$ is a diagram modeled on a (finite) small category $D$, then
    \begin{equation*}
      \xymatrix@R=5ex@C=4em{
      \CoLimOf{F} & \Ker{ \CoLimOf{q\Comp F}} \ar@{{ |>}->}[r] &
      \CoLimOf{X\Comp F} \ar@{-{ >>}}[r]^-{ \CoLimOf{q\Comp F}} &
      \CoLimOf{Q\Comp F}
      }
    \end{equation*}
    \item $C$-limits in $\SESCat{X}$ exist. If $F\from D\to \SESCat{X}$ is a diagram modeled on a (finite) small category $D$, then
    \begin{equation*}
      \xymatrix@R=5ex@C=4em{
      \LimOf{F} & \LimOf{K\Comp F} \ar@{{ |>}->}[r]^-{\LimOf{\kappa\Comp F}} &
      \LimOf{X\Comp F} \ar@{-{ >>}}[r] &
      \CoKer{\LimOf{\kappa\Comp F}}
      }
    \end{equation*}
  \end{thmlist}
\end{corollary}
\begin{proof}
  These claims may be verified directly. However, by (\ref{thm:(Co)Limits-NMono(X)/NEpi(X)}), limits in $\NMonoCat{X}$ are computed pointwise, and colimits in $\NEpiCat{X}$ are computed pointwise.  Via the equivalences between $\NMonoCat{X}$, $\NEpiCat{X}$, and $\SESCat{X}$, the computation of (co)limits in $\SESCat{X}$ is encapsulated in these diagrams:
  \begin{equation*}
    \xymatrix@R=5ex@C=5em{
    \SESCat{X} \ar[r]^-{\Lim\ \text{in}\ \SESCat{X}} \ar[d]_{\simeq}^{K} &
    \SESCat{X} &
    \SESCat{X} \ar[r]^-{\CoLim\ \text{in}\ \SESCat{X}} \ar[d]_{\simeq}^{C} &
    \SESCat{X} \\
    \NMonoCat{X} \ar[r]_-{\Lim\ \text{in}\ \NMonoCat{X}} &
    \NMonoCat{X} \ar[u]_{\simeq}^{\KerFunc} &
    \NEpiCat{X} \ar[r]_-{\CoLim\ \text{in}\ \NEpiCat{X}} &
    \NEpiCat{X} \ar[u]_{\simeq}^{\CoKerFunc}
    }
  \end{equation*}
  This is what we wanted to show.
\end{proof}

We conclude that, while $\SESCat{X}$ is a category of diagrams, it is not equivalent to a functor category because (co-)limits in $\SESCat{X}$ are \emph{not necessarily} computed object-wise. For example:

\begin{corollary}[(Co)kernels in $\SESCat{X}$\ZExactTag]
  \label{thm:(Co)Kernels-SES(C)}%
  The (co)kernel of a morphism $(\mu,\xi,\eta)$ in $\SESCat{X}$ is computed as in this commutative diagram with horizontal short exact sequences, and where $\lambda$ is the kernel, and $\sigma$ the cokernel of $\xi$:
  \begin{equation}\label{eq:KC}
    \begin{vcenter}{
      \xymatrix@R=5ex@C=4em{
      & M \ar@{{ |>}->}[r]^-{m} \ar@{{ |>}->}[d]_{\mu} \PullLU{rd} &
      N \ar@{-{ >>}}@[blue][r] \ar@{{ |>}->}[d]_{\nu} &
      {\color{blue} \CoKer{m}} \ar[d]^{\zeta} \\
      \varepsilon \ar[d]_{(\kappa,\xi,\rho)} &
      K \ar@{{ |>}->}[r]^-{k} \ar[d]_{\kappa} &
      X \ar@{-{ >>}}[r]^-{q} \ar[d]_{\xi} &
      Q \ar[d]^{\rho} \\
      \varphi &
      L \ar@{{ |>}->}[r]_-{l} \ar[d]_{\eta} &
      Y \ar@{-{ >>}}[r]_-{r} \ar@{-{ >>}}[d]_{\sigma} \PushRD{rd} &
      R \ar@{-{ >>}}[d]^{\tau} \\
      & {\color{red} \Ker{t}} \ar@[red]@{{ |>}->}[r] &
      S \ar@{-{ >>}}[r]_-{t} & T
      }
      }
    \end{vcenter}
  \end{equation}
  In this situation, the top left square is a pullback, and the bottom right square is a pushout.
\end{corollary}
\begin{proof}
  The computation of $\KerMap{\kappa,\xi,\rho}=(\mu,\nu,\zeta)$ and of $\CoKerMap{\kappa,\xi,\rho}=(\eta,\sigma,\tau)$ is given by (\ref{thm:SESCat(X)(Finite)BiComplete}). The pullback/pushout claims follow from (\ref{thm:KernelsInNMono(X)-CoKernelsInNEpi(X)}).
\end{proof}

\begin{corollary}[Monos/epis in $\NMonoCat{X}$, $\SESCat{X}$, $\NEpiCat{X}$\ZExactTag]
  \label{thm:NMono(X)/SES(X)-Monos/Epis}
  The following hold for a morphism in $\NMonoCat{X}$, respectively $\SESCat{X}$, respectively $\NEpiCat{X}$:
  \begin{equation*}
    \xymatrix@R=5ex@C=4em{
    K \ar@{{ |>}->}[r]^-{k} \ar[d]_{\kappa} &
    X \ar[d]^{\xi} &
    K \ar@{{ |>}->}[r]^-{k} \ar[d]_{\kappa} &
    X \ar@{-{ >>}}[r]^-{q} \ar[d]^{\xi} &
    Q \ar[d]^{\rho} &
    X \ar@{-{ >>}}[r]^-{q} \ar[d]_{\xi} &
    Q \ar[d]^{\rho} \\
    L \ar@{{ |>}->}[r]_-{l} &
    Y &
    L \ar@{{ |>}->}[r]_-{l} &
    Y \ar@{-{ >>}}[r]_-{r} &
    R &
    Y \ar@{-{ >>}}[r]_-{r} &
    R
    }
  \end{equation*}
  \begin{thmlist}
    \item Any of the three morphisms above is a monomorphism in its category if and only if $\xi$ is a monomorphism.
    \item Any of the three morphisms above is a epimorphism in its category if and only if $\xi$ is a epimorphism.
  \end{thmlist}
\end{corollary}
\begin{proof}
  (i) We prove the claim in $\SESCat{X}$, and infer the claim in $\NMonoCat{X}$, respectively in $\NEpiCat{X}$ from the equivalences between these three categories. - If $(\kappa,\xi,\rho)$ is a monomorphism in $\SESCat{X}$ we infer that $\xi$ is monic from this situation:
  \begin{equation*}
    \xymatrix@R=5ex@C=4em{
    0 \ar@{{ |>}->}[r] \ar[d] &
    A \ar@{=}[r] \ar@<-0.5ex>[d]_{a} \ar@<0.5ex>[d]^{b} &
    A \ar@<-0.5ex>[d]_{qa} \ar@<0.5ex>[d]^{qb}\\
    K \ar@{{ |>}->}[r]^-{k} \ar[d]_{\kappa} &
    X \ar@{-{ >>}}[r]^-{q} \ar[d]^{\xi} &
    Q \ar[d]^{\rho} \\
    L \ar@{{ |>}->}[r]_-{l} &
    Y \ar@{-{ >>}}[r]_-{r} &
    R
    }
  \end{equation*}
  Suppose $\xi a=\xi b$. Then we obtain two morphisms into the middle row which, after composing with the monomorphism $(\kappa,\xi, \rho)$ are equal. Thus $a=b$, implying that $\xi$ is monic.

  Conversely, suppose $\xi$ is a monomorphism in $\Ctgry{X}$. To show that $(\kappa,\xi,\rho)$ is a monomorphism in $\SESCat{X}$ consider the situation below.
  \begin{equation*}
    \xymatrix@R=5ex@C=5em{
    M \ar@{{ |>}->}[r]^-{m} \ar@<-0.5ex>[d]_{\mu} \ar@<0.5ex>[d]^{\mu'} &
    Z \ar@{-{ >>}}[r]^-{s} \ar@<-0.5ex>[d]_{\zeta} \ar@<0.5ex>[d]^{\zeta'} &
    S \ar@<-0.5ex>[d]_{\gamma} \ar@<0.5ex>[d]^{\gamma'}\\
    K \ar@{{ |>}->}[r]^-{k} \ar[d]_{\kappa} &
    X \ar@{-{ >>}}[r]^-{q} \ar[d]^{\xi} &
    Q \ar[d]^{\rho} \\
    L \ar@{{ |>}->}[r]_-{l} &
    Y \ar@{-{ >>}}[r]_-{r} &
    R
    }
  \end{equation*}
  If $(\kappa,\xi,\rho)\Comp (\mu,\zeta,\gamma) = (\kappa,\xi,\rho)\Comp (\mu',\zeta',\gamma')$, then $\xi\zeta = \xi\zeta'$, implying that $\zeta=\zeta'$ since $\xi$ is monic. It follows that
  \begin{equation*}
    \gamma s = q\zeta = q\zeta' = \gamma's
  \end{equation*}
  So, $\gamma=\gamma'$ via the epic property of the normal epimorphism $b$. Similarly,
  \begin{equation*}
    k\mu = \zeta m = \zeta' m = k\mu'
  \end{equation*}
  So, $\mu=\mu'$ via the monic property of the normal monomorphism $k$. - The proof of (ii) is similar.
\end{proof}

\begin{corollary}[Adjunction between $\ArrowCat{\Ctgry{X}}$ and $\SESCat{X}$\ZExactTag]
  \label{thm:(Co)ReflectArrowToSES}
  Associated to a \BCPtdCat\  $\SACtgry{X}$ are these adjunctions
  \begin{equation*}
    \xymatrix@C=4em{\ArrowCat{\Ctgry{X}} \ar@<1ex>[r] \ar@{}[r]|-{\bot} & \SESCat{X} \ar@<1ex>[l]^-{G\Comp K}}
    \qquad\qquad
    \xymatrix@C=4em{\ArrowCat{\Ctgry{X}} \ar@<1ex>[r] \ar@{}[r]|-{\top} & \SESCat{X} \ar@<1ex>[l]^-{G'\Comp C}}
  \end{equation*}
  Given a morphism $f\from X\to Y$ in $\SACtgry{X}$, the following hold:
  \begin{enumerate}[(i)]
    \item The left adjoint of the adjunction on the left sends $f$ to the short exact sequence
          \begin{equation*}
            \xymatrix@R=5ex@C=3em{
            \Ker{\CoKer{f}} \ar@{{ |>}->}[r] &
            Y \ar@{-{ >>}}[r] &
            \CoKer{f}
            }
          \end{equation*}
    \item The right adjoint of the adjunction on the right sends $f$ to the short exact sequence
          \begin{equation*}
            \xymatrix@R=5ex@C=3em{
            \Ker{f} \ar@{{ |>}->}[r] &
            X \ar@{-{ >>}}[r] &
            \CoKer{\Ker{f}}
            }
          \end{equation*}
  \end{enumerate}
\end{corollary}

Since the category of short exact sequences in a \ZExact\ category is again a \ZExact\ category, we may consider short exact sequences of short exact sequences; that is we may consider the category  $\SESCat{\SESCat{X}}$. Due to the particular way in which we compute (co)limits in $\SESCat{X}$ unexpected phenomena may occur. For example:

\begin{example}[Short exact sequence of short exact sequences in $\TopGrps$]
  \label{exa:SESofSES's-TopGrps}
  Based on \cite[p.~43]{MMClementino2021-TopAlgs}:\quad Consider this commutative diagram of topological groups:
  \begin{equation*}
    \xymatrix@R=5ex@C=4em{
    0 \ar[r] \ar[d] &
    R^{\delta} \ar@{=}[r] \ar@{{ |>}->}[d]^{\PrdctMapInto{0,\IdMap}} &
    R^{\delta} \ar[d]^{h} \\
    R^{\delta} \ar@{{ |>}->}[r]^-{\PrdctMapInto{\IdMap,0}} \ar[d]_{h} &
    X \ar@{-{ >>}}[r]_-{\PrjctnOnto{2}} \ar@{-{ >>}}[d]_-{\PrjctnOnto{1}} &
    R \ar@{-{ >>}}[d]\\
    R \ar@{=}[r] &
    R \ar@{-{ >>}}[r] &
    0
    }
  \end{equation*}
  Here $R$ denotes the additive group of real numbers with the Euclidean topology, and $R^{\delta}$ the same group with the discrete topology. Let $h\from R^{\delta}\to R$ denote the continuous bijection given by the identity map on underlying sets. Further, let $X$ be the topological space obtained by refining the product topology on $\Prdct{R}{R}$ by the sets
  \begin{equation*}
    U_s \DefEq \SetSlct{(x,y)\in \Prdct{R}{R}}{x-y=s}
  \end{equation*}
  Then addition yields the structure of a topological group on $X$. Moreover, the functions $\PrjctnOnto{1},\PrjctnOnto{2}\from X\to R$ which are projections on underlying sets are normal epimorphisms whose kernels are copies of $R^{\delta}$.

  So, the rows of the above diagram are short exact sequences in the category $\TopGrps$ of topological groups. Vertically, we have a short exact sequence of short exact sequences. Remarkable though is that the vertical sequences on the left and on the right fail to be short exact.
\end{example}

\begin{exercises}

\begin{exercise}[Non-pointwise mystery in short exact sequence of short exact sequences\ANKTag]
  \label{exe:SES(SES)-Mystery}%
  For a short exact sequence of short exact sequences in a pointed category %
  \index[acr]{a!$\EuRoman{A{\kern-0.2ex}N{\kern-0.15ex}K}$\IndSep exercise asking a question whose answer we do not know}%
  \begin{equation*}
    \vcenter{\xymatrix@R=5ex@C=3em{
    M  \ar@{{ |>}->}[r]^-{a} \ar@{{ |>}->}[d]_-{u} &
    L \ar@{-{ >>}}[r]^-{d} \ar@{{ |>}->}[d]^-{v} &
    I \ar[d]^-{w} \\
    K \ar@{{ |>}->}[r]_-{b} \ar[d]_-{x} &
    X \ar@{-{ >>}}[r]^-{e} \ar@{-{ >>}}[d]_-{y} &
    Q \ar@{-{ >>}}[d]^-{z} \\
    J \ar@{{ |>}->}[r]_-{c} &
    R \ar@{-{ >>}}[r]_-{f} &
    S }}
  \end{equation*}
  decide if $w$ is always a monomorphism, and if $x$ is always an epimorphism.
\end{exercise}
\end{exercises}
\section[The Normal Short \texorpdfstring{$5$}{5}-Lemma]{The Normal Short \texorpdfstring{$5$}{5}-Lemma}
\label{sec:SESMaps}%

Here, we collect basic properties of morphisms of short exact sequences. This development results in a version of the (Short) $5$-Lemma which holds in \ZExact\ categories.

\begin{proposition}[Morphism of short exact sequences: properties\ZExactTag]
  \label{thm:MorphismSESs-Props-ANN}
  In a morphism of short exact sequences assume that $\xi$ is a normal map.
  \begin{equation*}
    \xymatrix@R=5ex@C=4em{
    K \ar@{{ |>}->}[r]^{k} \ar[d]_{\kappa} &
    X \ar@{-{ >>}}[r]^{q} \ar[d]_{\xi} &
    Q \ar[d]^{\rho} \\
    L \ar@{{ |>}->}[r]_{l} &
    Y \ar@{-{ >>}}[r]_{r} &
    R
    }
  \end{equation*}
  Then the following hold:
  \begin{thmlist}
    \item If $\kappa$ and $\rho$ are epimorphisms, then $\xi$ is a normal epimorphism.
    \item If $\kappa$ and $\rho$ are monomorphisms, then $\xi$ is a normal monomorphism.
  \end{thmlist}
\end{proposition}
\begin{proof}
  (i)\quad Augment the given morphism of short exact sequences with the cokernels of $\kappa$, $\xi$, and $\rho$ to obtain this commutative diagram:
  \begin{equation*}
    \xymatrix@R=5ex@C=4em{
    K \ar@{{ |>}->}[r]^{k} \ar@{-{>>}}[d]_{\kappa} &
    X \ar@{-{ >>}}[r]^{q} \ar[d]_{\xi} &
    Q \ar@{-{>>}}[d]^{\rho} \\
    L \ar@{{ |>}->}[r]_{l} \ar@{-{ >>}}[d]  &
    Y \ar@{-{ >>}}[r]_{r} \ar@{-{ >>}}[d] &
    R \ar@{-{ >>}}[d] \\
    \CoKer{\kappa}=\ZeroObject \ar@{{ |>}->}[r] &
    \CoKer{\xi} \ar@{-{ >>}}[r]_-{s} &
    \CoKer{\rho}=\ZeroObject
    }
  \end{equation*}
  The bottom row happens to be short exact because (a) the map $\ZeroObject\to\CoKer{\xi}$ is a normal monomorphism, and (b) cokernels commute to the effect that
  \begin{equation*}
    \CoKer{\rho}=\CoKer{\CoKer{q}\to \CoKer{r}} \cong \CoKer{\CoKer{\kappa}\to \CoKer{\xi}}=\CoKer{\ZeroMap}
  \end{equation*}
  Thus $s$ is an isomorphism, implying that $\CoKer{\xi}=\ZeroObject$.

  The normal map $\xi$ is a composite $\xi=me$, where $e$ is a normal epimorphism, and  $m$ is a normal monomorphism. With (\ref{thm:NormalEpi-Props}), we see that $\CoKer{\xi}=\CoKer{m}=\ZeroObject$. Thus $m$ is an isomorphism.

  (ii)\quad Coaugmenting the given morphism of short exact sequences with the kernels of $\kappa$, $\xi$, and $\rho$, we start an argument which is dual to the one given in (i), and the claim follows.
\end{proof}

\begin{corollary}[Normal Short $5$-Lemma\ZExactTag]
  \label{thm:Short-5-Primordial}
  \label{thm:Short-5-Normal}%
  Consider a morphism of short exact sequences in which $\xi$ is a normal map. %
  \index{Short 5-Lemma!}%
  \begin{equation*}
    \xymatrix@R=5ex@C=3em{
    K \ar@{{ |>}->}[r]^{k} \ar@{{ >}-{>>}}[d]_{\kappa} &
    X \ar@{-{ >>}}[r]^{q} \ar[d]_{\xi} &
    Q \ar@{{ >}-{>>}}[d]^{\rho} \\
    L \ar@{{ |>}->}[r]_{l} &
    Y \ar@{-{ >>}}[r]_{r} &
    R
    }
  \end{equation*}
  If $\kappa$ and $\rho$ are both monic and epic, then $\kappa$, $\xi$, $\rho$ are isomorphisms.
\end{corollary}
\begin{proof}
  With (\ref{thm:MorphismSESs-Props-ANN}) we see that $\xi$ is an isomorphism. Consequently, the epimorphism $\kappa$ is also a normal monomorphism (\ref{thm:NormalMono-Props}), and the monomorphism $\rho$ is also a normal epimorphism (\ref{thm:NormalEpi-Props}). So, all of $\kappa$, $\xi$, $\rho$ are isomorphisms.
\end{proof}

\begin{corollary}[Normal $5$-Lemma\ZExactTag]
  \label{thm:5-Lemma-Primordial}
  \label{thm:5-Lemma-Normal}
  Consider a morphism of exact sequences. %
  \index{$5$-Lemma}
  \begin{equation*}
    \xymatrix@R=5ex@C=3em{
    \DiagObj \ar@{->>}[d]_{a} \ar[r] &
    \DiagObj \ar[d]_{b}^{\cong} \ar[r] &
    \DiagObj \ar[d]_{c} \ar[r] &
    \DiagObj \ar[d]^{\cong}_{d} \ar[r] &
    \DiagObj \ar@{{ >}->}[d]_{e} \\
    \DiagObj \ar[r] &
    \DiagObj \ar[r] &
    \DiagObj \ar[r] &
    \DiagObj \ar[r] &
    \DiagObj
    }
  \end{equation*}
  If $c$ is a normal map, $a$ is an epimorphism, $e$ a monomorphism, and $b$, $d$ isomorphisms, then $c$ is an isomorphism.
\end{corollary}
\begin{proof}
  Via pushout and pullback recognition, we reduce the proof of the Normal $5$-Lemma to an application of the Normal Short $5$-Lemma (\ref{thm:Short-5-Normal}). As the rows are exact we obtain normal factorizations at selected points as shown below.
  
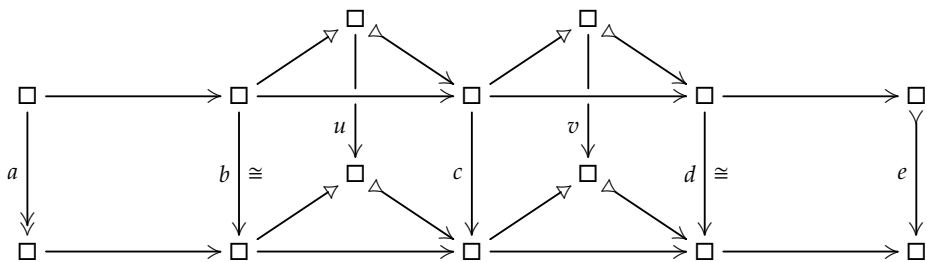
\begin{equation*}
    \xymatrix@R=3ex@C=2.5em{
    &&& \DiagObj \ar@{{ |>}->}[dr] \ar[dd]|\hole_(.7){u} &&
    \DiagObj \ar@{{ |>}->}[dr] \ar[dd]|\hole_(.7){v} \\
    \DiagObj \ar@{ ->>}[dd]_{a} \ar[rr] &&
    \DiagObj \ar[dd]_{b}^{\cong} \ar[rr] \ar@{-{ >>}}[ru] &&
    \DiagObj \ar[dd]_{c} \ar[rr] \ar@{-{ >>}}[ru] &&
    \DiagObj \ar[dd]^{\cong}_{d} \ar[rr] &&
    \DiagObj \ar@{{ >}->}[dd]_{e} \\
    &&& \DiagObj \ar@{{ |>}->}[dr] &&
    \DiagObj \ar@{{ |>}->}[dr] \\
    \DiagObj \ar[rr] &&
    \DiagObj \ar[rr] \ar@{-{ >>}}[ru]&&
    \DiagObj \ar[rr] \ar@{-{ >>}}[ru] &&
    \DiagObj \ar[rr] &&
    \DiagObj
    }
  \end{equation*}
  The square with vertical maps $b$ and $u$ is a pushout because $a$ is an epimorphism \eqref{thm:PushoutRecognize-Categorical}, implying that $u$ is an isomorphism. The square with vertical maps $v$ and $d$ is a pullback because $e$ is a monomorphism \eqref{thm:PullbackRecognition-KernelSide-1}, implying that $v$ is an isomorphism. Now the Normal Short $5$-Lemma (\ref{thm:Short-5-Normal}) tells us that $c$ is an isomorphism.
\end{proof}

\begin{subordinate}[Comments]{}
  In the Normal Short $5$-Lemma and the Normal $5$-Lemma, the requirement that the vertical middle maps~$\xi$, respectively $c$, be normal is very strong. We will see that both results hold in homological categories, without assuming that $\xi$, respectively $c$, are normal maps. On the other hand, the result encompasses the full Short $5$-Lemma in the context of any p-exact category; in particular, we regain the result for abelian categories.

  The earliest version of the Normal Short $5$-Lemma we know of appears as Theorem~6.2 of~\cite{ZJanelidze-Snake}, where it is shown to characterize subtractive categories amongst pointed regular categories. Note, however, that the structural axioms at the foundation of~\cite{ZJanelidze-Snake} are incompatibly different from ours. To begin, the notion of exactness is different, and as a result the nature of the Short $5$-Lemma is different. In fact, the concept of short exact sequence used there is weaker than ours, and the image of the middle map $\xi$ in the diagram is only assumed to be an ideal (i.e., a regular image of a normal monomorphism along a regular epimorphism), rather than a normal monomorphism. Hence the two versions of a Normal Short $5$-Lemma are not immediately comparable.
\end{subordinate}

\bigskip

\begin{exercises}

\begin{exercise}[Short $5$-Lemma in $\SetsBsd$]
  With the notation of (\ref{thm:Short-5-Primordial}), show that the primordial Short $5$-Lemma holds in the category $\SetsBsd$ of pointed sets, without assuming  that $\xi$ be a normal map.
\end{exercise}
\end{exercises}
\chapter{Homology in Di-Exact Categories}
\label{chap:HomologicalFoundations}
\label{chap:Di-ExactCats}

We build on the preparations in the previous chapter and introduce the concept of a chain complex, along with the concept of homology as a measure of the failure of exactness of a chain complex. Actually, contrary to what happens in an abelian category, the analysis presented in Section \ref{sec:Homology} shows that there are two canonical measures of the failure of a chain complex to be exact in a given position. While the constructions of these two measures are dual to one another, they are \emph{not necessarily isomorphic}.

Thus, we face two challenges: (a) determine under which conditions the two dual constructions of homology are naturally isomorphic, and (b) provide standard tools to be able to compute effectively with homological invariants; tools such as:
\begin{ulist}
  \item The Snake Lemma;
  \item the $(\Prdct{3}{3})$-Lemma;
  \item the (Short) $5$-Lemma.
\end{ulist}
As explained in the introduction to this part of the work (Chapter~\ref{chap:IntroPart-I}), whenever these tools are available in a given \ZExact\ category $\Ctgry{X}$, they belong to the self-dual axis of $\Ctgry{X}$. Accordingly, we look for least demanding self-dual structural axioms which ensure the validity of (a) and (b).

Perhaps surprisingly, the structural axioms we distil in Section \ref{sec:DiExtensions} are all related to situations which permit the construction of a \Defn{di-extension}, that is a commutative diagram in which every row and column is a short exact sequence:
\begin{equation*}
  \xymatrix@R=5ex@C=4em{
  M \ar@{{ |>}->}[r]^-{a}	\ar@{{ |>}->}[d]_{u} &
  L \ar@{{ |>}.>}[d]_{\kappa} \ar@{-{ >>}}[r]^-{d} &
  \DiagObj  \ar@{{ |>}.>}[d]^{w} \\
  K \ar@{{ |>}->}[r]^-{m} \ar@{.{ >>}}[d]_{\varepsilon} &
  X \ar@{-{ >>}}[r]_-{\pi} \ar@{-{ >>}}[d]^{e} &
  Q \ar@{-{ >>}}[d]^{z} \\
  \DiagObj	\ar@{{ |>}.>}[r]_-{\mu} &
  R \ar@{-{ >>}}[r]_{f} &
  S
  }
\end{equation*}
For example, the two dual constructions of homology are canonically isomorphic if and only if, whenever the composite $em=\ZeroMap$, then the composite $\pi\kappa$ is a normal map, which happens precisely when the diagram is a di-extension. ---Indeed, in this situation the top left square is a pullback with $M=K$, and the bottom right square is a pushout with $R=S$. So, the top right corner simultaneously represents both the familiar homology of the chain complex $K \XRA{m} X \XRA{e} R$, as well as its dually constructed homology.

We say that a \ZExact\ category $\Ctgry{X}$ is \Defn{homologically self-dual}, \HSDInline\ for short, if it satisfies the condition just described. We refer to the condition itself as the \HSDInline-condition\footnote{The significance of the \HSDInline-condition has already been recognized by Grandis in \cite{Grandis-HA2} in settings which require normal monomorphisms and normal epimorphisms to be closed under composition; there it is called the \emph{homology axiom}.}. It turns out that the scope of the \HSDInline-condition extends beyond immediate homological concerns. For example, we show in Section \ref{sec:HomologicalSelfDuality} that it is equivalent to the validity of the Third Isomorphism Theorem, as well as the validity of the Pure Snake Lemma.

In Sections \ref{sec:DinversionPreservesNormalMaps}, \ref{sec:Sub-di-exact Categories}, and \ref{sec:DiExactCats} we introduce increasingly more demanding self-dual structural axioms we ask a \ZExact\ category to satisfy. For an outline, we refer to the diagram above.

The weakest of these axioms, discussed in Section \ref{sec:DinversionPreservesNormalMaps},  requires $em$ to be a normal map if and only if $\pi\kappa$ is a normal map. This is equivalent to the validity of the border cases of the $(\Prdct{3}{3})$-Lemma.

The second of these axioms, discussed in Section \ref{sec:Sub-di-exact Categories}, is a convenient self-dual requirement which ensures that the Snake Lemma is true; see Section \ref{sec:Snakes}.

The strongest of these axioms, discussed in Section \ref{sec:DiExactCats} plays a special role. It requires every antinormal map to be a normal map. This means that every composite of the kind $em$ above admits a factorization $\mu\varepsilon$. We call this the \ANNInline-condition (for `antinormal composites are normal'), and a \ZExact\ category satisfying the \ANNInline-condition is called \Defn{di-exact}.

On the one hand, the \ANNInline-condition characterizes the separation between homological categories, Chapter \ref{chap:HomologicalCats}, and semiabelian categories, Chapter \ref{chap:SACats}. On the other, in a di-exact category, the Snake Lemma holds, the border cases of the $(\Prdct{3}{3})$-Lemma hold, and a primordial version of the (Short) $5$-Lemma holds. In addition, the remarkably subtle Theorem \ref{thm:(Co)Ker(ProperMapLESs)} holds.

We close this chapter with a discussion of normal pushouts and normal pullbacks in Section \ref{sec:NormalPushouts/Pullbacks}, and of higher extensions in Section \ref{sec:HigherExtensions}. The latter results in boundary cases of higher-dimensional versions of the $(\Prdct{3}{3})$-Lemma.

\newpage

\begin{center}
  \textbf{Leitfaden for Chapter \ref{chap:HomologicalFoundations}}
\end{center}

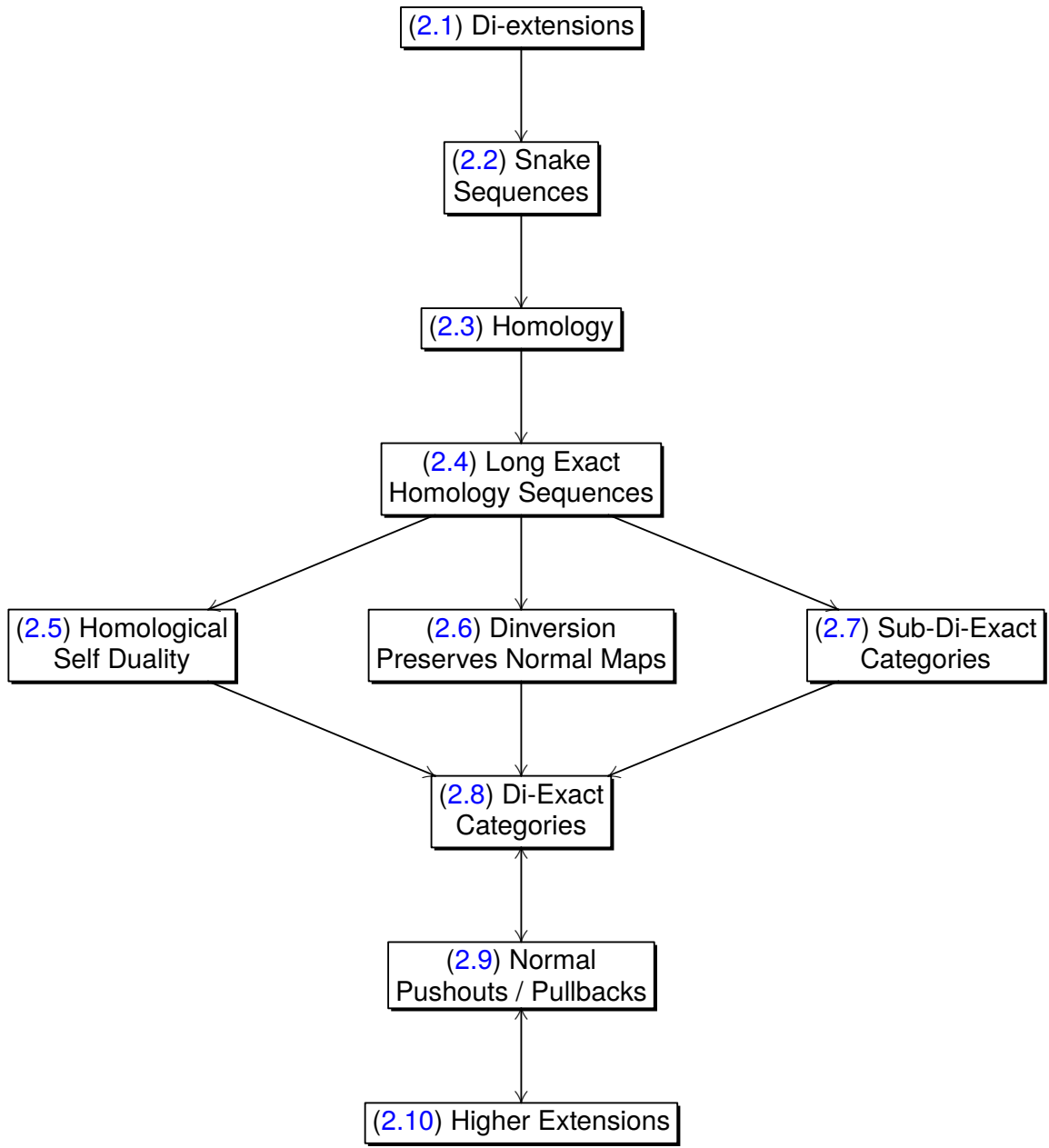
\bigskip

\begin{equation*}
  \xymatrix@R=7ex@C=4.5em{
  & *+[F-,]{\txt{\sffamily (\ref{sec:DiExtensions}) Di-extensions}} \ar[d] \\
  & *+[F-,]{\txt{\sffamily (\ref{sec:Snakes}) Snake \\ \sffamily Sequences}} \ar[d] \\
  & *+[F-,]{\txt{\sffamily (\ref{sec:Homology}) Homology}} \ar[d] \\
  & *+[F-,]{\txt{\sffamily (\ref{sec:LES-Homology}) Long Exact \\ \sffamily Homology Sequences}} \ar[d] \ar[dl] \ar[dr] \\
  *+[F-,]{\txt{\sffamily (\ref{sec:HomologicalSelfDuality}) Homological \\ \sffamily Self Duality}} \ar[dr] &
  *+[F-,]{\txt{\sffamily (\ref{sec:DinversionPreservesNormalMaps}) Dinversion \\ \sffamily Preserves Normal Maps}} \ar[d] &
  *+[F-,]{\txt{\sffamily (\ref{sec:Sub-di-exact Categories}) Sub-Di-Exact \\ \sffamily Categories}} \ar[dl] \\
  & *+[F-,]{\txt{\sffamily (\ref{sec:DiExactCats}) Di-Exact \\ \sffamily Categories}} \ar@{<->}[d] \\
  & *+[F-,]{\txt{\sffamily (\ref{sec:NormalPushouts/Pullbacks}) Normal  \\ \sffamily Pushouts / Pullbacks}} \ar@{<->}[d] \\
  & *+[F-,]{\txt{\sffamily (\ref{sec:HigherExtensions}) Higher Extensions}} \\
  }
\end{equation*}
\section[Di-Extensions]{Di-Extensions}
\label{sec:DiExtensions}

In Section \ref{sec:CatSESs-I}, we observed in Example \ref{exa:SESofSES's-TopGrps} that a short exact sequence in the category $\SESCat{X}$ of short exact sequences is a $(\Prdct{3}{3})$-diagram which might have a column/row which fails to be a short exact sequence in the underlying category. Here we single out a special class of short exact sequences of short exact sequences which are free of this phenomenon:

\begin{definition}[$(\Prdct{3}{3})$-diagram / di-extension\ZExactTag]
  \label{def:3x3Diagram-DoubleExtension}%
  In any category, a \Defn{$(\Prdct{3}{3})$-diagram} is a commutative diagram of the form %
  \index{$(\Prdct{3}{3})$-diagram}\index{di-extension}\index{extension!di-}
  \begin{equation}\label{eq:3x3}
    \vcenter{
    \xymatrix@R=5ex@C=3em{
    M  \ar[r]^-{a} \ar[d]_-{u} &
    L \ar[r]^-{d} \ar[d]^-{v} &
    I \ar[d]^-{w} \\
    K \ar[r]_-{b} \ar[d]_-{x} &
    X \ar[r]^-{e} \ar[d]_-{y} &
    Q \ar[d]^-{z} \\
    J \ar[r]_-{c} &
    R \ar[r]_-{f} &
    S }}
  \end{equation}
  In a \ZExact\ category $\Ctgry{X}$, such a diagram is called a \Defn{di-extension of $M$ and over $S$} if each of its rows and columns is a short exact sequence in the underlying category $\Ctgry{X}$.
\end{definition}

We develop methods for the construction of di-extensions. Perhaps surprisingly, the conditions we encounter match exactly those conditions  which a \ZExact\ category must satisfy so that diagram lemmas, such as the Snake Lemma, familiar from homological algebra hold. Here is an outline. If \eqref{eq:3x3} is a di-extension, then the square $K\rightrightarrows R$ provides the normal decomposition $cx$ of the map $K\to R$. This map is also antinormally decomposed as $yb$, i.e., a normal monomorphism followed by a normal epimorphism. - The same kind of structure is associated to the square $L\rightrightarrows Q$.

We will see that every di-extension centered upon $X$ can be constructed from a suitable antinormal pair $(y,b)$, or equivalently, the antinormal pair $(e,v)$. In fact, in a di-extension the antinormal pairs $(y,b)$ and $(e,v)$ determine each other, up to isomorphism, through dinversion ($=$ di-inversion), that is the normal subobject / quotient inversion (\ref{term:NormalSubobject/QuotientObjectInversion}) applied to $y$ and $b$, respectively $e$ and $v$.

\begin{definition}[Antinormal (de)composition\ZExactTag]
  \label{def:Antinormal(De)Composition}%
  A composition of morphisms $\varepsilon \mu$ is called \Defn{antinormal} if $\mu$ is a normal monomorphism, and $\varepsilon$ is a normal epimorphism. An \Defn{antinormal decomposition} of a given morphism $f$ is given by an identity $f=\varepsilon\mu$, where $\varepsilon\mu$ is an antinormal composition. A morphism is said to be \Defn{antinormal} if it admits an antinormal decomposition. %
  \index{antinormal!morphism}\index{antinormal!(de)composition}%
\end{definition}

Let us clarify: A normal map $f$ has an essentially unique factorization $f=me$, with $e$ a normal epimorphism and $m$ a normal monomorphism. On the other hand, antinormal decompositions, if they exist, need not be unique. For example, every object $X$ yields an antinormal decomposition $\ZeroObject \NMono X \NEpi \ZeroObject$ of the zero map. Thus, being normal is a property which a given morphism might or might not have. However, an antinormal decomposition of $f=\varepsilon\mu$ is a choice of additional structure associated with $f$; whether or not $f$ admits such a decomposition is indeed a property of $f$, but many such decompositions may exist for this morphism. In response, we refer to a pair $(\varepsilon,\mu)$ consisting of a normal epimorphism $\varepsilon$ and a normal monomorphism $\mu$ as an antinormal pair. By design, it provides an antinormal decomposition of $f\DefEq \varepsilon\mu$. %
\index{antinormal!pair}%

We will frequently encounter an antinormal pair whose components compose to a normal map $f$. So, we spell out explicitly how $f$ may be factored as a normal map.

\begin{lemma}[Factoring an normal map via an antinormal decomposition\ZExactTag]
  \label{thm:NormalMap-FactorViaAntinormalPair}%
  If an antinormal decomposition of $f$ is given by $X \overset{\mu}{\NMono} S \overset{\varepsilon}{\NEpi} Y$ , then the following are equivalent.
  \begin{tfae}
    \item $f=me$ is a normal factorization of $f$.
    \item $f$ is a normal map, and $e=\CoKerMap{\bar{k}}$, where $\bar{k}$ comes from the pullback of $\mu$ along $k\DefEq\KerMap{\varepsilon}$ on the left below.
    
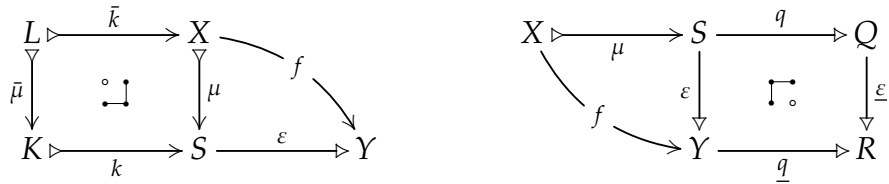
\begin{equation*}
      \xymatrix@R=5ex@C=4em{
      L \ar@{{ |>}->}[r]^-{\bar{k}} \ar@{{ |>}->}[d]_{\bar{\mu}} \PullLU{rd} &
      X \ar@{{ |>}->}[d]^{\mu} \ar@/^2ex/[rd]|-{\ f\ }&&
      X \ar@{{ |>}->}[r]_-{\mu} \ar@/_2ex/[rd]|-{\ f\ }&
      S \ar@{-{ >>}}[r]^-{q} \ar@{-{ >>}}[d]_{\varepsilon} \PushRD{rd} &
      Q \ar@{-{ >>}}[d]^{\underline{\varepsilon}}  \\
      K \ar@{{ |>}->}[r]_-{k} &
      S \ar@{-{ >>}}[r]^-{\varepsilon} &
      Y &&
      Y \ar@{-{ >>}}[r]_-{\underline{q}} &
      R
      }
    \end{equation*}
    \item $f$ is a normal map, and $m=\KerMap{\underline{q}}$, where $\underline{q}$ comes from the pushout of $\varepsilon$ along $q\DefEq \CoKerMap{\mu}$ on the right above.
  \end{tfae}
\end{lemma}
\begin{proof}
  (I) $\Leftrightarrow$ (II)\quad We know that $f$ is a normal map if and only if it admits a factorization $f=me$ with $m$ a normal monomorphism and $e$ a normal epimorphism. In this situation $m$ is uniquely determined by factoring $f$ through $e\DefEq \CoKerMap{\KerMap{f}}$. This latter map is the pullback of $\KerMap{\varepsilon}$ along $\mu$; i.e., it is the map $\bar{k}$ in the pullback diagram under (II), which is in fact just the kernel of $f$. - The equivalence of (I) and (III) follows by dual reasoning.
\end{proof}

Every antinormal pair $(\varepsilon,\mu)$ provides initial data for a double extension via the following construction:
\begin{enumerate}
  \item Construct the inverse antinormal pair $(\pi,\kappa)$ of $(\varepsilon,\mu)$ by taking $\pi\DefEq \CoKerMap{\mu}$ and $\kappa\DefEq \KerMap{\varepsilon}$
  \item Construct $M$ as the pullback as indicated (\ref{thm:Pullback/Pushout-Existence}). Then the maps $a$ and $u$ are normal monomorphisms and, with $\omega\DefEq \pi\kappa$,  $a=\KerMap{\omega}$, and $u=\KerMap{\alpha}$.
  \item Construct $S$ as the pushout as indicated (\ref{thm:Pullback/Pushout-Existence}). Then the maps $f$ and $z$ are normal epimorphisms and, with $\alpha\DefEq \varepsilon\mu$, $f=\CoKerMap{\alpha}$, and $z=\CoKerMap{\omega}$.
\end{enumerate}
\stepcounter{theorem}%
\begin{equation}
  \label{eq:Dinversion}
  \vcenter{
  \xymatrix@R=6ex@C=4.5em{
  M \ar@{{ |>}->}[r]^-{a}	\ar@{{ |>}->}[d]_{u} &
  L \ar@{{ |>}->}[d]_{\kappa} \ar@{-->}[rd]|-{\ \omega\ } \ar@{.>}[r]^-{d} &
  \DiagObj  \ar@{.>}[d]^{w} \\
  K \ar@{{ |>}->}[r]^-{\mu} \ar@{-->}[rd]|-{\ \alpha\ } \ar@{.>}[d]_{x} &
  X \ar@{-{ >>}}[r]_-{\pi} \ar@{-{ >>}}[d]^{\varepsilon} \PullLU{lu} \PushRD{rd} &
  Q \ar@{-{ >>}}[d]^{z} \\
  \DiagObj	\ar@{.>}[r]_-{c} &
  R \ar@{-{ >>}}[r]_{f} &
  S
  }}
\end{equation}%
\index{dinversion}%
By design, the maps $\alpha$ and $\omega$ come with antinormal decompositions $\alpha\DefEq \varepsilon\mu$ and $\omega\DefEq qk$. We are left with the question if $d\DefEq \CoKerMap{a}$ and $w\DefEq \KerMap{z}$ form the normal factorization of $\omega$, and if $x\DefEq\CoKerMap{u}$ and $c\DefEq\KerMap{f}$ form the normal factorization of $\alpha$. - The following lemma confirms that Diagram \ref{eq:Dinversion} admits an expansion to a di-extension if and only if both, $\alpha$ and $\beta$ are normal maps.

\begin{lemma}[Di-extension from antinormal pair\ZExactTag]
  \label{thm:DoubleExtensionFromAntiNormalPair}
  \label{thm:DiExtensionFromAntiNormalPair}%
  Diagram \ref{eq:Dinversion} admits an expansion to a di-extension if and only if $\alpha$ and $\omega$ are normal maps. %
  \index{di-extension!from normal antinormal pairs}%
\end{lemma}
\begin{proof}
  If the maps $\alpha$ and $\omega$ are normal, then (\ref{thm:NormalMap-FactorViaAntinormalPair}) tells us that their normal factorizations are constructed as
  \begin{equation*}
    d\DefEq\CoKerMap{a},\quad  w\DefEq \KerMap{z},\quad x\DefEq \CoKerMap{u},\quad c\DefEq \KerMap{f}
  \end{equation*}
  Thus all rows and columns of the resulting $(\Prdct{3}{3})$-diagram are short exact sequences; i.e., we obtain a di-extension.

  Conversely, if  $(\varepsilon,\mu)$ together with its antinormal inversion $(\pi,\kappa)$ expand to a di-extension, then this di- extension provides normal factorizations of $\alpha$ and $\omega$.
\end{proof}

In general, if a map $f$ admits an antinormal decomposition $f=\varepsilon\mu$,  then $f$ need not be a normal map. So, in general, it makes no sense to ask for a functorial construction from antinormal pairs to di-extensions.  However:

\begin{lemma}[Map of di-extensions from map of antinormal pairs\ZExactTag]
  \label{thm:MorphismDoubleExtensionsFromMorphismAntinormalPairs}
  \label{thm:MorphismDiExtensionsFromMorphismAntinormalPairs}%
  Using the notation \eqref{eq:Dinversion}, suppose the di-extensions (E1) and (E2) contain the antinormal pairs $(\varepsilon,\mu)$ and $(\varepsilon',\mu')$ in their bottom left corners. Then every morphism $(\alpha,\beta,\gamma)\from (\varepsilon,\mu)\to (\varepsilon',\mu')$ of antinormal pairs is the unique restriction of a morphism $(E1)\to(E2)$ of di-extensions.
  
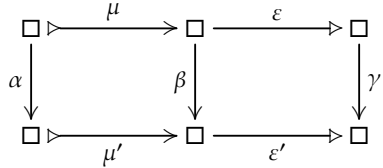
\begin{equation*}
    \xymatrix@R=5ex@C=4em{
    \DiagObj \ar@{{ |>}->}[r]^-{\mu} \ar[d]_{\alpha} &
    \DiagObj \ar@{-{ >>}}[r]^-{\varepsilon} \ar[d]_{\beta}  &
    \DiagObj  \ar[d]^{\gamma} \\
    \DiagObj \ar@{{ |>}->}[r]_-{\mu'} &
    \DiagObj \ar@{-{ >>}}[r]_-{\varepsilon'} &
    \DiagObj
    }
  \end{equation*}
\end{lemma}

Based on (\ref{thm:DiExtensionFromAntiNormalPair}), we single out several settings in which an antinormal pair $(\varepsilon,\mu)$ yields a double extension.

\begin{corollary}[Di-extension from antinormal pair\ZExactTag]
  \label{thm:DoubleExtensionFromAN->N.Pairs}
  \label{thm:DiExtensionFromAN->N.Pairs}%
  Let $\Ctgry{X}$ be a \ZExact\ category in which every antinormal map is a normal map. Then every antinormal pair yields a di-extension. \NoProof
\end{corollary}

\begin{definition}[Di-exact category\ZExactTag]
  \label{def:DiExactCat}
  Of a \ZExact\ category $\Ctgry{X}$, we say that it satisfies the \ANNInline-condition if every antinormal map is normal.  If so, we call $\Ctgry{X}$ \Defn{di-exact}. %
  \index{di-exact category}\index{category!di-exact}%
\end{definition}

In other words, a category $\Ctgry{X}$ is di-exact if and only if the following axioms hold.
\begin{enumerate}
  \item $\Ctgry{X}$ has a zero object; see (\ref{def:0-Object}).
  \item Any morphism in $\Ctgry{X}$ has a (functorially chosen) kernel and a cokernel; see (\ref{sec:Kernel/CoKernel}).
  \item (ANN): Every antinormal map is normal; whenever a morphism factors as a normal monomorphism $\mu$ followed by a normal epimorphism $\epsilon$
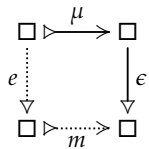
        \[
          \xymatrix{\DiagObj \ar@{{ |>}->}[r]^-{\mu}   \ar@{.{ >>}}[d]_-e  & \DiagObj \ar@{-{ >>}}[d]^-{\epsilon} \\
          \DiagObj  \ar@{{ |>}.>}[r]_m& \DiagObj}
        \]
        it may be written as a composite of a normal epimorphism $e$ followed by a monomorphism $m$.
\end{enumerate}

We use the label `{\color{Cerulean} $\EuRoman{DEx}$}' to identify a di-exact category. %
\index[acr]{d!{\color{Cerulean} $\EuRoman{DEx}$}\IndSep category in which antinormal composites are normal}%

Less demanding is the condition that dinversion preserves normal maps. This means that, dinversion applied to antinormal pair which happens to compose to a normal map yields an antinormal pair which composes to a normal map. With the notation of Diagram \eqref{eq:Dinversion}, this means that $\alpha$ is a normal map if and only if $\omega$ is a normal map. We use the acronym \DPNInline\ for this condition, and we use the label `{\color{Cerulean} $\EuRoman{DPN}$}' to identify \ZExact\ categories which satisfy the \DPNInline-condition. %
\index[acr]{d!{\color{Cerulean} $\EuRoman{DPN}$}\IndSep category in which dinversion preserves normal maps}

\begin{corollary}[Di-extension from normal antinormal pair\DPNTag]
  \label{thm:DiExtensionFromDPN}%
  Let $\Ctgry{X}$ be a \ZExact\ category in which dinversion preserves normal maps. Then every antinormal pair $(\varepsilon,\mu)$ with normal composite $\alpha\DefEq \varepsilon\mu$ provides initial data for a di-extension. \NoProof
\end{corollary}

Of particular interest will be antinormal pairs $(\varepsilon,\mu)$ for which $\varepsilon\mu=\ZeroMap$. We say that a \ZExact\ category is \Defn{homologically self-dual} if dinversion turns every antinormal decomposition of $\ZeroMap$ into an antinormal pair with a normal composite. With the notation of Diagram \eqref{eq:Dinversion}, this means that whenever $\alpha=\ZeroMap$ then  $\omega$ is a normal map. We use the label `{\color{Cerulean} $\EuRoman{HSD}$}' to identify homologically self-dual categories.

\begin{corollary}[Di-extension from antinormal decomposition of $\ZeroMap$\HSDTag]
  \label{thm:DiExtensionFromHSD}%
  In a homologically self-dual category every antinormal pair $(\varepsilon,\mu)$ with $\mu\varepsilon=\ZeroMap$ provides initial data for a di-extension. \NoProof
\end{corollary}

Here is a summary of the structural axioms introduced above:

\begin{terminology}[Di-extensive conditions]
  \label{term:DiExtensiveConditions}%
  Given a \ZExact\ category $\Ctgry{X}$, we introduce the following conditions:
  \begin{enumerate}[(i)]
    \item \label{Ax:0->DPNnverseIsNormal}
          \label{Ax:HSD}%
          \emph{Homological self-duality \HSDInline} holds in $\Ctgry{X}$ if  dinversion turns every antinormal decomposition of the zero map into an antinormal pair with normal composite. - A \ZExact\ category in which the \HSDInline-condition is satisfied is called \Defn{homologically self-dual}. %
          \index[acr]{h!\HSDInline\IndSep homological self-duality property}
    \item \label{Ax:NormalMaps->NANClosure}
          \label{Ax:DinversionPreservesANN}%
          \Defn{Dinversion preserves normal maps \DPNInline} if dinversion turns every antinormal composite which is normal into an antinormal map which is normal. %
          \index[acr]{d!\DPNInline\IndSep property that dinversion preserves normal maps}
    \item \label{Ax:Subnormal+Antinormal->Normal}
          \label{Ax:SubDiExact}%
          \Defn{Sub-di-exactness}\quad $\Ctgry{X}$ is \Defn{sub-di-exact} if every (co)subnormal map in $\Ctgry{X}$ which admits an antinormal decomposition is normal, and dinversion preserves normal antinormal maps.
    \item \label{Ax:Antinormal->Normal}
          \label{Ax:DiExact}%
          \emph{Every antinormal composite is a normal map \ANNInline} - A \ZExact\ category in which the \ANNInline-condition is satisfied is called \Defn{di-exact}. \DExTag%
          \index{di-exact category}\index{category!di-exact}%
          \index[acr]{a!\ANNInline\IndSep property that antinormal maps are normal}%
    \item \label{Ax:ALLNormal} \Defn{All maps are normal}.
  \end{enumerate}
\end{terminology}

Evidently, we have (\ref{Ax:ALLNormal}) $\implies$ (\ref{Ax:DiExact}) $\implies$ (\ref{Ax:SubDiExact}) $\implies$ (\ref{Ax:DinversionPreservesANN}) $\implies$ (\ref{Ax:HSD}). - In Table~\ref{fig:DiExtensionTypes} we graphically summarize the kinds of di-extensions which are associated to the types of antinormal pairs we discussed.

\begin{table}[H]
  \stepcounter{theorem}
  \begin{center}
    \resizebox{.95\textwidth}{!}{%
      \begin{tabular}{m{.3\textwidth}  m{.3\textwidth} | m{.3\textwidth}  m{.3\textwidth}} %
        \toprule                                                                            \\ %
        $\xymatrix@R=5ex@C=4em{
        \DiagObj \ar@{=}[r] \ar@{=}[d]  \PullLU{rd}                                       &
        \DiagObj \ar@{{ |>}->}[d] _{\mu} \ar@{-{ >>}}[r] \BiCart{rd}                      &
        0 \ar@{{ |>}->}[d]                                                                  \\
        \DiagObj \ar@{{ |>}->}[r]^-{\mu} \ar@{-{ >>}}[d] \BiCart{rd}                      &
        \DiagObj \ar@{-{ >>}}[r]_-{\varepsilon} \ar@{-{ >>}}[d]^{\varepsilon} \PushRD{rd} &
        \DiagObj \ar@{=}[d]                                                                 \\
        0 \ar[r]                                                                          &
        \DiagObj \ar@{=}[r]                                                               &
          \DiagObj
          }$
                                                                                          &
        $\begin{array}{c}
             \DiagObj \overset{\mu}{\NMono} \DiagObj \overset{\varepsilon}{\NEpi} \DiagObj \\
             \text{is short exact \ZExactTag}
           \end{array}$
                                                                                          &
        $\xymatrix@R=5ex@C=4em{
        0 \ar@{{ |>}->}[r] \ar@{{ |>}->}[d]  \PullLU{rd}                                  &
        \DiagObj \ar@{=}[r] \ar@{{ |>}->}[d] _{\mu}                                       &
        \DiagObj \ar@{=}[d]                                                                 \\
        \DiagObj \ar@{{ |>}->}[r] \ar@{=}[d]                                              &
        \DiagObj \ar@{-{ >>}}[r]_-{\varepsilon} \ar@{-{ >>}}[d] \PushRD{rd}               &
        \DiagObj \ar@{-{ >>}}[d]                                                            \\
        \DiagObj \ar@{=}[r]                                                               &
        \DiagObj \ar@{-{ >>}}[r]                                                          &
          0
          }$
                                                                                          &
        $\varepsilon\mu=\IdMap$ \DPNTag                                                     \\
        \midrule
        $\xymatrix@R=5ex@C=4em{
        \DiagObj \ar@{{ |>}->}[r] \ar@{=}[d] \PullLU{rd}                                  &
        \DiagObj \ar@{-{ >>}}[r] \ar@{{ |>}->}[d]  \BiCart{rd}                            &
        \DiagObj \ar@{{ |>}->}[d]                                                           \\
        \DiagObj \ar@{{ |>}->}[r]_{\mu} \ar@{-{ >>}}[d]                                   &
        \DiagObj \ar@{-{ >>}}[d]_-{\varepsilon} \ar@{-{ >>}}[r] \PushRD{rd}               &
        \DiagObj \ar@{-{ >>}}[d]                                                            \\
        0 \ar@{{ |>}->}[r]                                                                &
        \DiagObj \ar@{=}[r]                                                               &
          \DiagObj
          }$
                                                                                          &
        $\begin{array}{c}\varepsilon\mu=\ZeroMap \\ \text{\HSDTag, variant 1} \end{array}$
                                                                                          &
        $\xymatrix@R=5ex@C=4em{
        \DiagObj \ar@{{ |>}->}[d] \ar@{=}[r]  \PullLU{rd}                                 &
        \DiagObj \ar@{{ |>}->}[d]^{\mu} \ar@{-{ >>}}[r]                                   &
        0 \ar@{{ |>}->}[d]                                                                  \\
        \DiagObj \ar@{{ |>}->}[r] \ar@{-{ >>}}[d] \BiCart{rd}                             &
        \DiagObj \ar@{-{ >>}}[r]^-{\varepsilon} \ar@{-{ >>}}[d] \PushRD{rd}               &
        \DiagObj \ar@{=}[d]                                                                 \\
        \DiagObj \ar@{{ |>}->}[r]                                                         &
        \DiagObj \ar@{-{ >>}}[r]                                                          &
          \DiagObj
          }$
                                                                                          &
        $\begin{array}{c}\varepsilon\mu=\ZeroMap \\ \text{\HSDTag, variant 2} \end{array}$  \\
        \midrule
        $\xymatrix@R=5ex@C=4em{
        \DiagObj \ar@{{ |>}->}[r] \ar@{{ |>}->}[d] \PullLU{rd}                            &
        \DiagObj \ar@{-{ >>}}[r] \ar@{{ |>}->}[d] _{\mu}                                  &
        \DiagObj \ar@{{ |>}->}[d]                                                           \\
        \DiagObj \ar@{{ |>}->}[r] \ar@{-{ >>}}[d]                                         &
        \DiagObj \ar@{-{ >>}}[r]^-{\varepsilon} \ar@{-{ >>}}[d] \PushRD{rd}               &
        \DiagObj \ar@{-{ >>}}[d]                                                            \\
        \DiagObj \ar@{{ |>}->}[r]                                                         &
        \DiagObj \ar@{-{ >>}}[r]                                                          &
          \DiagObj
          }$
                                                                                          &
        $\begin{array}{c}
             \varepsilon\mu\neq \ZeroMap,\IdMap\quad \text{and} \\
             \varepsilon\mu\ \text{normal \DPNTag}
           \end{array}$
                                                                                          &
        $\xymatrix@R=5ex@C=4em{
        \DiagObj \ar@{{ |>}->}[r] \ar@{{ |>}->}[d] \PullLU{rd}                            &
        \DiagObj \ar@{-{ >>}}[r] \ar@{{ |>}->}[d] _{\mu}                                  &
        \DiagObj \ar@{{ |>}->}[d]                                                           \\
        \DiagObj \ar@{{ |>}->}[r] \ar@{-{ >>}}[d]                                         &
        \DiagObj \ar@{-{ >>}}[r]^-{\varepsilon} \ar@{-{ >>}}[d] \PushRD{rd}               &
        \DiagObj \ar@{-{ >>}}[d]                                                            \\
        \DiagObj \ar@{{ |>}->}[r]                                                         &
        \DiagObj \ar@{-{ >>}}[r]                                                          &
          \DiagObj
          }$
                                                                                          &
        $\varepsilon\mu\neq \ZeroMap,\IdMap$ \DExTag
        \\
        \bottomrule
      \end{tabular}
    }
  \end{center}
  \caption{Di-extension types}
  \label{fig:DiExtensionTypes}
\end{table}

Next, we clarify the relationship between the category $\ANPCat{X}$ of antinormal pairs in a \ZExact\ category $\Ctgry{X}$ and the category $\DExCat{X}$ of di-extensions in $\Ctgry{X}$. To this end, let $\varphi\from \Ord{2}\to \Ord{2}^2$ be the functor with the following description: %
\index[not]{a!$\ANPCat{X}$\IndSep category of antinormal pairs in $\Ctgry{X}$}%
\index[not]{d!$\DExCat{X}$\IndSep category of di-extensions in $\Ctgry{X}$}%

\begin{center}
  \begin{minipage}[m]{3.5cm}
    $\left[0 \XRA{a_{0}} 1 \XRA{a_{1}} 2\right]$
  \end{minipage}\qquad $\XRA{\varphi}$ \qquad
  \begin{minipage}[m]{7cm}
    \begin{equation*}
      \xymatrix@R=5ex@C=4em{
      (0,0) \ar[r] \ar[d] &
      (1,0) \ar[r] \ar[d] &
      (2,0) \ar[d] \\
      (0,1) \ar[r]^-{\varphi(a_{0})\DefEq a_{0}\times \IdMapOn{1}} \ar[d] &
      (1,1) \ar[r] \ar[d]^{\varphi(a_{1})\DefEq \IdMapOn{1}\prdct a_{1}} &
      (2,1) \ar[d] \\
      (0,2) \ar[r] &
      (1,2) \ar[r] &
      (2,2)
      }
    \end{equation*}
  \end{minipage}
\end{center}
Thus, if $X\from \Ord{2}^2\to \Ctgry{X}$ is a di-extension then, using the notation from \eqref{eq:Dinversion}, $X\varphi$ is the antinormal diagram $\alpha=K\XRA{\mu} X \XRA{\varepsilon} R$.

\begin{proposition}[$\varphi^{\ast}$ is equivalence of categories\DExTag]
  \label{thm:ANP(X)<->DEx(X)}%
  The functor $\varphi^{\ast} \from \DExCat{X}\to \ANPCat{X}$ is an equivalence of categories. Using the notation from \eqref{eq:Dinversion}, a pseudo-inverse $\psi\from \ANPCat{X}\to \DExCat{X}$ comes from applying this sequence of functorial constructions to a given antinormal pair $(\varepsilon,\mu)$:
  \begin{enumerate}
    \item Form the antinormal inverse $(\pi\DefEq \CoKerMap{\mu},\kappa\DefEq \Ker{\varepsilon})$.
    \item Form the pullback of $\mu$ and $\kappa$ to obtain the top left square of \eqref{eq:Dinversion}.
    \item Form the pushout of $\varepsilon$ and $\pi$ to obtain the bottom right square of \eqref{eq:Dinversion}.
    \item Form the normal factorization $\alpha=cx$ by putting $x\DefEq \CoKerMap{u}$. Then $\alpha$ factors uniquely and functorially through $x$ via the normal monomorphism $c$.
    \item Form the normal factorization $\omega=wd$ by putting $w\DefEq \Ker{z}$. Then $\omega$ factors uniquely and functorially through $w$ via the normal epimorphism $d$.
  \end{enumerate}
\end{proposition}
\begin{proof}
  The proposed construction yields a di-extension by (\ref{thm:DoubleExtensionFromAntiNormalPair}). The construction is functorial by (\ref{thm:MorphismDoubleExtensionsFromMorphismAntinormalPairs}).
\end{proof}

To compute (co)limits in $\DExCat{X}$, we first clarify how to compute (co)limits in $\ANPCat{X}$.

\begin{proposition}[(Co)limits in $\ANPCat{X}$\ZExactTag]
  \label{thm:ANP(X)-(Co)Limits}%
  The colimit of $J$-diagram $F\from J\to \ANPCat{X}$ may be computed as
  \begin{equation*}
    \CoLimOf{F} = \Biggl(\DiagObj \XRA{\KKerMap{\CoLimOf{\CoKerFunc\Comp F|[0\XRA{a_{1}} 1]}}} \CoLimOf{F|1} \XRA{\CoLimOf{F|a_{1}}} \CoLimOf{F|2}\Biggr)
  \end{equation*}
  The limit of $F$ may be computed as
  \begin{equation*}
    \LimOf{F} = \Biggl(\LimOf{F|0} \XRA{\LimOf{F|a_{0}}} \LimOf{F|1} \XRA{\CCoKerMap{\LimOf{\KerFunc\Comp F|[1\XRA{a_{1}} 2]}} } \DiagObj\Biggr)
  \end{equation*}
  Hence $\ANPCat{X}$ has all (co)limits $\Ctgry{X}$ has. In particular, if $\Ctgry{X}$ is \ZExact, then so is $\ANPCat{X}$.
\end{proposition}
\begin{proof}
  Let us consider the proposed computation of colimits in $\ANPCat{X}$. It is the computation of the colimit in $\NMonoCat{X}$ of the underlying $J$-diagram of normal monomorphisms, composed with the computation of the colimit in $\NEpiCat{X}$ of the underlying diagram of normal epimorphisms. We only need to observe that these two colimits yield composable morphisms in $\EuScript{X}$. Indeed, the codomain of the colimit in $\NMonoCat{X}$ is $\CoLimOf{F|1}$, which is exactly the same as the domain of the colimit in $\NEpiCat{X}$.

  Dual reasoning yields the limit computation of $F$.
\end{proof}

\begin{corollary}[$\DExCat{X}$ is a \ZExact\ category\DExTag]
  \label{thm:DEx(X)Pointed}%
  The category $\DExCat{X}$ associated to a di-exact category $\Ctgry{X}$ is a \ZExact\ category.
\end{corollary}
\begin{proof}
  We only need to show that the category $\DExCat{X}$ has kernels and cokernels. From (\ref{thm:ANP(X)<->DEx(X)}) we know that $\DExCat{X}$ is equivalent to $\ANPCat{X}$. In (\ref{thm:ANP(X)-(Co)Limits}), we reduce the computation of (co)kernels in $\ANPCat{X}$ to the computation of (co)kernels in $\NMonoCat{X}$ and in $\NEpiCat{X}$, and these are computed in  (\ref{thm:(Co)Limits-NMono(X)/NEpi(X)}). So, the claim follows.
\end{proof}

We observe that the requirement that antinormal maps be normal is equivalent to the following condition.

\begin{definition}[Normal image preservation of normal maps\ZExactTag]
  \label{def:NICNormalMaps}
  In a \ZExact\ category \Defn{normal images preserve normal maps} if, for every morphism of normal factorizations in which $f$ and $g$ are normal maps, the map $g$ is normal as well.
  
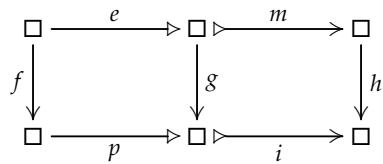
\begin{equation*}
    \xymatrix@R=5ex@C=4em{
    \DiagObj \ar[d]_f \ar@{-{ >>}}[r]^-{e} &
    \DiagObj \ar[d]^{g} \ar@{{ |>}->}[r]^-{m} &
    \DiagObj \ar[d]^h \\
    \DiagObj  \ar@{-{ >>}}[r]_-{p} &
    \DiagObj \ar@{{ |>}->}[r]_-{i} &
    \DiagObj
    }
  \end{equation*}
  \index{normal! map preserved by normal images}
\end{definition}

Keeping in mind that the map $g$ in (\ref{def:NICNormalMaps}) is uniquely determined by the remaining maps, we conclude that, in a category where normal images preserve normal maps, every commutative square of normal maps, as on the left below, has normal factorizations which form the commutative diagram on the right below.
\begin{center}
  \begin{minipage}[m]{4cm}
    \begin{equation*}
      \xymatrix@R=5ex@C=4em{
      \square \ar[r]^-{a} \ar[d]_{\alpha} &
      \square \ar[d]^{\beta} \\
      \square \ar[r]_-{b} &
      \square
      }
    \end{equation*}
  \end{minipage}\qquad\qquad
  \begin{minipage}[m]{6cm}
    \begin{equation*}
      \xymatrix@R=5ex@C=4em{
      \square \ar@/^2ex/[rr]^-{a} \ar@/_2ex/[dd]_{\alpha} \ar@{-{ >>}}[r] \ar@{-{ >>}}[d] &
      \square \ar@{{ |>}->}[r] \ar@{-{ >>}}[d] &
      \square \ar@/^2ex/[dd]^{\beta} \ar@{-{ >>}}[d]  \\
      \square \ar@{-{ >>}}[r] \ar@{{ |>}->}[d] &
      \square \ar@{{ |>}->}[r] \ar@{{ |>}->}[d] &
      \square \ar@{{ |>}->}[d] \\
      \square \ar@/_2ex/[rr]_-{b} \ar@{-{ >>}}[r] &
      \square \ar@{{ |>}->}[r] &
      \square
      }
    \end{equation*}
  \end{minipage}
\end{center}

\begin{proposition}[Antinormal composites and normal images\ZExactTag]
  \label{thm:ANN<->NMNIC}
  \label{thm:ANN<->NIPN}%
  In a \ZExact\ category, the following conditions are equivalent:
  \begin{tfae}
    \item Antinormal composites are normal maps.
    \item Normal images preserve normal maps.
  \end{tfae}
\end{proposition}
\begin{proof}
  (II) $\Rightarrow$ (I)\quad Suppose $g$ admits a factorization $g = \mu \varepsilon$, where $\varepsilon$ is a normal epimorphism, and $\mu$ is a normal monomorphism. Then take $f\DefEq m\DefEq \mu$, $p\DefEq h \DefEq \varepsilon$, and $e$ and $i$  identities. Now (II) implies that $g$ is normal.

  (I) $\Rightarrow$ (II)\quad For the converse, in the diagram on the left, we use the normal factorizations $f=f_2f_1$ and $h=h_2h_1$, along with the decomposition $g=g_{2}\dot{g} g_{1}$ from (\ref{thm:NormalEpi/MonoDecompositions-Relationship}).
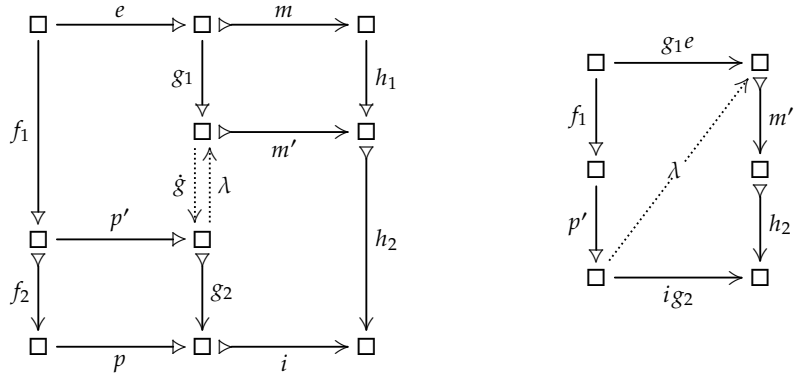
  \begin{center}
    \begin{minipage}[m]{6cm}
      \begin{equation*}
        \xymatrix@R=5ex@C=4em{
        \DiagObj \ar@{-{ >>}}[dd]_{f_1} \ar@{-{ >>}}[r]^-{e} &
        \DiagObj \ar@{-{ >>}}[d]_{g_1} \ar@{{ |>}->}[r]^-{m} &
        \DiagObj \ar@{-{ >>}}[d]^{h_1} \\
        & \DiagObj \ar@{{ |>}->}[r]_-{m'} \ar@{.>}@<-0.5ex>[d]_-{\dot{g}} &
        \DiagObj \ar@{{ |>}->}[dd]^{h_2} \\
        \DiagObj \ar@{{ |>}->}[d]_{f_2} \ar@{-{ >>}}[r]^-{p'} &
        \DiagObj \ar@{{ |>}->}[d]^{g_2} \ar@{.>}@<-0.5ex>[u]_-{\lambda} & 		\\
        \DiagObj  \ar@{-{ >>}}[r]_-{p} &
        \DiagObj \ar@{{ |>}->}[r]_-{i} &
        \DiagObj
        }
      \end{equation*}
    \end{minipage}\qquad
    \begin{minipage}[m]{5cm}
      \begin{equation*}
        \xymatrix@R=5ex@C=4em{
        \DiagObj \ar[r]^-{g_{1}e} \ar@{-{ >>}}[d]_{f_{1}}&
        \DiagObj \ar@{{ |>}->}[d]^{m'} \\
        \DiagObj \ar@{-{ >>}}[d]_{p'} &
        \DiagObj \ar@{{ |>}->}[d]^{h_{2}} \\
        \DiagObj \ar[r]_-{ig_{2}} \ar@{.>}[ruu]|-{\lambda\ } &
        \DiagObj
        }
      \end{equation*}
    \end{minipage}
  \end{center}
  We claim that $\dot{g}$ has an inverse given by the map $\lambda$ which comes from the orthogonality between the strong epimorphism $p'f$ and the monomorphism $h_{2}m'$. This follows from the computations
  \begin{align*}
    h_{2} m' \lambda \dot{g} g_{1}e  & = i g_{2} \dot{g} g_{1}e = h_{2} m' g_{1} e & \Rightarrow\quad \lambda\dot{g} =\IdMap   \\
    i g_{2} \dot{g} \lambda p' f_{1} & = i g_{2}\dot{g} g_{1} e = i g_{2} p' f_{1} & \Rightarrow\quad \dot{g} \lambda = \IdMap
  \end{align*}
  This implies that $g$ is a normal map, and the proof that (I) $\Leftrightarrow$ (II) is complete.
\end{proof}

\bigskip

\begin{subordinate}{}

  \begin{subsubordinate}{On di-exact categories}
    Of a di-exact category, (\ref{term:DiExtensiveConditions}.\ref{Ax:Antinormal->Normal}), we ask that every map which admits an antinormal decomposition be normal. It plays a special role in what follows for the following reasons. We will introduce homological categories and semiabelian categories. To, the set of structural axioms characterizing homological categories, we only need to add one more axiom to characterize semiabelian categories, namely (\ref{term:DiExtensiveConditions}.\ref{Ax:Antinormal->Normal}). Condition (\ref{Ax:ALLNormal}) is stronger and holds, for instance, for abelian categories.

    It is, therefore quite surprising that in di-exact categories most of the basic tool set is available so that effective computations with homology can be carried out. Some of these tools are even available under significantly weaker conditions. So, while we uniformly adopt di-exact categories as a working environment, we will point out where results hold under any of the weaker assumptions listed in (\ref{term:DiExtensiveConditions}). This broadens the scope of those results, in particular toward the inclusion of certain topological varieties. We will discuss this in detail.
  \end{subsubordinate}

  \begin{subsubordinate}{On homological self-duality}
    At first, this condition may look quite specialized. However, we see in Section \ref{sec:HomologicalSelfDuality} that its scope is surprisingly diverse. The terminology is motivated by the fact that, in a homologically self-dual category two dual, but a priori different constructions of `homology' are naturally equivalent; see the Notes in Section \ref{sec:Homology}.
  \end{subsubordinate}

  \begin{subsubordinate}{On the dinversion of antinormal pairs}
    Associated to a \ZExact\ category $\Ctgry{X}$ is the category $\ANPCat{X}$ of antinormal pairs $(\varepsilon,\mu)$. With functorial kernels and cokernels, `dinversion' is a functor
    \begin{equation*}
      B\from \ANPCat{X} \longrightarrow \ANPCat{X}
    \end{equation*}
    which is its own inverse up to unique equivalence. The up-to-equivalence fixed points of $B$ are antinormal pairs which concatenate to short exact sequences.
  \end{subsubordinate}
\end{subordinate}

\begin{exercises}

\begin{exercise}[Constructing di-extensions\ZExactTag]
  \label{thm:DoubleExtensionsConstruct}
  \label{thm:DiExtensionsConstruct}
  Show that each of the following yields a double extension.
  \begin{thmlist}
    \item \label{thm:DiExtensionsConstruct-EqualNormalMonos/Epis}%
    Extending the kernel pair of a normal monomorphism or the cokernel pair of a normal epimorphism yields a di-extension with $0$'s in the antidiagonal corners.
    \item \label{thm:DiExtensionsConstruct-ZRB}%
    Extending a bicartesian square of normal monomorphisms yields a di-extension with a $\ZeroObject$ in the bottom right corner.
    \item \label{thm:DiExtensionsConstruct-ZLT}%
    Extending a bicartesian square of normal epimorphisms yields a di-extension with a $\ZeroObject$ in the top left corner.
    \item \label{thm:DiExtensionsConstruct-ZLB}%
    Extending a bicartesian square in the top right whose edges are normal monomorphisms, respectively normal epimorphisms, as shown on the left below, yields a di-extension with a $\ZeroObject$ in the bottom left corner.
    
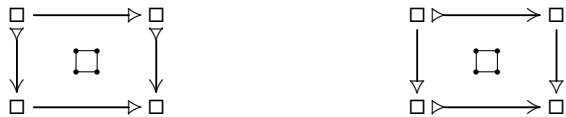
\begin{equation*}
      \xymatrix@R=5ex@C=4em{
      \DiagObj \ar@{-{ >>}}[r] \ar@{{ |>}->}[d] \BiCart{rd} &
      \DiagObj \ar@{{ |>}->}[d] &&
      \DiagObj \ar@{{ |>}->}[r] \ar@{-{ >>}}[d] \BiCart{rd} &
      \DiagObj \ar@{-{ >>}}[d] \\
      \DiagObj \ar@{-{ >>}}[r]  &
      \DiagObj &&
      \DiagObj \ar@{{ |>}->}[r] &
      \DiagObj
      }
    \end{equation*}
    \item \label{thm:DiExtensionsConstruct-ZRT}%
    Extending a bicartesian square in the bottom left whose edges are normal monomorphisms, respectively normal epimorphisms, as shown on the right above yields a di-extension with a $\ZeroObject$ in the top right corner.
  \end{thmlist}
\end{exercise}

\begin{exercise}[Dinversion in $\SetsBsd$]
  \label{exe:Dinversion-Sets_*}
  Show that, in the category $\SetsBsd$ of pointed sets, dinversion preserves normal maps.
\end{exercise}

\begin{exercise}[Antinormal pairs with prescribed properties for $\varepsilon\mu$\ZExactTag]
  \label{exe:AntinormalPairs-emPrescribedProps}%
  In a \ZExact\ category $\Ctgry{X}$, let $(\varepsilon,\mu)$ be an antinormal pair. Using the notation of \eqref{eq:Dinversion} prove the following:
  \begin{thmlist}
    \item If $\varepsilon\mu$ is a normal monomorphism which turns into a normal map under antinormal inversion, then the map $\Ker{\varepsilon}\to \CoKer{\mu}$ is a normal monomorphism as well, and the resulting di-extension has a $\ZeroObject$ in the top left corner.
    \item If $\varepsilon\mu$ is a normal epimorphism which which turns into a normal map under antinormal inversion, then the map $\Ker{\varepsilon}\to \CoKer{\mu}$ is a normal epimorphism as well, and the resulting di-extension has a $\ZeroObject$ in the bottom right corner.
    \item If $\varepsilon\mu$ is an isomorphism which turns into a normal map under antinormal inversion, then the map $\Ker{\varepsilon}\to \CoKer{\mu}$ is an isomorphism as well, and the resulting di-extension has a $\ZeroObject$ in the top left and the bottom right corners.
    \item Decide whether there are pointed categories in which an antinormal pair $(\varepsilon,\mu)$, with $\varepsilon\mu=\IdMap$, need not yield a double extension. \ANKTag
  \end{thmlist}
\end{exercise}

\begin{exercise}[Dinversion of isomorphisms\ANKTag]
  \label{exe:DinversionOfIso}
  For a \ZExact\ category, determine the significance of the condition that dinversion turns every isomorphism into an isomorphism. (Compare with (\ref{thm:ProductRecognition}) and the condition that dinversion preserves all normal maps).
\end{exercise}

\begin{exercise}[Di-extensions in $\Monoids$ and $\Magmas$\ANKTag]
  \label{exe:DoubleExt-Mon/Mag}
  In each of the categories $\Magmas$ of unital magmas or $\Monoids$ of monoids decide if every short exact sequence of short exact sequences is a di-extension.
\end{exercise}

\begin{exercise}[Pushout of normal epimorphisms\ANKTag]
  \label{exe:PushoutNormalEpis-EpicKernelMap}%
  In a \ZExact\ category, consider the morphism of short exact sequences below.
  \begin{equation*}
    \xymatrix@R=5ex@C=4em{
    M \ar@{{ |>}->}[r] \ar[d]_{m} &
    X \ar@{-{ >>}}[r] \ar@{-{ >>}}[d] \ar@{}[rd]|-{\text{(R)}}&
    Q \ar@{-{ >>}}[d] \\
    N \ar@{{ |>}->}[r] &
    Y \ar@{-{ >>}}[r] &
    R
    }
  \end{equation*}
  Assume that all maps in square (R) are normal epimorphisms, and that (R) is a pushout. Is the map $m$ on the left always an epimorphism? - Also formulate the dual of this question, and try to answer it.
\end{exercise}

\begin{exercise}[Antinormal inversion in (topological) varieties\ANKTag]
  \label{exe:Dinversion-TopologicalVarieties}%
  For a variety of algebras $\Ctgry{V}$ and the associated topological variety $\bar{\Ctgry{V}}$ do the following:
  \begin{thmlist}
    \item If $\bar{\Ctgry{V}}$ satisfies any of the conditions in (\ref{term:DiExtensiveConditions}) show that so does $\Ctgry{V}$.
    \item Determine those varieties of algebras in which dinversion preserves normal maps for which dinversion in $\bar{\EuScript{V}}$ also preserves normal maps.
    \item Determine those di-exact varieties  for which  $\bar{\Ctgry{V}}$ is di-exact as well.
  \end{thmlist}
\end{exercise}

\begin{exercise}[Antinormal inversion in varieties\ANKTag]
  \label{exe:AntiNormalInversion-Varieties-Que}%
  Determine in which \ZExact\ varieties dinversion preserves normal maps.
\end{exercise}

\begin{exercise}[Special cases of $(\Prdct{3}{3})$ normal factorizations of normal maps\DExTag]
  \label{exe:(3x3)NormalFactorizations-SpecialCases}
  Given the commutative square below, show the following:
  \begin{equation*}
    \xymatrix@R=5ex@C=4em{
    \square \ar[r]^-{a} \ar[d]_{\alpha} &
    \square \ar[d]^{\beta} \\
    \square \ar[r]_-{b} &
    \square
    }
  \end{equation*}
  \begin{thmlist}
    \item Suppose $a$ and $b$ are normal epimorphisms. If $\alpha$ is a normal map, then $\beta$ is a normal map.
    \item Suppose $a$ and $b$ are normal monomorphisms. If $\beta$ is a normal map, then $\alpha$ is a normal map.
  \end{thmlist}
\end{exercise}

Given a map $f$ in a pointed category, it is not at all clear whether it admits an antinormal decomposition or not. The next exercises provide constructions of an antinormal decomposition in special cases.

\begin{exercise}[Maps in $\AbGrps$ are antinormal]
  \label{exe:AntiNormal-Ab}%
  Show that every morphism in the category $\AbGrps$ of abelian groups (more generally: in an abelian category) admits antinormal decompositions. %
  \index{antinormal!decomposition in abelian category}%
\end{exercise}

\begin{exercise}[Precrossed module admits antinormal decomposition - $\Grps$]
  \label{exe:CrossedModule->AntinormalDecomposition}%
  A \Defn{precrossed module (of groups)} is a group homomorphism $\del\from T\to G$, called the \Defn{boundary map}, together with a left group action $\alpha\from \Prdct{G}{T}\to T$ of $G$ on $T$; $g. t\DefEq \alpha(g,t)$ such that
  \begin{equation*}
    \del(g. t)=g \del(t)g^{-1}
  \end{equation*}
  This condition says that $\del$ is equivariant with respect to $\alpha$ and the conjugation action $\gamma$ of $G$ on itself.
  \begin{thmlist}
    \item The map $\del$ of a precrossed module condition is a normal map.
    \item The map $\mu\from G\rtimes_{\gamma} G\to G$, $\mu(x,g)\DefEq xg$ is a group homomorphism.
    \item An antinormal decomposition of $\del$ is given by the commutative diagram
    
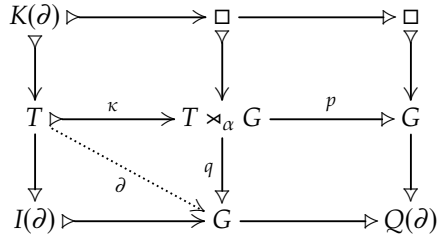
\begin{equation*}
      \xymatrix@R=5ex@C=4em{
      \Ker{\del} \ar@{{ |>}->}[r] \ar@{{ |>}->}[d] &
      \DiagObj \ar@{-{ >>}}[r] \ar@{{ |>}->}[d] &
      \DiagObj \ar@{{ |>}->}[d] \\
      T \ar@{.>}[rd]_-{\del} \ar@{{ |>}->}[r]^-{\kappa} \ar@{-{ >>}}[d] &
      T\rtimes_{\alpha} G \ar@{-{ >>}}[d]_{q} \ar@{-{ >>}}[r]^-{p} &
      G \ar@{-{ >>}}[d] \\
      \Img{\del} \ar@{{ |>}->}[r] &
      G \ar@{-{ >>}}[r] &
      \CoKer{\del}
      }
    \end{equation*}
    in which $\kappa(t)\DefEq (t,1)$, $q\DefEq \mu\Comp (\del\rtimes \IdMapOn{G})$, and $p(t,g)\DefEq g$.
    \item Show that a section of $q$ is given by $s\from G\to T\rtimes_{\alpha} G$.
    \item Show that, as a sectioned decomposition, $q\kappa$ is initial among all antinormal decompositions $\del = rk$ for which $r$ admits a section, say $\sigma$.
  \end{thmlist}
\end{exercise}

\begin{exercise}[Are normal maps antinormal?\ANKTag]
  \label{exe:NormalMapAntiNormal?}
  Determine if there is a pointed category in which there exists a normal map which does not admit an antinormal decomposition. If such a category exists, develop sufficient conditions to identify pointed categories in which every normal map admits an antinormal decomposition.
\end{exercise}
\end{exercises}
\section[Snake Sequences]{Snake Sequences}
\label{sec:Snakes}

We show that various versions of the Snake Lemma hold in di-exact categories; see (\ref{term:DiExtensiveConditions}). We begin with the Pure Snake Lemma and use it as a stepping stone to prove the Classical Snake Lemma and variations thereof.

\begin{lemma}[Pure Snake Lemma\DExTag]
  \label{thm:SnakeLemma-Pure}%
  Consider this morphism of short exact sequences:
  \index{Snake Lemma!pure version}%
  \begin{equation}\label{fig:PureSnakeDiagram}
    \vcenter{
    \xymatrix@!0@R=7ex@C=5em{
    K \ar[d]_-{a} \ar@{{ |>}->}[r]^-{\kappa} &
    X \ar@{=}[d] \ar@{-{ >>}}[r]^-{\xi} &
    Q \ar[d]^-{\rho} \\
    L \ar@{{ |>}->}[r]_-{\xi} &
    X \ar@{-{ >>}}[r]_-{\eta} &
    R
    }
    }
  \end{equation}
  Then $\kappa$ is a normal monomorphism, and $\rho$ is a normal epimorphism. Moreover, the normal factorization of the antinormal composite $t\DefEq ql$ yields the exact sequence below.
  \begin{equation*}
    \xymatrix@R=5ex@C=3em{
    0 \longrightarrow A \ar@{{ |>}->}[r]^-{k} &
    X \ar[r]^-{\beta\xi} \ar@{-{ >>}}[d] &
    Q \ar@{-{ >>}}[r]^-{\rho} &
    R \longrightarrow 0 \\
    & \CoKer{\kappa} \ar[r]^-{\partial^{-1}}_-{\cong} &
    \Ker{\rho} \ar@{{ |>}->}[u]
    }
  \end{equation*}
  This sequence is functorial with respect to morphisms of diagrams of the form \eqref{fig:PureSnakeDiagram}.
\end{lemma}
\begin{proof}
  That $\kappa$ and $\rho$ are normal monomorphism and normal epimorphism, respectively, follows via commutativity and (\ref{thm:NormalMono-Props}), (\ref{thm:NormalEpi-Props}). Next, consider the following rearrangement of the information from the given morphism of short exact sequences:
  \begin{equation*}
    \xymatrix@!0@R=8ex@C=5em{
    A \ar@{=}[d] \ar@{{ |>}->}[r]^-{\kappa} &
    X \ar@{.>}[rd]|-{t} \ar@{{ |>}->}[d]_-{l} \ar@{.{ >>}}[r]^-{\tilde{q}} &
    \Ker{\rho} \ar@{{ |>}->}[d]^-{\mu} \\
    A \ar@{{ |>}->}[r]_-{k} \ar@{-{ >>}}[d] &
    B  \ar@{-{ >>}}[r]_-{q} \ar@{-{ >>}}[d]_{r} &
    C \ar@{-{ >>}}[d]^{\rho} \\
    0 \ar@{{ |>}->}[r] &
    R \ar@{=}[r] &
    R
    }
  \end{equation*}
  In a homologically self-dual category, the antinormally decomposed map $t=ql$ has the normal factorization presented on the top right; see (\ref{fig:DiExtensionTypes}). This proves the claim. Functoriality follows from the uniqueness property of the normal factorization of $t$.
\end{proof}

Repeated use of the Pure Snake Lemma yields the Classical Snake Lemma. It plays a pivotal role in homological algebra.

\begin{theorem}[Classical Snake Lemma\DExTag]
  \label{thm:SnakeLemma-Classical}%
  Suppose the maps $\kappa$, $\xi$, $\rho$ in the following morphism of short exact sequences are normal.
  \index{Snake Lemma!classical - conditionally extensive cats}%
  \begin{equation*}
    \xymatrix@R=5ex@C=3em{
    K \ar[d]_{\kappa} \ar@{{ |>}->}[r]^-{k} &
    X \ar[d]_{\xi} \ar@{-{ >>}}[r]^-{q} &
    Q \ar[d]^{\rho} \\
    L \ar@{{ |>}->}[r]_-{l} &
    Y \ar@{-{ >>}}[r]_-{r} &
    R
    }
  \end{equation*}
  Then the kernels and cokernels of $\kappa$, $\xi$, $\rho$ form a six-term exact sequence
  \begin{equation*}
    \xymatrix@R=5ex@C=2.2em{
    0 \longrightarrow \Ker{\kappa} \ar@{{ |>}->}[r]^-{k^{\ast}} &
    \Ker{\xi} \ar[r]^-{q^{\ast}} &
    \Ker{\rho} \ar[r]^-{\partial} &
    \CoKer{\kappa} \ar[r]^-{l_{\ast}} &
    \CoKer{\xi} \ar@{-{ >>}}[r]^-{r_{\ast}} &
    \CoKer{\rho} \longrightarrow 0
    }
  \end{equation*}
  The six-term exact sequence depends functorially on morphisms of the underlying diagrams of short exact sequences.
\end{theorem}
\begin{proof}
  Compare \cite[p.~47f]{Borceux-Semiab}\cite[p.~297f]{FBorceuxDBourn2004}. We reduce the proof of the classical Snake Lemma to three applications of the Pure Snake Lemma.

  \emph{Step 1: Comparison with the extension of kernels}\quad In the commutative diagram below, via (\ref{thm:NormalMap-FactorViaAntinormalPair}), the composite $K\to \Img{\kappa}\to \Img{\xi}$ is the normal factorization of the antinormal composite $K\to X\to \Img{\xi}$. By (\ref{thm:DiExtensionFromAN->N.Pairs}), we obtain a di-extension of $\Ker{\kappa}$ over $S$.
  \begin{equation*}
    \xymatrix@!0@R=6.5ex@C=3em{
    \Ker{\kappa} \ar@{{ |>}->}[rr]^-{k^{\ast}} \ar@{{ |>}->}[d] \PullLU{rrd} &&
    \Ker{\xi} \ar@{{ |>}->}[d] \ar@{-{ >>}}[rr] \ar@/^2ex/[rrrr]^-{q^{\ast}} &&
    J \ar@{{ |>}->}[rr] \ar@{{ |>}->}[d] &&
    \Ker{\rho} \ar@{-{ >>}}[rr] \ar@{{ |>}->}[d] &&
    \CoKer{q^{\ast}} \\
    K \ar@{{ |>}->}[rr]_-{\alpha} \ar[dd]_-{\kappa} \ar@{-{ >>}}[rd] &&
    X \ar[dd]_(.3){\xi} \ar@{-{ >>}}[rd] \ar@{-{ >>}}[rr]^-{q} &&
    Q \ar[dd]_(.3){\rho} \ar@{=}[rr] \ar@{-{ >>}}[rd] &&
    Q \ar[dd]_(.3){\rho} \ar@{-{ >>}}[rd] \\
    & \Img{\kappa} \ar@{{ |>}->}[ld] \ar@{{ |>}->}[rr]|\hole &&
    \Img{\xi} \ar@{{ |>}->}[ld] \ar@{-{ >>}}[rr]|\hole &&
    S \ar[dl]^(.4){\mu} \ar@{-{ >>}}[rr]|\hole_(.3){\tau} &&
    \Img{\rho} \ar@{{ |>}->}[ld] \\
    L \ar@{-{ >>}}[d] \ar@{{ |>}->}[rr] &&
    Y \ar@{-{ >>}}[d] \ar@{-{ >>}}[rr] &&
    R \ar@{-{ >>}}[d] \ar@{=}[rr] &&
    R \ar@{-{ >>}}[d] \\
    \CoKer{\kappa} \ar[rr]_-{l_{\ast}} &&
    \CoKer{\xi} \ar@{-{ >>}}[rr]_-{r_{\ast}} &&
    \CoKer{\rho} \ar@{=}[rr] &&
    \CoKer{\rho}
    }
  \end{equation*}
  Working with the universal properties of the kernels and cokernels, combined with their monomorphic and epimorphic properties, shows that the supplemented maps exist and are unique rendering the entire diagram commutative. So, $q^{\ast}$ is a normal map, and the Pure Snake Lemma yields a functorial isomorphism $\CoKer{q^{\ast}}\to \Ker{\tau}$.

  \emph{Step 2: Comparison with the extension of cokernels}\quad The dual construction yields a di-extension of $T$ over $\CoKer{\rho}$.
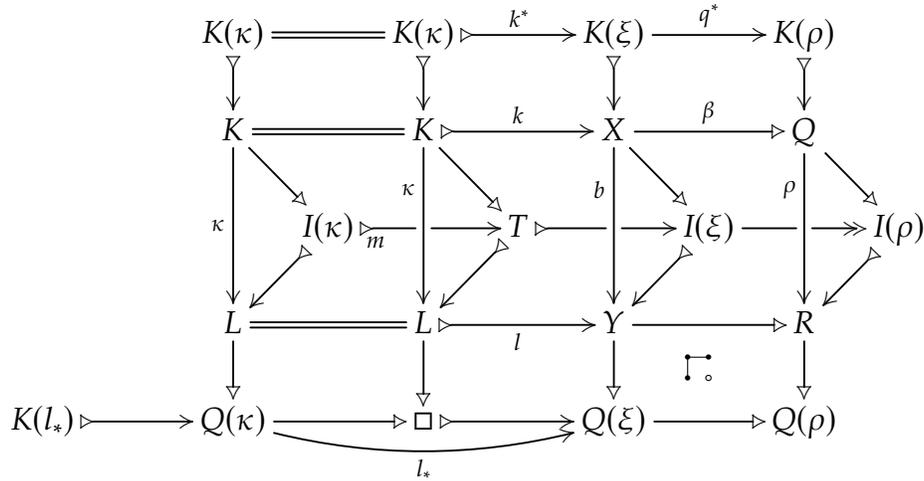
  \begin{equation*}
    \xymatrix@!0@R=6.5ex@C=3em{
    && \Ker{\kappa} \ar@{{ |>}->}[d] \ar@{=}[rr] &&
    \Ker{\kappa} \ar@{{ |>}->}[d] \ar@{{ |>}->}[rr]^-{k^{\ast}} &&
    \Ker{\xi} \ar[rr]^-{q^{\ast}} \ar@{{ |>}->}[d] &&
    \Ker{\rho} \ar@{{ |>}->}[d] \\
    && K \ar@{=}[rr] \ar[dd]_-{\kappa} \ar@{-{ >>}}[rd] &&
    K \ar[dd]_(.3){\kappa} \ar@{-{ >>}}[rd] \ar@{{ |>}->}[rr]^-{k} &&
    X \ar@{-{ >>}}[rr]^-{\beta} \ar[dd]_(.3){b} \ar@{-{ >>}}[rd] &&
    Q \ar[dd]_(.3){\rho} \ar@{-{ >>}}[rd] \\
    &&& \Img{\kappa} \ar@{{ |>}->}[ld] \ar@{{ |>}->}[rr]|\hole_(.25){m} &&
    T \ar@{{ |>}->}[ld] \ar@{{ |>}->}[rr]|\hole &&
    \Img{\xi} \ar@{{ |>}->}[dl] \ar@{-{>>}}[rr]|\hole &&
    \Img{\rho} \ar@{{ |>}->}[ld] \\
    && L \ar@{=}[rr] \ar@{-{ >>}}[d] &&
    L \ar@{-{ >>}}[d] \ar@{{ |>}->}[rr]_-{l} &&
    Y \ar@{-{ >>}}[d] \ar@{-{ >>}}[rr] \PushRD{rrd} &&
    R \ar@{-{ >>}}[d] \\
    \Ker{l_{\ast}} \ar@{{ |>}->}[rr] &&
    \CoKer{\kappa} \ar@{-{ >>}}[rr] \ar@/_2ex/[rrrr]_-{l_{\ast}} &&
    \DiagObj \ar@{{ |>}->}[rr] &&
    \CoKer{\xi} \ar@{-{ >>}}[rr] &&
    \CoKer{\rho}
    }
  \end{equation*}
  Working with the universal properties of the kernels and cokernels shows that the supplemented maps exist and are unique rendering the entire diagram commutative (the pullback property of the square with corners $T$ and $Y$ is essential here). The Pure Snake Lemma yields a functorial isomorphism $\Ker{l_{\ast}}\to \CoKer{m}$.

  \emph{Step 3: Connecting steps 1 and 2}\quad From steps 1 and 2, we obtain the diagram of short exact sequences below.
  
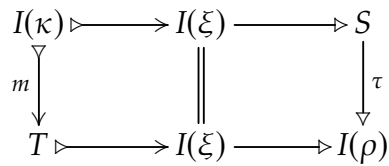
\begin{equation*}
    \xymatrix@R=5ex@C=3em{
    \Img{\kappa} \ar@{{ |>}->}[d]_{m} \ar@{{ |>}->}[r] &
    \Img{\xi} \ar@{=}[d] \ar@{-{ >>}}[r] &
    S \ar@{-{ >>}}[d]^-{\tau} \\
    T \ar@{{ |>}->}[r] &
    \Img{\xi} \ar@{-{ >>}}[r] &
    \Img{\rho}
    }
  \end{equation*}
  Combining the epimorphic property of the cokernel $K \NEpi \Img{\kappa}$ with the monomorphic property of the kernel $\Img{\xi}\NMono Y$, we see that the left hand square commutes. Dually, combining the epimorphic property of the cokernel $B \NEpi \Img{b}$ with the monomorphic property of the kernel $\Img{\rho}\NMono R$, we see that the right hand square commutes. The Pure Snake Lemma yields a functorial isomorphism $\Ker{\tau}\to \CoKer{m}$.

  \emph{Step 4: Assembly}\quad The functorial isomorphisms we found in steps 1, 2, 3, may be composed as in the diagram below to yield the desired isomorphism $\partial$:
  \begin{equation*}
    \xymatrix@R=5ex@C=3em{
    \CoKer{q^{\ast}} \ar[d]^{\cong}_{\text{Step 1}} \ar[r]^-{\partial}_-{\cong} &
    \Ker{\xi_{\ast}} \\
    \Ker{\tau} \ar[r]^-{\cong}_-{\text{Step 3}} &
    \CoKer{m} \ar[u]^-{\cong}_-{\text{Step 2}}
    }
  \end{equation*}
  The arguments given are functorial with respect to morphisms of diagrams, as in the statement of the Classical Snake Lemma. So, its proof is complete.

  To extend the scope of the Snake Lemma, we explain how the argument given above can be refined so that the hypotheses can be weakened to: ``(co)subnormal maps in $\Ctgry{X}$ which admit an antinormal decomposition are normal'', combined with the assumption that dinversion preserves normal maps in $\Ctgry{X}$; see (\ref{term:DiExtensiveConditions}). This refinement affects step 1 and, by duality, step 2. Here is an excerpt of the relevant part of the diagram used in step 1.
  \begin{equation*}
    \xymatrix@!0@R=6.5ex@C=3em{
    \Ker{\kappa} \ar@{{ |>}->}[rr]^-{k^{\ast}} \ar@{{ |>}->}[d] \PullLU{rrd} &&
    \Ker{\xi} \ar@{{ |>}->}[d] \\
    K \ar@{{ |>}->}[rr]_-{k} \ar@/_4ex/[ddd]_-{\kappa} \ar@{-{ >>}}[rd] &&
    X \ar@{-{ >>}}[rd] \\
    & \Img{\kappa} \ar@{{ |>}->}@/_2ex/[ldd] \ar@{{ >}->}[rr] \ar@{{ |>}.>}[d]^{\gamma} &&
    \Img{\xi} \ar@{=}[d] \\
    & P \ar@{{ |>}->}[rr] \ar@{{ |>}->}[ld] \PullLU{rd} &&
    \Img{\xi} \ar@{{ |>}->}[ld] \\
    L  \ar@{{ |>}->}[rr] &&
    Y
    }
  \end{equation*}
  In the above diagram we used (\ref{thm:Pullback/Pushout-Existence}) to insert the pullback of $L\NMono Y$ along $\Img{\xi} \NMono Y$ into the square $\Img{\kappa}\rightrightarrows Y$. The comparison map $\gamma$ is a normal monomorphism by (\ref{thm:NormalMono-Props}). Thus we have a factorization of the antinormal composite $K \NMono X \NEpi \Img{\xi}$ as a subnormal map. In a sub-di-exact category, such a map is normal, and dinversion turns it into a normal map. Thus, we obtain the di- extension we were working with in step 1.

  Dually, in step 2, we find that the antinormal composite $\Img{\xi}\NMono Y \NEpi R$ is cosubnormal, hence normal.
\end{proof}

\begin{terminology}
  The morphism $\partial$ in the statement of the Classical Snake Lemma (\ref{thm:SnakeLemma-Classical}) is also called the \Defn{snake map/morphism} or the \Defn{connecting map/morphism}. %
  \index{connecting map!in Snake Lemma}%
  \index{snake!map/morphism}%
\end{terminology}

\begin{corollary}[Relaxed Snake Lemma\DExTag]
  \label{thm:SnakeLemma-Relaxed}%
  \cite[p.~47f]{Borceux-Semiab}\cite[p.~297f]{FBorceuxDBourn2004}\quad Consider the commutative diagram below. %
  \index{Snake Lemma!relaxed version}%
  \begin{equation*}
    \xymatrix@R=5ex@C=3em{
    & K \ar[d]_{\kappa} \ar[r]^-{k} &
    X \ar[d]_{\xi} \ar@{-{ >>}}[r]^-{q} &
    Q \ar[d]^{\rho} \ar[r] &
    0 \\
    0 \ar[r] &
    L \ar@{{ |>}->}[r]_-{l} &
    Y \ar[r]_-{r} &
    R \\
    }
  \end{equation*}
  Assume that the maps $k$, $\kappa$, $\xi$, $\rho$, $r$ are normal, the top sequence is exact, and the bottom sequence is coexact. Then there is this functorial six-term exact sequence:
  \begin{equation*}
    \xymatrix@R=5ex@C=2em{
    \Ker{\kappa} \ar@{.>}[r]^-{k^{\ast}} &
    \Ker{\xi} \ar@{.>}[r]^-{q^{\ast}} &
    \Ker{\rho} \ar@{.>}[r]^-{\partial} &
    \CoKer{\kappa} \ar@{.>}[r]^-{l_{\ast}} &
    \CoKer{\xi} \ar@{.>}[r]^-{r_{\ast}} &
    \CoKer{\rho}
    }
  \end{equation*}
\end{corollary}
\begin{proof}
  We reduce the given situation to the Classical Snake Lemma \ref{thm:SnakeLemma-Classical}:
  \begin{equation*}
    \xymatrix@R=5ex@C=3em{
    \Ker{\kappa} \ar@{{ |>}->}[d] \ar@{-{ >>}}[r] \PullLU{rd} \ar@/^2ex/[rr]^-{k^{\ast}} &
    \Ker{\kappa'} \ar@{{ |>}->}[r] \ar@{{ |>}->}[d] &
    \Ker{\xi} \ar@{{ |>}->}[d] \ar[r] \ar@{-{ >>}}@/^2ex/[rr]^-{q^{\ast}} &
    \Ker{\rho'} \ar@{{ |>}->}[d] \ar[r]_-{\cong} &
    \Ker{\rho} \ar@{{ |>}->}[d] \\
    K \ar[d]_-{\kappa} \ar@{-{ >>}}[r]^-{u} &
    K' \ar@{{ |>}->}[r]^-{k'} \ar[d]_-{\kappa'} &
    X \ar[d]_{b} \ar@{-{ >>}}[r]^-{q} &
    Q \ar@{=}[r] \ar[d]_-{\rho'} &
    Q \ar[d]^-{\rho} \\
    L \ar@{=}[r] \ar@{-{ >>}}[d] &
    L \ar@{{ |>}->}[r]_-{l} \ar@{-{ >>}}[d] &
    Y \ar@{-{ >>}}[r]_-{r'} \ar@{-{ >>}}[d] &
    R' \ar@{-{ >>}}[d] \ar@{{ |>}->}[r]_-{v} &
    R \ar@{-{ >>}}[d] \\
    \CoKer{\kappa} \ar[r]^-{\cong} \ar@{{ |>}->}@/_2ex/[rr]_-{l_{\ast}}&
    \CoKer{\kappa'} \ar[r] &
    \CoKer{\xi} \ar@{-{ >>}}[r] \ar@/_2ex/[rr]_-{r_{\ast}}&
    \CoKer{\rho'} \ar@{{ >}->}[r]^-{x} &
    \CoKer{\rho}
    }
  \end{equation*}
  The composite $k'u$ is the normal factorization of the normal map $k$. So, $k' = \Ker{q}$. The universal property of the kernel $l$ yields $k'$ unique with $l \kappa'u = l \kappa$. We conclude $\kappa'u=\kappa$ via the monomorphic property of $l$.  With (\ref{thm:FactorNormalQuotientObject->Normal}), we see that $\kappa'$ is a normal map.

  The composite $vr'$ is the normal factorization of $r$, coming from the coexactness in position~$Y$. Via the cokernel $q$ we obtain $\rho'$, unique with $\rho' q=r' \xi$. We conclude $v\rho'=\rho$ from the epimorphic property of $\beta$. With (\ref{thm:FactorNormalThroughSubObject->Normal}), we see that $\rho'$ is a normal map.

  Thus we may apply the classical Snake Lemma to the morphism $(\kappa',\xi,\rho')$ of short exact sequences, to obtain the following exact sequence:
  \begin{equation*}
    0 \to \Ker{\kappa'} \NMono \Ker{\xi} \longrightarrow \Ker{\rho'} \XRA{\partial} \CoKer{\kappa'} \XRA{x} \CoKer{\xi} \NEpi \CoKer{\rho'} \to 0
  \end{equation*}
  Using the normal factorizations of $\kappa$ and $\kappa'$ we arrive at the following commutative diagram:
  
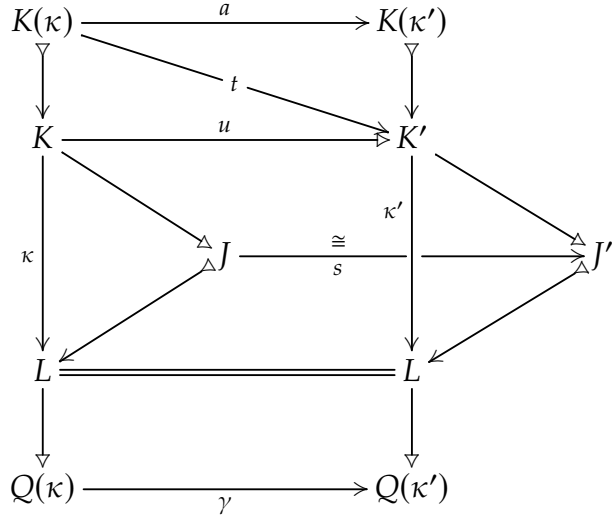
\begin{equation*}
    \xymatrix@R=5ex@C=4em{
    \Ker{\kappa} \ar[rr]^-{a} \ar@{{ |>}->}[d] \ar[rrd]|-{\ t\ }&&
    \Ker{\kappa'} \ar@{{ |>}->}[d] \\
    K \ar@{-{ >>}}[rr]^-{u} \ar@{-{ >>}}[rd] \ar[dd]_{\kappa} &&
    K' \ar@{-{ >>}}[rd] \ar[dd]_(0.3){\kappa'} \\
    & J \ar[rr]|\hole_(.3){s}^(0.3){\cong} \ar@{{ |>}->}[ld] &&
    J' \ar@{{ |>}->}[ld]  \\
    L \ar@{=}[rr] \ar@{-{ >>}}[d] &&
    L \ar@{-{ >>}}[d] \\
    \CoKer{\kappa} \ar[rr]_-{\gamma} &&
    \CoKer{\kappa'}
    }
  \end{equation*}
  The map $s$ is an isomorphism because it is a normal monomorphism by commutativity of the square $J\rightrightarrows L$, and it is an epimorphism by commutativity of the square $K\rightrightarrows J'$. So, the map $\gamma$ is an isomorphism as well. Via the recognition criterion for homologically self-dual categories (\ref{thm:HomologicalSelfDuality-Recognize-II}), $a$ is a normal epimorphism. Thus $\Ker{\kappa}\to \Ker{\xi}\to \Ker{\rho}$ is exact at $\Ker{\xi}$, and $\CoKer{\kappa'}$ may be replaced by $\CoKer{\kappa}$ via the isomorphism we found.

  Dual argumentation completes the validation of the relaxed snake sequence. Checking that all constructions involved are functorial completes the proof.
\end{proof}

For later use, we develop a number generalizations of the Pure Snake Lemma \ref{thm:SnakeLemma-Pure}. We start from a commutative diagram in which the composites $\beta\alpha$ and $\eta\xi$ are $\ZeroMap$.
\begin{equation}\label{fig:PureSnakeDiagram-Generalized}
  \vcenter{
  \xymatrix@!0@R=7ex@C=5em{
  A \ar[d]_-{a} \ar[r]^-{\alpha} &
  B \ar@{=}[d] \ar@{-{ >>}}[r]^-{\beta} &
  C \ar[d]^-{c}  \\
  X \ar@{{ |>}->}[r]_-{\xi} &
  B \ar[r]_-{\eta} &
  Z
  }
  }
\end{equation}
We further assume that $\beta$ is a normal epimorphism and $\xi$ is a normal monomorphism.
The above diagram has the following expansion which we will refer to.
\begin{equation*}
  \xymatrix@!0@R=7ex@C=7em{
  &&&
  \Ker{\rho} \ar@{{ |>}->}[d] \ar@{{ |>}->}[r]^{v} &
  \Ker{c} \ar@{{ |>}->}[d] \\
  A \ar[d]_{a} \ar@{.>}[r]_-{\alpha'} \ar@/^3ex/[rr]^-{\alpha}&
  \Ker{\beta} \ar@{{ |>}->}[d]_-{\kappa} \ar@{{ |>}->}[r] &
  B \ar@{=}[d] \ar@{-{ >>}}[r]^-{\beta} &
  C \ar@{-{ >>}}[d]^-{\rho}  \ar@{=}[r] &
  C \ar[d]^{c} \\
  X \ar@{=}[r] \ar@{-{ >>}}[d] &
  X \ar@{{ |>}->}[r]_-{\xi} \ar@{-{ >>}}[d]  &
  B \ar@{-{ >>}}[r] \ar@/_3ex/[rr]_-{\eta} &
  \CoKer{\xi} \ar@{.>}[r]^-{\eta'}&
  Z \\
  \CoKer{a} \ar@{-{ >>}}[r]_-{u} &
  \CoKer{\kappa}
  }
\end{equation*}
\begin{corollary}[Generalized pure snake I\DExTag]
  \label{thm:PureSnake-Generalized-I}
  In the setup of diagram \eqref{fig:PureSnakeDiagram-Generalized}, the following hold and are functorial with respect to morphisms of diagrams of type \eqref{fig:PureSnakeDiagram-Generalized}:
  \begin{thmlist}
    \item The top row in the diagram below is exact, and the map $d$ is normal.
    
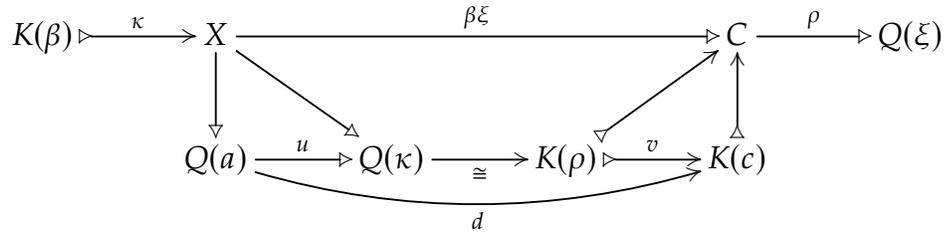
\begin{equation*}
      \xymatrix@R=5ex@C=3em{
      \Ker{\beta} \ar@{{ |>}->}[r]^-{\kappa} &
      X \ar@{-{ >>}}[rrr]^-{\beta\xi} \ar@{-{ >>}}[d] \ar@{-{ >>}}[dr] &&&
      C \ar@{-{ >>}}[r]^-{\rho} &
      \CoKer{\xi} \\
      & \CoKer{a} \ar@{-{ >>}}[r]^-{u} \ar@/_3ex/[rrr]_-{d} &
      \CoKer{\kappa} \ar[r]_-{\cong} &
      \Ker{\rho} \ar@{{ |>}->}[r]^-{v} \ar@{{ |>}->}[ru] &
      \Ker{c} \ar@{{ |>}->}[u]
      }
    \end{equation*}
    \item If the top row in \eqref{fig:PureSnakeDiagram-Generalized} is exact, then $a$ is a normal map, $u$ is an isomorphism, and $d$ is a normal monomorphism. Moreover, the sequence below is exact.
    \begin{equation*}
      \xymatrix@R=5ex@C=3em{
      \Ker{\beta} \ar@{{ |>}->}[r]^-{\kappa} &
      X \ar[r] &
      \Ker{c} \ar[r] &
      \CoKer{\xi} \ar[r]^-{\eta'} &
      Z
      }
    \end{equation*}
    \item If the bottom row in \eqref{fig:PureSnakeDiagram-Generalized} is coexact, then $c$ is a normal map, $v$ is an isomorphism, and $d$ is a normal epimorphism. Moreover, the sequence below is coexact.
    \begin{equation*}
      \xymatrix@R=5ex@C=4em{
      A \ar[r]^-{\alpha'} &
      \Ker{\beta} \ar[r] &
      \CoKer{a} \ar[r] &
      C \ar@{-{ >>}}[r]^-{\rho} &
      \CoKer{\xi}
      }
    \end{equation*}
  \end{thmlist}
\end{corollary}
\begin{proof}
  Upon confirming the normal monomorphism and epimorphism properties of the maps $\kappa$ and $\rho$, we see that the Pure Snake Lemma (\ref{thm:SnakeLemma-Pure}) is applicable to the morphism $(\kappa,\IdMapOn{B},\rho)$, and claim (i) follows since $\CoKer{\kappa}\cong\Ker{\rho}$ while $u$ is a normal epimorphism and $v$ a normal monomorphism.

  (ii)\quad If the top row in \eqref{fig:PureSnakeDiagram-Generalized} is exact then $\alpha'$ is a normal epimorphism. So $a=\kappa \alpha'$ is a normal map, and (\ref{thm:PushoutRecognize-Categorical}) the bottom left square is a pushout. It follows that $u$ is an isomorphism  and, hence, $d$ is a normal monomorphism.

  Now, we observe that we can apply (i) to the morphism $(v,\IdMapOn{C},\eta')$ and splice the result into the exact sequence from (i). The result is:
  \begin{equation*}
    \xymatrix@!@R=-4ex@C=-2em{
    \Ker{\beta} \ar@{{ |>}->}[r]^-{\kappa} &
    X \ar[rr] \ar@{-{ >>}}[rd] &&
    \Ker{c} \ar[rr] \ar@{-{ >>}}[rd]&&
    \CoKer{\xi} \ar[rr]^{\eta'} &&
    Z \\
    && \CoKer{\kappa}\cong \Ker{\rho} \ar@{{ |>}->}[ru]_{v} &&
    \CoKer{v}\cong\Ker{\eta'} \ar@{{ |>}->}[ru]
    }
  \end{equation*}
  So, the top row is exact as claimed.

  (iii) The argument is dual to that of (ii).
\end{proof}

\begin{corollary}[Relaxed Pure Snake Lemma\DExTag]
  \label{thm:SnakeLemma-RelaxedPure}%
  In the setup of diagram \eqref{fig:PureSnakeDiagram-Generalized}, assume that the top row is exact and that the bottom row is coexact. Then $a$ and $c$ are normal maps, and the normal factorization of the antinormal composite $t\DefEq\beta\xi$ yields
  \index{Snake Lemma!pure version, relaxed}%
  \begin{equation*}
    \xymatrix@R=5ex@C=3em{
    A \ar[r]^-{a} &
    X \ar[r]^-{\beta\xi} \ar@{-{ >>}}[d] &
    C \ar[r]^-{c} &
    Z \\
    & \CoKer{a} \ar[r]_-{\cong}^-{d} &
    \Ker{c} \ar@{{ |>}->}[u]
    }
  \end{equation*}
  as a functorial exact and coexact sequence. \NoProof%
\end{corollary}

\begin{corollary}[\DExTag]
  \label{thm:PureSnake-Generalized-II}
  In the setup of diagram \eqref{fig:PureSnakeDiagram-Generalized}, suppose at the following two conditions hold.
  \begin{thmlist}
    \item $\beta=\CoKerMap{\alpha}$ and $\xi=\KerMap{\eta}$.
    \item The maps $a$ and $c$ are normal.
  \end{thmlist}
  Then the sequence below is exact and coexact, and is functorial with respect to morphisms of diagrams of type \eqref{fig:PureSnakeDiagram-Generalized}.
  \begin{equation*}
    A \XRA{a'} \Ker{\beta} \longrightarrow \CoKer{a} \longrightarrow \Ker{c} \longrightarrow \CoKer{\xi} \XRA{\eta'} Z
  \end{equation*}
\end{corollary}
\begin{proof}
  We refine the expansion of \eqref{fig:PureSnakeDiagram-Generalized} as follows.
  \begin{equation*}
    \xymatrix@!0@R=7ex@C=7em{
    A \ar@/_{3ex}/[dd]_{a} \ar@{-{ >>}}[d]^(0.4){a'} \ar@{.>}[r]_-{\alpha'} \ar@/^3ex/[rr]^-{\alpha}&
    \Ker{\beta} \ar@{=}[d] \ar@{{ |>}->}[r]_-{\alpha''} &
    B \ar@{=}[d] &
    \Ker{\rho}  \ar@{{ |>}->}[r]^{v} \ar@{{ |>}->}[d]&
    \Ker{c} \ar@{{ |>}->}[d]^{k} \\
    \Img{a} \ar@{{ |>}->}[d]^(0.4){a''} \ar@{{ |>}->}[r] &
    \Ker{\beta} \ar@{{ |>}->}[d]_-{\kappa} \ar@{{ |>}->}[r]^-{\alpha''} &
    B \ar@{=}[d] \ar@{-{ >>}}[r]^-{\beta} &
    C \ar@{=}[r] \ar@{-{ >>}}[d]^(.4){\rho}&
    C \ar@{-{ >>}}[d]_(0.4){c'} \ar@/^3ex/[dd]^{c} \\
    X \ar@{=}[r] \ar@{-{ >>}}[d]_{r} &
    X \ar@{{ |>}->}[r]_-{\xi} \ar@{-{ >>}}[d]  &
    B \ar@{-{ >>}}[r]  \ar@{=}[d] &
    \CoKer{\xi} \ar@{-{ >>}}[r] \ar@{=}[d] &
    \Img{c} \ar@{{ |>}->}[d]_{c''} \\
    \CoKer{a} \ar@{-{ >>}}[r]_-{u} &
    \CoKer{\kappa} &
    B \ar@{-{ >>}}[r]^-{\eta''} \ar@/_3ex/[rr]_-{\eta}&
    \CoKer{\xi} \ar@{.>}[r]^-{\eta'} &
    Z
    }
  \end{equation*}
  Upon confirming that the diagram commutes and that the maps have the indicated properties, a $3$-fold application of the Pure Snake Lemma yields the following exact and coexact sequence which is spliced from short exact sequences.
  \begin{equation*}
    \xymatrix@R=6ex@!C=1.5em{
    A \ar[rr]^-{\alpha'} \ar@{-{ >>}}[rd]^{a'} &&
    \Ker{\beta} \ar[rr]^{r\kappa} \ar@{-{ >>}}[rd] &&
    \CoKer{a} \ar[rr]^{vu} \ar@{-{ >>}}[rd]^{u} &&
    \Ker{c} \ar[rr]^{\rho k} \ar@{-{ >>}}[rd] &&
    \CoKer{\xi} \ar[rr]^-{\eta'} \ar@{-{ >>}}[rd] &&
    Z \\
    & \Img{a} \ar@{{ |>}->}[ru] &&
    \CoKer{\alpha'}\cong \Ker{u} \ar@{{ |>}->}[ru] &&
    \CoKer{\kappa}\cong \Ker{\rho} \ar@{{ |>}->}[ru]^{v}&&
    \CoKer{v}\cong \Ker{\eta'} \ar@{-{ >>}}[ru] &&
    \Img{c} \ar@{{ |>}->}[ru]
    }
  \end{equation*}
  This proves the claim.
\end{proof}

\begin{subordinate}[Comments]{On minimal conditions}
  In the proof of (\ref{thm:SnakeLemma-Pure}), we actually only assume that we are working in a homologically self-dual category (\ref{term:DiExtensiveConditions}). Its generalizations (\ref{thm:PureSnake-Generalized-I}),  (\ref{thm:SnakeLemma-RelaxedPure}), and (\ref{thm:PureSnake-Generalized-II}) rely only on properties of \ZExact\ categories and the validity of the Pure Snake Lemma. So, they also hold under this weakened assumption.

  The proof of the Classical Snake Lemma (\ref{thm:SnakeLemma-Classical}) involves four steps. As presented, Steps 1 and 2 are valid with a straightforward  argument in di-exact categories. However, more subtly, we also presented an argument which is valid  \ZExact\ categories in which antinormal composites which are (co)subnormal are actually normal and that normal maps are preserved by dinversion; that is in sub-di-exact categories (\ref{term:DiExtensiveConditions}). Step 3 only relies on homological self-duality. So, the Classical Snake Lemma and its corollary (\ref{thm:SnakeLemma-Relaxed}) hold in sub-di-exact categories as well. We present a further analysis of sub-di-exact categories in Section \ref{sec:Sub-di-exact Categories}.
\end{subordinate}

\begin{exercises}
\begin{exercise}[Generalized pure snake III\DExTag]
  \label{thm:SnakeLemma-Pure-Expansion}
  In the setup of diagram \eqref{fig:PureSnakeDiagram-Generalized}, suppose at least one of the following two conditions hold.
  \begin{thmlist}
    \item The map $a$ is normal, and the bottom row is coexact.
    \item The map $c$ is normal, and the top row is exact.
  \end{thmlist}
  Then the sequence below is exact and coexact, and is functorial with respect to morphisms of diagrams of type \eqref{fig:PureSnakeDiagram-Generalized}.
  \begin{equation*}
    A \XRA{a} X \XRA{\beta\xi} C \XRA{c} Z
  \end{equation*}
\end{exercise}
\end{exercises}
\section[Homology]{Homology}
\label{sec:Homology}

We consider chain complexes in \ZExact\ categories, and we introduce two types of homology as a measure of failure of exactness. Technically, these two definitions of homology are dual to one another. In a di-exact category, and under mild conditions required of the chain complex, these two constructions of homology are related as is explained in Lemma \ref{thm:Homology-CoKerHomology->KerHomology}. For normal chain complexes, this relationship turns into an isomorphism.

Subsequently, we confirm that either type of homology serves its purpose as expected. Under the conditions formulated in (\ref{thm:H^c(C)=0<->Exact}) and (\ref{thm:H^k(C)=0<=>Exact}), homology vanishes if and only if the chain complex is exact.

\begin{definition}[Chain complex\ZExactTag]%
  \label{def:ChainComplex}%
  A \Defn{chain complex} is a $\ZNr$-indexed sequence of maps %
  \index{chain complex}%
  \begin{equation*}
    \xymatrix@R=5ex@C=3em{
    \cdots \ar[r] &
    C_{n+1} \ar[r]^-{d_{n+1}} &
    C_n \ar[r]^-{d_n} &
    C_{n-1} \ar[r] &
    \cdots
    }
  \end{equation*}
  such that $d_n\Comp d_{n+1} = 0$  for each  $n\in \ZNr$.
\end{definition}

\begin{notation}[Chain complex]%
  \label{not:ChainComplex}%
  We write  $ (C,d^C)$  for a chain complex, and shorten the notation to  $(C,d)$  or even $C$, if there is no risk of confusion. The family of maps $d^C$ is called the \Defn{differential} or the \Defn{boundary operator} of the chain complex.%
  \index[not]{c!$(C,d^C)$\IndSep chain complex $C$ with differential $d^C$}%
  \index{chain complex!differential}\index{chain complex!boundary operator}\index{differential of a chain complex}%
  \index{boundary operator of a chain complex}%
\end{notation}

Thus every long exact sequence is a chain complex. The converse is far from true. To see how special exact chain complexes are, let us analyze consequences of the condition $d_n\Comp d_{n+1} = 0$. For every $n\in \ZNr$, we find this decomposition diagram:
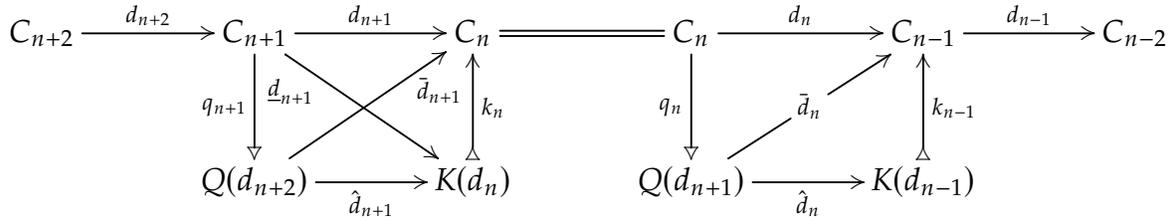
\begin{figure}[H]
  \stepcounter{theorem}
  \begin{equation*}
    \begin{vcenter}{\label{eq:ChainMapDecomp}
      \xymatrix@R=7ex@C=3.5em{
      C_{n+2} \ar[r]^-{d_{n+2}} &
      C_{n+1} \ar[r]^-{d_{n+1}} \ar@{-{ >>}}[d]_{q_{n+1}} \ar[rd]_(0.25){\ \ \underline{d}_{n+1} } &
      C_{n} \ar@{=}[r] &
      C_{n}\ar[r]^-{d_{n}} \ar@{-{ >>}}[d]_{q_{n}}&
      C_{n-1} \ar[r]^-{d_{n-1}} &
      C_{n-2}\\
      & \CoKer{d_{n+2}} \ar[r]_-{\hat{d}_{n+1}} \ar[ru]_(.75){\bar{d}_{n+1}} &
      \Ker{d_{n}} \ar@{{ |>}->}[u]_{k_n} &
      \CoKer{d_{n+1}} \ar[r]_-{\hat{d}_{n}} \ar|-{\ \bar{d}_{n}\ }[ru]&
      \Ker{d_{n-1}} \ar@{{ |>}->}[u]_{k_{n-1}}
      }
      }\end{vcenter}
  \end{equation*}
  \caption{Decomposition diagram of $d$}
\end{figure}
Using that $k_n$ is a kernel, we find a unique map $\underline{d}_{n+1}$ with $d_{n+1}= k_n\Comp \underline{d}_{n+1}$. Using that $q_{n+1}$ is a cokernel, we find a unique $\bar{d}_{n+1}$ with $d_{n+1}=\bar{d}_{n+1}\Comp q_{n+1}$.  Then:
\begin{enumerate}[(a)]
  \item $d_n\Comp \bar{d}_{n+1}=\ZeroMap$; so there is $r_{n+1}\from \CoKer{d_{n+2}}\to \Ker{d_n}$, unique with $\bar{d}_{n+1}= k_nr_{n+1}$.
  \item $\underline{d}_{n+1}\Comp d_{n+2}=\ZeroMap$; so there is $s_{n+1}\from \CoKer{d_{n+2}}\to \Ker{d_n}$, unique with $\underline{d}_{n+1}= s_{n+1}q_{n+1}$.
\end{enumerate}
We have $r_{n+1}=s_{n+1}$ because
\begin{equation*}
  k_ns_{n+1}q_{n+1} = k_{n}\underline{d}_{n+1}=d_{n+1} =\bar{d}_{n+1}q_{n+1}=k_{n}r_{n+1}q_{n+1}
\end{equation*}
Thus $\hat{d}_{n+1}\DefEq r_{n+1}=s_{n+1}$ is the unique map rendering the entire diagram commutative.

In general, there is no reason why $\underline{d}_{n+1}$ should be a normal epimorphism and, hence, there is no reason why $C$ should be exact in position $n$. Similarly, there is no reason why $\bar{d}_{n}$ should be a normal monomorphism and, hence, there is no reason why $C$ should be coexact in position $n$. Each of these observations leads to a measure of failure of exactness for chain complexes. The resulting discussion involves the following settings:

\begin{definition}[Types of chain complexes\ZExactTag]
  \label{def:ChainComplex-Types}%
  Of a chain complex $(C,d)$, we distinguish the following properties:
  \begin{thmlist}
    \item \label{def:ChainComplex-Normal}%
    $(C,d)$ is a \Defn{normal chain complex} if $d_{n}$ is a normal map for each $n\in\ZNr$. %
    \index{normal!chain complex}\index{chain complex!normal}%
    \index[not]{z!$\ZNr$\IndSep number system of integers}%
    \item \label{def:ChainComplex-Subnormal}%
    $(C,d)$ is a \Defn{subnormal chain complex} if $\underline{d}_{n}$ is a normal map for each $n\in\ZNr$. %
    \index{subnormal!chain complex}\index{chain complex!subnormal}%
    \item \label{def:ChainComplex-CoSubnormal}%
    $(C,d)$ is a \Defn{cosubnormal chain complex} if $\bar{d}_{n}$ is a normal map for each $n\in\ZNr$. %
    \index{cosubnormal!chain complex}\index{chain complex!cosubnormal}%
    \item \label{def:ChainComplex-WeaklyNormal}%
    $(C,d)$ is a \Defn{weakly normal chain complex} if $\hat{d}_{n}$ is a normal map for each $n\in\ZNr$. %
    \index{weakly normal chain complex}\index{chain complex!weakly normal}%
  \end{thmlist}
\end{definition}

We clarify the relationship between these types of chain complexes via the following three lemmas.

\begin{lemma}[Recognizing normality of $d_{n+1}$\ZExactTag]
  \label{thm:DifferentialNormal-Recognize}%
  Using the notation of the decomposition diagram \eqref{eq:ChainMapDecomp}, the following are equivalent.
  \begin{tfae}
    \item $d_{n+1}$ is normal.
    \item $\underline{d}_{n+1}=me$ is normal, and the composite $k_nm$ is a normal monomorphism.
    \item $\bar{d}_{n+1}=\nu \varepsilon$ is normal, and the composite $\varepsilon q_{n+1}$ is a normal epimorphism.\NoProof
  \end{tfae}
\end{lemma}

\begin{lemma}[Recognizing normality of $\hat{d}_{n+1}$\ZExactTag]
  \label{thm:DifferentialHatNormal-Recognize}%
  Using the notation of the decomposition diagram \eqref{eq:ChainMapDecomp}, the following are equivalent.
  \begin{tfae}
    \item $\hat{d}_{n+1}$ is normal.
    \item The initial normal mono factorization $\underline{d}_{n+1}=\mu u$ of $\underline{d}_{n+1}$ admits a cosubnormal refinement $\underline{d}_{n+1}=\mu \varepsilon q_{n+1}$ .
    \item The terminal normal epi factorization $\bar{d}_{n+1}=ve$ of $\bar{d}_{n+1}$ admits a subnormal refinement $\bar{d}_{n+1}=k_{n}me$. \NoProof
  \end{tfae}
\end{lemma}

\begin{lemma}[Cosubnormal and subnormal $\Rightarrow$ Normal\ZExactTag]
  \label{thm:DifferentialSubnormal+Cosubnormal->Normal}%
  A chain complex which is subnormal and cosubnormal is a normal chain complex.
\end{lemma}
\begin{proof}
  This follows from Proposition \ref{thm:Subnormal+Cosubnormal->Normal}.
\end{proof}

The following diagram summarizes the relationships between the types of chain complexes in (\ref{def:ChainComplex-Types}).
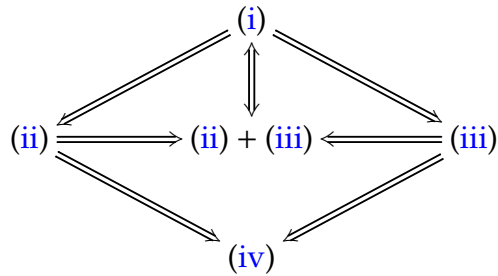
\begin{equation*}
  \xymatrix@!0@R=8ex@C=7em{
  & \text{(\ref{def:ChainComplex-Normal})} \ar@{=>}[ld] \ar@{=>}[rd] \ar@{<=>}[d] \\
  \text{(\ref{def:ChainComplex-Subnormal})} \ar@{=>}[rd] \ar@{=>}[r]&
  \text{(\ref{def:ChainComplex-Subnormal})} + \text{(\ref{def:ChainComplex-CoSubnormal})} & \text{(\ref{def:ChainComplex-CoSubnormal})} \ar@{=>}[ld] \ar@{=>}[l]\\
  &
  \text{(\ref{def:ChainComplex-WeaklyNormal})}
  }
\end{equation*}

For weakly normal chain complexes we introduce the following measures of `failure of exactness':

\begin{definition}[Homology of a chain complex\ZExactTag]%
  \label{def:Homology(ChainComplex)}%
  Let $(C,d)$ denote a weakly normal chain complex. Using the notation of the decomposition diagram \eqref{eq:ChainMapDecomp}, the \Defn{cokernel homology} of  $(C,d)$  in position  $n\in \ZNr$  is %
  \index[not]{h!$\HmlgyCoKer{n}{C}$\IndSep cokernel construction of homology}%
  \begin{equation*}
    \HmlgyCoKer{n}{C} \DefEq \CoKer{\hat{d}_{n+1}} \cong \CoKer{ \underline{d}_{n+1}\from C_{n+1}\to \Ker{d_n} }
  \end{equation*}
  The \Defn{kernel homology} of $(C,d)$ in position $n$ is %
  \index[not]{h!$\HmlgyKer{n}{C}$\IndSep kernel construction of homology}
  \begin{equation*}
    \HmlgyKer{n}{C} \DefEq \Ker{\hat{d}_{n}} \cong \KerMap{\bar{d}_{n}}
  \end{equation*}
\end{definition}

The two constructions of homology are related as follows.

\begin{lemma}[Relationship: $\HmlgyCoKer{n}{C}$ and $\HmlgyKer{n}{C}$\DExTag]%
  \label{thm:Homology-CoKerHomology->KerHomology}%
  For a weakly normal chain complex $(C,d)$, the objects $\HmlgyCoKer{n}{C}$ and $\HmlgyKer{n}{C}$ are related via this functorial exact and coexact sequence: %
  \begin{equation*}
    \xymatrix@R=5ex@C=3em{
    \CoKer{d_{n+2}} \ar[r]^-{\bar{d}'_{n+1}} &
    \Ker{q_n} \ar[r] &
    \HmlgyCoKer{n}{C} \ar[r]^-{\lambda_n(C)} &
    \HmlgyKer{n}{C}  \ar[r] &
    \CoKer{k_{n}} \ar[r]^-{\underline{d}'_{n}} &
    \Ker{d_{n-1}}
    }
  \end{equation*}
  where $\bar{d}_{n+1} = \KerMap{q_{n}}\Comp \bar{d}'_{n+1}$ and $\underline{d}_{n} = \underline{d}'_{n}\Comp \CoKerMap{k_{n}}$. The comparison map $\lambda_n(C)$ has the following properties:
  \begin{thmlist}
    \item $\lambda_{n}(C)$ is a normal monomorphism if and only if $d_{n+1}$ is cosubnormal.
    \item $\lambda_{n}(C)$ is a normal epimorphism if and only if $d_{n}$ is subnormal.
  \end{thmlist}
\end{lemma}
\begin{proof}
  We relate the kernel and cokernel variants of homology via this diagram:
  \begin{equation*}
    \xymatrix@R=4.5ex@C=4em{
    C_{n+2} \ar[r]^-{d_{n+2}} &
    C_{n+1} \ar@{-{ >>}}[d]_{q_{n+1}} \ar[r]^-{d_{n+1}} &
    C_{n} \ar@{=}[d] &
    \HmlgyKer{n}{C} \ar@{{ |>}->}[d] \\
    & \CoKer{d_{n+2}} \ar[r]^-{\bar{d}_{n+1}} \ar[d]_{\hat{d}_{n+1}} &
    C_n \ar@{-{ >>}}[r]^-{q_n} \ar@{=}[d] &
    \CoKer{d_{n+1}} \ar[d]^{\hat{d}_{n}} \\
    & \Ker{d_{n}} \ar@{-{ >>}}[d] \ar@{{ |>}->}[r]_-{k_n} &
    C_{n} \ar[r]_-{\underline{d}_{n}} \ar@{=}[d]  &
    \Ker{d_{n-1}} \ar@{{ |>}->}[d]^{k_{n-1}} \\
    & \HmlgyCoKer{n}{C} &
    C_{n} \ar[r]_-{d_{n}} &
    C_{n-1} \ar[r]_-{d_{n-1}} &
    C_{n-2}
    }
  \end{equation*}
  By assumption, the maps $\hat{d}_{k}$ are normal. Thus the Generalized Pure Snake Lemma \ref{thm:PureSnake-Generalized-II} yields the claimed exact sequence and its functoriality.

  (i) By exactness, the map $\lambda_{n}(C)$ is a normal monomorphism if and only if $\bar{d}'_{n+1}$ is a normal epimorphism. This happens if and only if $\bar{d}_{n+1}$ is a normal map which, in turn, is equivalent to $d_{n+1}$ being cosubnormal.

  (ii) Similarly, $\lambda_{n}(C)$ is a normal epimorphism if and only if $\underline{d}'_{n}$ is a normal monomorphism. This happens if and only if $\underline{d}_{n}$ is a normal map which, in turn, is equivalent to $d_{n}$ being subnormal.
\end{proof}

\begin{corollary}[Homology of a types of chain complexes\DExTag]
  \label{thm:Homology-TypesChainComplexes}
  For a weakly normal chain complex $(C,d)$ the following hold
  \begin{thmlist}
    \item If $(C,d)$ is subnormal, then $\lambda_{n}(C)\from \HmlgyCoKer{n}{C}\to \HmlgyKer{n}{C}$ is a normal epimorphism for all $n\in \ZNr$.
    \item If $(C,d)$ is cosubnormal, then $\lambda_{n}(C)\from \HmlgyCoKer{n}{C}\to \HmlgyKer{n}{C}$ is a normal monomorphism for all $n\in \ZNr$.
    \item If $(C,d)$ is normal, then $\lambda_{n}(C)\from \HmlgyCoKer{n}{C}\to \HmlgyKer{n}{C}$ is an isomorphism for all $n\in \ZNr$.
  \end{thmlist}
\end{corollary}
\begin{proof}
  (i) and (ii) follow from (\ref{thm:Homology-CoKerHomology->KerHomology}). This implies (iii) because a map that is both subnormal and cosubnormal is normal (\ref{thm:Subnormal+Cosubnormal->Normal}).
\end{proof}

\begin{corollary}[Homology of a normal chain complex\DExTag]
  \label{thm:Homology-NormalChainComplex}%
  For a normal chain complex $(C,d)$ in a di-exact category, $\lambda_{n}(C)\from \HmlgyCoKer{n}{C} \cong \HmlgyKer{n}{C}$ is an isomorphism for every $\in \ZNr$. \NoProof
\end{corollary}

In view of Corollary \ref{thm:Homology-NormalChainComplex}, there is no need to distinguish between the cokernel vs. the kernel based approach to homology. Instead, as in abelian categories, we use the isomorphisms $\lambda_{n}(C)$ to define the homology of the complex:

\begin{notation}[Homology of a normal chain complex\DExTag]
  Given a normal chain complex $(C,d)$ in a di-exact category, we define
  \begin{equation*}
    \Hmlgy{n}{C} \DefEq \HmlgyCoKer{n}{C} \cong \HmlgyKer{n}{C}
  \end{equation*}
\end{notation}

We confirm that either of the two homology variants we just introduced vanishes if and only if the chain complex in question is exact in the appropriate position.

\begin{proposition}[$\HmlgyCoKer{n}{C}$ vanishes if and only if chain complex exact\ZExactTag]
  \label{thm:H^c(C)=0<->Exact}
  For a weakly normal chain complex $(C,d)$ in a pointed category $\EuScript{X}$ the following hold.
  \begin{thmlist}
    \item If $(C,d)$ is exact in position $n\in \ZNr$, then $\HmlgyCoKer{n}{C}=\ZeroObject$.
    \item If $\underline{d}_{n}$ is normal, and $\HmlgyCoKer{n}{C}=\ZeroObject$, then $\HmlgyKer{n}{C}=\ZeroObject$ as well, and $C$ is exact in position $n$.
  \end{thmlist}
  Hence a subnormal chain complex $(C,d)$ is exact in position $n\in \ZNr$ if and only if $\HmlgyCoKer{n}{C}=\HmlgyKer{n}{C}=\ZeroObject$.
\end{proposition}
\begin{proof}
  We reason via the commutative diagram below.
  \begin{equation*}
    \xymatrix@R=5ex@C=4em{
    C_{n+1} \ar@{-{ >>}}[r]_(0.5){q_{n+1}} \ar@{-{ >>}}[d]_{\rho_{n+1}} \ar[rd]|-{\underline{d}_{n+1}} \ar@/^3ex/[rr]^-{d_{n+1}} &
    \CoKer{d_{n+2}} \ar[r]_-{\bar{d}_{n+1}} \ar[d]^{\hat{d}_{n+1}} &
    C_n \ar@{-{ >>}}[r]^-{q_n} \ar@{=}[d] &
    \CoKer{d_{n+1}} \ar[d]^{\bar{d}_{n}} \\
    I \ar@{{ |>}->}[r]_-{e_{n+1}} &
    \Ker{d_n} \ar@{{ |>}->}[r]_-{k_n} &
    C_n \ar[r]_-{d_n} &
    C_{n-1}
    }
  \end{equation*}
  (i)\quad If $C$ is exact in position $n$, then $\underline{d}_{n+1}$ is a normal epimorphism. So $\HmlgyCoKer{n}{C}=\CoKer{\underline{d}_{n+1}}=0$.

  (ii)\quad If $\HmlgyCoKer{n}{C}=\CoKer{\hat{d}_{n+1}}=0$, and $\underline{d}_{n+1}=e_{n+1}\rho_{n+1}$ is normal, then $e_{n+1}=\KerMap{\ZeroMap}$ may be chosen to be $\IdMapOn{\Ker{d_n}}$. But then $\underline{d}_{n+1}=\rho_{n+1}$, implying that $C$ is exact in position $n$.
\end{proof}

Dually:

\begin{proposition}[$\HmlgyKer{n}{C}$ vanishes if and only if chain complex coexact\ZExactTag]
  \label{thm:H^k(C)=0<=>Exact}%
  For a weakly normal chain complex $(C,d)$ in a pointed category $\EuScript{X}$ the following hold.
  \begin{thmlist}
    \item If $(C,d)$ is coexact in position $n\in \ZNr$, then $\HmlgyKer{n}{C}=\ZeroObject$.
    \item If $\bar{d}_{n}$ is normal, and $\HmlgyKer{n}{C}=\ZeroObject$, then $C$ is coexact in position $n$.
  \end{thmlist}
  A cosubnormal chain complex $(C,d)$ is coexact in position $n\in \ZNr$ if and only if $\HmlgyKer{n}{C}=\ZeroObject$.\NoProof
\end{proposition}

Hence:

\begin{corollary}[Simultaneous vanishing of $\HmlgyKer{n}{C}$ and $\HmlgyCoKer{n}{C}$\ZExactTag]
  \label{thm:Exact<=>H_n=0}%
  For a normal chain complex $C$ and $n\in \ZNr$, we have $\HmlgyKer{n}{C}=\ZeroObject=\HmlgyCoKer{n}{C}$ if and only if $C$ is exact in position $n$.\NoProof
\end{corollary}

Lemma \ref{thm:Homology-CoKerHomology->KerHomology} leaves open the possibility that, in a subnormal chain complex, $\HmlgyKer{n}{C}$ can vanish while $\HmlgyCoKer{n}{C}\neq 0$. In this case, the property `$\HmlgyKer{n}{C}=0$' does not characterize exactness of $C$ in position $n$. Examples show that this actually happens:

\begin{example}[Non-subnormal chain complex]
  \label{exa:Non-Subnormal-ChainComplex-Different-H-and-K}%
  In the category $\Grps$ of groups, consider the following chain complex $C$, in which $C_0=\langle y\rangle$ is the free group on the single generator $y$:
  \begin{equation*}
    \xymatrix{
    {\cdots} \ar[r] &
    0 \ar[r] &
    {\langle x \rangle} \ar[r]^-{d_{2}} &
    {\langle x,y \rangle} \ar@{-{ >>}}[r]^-{d_{1}} &
    {\langle y \rangle} \ar[r] &
    0 \ar[r] &
    {\cdots.}
    }
  \end{equation*}
  Further, $d_{2}$ is the inclusion which maps $x$ to $x$, and $d_{1}$ is the quotient map which sends $x$ to $1$ and $y$ to $y$. Then we find that $\HmlgyKer{1}{C}=0$, while
  \begin{equation*}
    \HmlgyCoKer{1}{C}\cong\langle y^{n}xy^{-n},\,n\in \ZNr - \Set{0}\rangle
  \end{equation*}
  The complex $C$ is not subnormal, though.
\end{example}

\begin{example}[Subnormal, non-normal chain complex]
  \label{exa:SubNormal-NonNormal-ChainComplex-Different-H-and-K}
  In the category $\Grps$ of groups we now construct a subnormal chain complex $C$ which is not exact but has $\HmlgyKer{n}{C}=0$: Choose a group $M$ with normal subgroup $L$ containing a normal subgroup $K$ whose normal closure in $M$ is $L$. Then this diagram is a subnormal chain complex.
  \[
    \xymatrix{{\cdots} \ar[r] & 0 \ar[r] & K \ar[r]^-{d_{2}} & M \ar@{-{ >>}}[r]^-{d_{1}} & M/L \ar[r] & 0 \ar[r] & {\cdots,}}
  \]
  If $M$ sits in position $0$, then $\HmlgyKer{0}{C}=0$ while $\HmlgyCoKer{n}{C}\cong L/K\neq 0$. In particular, the chain complex is not exact at $M$.

  For a concrete example, choose $M/L=C_2=\Set{\pm 1}$ a $2$-element group acting on $L\DefEq \ZNr \prdct \ZNr$ by interchanging coordinates. Then take $K\DefEq \ZNr\prdct 0$.
\end{example}

\begin{definition}[Morphism of chain complexes\ZExactTag]%
  \label{def:Morphism-OfChainComplexes}%
  A \Defn{morphism $f\from C\to D$ of chain complexes} is given by a commutative diagram of the form: %
  \index{morphism!of chain complexes}\index{chain complex!morphism}
  \begin{equation*}
    \xymatrix@R=5ex@C=2em{
    C \ar[d]_{f} && \cdots \ar[r] &
    C_{n+1} \ar[rr]^-{d^{C}_{n+1}} \ar[d]^{f_{n+1}} &&
    C_n \ar[rr]^-{d^{C}_{n}} \ar[d]^{f_n} &&
    C_{n-1} \ar[d]^{f_{n-1}} \ar[r] &
    \cdots \\
    D && \cdots \ar[r] &
    D_{n+1} \ar[rr]_-{d^{D}_{n+1}} &&
    D_n \ar[rr]_-{d^{D}_{n}} &&
    D_{n-1} \ar[r]_-{d^{D}_{n-1}} &
    \cdots
    }
  \end{equation*}
\end{definition}

The chain complexes in a pointed category $\EuScript{C}$ form a category containing as full subcategories `subnormal chain complexes', `cosubnormal chain complexes' and `normal chain complexes'.

\begin{proposition}[$\HmlgyCoKer{n}{-}$ and $\HmlgyKer{n}{-}$ are functorial\ZExactTag]
  \label{thm:HomologyFunctorial}%
  A morphism $f\from A\to B$ of weakly normal chain complexes induces, for each $n\in \ZNr$, a commutative square of homology objects: %
  \index{homology!functorial}%
  \begin{equation*}
    \xymatrix@R=5ex@C=4em{
    \HmlgyCoKer{n}{A} \ar[r]^-{\lambda_n(A)} \ar[d]_{\HmlgyCoKer{n}{f}} &
    \HmlgyKer{n}{A} \ar[d]^{\HmlgyKer{n}{f}} \\
    \HmlgyCoKer{n}{B} \ar[r]_-{\lambda_n(A)} &
    \HmlgyKer{n}{B}
    }
  \end{equation*}
  The above diagram is functorial with respect to morphisms of weakly normal chain complexes.
\end{proposition}
\begin{proof}
  The diagram we used to define $\HmlgyCoKer{n}{}$ and $\HmlgyKer{n}{}$, see (\ref{def:Homology(ChainComplex)}) is functorial with respect to morphisms of chain complexes. So, the claim follows.
\end{proof}

\begin{subordinate}{}

  \begin{subsubordinate}{On Minimal Conditions}
    Both, lemma (\ref{thm:Homology-CoKerHomology->KerHomology}), which explains the relationship between the kernel and the cokernel constructions of homology, and its corollary (\ref{thm:Homology-TypesChainComplexes}) rely only on the Pure Snake Lemma. So, both are valid in homologically self-dual categories.
  \end{subsubordinate}

  \begin{subsubordinate}{On the Need for Both Types of Homology}
    In view of Corollary \ref{thm:Homology-NormalChainComplex}, one might be tempted to conclude that considering both, the kernel definition as well as the cokernel definition of homology is redundant if we agree to work with normal chain complexes only. This, however, is not so because constructions that we might want to perform force us to consider both kinds of homology in conjunction; see in particular Theorem \ref{thm:(Co)Ker(ProperMapLESs)}.
  \end{subsubordinate}

\end{subordinate}

\begin{exercises}

\begin{exercise}[(Non-)cancellation in a product\ZExactTag]
  \label{exe:(Non)CancellationInProduct}
  Give an example of a pointed category with non-zero objects $X$ and $Y$ such $\Prdct{X}{Y}\cong Y$. Then show the following: In every pointed category, for objects $X$ and $Y$ which admit a product $\Prdct{X}{Y}$, if $\PrjctnOnto{Y}\from \Prdct{X}{Y}\to Y$ is a monomorphism, then $X\cong \ZeroObject$.
\end{exercise}
\end{exercises}
\newpage
\section[The Long Exact Homology Sequence]{The Long Exact Homology Sequence}
\label{sec:LES-Homology}%

We provide effective tools for computing with homology. Most importantly, this includes the long exact sequence of homology objects associated with a short exact sequence of normal chain complexes; see (\ref{thm:LES-Homology}).

Next, we analyze the chain complexes obtained by taking  pointwise kernels and cokernels of a normal morphism of long exact sequences, Theorem \ref{thm:(Co)Ker(ProperMapLESs)}. We find that the chain complex of pointwise kernels is subnormal, while the chain complex of cokernels is cosubnormal. Then we show that the homology of the cokernel complex agrees with the homology of the kernel complex up to a dimension shift (\ref{thm:(Co)Ker(ProperMapLESs)}). In abelian categories this is well known and takes a much simpler form. Under the present minimal assumptions, it confirms that it is essential to introduce the fine grained cokernel and kernel homologies (\ref{def:Homology(ChainComplex)}), we well as various types of chain complexes (\ref{def:ChainComplex}) .

Theorem \ref{thm:(Co)Ker(ProperMapLESs)} provides homological methods for a whole range of results familiar from abelian categories. For example:
\begin{ulist}
  \item We find conditions under which a normal morphism of exact chain complexes has an exact kernel complex, respectively an exact cokernel complex; see (\ref{thm:LESMor-KerExact<->CoKerExact}).
  \item For a sequence $X\XRA{f} Y \XRA{g} Z$ of normal maps, we find its associated $6$-term exact sequence of kernels and cokernels; see (\ref{thm:6TermES-From-Composite}).
  \item We obtain further variations of the Snake Lemma in (\ref{exe:SnakeLemma-Relaxed-Strong}) and in (\ref{exe:SnakeLemma-Extended}).
\end{ulist}

\begin{definition}[Exactness of composite of morphisms of chain complexes\ZExactTag]
  \label{def:Exactness-MapsOfChainComplexes}%
  A diagram $A\to B\to C$ of chain complexes is called \Defn{short exact} if, for each $n\in \ZNr$, the sequence $A_n\to B_n\to C_n$ is short exact. %
  \index{short exact sequence!of chain complexes}%
\end{definition}

\begin{theorem}[Long exact homology sequence\DExTag]%
  \label{thm:LES-Homology}%
  A short exact sequence $A \overset{f}{\NMono} B \overset{g}{\NEpi} C$ of  normal chain complexes functorially induces a long exact sequence of homology objects:
  \begin{equation*}
    \cdots \to \Hmlgy{n+1}{C} \XRA{\partial_{n+1}} \Hmlgy{n}{A} \XRA{f_{\ast}} \Hmlgy{n}{B} \XRA{g_{\ast}} \Hmlgy{n}{C} \XRA{\partial_n} \Hmlgy{n-1}{A} \longrightarrow \cdots
  \end{equation*}
\end{theorem}
\begin{proof}
  \emph{Step 1}\quad Each position $k$ of the short exact sequence of chain complexes yields a morphism of short exact sequences in $\SACtgry{X}$:
  \begin{equation*}
    \xymatrix@R=5ex@C=3em{
    \Ker{d^{A}_{k+1}} \ar@{{ |>}->}[d] \ar@{{ |>}->}[r] &
    \Ker{d^{B}_{k+1}} \ar@{{ |>}->}[d] \ar[r] &
    \Ker{d^{C}_{k+1}} \ar@{{ |>}->}[d] \\
    A_{k+1} \ar[d]_{d^{A}_{k+1}} \ar@{{ |>}->}[r] &
    B_{k+1} \ar[d]_{d^{B}_{k+1}} \ar@{-{ >>}}[r] &
    C_{k+1} \ar[d]^{d^{C}_{k+1}} \\
    A_k \ar@{-{ >>}}[d] \ar@{{ |>}->}[r] &
    B_k \ar@{-{ >>}}[r] \ar@{-{ >>}}[d] &
    C_k \ar@{-{ >>}}[d] \\
    \CoKer{d^{A}_{k+1}} \ar[r] &
    \CoKer{d^{B}_{k+1}} \ar@{-{ >>}}[r] &
    \CoKer{d^{C}_{k+1}}
    }
  \end{equation*}
  The top sequence is coexact, and the bottom sequence is exact by the Snake Lemma  \ref{thm:SnakeLemma-Classical}.

  \emph{Step 2}\quad Choosing  $k=n+1$  and  $n-1$, we form the following commutative diagram.
  \begin{equation*}
    \xymatrix@R=5ex@C=3em{
    \HmlgyKer{n}{A} \ar[r] \ar@{{ |>}->}[d] &
    \HmlgyKer{n}{B} \ar[r] \ar@{{ |>}->}[d] &
    \HmlgyKer{n}{C} \ar@{{ |>}->}[d] \\
    \CoKer{d^{A}_{n+1}} \ar[r] \ar[d]_{\hat{d}^{A}_{n+1}} &
    \CoKer{d^{B}_{n+1}} \ar@{-{ >>}}[r] \ar[d]^{\hat{d}^{B}_{n+1}} &
    \CoKer{d^{C}_{n+1}} \ar[d]_{\hat{d}^{C}_{n+1}} \\
    \Ker{d^{A}_{n-1}} \ar@{{ |>}->}[r] \ar@{-{ >>}}[d] &
    \Ker{d^{B}_{n-1}} \ar[r] \ar@{-{ >>}}[d] &
    \Ker{d^{C}_{n-1}} \ar@{-{ >>}}[d] \\
    \HmlgyCoKer{n-1}{A} \ar[r] &
    \HmlgyCoKer{n-1}{B} \ar[r] &
    \HmlgyCoKer{n-1}{C}
    }
  \end{equation*}
  The maps $\hat{d}^{A}_{n},\hat{d}^{B}_{n},\hat{d}^{C}_{n}$ are normal. So we apply the Relaxed Snake Lemma \ref{thm:SnakeLemma-Relaxed} to obtain a segment of the claimed exact sequence of homology objects. To complete the proof, all that remains is to identify kernel/cokernel based homology objects via (\ref{thm:Homology-CoKerHomology->KerHomology}). Since the Relaxed Snake Lemma as well as the isomorphism of cokernel and kernel defined homology is functorial, so is the long exact homology sequence.
\end{proof}

\begin{terminology}[Long exact homology sequence]
  \label{def:LongExactHomologySequence}%
  In the setting of (\ref{thm:LES-Homology}), the sequence
  $$
    \cdots \to \Hmlgy{n+1}{C} \XRA{\partial_{n+1}} \Hmlgy{n}{A} \XRA{f_{\ast}} \Hmlgy{n}{B} \XRA{g_{\ast}} \Hmlgy{n}{C} \XRA{\partial_n} \Hmlgy{n-1}{A} \longrightarrow \cdots
  $$
  is called the \Defn{long exact homology sequence} of the given short exact sequence of normal chain complexes. The maps $\partial_{n}$ are known as connecting homomorphisms%
  \index{long exact homology sequence}
\end{terminology}

\begin{theorem}[(Co-)kernels of a normal morphism of long exact sequences\DExTag]
  \label{thm:(Co)Ker(ProperMapLESs)}%
  Consider a morphism of long exact sequences:
  \begin{equation*}
    \xymatrix@R=5ex@C=3em{
    \cdots \ar[r] &
    A_{n+2} \ar[r]^-{d^{A}_{n+2}} \ar[d]_{\alpha_{n+2}} &
    A_{n+1} \ar[r]^-{d^{A}_{n+1}} \ar[d]_{\alpha_{n+1}} &
    A_{n} \ar[r]^-{d^{A}_{n}} \ar[d]_{\alpha_{n}} &
    A_{n-1} \ar[r]^-{d^{A}_{n-1}} \ar[d]_{\alpha_{n-1}} &
    A_{n-2} \ar[r] \ar[d]_{\alpha_{n-2}} &
    \cdots \\
    \cdots \ar[r] &
    B_{n+2} \ar[r]_-{d^{B}_{n+2}} &
    B_{n+1} \ar[r]_-{d^{B}_{n+1}} &
    B_{n} \ar[r]_-{d^{B}_{n}} &
    B_{n-1} \ar[r]_-{d^{B}_{n-1}} &
    B_{n-2} \ar[r] &
    \cdots
    }
  \end{equation*}
  If each $\alpha_{n}$ is a normal map, then the kernel objects $K_n\DefEq \Ker{\alpha_n}$ form a subnormal chain complex, while the cokernel terms $C_n\DefEq \CoKer{\alpha_n}$ form a cosubnormal chain complex. The homologies of these two chain complexes are linked via an isomorphism which is functorial with respect to morphisms of normal morphisms of long exact sequences:
  \begin{equation*}
    \HmlgyCoKer{n-1}{K} \cong \HmlgyKer{n+1}{C}
  \end{equation*}
\end{theorem}
\begin{proof}
  \emph{Step 1: Long exact sequences are spliced short exact ones}\quad By Definition \ref{def:LES}, each of the long exact sequences $(A,d^A)$ and $(B,d^B)$ is spliced together from short exact sequences as shown. Since normal factorizations are functorial (\ref{thm:NormalFactorization-Functorial}), we obtain this commutative diagram.
  \begin{equation*}
    \xymatrix@!0@R=6ex@C=5em{
    A_{n+2} \ar[r]^-{d^{A}_{n+2}} \ar[dd]_{\alpha_{n+2}} &
    A_{n+1} \ar[rr]^-{d^{A}_{n+1}} \ar[dd]_{\alpha_{n+1}} \ar@{-{ >>}}[dr]_{\varepsilon^{A}_{n+1}} &&
    A_{n} \ar[rr]^-{d^{A}_{n}} \ar[dd]^{\alpha_{n}} \ar@{-{ >>}}[dr] &&
    A_{n-1} \ar[r]^-{d^{A}_{n-1}} \ar[dd]_{\alpha_{n-1}} &
    A_{n-2} \ar[dd]_{\alpha_{n-2}} \\
    && I_{n+1} \ar@{{ |>}->}[ru]_{\mu^{A}_{n+1}} \ar[dd]|\hole_(.3){\gamma_{n+1}} &&
    I_n \ar@{{ |>}->}[ru] \ar[dd]|\hole_(0.3){\gamma_{n}} \\
    B_{n+2} \ar[r]_-{d^{B}_{n+2}} &
    B_{n+1} \ar[rr]^(.75){d^{B}_{n+1}} \ar@{-{ >>}}[rd]_{\varepsilon^{B}_{n+1}} &&
    B_{n} \ar[rr]^(.7){d^{B}_{n}} \ar@{-{ >>}}[rd] &&
    B_{n-1} \ar[r]_-{d^{B}_{n-1}} &
    B_{n-2} \\
    && J_{n+1} \ar@{{ |>}->}[ru]_{\mu^{B}_{n+1}} &&
    J_{n} \ar@{{ |>}->}[ru]
    }
  \end{equation*}
  Here  $I_{n+1}\DefEq \CoKer{d^{A}_{n+2}}\cong \Ker{d^{A}_{n}}$ and $J_{n+1}\DefEq \CoKer{d^{B}_{n+2}}\cong \Ker{d^{B}_{n}}$. The maps $\gamma_{n}$ are normal by  (\ref{thm:ANN<->NMNIC}). Since kernels and cokernels are functorial,  Step 1 depends functorially on morphisms of normal morphisms of long exact sequences.

  To continue the argument, we insert the $(\Prdct{3}{3})$-diagram of normal factorizations from (\ref{def:NICNormalMaps}) into the squares $A_{n}\rightrightarrows B_{n-1}$  of normal maps.
  \begin{equation*}
    \resizebox{.95\textwidth}{!}{$
      \xymatrix@R=5ex@C=4em{
      && K_{n} \ar[r]_-{\underline{k}_{n}} \ar@{{ |>}->}[d]_{\kappa_{n}} \ar@/^2ex/[rr]^-{k_{n}}&
      L_{n} \ar@{{ |>}->}[r]_-{\lambda_{n}} \ar@{{ |>}->}[d]_{\sigma_{n}} \PullLU{rd} &
      K_{n-1} \ar@/^2ex/[rr]^-{k_{n-1}} \ar@{{ |>}->}[d]_{\kappa_{n-1}} \ar[r]_-{\underline{k}_{n-1}} &
      L_{n-1} \ar@{{ |>}->}[r]_-{\lambda_{n-1}} \ar@{{ |>}->}[d]_{\sigma_{n-1}} &
      K_{n-2} \ar@{{ |>}->}[d]_{\kappa_{n-2}}\\
      && A_{n} \ar@{-{ >>}}[r]_-{\varepsilon^{A}_{n}} \ar@{-{ >>}}[d]^{e_{n}} \ar@/_2ex/[dd]_{\alpha_{n}} &
      I_{n} \ar@{{ |>}->}[r]_-{\mu^{A}_{n}} \ar@{-{ >>}}[d] &
      A_{n-1} \ar@{-{ >>}}[r]_-{\varepsilon^{A}_{n-1}} \ar@{-{ >>}}[d]_{e_{n-1}}\ar@/^2ex/[dd]^-{\alpha_{n-1}} &
      I_{n-1} \ar@{{ |>}->}[r]_-{\mu^{A}_{n-1}} &
      A_{n-2} \\
      && \DiagObj \ar@{-{ >>}}[r] \ar@{{ |>}->}[d]^{\mu_{n}} &
      \DiagObj \ar@{{ |>}->}[r] \ar@{{ |>}->}[d]&
      \DiagObj \ar@{{ |>}->}[d]_{\mu_{n-1}} \\
      B_{n+1} \ar@{-{ >>}}[r]^-{\varepsilon^{B}_{n+1}} \ar@{-{ >>}}[d]_{\varepsilon_{n+1}} &
      J_{n+1}\ar@{{ |>}->}[r]^-{\mu^{B}_{n+1}} \ar[d]^{\tau_{n+1}} &
      B_{n} \ar@{-{ >>}}[r]^-{\varepsilon^{B}_{n}} \ar@{-{ >>}}[d]_{\varepsilon_{n}} \PushRD{rd} &
      J_{n} \ar@{{ |>}->}[r]^{\mu^{B}_{n}} \ar@{-{ >>}}[d]^{\tau_{n}} &
      B_{n-1} \ar@{-{ >>}}[d]^{\varepsilon_{n-1}}\\
      C_{n+1} \ar@{-{ >>}}[r]^-{\rho_{n+1}} \ar@/_2ex/[rr]_-{c_{n+1}} &
      T_{n+1} \ar[r]^-{\bar{c}_{n+1}} &
      C_{n} \ar@{-{ >>}}[r]^-{\rho_{n}} \ar@/_2ex/[rr]_-{c_{n}} &
      T_{n} \ar[r]^-{\bar{c}_{n}} &
      C_{n-1}
      }
    $}
  \end{equation*}

  \emph{Step 2: Computation of the homology of the subnormal kernel sequence}\quad By definition, $\sigma_{n}\DefEq \KerMap{\gamma_{n}}$. This yields the factorization $k_n=\lambda_n \bar{k}_n$. The square $L_{n}\rightrightarrows A_{n-1}$ is a pullback (\ref{thm:PullbackRecognition-KernelSide-1}). So $\lambda_{n}$ is a normal monomorphism by (\ref{thm:PullbackPreservesNormalMonos}). It is the kernel of the composite $\varepsilon^{A}_{n-1}\kappa_{n-1}$ by (\ref{thm:NormalMono-Props})  and, hence, the kernel of $k_{n-1}$.

  Further, the antinormal composite $\varepsilon^{A}_{n}\kappa_{n}$ also has a normal factorization $\sigma_{n} \underline{k}_{n} = me
  $, with $e$ a normal epimorphism, and $m$ a normal monomorphism. So $\underline{k}_{n}$ is a normal map by (\ref{thm:FactorNormalThroughSubObject->Normal}). Since $k_{n}$ factors through $\KerMap{k_{n-1}}$, $(K,k)$ is a chain complex. This chain complex is subnormal by what we just showed. It depends functorially on morphisms of normal morphisms between long exact sequences. - By definition:
  \begin{equation*}
    \HmlgyCoKer{n-1}{K}=\CoKer{\underline{k}_{n}}
  \end{equation*}
  \emph{Step 3: Computation of the homology of the cosubnormal cokernel complex}\quad By definition, $\tau_{n}\DefEq \CoKerMap{\gamma_{n}}$. This yields the factorization $c_{n}=\bar{c}_{n}\rho_{n}$. The square $B_{n}\rightrightarrows T_{n}$ is a pushout by (\ref{thm:PushoutRecognize-Categorical}). So $\rho_{n}$ is a normal epimorphism by (\ref{thm:PushoutPreservesNormalEpis}). It is the cokernel of the composite $\varepsilon_{n}\mu^{B}_{n+1}$ and, hence, the cokernel of $c_{n+1}$.

  Further, the antinormal composite $\varepsilon_{n}\mu^{B}_{n+1}$ has a normal factorization via which we see that $\bar{c}_{n+1}$ is a normal map. This renders $c_{n+1}$ a cosubnormal map. Since $c_{n}$ factors through $\CoKerMap{c_{n+1}}$, $(C,c)$ is a chain complex which is cosubnormal by what we just showed. It depends functorially on morphisms of normal morphisms between long exact sequences. - By definition:
  \begin{equation*}
    \HmlgyKer{n+1}{C} = \Ker{\bar{c}_{n+1}}
  \end{equation*}
  \emph{Step 4: $\HmlgyCoKer{n-1}{K}\cong \HmlgyKer{n+1}{C}$}\quad The Snake Lemma \ref{thm:SnakeLemma-Classical} yields:
  \begin{equation*}
    \xymatrix@R=5ex@C=4em{
    & L_{n+1} \ar@{{ |>}->}[r]^-{\lambda_{n+1}} \ar@{{ |>}->}[d]_{\sigma_{n+1}} &
    K_{n} \ar[r]^-{\underline{k}_{n}} \ar@{{ |>}->}[d]_{\kappa_{n}} &
    L_{n} \ar@{-{ >>}}[r] \ar@{{ |>}->}[d]^{\sigma_{n}} &
    \HmlgyCoKer{n-1}{K} \\
    & I_{n+1} \ar@{{ |>}->}[r] \ar[d]_{\gamma_{n+1}} &
    A_{n} \ar@{-{ >>}}[r] \ar[d]_{\alpha_{n}} &
    I_{n} \ar[d]_{\gamma_{n}} \\
    & J_{n+1} \ar@{{ |>}->}[r] \ar@{-{ >>}}[d]_{\tau_{n+1}} &
    B_{n} \ar@{-{ >>}}[r] \ar@{-{ >>}}[d]_{\varepsilon_{n}} &
    J_{n} \ar@{-{ >>}}[d]_{\tau_{n}} \\
    \HmlgyKer{n+1}{C} \ar@{{ |>}->}[r] &
    T_{n+1} \ar[r]_-{\bar{c}_{n+1}} &
    C_{n} \ar@{-{ >>}}[r]_{\rho_{n}} &
    T_{n}
    }
  \end{equation*}
  From the proof of the Snake Lemma we obtain the desired isomorphism
  \begin{equation*}
    \HmlgyCoKer{n-1}{K} \cong \HmlgyKer{n+1}{C}
  \end{equation*}
  This completes the proof.
\end{proof}

\begin{corollary}[Kernel sequence exact $\Leftrightarrow$ cokernel sequence exact\DExTag]
  \label{thm:LESMor-KerExact<->CoKerExact}%
  A normal morphism of long exact sequences has an exact sequence of pointwise kernels if and only if it has an exact sequence of pointwise cokernels.
\end{corollary}
\begin{proof}
  With (\ref{thm:H^c(C)=0<->Exact}) we see that $(K,k)$ is exact if and only if $\HmlgyCoKer{n}{K}=\ZeroObject$, for all $n\in \ZNr$. By (\ref{thm:(Co)Ker(ProperMapLESs)}), this happens if and only if $\HmlgyKer{n}{C}=\ZeroObject$, for all $n\in \ZNr$. With (\ref{thm:H^k(C)=0<=>Exact}), we see that this is equivalent to the exactness of $(C,c)$.
\end{proof}

\begin{corollary}[Positionwise (co-)kernel is (co-)kernel\DExTag]
  \label{thm:LES-Positionwise(Co-)KerIs(Co-)Ker}%
  For a normal morphism  $f\from A\to B$  of long exact sequences the following hold:
  \begin{thmlist}
    \item If each $f_n$ is a monomorphism, then $f$ is a normal monomorphism in the category of long exact sequences.
    \item If each $f_n$ is an epimorphism, then $f$ is a normal epimorphism in the category of long exact sequences. \NoProof
  \end{thmlist}
\end{corollary}

\begin{corollary}[$6$-term exact sequence from composite\DExTag]
  \label{thm:6TermES-From-Composite}%
  In the composite $X\XRA{f} Y \XRA{g} Z$, if $f$, $g$, and $g\Comp f$ are normal maps, then the sequence below is exact. It is functorial with respect to morphisms of such composition diagrams.
  \begin{equation*}
    \xymatrix@R=5ex@C=2em{
    \Ker{f} \ar@{{ |>}->}[r]^-{u} &
    \Ker{gf} \ar[r]^-{\tilde{f}} &
    \Ker{g} \ar[r] &
    \CoKer{f} \ar[r]^-{v} &
    \CoKer{gf} \ar@{-{ >>}}[r]^-{\bar{g}} &
    \CoKer{g}
    }
  \end{equation*}
  Here $u$ is the universal map of kernels, and $v$ is the universal map of cokernels. The maps $\tilde{f}$ and $\bar{g}$ are explained in the proof.
\end{corollary}
\begin{proof}
  We work with the commutative diagram below.
  \begin{equation*}
    \xymatrix@R=5ex@C=2em{
    0 \ar[r] \ar[d] &
    0 \ar[r] \ar[d] &
    \Ker{g} \ar@{{ |>}->}[d] \ar[r] &
    \Ker{v} \ar@{{ |>}->}[d] \ar[r] &
    0 \\
    \Ker{f} \ar@{{ |>}->}[r]^-{\kappa} \ar@{{ |>}->}[d]_{u} &
    X \ar[r]^-{f} \ar@{=}[d] &
    Y \ar@{-{ >>}}[r] \ar[d]^{g} &
    \CoKer{f} \ar[d]^{v} \\
    \Ker{gf} \ar@{{ |>}->}[r]_-{\kappa'} \ar@{-{ >>}}[d] &
    X \ar[r]^-{gf} \ar@{-{ >>}}[d] &
    Z \ar@{-{ >>}}[r] \ar@{-{ >>}}[d] &
    \CoKer{gf} \ar@{-{ >>}}[d] \\
    \CoKer{u} \ar[r] &
    0 \ar[r] &
    \CoKer{g} \ar@{-{ >>}}[r] &
    \CoKer{v}
    }
  \end{equation*}
  Note that the rows in the middle are exact because $f$ and $gf$ are normal. The quadruple $(u,\IdMapOn{X},g,v)$ forms a normal morphism of long exact sequences because (a) $u$ is a normal monomorphism by (\ref{thm:MonomorphismCancellationInKernel}); (b) $g$ is normal by assumption; and (c) $v$ is normal by (\ref{thm:ANN<->NMNIC}).

  Both, the row of cokernels and the row of kernels are normal chain complexes. In all positions but $\Ker{v}$ this is obvious. In position $\Ker{v}$ (\ref{thm:(Co)Ker(ProperMapLESs)}) tells us that the map $\Ker{g}\to \Ker{v}$ is subnormal. With $\KerMap{\Ker{v}\to \ZeroObject}= \IdMapOn{\Ker{v}}$, it follows that it is normal. Thus (\ref{thm:Exact<=>H_n=0}): either chain complex is exact in a given position if and only if its homology in that position vanishes.

  \emph{Step 1\quad $\CoKer{g} \to \CoKer{v}$ is an isomorphism}\quad Indeed, it is a normal epimorphism by (\ref{thm:NormalEpi-Props}). Further, its kernel $K$ is the homology of the bottom chain complex in position $\CoKer{g}$. By (\ref{thm:(Co)Ker(ProperMapLESs)}) this kernel vanishes. By (\ref{thm:IsomorphismRecognition}), the horizontal bottom right map is an isomorphism.

  \emph{Step 2\quad $\Ker{g}\to \Ker{v}$ is a normal epimorphism}\quad Again by \eqref{thm:(Co)Ker(ProperMapLESs)} we have that $0$ is the homology of the normal chain complex $0\to \Ker{g}\to \Ker{v}\to 0$ in position $\Ker{v}$. So, the normal map in question is a normal epimorphism.

  \emph{Step 3\quad $\CoKer{u}\cong \Ker{\Ker{g}\to\Ker{v}}$}\quad First, $\CoKer{u}$ is the homology of the bottom chain complex in the leftmost position. Second, $\Ker{\Ker{g}\to\Ker{v}}$ is the homology of the top chain complex in position $\Ker{g}$. By \eqref{thm:(Co)Ker(ProperMapLESs)}, both terms are isomorphic.

  \emph{Step 4}\quad We assemble the available information to obtain the $6$-term exact sequence:
  \begin{equation*}
    \xymatrix@R=5ex@!C=1.3em{
    \Ker{f} \ar@{{ |>}->}[rr]^-{u} &&
    \Ker{gf} \ar[rr]^-{\tilde{f}} \ar@{-{ >>}}[rd] &&
    \Ker{g} \ar[rr] \ar@{-{ >>}}[rd] &&
    \CoKer{f} \ar[rr]^-{v} &&
    \CoKer{gf} \ar@{-{ >>}}[rr] \ar@{-{ >>}}[rrd]_{\bar{g}} &&
    \CoKer{v} \ar[d]^{\cong} \\
    &&& \CoKer{u} \ar@{{ |>}->}[ru] &&
    \Ker{v}\ar@{{ |>}->}[ru] &&& &&
    \CoKer{g}
    }
  \end{equation*}
  This was to be shown. Since all the constructions and arguments used are functorial, so is the $6$-term exact sequence.
\end{proof}

\begin{subordinate}{On Minimal Conditions}

  The proof of theorem (\ref{thm:LES-Homology}) relies on the relaxed Snake Lemma (\ref{thm:SnakeLemma-Relaxed}), which is valid in sub-di-exact categories. Thus (\ref{thm:LES-Homology}) holds in such categories as well.

\end{subordinate}

\begin{exercises}

\begin{exercise}[Strong relaxed Snake Lemma\DExTag]
  \label{exe:SnakeLemma-Relaxed-Strong}%
  In the morphism of exact sequences below, assume that $a$, $b$, $c$ are normal maps.%
  \index{Snake Lemma!strong relaxed version}%
  \begin{equation*}
    \xymatrix@R=5ex@C=3em{
    \Ker{\alpha} \ar@{{ |>}->}[r] \ar[d] &
    A \ar[d]_{a} \ar[r]^-{\alpha} &
    B \ar[d]_{b} \ar@{-{ >>}}[r]^-{\beta} &
    C \ar[d]^{c} \ar[r] &
    0 \ar[d] \\
    0 \ar[r] &
    X \ar@{{ |>}->}[r]_-{\xi} &
    Y \ar[r]_-{\eta} &
    Z \ar@{-{ >>}}[r] &
    \CoKer{\eta}
    }
  \end{equation*}
  Then the sequence below is exact and functorial.
  \begin{equation*}
    \xymatrix@R=5ex@C=3em{
    0 \ar[r] &
    \Ker{\alpha}\ar[r] &
    \Ker{a} \ar[r] &
    \Ker{b} \ar[r] &
    \Ker{c} \ar`d[l]`[dlll][dlll] \\
    & \CoKer{a}\ar[r] &
    \CoKer{b}\ar[r] &
    \CoKer{c} \ar[r] &
    \CoKer{\eta} \ar[r] &
    0 .
    }
  \end{equation*}
\end{exercise}

\begin{exercise}[Extended Snake Lemmas\DExTag]
  \label{exe:SnakeLemma-Extended}
  Consider the commutative diagram with exact rows below.
  \begin{equation*}
    \xymatrix@R=5ex@C=4em{
    & A \ar[r]^-{a} \ar[d]_{\alpha} &
    B \ar[r]^-{b} \ar[d]_{\beta} &
    C \ar@{-{ >>}}[r] \ar[d]_{\gamma} &
    D \ar[r] \ar[d]^{\delta} &
    0 \\
    0 \ar[r] &
    W \ar[r]_-{w} &
    X \ar[r]_-{x} &
    Y \ar[r]_-{y} &
    Z
    }
  \end{equation*}
  If the vertical maps are normal, then the following hold:
  \begin{enumerate}[(i)]
    \item If $\CoKer{\alpha}=0$, then the sequence below is exact and functorial with respect to morphisms of such diagrams:
          {\small
          \begin{equation*}
            \Ker{\alpha} \longrightarrow \Ker{\beta} \longrightarrow\Ker{\gamma} \longrightarrow\Ker{\delta} \XRA{\partial} \CoKer{\beta} \longrightarrow \CoKer{\gamma} \longrightarrow \CoKer{\delta}
          \end{equation*} }
    \item If $\Ker{\delta}=0$, then the sequence below is exact and functorial with respect to morphisms of such diagrams:
          {\small
          \begin{equation*}
            \Ker{\alpha} \longrightarrow \Ker{\beta} \longrightarrow \Ker{\gamma} \XRA{\partial} \CoKer{\alpha} \longrightarrow \CoKer{\beta} \longrightarrow \CoKer{\gamma} \longrightarrow \CoKer{\delta}
          \end{equation*} }
  \end{enumerate}
\end{exercise}

\begin{exercise}[Kernel/cokernel sequences from a square of proper maps\DExTag]
  \label{exe:Ker/CoKer-SequencesFromSquareOfNormalMaps}%
  Consider the morphism of exact sequences constructed from the center square of normal maps.
  \begin{equation*}
    \xymatrix@R=5ex@C=4em{
    \Ker{a} \ar@{{ |>}->}[r] \ar[d]_{\Ker{\alpha,\beta}} &
    A \ar[r]^-{a} \ar[d]_{\alpha} &
    B \ar@{-{ >>}}[r] \ar[d]^{\beta} &
    \CoKer{a} \ar[d]^{\CoKer{\alpha,\beta}} \\
    \Ker{x} \ar@{{ |>}->}[r] &
    X \ar[r]_-{x} &
    Y \ar@{-{ >>}}[r] &
    \CoKer{x}
    }
  \end{equation*}
  If $\CoKer{\Ker{\alpha,\beta}}=0$, show that the sequence below is exact:
  \begin{equation*}
    \xymatrix@R=5ex@C=3em{
    0\ar[r] &
    \Ker{\Ker{\alpha,\beta}} \ar@{{ |>}->}[r] &
    \Ker{\alpha} \ar[r] &
    \Ker{\beta} \ar[r] &
    \Ker{\CoKer{\alpha,\beta}} \ar`d[l]`[dlll][dlll] \\
    & \CoKer{\alpha} \ar[r] &
    \CoKer{\beta} \ar@{-{ >>}}[r] &
    \CoKer{\CoKer{\alpha,\beta}} \ar[r] &
    0
    }
  \end{equation*}
  Further, if $\Ker{\CoKer{\alpha,\beta}}=0$, show that the sequence below is exact:
  \begin{equation*}
    \xymatrix@R=5ex@C=3em{
    0 \ar[r] &
    \Ker{\Ker{\alpha,\beta}} \ar@{{ |>}->}[r] &
    \Ker{\alpha} \ar[r] &
    \Ker{\beta} \ar`d[l]`[dll][dll] \\
    & \CoKer{\Ker{\alpha,\beta}} \ar[r] &
    \CoKer{\alpha} \ar[r] &
    \CoKer{\beta} \ar@{-{ >>}}[r] &
    \CoKer{\CoKer{\alpha,\beta}} \ar[r] &
    0
    }
  \end{equation*}
  Now use the symmetry in the center square of the above diagram to develop $7$-term exact sequences involving the kernels and cokernels of $a$ and $x$.
\end{exercise}
\end{exercises}
\newpage
\section[Homological Self Duality]{Homological Self Duality}
\label{sec:HomologicalSelfDuality}

In (\ref{term:DiExtensiveConditions}\ref{Ax:0->DPNnverseIsNormal}), we said that a \ZExact\ category is homologically self-dual if dinversion turns every antinormal decomposition of the zero map into a normal map and, hence, yields a di-extension. Here, we show that this condition has a much wider scope: For any normal chain complex $C$, it is equivalent to the condition that the cokernel construction of homology $\HmlgyCoKer{n}{C}$ and the kernel construction $\HmlgyKer{n}{C}$ coincide. Perhaps a bit surprisingly, it is also equivalent to the validity of the Pure Snake Lemma, as well as the validity of the Third Isomorphism Theorem.

We begin by examining properties of certain morphisms of short exact sequence which are directly related to the manner in which homological self-duality can be used to construct di-extensions as in Figure \ref{fig:DiExtensionTypes}.

\begin{definition}[Totally normal sequence of mono/epimorphisms\ZExactTag]
  \label{def:TotallyNormalSequenceMorphisms}
  A sequence of monomorphisms (resp.\ epimorphisms) $A\to B\to \cdots \to Z$ is \Defn{totally normal} if all composites of these morphisms are normal monomorphisms (resp.\ normal epimorphisms). %
  \index{totally!normal sequence of monomorphisms}\index{totally!normal sequence of epimorphisms}
\end{definition}

With (\ref{thm:PushoutRecognize-Categorical}), let us observe that every totally normal sequence of monomorphisms yields a special morphism of short exact sequences in which the right hand square is a pushout:
\stepcounter{theorem}
\begin{equation}\label{fig:SESs-TotallyNormalSeqMonos}
  \vcenter{
  \xymatrix@R=5ex@C=4em{
  A \ar@{{ |>}->}[r]^-{\alpha} \ar@{=}[d] &
  B \ar@{-{ >>}}[r]^-{q\DefEq \CoKerMap{\alpha}} \ar@{{ |>}->}[d]_{\beta} \PushRD{rd}&
  Q \ar[d]^{\gamma} \\
  A \ar@{{ |>}->}[r]_{\beta\alpha} &
  C \ar@{-{ >>}}[r]_-{r\DefEq \CoKerMap{\beta\alpha}} &
  R
  }
  }
\end{equation}
While pushouts preserve (normal) epimorphisms (\ref{thm:PushoutPreservesNormalEpis}), in general, they fail to preserve (normal) monomorphisms.

\begin{lemma}[Properties of Diagram \eqref{fig:SESs-TotallyNormalSeqMonos} \ZExactTag]
  \label{thm:SES-TNSM-Diagram-Props}%
  \label{thm:Pullback/Pushout-Recognition-HSD} 
  In diagram \eqref{fig:SESs-TotallyNormalSeqMonos} the left hand square is a pullback, and so $\alpha=\KerMap{r\beta}$. If $\gamma$ is a monomorphism, then the square on the right is bicartesian.
\end{lemma}
\begin{proof}
  To see that the left hand square is a pullback, consider maps $a\from T\to A$ and $b\from T\to B$ with $\beta b=\beta \alpha a$. The monic property of $\beta$ implies that $b=\alpha a$. Thus $a$ is the unique filler required for the pullback property to hold. With (\ref{thm:NormalMono-Props}.\ref{thm:Kernel(gf)}), we see that $\alpha=\KerMap{r\beta}$.

  Now, suppose $\gamma$ is a monomorphism. Since the right hand square is a pushout, the cokernels of $\beta$ and $\gamma$ are related by an isomorphism:
  \begin{equation*}
    \xymatrix@R=5ex@C=4em{
    A \ar@{{ |>}->}[r]^-{\alpha} \ar@{=}[d] &
    B \ar@{-{ >>}}[r]^-{q\DefEq \CoKerMap{\alpha}} \ar@{{ |>}->}[d]_{\beta} \PushRD{rd}&
    Q \ar[d]^{\gamma} \\
    A \ar@{{ |>}->}[r]_{\beta\alpha} &
    C \ar@{-{ >>}}[r]_-{r\DefEq \CoKerMap{\beta\alpha}} \ar@{-{ >>}}[d] &
    R \ar@{-{ >>}}[d] \\
    & S \ar@{=}[r] &
    S
    }
  \end{equation*}
  With (\ref{thm:PullbackRecognition-KernelSide-1}) we see that the right hand square is a pullback.
\end{proof}

Dually, every totally normal sequence of epimorphisms yields a special morphism of short exact sequences in which the left hand square is a pullback:
\stepcounter{theorem}
\begin{equation}\label{fig:SESs-TotallyNormalSeqEpis}
  \vcenter{
  \xymatrix@R=5ex@C=4em{
  K \ar@{{ |>}->}[r]^-{k\DefEq \KerMap{\eta\xi}} \ar[d]_{\kappa} \PullLU{rd} &
  X \ar@{-{ >>}}[r]^-{\eta\xi} \ar@{-{ >>}}[d]^{\xi}&
  Z \ar@{=}[d] \\
  L \ar@{{ |>}->}[r]_-{l\DefEq \KerMap{\eta}} &
  Y \ar@{-{ >>}}[r]_-{\eta} &
  Z
  }
  }
\end{equation}
While pullbacks preserve (normal) monomorphisms (\ref{thm:PullbackPreservesNormalMonos}), in general, they fail to preserve (normal) epimorphisms.

\begin{lemma}[Properties of Diagram \eqref{fig:SESs-TotallyNormalSeqEpis}\ZExactTag]
  \label{thm:SES-TNSE-Diagram-Props}
  In diagram \eqref{fig:SESs-TotallyNormalSeqEpis} the left hand square is a pushout, and so $\eta=\CoKerMap{\xi k}$. If $\kappa$ is an epimorphism, then the square on the left is bicartesian. \NoProof
\end{lemma}

\begin{proposition}[Dinversions of $\ZeroMap$ and \eqref{fig:SESs-TotallyNormalSeqMonos} / \eqref{fig:SESs-TotallyNormalSeqEpis}\ZExactTag]
  \label{thm:Dinversion0TotSeqMonos/TotSeqEpis}%
  In a \ZExact\ category, the following conditions are equivalent.
  \begin{tfae}
    \item \label{thm:Dinversion0TotSeqMonos/TotSeqEpis-TotMonos}%
    For every totally normal sequence of monomorphisms and its morphism of short exact sequences as in \eqref{fig:SESs-TotallyNormalSeqMonos}, the normal monomorphism $\beta$ pushes forward to a normal monomorphism $\gamma$.
    \item \label{thm:Dinversion0TotSeqMonos/TotSeqEpis-HSD}%
    Every dinversion of $\ZeroMap$ is a normal map. {\HSDTag}
    \item \label{thm:Dinversion0TotSeqMonos/TotSeqEpis-TotEpis}%
    For every totally normal sequence of epimorphisms and its morphism of short exact sequences, as in \eqref{fig:SESs-TotallyNormalSeqEpis}, the normal epimorphism $\xi$ pulls back to a normal epimorphism $\kappa$.
  \end{tfae}
\end{proposition}
\begin{proof}
  (\ref{thm:Dinversion0TotSeqMonos/TotSeqEpis-TotMonos}) $\Rightarrow$ (\ref{thm:Dinversion0TotSeqMonos/TotSeqEpis-HSD})\quad Given an antinormal decomposition $\ZeroMap = em$ of $\ZeroMap$, form the pullback on the top left:
  \begin{equation*}
    \xymatrix@R=5ex@C=4em{
    \DiagObj \ar@{{ |>}->}[r]^-{k} \ar@{=}[d] &
    \Ker{e} \ar@{-{ >>}}[r]^-{q} \ar@{{ |>}->}[d]^{\kappa} &
    \CoKer{k} \ar@{.>}[d]^{\underline{\kappa}} \\
    \DiagObj \ar@{{ |>}->}[r]_-{m} \ar@{-{ >>}}[d] &
    \DiagObj \ar@{-{ >>}}[r]_-{\varepsilon} \ar@{-{ >>}}[d]_{e} &
    \CoKer{m} \\
    \ZeroObject \ar@{{ |>}->}[r] &
    \DiagObj
    }
  \end{equation*}
  Then the universal map $\underline{\kappa}\from\CoKer{k}\to \CoKer{m}$ is a normal monomorphism by hypothesis. So the dinversion $\varepsilon\kappa$ of $em$ is a normal map.

  (\ref{thm:Dinversion0TotSeqMonos/TotSeqEpis-HSD}) $\Rightarrow$ (\ref{thm:Dinversion0TotSeqMonos/TotSeqEpis-TotEpis})\quad Given a diagram as in \eqref{fig:SESs-TotallyNormalSeqEpis}, consider the construction below.
  \begin{equation*}
    \xymatrix@R=5ex@C=4em{
    & K \ar@{{ |>}->}[d]_{k} \ar@{.>}[r]^-{\kappa} &
    L \ar@{{ |>}->}[d]^{l} \\
    \Ker{\xi} \ar@{{ |>}->}[r] \ar@{-{ >>}}[d] &
    X \ar@{-{ >>}}[r]_-{\xi} \ar@{-{ >>}}[d]_{\eta\xi} &
    Y \ar@{-{ >>}}[d]^{\eta} \\
    \ZeroObject \ar@{{ |>}->}[r] &
    Z \ar@{=}[r] &
    Z
    }
  \end{equation*}
  Since the dinversion $\xi k$ of $\ZeroMap$ is normal, we know that the universal map $K\to \Ker{\CoKerMap{\xi k}}$ is a normal monomorphism. However, we know that $\eta=\CoKerMap{\xi k}$ from (\ref{thm:SES-TNSE-Diagram-Props}) that. So, $\Ker{\CoKerMap{\xi k}}=L$, which implies the claim.

  (\ref{thm:Dinversion0TotSeqMonos/TotSeqEpis-TotEpis}) $\Rightarrow$ (\ref{thm:Dinversion0TotSeqMonos/TotSeqEpis-TotMonos})\quad Given a diagram as in \eqref{fig:SESs-TotallyNormalSeqMonos}, we obtain the square on the bottom right from the pushout property of the square $B\rightrightarrows R$:
  \begin{equation*}
    \xymatrix@R=5ex@C=4em{
    A \ar@{{ |>}->}[r]^-{\alpha} \ar@{=}[d] &
    B \ar@{-{ >>}}[r]^-{q\DefEq \CoKerMap{\alpha}} \ar@{{ |>}->}[d]_{\beta} \PushRD{rd}&
    Q \ar[d]^{\gamma} \\
    A \ar@{{ |>}->}[r]_{\beta\alpha} &
    C \ar@{-{ >>}}[r]_-{r\DefEq \CoKerMap{\beta\alpha}} \ar@{-{ >>}}[d] &
    R \ar@{-{ >>}}[d]^{t} \\
    & S \ar@{=}[r] &
    S
    }
  \end{equation*}
  The normal epimorphism $r$ pulls back to a normal epimorphism $s\from B\to \Ker{t}$. If $k=\KerMap{t}$, then we find:
  \begin{equation*}
    \KerMap{s} = \KerMap{ks} = \KerMap{r\beta} \overset{\text{(\ref{thm:SES-TNSM-Diagram-Props})} }{=} \alpha
  \end{equation*}
  But then $q=s$, and $\gamma=k$ is a normal monomorphism.
\end{proof}

\begin{proposition}[Criteria for homological self-duality\ZExactTag]
  \label{thm:HomologicalSelfDuality-Recognize}
  In a \ZExact\ category the following are equivalent.
  \begin{tfae}
    \item \label{thm:HomologicalSelfDuality-AxHSD}%
    Every dinversion of $\ZeroMap$ is a normal map {\HSDTag}
    \item \label{thm:HomologicalSelfDuality-Recognize-3rdIso}%
    \emph{The Third Isomorphism Property holds}\quad Every totally normal sequence of monomorphisms $X\NMono Y\NMono Z$ yields a short exact sequence \ $(Y/X)\NMono (Z/X)\NEpi (Z/Y)$. %
    \index{Isomorphism Theorem!condition for Third}%
    \item \label{thm:HomologicalSelfDuality-Recognize-PureSnake}%
    \emph{Pure Snake Condition}\quad For every morphism of short exact sequences
    \begin{equation*}
      \xymatrix@R=5ex@C=4em{
      A \ar@{{ |>}->}[r]^-{\alpha} \ar@{{ |>}->}[d]_{a} &
      B \ar@{-{ >>}}[r]^-{\beta} \ar@{=}[d] &
      C \ar@{-{ >>}}[d]^{c} \\
      X \ar@{{ |>}->}[r]_-{\xi} &
      Y \ar@{-{ >>}}[r]_-{\eta} &
      Z
      }
    \end{equation*}
    there is a functorial isomorphism $\CoKer{a}\cong \Ker{c}$. %
    \index{pure snake!condition}
    \item \label{thm:HomologicalSelfDuality-Recognize-HomologicalSelfD}%
    \emph{Homology is self-dual}\quad For every normal chain complex $(C,d)$, and every $n\in \ZNr$, $\HmlgyCoKer{n}{C} \cong \HmlgyKer{n}{C}$. %
    \index{homological!self-duality}%
  \end{tfae}
\end{proposition}
\begin{proof}
  (\ref{thm:HomologicalSelfDuality-AxHSD}) $\implies$ (\ref{thm:HomologicalSelfDuality-Recognize-3rdIso})\quad We know from Proposition  (\ref{thm:Dinversion0TotSeqMonos/TotSeqEpis}) that a category is homologically self-dual if and only if it satisfies condition (\ref{thm:Dinversion0TotSeqMonos/TotSeqEpis}.\ref{thm:Dinversion0TotSeqMonos/TotSeqEpis-TotMonos}). That the latter condition implies (\ref{thm:HomologicalSelfDuality-Recognize-3rdIso}) is a special case of Exercise \ref{exe:FactoringCoKer}.

  (\ref{thm:HomologicalSelfDuality-Recognize-3rdIso}) $\implies$ (\ref{thm:HomologicalSelfDuality-Recognize-PureSnake})\quad We observe that $A \overset{a}{\NMono} X \overset{\xi}{\NMono} Y$ is a totally normal sequence of monomorphisms. Thus we obtain the short exact sequence
  \begin{equation*}
    \xymatrix@R=5ex@C=4em{
    \CoKer{a}=(X/A) \ar@{{ |>}->}[r] &
    \CoKer{\alpha} = (B/A) \ar@{-{ >>}}[r] &
    \CoKer{\xi} = (Y/X)
    }
  \end{equation*}
  This means that $\CoKer{a}\cong \Ker{c}$, as required.

  (\ref{thm:HomologicalSelfDuality-Recognize-PureSnake}) $\implies$ (\ref{thm:HomologicalSelfDuality-Recognize-HomologicalSelfD})\quad The construction of the (co)kernel homology of $(C,d)$ in position $n$ fits into this pure snake diagram, which implies the claim.
  
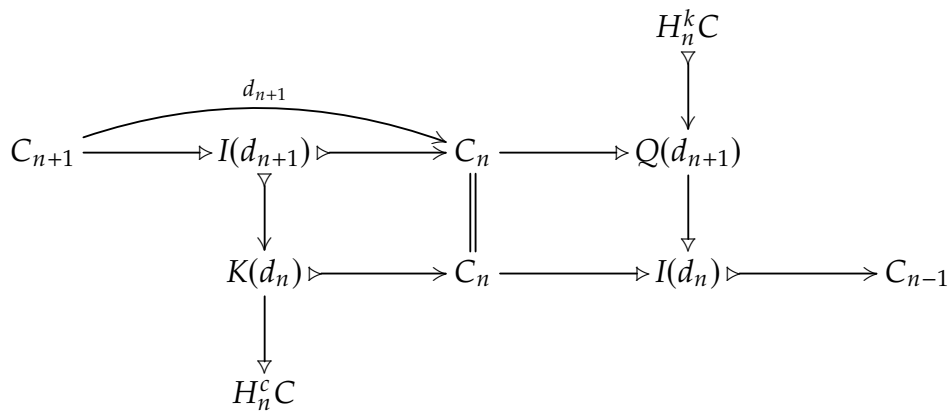
\begin{equation*}
    \xymatrix@R=5ex@C=4em{
    &&& \HmlgyKer{n}{C} \ar@{{ |>}->}[d] \\
    C_{n+1} \ar@{-{ >>}}[r] \ar@/^3ex/[rr]^-{d_{n+1}}&
    \Img{d_{n+1}} \ar@{{ |>}->}[r] \ar@{{ |>}->}[d] &
    C_{n} \ar@{-{ >>}}[r] \ar@{=}[d] &
    \CoKer{d_{n+1}} \ar@{-{ >>}}[d] \\
    & \Ker{d_{n}} \ar@{{ |>}->}[r] \ar@{-{ >>}}[d]&
    C_{n} \ar@{-{ >>}}[r] &
    \Img{d_{n}} \ar@{{ |>}->}[r] &
    C_{n-1} \\
    & \HmlgyCoKer{n}{C}
    }
  \end{equation*}

  (\ref{thm:HomologicalSelfDuality-Recognize-HomologicalSelfD}) $\Leftrightarrow$ (\ref{thm:HomologicalSelfDuality-AxHSD})\quad An antinormal pair $(\varepsilon,\mu)$ determines a normal chain complex $C$, as on the left below, if and only if   $\varepsilon\mu=\ZeroObject$. Construct $\kappa$ and $\pi$ by dinversion, then the top left square as a pullback, and the bottom right square as a pushout yields the diagram on the right.
  \begin{center}
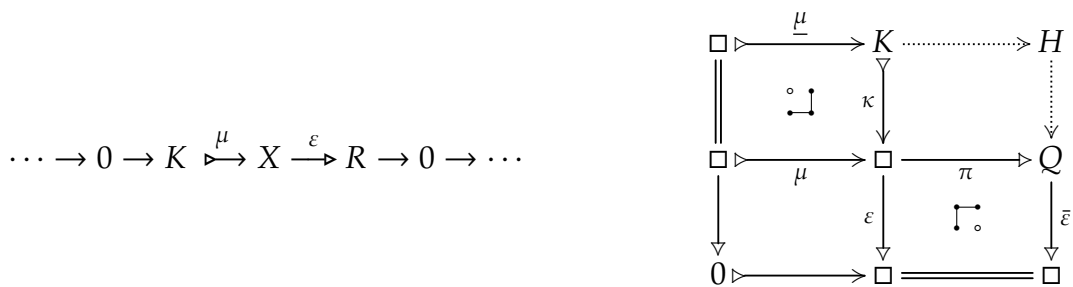
    \begin{minipage}[m]{7cm}
      \begin{equation*}
        \cdots\to 0 \to K \overset{\mu}{\NMono} X \overset{\varepsilon}{\NEpi} R \to 0\to \cdots
      \end{equation*}
    \end{minipage} \qquad\qquad
    \begin{minipage}[m]{6cm}
      \begin{equation*}
        \xymatrix@R=5ex@C=4em{
        \DiagObj \ar@{{ |>}->}[r]^-{\underline{\mu}} \ar@{=}[d] \PullLU{rd} &
        K \ar@{.>}[r] \ar@{{ |>}->}[d]_{\kappa} &
        H \ar@{.>}[d] \\
        \DiagObj \ar@{{ |>}->}[r]_{\mu} \ar@{-{ >>}}[d] &
        \DiagObj \ar@{-{ >>}}[d]_-{\varepsilon} \ar@{-{ >>}}[r]_-{\pi} \PushRD{rd}               &
        Q \ar@{-{ >>}}[d]^{\bar{\varepsilon}} \\
        0 \ar@{{ |>}->}[r]                                                               &
        \DiagObj \ar@{=}[r]                                                               &
        \DiagObj
        }
      \end{equation*}
    \end{minipage}
  \end{center}
  Placing $X$ in position $0$ of the chain complex $C$, observe the following:
  \begin{enumerate}
    \item With $\kappa=\KerMap{\varepsilon}$, we see that $\bar{\mu}$ is the unique factorization of $\mu$ through $K$. So, $\CoKer{\underline{\mu}}= \HmlgyCoKer{0}{C}$ is the cokernel constructed homology of $C$ in position $X$.
    \item With $\pi=\CoKerMap{\mu}$, we see that $\bar{\varepsilon}$ is the unique factorization of $\varepsilon$ through $Q$. So, $\Ker{\bar{\varepsilon}}=\HmlgyKer{0}{C}$ is the kernel constructed homology of $C$ in position $X$.
  \end{enumerate}
  Thus, if the antinormal pair $(\pi,\kappa)$ composes to a normal map, then its normal factorization consists of $\CoKer{\underline{\mu}}$ and $\Ker{\bar{\varepsilon}}$. So, we obtain $\HmlgyCoKer{0}{C}\cong H\cong \HmlgyKer{0}{C}$.

  On the other hand, if the two homology constructions are naturally isomorphic, then they provide a normal factorization of the composite $\pi\kappa$. - This completes the proof.
\end{proof}

The following paraphrasing of Proposition~\ref{thm:HomologicalSelfDuality-Recognize} makes the relationship between homological self-duality and di-extensions explicit:

\begin{proposition}[Homological self-duality and the $(\Prdct{3}{3})$-Lemma\ZExactTag]
  \label{thm:RR33}
  For a \ZExact\ category $\Ctgry{X}$ the following are equivalent:
  \begin{tfae}
    \item $\Ctgry{X}$ is homologically self-dual.
    \item Every dinversion of $\ZeroMap$ is a normal map. {\HSDTag}
    \item \label{thm:HomologicalSelfDuality-Recognize-Right-3x3}%
    In every $(\Prdct{3}{3})$-diagram, with a $\ZeroObject$ in the bottom left corner, in which the rows and the left and middle column are short exact sequences
    \begin{equation*}
      \vcenter{
      \xymatrix@R=5ex@C=3em{
      A \ar@{{ |>}->}[r]^-{\alpha} \ar@{=}[d] &
      B \ar@{-{ >>}}[r]^-{c} \ar@{{ |>}->}[d]_{\beta} &
      Q \ar[d]^{\gamma} \\
      A \ar@{-{ >>}}[d] \ar@{{ |>}->}[r]_{\beta\alpha} &
      C \ar@{-{ >>}}[d]_-b \ar@{-{ >>}}[r]_-{c'} &
      R \ar[d]\\
      0 \ar@{{ |>}->}[r] & S \ar@{=}[r] & S
      }}
    \end{equation*}
    it follows that the right hand column is a short exact sequence.\NoProof
  \end{tfae}
\end{proposition}

Condition (\ref{thm:RR33}.\ref{thm:HomologicalSelfDuality-Recognize-Right-3x3}) has an equivalent sibling, obtained by mirroring the diagram over the $AS$-diagonal.

\begin{corollary}[More criteria for homological self-duality\ZExactTag]
  \label{thm:HomologicalSelfDuality-Recognize-II}%
  A \ZExact\ category is homologically self-dual if and only if any of the conditions below is satisfied.
  \begin{tfae}
    \setcounter{enumi}{5}
    \item \label{thm:HomologicalSelfDuality-Recognize-TotPull}%
    In a morphism of short exact sequences as in \eqref{fig:SESs-TotallyNormalSeqEpis}, the normal epimorphism $\xi$ pulls back to a normal epimorphism $\kappa$.
    \item \label{thm:HomologicalSelfDuality-Recognize-3rdIso-Op}%
    Every totally normal sequence of epimorphisms $X\NEpi Y\NEpi Z$ yields a short exact sequence \ $\Ker{X\NEpi Y}\NMono \Ker{X\NEpi Z}\NEpi \Ker{Y\NEpi Z}$. %
    \index{Isomorphism Theorem!condition for Third}%
  \end{tfae}
\end{corollary}
\begin{proof}
  These claims may be verified directly. Alternatively, (\ref{thm:HomologicalSelfDuality-Recognize-TotPull}) is a paraphrasing of  (\ref{thm:Dinversion0TotSeqMonos/TotSeqEpis}.\ref{thm:Dinversion0TotSeqMonos/TotSeqEpis-TotEpis}), and (\ref{thm:HomologicalSelfDuality-Recognize-3rdIso-Op}) is dual to (\ref{thm:HomologicalSelfDuality-Recognize}.\ref{thm:HomologicalSelfDuality-Recognize-3rdIso}).
\end{proof}

\begin{corollary}[Cobase change of a normal epimorphism]
  \label{thm:CobaseChange-NormalEpi}%
  In a homologically self-dual category, the pushout of a short exact sequence $K\NMono X \NEpi Y$ along a normal epimorphism $h\from K\NEpi L$ whose kernel $m$ is normal in $X$ yields a short exact sequence as in the diagram below.
  \begin{equation*}
    \xymatrix@R=5ex@C=3em{
    K \PushRD{rd} \ar@{{ |>}->}[r]^-{k} \ar@{-{ >>}}[d]_h &
    X \ar@{-{ >>}}[r]^-{f} \ar@{-{ >>}}[d]^{h'} &
    Y \ar@{=}[d] \\
    L \ar@{{ |>}->}[r]_-{k'} &
    A \ar@{-{ >>}}[r]_-{f'} &
    Y
    }
  \end{equation*}
\end{corollary}
\begin{proof}
  The sequence $k\Comp m$ is a totally normal sequence of monomorphisms. So, the claim is equivalent to (\ref{thm:Dinversion0TotSeqMonos/TotSeqEpis}.\ref{thm:Dinversion0TotSeqMonos/TotSeqEpis-TotMonos}).
\end{proof}

\begin{proposition}[Homological self-duality of $\NMonoCat{X}$, $\SESCat{X}$, $\NEpiCat{X}$\DExTag]
  \label{thm:ANN->HSDOf-NM(X),SES(X),NE(X)}%
  For a di-exact category $\Ctgry{X}$, the categories $\NMonoCat{X}$, $\SESCat{X}$, and $\NEpiCat{X}$ are homologically self-dual. %
  \index[not]{s!$\SESCat{X}$\IndSep category of short exact sequences in $\Ctgry{X}$}%
  \index[not]{n!$\NMonoCat{X}$\IndSep category of normal monomorphisms in $\Ctgry{X}$}%
  \index[not]{n!$\NEpiCat{X}$\IndSep category of normal epimorphisms in $\Ctgry{X}$}%
\end{proposition}
\begin{proof}
  We know that the categories $\NMonoCat{X}$, $\SESCat{X}$, and $\NEpiCat{X}$ are equivalent. So, it suffices to prove that $\NMonoCat{X}$ is homologically self-dual. We verify that criterion (\ref{thm:Dinversion0TotSeqMonos/TotSeqEpis}.\ref{thm:Dinversion0TotSeqMonos/TotSeqEpis-TotMonos}) holds. The top four rows of the diagram below present a the situation of (\ref{fig:SESs-TotallyNormalSeqMonos}) in $\NMonoCat{X}$. Here we used that, by Proposition~\ref{thm:ANN<->PointwiseSES}, a short exact sequence in $\NMonoCat{X}$ is pointwise short exact. In the di-exact category $\Ctgry{X}$, the $(\Prdct{3}{3})$-diagrams in the front and back are di-extensions, and $\beta''$ is a normal monomorphism.
  \begin{equation*}
    \xymatrix@!0@C=5em@R=6ex{
    & A \ar@{=}[ld]  \ar@{{ |>}->}[rr]^-{\alpha} \ar@{{ |>}->}[dd]|\hole^(.2){a} &&
    B \ar@{{ |>}->}[ld]_-{\beta} \ar@{-{ >>}}[rr]^-{\CoKerMap{\alpha}} \ar@{{ |>}->}[dd]|\hole^(.3){b} &&
    Q \ar@{{ |>}->}[dd]^{q} \ar@{.>}[ld]^{\gamma} \\
    A \ar@{{ |>}->}[rr]^(.7){\beta\alpha} \ar@{{ |>}->}[dd]_{a} &&
    C \ar@{-{ >>}}[rr]_(.7){\CoKerMap{\beta\alpha}} \ar@{{ |>}->}[dd]_(0.3){c} &&
    R \ar@{{ |>}->}[dd]^(.3){r} \\
    & A' \ar@{=}[ld]  \ar@{{ |>}->}[rr]|(0.52)\hole_(.3){\alpha'} \ar@{-{ >>}}[dd]|\hole &&
    B'  \ar@{-{ >>}}[rr]|(.49)\hole \ar@{{ |>}->}[ld]_-{\beta'} \ar@{-{ >>}}[dd]|\hole^(.3){b'}&&
    Q' \ar@{.>}[ld]^{\gamma'} \ar@{-{ >>}}[dd]^{q'}\\
    A' \ar@{{ |>}->}[rr]^(0.75){\beta'\alpha'} \ar@{-{ >>}}[dd] &&
    C' \ar@{-{ >>}}[rr] \ar@{-{ >>}}[dd]^(0.25){c'} &&
    R' \ar@{-{ >>}}[dd]^(.3){r'} \\
    & \CoKer{a} \ar@{{ |>}->}[rr]|(0.52)\hole \ar@{=}[ld] &&
    B'' \ar@{{ |>}->}[ld]_-{\beta''} \ar@{-{ >>}}[rr]|\hole &&
    Q'' \ar@{.>}[ld]^{\gamma''} \\
    \CoKer{a} \ar@{{ |>}->}[rr] &&
    C''	\ar@{-{ >>}}[rr] &&
    R''
    }
  \end{equation*}
  We must show that $(\gamma,\gamma')$ is a normal monomorphism in $\NMonoCat{X}$ or, equivalently, that the square $R\rightrightarrows Q'$ is a pullback of normal monomorphisms; see (\ref{thm:KernelsInNMono(X)-CoKernelsInNEpi(X)}). To see this, note first that $\gamma$, $\gamma'$, $\gamma''$ are normal monomorphisms by homological self-duality in $\Ctgry{X}$.  Then the square $R\rightrightarrows Q'$ is a pullback by Proposition~\ref{thm:PullbackRecognition-KernelSide-1}.
\end{proof}

For identifying homologically self-dual categories, the following lemma provides a useful tool. We rely on background from Section \ref{sec:NormalCats}.

\begin{lemma}[Criterion for homological self-duality]
  \label{thm:TNSM-Mono->NormalMono}
  In a \ZExact\ regular category, if the morphism $\gamma$ in diagram \eqref{fig:SESs-TotallyNormalSeqMonos} is a monomorphism, then it is a normal monomorphism.
\end{lemma}
\begin{proof}
  We extend diagram \eqref{fig:SESs-TotallyNormalSeqMonos} as follows:
  
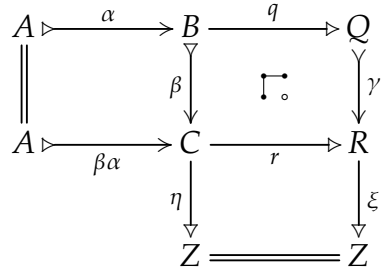
\begin{equation*}
    \xymatrix@R=5ex@C=4em{
    A \ar@{{ |>}->}[r]^-{\alpha} \ar@{=}[d] &
    B \ar@{-{ >>}}[r]^-{q} \ar@{{ |>}->}[d]_{\beta} \PushRD{rd}&
    Q \ar@{{ >}->}[d]^{\gamma} \\
    A \ar@{{ |>}->}[r]_{\beta\alpha} &
    C \ar@{-{ >>}}[r]_-{r} \ar@{-{ >>}}[d]_{\eta} &
    R \ar@{-{ >>}}[d]^{\xi} \\
    & Z \ar@{=}[r] &
    Z
    }
  \end{equation*}
  Thus $\eta\from  C\to Z$ and $\xi\from R\to Z$ represent cokernels of $\beta$  and $\gamma$, respectively. Since the top right square is a pushout, the bottom right square results via (\ref{thm:Pushout->IsoOfCoKers}). We show that $\gamma=\KerMap{\xi}$ by checking the universal property.

  Consider $x\from X\to R$ such that $\xi\Comp x=0$. The normal epimorphism $\CoKerMap{\beta\alpha}$ pulls back along $x$ to a regular epimorphism $\bar{c}$.
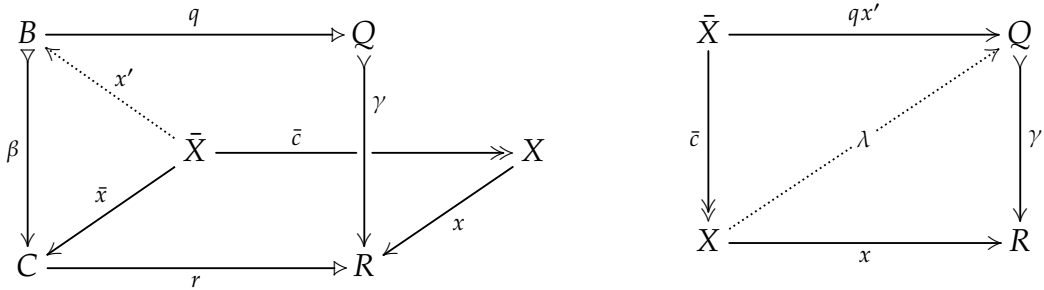
  \begin{equation*}
    \xymatrix@R=5ex@C=4em{
    B \ar@{-{ >>}}[rr]^-{q} \ar@{{ |>}->}[dd]_{\beta} &&
    Q \ar@{{ >}->}[dd] ^(.3){\gamma} \\
    & \bar{X} \ar@{-{>>}}[rr]|\hole^(.3){\bar{c}} \ar[ld]_{\bar{x}} \ar@{.>}[lu]_{x'} &&
    X \ar[ld]^{x} \\
    C \ar@{-{ >>}}[rr]_-{r} &&
    R
    }\qquad\qquad \xymatrix@R=5ex@C=4em{
    \bar{X} \ar[rr]^-{q x'} \ar@{-{>>}}[dd]_{\bar{c}} &&
    Q \ar@{{ >}->}[dd]^{\gamma} \\ \\
    X \ar[rr]_-{x} \ar@{.>}[rruu]|-{\ \lambda\ }&&
    R
    }
  \end{equation*}
  The computation $\eta\bar{x}=\xi r\bar{x}=\xi x\bar{c}=\ZeroMap$ yields a map $x'\from \bar{X}\to B$ unique with $\bar{x}=\beta x'$. From these data, we compose the square on the right. It commutes because
  \begin{equation*}
    x\bar{c} = r\bar{x} = r\beta x' = \gamma qx'.
  \end{equation*}
  The regular epimorphism $\bar{c}$ is strong by (\ref{sec:Surjectivity}), and so the square on the right has a unique filler $\lambda$ rendering it commutative. This implies that $\gamma$ is normal.
\end{proof}

We close this section by showing that, in the category of commutative monoids, homology is self-dual, even though this category is not normal, as follows from (\ref{thm:Mono<->0-Kernel}) combined with (\ref{exa:Kernels/CoKernels-Examples}). The argument relies on an explicit computation of the kernel pair of a normal epimorphism. This computation is of independent interest; see Lemma \ref{thm:FibersNormalEpiMorphismCommutativeMonoids}.

\begin{lemma}[Fibers of a normal epimorphism of commutative monoids]
  \label{thm:FibersNormalEpiMorphismCommutativeMonoids}
  Let $f\from B\to C$ be a normal epimorphism of commutative monoids, and let $K$	be its kernel. Then the kernel pair of $f$ consists of all pairs $(b,b')$ for which there exist $k$, $k'\in K$ such that $b+k=b'+k'$.
\end{lemma}
\begin{proof}
  We know that the kernel pair $\KrnlPr{f}$ is the least internal equivalence relation on $B$ whose equivalence class of $\ZeroObject$ is $K$.

  Let $R\subseteq \Prdct{B}{B}$ denote the set of all pairs $(b,b')$ for which there exist $k$, $k'\in K$ such that $b+k=b'+k'$. For simplicity, write $b\sim b'$ whenever $(b,b')\in R$. We show first that $R$ is an internal equivalence relation on $B$ containing $K$ as its zero-class. Indeed, the relation $R$ is reflexive because, for every $b\in B$, $b+0=b+0$. It is symmetric because so is $=$. To see that $R$ is transitive, suppose $b\sim b'$ and $b'\sim b''$. Then there exist $k,k',l,l'\in K$ such that
  \begin{equation*}
    b+k=b'+k'  \quad \text{and}\quad b'+l'=b''+l''  \quad \text{and so}\quad b+k+l'=b'+k'+l'=b'+l'+k'=b''+l''+k'
  \end{equation*}
  Finally, $R$ is an internal relation because if $b+k=b'+k'$ and $c+l=c'+l'$, then $b+c+k+l=b'+c'+k'+l'$, which means that $(b+c)\sim (b'+c')$. By design, $(\Prdct{K}{\Set{\ZeroObject}})\union (\Prdct{\Set{\ZeroObject}}{K})\subseteq R$. Furthermore, if $b\sim 0$, then $b+k=k'$ for some $k$, $k'\in K$, which implies that $b\in K$. Hence $K$ is the zero-class of $R$. - So, $\KrnlPr{f}\subseteq R$.

  Next, we show that $R$ is the \emph{least} internal equivalence relation on $B$ of which $K$ is the zero-class. Let $S\subseteq \Prdct{B}{B}$ be such an internal equivalence relation. Then, for every $k\in K$, $0Sk$ and $kS0$. So, since $S$ is internal,
  \begin{equation*}
    bS(b+k)  \qquad \text{and}\qquad (b+k)Sb\quad \text{for  every}\ b\in B,\ k\in K
  \end{equation*}
  Now, suppose $(b,b')\in R$. Then $b+k=b'+k'$ for some $k,k'\in K$. As $S$ is reflexive: $(b+k)S(b'+k')$ and, hence,
  \begin{equation*}
    bS(b+k)S(b'+k')Sb'
  \end{equation*}
  So, $(b,b')\in S$ by transitivity. Hence $R\subseteq \KrnlPr{f}$. This was to be shown.
\end{proof}

\begin{theorem}[$\CMon$ is homologically self-dual]
  \label{thm:CMonIsHSD}%
  The category $\CMon$ of commutative monoids is homologically self-dual.
\end{theorem}
\begin{proof}
  We show that in diagram \eqref{fig:SESs-TotallyNormalSeqMonos}, the morphism $\gamma$ is a normal monomorphism. By (\ref{thm:TNSM-Mono->NormalMono}) it suffices to show that $\gamma$ is monic.

  Write $c\DefEq \CoKerMap{\alpha}$ and $c'\DefEq \CoKerMap{\beta \alpha}$. If $\gamma(q)=\gamma(q')$ for $q$, $q'\in Q$, choose $b$, $b'\in B$ with $c(b)=q$ and $c(b')=q'$. Then $c'\beta(b) =\gamma(q)=\gamma(q')= c'\beta(b')$. By (\ref{thm:FibersNormalEpiMorphismCommutativeMonoids}) , there are $a$, $a'\in A$ such that
  \begin{equation*}
    \beta(b+\alpha(a)) = \beta(b)+\beta\alpha(a)=\beta(b')+\beta\alpha(a') = \beta(b' + \alpha(a'))
  \end{equation*}
  Since $\beta$ is monic, $b+\alpha(a)=b' + \alpha(a')$. But then
  \begin{equation*}
    q = c(b) = c(b+\alpha(a)) = c(b'+\alpha(a')) = c(b') =q'
  \end{equation*}
  So, $\gamma$ is monic, and the proof is complete.
\end{proof}

\begin{example}[z-Exact but not \HSDInline]
  \label{exa:NotHSD}%
  The category $\SESCat{\CMon}$ of short exact sequences of commutative monoids is z-exact. However, $\SESCat{\CMon}$ is not homologically self dual. This follows from a counterexample constructed by F.~Afsa.
\end{example}

\begin{example}[Commutative Hopf algebras are \HSDInline]\label{exa:CocommHopf}
  The category $\CHopfK$ of commutative Hopf algebras over a field $\mathbb{K}$ is homologically self-dual. Note that by~\cite{GM-VdL1} we know it is a co-protomodular category (see Section~\ref{sec:Protomodular-SEpi(X)->X}). We do, however, not know whether it is co-regular as well, so we cannot use that it is co-homological to obtain the result. We thus give a direct proof.

  In \eqref{fig:SESs-TotallyNormalSeqEpis} we write $I$ for the cokernel of the kernel of $\kappa$; it is the image of $\xi k$ as a morphism of vector spaces. We need to prove that the induced map $I\to L$ is an isomorphism. A~direct calculation shows that the Hopf algebra $I$ includes into $Y$ as a normal subHopf algebra in the sense of~\cite[Definition~3.4.1]{Montgomery}. In the commutative case, it thus becomes the kernel of its cokernel---see~\cite[Theorem~3.4.6]{Montgomery}. This finishes the proof that \HSDInline\ holds for commutative Hopf algebras.
\end{example}

For the sake of completeness, we recall that the following results which appeared in \eqref{sec:SESMaps} are actually valid in any homologically self-dual category.

\begin{proposition}[Primordial Short $5$-Lemma\HSDTag]
  \label{thm:Short-5-Primordial-bis}%
  Consider a morphism of short exact sequences in which $b$ is a normal map. %
  \index{Short 5-Lemma!}%
  \begin{equation*}
    \xymatrix@R=5ex@C=3em{
    A \ar@{{ |>}->}[r]^{k} \ar@{{ >}-{>>}}[d]_{a} &
    B \ar@{-{ >>}}[r]^{f} \ar[d]_{b} &
    C \ar@{{ >}-{>>}}[d]^{c} \\
    X \ar@{{ |>}->}[r]_{l} &
    Y \ar@{-{ >>}}[r]_{g} &
    Z
    }
  \end{equation*}
  If $a$ and $c$ are both monic and epic, then $a$, $b$, $c$ are isomorphisms.\NoProof
\end{proposition}

\begin{proposition}[Primordial $5$-Lemma\HSDTag]
  \label{thm:5-Lemma-Primordial-bis}%
  Consider a morphism of exact sequences. %
  \index{$5$-Lemma}
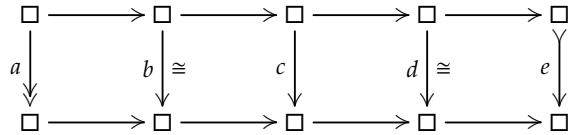
  \begin{equation*}
    \xymatrix@R=5ex@C=3em{
    \DiagObj \ar@{->>}[d]_{a} \ar[r] &
    \DiagObj \ar[d]_{b}^{\cong} \ar[r] &
    \DiagObj \ar[d]_{c} \ar[r] &
    \DiagObj \ar[d]^{\cong}_{d} \ar[r] &
    \DiagObj \ar@{{ >}->}[d]_{e} \\
    \DiagObj \ar[r] &
    \DiagObj \ar[r] &
    \DiagObj \ar[r] &
    \DiagObj \ar[r] &
    \DiagObj
    }
  \end{equation*}
  If $c$ is a normal map, $a$ is an epimorphism, $e$ a monomorphism, and $b$, $d$ isomorphisms, then $c$ is an isomorphism.\NoProof
\end{proposition}

\begin{subordinate}{On Homological self-duality in a wider context}

  \begin{subsubordinate}{On self-duality}
    A pointed category $\Ctgry{X}$ is homologically self-dual if and only if its opposite $\Ctgry{X}^{\op}$ is homologically self-dual.
  \end{subsubordinate}

  \begin{subsubordinate}{Why homologically self-dual categories only allow for a partial deployment of homological methods}
    From what we have seen about exactness and homology in homologically self-dual categories, we ask: To which extent can homological algebra tools and results known from abelian categories be adapted to homologically self-dual categories? - Answer: Very little beyond basic definitions and constructions. Important tools for computation, such as the full (Short) $5$-Lemma, the $(\Prdct{3}{3})$-Lemma, and others need not be available in every homologically self-dual category. This imposes severe limits on our capability to compute effectively. In the following sections we consider incrementally more demanding foundational axioms with the aim of being able to prove stronger results.
  \end{subsubordinate}

  \begin{subsubordinate}{On examples of homologically self-dual varieties}
    By (\ref{thm:NormalCat->HomologicallySelfDual}), for a pointed variety of algebras to satisfy the \HSDInline-condition, it suffices that each surjective algebra morphism is a cokernel. For example, in a pointed variety surjective morphisms are cokernels whenever it has a forget functor to the category $\Grps$ of groups.

    Thus, in particular, normal varieties are homologically self-dual. On the other hand, there do exist varieties of algebras which satisfy the \HSDInline-condition but aren't normal---the category $\CMon$ of commutative monoids is one, see (\ref{thm:CMonIsHSD}) and (\ref{exa:NxN-add->N}).
  \end{subsubordinate}
\end{subordinate}

\begin{exercises}

\begin{exercise}[Special property of \eqref{fig:SESs-TotallyNormalSeqMonos}-diagram\ZExactTag]
  \label{exe:SES-TNSM-DiagramProp}
  In the morphism of short exact sequences associated to a totally normal sequence of monomorphisms \eqref{fig:SESs-TotallyNormalSeqMonos}, show that the left hand square is a pullback. Then conclude that  $\alpha=\KerMap{\CoKerMap{\beta\alpha}\Comp \beta}$.
\end{exercise}

\begin{exercise}[Special property of \eqref{fig:SESs-TotallyNormalSeqEpis}-diagram\ZExactTag]
  \label{exe:SES-TNSE-DiagramProp}
  In the morphism of short exact sequences associated to a totally normal sequence of epimorphisms \eqref{fig:SESs-TotallyNormalSeqEpis}, show that the right hand square is a pushout. Then conclude that $\eta=\CoKerMap{\xi\Comp \KerMap{\eta\xi}}$.
\end{exercise}

\begin{exercise}[$\SetsBsd$ is homologically self-dual]
  \label{exe:Set_*HomologicallySelfDual}%
  Show that the category $\SetsBsd$ of pointed sets is homologically self-dual. - Hint: Recall (\ref{thm:NormalEpi/Mono-Set_*}).
\end{exercise}

\begin{exercise}[$\TopsBsd$ is homologically self-dual]
  \label{exe:Top_*IsHSD}
  Show that the category $\TopsBsd$ of based topological spaces is homologically self-dual.
\end{exercise}

With Exercise \ref{exe:Top_*IsHSD} we have a category which is homologically self-dual but is not regular; compare Lemma \ref{thm:TNSM-Mono->NormalMono}.

\begin{exercise}[Pushout property in (\ref{fig:SESs-TotallyNormalSeqEpis})]
  \label{exe:Pushout-SES-TNSE}
  In the situation of diagram (\ref{fig:SESs-TotallyNormalSeqEpis}) assume that $\kappa$ is an epimorphism. Show that the square on the left is a pushout.
\end{exercise}

\begin{exercise}[Pushout recognition under homological self-duality]
  \label{exe:Pushout-Recognition-HSD}
  In a homologically self-dual category, consider a morphism of short exact sequences of the form
  \begin{equation*}
    \xymatrix@R=5ex@C=4em{
    K \ar@{{ |>}->}[r] \ar[d]_{k} &
    X \ar@{-{ >>}}[r] \ar@{-{ >>}}[d]^{\xi}&
    Z \ar@{=}[d] \\
    L \ar@{{ |>}->}[r] &
    Y \ar@{-{ >>}}[r] &
    Z
    }
  \end{equation*}
  where $\xi$ is a normal epimorphism. Show that the map $k$ is a normal epimorphism and that the square on the left is both a pushout and a pullback.
\end{exercise}

\begin{exercise}[Totally normal sequence of mono/epi morphisms - Construction\ZExactTag]
  \label{thm:TotallyNormalMonos/Epis-(Co)KernelsComposite}
  In a pointed category $\Ctgry{X}$, show that an arbitrary sequence $A \XRA{a} B \XRA{b} C$ of morphisms yields
  \begin{thmlist}
    \item a totally normal sequence of monomorphisms $\Ker{a}\NMono \Ker{ba} \NMono A$, and
    \item a totally normal sequence of epimorphisms $C\NEpi \CoKer{ba}\NEpi \CoKer{b}$.
  \end{thmlist}
\end{exercise}

\begin{exercise}[Image factorization of (normal epi)$\Comp$(normal mono)\HSDTag]
  \label{exe:Hofmann}%
  For an antinormal pair $(h,k)$ in a homologically self-dual category the following are equivalent:
  \begin{tfae}
    \item \label{thm:Hofmann-I}%
    $\Ker{h}$ is a subobject of $K$.
    \item \label{thm:Hofmann-II}%
    The composite $h\Comp k$ admits a normal epi / mono image factorization so that the diagram below is a pushout.
    \begin{equation*}
      \xymatrix@R=5ex@C=3em{
      K \ar@{.{ >>}}[r]^-{h'} \ar@{{ |>}->}[d]_-{k} \PushRD{rd} &
      L \ar@{{ >}.>}[d]^-{k'} \\
      X \ar@{-{ >>}}[r]_-{h} &
      Q.
      }
    \end{equation*}
  \end{tfae}
  In this situation $k'$ is a normal map.
\end{exercise}

\begin{exercise}[Normality preserved by certain pullbacks\HSDTag]
  \label{exe:MorphismOfSESs-Properness-HSD}
  In a homologically self-dual category, consider a morphism of short exact sequences in which $c$ is a normal monomorphism. %
  \begin{equation*}
    \xymatrix@R=5ex@C=4em{
    A \ar@{{ |>}->}[r]^{\alpha} \ar[d]_{a} \PullLU{rd} &
    B \ar@{-{ >>}}[r]^{\beta} \ar[d]_{b} &
    C \ar@{{ |>}->}[d]^{c} \\
    X \ar@{{ |>}->}[r]_{\xi} &
    Y \ar@{-{ >>}}[r]_{\eta} &
    Z
    }
  \end{equation*}
  Then the square on the left is a pullback. Moreover,
  \begin{thmlist}
    \item if $b$ is a normal monomorphism, then $a$ is a normal monomorphism;
    \item if $b$ is a normal epimorphism, then $a$ is a normal epimorphism;
    \item if $b$ is a normal map, then $a$ is a normal map.
  \end{thmlist}
\end{exercise}

Dually:

\begin{exercise}[Normal maps preserved by certain pushouts\HSDTag]
  \label{exe:MorphismOfSESs-Properness-HSD-Dual}
  In a homologically self-dual category, consider a morphism of short exact sequences in which $a$ is a normal epimorphism. %
  \begin{equation*}
    \xymatrix@R=5ex@C=4em{
    A \ar@{{ |>}->}[r]^{\alpha} \ar@{-{ >>}}[d]_{a} &
    B \ar@{-{ >>}}[r]^{\beta} \ar[d]_{b} \PushRD{rd} &
    C \ar[d]^{c} \\
    X \ar@{{ |>}->}[r]_{\xi} &
    Y \ar@{-{ >>}}[r]_{\eta} &
    Z
    }
  \end{equation*}
  Then the square on the right is a pushout. Moreover,
  \begin{thmlist}
    \item if $b$ is a normal epimorphism, then $c$ is a normal epimorphism;
    \item if $b$ is a normal monomorphism, then $c$ is a normal monomorphism;
    \item if $b$ is a normal map, then $c$ is a normal map. \NoProof
  \end{thmlist}
\end{exercise}

\begin{exercise}[Pullback of a short exact sequence along a normal monomorphism\HSDTag]
  \label{exe:PullbackSES-AlongNormalMono}
  Compare (\ref{thm:CobaseChange-NormalEpi})\quad In a homologically self-dual category, pulling back a short exact sequence $A\NMono C\NEpi R$ along a normal monomorphism $\gamma\from Q\NMono R$ yields a commutative diagram
  \begin{equation*}
    \xymatrix@R=5ex@C=4em{
    A \ar@{{ |>}->}[r]^-{a} \ar@{=}[d] &
    B \ar[r]^-{f} \ar@{{ |>}->}[d]_{\beta} \PullLU{rd}&
    Q \ar@{{ |>}->}[d]^{\gamma} \\
    A \ar@{{ |>}->}[r]_{\alpha} &
    C \ar@{-{ >>}}[r]_-{c} &
    R
    }
  \end{equation*}
  in which $f$ is normal, and $\CoKer{f}\cong \Ker{\CoKer{\beta}\to \CoKer{\gamma}}$.
\end{exercise}

\begin{exercise}[Pushout of a short exact sequence along a normal epimorphism\HSDTag]
  \label{exe:PushOutSES-AlongNormalEpi}%
  Formulate the dual of Exercise \ref{exe:PullbackSES-AlongNormalMono} and prove it.
\end{exercise}

\begin{exercise}[Normal map essential in Short $5$-Lemma]
  \label{exe:Short-5-HSD-ProperEssential}%
  In the category $\CMon$ of commutative monoids, consider the morphism of short exact sequences.
  \begin{equation*}
    \xymatrix@R=5ex@C=4em{
    \NNr \ar@{{ |>}->}[r]^-{\PrdctMapInto{\IdMap,0}} \ar@{=}[d] &
    \Prdct{\NNr}{\NNr} \ar@{-{ >>}}[r]^-{\PrjctnOnto{2}} \ar[d]^{\PrdctMapInto{+,\PrjctnOnto{2}}}_{b} &
    \NNr \ar@{=}[d] \\
    \NNr \ar@{{ |>}->}[r]_-{\PrdctMapInto{\IdMap,0}} &
    \Prdct{\NNr}{\NNr} \ar@{-{ >>}}[r]_-{\PrjctnOnto{2}} &
    \NNr
    }
  \end{equation*}
  The coordinate maps of $b$ are 'addition' and 'projection onto the second coordinate'. Show that $b$ is injective, not surjective, hence is not an isomorphism.

  Then explain why this morphism of short exact sequences is not a counterexample to the primordial Short $5$-Lemma (\ref{thm:Short-5-Primordial}). Conclude further that the subobject $S\DefEq\SetSlct{(m,n)}{n\geq m}<\Prdct{\NNr}{\NNr}$ is not a normal subobject of $\Prdct{\NNr}{\NNr}$.
\end{exercise}

\begin{exercise}[Linear categories \ANKTag]
  In view of (\ref{thm:CMonIsHSD}) and (\ref{thm:LinearThenCMon}), determine if linear categories are homologically self-dual.
\end{exercise}

\begin{exercise}[Hopf algebras \ANKTag]
  Is the category of Hopf algebras over a field homologically self-dual?
\end{exercise}
\end{exercises}
\section[Preservation of Normal Maps by Dinversion]{Preservation of Normal Maps by Dinversion}
\label{sec:DinversionPreservesNormalMaps}%

In Section \ref{sec:DiExtensions}, we introduced the \DPNInline-condition which requires that \emph{dinversion preserves normal morphisms}. In (\ref{thm:DiExtensionFromDPN}), we explained its role in the construction of di-extensions. Here we introduce di-extensive pushouts and di-extensive pullbacks as conditions which are equivalent to the \DPNInline-condition; see (\ref{def:DiExtensivePushOut/PullBack}).

Then we explain how this development is related to the classical perspective of the `border cases' of the $(\Prdct{3}{3})$-Lemma. Facts surrounding the $(\Prdct{3}{3})$-Lemma are well documented in existing literature. This connection enables us to identify \ZExact\ categories which satisfy the \DPNInline-condition.

Throughout, we take advantage of the fact that di-extensions exhibit notational symmetry across the axis joining the initial and the terminal corners. It cuts the number of cases that need to be considered in half.

\begin{definition}[Di-extensive pushout / pullback\ZExactTag]
  \label{def:SymmetricallyDoubleExtensivePushOut/PullBack}
  \label{def:SymmetricallyDiExtensivePushOut/PullBack}
  \label{def:DiExtensivePushOut/PullBack}%
  A pullback of normal monomorphisms, left below, is \Defn{di-extensive} if, whenever one of $\gamma$ or $c$ is a normal monomorphism, then so is the other. %
  \index{di-extensive!normal pullback}\index{di-extensive!normal pushout}\index{normal pullback!di-extensive}\index{normal pushout!di-extensive}%
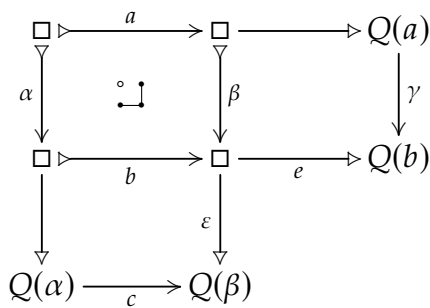
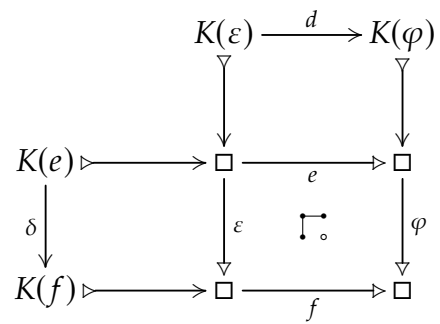
  \begin{equation*}
    \xymatrix@R=3ex@!@C=3em{
    \DiagObj \ar@{{ |>}->}[r]^-{a} \PullLU{rd} \ar@{{ |>}->}[d]_{\alpha} &
    \DiagObj \ar@{-{ >>}}[r] \ar@{{ |>}->}[d]^{\beta} &
    \CoKer{a} \ar[d]^{\gamma} &&
    & \Ker{\varepsilon} \ar[r]^-{d} \ar@{{ |>}->}[d] &
    \Ker{\varphi} \ar@{{ |>}->}[d] \\
    \DiagObj \ar@{{ |>}->}[r]_-{b} \ar@{-{ >>}}[d] &
    \DiagObj \ar@{-{ >>}}[r]_-{e} \ar@{-{ >>}}[d]_{\varepsilon} &
    \CoKer{b} &&
    \Ker{e} \ar@{{ |>}->}[r] \ar[d]_{\delta} &
    \DiagObj \ar@{-{ >>}}[r]_-{e} \ar@{-{ >>}}[d]^{\varepsilon} \PushRD{rd} &
    \DiagObj \ar@{-{ >>}}[d]^{\varphi} \\
    \CoKer{\alpha} \ar[r]_-{c} &
    \CoKer{\beta} &&&
    \Ker{f} \ar@{{ |>}->}[r] &
    \DiagObj \ar@{-{ >>}}[r]_-{f} &
    \DiagObj
    }
  \end{equation*}
  A pushout of normal epimorphisms, right above, is \Defn{di-extensive} if, whenever one of $d$ or $\delta$ is a normal monomorphism, then so is the other.
\end{definition}

\begin{proposition}[Dinversion and di-extensive pushouts / pullbacks\ZExactTag]
  \label{thm:AntiNormalInversion-DoubleExtensivePush/Pull}
  \label{thm:Dinversion-DiExtensivePush/Pull}%
  In a \ZExact\ category the following are equivalent.
  \begin{tfae}
    \item \label{thm:Dinversion-DiExtensivePush/Pull-Pull}%
    Pullbacks of normal monomorphisms are di-extensive.
    \item \label{thm:Dinversion-DiExtensivePush/Pull-DPN->N}%
    Dinversion preserves normal maps.
    \item \label{thm:Dinversion-DiExtensivePush/Pull-Push}%
    Pushouts of normal epimorphisms are di-extensive.
  \end{tfae}
\end{proposition}
\begin{proof}
  The key to observe that, whenever an antinormal composite $\varepsilon\mu$ or $\pi\kappa$ is a normal map, then (\ref{thm:PullbackRecognition-KernelSide-1}) and (\ref{thm:PushoutRecognize-Categorical}) yield pullback and pushout squares as indicated.
  
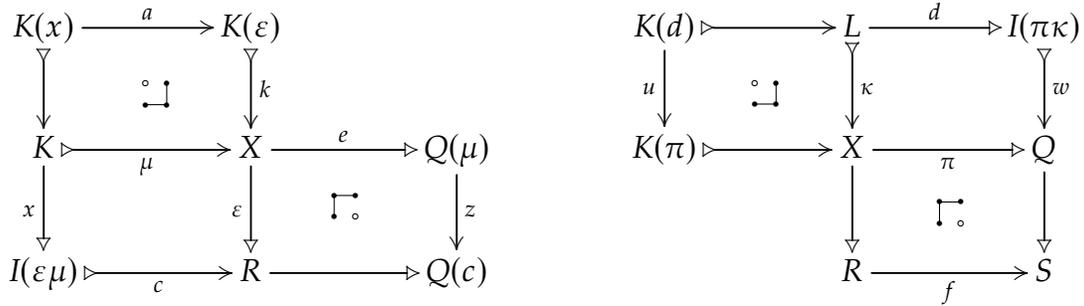
\begin{equation*}
    \xymatrix@R=5ex@C=4em{
    \Ker{x} \ar[r]^-{a} \ar@{{ |>}->}[d] \PullLU{rd} &
    \Ker{\varepsilon} \ar@{{ |>}->}[d]^{k} &&
    \Ker{d} \ar@{{ |>}->}[r] \ar[d]_{u} \PullLU{rd} &
    L \ar@{-{ >>}}[r]^-{d} \ar@{{ |>}->}[d]^{\kappa} &
    \Img{\pi\kappa} \ar@{{ |>}->}[d]^{w} \\
    K \ar@{{ |>}->}[r]_-{\mu} \ar@{-{ >>}}[d]_{x} &
    X \ar@{-{ >>}}[d]_{\varepsilon} \ar@{-{ >>}}[r]^-{e} \PushRD{rd} &
    \CoKer{\mu} \ar[d]^{z} &
    \Ker{\pi} \ar@{{ |>}->}[r] &
    X \ar@{-{ >>}}[r]_-{\pi} \ar@{-{ >>}}[d] \PushRD{rd} &
    Q \ar@{-{ >>}}[d] \\
    \Img{\varepsilon\mu} \ar@{{ |>}->}[r]_-{c} &
    R \ar@{-{ >>}}[r] &
    \CoKer{c} &&
    R \ar[r]_-{f} &
    S
    }
  \end{equation*}
  So, in the diagram on the left, $a$ is a normal monomorphism and $z$ is a normal epimorphism. Similarly, in the diagram on the right, $u$ is a normal monomorphism, and $f$ is a normal epimorphism.

  (\ref{thm:Dinversion-DiExtensivePush/Pull-Pull}) $\Leftrightarrow$ (\ref{thm:Dinversion-DiExtensivePush/Pull-DPN->N})\quad If $\varepsilon\mu$ is normal, we want to show that $ek$ is normal as well. This is so because, by (\ref{thm:Dinversion-DiExtensivePush/Pull-Pull}), $ek$ factors through $\CoKerMap{a}$ with a normal monomorphism.

  Conversely, if in the diagram on the left in (\ref{def:SymmetricallyDoubleExtensivePushOut/PullBack}) the map $c$ is a normal monomorphism, then we must show that $\gamma$ is a normal monomorphism as well. By assumption $e\beta$ is a normal map. Since $a$ is the pullback of $b=\KerMap{e}$ along $\beta$, we have $a=\KerMap{e\beta}$. So, $\CoKer{a}=\Img{e\beta}$, and so $\gamma$ is a normal monomorphism.

  (\ref{thm:Dinversion-DiExtensivePush/Pull-DPN->N}) $\Leftrightarrow$ (\ref{thm:Dinversion-DiExtensivePush/Pull-Push})\quad The argument is dual to the one for (\ref{thm:Dinversion-DiExtensivePush/Pull-Pull}) $\Leftrightarrow$ (\ref{thm:Dinversion-DiExtensivePush/Pull-DPN->N}).
\end{proof}

\begin{proposition}[(DPN) and border cases of the $(\Prdct{3}{3})$-Lemma\ZExactTag]
  \label{thm:DPN-PreservationNormalMaps-Border(3x3)}
  \label{thm:Dinversion-PreservationNormalMaps-Border(3x3)}%
  In a \ZExact\ category, the following are equivalent. %
  \index{$(\Prdct{3}{3})$-Lemma!border cases}%
  \begin{tfae}
    \item Dinversion preserves normal maps.
    \item\label{thm:DPNPNvs3x3-Right3x3}%
    Every $(\Prdct{3}{3})$-diagram in which all rows, together with the left and the middle column are short exact, is a di-extension: the right column is exact as well.
    \item\label{thm:DPNPNvs3x3-Left3x3}%
    Every $(\Prdct{3}{3})$-diagram in which all rows, together with the right and the middle column are short exact, is a di-extension: the left column is exact as well.
    \item\label{thm:DPNPNvs3x3-Top3x3}%
    Every $(\Prdct{3}{3})$-diagram in which all columns, together with  the bottom and the middle rows are short exact, is a di-extension: the top row is exact as well.
    \item\label{thm:DPNPNvs3x3-Bottom3x3}%
    Every $(\Prdct{3}{3})$-diagram in which all columns, together with  the top and the middle rows are short exact, is a di-extension: the bottom row is exact as well.
  \end{tfae}
\end{proposition}
\begin{proof}
  With Lemma \ref{thm:DiExtensionFromAntiNormalPair} we see that (i) is equivalent to each of (ii), (iii), (iv), and (v).
\end{proof}

We close this section by reflecting upon The Classical Snake Lemma (\ref{thm:SnakeLemma-Classical}). Let us say that a \ZExact\ category $\Ctgry{X}$ satisfies the \Defn{Snake Condition} if the conclusion of the Snake Lemma holds in $\Ctgry{X}$. Then we know already:
\begin{ulist}
  \item The Snake Condition implies that dinversion preserves normal maps (Exercise~\ref{exe:SnakeCondImplies3x3}).
  \item Di-exact categories satisfy the Snake Condition.
  \item Sub-di-exact categories satisfy the Snake Condition.
\end{ulist}
We observe that the Snake Condition is self-dual, as are the items listed above. In (\ref{thm:SnakeLemma-Classical-DPN-NormalCat}), we verify the validity of the Snake condition in categories in which dinversion preserves normal maps, and which also satisfy either of the non-self-dual conditions that normal monomorphisms, respectively normal epimorphisms, are closed under composition; (NE/MC).

\begin{example}[\HSDInline\ but not \DPNInline]
  \label{exa:HSD-not-DPN}%
  The category $\CMon$ is homologically self dual. But dinversion does not preserve all of its normal maps.
\end{example}

\begin{subordinate}{}

  \begin{subsubordinate}{Remark on the Middle Cases of the $(\Prdct{3}{3})$-Lemma}
    The nature of the `middle cases' of the $(\Prdct{3}{3})$-Lemma differs significantly from that of its `border cases'. No antinormal pair is present as initial data. Further, an additional requirement is needed: The row/column which is not known to be exact must compose to $\ZeroMap$. As far as we know, the environment of homological categories (\ref{sec:HomologicalCats-Axioms}) provides minimal practical conditions under which the middle cases of the $(\Prdct{3}{3})$-Lemma hold.
  \end{subsubordinate}
\end{subordinate}

\begin{exercises}

\begin{exercise}
  \label{exe:SnakeCondImplies3x3}%
  Use (\ref{thm:Dinversion-PreservationNormalMaps-Border(3x3)}) to show that in any \ZExact\ category in which the Snake Condition holds, dinversion preserves normal maps.
\end{exercise}

\begin{exercise}[\ANKTag]
  \label{exe:CatOfTopAlgs-DPNN}
  Determine for which categories of topological algebras dinversion preserves normal maps.
\end{exercise}
\end{exercises}
\section[Sub-Di-Exact Categories]{Sub-Di-Exact Categories}
\label{sec:Sub-di-exact Categories}%

In this section we discuss in more detail sub-di-exact categories (\ref{term:DiExtensiveConditions}), that is: \ZExact\ categories satisfying  the following two conditions:
\begin{ulist}
  \item Every (co)subnormal map that admits an antinormal decomposition is normal, and
  \item dinversion preserves normal maps.
\end{ulist}

At first sight, this condition may look a bit esoteric. However, it is exactly the condition which is needed in the proof of the Snake Lemma we presented in Section~\ref{sec:Snakes}. Further, it is also relevant in topological varieties. For example, in the category of topological groups, a map which admits an antinormal decomposition need not be a normal map. However, antinormal (co)subnormal maps are normal.

To facilitate identifying categories satisfying the two properties stated above, we reformulate them in terms of preservation of normal epi/monomorphisms under normal image factorization; see (\ref{def:NormalImageFactorizationPreservesNormalEpis/Monos}). For context, recall from  (\ref{thm:ANN<->NIPN}) that antinormal maps are normal if and only if normal images preserve normal maps. Restricting this condition to the preservation of normal epimorphisms as well as normal monomorphisms, yields the following strictly weaker condition:

\begin{definition}[Normal image preservation of normal epis / monos\ZExactTag]
  \label{def:NIC}
  \label{def:NormalImageFactorizationPreservesNormalEpis/Monos}%
  In a pointed category \Defn{normal images preserve normal monomorphisms and normal epimorphisms} if every morphism of normal factorizations, as in the commutative diagram below, satisfies the following conditions: %
  \index{normal!images preserve normal monos / epis}%
  \begin{equation}\label{eq:NIC}%
    \vcenter{
    \xymatrix@R=5ex@C=4em{
    \DiagObj \ar[d]_f \ar@{-{ >>}}[r]^-{e} &
    \DiagObj \ar[d]^{g} \ar@{{ |>}->}[r]^-{m} &
    \DiagObj \ar[d]^h \\
    \DiagObj  \ar@{-{ >>}}[r]_-{p} &
    \DiagObj \ar@{{ |>}->}[r]_-{i} &
    \DiagObj
    }
    }
  \end{equation}
  If the morphisms $f$ and $h$ are normal monomorphisms, then $g$ is a normal monomorphisms, and if $f$ and $h$ are normal epimorphisms, then $g$ is a normal epimorphism.
\end{definition}

\begin{proposition}[Normal images and normality of (co)subnormal maps\ZExactTag]
  \label{thm:NIC}
  For an \ZExact\ category $\Ctgry{X}$ the following conditions are equivalent:
  \begin{tfae}
    \item For each antinormal decomposition of a (co)subnormal map the construction \eqref{eq:Dinversion} yields a di-extension.
    \item $\Ctgry{X}$ is sub-di-exact.
    \item Normal images preserve normal monomorphisms and normal epimorphisms, and dinversion preserves normal maps.
  \end{tfae}
\end{proposition}
\begin{proof}
  (I) $\Leftrightarrow$ (II) follows from Lemma~\ref{thm:DiExtensionFromAntiNormalPair}, because every normal map is both subnormal and cosubnormal.

  (II) $\Rightarrow$ (III)  In the diagram of Definition~\ref{def:NIC}, assume that $f$ and $h$ are normal monomorphisms. The pullback of $i$ and $h$ as in the diagram below induces a factorization $l$ of $g$ through $k$.
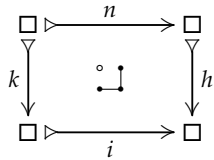
  \begin{equation*}
    \xymatrix@R=5ex@C=4em{
    \DiagObj \PullLU{rd} \ar@{{ |>}->}[d]_{k} \ar@{{ |>}->}[r]^-{n} &
    \DiagObj \ar@{{ |>}->}[d]^h \\
    \DiagObj \ar@{{ |>}->}[r]_-{i} &
    \DiagObj
    }
  \end{equation*}
  We see that $l$ is a normal monomorphism, because so are $m$ and $n$. Now $pf=kle$ is both antinormal and subnormal. So it is normal by assumption. It follows that $g=kl$ is a normal monomorphism. - The case where $f$ and $h$ are normal epimorphisms follows by duality.

  (III) $\Rightarrow$ (II) To see an antinormal subnormal map is normal, we relate it to diagram \eqref{eq:NIC} as follows. Let $pf$ represent the antinormal decomposition. Choose $i$ to be the identity, and let $hme$ represent the subnormal decomposition. Then $g$ is a normal monomorphism because, by the assumption, normal images preserve normal monomorphisms. So, $pf$ is a normal map. - That antinormal cosubnormal maps are normal follows by duality.
\end{proof}

In light of (\ref{thm:NIC}), we may reformulate the di-extensive conditions formulated in (\ref{term:DiExtensiveConditions}) as follows:

\begin{proposition}[Paraphrasing of structural axioms]
  \label{thm:DiExtensiveConditions}%
  Given an \ZExact\ category $\Ctgry{X}$, consider the following conditions:
  \begin{enumerate}[(i)]
    \item \label{Ax:AN0->DiExtension}%
          $\Ctgry{X}$ is \emph{homologically self-dual}: Whenever an antinormal pair composes to the zero map, then it provides initial data for a di-extension.
    \item \label{Ax:DinversionPreservesNormalMaps}%
          \emph{Dinversion preserves normal maps}: whenever an antinormal pair composes to a normal map, then it provides initial data for a di-extension.
    \item \label{Ax:CoSubnormalGenerates}%
          \emph{$\Ctgry{X}$ is sub-di-exact}: Whenever an antinormal pair in $\Ctgry{X}$ composes to a subnormal or a cosubnormal map, then it provides initial data for a di-extension.
    \item \label{Ax:AntinormalGenerates}%
          $\Ctgry{X}$ is \Defn{di-exact}: Antinormal composites are normal maps; every antinormal pair in $\Ctgry{X}$ generates double extension.
    \item \label{Ax:ALLNormal'} $\Ctgry{X}$ is \Defn{Puppe-exact} or \Defn{p-exact}: All maps are normal. %
          \index{Puppe!exact category}\index{p-exact category}%
  \end{enumerate}
  The implications (\ref{Ax:ALLNormal'}) $\implies$ (\ref{Ax:AntinormalGenerates}) $\implies$ (\ref{Ax:CoSubnormalGenerates}) $\implies$ (\ref{Ax:DinversionPreservesNormalMaps}) $\implies$ (\ref{Ax:AN0->DiExtension}) hold. \NoProof
\end{proposition}

\begin{subordinate}{}

  \begin{subsubordinate}{On the notion of p-exact categories}
    In (\ref{thm:DiExtensiveConditions}.\ref{Ax:ALLNormal}) we defined a category to be Puppe-exact if it is \ZExact\ and all morphisms are normal. We should note that, equivalently, a category is Puppe-exact if it has a zero object and all morphisms are normal; see Exercise \ref{exe:PuppeExactCat-EquivalentCharacterization}.
  \end{subsubordinate}

  \begin{subsubordinate}{Examples of sub-di-exact categories}
    We now explain that Axiom~(\ref{thm:DiExtensiveConditions}.\ref{Ax:CoSubnormalGenerates}) holds in the category $\TopGrps$ of topological groups. (We know that (\ref{thm:DiExtensiveConditions}.\ref{Ax:AntinormalGenerates}) doesn't hold: this is related to the fact that $\TopGrps$ is not semiabelian.)

    \begin{proposition}[Normal images of normal epis / monos in top varieties]
      \label{thm:NIC Topological Algebras}%
      If $\Ctgry{V}$ is a pointed variety of algebras in which normal monomorphisms (respectively, normal epimorphisms) are preserved under normal images in the sense of \eqref{def:NormalImageFactorizationPreservesNormalEpis/Monos}, then they are also preserved in the category $\Ctgry{V}(\Tops)$ of $\Ctgry{V}$-algebras in $\Tops$. %
      \index{topological!variety}%
    \end{proposition}
    \begin{proof}
      Suppose $f$ and $h$ are normal monomorphisms in $\Ctgry{V}(\Tops)$. Then we know that the unique map $g$ is a normal monomorphism in $\Ctgry{V}$, and we must show that it is a subspace inclusion. This is so because (a) a composite of two subspace inclusions is a subspace inclusion, and (b) the analogue of (\ref{thm:NormalMono-Props}.iv) holds for subspace inclusions.

      This argument dualizes to normal epimorphisms. Essential here is (a) a composite of two quotient maps is again a quotient map, and (b) the analogue of (\ref{thm:NormalEpi-Props}.iii) holds for quotient maps.
    \end{proof}
  \end{subsubordinate}

\end{subordinate}

\begin{exercises}

\begin{exercise}[Preservation of normal epis / monos in selected categories]
  \label{exe:NormalEpi/MonoPreservationSelectedCats}%
  For each of the following categories determine of normal images preserve normal epimorphisms and normal monomorphisms.
  \begin{thmlist}
    \item $\TopGrps$, the category of topological groups.
    \item Pointed topological spaces.
  \end{thmlist}
\end{exercise}

\begin{exercise}[Topological monoids and di-extension]
  \label{exe:TopMon/DiExtension}
  In the category $\TopMonoids$ of topological monoids do the following: %
  \index{topological!monoids}%
  \begin{thmlist}
    \item Determine if $\TopMonoids$ is homologically self-dual.
    \item Determine if dinversion in $\TopMonoids$ preserves normal maps.
  \end{thmlist}
\end{exercise}

\begin{exercise}[Subtractive varieties and di-extension]
  \label{exe:SubtractionVarieties-DiExtension}%
  For an arbitrary subtractive variety $\Ctgry{V}$ in which $x-y=\ZeroObject$ implies $x=y$ do the following: %
  \index{subtractive!variety}%
  \begin{thmlist}
    \item Determine if dinversion in $\Ctgry{V}$ preserves normal maps.
    \item Determine if $\Ctgry{V}$ is sub-di-exact.
  \end{thmlist}
\end{exercise}

\begin{exercise}[Equivalent characterization of Puppe-exact categories]%
  \label{exe:PuppeExactCat-EquivalentCharacterization}%
  Show that a category is Puppe-exact in the sense of Proposition~\ref{thm:DiExtensiveConditions} if and only if it is a category with a zero object in which every morphism is a normal map.
\end{exercise}
\end{exercises}
\section{Di-Exact Categories}
\label{sec:DiExactCats}%

In a di-exact category every antinormal map is normal; see (\ref{def:DiExactCat}). This means that, whenever a morphism $f$ admits a decomposition $f=em$ in which $m$ is a normal monomorphism and $e$ a normal epimorphism, the morphism $f$ is normal as well. In previous sections we saw that di-exact categories support setting up a foundation for the effective use of homological invariants. From (\ref{thm:ANN<->NIPN}), we already know that this condition is satisfied if and only if normal images preserve normal maps. Here, we present another recognition criterion, for this condition in terms of pullbacks or pushouts with a special property:

\begin{definition}[Di-extensive pullbacks / pushouts\ZExactTag]
  \label{def:DoubleExtensivePull/Push}
  \label{def:DiExtensivePull/Push}%
  A pullback of normal monomorphisms is \Defn{di-extensive} if, as on the left below, the connecting maps between their cokernels are normal monomorphisms as well. %
  \index{di-extensive!pushout}\index{di-extensive!pullback}\index{pushout!di-extensive}\index{pullback!di-extensive}%
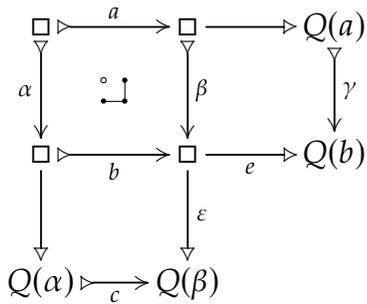
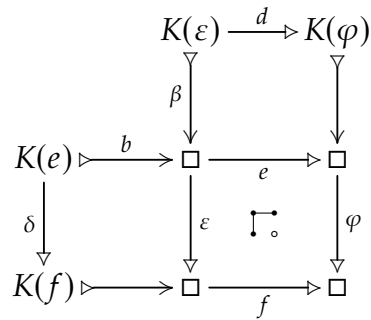
  \begin{equation*}
    \xymatrix@R=3ex@!@C=2em{
    \DiagObj \ar@{{ |>}->}[r]^-{a} \PullLU{rd} \ar@{{ |>}->}[d]_{\alpha} &
    \DiagObj \ar@{-{ >>}}[r] \ar@{{ |>}->}[d]^{\beta} &
    \CoKer{a} \ar@{{ |>}->}[d]^{\gamma} &&
    & \Ker{\varepsilon} \ar@{-{ >>}}[r]^-{d} \ar@{{ |>}->}[d]_{\beta} &
    \Ker{\varphi} \ar@{{ |>}->}[d] \\
    \DiagObj \ar@{{ |>}->}[r]_-{b} \ar@{-{ >>}}[d] &
    \DiagObj \ar@{-{ >>}}[r]_-{e} \ar@{-{ >>}}[d]^-{\varepsilon} &
    \CoKer{b} &&
    \Ker{e} \ar@{{ |>}->}[r]^-{b} \ar@{-{ >>}}[d]_{\delta} &
    \DiagObj \ar@{-{ >>}}[r]_-{e} \ar@{-{ >>}}[d]^{\varepsilon} \PushRD{rd} &
    \DiagObj \ar@{-{ >>}}[d]^{\varphi} \\
    \CoKer{\alpha} \ar@{{ |>}->}[r]_-{c} &
    \CoKer{\beta} &&&
    \Ker{f} \ar@{{ |>}->}[r] &
    \DiagObj \ar@{-{ >>}}[r]_-{f} &
    \DiagObj
    }
  \end{equation*}
  A pushout of normal epimorphisms is di-extensive if, as on the right above, the connecting maps between their kernels are normal epimorphisms as well.
\end{definition}

Here is the motivation for introducing di-extensive pullbacks, respectively di-extensive pushouts.

\begin{proposition}[Di-exact categories and di-extensive pullback / pushouts\ZExactTag]
  \label{thm:DPN-DiExtensivePull/Push}%
  In a \ZExact\ category, the following conditions are equivalent:
  \begin{tfae}
    \item \label{thm:DPN-DiExtensivePull/Push-Pull}%
    Pullbacks of normal monomorphisms are di-extensive.
    \item \label{thm:DPN-DiExtensivePull/Push-ANN}%
    Antinormal composites are normal maps.
    \item \label{thm:DPN-DiExtensivePull/Push-Push}%
    Pushouts of normal epimorphisms are di-extensive.
  \end{tfae}
\end{proposition}
\begin{proof}
  (\ref{thm:DPN-DiExtensivePull/Push-Pull}) $\Leftrightarrow$ (\ref{thm:DPN-DiExtensivePull/Push-Push})\quad With (\ref{thm:DiExtensionFromAntiNormalPair}) we see that the di-extensive pullback on the left above yields a di-extension whose terminal vertex is the pushout of $e$ and $\varepsilon$. Conversely, the di-extensive pushout on the right above yields a di-extension whose initial vertex is the pullback of $b$ and $\beta$.

  Further, with (\ref{thm:DiExtensionFromAN->N.Pairs}), we see that (\ref{thm:DPN-DiExtensivePull/Push-Pull}) $\Leftrightarrow$  (\ref{thm:DPN-DiExtensivePull/Push-ANN}) $\Leftrightarrow$ (\ref{thm:DPN-DiExtensivePull/Push-Push}).
\end{proof}

\begin{corollary}[Normal mono/epi in $\SESCat{X}$\DExTag]
  \label{thm:NormalEpi/Mono-SES(X)-SACat}
  \label{thm:NormalEpi/Mono-SES(X)-DiExact}%
  In the category $\SESCat{X}$ of a di-exact category $\SACtgry{X}$ the following hold for a morphism $F=(\mu,\xi,\eta)$ of short exact sequences:
  \begin{thmlist}
    \item $F$ is a normal monomorphism if and only if $\mu,\xi,\eta$ are normal monomorphisms in $\Ctgry{X}$. This happens if and only if  $(\mu,\xi)$ belong to a pullback square of normal monomorphisms.
    \item $F$ is a normal epimorphism if and only if $\mu,\xi,\eta$ are normal epimorphisms in $\Ctgry{X}$. This happens if and only if  $(\xi,\eta)$ belong to a pushout square of normal epimorphisms.
  \end{thmlist}
\end{corollary}
\begin{proof}
  (i)\quad $F$ is a normal monomorphism in $\SESCat{X}$ if and only if $F$ followed by its cokernel is a short exact sequence, hence a di-extension in $\SESCat{X}$. With the analysis related to (\ref{eq:Dinversion}), we see that this happens if and only if $\mu,\xi,\eta$ are normal monomorphisms in $\Ctgry{X}$. By (\ref{thm:DPN-DiExtensivePull/Push}),  this happens if and only if  $(\mu,\xi)$ belong to a pullback square of normal monomorphisms.

  (ii)\quad The proof is dual to that of (i).
\end{proof}

This means that in a di-exact category, the category $\DExCat{X}$ of di-extensions in $\Ctgry{X}$ is equivalent to both, the category of pushouts of normal epimorphisms in $\Ctgry{X}$, and the category of pullbacks of normal monomorphisms in $\Ctgry{X}$.

\begin{corollary}[Comparison $\DExCat{X}$ - $\NMonoCat{\NMonoCat{X}}$ - $\NEpiCat{\NEpiCat{X}}$ \DExTag]
  \label{thm:DEx(X)-NM(NM(X))-NE(NE(X))}%
  Given a di-exact category $\Ctgry{X}$, the following categories associated to $\Ctgry{X}$ are equivalent: %
  \index{di-exact category!equivalent siblings}%
  \begin{thmlist}
    \item the category $\DExCat{X}$ of di-extensions in $\Ctgry{X}$;
    \item the category $\NMonoCat{\NMonoCat{X}}$ of normal monomorphisms in the category of normal monomorphisms in $\Ctgry{X}$;
    \item the category $\NEpiCat{\NEpiCat{X}}$ of normal epimorphisms in the category of normal epimorphisms in $\Ctgry{X}$.
  \end{thmlist}
\end{corollary}
\begin{proof}
  The functors $m,e\from \Ord{1}^2\to \Ord{2}^2$ induce restrictions of $M\from \DExCat{X}\to \NMonoCat{\NMonoCat{X}}$ and $E\from \DExCat{X}\to \NEpiCat{\NEpiCat{X}}$. An inverse to each comes from (\ref{thm:KernelsInNMono(X)-CoKernelsInNEpi(X)}), which says:
  \begin{ulist}
    \item a normal monomorphism in $\NMonoCat{X}$ is a pullback diagram of normal monomorphisms in $\Ctgry{X}$; and
    \item a normal epimorphism in $\NEpiCat{X}$ is a pushout diagram of normal epis in $\Ctgry{X}$.
  \end{ulist}
  With (\ref{thm:DPN-DiExtensivePull/Push}) the claim follows.
\end{proof}

Since $\NMonoCat{X}\simeq\SESCat{X}\simeq\NEpiCat{X}$, the category of short exact sequences in $\Ctgry{X}$ provides yet another viewpoint:

\begin{proposition}[Antinormal composites and short exact sequences in $\SESCat{X}$\ZExactTag]
  \label{thm:ANN<->PointwiseSES}%
  In a \ZExact\ category $\Ctgry{X}$ the following conditions are equivalent:
  \begin{tfae}
    \item \label{thm:ANN<->PointwiseSES-DiExact}%
    $\Ctgry{X}$ is di-exact.
    \item \label{thm:ANN<->PointwiseSES-SESInSES(X)}%
    Every short exact sequence in $\SESCat{X}$  is a di-extension in $\Ctgry{X}$.
    \item \label{thm:ANN<->PointwiseSES-(Co)Ker(Normal)}%
    Cokernels and kernels of normal maps in any of the categories $\SESCat{X}$, $\NEpiCat{X}$, and $\NMonoCat{X}$ are computed object-wise.
  \end{tfae}
\end{proposition}
\begin{proof}
  (\ref{thm:ANN<->PointwiseSES-DiExact})  $\Rightarrow$ (\ref{thm:ANN<->PointwiseSES-SESInSES(X)})\quad By (\ref{thm:(Co)Kernels-SES(C)}), a short exact sequence in $\SESCat{X}$ is given by a commutative ($\Prdct{3}{3})$-diagram with the indicated features:
  \begin{equation*}
    \xymatrix@R=5ex@C=4em{
    \DiagObj \ar@{{ |>}->}[r]^-{a} \PullLU{rd} \ar@{{ |>}->}[d]_{u} &
    \DiagObj  \ar@{{ |>}->}[d]^{v} \ar[r]^-{d} &
    \DiagObj  \ar@{{ |>}->}[d]^{w} \\
    \DiagObj \ar@{{ |>}->}[r]_-{b}  \ar@{-{ >>}}[d]_{x} &
    \DiagObj \ar@{-{ >>}}[r]^-{e} \ar@{-{ >>}}[d]_{y} \PushRD{rd} &
    \DiagObj \ar@{-{ >>}}[d]^{z} \\
    \DiagObj \ar[r]_-{c} &
    \DiagObj \ar@{-{ >>}}[r]_-{f} &
    \DiagObj
    }
  \end{equation*}
  If $\Ctgry{X}$ is di-exact, then $ev$ is normal. Since $z=\CoKer{ev}$, $w=\KerMap{z}$ is the normal mono part of the normal factorization of $ev$. Thus $d=\CoKer{\Ker{ev}}=\CoKer{a}$ is a normal epimorphism and the top row is short exact. Similar reasoning shows that $c=\KerMap{f}$ renders the bottom sequence short exact. Thus the given short exact sequence in $\SESCat{X}$ is a di-extension.

  (\ref{thm:ANN<->PointwiseSES-DiExact}) $\Rightarrow $ (\ref{thm:ANN<->PointwiseSES-(Co)Ker(Normal)})\quad A normal map $F=(\alpha,\beta,\gamma)$ in $\SESCat{X}$ is a normal epimorphism $E=(\varepsilon_{1},\varepsilon_{2},\varepsilon_{3})$, followed by a normal monomorphism $M=(\mu_{1},\mu_{2},\mu_{3})$. If $\Ctgry{X}$ is di-exact, then we have $\KerMap{F}=\KerMap{E}$, which is computed object-wise by (\ref{thm:NormalEpi/Mono-SES(X)-DiExact}). Similarly, $\CoKerMap{F}=\CoKerMap{M}$ which is also computed object-wise by (\ref{thm:NormalEpi/Mono-SES(X)-DiExact}). This implies the claim for $\NEpiCat{X}$ and for $\NMonoCat{X}$ because the equivalences in Section \ref{sec:CatSESs-I} preserve and reflect limits and colimits.

  (\ref{thm:ANN<->PointwiseSES-(Co)Ker(Normal)}) $\Rightarrow $ (\ref{thm:ANN<->PointwiseSES-SESInSES(X)})\quad Without loss of generality, we argue based on the assumption that cokernels and kernels in $\SESCat{X}$ are computed object-wise. Given a short exact sequence in $\SESCat{X}$, as in the diagram above, $(d,e,f)=\CoKerMap{a,b,c}$ is computed by taking object-wise cokernels. Similarly, $(a,b,c)=\KerMap{d,e,f}$ is computed by taking object-wise kernels. Thus the short exact sequence is a di-extension.

  (\ref{thm:ANN<->PointwiseSES-SESInSES(X)}) $\Rightarrow$ (\ref{thm:ANN<->PointwiseSES-DiExact})\quad Given an antinormal composite $yb$, we show that it is a normal map using the diagram above: Let $v\DefEq \KerMap{y}$, then construct the top left square as a pullback. Put $x\DefEq \CoKerMap{u}$, and let $c$ denote the unique factorization. Then the morphism $(a,b)\from u\to v$ is a normal monomorphism in $\NMonoCat{X}$. Thus, the morphism $(a,b,c)$ of short exact sequences is a normal monomorphism as well. Combined with its cokernel, it is a short exact sequence in $\SESCat{X}$. By assumption, this is a di-extension. So, $c$ is a normal monomorphism, and $cx$ is the normal factorization of $yb$.
\end{proof}

\begin{corollary}[Normal epis in $\NEpiCat{X}$ and $\NMonoCat{X}$ for di-exact $\SACtgry{X}$\DExTag]
  \label{thm:CoKernelsInNEpi(X)-DiExact}%
  In a di-exact category $\Ctgry{X}$ the following hold:
  \begin{enumerate}[(i)]
    \item In $\NEpiCat{X}$ a map $(a,b)$ from a normal epimorphism $x$ to another such $y$ is a normal epimorphism in $\NEpiCat{X}$ if and only if $(a,b)$ and $(x,y)$ form a pushout square of normal epimorphisms in $\Ctgry{X}$.
    \item In $\NMonoCat{X}$ a map $(r,s)$ from a normal monomorphism $u$ to another such $v$ is a normal epimorphism in $\NMonoCat{X}$ if and only if  both $r$ and $s$ are normal epimorphisms in $\Ctgry{X}$. \NoProof
  \end{enumerate}
\end{corollary}

Dually:

\begin{corollary}[Normal monos in $\NEpiCat{X}$ and $\NMonoCat{X}$ for di-exact $\SACtgry{X}$\DExTag]
  \label{thm:NormalMonosInNEpi(X)-Semiabelian}
  In a di-exact category $\Ctgry{X}$ the following hold:
  \begin{enumerate}[(i)]
    \item In $\NMonoCat{X}$ a map $(a,b)$ from a normal monomorphism $x$ to another such $y$ is a normal monomorphism in $\NMonoCat{X}$ if and only if $(a,b)$ and $(x,y)$ form a pullback square of normal monomorphisms in $\Ctgry{X}$.
    \item In $\NEpiCat{X}$ a map $(r,s)$ from a normal epimorphism $u$ to another such $v$ is a normal monomorphism $\NEpiCat{X}$ if and only if  both $r$ and $s$ are normal epimorphisms in $\Ctgry{X}$.
  \end{enumerate}
\end{corollary}

\begin{exercises}

\begin{exercise}[Pushout of normal epimorphisms]
  \label{exe:PushoutNormalEpis-NEpicKernelMap}%
  In a \ZExact\ category, consider the morphism of short exact sequences below.
  
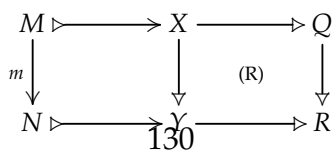
\begin{equation*}
    \xymatrix@R=5ex@C=4em{
    M \ar@{{ |>}->}[r] \ar[d]_{m} &
    X \ar@{-{ >>}}[r] \ar@{-{ >>}}[d] \ar@{}[rd]|-{\text{(R)}}&
    Q \ar@{-{ >>}}[d] \\
    N \ar@{{ |>}->}[r] &
    Y \ar@{-{ >>}}[r] &
    R
    }
  \end{equation*}
  Assume that all maps in square (R) are normal epimorphisms, and that (R) is a pushout. Is the map $m$ on the left always a normal epimorphism? - Also formulate the dual of this question, and try to answer it.
\end{exercise}

\begin{exercise}[Recursiveness of di-extension conditions\ANKTag]
  \label{exe:DiExtensionConditionsRecursive?}%
  In \ZExact\ category $\Ctgry{X}$, do the following:
  \begin{thmlist}
    \item If $\Ctgry{X}$ is homologically self-dual, determine if $\SESCat{X}$ is homologically self-dual.
    \item If dinversion in $\Ctgry{X}$ preserves normal maps, determine if it does so in $\SESCat{X}$.
    \item If $\Ctgry{X}$ is di-exact, determine if $\SESCatn{X}$ is di-exact.
  \end{thmlist}
  We suspect that the answer is `no' at least for item (iii).
\end{exercise}
\end{exercises}
\section[Normal Pushouts and Normal Pullbacks]{Normal Pushouts and Normal Pullbacks}
\label{sec:NormalPushouts/Pullbacks}%

We introduce the concept of normal pushout and explain how it refines the notion of di-extensive pushout from Section \ref{sec:DiExactCats}. Dually, the notion of normal pullback refines the notion of di-extensive pullback. - Starting from a commuting square (S), we construct the pullback $P\rightrightarrows Y$ and the pushout $U\rightrightarrows Q$.
\begin{equation*}
  \xymatrix@R=5ex@C=3em{
  P \ar[r]^-{\pi_V} \ar[d]_-{\pi_X} \PullLU{rd} &
  V \ar[d]^-v \\
  X \ar[r]_-f &
  Y
  }\qquad
  \xymatrix@R=5ex@C=3em{
  U \ar[r]^-g \ar[d]_-u \ar@{}[dr]|-{\text{(S)}} &
  V \ar[d]^-v \\
  X \ar[r]_-f &
  Y
  }\qquad
  \xymatrix@R=5ex@C=3em{
  U \ar[r]^-g \ar[d]_-u \PushRD{rd} &
  V \ar[d]^-{q_V} \\
  X \ar[r]_-{q_X} &
  Q
  }
\end{equation*}
We then have comparison maps $\PrdctMapInto{u,g}\from U\to P$ and $\SumMapOutOf{f,v}\from Q \to Y$.

\begin{definition}[Seminormal/normal pushout\ZExactTag]
  \label{def:NormalPushout}%
  \index{normal pushout}\index{seminormal!pushout}%
  \index{pushout!normal}\index{pushout!seminormal}%
  The square (S) is called
  \begin{enumerate}
    \item a \Defn{normal pushout}, if $u$, $v$, $f$, $g$ and $\PrdctMapInto{u,g}$ are normal epimorphisms; %
          \index{normal!pushout}\index{pushout!normal}%
    \item a \Defn{seminormal pushout}, if $\PrdctMapInto{u,g}$ and at least one of the pairs ($f$ and $g$) or ($u$ and $v$) are normal epimorphisms. %
          \index{seminormal!pushout}\index{pushout!seminormal}
  \end{enumerate}
\end{definition}

Dually:

\begin{definition}[Seminormal/normal pullback]
  \label{def:NormalPullBack}%
  \index{normal pullback}\index{seminormal!pullback}%
  \index{pullback!normal}\index{pullback!seminormal}%
  The square (S) is called
  \begin{enumerate}
    \item a \Defn{normal pullback}, if $u$, $v$, $f$, $g$ and $\SumMapOutOf{f,v}$ are normal monomorphisms; %
          \index{normal!pushout}\index{pushout!normal}%
    \item a \Defn{seminormal pullback}, if $\SumMapOutOf{f,v}$ and at least one of the pairs ($f$ and $g$) or ($u$ and $v$) are normal monomorphisms. %
          \index{seminormal!pushout}\index{pushout!seminormal}
  \end{enumerate}
\end{definition}

\begin{proposition}[Homological self-duality / seminormal pushout property\ZExactTag]
  \label{thm:HDS<->SemiNormalPushoutProp}%
  A \ZExact\ category $\Ctgry{X}$ is homologically self-dual if and only if whenever in a morphism of short exact sequences %
  \index{seminormal!pushout and HSD}%
  \begin{equation*}
    \xymatrix@R=5ex@C=4em{
    K \ar@{{ |>}->}[r] \ar[d]_{\kappa} &
    X \ar@{-{ >>}}[r] \ar[d]_-{\xi} \ar@{}[rd]|-{\text{(R)}}&
    Q \ar[d]^-{\rho} \\
    L \ar@{{ |>}->}[r] &
    Y \ar@{-{ >>}}[r] &
    R
    }
  \end{equation*}
  the square (R) is a seminormal pushout, then $\kappa$ is a normal epimorphism.
\end{proposition}
\begin{proof}
  Suppose $\Ctgry{X}$ is homologically self-dual. In the commutative diagram below, the bottom right square is constructed as the pullback of $q$ along $v$.
  \begin{equation*}
    \xymatrix@R=5ex@C=3em{
    K \ar@{{ |>}->}[rr]^-{k} \ar[dd]_{\kappa} \ar[dr]_{\kappa} &&
    X \ar@{-{ >>}}[rr] \ar[dd]^(0.35){\xi} \ar@{-{ >>}}[rd]^{\hat{\xi}} &&
    V \ar[dd]_(.3){\rho} \ar@{=}[rd] \\
    & L \ar@{{ |>}->}[rr]|\hole^(.3){\hat{k}} &&
    P \ar[ld] \ar[rr]|\hole_(.25){\bar{r}} \ar[ld]^(0.4){\bar{\rho}} &&
    V \ar[dl]^{\rho} \\
    L \ar@{{ |>}->}[rr]_-{l} \ar@{=}[ru] &&
    X \ar@{-{ >>}}[rr]_-{r} &&
    Y
    }
  \end{equation*}
  If the square \ (R) is a seminormal pushout, then the comparison map $\hat{\xi}$ is a normal epimorphism. So, $\bar{r}$ is a normal epimorphism by (\ref{thm:NormalEpi-Props}.\ref{thm:NormalEpiComposite}). Thus the middle row is short exact. So, $\kappa$ is a normal epimorphism by the recognition criterion (\ref{thm:HomologicalSelfDuality-Recognize-II}) for homologically self-dual categories.

  Conversely, suppose $\rho$ is an isomorphism and $\xi$ a normal epimorphism. According to the recognition criterion (\ref{thm:HomologicalSelfDuality-Recognize-II}) for homologically self-dual categories we need to show that $\kappa$ is a normal epimorphism. To see this, construct the square $P\rightrightarrows Y$ as a pullback. So, the isomorphism $\rho$ pulls back to the isomorphism $\bar{\rho}$. Then $\hat{\xi}=\bar{\rho}^{-1}\xi$ is a normal epimorphism. This means that the square $X\rightrightarrows R$ is a normal pushout. By hypothesis, $\kappa$ is a normal epimorphism.
\end{proof}

We need to confirm that a seminormal pushout actually has the universal property expected of a pushout.

\begin{corollary}[Seminormal pushout $\implies$ pushout\HSDTag]
  \label{thm:SemiNormalPushoutIsPushout}%
  In a homologically self-dual category, every seminormal pushout is a pushout.
\end{corollary}
\begin{proof}
  Given a seminormal pushout square, consider the kernels of the pair of opposing normal epimorphisms, as in (\ref{thm:HDS<->SemiNormalPushoutProp}). Then $k$ is a (normal) epimorphism. So, the square (R) is a pushout by (\ref{thm:PushoutRecognize-Categorical}).
\end{proof}

\begin{corollary}[Normal pushout is di-extensive\HSDTag]
  \label{thm:NormalPushout->DoubleExtensive}
  \label{thm:NormalPushout->DiExtensive}%
  In a homologically self-dual category every normal pushout square (R) is a di-extensive pushout of normal epimorphisms.
\end{corollary}
\begin{proof}
  The square (R) is a pushout by (\ref{thm:SemiNormalPushoutIsPushout}). Moreover, the induced maps between kernels of the normal epimorphisms are normal epimorphisms by (\ref{thm:HDS<->SemiNormalPushoutProp}). So, (R) is a di-extensive pushout of normal epimorphisms.
\end{proof}

Dually:

\begin{proposition}[Homological self-duality / seminormal pullback property\ZExactTag]
  \label{thm:HDS<->SemiNormalPullBackProp}%
  A \ZExact\ category $\EuScript{X}$ is homologically self-dual if and only if, whenever in a morphism of short exact sequences the square (L) is a seminormal pullback,%
  \index{seminormal!pushout and HSD}%
  \begin{equation*}
    \xymatrix@R=5ex@C=4em{
    K \ar@{{ |>}->}[r] \ar[d]_{\kappa} \ar@{}[rd]|-{\text{(L)}} &
    X \ar@{-{ >>}}[r] \ar[d]_-{\xi} &
    Q \ar[d]^-{\rho} \\
    L \ar@{{ |>}->}[r] &
    Y \ar@{-{ >>}}[r] &
    R
    }
  \end{equation*}
  then $\rho$ is a normal monomorphism. \NoProof
\end{proposition}

\begin{corollary}[Seminormal pullback $\implies$ pullback\HSDTag]
  \label{thm:SemiNormalPullBackIsPullBack}%
  In a homologically self-dual category, every seminormal pullback is a pullback. %
  \index{pullback}\index{seminormal!pullback}%
  \NoProof
\end{corollary}

\begin{corollary}[Normal pullback is di-extensive\HSDTag]
  \label{thm:NormalPullback->DoubleExtensive}
  \label{thm:NormalPullback->DiExtensive}%
  In a homologically self-dual category every normal pullback square (L) is a di-extensive pullback of normal monomorphisms. %
  \index{di-extensive!pullback}%
  \NoProof
\end{corollary}

In a normal category pullbacks preserve normal epimorphisms. Because of this property, we obtain a convenient criterion for identifying seminormal pushouts. As a consequence, a pushout of normal epimorphisms is double extensive if and only if it is a normal pushout. Thus, in a normal category antinormal inversion preserves normal maps. Therefore, much of what we developed earlier in this Chapter is satisfied in normal categories.

\begin{example}[A non-seminormal pushout]
  Pick any two non-zero objects $A$ and $B$ in the category $\Grps$ of groups. If $\InclsnOf{A}\from A\to A+B$ is the canonical inclusion into the coproduct, then we obtain the morphism of short exact sequences whose right hand square is a pushout:
  \begin{equation*}
    \xymatrix@R=5ex@C=4em{
    A \ar@{=}[r] \ar[d]_{i} &
    A \ar@{-{ >>}}[r] \ar[d]_{\InclsnOf{A}} \PushRD{rd} &
    0  \ar[d] \\
    \bar{A} \ar@{{ |>}->}[r] &
    A+B \ar@{-{ >>}}[r]_-{\SumMapOutOf{\ZeroMap,\IdMapOn{B}}} &
    B
    }
  \end{equation*}
  The comparison map of kernels is the inclusion of $A$ into its normal closure in $A+B$, which is not a normal epimorphism. So, the pushout on the right can not be a seminormal pushout by (\ref{thm:HDS<->SemiNormalPushoutProp}).
\end{example}

\begin{subordinate}{}

  \begin{subsubordinate}{Origin of `normal pushout'}
    \label{rem:RegularPushout-ConceptVariations}%
    Normal pushouts were preceded by regular pushouts, introduced in~\cite{DBourn2003}. The more general variant of a semiregular pushout appears in \cite[3.5]{VdLinden:Simp}. We want to work with exact sequences and with homological invariants. We introduced those concepts in \ZExact\ categories, and we saw that normal epimorphisms and normal monomorphisms are the central concepts for our purposes. Adapting the notion of (semi)regular pushout resulted in (semi)normal pushout. - In a normal category, both concepts coincide.
  \end{subsubordinate}
\end{subordinate}

\begin{exercises}

\begin{exercise}[Normal pullbacks in every \ZExact\ category\ZExactTag]
  \label{exe:NormalPullback-Ptd}%
  If $\kappa\from K\NMono X$ is a normal monomorphism, show that the pullback of $\kappa$ along itself is a normal pullback. %
  \index{normal!pullback}%
\end{exercise}

\begin{exercise}[Pullback of normal monomorphisms in $\SetsBsd$]
  \label{exe:PullbackNormalMonosSet_*}%
  Show that, in the category $\SetsBsd$ of pointed sets, every pullback of normal monomorphisms is a normal pullback.
\end{exercise}

\begin{exercise}[Pullback of normal monomorphisms of groups need not be normal]
  \label{exe:NormalPullbackNormalMono-Grp}
  Show that, in the category $\Grps$ of groups, the pullback of two normal monomorphisms need not be a normal pullback. %
  \index{normal!pullback in $\Grps$}
\end{exercise}
\end{exercises}
\newpage
\section{Higher Extensions}
\label{sec:HigherExtensions}

To define what an $n$-fold extension in a \ZExact\ category $\Ctgry{X}$ is, we adopt the following notational conventions.

For $n\geq 0$, let $\Ord{n}$ denote the category $0\to \cdots \to n$. In particular, $\Ord{2} = 0\to 1\to 2$. For $k_{1},\dots ,k_{n-1},t\in \Set{0,1,2}$, define the $n$-tuple $d_{i}^{n}(k_{1},\dots ,k_{n-1};t)$ by
\begin{equation*}
  d_{i}^{n}(k_{1},\dots ,k_{n-1};t) \DefEq
  \begin{cases}
    (t,k_{1},\dots ,k_{n-1})                       & \text{if}\ \ i=1    \\
    (k_{1},\dots , k_{i-1},t,k_{i},\dots ,k_{n-1}) & \text{if}\ \ 1< i<n \\
    (k_{1},\dots ,k_{n-1},t)                       & \text{if}\ \ k=n
  \end{cases}
\end{equation*}
In a \ZExact\ category $\Ctgry{X}$ and $n\geq 1$, an \Defn{$n$-fold extension} is a diagram $E$ modeled on $\Ord{2}^{n}$ in which every sequence of the form %
\index{higher extension}\index{$n$-fold extension}%
\begin{equation*}
  E(d^{n}_{i}(k_{1},\dots ,k_{n-1};0)) \longrightarrow E(d^{n}_{i}(k_{1},\dots ,k_{n-1};1)) \longrightarrow E(d^{n}_{i}(k_{1},\dots ,k_{n-1};2))
\end{equation*}
is short exact. The $n$-fold extensions in $\Ctgry{X}$ form the full subcategory $\HExCat{n}{X}$ of the category of $\Ord{2}^{n}$-diagrams in $\Ctgry{X}$. For example, $\HExCat{1}{X}=\SESCat{X}$, and $\HExCat{2}{X} = \DExCat{X}$. %
\index[not]{e!$\HExCat{n}{X}$\IndSep category of $n$-fold extensions in $\Ctgry{X}$}%

Assuming that antinormal maps in $\Ctgry{X}$ are normal (\ref{term:DiExtensiveConditions}), we show how to construct $n$-fold extensions inductively from suitable morphisms of $(n-1)$-extensions. This construction may be re-interpreted as the border cases of a higher order $(\Prdct{3}{3})$-Lemma.

\begin{proposition}[Sources of di-extensions\DExTag]
  \label{thm:DoubleExtensionSource}
  \label{thm:DiExtensionSource}
  In a di-exact category $\Ctgry{X}$, the following hold for a morphism $F=(\mu,\xi,\eta)$ of short exact sequences
  \begin{equation}
    \vcenter{
    \xymatrix@R=5ex@C=4em{
    \varepsilon \ar[d]_{(\mu,\xi,\eta)} &
    M \ar@{{ |>}->}[r]^-m \ar[d]_{\mu} &
    X \ar@{-{ >>}}[r]^-q \ar[d]_{\xi} &
    Q \ar[d]^{\eta} \\
    \varphi &
    N \ar@{{ |>}->}[r]_-n &
    Y \ar@{-{ >>}}[r]_-r &
    R
    }}
  \end{equation}
  \begin{thmlist}
    \item \label{thm:DiExtensionSource-NEpis}%
    If each of $\mu,\xi,\eta$ is a normal epimorphism then the kernels of $\mu,\xi,\eta$ extend $F$ to a di-extension which is functorial with respect to morphisms of such diagrams.
    \item \label{thm:DiExtensionSource-NMonos}%
    If each of $\mu,\xi,\eta$ is a normal monomorphism then the cokernels of $\mu,\xi,\eta$ extend $F$ to a di-extension  which is functorial with respect to morphisms of such diagrams.
  \end{thmlist}
\end{proposition}
\begin{proof}
  Both claims follow using  one of the `border cases' of the $(\Prdct{3}{3})$-Lemma, as in (\ref{thm:DPN-PreservationNormalMaps-Border(3x3)}). We note  in passing that the conclusion of this proposition holds under the weaker hypothesis of `dinversion preserves normal maps', as in (\ref{thm:DiExtensionFromAntiNormalPair}).
\end{proof}

We use (\ref{thm:DiExtensionSource}) inductively to construct $n$-fold extensions from $(n-1)$-extensions. Equivalently, we formulate versions of border cases of the $(3^n)$-lemma for $(n-1)$-fold extensions.

\begin{theorem}[Recursive construction of $n$-fold extensions\DExTag]
  \label{thm:n-FoldExtension-Recursion}%
  For a morphism
  \begin{equation*}
    F=\SetSlct{f_{\mathbf{i}}}{\mathbf{i}=(i_{1},\dots ,i_{n})\in \Set{0,1,2}^{n}}\from X\to Y
  \end{equation*}
  of $n$-fold extensions the following hold:
  \begin{thmlist}
    \item If each $f_{\mathbf{i}}$ is a normal epimorphism, then adjoining the kernels of the $f_{\mathbf{i}}$ yields an $(n+1)$-fold extension.
    \item If each $f_{\mathbf{i}}$ is a normal monomorphism, then adjoining the cokernels of the $f_{\mathbf{i}}$ yields an $(n+1)$-fold extension.
  \end{thmlist}
\end{theorem}
\begin{proof}
  (i)\quad Let $K$ denote the $\Ord{2}^{n}$-diagram constructed from the kernels of the normal epimorphisms $f_{\mathbf{i}}$. It only remains to show that each of the top sequences in the $(\Prdct{3}{3})$-diagram below is short exact:
  \begin{equation*}
    \xymatrix@R=5ex@C=4em{
    K(d^{n}_{i}(k_{1},\dots ,k_{n-1};0)) \ar[r] \ar@{{ |>}->}[d] &
    K(d^{n}_{i}(k_{1},\dots ,k_{n-1};1)) \ar[r] \ar@{{ |>}->}[d] &
    K(d^{n}_{i}(k_{1},\dots ,k_{n-1};2)) \ar@{{ |>}->}[d] \\
    X(d^{n}_{i}(k_{1},\dots ,k_{n-1};0)) \ar@{{ |>}->}[r] \ar@{-{ >>}}[d] &
    X(d^{n}_{i}(k_{1},\dots ,k_{n-1};1)) \ar@{-{ >>}}[r] \ar@{-{ >>}}[d] &
    X(d^{n}_{i}(k_{1},\dots ,k_{n-1};2)) \ar@{-{ >>}}[d] \\
    Y(d^{n}_{i}(k_{1},\dots ,k_{n-1};0)) \ar@{{ |>}->}[r] &
    Y(d^{n}_{i}(k_{1},\dots ,k_{n-1};1)) \ar@{-{ >>}}[r] &
    Y(d^{n}_{i}(k_{1},\dots ,k_{n-1};2))
    }
  \end{equation*}
  This, however, is the case by (\ref{thm:DoubleExtensionSource}). The proof of (ii) is similar.
\end{proof}

\begin{corollary}[$3^{n}$-lemma\DExTag]
  \label{thm:3^nLemma}%
  Let $X$ be a $\Ord{2}^{n}$-diagram, $n\geq 2$, with sequences
  \begin{equation*}
    X(d^{n}_{i}(k_{1},\dots ,k_{n-1};0)) \longrightarrow X(d^{n}_{i}(k_{1},\dots ,k_{n-1};1)) \longrightarrow
    X(d^{n}_{i}(k_{1},\dots ,k_{n-1};2))
  \end{equation*}
  Then the following hold:
  \begin{thmlist}
    \item If all of the above sequences are known to be short exact, except possibly the ones with $i=1$, then $X$ is an $n$-fold extension.
    \item If all of the above sequences are known to be short exact, except possibly the ones with $i=n$, then $X$ is an $n$-fold extension. \NoProof
  \end{thmlist}
\end{corollary}

\begin{proposition}[Normality of antinormal subnormal maps is recursive w.r.t. $\SESCat{X}$]
  \label{thm:NormalityAntiNormalSubNormal-Recursive}
  If antinormal subnormal maps $\Ctgry{X}$ are normal, then this same property holds in each of the equivalent categories $\NMonoCat{X}$, $\SESCat{X}$, and $\NEpiCat{X}$.
\end{proposition}
\begin{proof}
  We prove the claim for $\NMonoCat{X}$. A subnormal antinormal map in $\NMonoCat{X}$ is given by this kind of commutative diagram:
  
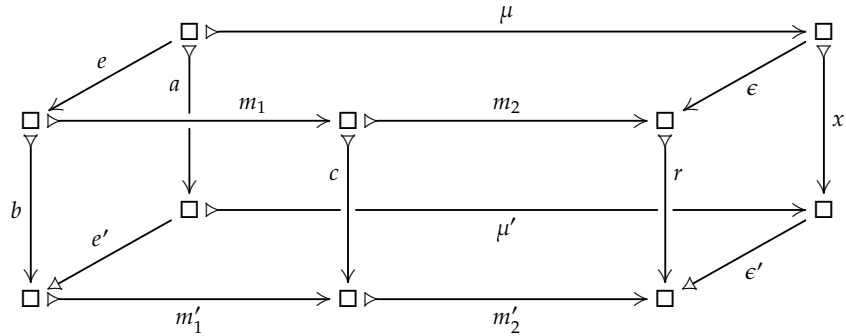
\begin{equation*}
    \xymatrix@!0@C=5em@R=6ex{
    & \DiagObj \ar[ld]_e \ar@{{ |>}->}[rrrr]^-{\mu} \ar@{{ |>}->}[dd]|\hole_(.3){a} &&&&
    \DiagObj \ar@{{ |>}->}[dd]^{x} \ar@{->}[ld]^{\epsilon} \\
    \DiagObj \ar@{{ |>}->}[rr]^(.7){m_1} \ar@{{ |>}->}[dd]_{b} &&
    \DiagObj \ar@{{ |>}->}[rr]^{m_2} \ar@{{ |>}->}[dd]_(0.3){c} &&
    \DiagObj \ar@{{ |>}->}[dd]^(.3){r} \\
    & \DiagObj  \ar@{-{ >>}}[ld]_-{e'}  \ar@{{ |>}->}[rrrr]|(0.25)\hole_(.5){\mu'}|(.75){\hole} && &&
    \DiagObj  \ar@{-{ >>}}[ld]^{\epsilon'}\\
    \DiagObj  \ar@{{ |>}->}[rr]_-{m_1'} &&
    \DiagObj  \ar@{{ |>}->}[rr]_{m_2'}  &&
    \DiagObj
    }
  \end{equation*}
  Consequently, both front squares and the back square are pullbacks of normal monomorphisms in $\Ctgry{X}$. The left and right squares depict normal epimorphisms in $\NMonoCat{X}$. So, $e'$ and $\epsilon'$ are normal epic in $\Ctgry{X}$. By the assumption on $\Ctgry{X}$, the composite $m_2'm_1'$ is a normal monomorphism in $\Ctgry{X}$. Hence so is its pullback $m_2m_1$. It follows that the front square is a pullback of normal monomorphisms in $\Ctgry{X}$, or in other words, a normal monomorphism in $\NMonoCat{X}$.
\end{proof}

\begin{subordinate}{On higher extensions}

  \begin{subsubordinate}{On weakening hypotheses to \DPNTag}
    Reviewing the proof of Proposition~\ref{thm:DoubleExtensionSource}, we see that normal maps in $\EuScript{X}$ are closed under dinversion  (\ref{thm:DiExtensionFromDPN}) if and only if all border cases of the $(\Prdct{3}{3})$-Lemma hold. That is: Suppose in a $(\Prdct{3}{3})$-diagram all rows and columns, except possibly for one of the four border sequences are known to be short exact, then the diagram is a di-extension.
  \end{subsubordinate}

  \begin{subsubordinate}{On viewpoints of higher extensions}
    Every $n$-fold extension may be viewed as a short exact sequence of $(n-1)$-fold extensions in $n$ `orthogonal' ways. Thus, we have two ways of distinguishing a higher extension in the category of short exact sequences of $(n-1)$-fold extensions:
    \begin{ulist}
      \item An $n$-fold extension is a short exact sequence in the category of $(n-1)$-fold extensions in all $n$ orthogonal ways.
      \item An $n$-fold extension is a short exact sequence of $(n-1)$-fold extensions which forgets to a short exact sequence in the category of $(n-1)$-fold diagrams.
    \end{ulist}
    Both points of view complement each other. While the first is appealing because of its fully inductive nature, the second benefits from pointwise limits and colimits in the underlying category.
  \end{subsubordinate}
\end{subordinate}
\chapter[Normal Categories]{Normal Categories}
\label{chap:NormalCategories}

The term `normal category' was first used by Mitchell \cite[I.14]{BMitchell1965-Cats} to refer to a pointed category in which every monomorphism is a normal map. In \cite[Sec.~1]{ZJanelidze-Snake}, Z.~Janelidze repurposed this term to refer to a different concept: a category which is pointed and regular, and such that regular epimorphisms and normal epimorphisms coincide. In a finitely bicomplete context, this amounts to the existence of pullback-stable normal epi/mono factorizations which, in turn, corresponds to the condition that split epimorphisms are normal, while normal epimorphisms are preserved under pullbacks. This latter characterization is what we adopt as a definition: see (\ref{def:NormalCat}).

\bigskip
\begin{center}
  \textbf{Leitfaden for Chapter \ref{chap:NormalCategories}}
\end{center}

\bigskip

\begin{equation*}
  \xymatrix@R=6ex@C=4.5em{
  *+[F-,]{\txt{\sffamily (\ref{sec:NormalCats}) Normal Categories}} \ar@{<->}@<10ex>@/^5ex/[ddd] \ar[d] \\
  *+[F-,]{\txt{\sffamily (\ref{sec:ImageFactorizations}) Image Factorization}} \ar[d] \\
  *+[F-,]{\txt{\sffamily (\ref{sec:SnakeLemma-N+DPN}) Snake Lemma}} \ar[d] \\
  *+[F-,]{\txt{\sffamily (\ref{sec:AlternateNormal}) Alternate Characterizations \\ \sffamily of Normal Categories}}
  }
\end{equation*}
\section{Normal Categories}
\label{sec:NormalCats}%

In \cite[Sec.~1]{ZJanelidze-Snake}, Z.~Janelidze defined a normal category as being a pointed regular category where all regular epis are normal. Here, guided by the foundation developed earlier, we adapt his ideas as follows:

\begin{definition}[Normal category]
  \label{def:NormalCat}%
  A category $\Ctgry{X}$ is called \Defn{normal} if the following structural axioms are satisfied:  %
  \index{normal!category}\index{category!normal}%
  \begin{ulist}
    \item \emph{Zero object}: $\Ctgry{X}$ has a zero object; see (\ref{def:0-Object}).
    \item \emph{Functorially finitely complete}: For every functor $F\from J\to \Ctgry{X}$ whose domain is a finite category $J$, the limit $\LimOf{F}$ exists; see (\ref{sec:Limits-CoLimits}).
    \item \emph{Functorially finitely cocomplete}: For every functor $F\from J\to \Ctgry{X}$ whose domain is a finite category $J$, the colimit $\CoLimOf{F}$ exists; see (\ref{sec:Limits-CoLimits}).
    \item \emph{\PNEInline\ Pullbacks preserve normal epimorphisms}: The pullback $\bar{g}$ of a normal epimorphism $g$ along any morphism $f$ in $\Ctgry{X}$ is a normal epimorphism. %
    \index[acr]{p!\PNEInline\IndSep condition that pullbacks preserve normal epimorphisms}%
    \begin{equation*}
      \xymatrix@R=5ex@C=4em{
      P \ar@{-{ >>}}[d]_-{\bar{g}} \ar[r] \PullLU{rd} &
      Z \ar@{-{ >>}}[d]^-{g} \\
      X \ar[r]_-{f} &
      Y
      }
    \end{equation*}
    \item \emph{\AENInline\ Absolute epimorphisms are normal:} In $\Ctgry{X}$, every absolute epimorphism is a cokernel of its kernel. %
    \index{absolute!epimorphism}\index{epimorphism!absolute}%
    \index[acr]{a!\AENInline\IndSep condition that absolute epimorphisms are normal maps}%
    \index[acr]{n!{\color{Brown} $\EuRoman{N}$}\IndSep normal category}%
  \end{ulist}
\end{definition}

To comment on the structural axioms characterizing a normal category, we start in a finitely bicomplete pointed category. Thus every normal category is automatically \ZExact. Via (\ref{thm:CoKer=NormalEpi=RegEpi=EffectiveEpi}), the \PNEInline-condition is inspired by the fact that it is satisfied in every variety of algebras in which normal and regular epimorphisms coincide.

As to the \AENInline-condition, an immediate consequence is that every sectioned epimorphism determines a special kind of short exact sequence:

\begin{definition}[Split short exact sequence\ZExactTag]
  A \Defn{split short exact sequence} is given by a diagram
  \index{split!short exact sequence}%
  \begin{equation*}
    \xymatrix@R=5ex@C=2em{
    K \ar@{{ |>}->}[rr]^-{k} &&
    X \ar@{-{ >>}}@<.5ex>[rr]^-{q} &&
    Q \ar@<.5ex>[ll]^-{x}
    }
  \end{equation*}
  where $k=\KerMap{q}$, $q=\CoKerMap{k}$ and $qx=\IdMapOn{Q}$.
\end{definition}

\begin{proposition}[Sectioned epimorphism yields split short exact sequence\NTag]
  \label{thm:SplitEpi->SES}
  An absolute epimorphism $q\from X\to Q$, sectioned by $x\from Q\to X$, together with any choice of a morphism $k\from K\to X$ representing $\KerMap{q}$, yields a split short exact sequence. %
  \index{split!short exact sequence from split epi\NTag}
\end{proposition}
\begin{proof}
  Via the \AENInline-condition, this is so because every absolute epimorphism is the cokernel of its kernel.
\end{proof}

Further, we show in the next section that the \AENInline-condition plays a key role in establishing the existence of image factorizations. Here is an example which shows that in general, not every absolute epimorphism is a normal map.

\begin{example}[Addition $\NNr\prdct\NNr \to \NNr$]
  \label{exa:NxN-add->N}
  In the category $\Monoids$ of monoids (or in the category $\CMon$ of commutative monoids), consider the addition morphism together with a chosen section
  \begin{equation*}
    a\from \NNr \prdct \NNr \longrightarrow \NNr,\qquad a(m,n)\DefEq m+n,\qquad s(n)\DefEq (n,0)
  \end{equation*}
  Then $k\DefEq\KerMap{a}=0$. So, $a\neq \CoKerMap{k}$.
\end{example}

We turn to developing properties of normal categories. To begin, we show that the situation of Example~\ref{exa:NxN-add->N} cannot occur: a map is a monomorphism exactly when its kernel vanishes. We prove this in two steps:

\begin{lemma}[Monomorphism if and only if kernel is $\ZeroObject$]
  \label{thm:Proto-Mono<->0-Kernel}%
  In category with zero object assume that every map has a kernel pair, and that the \AENInline-condition holds. Then a morphism $f\from X\to Y$ is a monomorphism if and only if its kernel vanishes. %
  \index{monomorphism!recognition by $\ZeroObject$-kernel}%
\end{lemma}
\begin{proof}
  If $f$ is monic, then $\KerMap{f}=\ZeroMap$ by Proposition~\ref{thm:Ker(Mono)=0}. For the converse, consider the kernel pair diagram of $f$:
  \begin{equation*}
    \xymatrix@R=5ex@C=4em{
    \KrnlPr{f} \ar@<.5ex>[r]^-{\PrjctnOnto{1}} \ar@<-0.5ex>[d]_{\PrjctnOnto{2}} \PullLU{rd} &
    X \ar[d]^{f} \ar@<0.5ex>[l]^-{s} \\
    X \ar[r]_-{f} \ar@<-0.5ex>[u]_{s} \ar[r]_-{f} &
    Y
    }
  \end{equation*}
  Since $\Ker{f}=\ZeroObject$, by Proposition~\ref{thm:Pullback->IsoOfKernels} the kernels of $\PrjctnOnto{1}$ and $\PrjctnOnto{2}$ vanish as well. But these sectioned epimorphisms are normal epimorphisms by assumption. In particular (\ref{thm:Ker(CoKer)-CoKer(Ker)}), the morphism $\PrjctnOnto{1}$ is a cokernel of its kernel $\ZeroObject$.  So, it is an isomorphism. But then $f$ is a monomorphism by (\ref{thm:KernelPair-Monos}).
\end{proof}

\begin{corollary}[Monomorphism if and only if kernel is $\ZeroObject$\NTag]
  \label{thm:Mono<->0-Kernel}
  In a normal category, a morphism $f\from X\to Y$ is a monomorphism if and only if its kernel vanishes.
  \index{monomorphism!recognition by $0$-kernel}
\end{corollary}
\begin{proof}
  One of the structural axioms of a normal category requires sectioned epimorphisms to be normal. Thus the claim follows from  (\ref{thm:Proto-Mono<->0-Kernel}).
\end{proof}

Next, we relate normal categories to the development in Chapter \ref{chap:Di-ExactCats}. We show that every normal category is homologically self-dual and even satisfies the \DPNInline-condition; see Section \ref{sec:ImageFactorizations}. In what follows, we will refer to the notation
\stepcounter{theorem}
\begin{equation}\label{fig:NormalPushoutSES}
  \vcenter{
  \xymatrix@R=5ex@C=4em{
  K \ar@{{ |>}->}[r]^-{k} \ar[d]_{\kappa} \ar@{}[rd]|-{\text{(L)}} &
  X \ar@{-{ >>}}[r]^-{q} \ar[d]_{\xi} \ar@{}[rd]|-{\text{(R)}}&
  Q \ar[d]^{\rho} \\
  L \ar@{{ |>}->}[r]_-{l} &
  X \ar@{-{ >>}}[r]_-{r} &
  R
  }
  }
\end{equation}
for a morphism of short exact sequences.

\begin{lemma}[Normal categories are homologically self-dual\NTag]
  \label{thm:NormalCat->HomologicallySelfDual}%
  Every normal category $\Ctgry{X}$ is homologically self-dual. %
  \index{normal!category - homologically self-dual}\index{HSD!property of normal category}%
\end{lemma}
\begin{proof}
  We use the recognition criterion (\ref{thm:HomologicalSelfDuality-Recognize-II}.\ref{thm:HomologicalSelfDuality-Recognize-TotPull}). In Diagram \eqref{fig:NormalPushoutSES}, consider the case where $\rho$ is an isomorphism, and $\xi$ is a normal epimorphism. We must show that $\kappa$ is a normal epimorphism as well. To see this, note first that the square (L) is a pullback by (\ref{thm:PullbackRecognition-KernelSide-1}). Via the \PNEInline-condition, it pulls the normal epimorphism $\xi$ back to the normal epimorphism $\kappa$. So, $\Ctgry{X}$ is homologically self-dual.
\end{proof}

Among the normal categories, we find certain subtractive varieties:

\begin{definition}[Subtractive variety]
  \label{def:SubtractiveVariety}%
  A pointed variety is said to be \Defn{subtractive} if it has a binary operation `$-$' which satisfies $x-x=0$ and $x-0=x$. %
  \index{subtractive!variety}\index{variety!subtractive}%
\end{definition}

We adopted the terminology `subtractive variety' from \cite{ZJanelidze-Snake}. Alternatively, one may conceptualize such a variety as consisting of left unital magmas satisfying $x^2=1$ for all $x$. For example, every pointed variety whose algebras have an underlying group or loop structure is subtractive. In a group, define  $x-y\DefEq xy^{-1}$. In a loop define $x-y\DefEq x/y$ where $/$ is the division on the right.

The following result presents a sufficient condition for a subtractive variety to be a normal category.

\begin{proposition}[Condition for subtractive variety to be normal category]
  \label{thm:SubtractiveVariety->NormalCat}%
  Any subtractive variety in which $x-y=0$ implies $x=y$ is a normal category.
\end{proposition}
\begin{proof}
  If $f\from X\to Y$ is a regular epimorphism, given as the coequalizer of $r_1$, $r_2\from R\to X$, then $f$ is a cokernel of its kernel $\kappa\from {\Ker{f}\to X}$. If $g\from X\to Z$ satisfies $g\kappa=0$, we claim that $gr_{1}=gr_{2}$. Indeed, for each $x\in R$, we have $f(r_1(x))=f(r_2(x))$ so that $f(r_1(x)-r_2(x))=0$ and $(r_1(x)-r_2(x))\in \Ker{f}$. Then $g(r_1(x)-r_2(x))=0$, and so $g(r_1(x))=g(r_2(x))$ as claimed. This implies that $g$ factors uniquely through the coequalizer~$f$ of $r_1$ and $r_2$. So, $f$ is a normal epimorphism and the \PNEInline-condition is satisfied.

  Finally, in a finitely complete category, absolute epimorphisms are regular. So, what we just proved also implies that the \AENInline-condition is satisfied.
\end{proof}

Thus all varieties of groups are subtractive and, hence, normal. More generally, a pointed variety which has a group structure amongst its operations is subtractive.

In a loop, $x/y=1$ implies $y=1\cdot y=(x/y)\cdot y=x$ so the variety of loops is a normal category as well.
\section[Image Factorization]{Image Factorization}%
\label{sec:ImageFactorizations}%

In Section \ref{sec:NormalDecompositionsFactorizations}, we showed that, in a \ZExact\ category, every morphism $f$ may be factored in a universal way as $f=ve$, with $v$ a normal epimorphism. We called this the normal-epi factorization of $f$, and know that it may be constructed thus: %
\stepcounter{theorem}%
\begin{equation}\label{eq:ImageFactorization}
  \vcenter{
  \xymatrix@R=5ex@C=4em{
  \Ker{f} \ar@{{ |>}->}[r]_-{k} &
  X \ar@{-{ >>}}[r]_-{e} \ar@/^3ex/[rr]^-{f} &
  \CoKer{\kappa}=I\ \ar[r]_-{v} &
  Y.
  }
  }
\end{equation}
In this diagram $k = \KerMap{f}$, and $e=\CoKerMap{k}$. %
\index{image!of a morphism in a normal category}%
\index{image!factorization in a normal category}%

We will show that, in a normal category,  the map $v$ in the normal epi-factorization  of $f$ is always a monomorphism. The means that every map $f$ has an image factorization; see (\ref{def:ImageFactorization}). So, when working in a normal category, we use the terms `normal epi factorization' and `image factorization' interchangeably. %
\index{image factorization}%

As an application, we show that in a normal category, normal, regular, effective, strong, and extremal epimorphisms coincide. This enables us to interpret image factorizations within the framework of orthogonal factorization systems; see (\ref{thm:OFS-(NormalEpis,Monos)}).

We turn to developing auxiliary results which we use to show that, in a homological category, every morphism admits an image factorization.

\begin{lemma}[$f$ factors through $\CoKer{\KerMap{f}}$ via map with trivial kernel\NTag]
  \label{thm:FactorizationThroughCoKer(Ker(f))ViaKer-0-Map}%
  In the factorization diagram \ref{eq:ImageFactorization}, the kernel $\Ker{v}$ of $v$ is $\ZeroObject$.
\end{lemma}
\begin{proof}
  Via the universal properties of kernels and cokernels the morphism $f$ yields the commutative diagram below.
  \begin{equation*}
    \xymatrix@R=5ex@C=4em{
    \Ker{f} \ar@{{ |>}->}[r]^-{k} \ar[d]_{\bar{e}} &
    X \ar[r]^-{f} \ar@{-{ >>}}[d]_{e} &
    Y \ar@{=}[d] \\
    \Ker{v} \ar@{{ |>}->}[r]_-{k'} &
    I \ar[r]_-{v} & Y
    }
  \end{equation*}
  As the diagram commutes,
  \begin{equation*}
    k' \bar{e} = ek = \ZeroMap = \ZeroMap \bar{e}
  \end{equation*}
  With Proposition~\ref{thm:PullbackRecognition-KernelSide-1} we see that the left hand square is a pullback. Via the \PNEInline-property, it preserves the normal epimorphism $e$; i.e., $\bar{e}$ is a normal epimorphism. Hence from the calculation above, it follows that $k'=\ZeroMap$.
\end{proof}

\begin{corollary}[Existence of image factorizations\NTag]
  \label{thm:NEM-Img-Fact-Existence}%
  In a normal category, in the normal epi-factorization $f=ev$ of a morphism $f$, the map $v$ is a monomorphism.
\end{corollary}
\begin{proof}
  From (\ref{thm:FactorizationThroughCoKer(Ker(f))ViaKer-0-Map}), we know that $\Ker{v}=\ZeroObject$. Then (\ref{thm:Mono<->0-Kernel}) tells us that $v$ is a monomorphism.
\end{proof}

This means that the normal epi-factorization of $f$ is an image factorization. - We turn to comparing various types of epimorphisms in a normal category.

\begin{corollary}[Extremal epimorphism is normal\NTag]
  \label{thm:ExtremalEpi->Cokernel}%
  In a normal category, every extremal epimorphism is a normal epimorphism.
\end{corollary}
\begin{proof}
  If $f\from X\to Y$ is an extremal epimorphism, consider the image factorization (\ref{thm:NEM-Img-Fact-Existence}) of $f$ through the cokernel of $k =\KerMap{f}$:
  \begin{equation*}
    \xymatrix@R=5ex@C=4em{
    && \CoKer{k} \ar@{ >->}[d]^-{m} \\
    \Ker{f} \ar@{{ |>}->}[r]_-{k} &
    X \ar[r]_-{f} \ar@{-{ >>}}[ru]^-{e} &
    Y
    }
  \end{equation*}
  As $f$ is an extremal epimorphism, $m$ is an isomorphism. So $f$, along with $e$, is a normal epimorphism.
\end{proof}

We can now establish the equivalence of certain types of epimorphism in normal categories displayed in Figure~\ref{fig:Epis}.

\begin{corollary}[Equivalence of various epimorphisms\NTag]
  \label{thm:CoKer=NormalEpi=RegEpi=EffectiveEpi}%
  \index{cokernel!same as strong epimorphism}%
  \index{cokernel!same as extremal epimorphism}%
  \index{strong epimorphism!same as cokernel}%
  \index{strong epimorphism!same as extremal epimorphism}%
  \index{extremal epimorphism!same as cokernel}%
  \index{extremal epimorphism!same as strong epimorphism}%
  For any morphism $ f\from{X\to Y}$, the following conditions are equivalent:
  \begin{tfae}
    \item $f$ is a normal epimorphism.
    \item $f$ is a regular epimorphism.
    \item $f$ is an effective epimorphism.
    \item $f$ is a strong epimorphism.
    \item $f$ is an extremal epimorphism.
  \end{tfae}
\end{corollary}
\begin{proof}
  In a category with finite limits, we have (I) $\implies$ (II) $\Leftrightarrow$ (III), and (II) $\implies$ (IV) $\implies$ (V). The implication (V) $\implies$ (I) from (\ref{thm:ExtremalEpi->Cokernel}) completes the proof.
\end{proof}

\begin{proposition}[Pullbacks preserve image factorizations and normal maps\NTag]
  \label{thm:PullbacksPreserveImageFactorizations}%
  Pull the map $u$ back along a map $f$ so as to obtain the map $\bar{u}$:
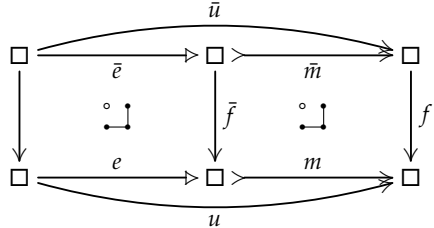
  \begin{equation*}
    \xymatrix@R=6ex@C=5em{
    \DiagObj \ar@{-{ >>}}[r]_-{\bar{e}} \ar@/^2ex/[rr]^-{\bar{u}} \ar[d] \PullLU{rd} &
    \DiagObj \ar@{{ >}->}[r]_-{\bar{m}} \ar[d]^-{\bar{f}} \PullLU{rd} &
    \DiagObj \ar[d]^-{f} \\
    \DiagObj \ar@/_2ex/[rr]_-{u} \ar@{-{ >>}}[r]^-{e} &
    \DiagObj \ar@{{ >}->}[r]^-{m} &
    \DiagObj
    }
  \end{equation*}
  If $u=me$ is an image factorization of $u$, then an image factorization of $\bar{u}$ is given by $\bar{u}=\bar{m}\bar{e}$, where $\bar{m}$ is the pullback of $m$ along $f$ and $\bar{e}$ is the pullback of $e$ along $\bar{f}$. Further, if $m$ is a normal monomorphism then $\bar{m}$ is a normal monomorphism as well.
\end{proposition}
\begin{proof}
  The map $f$ pulls the monomorphism $m$ back to the monomorphism $\bar{m}$ by (\ref{exe:BaseChange-Mono-Epi-Iso}). Then $\bar{f}$ pulls the normal epimorphism $e$ axiomatically back to the normal epimorphism~$\bar{e}$. Finally, the concatenation of these two pullbacks is again a pullback by (\ref{thm:Pullbacks,ConcatenatedSquares}). Thus $\bar{m}\Comp \bar{e}$ is the image factorization of the pullback $\bar{u}$ of $u$ along $f$. Finally, if $m$ happens to be a normal monomorphism, then $\bar{m}$ is a normal monomorphism as well by (\ref{thm:KernelFunctor-Props}).
\end{proof}

As a consequence, we see that, in a normal category, normal epimorphisms enjoy the following properties:

\begin{proposition}[Properties of normal epimorphisms\NTag]
  \label{thm:Cokernels-Props-Normal}
  \label{thm:NormalEpis-Props-Normal}
  For morphisms $ f\from{X\to Y}$, $g\from{Y\to Z}$, and $ f'\from{X'\to Y'}$ in a normal category, the following hold. %
  \index{normal!epimorphism - properties}
  \begin{thmlist}
    \item \label{thm:Cokernels-Props.gfNEpi->gNEpi}%
    If $ g\Comp f$ is a normal epimorphism, then so is $g$.
    \item \label{thm:Cokernels-Props.f,gNEpis->gfNEpi}%
    If $f$ and $g$ are normal epimorphisms, then so is $g\Comp f$. %
    \index{normal!epimorphism - composite}\index{composite!of normal epimorphisms}
    \item \label{thm:Cokernels-Props.f,g'NEpis->fxf'NEpi}%
    If $f$ and $f'$ are normal epimorphisms, then so is $f\times f'$. %
    \index{normal!epimorphism - product}\index{product!of normal epimorphisms}
  \end{thmlist}
\end{proposition}
\begin{proof}
  We use the fact (\ref{thm:CoKer=NormalEpi=RegEpi=EffectiveEpi}) that, under the stated hypotheses, normal epimorphisms, strong epimorphisms and extremal epimorphisms are equivalent concepts. In situation (\ref{thm:Cokernels-Props.gfNEpi->gNEpi}) it therefore suffices to show that $g$ is a extremal epimorphism. Consider a factorization of $g$ through a subobject $m$ of $Z$:
  \begin{equation*}
    \xymatrix@R=5ex@C=4em{
    && M \ar@{ >->}[d]^-{m} \\
    X \ar[r]_-{f} &
    Y \ar[r]_-{g} \ar[ru]^-{\gamma} &
    Z
    }
  \end{equation*}
  Then $ \gamma\Comp f$ is a factorization of the extremal epimorphism $ g\Comp f$ through $m$. So $m$ is an isomorphism, implying that $g$ is an extremal epimorphism.

  (\ref{thm:Cokernels-Props.f,gNEpis->gfNEpi}) It suffices to show that a composite of strong epimorphisms is a strong epimorphism. This is the content of \eqref{exe:StrongEpi-Props}.

  (\ref{thm:Cokernels-Props.f,g'NEpis->fxf'NEpi}) We use that normal epimorphisms are pullback-stable to see that in the two pullback squares
  \begin{equation*}
    \vcenter{
    \xymatrix@R=5ex@C=4em{X\times X' \ar[d]_-{\pi_{X'}} \ar@{-{ >>}}[r]^-{1_{X}\times f'} & X\times Y' \ar[d]^-{\pi_{Y'}}\\
    X' \ar@{-{ >>}}[r]_-{f'} & Y'}}
    \qquad\qquad
    \vcenter{
    \xymatrix@R=5ex@C=4em{X\times Y' \ar[d]_-{\pi_{X}} \ar@{-{ >>}}[r]^-{f\times 1_{Y'}} & Y\times Y' \ar[d]^-{\pi_{Y}}\\
    X \ar@{-{ >>}}[r]_-{f} & Y}}
  \end{equation*}
  the top horizontal arrows are normal epimorphisms. The result now follows from (ii), because $f\times f'=(f\times \IdMapOn{Y})\Comp(\IdMapOn{X}\times f')$.
\end{proof}

\begin{corollary}[Product of short exact sequences\NTag]
  \label{thm:Product-SESs}
  If the sequences $K\XRA{k} X \XRA{q} Q$ and $L\XRA{l} Y \XRA{r} R$ are short exact, then so is the sequence %
  \index{product!of short exact sequences}\index{short exact sequence!product}%
  \begin{equation*}
    \xymatrix@R=5ex@C=4em{
    \Prdct{K}{L} \ar@{{ |>}->}[r]^-{k\times l} &
    \Prdct{X}{Y} \ar@{-{ >>}}[r]^-{q\times r} &
    \Prdct{Q}{R}
    }
  \end{equation*}
\end{corollary}
\begin{proof}
  We know	from (\ref{thm:NormalEpis-Props-Normal}) that the product $\Prdct{q}{r}$ of the normal epimorphisms $q$ and $r$ is again a normal epimorphism. So $\Prdct{q}{r}$ is the cokernel of its kernel, by (\ref{thm:Ker(CoKer)-CoKer(Ker)}). However, $\KerMap{\Prdct{q}{r}} = \Prdct{\KerMap{q}}{\KerMap{r}} = \Prdct{k}{l}$ because limits commute with limits. This proves that a product of two short exact sequences is short exact.
\end{proof}

\begin{corollary}[Cobase change preserves extremal epimorphisms\NTag]
  \label{thm:Regular/Strong/ExtremalEpis-CoBaseChange}
  Pushouts along any map preserve regular/strong/extremal epimorphisms. %
  \index{cobase change!preserves regular epimorphisms}\index{cobase change!preserves strong epimorphisms}%
  \index{cobase change!preserves extremal epimorphisms}%
\end{corollary}
\begin{proof}
  From (\ref{thm:PushoutPreservesNormalEpis}) we know that cobase change preserves normal epimorphisms. From (\ref{thm:CoKer=NormalEpi=RegEpi=EffectiveEpi}) we know that in a normal category the concepts of normal epimorphism and regular / strong / extremal epimorphism coincide.
\end{proof}

\begin{corollary}[Recognizing normal epi via jointly extremal epis\NTag]
  \label{thm:JEE-then-cokernel-NE}%
  In the diagram below, if the maps $\lambda$ and $u$ are jointly extremally epimorphic, then $v$ is a normal epimorphism.
  \begin{equation*}
    \xymatrix@R=5ex@C=4em{
    K \ar[r]^-{\kappa} \ar[d]_{k} &
    U \ar@{-{ >>}}[r]^-{\CoKerMap{\kappa}} \ar[d]_{u} &
    V \ar[d]^{v} \\
    L \ar[r]_-{\lambda} &
    X \ar@{-{ >>}}[r]_-{\CoKerMap{\lambda}} &
    Y
    }
  \end{equation*}
\end{corollary}
\begin{proof}
  It follows from the hypotheses that $\CoKerMap{\lambda}\Comp \lambda=0$ and $\CoKerMap{\lambda}\Comp u$ are jointly extremally epimorphic. So, $\CoKerMap{\lambda}\Comp u=v\Comp \CoKerMap{\kappa}$ is an extremal epimorphism. Hence, by (\ref{thm:CoKer=NormalEpi=RegEpi=EffectiveEpi}), it is a normal epimorphism. By (\ref{thm:NormalEpi-Props}.\ref{thm:Cokernels-Props.gfNEpi->gNEpi}), $v$ is a normal epimorphism.
\end{proof}

\begin{proposition}[Image factorization of morphism out of a binary sum\NTag]
  \label{thm:ImageFactorizationCommutesBinarySums}
  For any two morphisms $f\from X\to Z$ and $g\from Y\to Z$:
  \begin{equation*}
    \Img{f}\join \Img{g} = \Img{\SumMapOutOf{f,g}\from X+Y\to Z}
  \end{equation*}
\end{proposition}
\begin{proof}
  Consider the commutative diagram below:
  \begin{equation*}
    \xymatrix@R=5ex@C=4em{
    && X+Y \ar[d] \\
    && S \ar@{{ >}->}[d] \\
    X \ar@{-{ >>}}[r] \ar[rruu]^{\InclsnOf{X}}&
    \Img{f} \ar@{{ >}->}[r] \ar@{{ >}.>}[ru] &
    \Img{\SumMapOutOf{f,g}} &
    \Img{g} \ar@{{ >}->}[l] \ar@{{ >}.>}[lu] &
    Y \ar@{-{ >>}}[l] \ar[lluu]_{\InclsnOf{Y}}\\
    }
  \end{equation*}
  The claim follows via (\ref{thm:ExtremalEpiCharacterizationOfJoin}) once we have shown that the canonical inclusions of $\Img{f}$ and $\Img{g}$ in $\Img{\SumMapOutOf{f,g}}$ are jointly extremally epimorphic. Indeed, given a subobject $S$ of $\Img{\SumMapOutOf{f,g}}$ which contains $\Img{f}$ and $\Img{g}$ as subobjects, the normal epimorphism $X+Y\NEpi \Img{\SumMapOutOf{f,g}}$ factors through $S$. With (\ref{thm:NormalEpis-Props-Normal}) it follows that the monomorphism $S\to \Img{\SumMapOutOf{f,g}}$ is a normal epimorphism as well and, hence (\ref{thm:IsomorphismRecognition}), an isomorphism. - This implies the claim.
\end{proof}

\begin{subordinate}{}

  \begin{subsubordinate}{Existence of image factorizations and orthogonal factorization systems}
    The existence of a \NEM factorization for \emph{every} morphism has a direct tie to the theory of orthogonal prefactorization systems; see Appendix~\ref{sec:FactorizationSystems}. Beginning with the class $\EuScript{M}$ of monomorphism in a category $\Ctgry{C}$, one obtains the class $\EuScript{E}$ of strong epimorphisms as its left orthogonal complement. In a normal category strong epimorphisms are normal. Here, then, is a reformulation of  what we found earlier in this section.

    \begin{theorem}[OFS from \NEM factorizations\NTag]
      \label{thm:OFS-(NormalEpis,Monos)}
      We let $\EuScript{E}$ be the class of normal epimorphisms, and $\EuScript{M}$ the class of monomorphisms. Then %
      \index{orthogonal!factorization system}\index{OFS}%
      \begin{thmlist}
        \item every morphism $f\from X\to Y$ in $\Ctgry{C}$ admits a factorization as $f=me$, where $m\in \EuScript{M}$ and $e\in\EuScript{E}$;
        \item $\EuScript{E}$ and $\EuScript{M}$ are each other's \Defn{orthogonal complement}\index{orthogonal!complement}: given any commutative square
        \begin{equation*}
          \xymatrix@R=5ex@C=4em{
          A \ar[d]_{e} \ar[r]^-{u} &
          U \ar[d]^{m} \\
          B \ar[r]_-{v} \ar@{.>}[ru]|-{\ \varphi\ } &
          V,
          }
        \end{equation*}
        \begin{enumerate}
          \item  $e\in \EuScript{E}$ if and only if for every $m\in \EuScript{M}$, a unique filler $\varphi$ exists, while
          \item $m\in \EuScript{M}$ if and only if for every $e\in \EuScript{E}$, a unique filler $\varphi$ exists.
        \end{enumerate}
        \index{left!orthogonal}\index{right!orthogonal}\index{left!orthogonal complement}\index{right!orthogonal complement} %
      \end{thmlist}
    \end{theorem}
    \begin{proof}
      In Lemma~\ref{thm:FactorizationThroughCoKer(Ker(f))ViaKer-0-Map} we show that $\Ker{m}=0$ and, with (\ref{thm:Mono<->0-Kernel}), we infer that $m$ is a monomorphism. Since by definition (\ref{def:StrongEpimorphism}), strong epimorphisms satisfy (a) and one of the implications in (b), it remains to show that every strong epimorphism is a normal epimorphism, and that these satisfy the other implication in (b). The first claim holds because every strong epimorphism is extremal and in a normal category, every extremal epimorphism is a normal epimorphism; see (\ref{thm:ExtremalEpi->Cokernel}). The second claim uses the factorization: suppose $m$ is such that for every $e\in \EuScript{E}$, the diagram admits a unique filler~$\varphi$. We may factor $m=ip$ where $m\in\EuScript{M}$ and $p\in \EuScript{E}$. Then a filler $\varphi$ exists for the square $m1_U=ip$, which shows that $p$ is an isomorphism, and thus $m$ is a monomorphism, because so is $i$.
    \end{proof}
  \end{subsubordinate}
\end{subordinate}

\begin{exercises}

\begin{exercise}[$\CoKer{f}=\ZeroMap$ does not imply $f$ is an epimorphism]
  \label{exe:CoKer(f)=0=/=>fEpi}%
  Find an example of a morphism $f\from X\to Y$ in a normal category where $\CoKer{f}=\ZeroMap$, but $f$ is not an epimorphism.
\end{exercise}

\begin{exercise}[$v\Comp u$ normal monomorphism does not imply that $u$ is normal]
  \label{exe:vuNormalMono=/=>uNormal}
  Find an example of a normal category in which there exists a composite $vu$ which is a normal monomorphism but $u$ is not a normal monomorphism; compare  (\ref{thm:NormalEpis-Props-Normal}) and recall (\ref{thm:NormalMono-Props}.\ref{thm:MonomorphismCancellationInKernel}).
\end{exercise}
\end{exercises}
\section[The Snake Lemma]{The Snake Lemma}
\label{sec:SnakeLemma-N+DPN}

The following proof of the Snake Lemma is inspired by the proof of the Snake Lemma in~\cite{FBorceuxDBourn2004}. It is valid in normal categories satisfying the condition \DPNInline. It is independent of the \ANNInline-condition which forms the basis of the proof given in (\ref{thm:SnakeLemma-Classical}).

\begin{theorem}[Classical Snake Lemma\NDPNTag]
  \label{thm:SnakeLemma-Classical-DPN-NormalCat}%
  In a normal category satisfying the \DPNInline-condition, suppose the maps $\kappa$, $\xi$ and $\rho$ in the following morphism of short exact sequences are normal.
  \index{Snake Lemma!classical in \DPNInline-normal cat}%
  \begin{equation*}
    \xymatrix@R=5ex@C=4em{
    K \ar[d]_{\kappa} \ar@{{ |>}->}[r]^-{k} &
    X \ar[d]_{\xi} \ar@{-{ >>}}[r]^-{q} &
    Q \ar[d]^{\rho} \\
    L\ar@{{ |>}->}[r]_-{l} &
    Y \ar@{-{ >>}}[r]_-{r} &
    R
    }
  \end{equation*}
  Then the kernels and cokernels of $\kappa$, $\xi$, $\rho$ form a six-term exact sequence
  \begin{equation*}
    \xymatrix@R=5ex@C=2.2em{
    0 \longrightarrow \Ker{\kappa} \ar@{{ |>}->}[r]^-{k^{\ast}} &
    \Ker{\xi} \ar[r]^-{q^{\ast}} &
    \Ker{\rho} \ar[r]^-{\partial} &
    \CoKer{\kappa} \ar[r]^-{l_{\ast}} &
    \CoKer{\xi} \ar@{-{ >>}}[r]^-{r_{\ast}} &
    \CoKer{\rho} \longrightarrow 0
    }
  \end{equation*}
  This six-term exact sequence depends functorially on morphisms of the underlying diagrams of short exact sequences.
\end{theorem}
\begin{proof}
  We construct the diagram below starting from the two columns on the right.
  \begin{equation*}
    \xymatrix@!0@R=6.5ex@C=3em{
    \Ker{\kappa} \ar@{{ |>}->}[rr]^-{k^{\ast}} \ar@{{ |>}->}[d] \PullLU{rrd} &&
    \Ker{\xi} \ar@{{ |>}->}[d] \ar@{=}[rr] &&
    \Ker{\xi} \ar[rr]^-{q^{\ast}} \ar@{{ |>}->}[d] &&
    \Ker{\rho} \ar@{{ |>}->}[d] \\
    K \ar@{{ |>}->}[rr] \ar[dd]_-{\kappa} \ar@{-{ >>}}[rd] &&
    \Ker{e} \BiCart{rrrd} \ar[dd] \ar@{{ |>}->}[rr] \ar@{-{ >>}}[rd] &&
    X \ar@{-{ >>}}[rrrd]|(.4){\ e\ }|(.66){\hole} \ar[dd]_(.75){\xi} \ar@{-{ >>}}[rd] \ar@{-{ >>}}[rr]^-{q} &&
    Q \ar[dd]_(.75){\rho} \ar@{-{ >>}}[rd] \\
    & \Img{\kappa} \ar@{{ |>}->}[ld] \ar@{{ |>}->}[rr]|\hole &&
    \Ker{\eta} \ar@{{ |>}->}[ld] \ar@{{ |>}->}[rr]|\hole &&
    \Img{\xi} \ar@{{ |>}->}[dl] \ar@{-{ >>}}[rr]|\hole_(.3){\eta} &&
    \Img{\rho} \ar@{{ |>}->}[ld] \\
    L \ar@{-{ >>}}[d] \ar@{=}[rr] &&
    L \ar@{-{ >>}}[d] \ar@{{ |>}->}[rr] &&
    Y \ar@{-{ >>}}[d] \ar@{-{ >>}}[rr] &&
    R \ar@{-{ >>}}[d] \\
    \CoKer{\kappa} \ar@{-{ >>}}[rr] \ar@/_2ex/[rrrr]_-{l_{\ast}} &&
    \Ker{r_{\ast}} \ar@{{ |>}->}[rr] &&
    \CoKer{\xi} \ar@{-{ >>}}[rr]_-{r_{\ast}} &&
    \CoKer{\rho}
    }
  \end{equation*}
  The map $e$ is a normal epimorphism, being a composite of two such. So, $\eta$ is a normal epimorphism as well. Via the \HSDInline-property (\ref{thm:NormalCat->HomologicallySelfDual}), the morphism of short exact sequences
  \begin{equation*}
    (\Ker{e}\to \Ker{\eta},X\NEpi \Img{\xi},\Img{\rho}=\Img{\rho})
  \end{equation*}
  yields the bi-cartesian square  $\Ker{e}\rightrightarrows \Img{\xi}$ in the middle and the identity $\Ker{b}=\Ker{b}$ in the upper row.

  From the morphism of short exact sequences $(\Ker{\eta}\to L,\Img{\xi}\NMono Y,\Img{\rho}\to R)$, we see that the square $\Ker{\eta}\rightrightarrows Y$ is a pullback. With the \DPNInline-condition we obtain the di-extension whose middle object is $Y$; see (\ref{thm:DiExtensionFromDPN}).

  The map $\Img{\kappa}\NMono L$ factors through $\Ker{\eta}$ via a normal monomorphism. The Pure Snake Lemma (\ref{thm:HomologicalSelfDuality-Recognize}) yields a functorial isomorphism
  \begin{equation*}
    S\DefEq \CoKer{\Img{\kappa}\to \Ker{\eta }}\cong \Ker{\CoKer{\kappa}\to \Ker{r_{\ast}}} = \Ker{l_{\ast}}
  \end{equation*}
  The upper-left part of the diagram extends to the following di-extension.
  \begin{equation*}
    \xymatrix@!0@=4em{
    \Ker{\kappa} \ar@{{ |>}->}[r]^-{k^{\ast}} \ar@{{ |>}->}[d]  &
    \Ker{l} \ar@{{ |>}->}[d] \ar@{-{ >>}}[r] \ar[rd]^{q^\ast } &
    J \ar@{{ |>}->}[d]\\
    K \ar@{{ |>}->}[r] \ar@{-{ >>}}[d] &
    \Ker{e}  \ar@{-{ >>}}[r] \ar@{-{ >>}}[d] &
    \Ker{\rho} \ar@{-{ >>}}[d]\\
    \Img{\kappa} \ar@{{ |>}->}[r] &
    \Ker{\eta} \ar@{-{ >>}}[r] &
    S
    }
  \end{equation*}
  The bottom short exact sequence has been constructed above. The horizontal short exact sequence in the middle comes from homological self-duality:
  \begin{equation*}
    \xymatrix@R=5ex@C=4em{
    K \ar@{{ |>}->}[r] \ar@{=}[d] &
    \Ker{e} \ar@{-{ >>}}[r] \ar@{{ |>}->}[d] &
    \Ker{\rho} \ar@{{ |>}->}[d] \\
    K \ar@{{ |>}->}[r] &
    X \ar@{-{ >>}}[r] \ar@{-{ >>}}[d]_{e} &
    Q \ar@{-{ >>}}[d] \\
    & \Img{\rho} \ar@{=}[r] &
    \Img{\rho}
    }
  \end{equation*}
  We have $k^{\ast}=\KerMap{q^{\ast}}$ because kernels commute with kernels. So, $q^\ast$ admits an antinormal decomposition. Given that its di-inverse is the normal map $K\to \Ker{\eta}$, $q^{\ast}$ is normal as well. Assembling the available information, we compute:
  \begin{equation*}
    \CoKer{q^{\ast}} \cong \CoKer{J\NMono\Ker{\rho}} \cong S \cong \Ker{l_{\ast}}
  \end{equation*}
  This yields the required $6$-term exact sequence. It is functorial because we assume that kernels and cokernels are functorial.
\end{proof}

\begin{remark}[Snake when normal monos are closed under composition]
  An analysis of the proof of (\ref{thm:SnakeLemma-Classical-DPN-NormalCat}) shows that it is valid in \ZExact\ categories which enjoy the  \DPNInline-property, and in which \emph{normal epimorphisms} are closed under composition. Therefore, dual reasoning shows that the Classical Snake Lemma is valid in \ZExact\ categories which enjoy the \DPNInline-property, and in which \emph{normal monomorphisms} are closed under composition.
\end{remark}

The Relaxed Snake Lemma  (\ref{thm:SnakeLemma-Relaxed}) is also valid under the assumptions stated in (\ref{thm:SnakeLemma-Classical-DPN-NormalCat}).

\begin{exercises}

\begin{exercise}[Classical Snake vs. border cases of $(\Prdct{3}{3})$-Lemma\ZExactTag]
  \label{exe:SnakeVs(3,3)-Lemma}%
  In a \ZExact\ category assume that normal monomorphisms are closed under composition or that normal epimorphisms are closed under composition. Show that the Classical Snake Lemma holds if and only if the border cases of the $(\Prdct{3}{3})$-Lemma (\ref{thm:Dinversion-PreservationNormalMaps-Border(3x3)}) hold.
\end{exercise}
\end{exercises}
\section[Alternate Characterizations of Normal Categories]{Alternate Characterizations of Normal Categories}
\label{sec:AlternateNormal}

We compare the concept of a normal category, as defined in Section \ref{sec:NormalCats} with the concept of a normal category as defined by Z.~Janelidze in \cite[Sec.~1]{ZJanelidze-Snake}. To express the distinction, we call a category \emph{Z.~Janelidze-normal} if it is pointed, regular, and regular epimorphisms are normal epimorphisms. In (\ref{thm:CharNormalCat}), we show that a finitely cocomplete Z.~Janelidze-normal category is the same thing as a normal category from~(\ref{sec:NormalCats}). %

In (\ref{thm:SubtractiveVariety-DinversionPreservesNormalMaps}) we show that the \DPNInline-condition is satisfied in all normal subtractive varieties~\cite{ZJanelidze-Snake}, i.e., in any subtractive variety where all surjective algebra maps are normal epimorphisms, see (\ref{def:NormalCat}). By (\ref{thm:SubtractiveVariety->NormalCat}), every pointed variety with a binary operation `$-$' such that $x-x=0$ and $x-0=x$, and where $x-y=0$ implies $x=y$ is normal subtractive.

As it turns out, in a normal context, the condition \DPNInline\ is characteristic of subtractivity: we proceed by showing that a normal category satisfies \DPNInline\ if and only if it is \emph{subtractive} in a suitable, purely categorical sense~\cite{ZJanelidze-Snake}. We approach this concept by means of the \emph{difference object} $D(X)$, which is the kernel $\delta_X\colon D(X)\to X+X$ of the codiagonal map $\FoldOn{X}\from X+X\to X$. The difference object enables us to define the \emph{formal difference} of two morphisms $f$, $g\from X\to Y$ as the composite $f-g=\SumMapOutOf{f,g}\Comp \delta_{X}\from D(X)\to Y$. In (\ref{thm:FormalDifference-Composition}), we show that the formal difference behaves well with respect to composition. In (\ref{thm:NormalCatViaSubtraction}), we characterize normal categories by means of the formal difference. In (\ref{thm:DPN-Then-Subtractive}), we show that any normal category in which the \DPNInline-property holds is subtractive. Conversely, in (\ref{thm:DPN-for-Sub-Cats}), we show that a normal subtractive category satisfies the \DPNInline-property.%

From the perspective of categorical algebra, Axiom \PNEInline\ of (\ref{def:NormalCat}) is at the core of the concept of a \emph{regular} category, introduced in~\cite{Barr-Grillet-vanOsdol}; see (\ref{def:RegularCategory}). In the presence of the other axioms characterizing a normal category, it implies that regular epimorphisms are stable under pullbacks.

\begin{definition}[Z.~Janelidze-normal category]
  \label{def:JanNormalCat}%
  \cite{ZJanelidze-Subtractive}\quad A \Defn{Z.~Janelidze-normal category} is a regular category with a zero object in which every regular epimorphism is a normal epimorphism.%
  \index{normal!category}\index{category!normal}\index{Z.~Janelidze-normal category}\index{category!Z.~Janelidze-normal}%
\end{definition}

In a Janelidze-normal category the \PNEInline-condition is satisfied, because normal epimorphisms and regular epimorphisms coincide.

\begin{theorem}[Characterization of normal categories]
  \label{thm:CharNormalCat}%
  For  a pointed category $\Ctgry{X}$, the following are equivalent:
  \begin{tfae}
    \item $\Ctgry{X}$ is a normal category, as per Definition \ref{def:NormalCat}.
    \item $\Ctgry{X}$ is finitely cocomplete and Z.~Janelidze-normal.
  \end{tfae}
\end{theorem}
\begin{proof}
  Suppose $\Ctgry{X}$ is normal. In $\Ctgry{X}$, normal epimorphisms and regular epimorphisms coincide by Corollary~\ref{thm:CoKer=NormalEpi=RegEpi=EffectiveEpi}. Then every morphism has an image factorization by  Corollary~\ref{thm:NEM-Img-Fact-Existence}. By (\ref{thm:PullbacksPreserveImageFactorizations}), such an image factorization is preserved by pullbacks because of the \PNEInline-condition. So, $\Ctgry{X}$ is Z.~Janelidze-normal.

  Conversely, it only remains to show that the \AENInline-condition is satisfied. This is so because an absolute epimorphism is always a regular epimorphism.
\end{proof}

We are now ready to characterize normal categories satisfying the \DPNInline-condition by means of a categorical version of the subtractivity condition. We first consider the varietal case, because as we shall see, the proof there can be mimicked in a purely categorical setting.

\begin{proposition}[Subtractive varieties and dinversion of normal maps]
  \label{thm:SubtractiveVariety:DPN}
  \label{thm:SubtractiveVariety-DinversionPreservesNormalMaps}%
  In a normal subtractive variety, dinversion preserves normal maps.
\end{proposition}
\begin{proof}
  We show that condition (\ref{thm:DPNPNvs3x3-Right3x3}) of Proposition~\ref{thm:Dinversion-PreservationNormalMaps-Border(3x3)} holds. We use the notation of \eqref{eq:3x3}. First, we see that the bottom right square is a pushout by (\ref{thm:PushoutRecognize-Categorical}). So, $z$ is a normal epimorphism, and it remains to show that $w$ is the kernel of $z$.

  To see this, note first that $zw=\ZeroObject$; i.e., the set image of $w$ is contained in $\Ker{z}$. Now let $q\in Q$ such that $z(q)=\ZeroObject$. Pick $x\in X$ with $e(x)=q$. Then $f(y(x))=z(e(x))=z(w)=\ZeroObject$. So, there is $j\in J$ with $c(j)=y(x)$. Pick $k\in K$ such that $x(k)=j$. Then $y$ sends $x-b(k)$ to $0\in R$. So, there is $l\in L$ with $v(l)=x-b(k)$. Now $w(d(l))=e(x-(b(k)))=q$.

  So, every element of $\Ker{z}$ belongs to the set image of $w$. Thus, the argument is complete once we show that $w$ is injective or, equivalently, $\Ker{w}=\ZeroMap$. In this last step, it is essential that a subtractive category is a normal category. So monomorphisms are exactly those maps whose kernel vanishes; see (\ref{thm:Mono<->0-Kernel}).

  Take $i\in I$ with $w(i)=0$. Let $l\in L$ with $d(l)=i$. Then $e(v(l))=0$ so that $k\in K$ exists with $b(k)=v(l)$. Now $x(k)$ is zero because $c$ is a monomorphism and $y(b(k))=y(v(l))=\ZeroObject$. Which gives us $m\in M$ with $u(m)=k$. Now $a(m)=l$ since $v$ is a monomorphism and $v(a(m))=b(u(m))=b(k)=v(l)$. Hence $i=d(l)=d(a(m))=0$.
\end{proof}

Given an object $X$ in a normal category $\Ctgry{X}$, the codiagonal map $\FoldOn{X} \DefEq \SumMapOutOf{\IdMapOn{X}$, $\IdMapOn{X}}\from X+X\to X$ is sectioned by the structure maps $\InclsnOf{1}$, $\InclsnOf{2}\from X\to X+X$. Thus $\FoldOn{X}$ is an absolute, hence normal, epimorphism which gives rise to the sectioned short exact sequence
\begin{equation*}
  \xymatrix@R=5ex@C=2em{
  D(X) \ar@{{ |>}->}[rr]^-{\delta_X} &&
  X+X \ar@{-{ >>}}@<.5ex>[rr]^-{\FoldOn{X}} &&
  X \ar@<.5ex>[ll]^-{\InclsnOf{1}}
  }
\end{equation*}
Using functorial (co)limits, we obtain a functor $D\from \Ctgry{X}\to \Ctgry{X}$.

\begin{definition}[Difference object\NTag]
  \label{def:DifferenceObject}%
  Given an object $X$ in a normal category, the associated \Defn{difference object} is $D(X)$. %
  \index{difference!object}\index[not]{d!$D(X)$\IndSep difference object of $X$}%
\end{definition}

\begin{example}[Difference objects in $\AbGrps$]
  \label{exa:DifferenceObject-Ab}%
  In the category $\AbGrps$, the coproduct $X+X$ is the biproduct $X\oplus X$, so that the kernel $\delta_X$ is represented by the map $\PrdctMapInto{\IdMapOn{X},-\IdMapOn{X}}\from X\to X\oplus X$. Hence the functor $D$ is naturally equivalent to the identity functor on $\AbGrps$.%
  \index{difference!object in $\AbGrps$}%
\end{example}

\begin{example}[Difference objects in $\Grps$]
  Given an element $x$ of a group $X$, we write $\bar{x}\DefEq \InclsnOf{1}(x) \in X+X$ and $\underline{x}\DefEq \InclsnOf{2}(x)\in X+X$ for the images of $x$ in the respective summands of $X+X$. Then the difference group $D(X)$ is the normal subgroup generated by elements of the shapes
  \begin{equation*}
    \bar{x}\underline{x}^{-1}  \qquad \text{and}\qquad \underline{x}(\bar{x})^{-1},\quad x\in X
  \end{equation*}
\end{example}
\begin{proof}
  An arbitrary element $w$ of $X+X$ may be written as $w=\overline{x_1}\underline{y_1}\cdots \overline{x_n}\underline{y_n}$. We argue by induction. If $n=1$, then $\FoldOn{X}(w)=1$ if and only if  $y_{1}=x_{1}^{-1}$. So, assume the claim is true for $1\leq n-1$. We notice that $\FoldOn{X}(w)=1$ if and only if $\FoldOn{X}\left( \underline{x_{1}}^{-1} w \underline{x_{1}}\right) = 1$. Then, we compute:
  \begin{equation*}
    \begin{array}{rcl}
      1                      & = & \FoldOn{X} \left( \underline{x_{1}}^{-1}\, \left[ \overline{x_1}\underline{y_1}\cdots  \overline{x_n}\underline{y_n}\right] \ \underline{x_{1}}\right)                                                          \\
      \Leftrightarrow\qquad1 & = & \FoldOn{X}\left( (\underline{x_{1}}^{-1}\overline{x_{1}})\ \left[  \underline{y_{1}}\ \overline{x_{2}} \underline{y_{2}}\cdots \overline{x_{n}}\underline{y_{n}}\ \underline{x_{1}} \right] \right)             \\
      \Leftrightarrow\qquad1 & = & \FoldOn{X}\left( (\underline{x_{1}}^{-1}\overline{x_{1}}) \right) \ \FoldOn{X}\left( \underline{y_{1}}\ \overline{x_{2}} \underline{y_{2}}\cdots \overline{x_{n}}\underline{y_{n}}\ \underline{x_{1}}   \right) \\
      \Leftrightarrow\qquad1 & = & \FoldOn{X}\left( \underline{y_{1}}\ \left[ \overline{x_{2}}\underline{y_{2}}\ \cdots\ \overline{x_{n}} \underline{y_{n}x_{1}y_{1}}\right]\, \underline{y_{1}}^{-1}\right)                                       \\
      \Leftrightarrow\qquad1 & = & \FoldOn{X}\left[ \overline{x_{2}}\underline{y_{2}}\ \cdots\ \overline{x_{n}} \ \left( \underline{y_{n}x_{1}y_{1}}\right) \right]
    \end{array}
  \end{equation*}
  Hence the claim follows by induction.
\end{proof}

The difference object enables us to formally define the difference of two morphisms outside an abelian environment:

\begin{definition}[Formal difference of two morphims]
  \label{def:DifferenceMorphisms}%
  Given two morphisms $f$, $g\from X\to Y$, their \Defn{formal difference} is the composite %
  \index{formal difference!of two morphisms}\index[not]{d!$f-g$\IndSep formal difference of morphisms}%
  \begin{equation*}
    f-g\DefEq \bigl( D(X) \XRA{\delta_{X}} X+X \XRA{\SumMapOutOf{f,g}}  Y\bigr)
  \end{equation*}
\end{definition}

\begin{example}[The difference of two morphisms of abelian groups]
  \label{exa:Difference-AbGroup-Maps}%
  In the category $\AbGrps$ of abelian groups, the difference object $D(A)$ of $A$ is functorially isomorphic to $A$. If now $f,g\from A\to B$, then this isomorphism identifies the formal difference $(f-g)=\SumMapOutOf{f,g}\Comp \delta_{A}\from D(A)\to B$ with the pointwise difference $f-g\from A\to B$.
\end{example}

\begin{example}[The difference of two group homomorphisms]
  \label{exa:Difference-Group-Maps}%
  In the category of groups, the pointwise difference of $f$ and $g$ is not a morphism. The group homomorphism $f-g\from D(X)\to Y$ amends this: it takes a generator $\overline{x}\underline{x}^{-1}$ and sends it to $f(x)g(x)^{-1}$.
\end{example}

\begin{proposition}[Formal difference maps and composition\HTag]
  \label{thm:FormalDifference-Composition}
  The formal difference construction of maps $f,g\from X\to Y$ has the following properties:
  \begin{thmlist}
    \item For $z\from Y\to Z$, \ $z\Comp (f-g) = (z\Comp f) - (z\Comp g)$
    \item For $a\from A\to X$, \ $(f-g)\Comp D(k) = (f\Comp a) - (g\Comp a)$
  \end{thmlist}
\end{proposition}
\begin{proof}
  We compute
  \begin{equation*}
    z\Comp (f-g)=z\Comp \SumMapOutOf{f,g}\Comp \delta_{X} = \SumMapOutOf{(z\Comp f),(z\Comp g)} \Comp \delta_{X}=(z\Comp f) - (z\Comp g)
  \end{equation*}
  Similarly,
  \begin{equation*}
    (f-g)\Comp D(a)=\SumMapOutOf{f,g}\Comp \delta_{X}\Comp D(k)=\SumMapOutOf{f,g}\Comp (a+a)\Comp \delta_A=(f\Comp a) - (g\Comp a)
  \end{equation*}
  This was to be shown.
\end{proof}

For the following proposition, let us observe that the difference object is well defined under the weaker assumptions stated there. We just need to be aware that the absolue epimorphism $\FoldOn{X}$ need not be normal.

\begin{proposition}
  \label{thm:NormalCatViaSubtraction}
  A pointed regular finitely bicomplete category is normal if and only if for any two given morphisms~$f,g\from{X\to Y}$,
  \begin{equation*}
    f-g = 0\quad\To\quad f=g
  \end{equation*}
\end{proposition}
\begin{proof}
  In a normal category, the morphism $\nabla_X$ is a normal epimorphism, so that $\SumMapOutOf{f,g}\delta_X=f-g = 0$ gives us a unique factorization $h\colon X\to Y$ of $\SumMapOutOf{f,g}$ over $\nabla_X$. This means that $\SumMapOutOf{f,g}=h\nabla_X=\SumMapOutOf{h,h}$ and thus $f=h=g$.

  For the converse, via (\ref{thm:CharNormalCat}) it suffices to prove that every regular epimorphism is normal. To see this, we show that, for a pair of parallel arrows $f$,~$g\from{X\to Y}$, their coequalizer is a cokernel of $f-g\from D(X)\to Y$. Consider the commutative diagram below, in which $e$ is a coequalizer of $f$ and $g$.
  \begin{equation*}
    \xymatrix@R=5ex@C=4em{
    D(X) \ar@{{ |>}->}[r]^-{\delta_{X}} \ar[rd]_{f-g} &
    X+X \ar[r]^-{\FoldOn{X}} \ar[d]^{\SumMapOutOf{f,g}} &
    X \ar@<0.5ex>[r]^-{f} \ar@<-0.5ex>[r]_-{g} &
    Y \ar[d]^{e} \\
    & Y \ar[rr]_-{e} &&
    Q
    }
  \end{equation*}
  If $h\from Y\to Z$ satisfies $h(f-g)=hf-hg=\ZeroMap$, then $hf=hg$ by assumption. So, the universal property of the coequalizer implies the universal property of the cokernel.
\end{proof}

The difference $\rho_X\DefEq \IdMapOn{X} -0\from D(X)\to X$ plays a special role. It appears in the next definition:

\begin{definition}[Normal subtractive category]
  \label{def:SubtractiveCat}%
  A normal category is \Defn{subtractive} if, for each object $X$, the canonical map
  \[
    \rho_{X}\DefEq \IdMapOn{X}- \ZeroMap\from D(X)\to X
  \]
  is a normal epimorphism. %
  \index{subtractive!normal category}\index{normal!subtractive category}%
\end{definition}

Thus, a subtractive category comes equipped with a normal epic natural transformation $\rho\from D\Rightarrow \IdMapOn{\Ctgry{X}}$.

\begin{lemma}[Condition for $f=g$\NTag]
  \label{thm:f=gRecognizeSubtractive}
  In a subtractive category, suppose $f$, $g\from X\to Y$ satisfy $(f-\ZeroMap)=(g-\ZeroMap) \from D(X)\to Y$, then $f=g$.
\end{lemma}
\begin{proof}
  We compute:
  \begin{equation*}
    f\Comp \rho_{X}=f\Comp (\IdMapOn{X}- \ZeroMap) = f-\ZeroMap=g-\ZeroMap=g\Comp (\IdMapOn{X}-\ZeroMap)=g\Comp \rho_{X}
  \end{equation*}
  so that $f=g$, because $\rho_{X}$ is an epimorphism.
\end{proof}

\begin{proposition}[\DPNInline\ implies subtractivity\NTag]
  \label{thm:DPN-Then-Subtractive}%
  Any normal category in which the \DPNInline-property holds is subtractive.
\end{proposition}
\begin{proof}
  We contemplate the di-extension
  \begin{equation*}
    \xymatrix@R=5ex@C=3em{
    \DiagObj \ar@{{ |>}->}[r] \ar@{{ |>}->}[d] &
    X\flat X \ar@{-{ >>}}[r]^-{\xi} \ar@{{ |>}->}[d]_{\kappa_{X}^{X}} &
    X \ar@{=}[d] \\
    D(X) \ar@{{ |>}->}[r]^-{\delta_X} \ar@{-{ >>}}[d]_-{\rho_X} &
    X+X \ar@{-{ >>}}[r]^-{\nabla_X} \ar@{-{ >>}}[d]_-{\SumMapOutOf{1_X,0} } &
    X \ar[d] \\
    X \ar@{=}[r] &
    X \ar[r] &
    \ZeroObject }
  \end{equation*}
  which comes from the top right antinormal composite $\nabla_X\Comp \kappa_X^X$: Note that $\xi$ is an absolute, hence normal, epimorphism. For it admits a section induced by $\iota_2\from X\to X+X$. This proves that $\rho_X$ is a normal epimorphism.
\end{proof}

It follows that any subtractive variety (\ref{def:SubtractiveVariety}) is a subtractive category in the sense of (\ref{def:SubtractiveCat}).

For the converse of (\ref{thm:DPN-Then-Subtractive}) we may essentially repeat proof of (\ref{thm:SubtractiveVariety:DPN}). This illustrates how in this context, certain varietal proofs involving elements admit a direct translation into a categorical proof.

\begin{proposition}[Subtractivity implies \DPNInline\NTag]
  \label{thm:DPN-for-Sub-Cats}%
  Any normal subtractive category satisfies the \DPNInline-property.
\end{proposition}
\begin{proof}
  We show that condition (\ref{thm:DPNPNvs3x3-Right3x3}) of Proposition~\ref{thm:Dinversion-PreservationNormalMaps-Border(3x3)} holds. We use the notation of \eqref{eq:3x3}. First, we see that the bottom right square is a pushout by (\ref{thm:PushoutRecognize-Categorical}). So, $z$ is a normal epimorphism, and it remains to show that $w$ is the kernel of $z$. To see this, note first that $z\Comp w=\ZeroMap$ by commutativity.

  First, we show that $w$ is a monomorphism or, equivalently, $\Ker{w}=\ZeroMap$; see (\ref{thm:Mono<->0-Kernel}). Consider $i\from Y_1\to I$ such that $w\Comp i=\ZeroMap$. Take the pullback $l\from Y_2\to L$ of $i$ along $d$, so that $d\Comp l=i\Comp p_1$ where the normal epimorphism $p_1\from Y_2\to Y_1$ is the pullback of $d$ along $i$. Then $e\Comp v\Comp l=\ZeroMap$ so that $k\from Y_2\to K$ exists with $b\Comp k=v\Comp l$. Now $cxk=\ZeroMap$ by commutativity, and so $x\Comp k=\ZeroMap$ because $c$ is a monomorphism. So, there is $m\from Y_2\to M$, unique with $u\Comp m=k$. Now $a\Comp m=l$ since $v$ is a monomorphism and $v\Comp a\Comp m=b\Comp u\Comp m=b\Comp k=v\Comp l$. Hence $i\Comp p=d\Comp l=d\Comp a\Comp m=\ZeroMap$. The epic property of $p$ implies that $i=\ZeroMap$.

  Now let $q\from Z_1\to Q$ be such that $z\Comp q=\ZeroMap$. Take the pullback $\bar{q}\from Z_2\to X$ of $q$ along~$e$, so that $e\Comp \bar{q}=q\Comp \bar{e}$, where the normal epimorphism $\bar{e}\from Z_2\to Z_1$ is the pullback of $e$ along~$q$. Then $f\Comp y\Comp \bar{q}=z\Comp e\Comp \bar{q}=z\Comp q\Comp \bar{e}=\ZeroMap$. So, there is $j\from Z_2\to J$, unique with $c\Comp j=y\Comp \bar{q}$. Pull back $j$ along $x$ to a map $\bar{j}\from Z_3\to K$ such that $x\Comp \bar{j}=j\bar{x}$ where $\bar{x}\from Z_3\to Z_2$ is the normal epic pullback of $x$ along $j$. Then $\bar{q}\Comp \bar{x}-b\Comp \bar{j}\from D(Z_3)\to X$ is such that $y\Comp (\bar{q}\Comp \bar{x } -b\Comp \bar{j})=\ZeroMap\from D(Z_3)\to R$. So, there is $l\from D(Z_3)\to L$, unique with $v\Comp l=\bar{q}\Comp \bar{x}-b\Comp \bar{j}$. Now,
  \begin{equation*}
    w\Comp d\Comp l=e\Comp v\Comp l=e\Comp (\bar{q}\Comp \bar{x}-b\Comp \bar{j}) = q \bar{e}-\ZeroMap =q\Comp \bar{e}\Comp(\IdMapOn{Z_{2}}-\ZeroMap)=q\Comp \bar{e}\Comp\rho_{Z_2}
  \end{equation*}
  so that---by orthogonality of the normal epimorphism $r_1\Comp\rho_{Z_2}$ and the monomorphism $w$, it is here that subtractivity is used---a unique morphism $\hat{q}\from Z_1\to I$ exists which satisfies $\hat{q}\Comp \bar{e}\Comp\rho_{Z_2}= d\Comp l$ and the needed $w\Comp \hat{q}=q$. - Thus $w=\KerMap{z}$.
\end{proof}

\begin{subordinate}{}
  \begin{subsubordinate}{Scope of regular categories}
    Regular categories provide a convenient environment for developing a calculus of internal relations or for the study of internal logic. It is part of the definition of a \emph{Barr exact} category, and as such one of the axioms defining an \emph{elementary topos}. Furthermore, all abelian categories and all varieties in the sense of universal algebra are regular categories (since they are actually Barr exact).
  \end{subsubordinate}

  \begin{subsubordinate}{Subtractive varieties}
    Normal subtractive varieties may be characterized as varieties which are normal subtractive categories. We keep this characterization (which combines results of~\cite{ZJanelidze-Snake} and \cite{ZJanelidze-Subtractive}) for a later version of the text which will include a chapter on varieties of algebras.
  \end{subsubordinate}

  \begin{subsubordinate}{Internal actions}
    Some of the exact sequences in (\ref{thm:DPN-Then-Subtractive}) play an important role in the theory of internal actions in the context of semiabelian categories~\cite{BJK,Bourn-Janelidze:Semidirect}---more about this subject in a later version of this document.
  \end{subsubordinate}

  \begin{subsubordinate}{Ideal-Determined Categories}\label{IdealDetermined}
    A category is \Defn{ideal determined} (\cite{JMU}, see also \cite{MM-NC}) if it is normal and di-exact. In particular, such a category satisfies the \DPNInline-property. So it is subtractive by (\ref{thm:DPN-Then-Subtractive}). An explicit characterization of ideal-determined varieties appears in \cite{JMU}.
  \end{subsubordinate}

  \begin{subsubordinate}{Formal difference of maps in an abelian category}
    Example \ref{exa:Difference-AbGroup-Maps} generalizes to abelian categories; see Exercise (\ref{exe:FormalDifferenceInAb}).
  \end{subsubordinate}

\end{subordinate}
\chapter{Homological Categories}
\label{chap:HomologicalCats}

In Section~\ref{sec:HomologicalCats-Axioms}, we introduce the concept of \Defn{homological category} which first appeared in work of Borceux and Bourn in \cite{FBorceuxDBourn2004}. Homological categories and normal categories from Chapter~\ref{chap:NormalCategories} are closely related: In the collection of structural axioms characterizing a normal category, we replace the \AENInline-condition, every split epimorphism is a normal map, by the assumption that, for every split epimorphism $p\from X\to Y$ with section $s\from Y\to X$, the morphisms $s$ and $\kappa=\KerMap{p}$ are jointly extremal-epimorphic:
\begin{equation*}
  \xymatrix@C=4em{
  \Ker{p} \ar@{{ |>}->}@[blue][r]^-{\color{blue} \kappa} &
  X \ar@<-.5ex>[r]_-{p} &
  Y \ar@[blue]@<-.5ex>[l]_-{\color{blue} s}
  }
\end{equation*}
We interpret this as saying that `kernel and section of $p$ generate $X$', and refer to this property as the \KSGInline-condition.

The \KSGInline-condition is extracted from \emph{Bourn protomodular} categories~\cite{DBourn1991}. It implies the \AENInline-condition, and it enables us to prove results which need not hold in a normal category, such as the (Short) $5$-Lemma, the `middle case' of the $(\Prdct{3}{3})$-Lemma, and more as presented in Sections \ref{sec:HomologicalCats-Axioms} to  \ref{sec:3x3-Lemma-Homological}:

\begin{ulist}
  \item Homological categories are normal; hence all of the epimorphism types `normal', `effective', `regular', `strong', and `extremal' coincide, and the Snake Lemma holds.
  \item Every morphism $f$ admits a normal epi factorization; that is $f=me$, where $e$ is a normal epimorphism, and $m$ is a monomorphism (not necessarily normal).
  \item The Short $5$-Lemma and the $5$-Lemma are true.
  \item The $(\Prdct{3}{3})$-Lemma holds in its full generality, with a formulation which matches the one familiar from abelian categories.
\end{ulist}

See Section~\ref{sec:Protomodular-SEpi(X)->X} for further insights on the original viewpoint on protomodularity in terms of the so-called \emph{fibration of points}.

\bigskip

\begin{center}
  \textbf{Leitfaden for Chapter \ref{chap:HomologicalCats}}
\end{center}

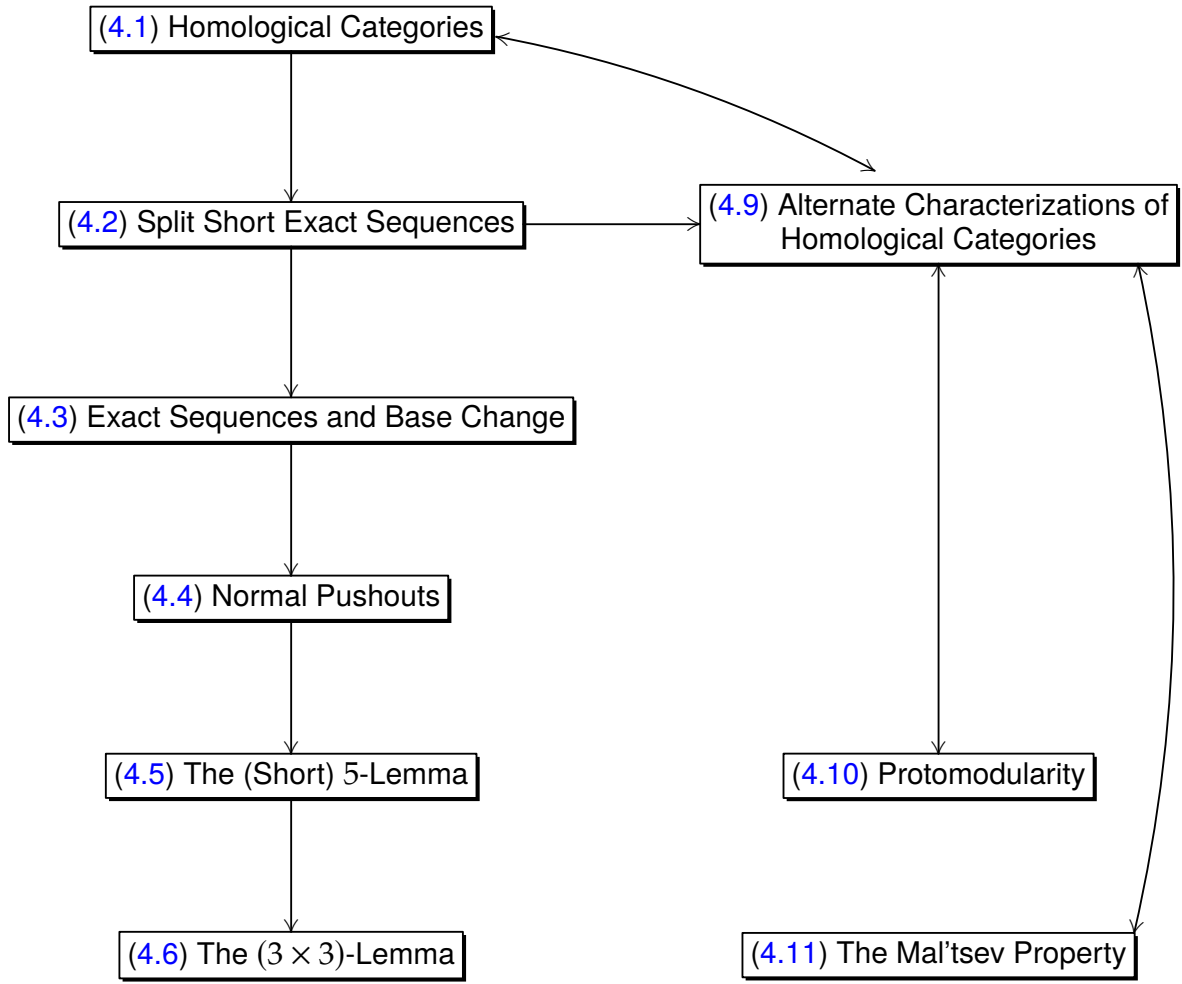
\bigskip

\begin{equation*}
  \xymatrix@R=9ex@C=4em{
  *+[F-,]{\txt{\sffamily (\ref{sec:HomologicalCats-Axioms}) Homological Categories}} \ar@{<->}@<1ex>@/^3ex/[rd] \ar[d] \\
  *+[F-,]{\txt{\sffamily (\ref{sec:SplitShortExactSequences}) Split Short Exact Sequences}} \ar[d] \ar[r] &
  *+[F-,]{\txt{\sffamily (\ref{sec:HomologicalCats-AlternateAxioms}) Alternate Characterizations of \\ \sffamily Homological Categories}} \ar@{<->}[ddd] \\
  *+[F-,]{\txt{\sffamily (\ref{sec:ExactSequences-BaseChange}) Exact Sequences and Base Change}} \ar[d] \\
  *+[F-,]{\txt{\sffamily (\ref{sec:NormalPushouts-Homological}) Normal Pushouts}} \ar[d] \\
  *+[F-,]{\txt{\sffamily (\ref{sec:(Short)5-Lemma}) The (Short) $5$-Lemma}} \ar[d] &
  *+[F-,]{\txt{\sffamily (\ref{sec:Protomodular-SEpi(X)->X}) Protomodularity}} \\
  *+[F-,]{\txt{\sffamily (\ref{sec:3x3-Lemma-Homological}) The $(\Prdct{3}{3})$-Lemma}} &
  *+[F-,]{\txt{\sffamily (\ref{sec:Maltsev}) The Mal'tsev Property}} \ar@{<->}@/_3ex/@<-6em>[uuuu]
  }
\end{equation*}
\newpage
\section[Homological Categories]{Axioms for Homological Categories}
\label{sec:HomologicalCats-Axioms}%

Originally, in~\cite{FBorceuxDBourn2004}, Borceux and Bourn defined a pointed category as being homological if it is regular and protomodular. Here, building on the foundation developed earlier, we arrive at the following presentation of their work:

\begin{definition}[Homological category]
  \label{def:HomologicalCategory}%
  A category $\Ctgry{X}$ is \Defn{homological} if it satisfies the following structural axioms. %
  \index{category!homological}\index{homological!category}%
  \index[acr]{h!{\color{Brown} $\EuRoman{H}$}\IndSep homological category}
  \begin{ulist}
    \item \label{SA-Ax:Pointed}%
    \emph{Zero object}: $\Ctgry{X}$ has a zero object; see (\ref{def:0-Object}).
    \item \label{SA-Ax:FinitelyComplete}%
    \emph{Functorially finitely complete}: For every functor $F\from J\to \Ctgry{X}$ whose domain is a finite category $J$, the limit $\LimOf{F}$ exists; see (\ref{sec:Limits-CoLimits}).
    \item \label{SA-Ax:FinitelyCoComplete}%
    \emph{Functorially finitely cocomplete}: For every functor $F\from J\to \Ctgry{X}$ whose domain is a finite category $J$, the colimit $\CoLimOf{F}$ exists; see (\ref{sec:Limits-CoLimits}).
    \item \label{SA-Ax:NormalEpis-BaseChangePreserved}%
    \emph{\PNEInline\ Pullbacks preserve normal epimorphisms}: The pullback $\bar{g}$ of a normal epimorphism $g$ along any morphism $f$ in $\Ctgry{X}$ is a normal epimorphism. %
    \index[acr]{p!\PNEInline\IndSep condition that pullbacks preserve normal epimorphisms}
    \begin{equation*}
      \xymatrix@R=5ex@C=4em{
      P \ar@{-{ >>}}[d]_-{\bar{g}} \ar[r] \PullLU{rd} &
      Z \ar@{-{ >>}}[d]^-{g} \\
      X \ar[r]_-{f} &
      Y
      }
    \end{equation*}
    \item \label{SA-Ax:KSG}%
    \emph{\KSGInline\ Kernel and section generate:} For a morphism $q\from X\to Q$ with any given section $x\from {Q\to X}$ (so that $qx=\OneMapOn{Q}$), %
    \index[acr]{k!\KSGInline\IndSep for $q\from X\to Q$ sectioned by $x$, images of $x$ and $\KerMap{q}$ generate $X$}%
    \begin{equation}
      \label{eq:KS}%
      \vcenter{
      \xymatrix@C=4em{
      \Ker{q} \ar@{{ |>}->}@[blue][r]^-{\color{blue} k} &
      X \ar@<0.5ex>[r]^-{q} &
      Q \ar@[blue]@<0.5ex>[l]^-{\color{blue} x}
      }
      }
    \end{equation}
    the morphisms $x$ and $k=\KerMap{q}$ are jointly extremal-epimorphic; see (\ref{def:ExtremalEpimorphism}).
  \end{ulist}
\end{definition}

\emph{Comparison with normal categories}\quad In transitioning from normal categories to homological ones, we swap out the condition \AENInline---absolute epimorphisms are normal---for the condition \KSGInline. In (\ref{thm:SplitEpi->Normal-Homological}), we show that every homological category is a normal category. Beyond this fact, the role of the \KSGInline-condition is not immediately obvious. Yet, it is far reaching. Given a short exact sequence, it is the key toward gaining information about its middle object from its end objects. Indeed, we shall prove (\ref{thm:Hofmann}) that, in the context of a normal category, the \KSGInline-condition is equivalent to the validity of the Short $5$-Lemma; see Theorem~\ref{thm:Short5}.

To see the difficulty one encounters if the \KSGInline-condition is not satisfied, consider the category $\SetsBsd$ of pointed sets. Here we encounter the phenomenon that a sectioned epimorphism $q$ as in Diagram \ref{eq:KS} need not be the cokernel of its kernel $k$, because $K$ combined with any section of $q$ fail to generate $X$. To see this in $\SetsBsd$, take a `huge' pointed set $(X,t)$ and let $q\from (X,t)\to (\Set{-1, 1},1)$ send $a\neq t$ to $-1$. Then $K\DefEq\Ker{q}=\Set{t}$, and so $\CoKerMap{\KerMap{q}} = \IdMapOn{X}\neq q$.

The following proposition implies that every  homological category is also a normal category.

\begin{proposition}[Split epimorphism is normal\HTag]
  \label{thm:SplitEpi->Normal-Homological}%
  In a \ZExact\ category which satisfies the \KSGInline-condition, every split epimorphism is a normal epimorphism.
\end{proposition}
\begin{proof}
  Given a split epimorphism $p$ with its kernel $k$ as in Diagram \ref{eq:KS}, we check that $p$ is a cokernel of $k$. Let $f\from X\to Z$ be such that $fk=0$. Then $fs\from Y\to Z$ is a factorization of $f$ through $p$. Indeed, $fsp=f$ because $k$ and $s$ are jointly epimorphic and $fsps=fs$ while $fspk=0=fk$. Since $p$ is an epimorphism, this factorization is unique.
\end{proof}

\begin{corollary}[Homological $\Rightarrow$ normal\HTag]
  \label{thm:HomologicalCat->NormalCat}%
  Every homological category $\Ctgry{X}$ enjoys the following properties.
  \begin{thmlist}
    \item $\Ctgry{X}$ is a normal category; see (\ref{def:NormalCat}).
    \item \HSDInline\ $\Ctgry{X}$ is homologically self-dual.
  \end{thmlist}
\end{corollary}
\begin{proof}
  (i)\quad With (\ref{thm:SplitEpi->Normal-Homological}) we see that the structural axioms (\ref{def:NormalCat}) characterizing a normal category are satisfied. This implies (ii) because we know already that every normal category is homologically self-dual;  (\ref{thm:NormalCat->HomologicallySelfDual}).
\end{proof}

\begin{example}[Normal but not homological]\label{exa:Normal-Not-Homological}
  The category $\OrdGrp$ of preordered groups is known to be normal but not homological~\cite{MMClementinoNMartinsFerreiraAMontoli2019-Preordered}: a preordered group being a group equipped with a preorder (a reflexive and transitive relation) for which the group operation is monotone; arrows are monotone group homomorphisms.
\end{example}

\begin{subordinate}{Borceux-Bourn homological vs homological}
  \begin{subsubordinate}{Finite (co)completeness}
    We require a homological category to be finitely bicomplete. This means that our definition of `homological category' differs slightly from earlier definitions such as the one in~\cite{FBorceuxDBourn2004}, where the only colimits required to exist are coequalizers of kernel pairs. More precisely, whenever a morphism $f\from X\to Y$ admits a kernel pair, the coequalizer of the kernel pair projections must exist. While it is clear that this condition holds in a finitely cocomplete category, there are categories in which the converse does not hold. For example, these two conditions do not guarantee the existence of cokernels or coproducts. This means that a category which merely admits this type of colimits need not even be \ZExact. Requiring a homological category to be finitely cocomplete is a matter of convenience.

    For clarity, we say that a category is \Defn{Borceux--Bourn homological} (or \Defn{BB-homological}) if it is pointed, finitely complete, admits coequalizers of kernel pairs, preserves normal epimorphisms under base change, and has the \KSGInline-property. We explain in (\ref{sec:BBHom}) that this agrees with the original definition of Borceux and Bourn.

    As we shall see in Chapter~\ref{chap:SACats}, in the context of a semiabelian category, this problem no longer occurs. Indeed, if a BB-homological category satisfies \ANNInline, then existence of binary coproducts entails finite cocompleteness. This is a non-trivial result\footnote{This result is omitted in the current rendition of the text, but will be available in a later version.}.
  \end{subsubordinate}
\end{subordinate}

\begin{exercises}
\begin{exercise}[Torsion-free groups]
  \label{exe:TorsionFreeGrps->Homological}%
  Show that the category $\GrpsTF$ of torsion-free groups has the following properties: %
  \index{homological category!torsion-free groups}\index{torsion-free!groups form homological cat}%
  \begin{thmlist}
    \item Given a group $G$, its torsion subgroup is $T(G)$ is a normal subgroup.
    \item The quotient map $G\NEpi G/T(G)$ determines a functor $\Grps\to \GrpsTF$ to the effect that $\GrpsTF$ is a reflective subcategory of $\Grps$.
    \item The kernel of a morphism $f\from G\to H$ of torsion free groups is its kernel in $\Grps$; i.e. the set of all those $g\in G$ with $f(g)=\ZeroObject$.
    \item The cokernel of a morphism $f\from G\to H$ of torsion-free groups is given by its cokernel in $\Grps$ modulo its torsion subgroup.
    \item $\GrpsTF$ is complete and limits may be computed in $\Grps$. - Conclude that the \PNEInline-property holds in $\GrpsTF$.
    \item The \KSGInline-property holds in $\GrpsTF$.
  \end{thmlist}
  Conclude that $\GrpsTF$ is a homological category.
\end{exercise}

\begin{exercise}[Torsion-free abelian groups]
  \label{exe:TorsionFreeAbGroups->Homological}%
  Show that the category $\AbGrpsTF$ is a homological category. %
  \index{homological category!torsion-free abelian groups}\index{torsion-free!abelian groups form homological cat}%
\end{exercise}
\end{exercises}
\newpage
\section[Split Short Exact Sequences]{Split Short Exact Sequences}
\label{sec:SplitShortExactSequences}

We know from (\ref{thm:SplitEpi->Normal-Homological}) that, in a homological category, every sectioned epimorphism is the cokernel of its kernel. Thus, as in a normal category, every sectioned epimorphism determines a split short exact sequence. Via the \KSGInline-property, kernel and section of such a split short exact sequence generate the middle object. As consequence, we show in (\ref{thm:SplitShort5}) that  a morphism of split short exact sequences is an isomorphism if and only if it restricts to an isomorphism of the end objects. This result is called the \Defn{Split Short $5$-Lemma}.

Using the Split Short $5$-Lemma, we establish a converse to (\ref{thm:Pullback->IsoOfKernels}). There, we showed that, in a pullback square, the kernels of two opposing maps are canonically isomorphic. Here, we show that, if a morphism of split short exact sequence involves an isomorphism of kernel objects, then the square of cokernels is a pullback; see (\ref{thm:PullbackRecognition-SplitSES}).

Next, we use this information to identify a product of two objects in a homological category. In (\ref{thm:ProductRecognition}) we show that an object $T$ admits the structure of a product of $A$ and $B$ if and only if $T$ fits into mutually split short exact sequences:
\begin{equation*}
  \xymatrix@R=5ex@C=4em{
  A \ar@{{ |>}->}@<+.5ex>[r]^-{i} &
  T \ar@{-{ >>}}@<+.5ex>[r]^-{p} \ar@{-{ >>}}@<+.5ex>[l]^-{q} &
  B \ar@{{ |>}->}@<+.5ex>[l]^-{j}
  }
\end{equation*}
As another application we show in (\ref{thm:Sum->ProductIsCokernel}) that the canonical map $A+B\to \Prdct{A}{B}$ is a normal epimorphism.

\begin{definition}[Morphism of sectioned/split short exact sequences\ZExactTag]
  \label{def:SplitSES-Morphism}%
  A \Defn{morphism of split short exact sequences} is given by a diagram whose rows are split short exact sequences: %
  \index{morphism!split short exact sequences}\index{split!short exact sequence - morphism}%
  \begin{equation}
    \label{fig:Morphism-SectionedSESs}
    \vcenter{\xymatrix@R=5ex@C=4em{
    K \ar@{{ |>}->}[r]^-{k} \ar[d]_-{\kappa} &
    X \ar@{-{ >>}}@<.5ex>[r]^-{q} \ar[d]_-{\xi} &
    Q \ar[d]^-{\rho} \ar@<.5ex>[l]^-{x} \\
    L \ar@{{ |>}->}[r]_-{l} &
    Y \ar@{-{ >>}}@<.5ex>[r]^-{r} &
    R \ar@<.5ex>[l]^-{y}
    }
    }
  \end{equation}
  such that the left hand square commutes, $\rho q=r\xi$, and $\xi x=y\rho$.
\end{definition}

We write $\SSESCat{X}$ for the category of split short exact sequences and their morphisms. %
\index[not]{s!$\SSESCat{X}$\IndSep category of split short exact sequences}%

Here is an example which demonstrates what can happen in a \ZExact\ category  which fails to satisfy the \KSGInline-Axiom.
\newpage
\begin{example}[Projection $\NNr\prdct\NNr \to \NNr$]
  \label{exa:NxN-project->N}
  In the category $\Monoids$ of monoids, consider the split short exact sequence
  \begin{equation*}
    \xymatrix@R=5ex@C=2em{
    \NNr \ar@{{ |>}->}[rr]^-{k=\PrdctMapInto{\IdMap,0}} &&
    \Prdct{\NNr}{\NNr} \ar@{-{ >>}}@<.5ex>[rr]^-{\PrjctnOnto{1}} &&
    \NNr \ar@<.5ex>[ll]^-{s=\Dgnl}
    }
  \end{equation*}
  where $\PrjctnOnto{1}(m,n)\DefEq m$, and $s(n)\DefEq (n,n)$. Still, $k$ and $s$ fail to generate $\Prdct{\NNr}{\NNr}$. This corresponds to the existence of factorizations of $k$ and $s$ through proper subobjects of $\Prdct{\NNr}{\NNr}$ such as, for $c\geq 2$ fixed,
  \begin{equation*}
    M = \SetSlct{(m,0)+(n,n) + (0,pc)}{m,n,p\in\NNr}
  \end{equation*}
\end{example}

\begin{example}[Split short exact sequence in $\SetsBsd$ vs.\ \KSGInline]
  More subtle yet is what happens in the category $\SetsBsd$ of pointed sets. On the one hand, in every \emph{split short exact sequence}, kernel and section generate the middle object. However, absolute epimorphisms need not be normal maps. So, the \KSGInline-property is not satisfied in $\SetsBsd$; see Exercise~\ref{exe:KSG-In-Set_*}.
\end{example}

\begin{proposition}[Split Short 5-Lemma\HTag]
  \label{thm:SplitShort5}%
  Consider a morphism of split short exact sequences as in \eqref{fig:Morphism-SectionedSESs}.  If $k$ and $\rho$ are isomorphisms, then $\xi$ is an isomorphism as well.
\end{proposition}
\begin{proof}
  In the given diagram, consider the diagram of kernels of $\kappa$, $\xi$, and $\rho$:
  \begin{equation*}
    \xymatrix@R=5ex@C=3.5em{
    \ZeroObject \ar@{{ |>}->}[r] &
    \Ker{\xi} \ar@<.5ex>@{-{ >>}}[r]^-{s} &
    \ZeroObject \ar@<.5ex>[l]
    }
  \end{equation*}
  We find that $\Ker{s}=\Ker{\kappa}=\ZeroObject$ because limits commute with limits. With (\ref{exa:Extension-By-0's}), we see that $\Ker{\xi}=\ZeroObject$.  So, the map $\xi$ is a monomorphism by Proposition~\ref{thm:Mono<->0-Kernel}. Now $l$ and $y$ factor through this monomorphism $\xi$ via $l'\kappa^{-1}$ and $x\rho^{-1}$. Since $l$ and $y$ are jointly extremal-epimorphic, we infer with the \KSGInline-property that $\xi$ is an isomorphism.
\end{proof}

\begin{proposition}[Pullback of sectioned epimorphism\HTag]
  \label{thm:PullbackOfSplitEpi-II}
  \label{thm:proto-object}	
  Consider the pullback of a sectioned epimorphism $p\from X\to Y$, with a section $s\from Y\to X$, along an arbitrary $f\from A\to Y$, as in (\ref{thm:PullbackOfSplitEpi}):
  \begin{equation*}
    \xymatrix@R=5ex@C=4em{
    P \ar[r]^-{\bar{f}} \ar@{-{ >>}}@<-.5ex>[d]_-{\bar{p}} \PullLU{rd} &
    X \ar@{-{ >>}}@<-.5ex>[d]_-{p} \\
    A \ar@<-.5ex>[u]_-{\overline{s}} \ar[r]_-{f} &
    Y \ar@<-.5ex>[u]_-s
    }
  \end{equation*}
  Then the morphisms $\bar{f}$ and $s$ are jointly extremal-epimorphic.
\end{proposition}
\begin{proof}
  In the commutative diagram below, assume that $m\from M\to X$ is a monomorphism.
  \begin{equation*}
    \xymatrix@R=5ex@C=3em{
    & M \ar@{ >->}[d]^- {m} \\
    P \ar[r]_-{\bar{f}} \ar[ur]^-{r'} &
    X &
    Y \ar[l]^-s \ar[ul]_-{s'}
    }
  \end{equation*}
  We want to show that $m$ is an isomorphism. To see this, note first that $pm\from {M\to Y}$ satisfies $ pms'=ps=\IdMapOn{Y}$. So $pm$ admits a section $s'$, which makes it an absolute epimorphism. Now expand the original diagram by taking kernels:
  \begin{equation*}
    \xymatrix@C=4em@R=6ex{
    \Ker{\bar{p}} \ar@{{ |>}->}[d] \ar[r]^{\tilde{f}}_-{\cong} &
    \Ker{p} \ar@{{ |>}->}[d] &
    \Ker{pm} \ar@{<-}`u[l]`[ll]_-{\tilde{r}'}[ll] \ar[l]_{\widetilde{m}} \ar@{{ |>}->}[d] \PullRU{ld} \\
    P \ar[r]^-{\bar{f}} \ar[d]_-{\bar{p}} \PullLU{rd} &
    X \ar@<-.5ex>[d]_-{p} &
    M \ar[l]_-{m} \ar@<-.5ex>[d]_-{pm} \\
    A \ar[r]_-{f} &
    Y \ar@<-.5ex>[u]_-s &
    Y \ar@<-.5ex>[u]_-{s'} \ar@{=}[l]
    }
  \end{equation*}
  The map $ \tilde{f}$ is an isomorphism by (\ref{thm:KernelFunctor-Props}). The top right square is a pullback \eqref{thm:PullbackRecognition-KernelSide-1}, and so pullback stability of monomorphisms (\ref{exe:BaseChange-Mono-Epi-Iso}) shows that $\widetilde{m}$ is a monomorphism. The morphism $\tilde{r}'$ induced by $r'$ via the top right pullback, however, is such that $\widetilde{m}\Comp\tilde{r}' =\tilde{f}$. Hence $\widetilde{m}$ is both a monomorphism and an absolute epimorphism, which makes it an isomorphism; see (\ref{exe:SectionableEpis}). With the split Short $5$-Lemma (\ref{thm:SplitShort5}), we see that $m$ is an isomorphism.
\end{proof}

Any pullback of a sectioned epimorphism along a sectioned epimorphism yields a sectioned epimorphism in the category of sectioned epimorphisms in $\Ctgry{X}$, as in (\ref{thm:AbsolutePush/Pull-SectionedMorInSEpi(X)}):
\stepcounter{theorem}
\begin{equation}\label{Diag:DoubleSplitEpi}
  \vcenter{
  \xymatrix@R=5ex@C=4em{
  P \ar@<-.5ex>@{-{ >>}}[r]_-{\bar{f}} \ar@{-{ >>}}@<-.5ex>[d]_-{\bar{p}} &
  X \ar@{-{ >>}}@<-.5ex>[d]_-{p} \ar@<-.5ex>@{{ >}->}[l]_-{\bar{t}} \\
  A \ar@<-.5ex>@{-{ >>}}[r]_-{f} \ar@<-.5ex>@{{ >}->}[u]_-{\bar{s}} &
  Y \ar@<-.5ex>@{{ >}->}[u]_-{s} \ar@<-.5ex>@{{ >}->}[l]_-{t}
  } 
  }
\end{equation}
This means that $\bar{f}\Comp \bar{t} = \IdMapOn{X}$, $p\Comp s= \IdMapOn{Y}$, $f\Comp t=\IdMapOn{Y}$ and $\bar{p}\Comp \bar{s}=\IdMapOn{A}$, and that the two possible composites from any corner to its opposite corner are equal. We find the next result, which as explained in (\ref{sec:Maltsev}) is related to the concept of a Mal'tsev category.

\begin{corollary}[Mal'tsev property\HTag]
  \label{thm:ProtoMaltsev}
  In a pointed, finitely complete category satisfying the \KSGInline-axiom, suppose diagram \eqref{Diag:DoubleSplitEpi} is a pullback of $p\from X\to Y$ along a morphism $f\from A\to Y$. Then, the sections $\bar{s}$ and $\bar{t}$ are jointly extremal-epimorphic.
\end{corollary}
\begin{proof}
  Via the \KSGInline-axiom, the maps $\bar{s}$ and $\Ker{\bar{p}}\to P$ are jointly extremally epimorphic. But the map $\Ker{\bar{p}}\cong\Ker{p} \to P$ factors through $X$ via $\bar{t}$. So, $\bar{s}$ and $\bar{t}$ are jointly extremally epimorphic.
\end{proof}

\begin{lemma}[Pullback recognition, split short exact version\HTag]
  \label{thm:PullbackRecognition-SplitSES}%
  Consider a morphism of split short exact sequences.
  \begin{equation*}
    \xymatrix@R=5ex@C=3em{
    K' \ar@{{ |>}->}[r]^-{\kappa'} \ar[d]_-{k}^{\cong} &
    X' \ar@<.5ex>@{-{ >>}}[r]^-{p'} \ar[d]_-{x} &
    Y' \ar[d]^-{y} \ar@<.5ex>[l]^-{s'} \\
    K \ar@{{ |>}->}[r]_-{\kappa} &
    X \ar@<.5ex>@{-{ >>}}[r]^-{p} &
    Y \ar@<.5ex>[l]^-{s}
    }
  \end{equation*}
  Then $k$ is an isomorphism if and only if the square $X'\rightrightarrows Y$ is a pullback.
\end{lemma}
\begin{proof}
  If the square on the right is a pullback, then $k$ is an isomorphism by (\ref{thm:Pullback->IsoOfKernels}). The converse follows via the Split Short 5-Lemma (\ref{thm:SplitShort5}): Pulling back $p$ along $y$ yields the commutative diagram below.
  \begin{equation*}
    \xymatrix@R=5ex@C=3em{
    K' \ar@{{ |>}->}[rr]^-{\kappa'} \ar[dd]_-{k}^{\cong} \ar[rd]_{\cong}^{k} &&
    X' \ar@<.5ex>@{-{>}}[rr]^-{p'} \ar[dd]_(.3){x} \ar@{.>}[rd]^{t} &&
    Y' \ar[dd]_(.25){y} \ar@<.5ex>[ll]^-{s'} \\
    & K \ar@{{ |>}->}[rr]|(.43)\hole \ar@{=}[dl] &&
    X\times_Y Y' \ar@<+.5ex>[rr]|(.57)\hole^(.35){\bar{p}} \ar[ld] &&
    Y' \ar@{=}[lu] \ar@{.>}@<.5ex>[ll]|(.44)\hole^(.65){\bar{s}} \ar[ld]^{y} \\
    K \ar@{{ |>}->}[rr]_-{\kappa} &&
    X \ar@<.5ex>@{-{>}}[rr]^-{p} &&
    Y \ar@<.5ex>[ll]^-{s}}
  \end{equation*}
  The universal property of $ X\times_Y Y'$ yields the comparison map $t$, and we claim that $t$ is an isomorphism. Indeed, $ \Ker{p}=K\cong\Ker{\bar{p}}$\*. A section $\bar{s}$ of $\bar{p}$ exists via the pullback property of the bottom right square. Now the Split Short 5-Lemma tells us that $t$ is an isomorphism. - This implies the claim.
\end{proof}

\begin{proposition}[Product properties\HTag]
  \label{thm:ProductRecognition}%
  In any homological category the following hold. %
  \index{product!recognition}
  \begin{enumerate}[(i)]
    \item \label{thm:ProductRecognition-ProjectionIsCoKer}%
          Any product projection $ \PrjctnOnto{A}\from \Prdct{A}{B}\to A$ is a normal epimorphism.
    \item \label{thm:ProductRecognition-2SESs}%
          The product of two objects $A$ and $B$ fits into two short exact sequences:
          \begin{equation*}
            \xymatrix@R=5ex@C=4em{
            A \ar@{{ |>}->}@<+.5ex>[r]^-{\PrdctMapInto{\IdMapOn{A},\ZeroMap}} &
            A\times B \ar@{-{ >>}}@<+.5ex>[r]^-{\PrjctnOnto{B}} \ar@{-{ >>}}@<+.5ex>[l]^-{\PrjctnOnto{A}} &
            B \ar@{{ |>}->}@<+.5ex>[l]^-{\PrdctMapInto{\ZeroMap , \IdMapOn{B}}}
            }
          \end{equation*}
    \item \label{thm:ProductRecognition-ViaDoubleSplit}%
          Two maps $q\from T\to A$ and $p\from T\to B$ represent $T$ as a product of $A$ and $B$ if and only if $p$ and $q$, together with their kernels, yield mutually sectioning short exact sequences such as
          \begin{equation*}
            \xymatrix@R=5ex@C=4em{
            A \ar@{{ |>}->}@<+.5ex>[r]^-{i} &
            T \ar@{-{ >>}}@<+.5ex>[r]^-{p} \ar@{-{ >>}}@<+.5ex>[l]^-{q} &
            B \ar@{{ |>}->}@<+.5ex>[l]^-{j}
            }
          \end{equation*}
          If so, then $(q,p)\from T\to \Prdct{A}{B}$ is an isomorphism.
  \end{enumerate}
\end{proposition}
\begin{proof}
  (\ref{thm:ProductRecognition-ProjectionIsCoKer})\quad $\PrjctnOnto{A}$ has a section $\PrdctMapInto{\ZeroMap,\IdMapOn{B}}$. So, the claim  follows from (\ref{thm:SplitEpi->Normal-Homological}).

  (\ref{thm:ProductRecognition-2SESs})\quad To see this, we show that $\KerMap{\PrjctnOnto{A}}=\PrdctMapInto{\ZeroMap , \IdMapOn{B}}$. Consider the commutative diagram below.
  \begin{equation*}
    \xymatrix@!0@R=10ex@C=6em{
    B \ar[r]^-{\PrdctMapInto{\ZeroMap , \IdMapOn{B}} } \ar[d] &
    \Prdct{A}{B} \ar[d]_-{\PrjctnOnto{A}} \ar[r]^-{\PrjctnOnto{B}} \ar[d] &
    B \ar[d] \ar@{<-}`u[l]`[ll]_-{1_{B}}[ll] \PullLU{ld} \\
    0 \ar[r] &
    A \ar[r] &
    0
    }
  \end{equation*}
  The square on the right is a pullback, as is the outer rectangle. By pullback cancellation \eqref{thm:Pullbacks,ConcatenatedSquares}, the left hand square is a pullback as well. With (\ref{thm:Ker/CoKerExist}) we see that $\PrdctMapInto{0,\IdMapOn{B}}=\KerMap{\PrjctnOnto{A}}$.

  (\ref{thm:ProductRecognition-ViaDoubleSplit})\quad If $p$ and $q$, together with their kernels, yield those mutually sectioning short exact sequences, then the diagram below commutes.
  \begin{equation*}
    \xymatrix@R=5ex@C=4em{
    A \ar@{{ |>}->}@<+.5ex>[r]^-{i} \ar@{=}[d] &
    T \ar@{-{ >>}}@<+.5ex>[r]^-{p}  \ar[d]^{(q,p)} \ar@{-{ >>}}@<+.5ex>[l]^-{q} &
    B \ar@{{ |>}->}@<+.5ex>[l]^-{j} \ar@{=}[d] \\
    A \ar@{{ |>}->}@<+.5ex>[r]^-{\PrdctMapInto{\IdMapOn{A},\ZeroMap}} &
    A\times B \ar@{-{ >>}}@<+.5ex>[r]^-{\PrjctnOnto{B}}  \ar@{-{ >>}}@<+.5ex>[l]^-{\PrjctnOnto{A}} &
    B \ar@{{ |>}->}@<+.5ex>[l]^-{\PrdctMapInto{\ZeroMap , \IdMapOn{B}}}
    }
  \end{equation*}
  By the Split Short $5$-Lemma \ref{thm:SplitShort5}, the map $(q,p)$ is an isomorphism. Via the commutativity of the diagram, we see that $(T,q,p)$ is a product cone.

  Conversely, if $(T,q,p)$ is a product cone, then the map $(q,p)$, induced by the product property of $\Prdct{A}{B}$ is an isomorphism.
  \begin{equation*}
    \xymatrix@R=5ex@C=4em{
    A \ar@{=}[d] &
    T \ar[r]^-{p}  \ar[d]^{(q,p)}_{\cong} \ar[l]^-{q} &
    B  \ar@{=}[d] \\
    A \ar@{{ |>}->}@<+.5ex>[r]^-{\PrdctMapInto{\IdMapOn{A},\ZeroMap}} &
    A\times B \ar@{-{ >>}}@<+.5ex>[r]^-{\PrjctnOnto{B}}  \ar@{-{ >>}}@<+.5ex>[l]^-{\PrjctnOnto{A}} &
    B \ar@{{ |>}->}@<+.5ex>[l]^-{\PrdctMapInto{\ZeroMap , \IdMapOn{B}}}
    }
  \end{equation*}
  It follows that $p$ and $q$ are normal epimorphisms. Further, $i\DefEq  (q,p)^{-1}\Comp (\IdMapOn{A},0) = \Ker{p}$ is a section for $q$, and $j\DefEq (q,p)^{-1}\Comp (0,\IdMapOn{B})= \Ker{q}$ is a section for $p$.
\end{proof}

\begin{proposition}[$A+B\to \Prdct{A}{B}$ is a normal epimorphism\HTag]
  \label{thm:Sum->ProductIsCokernel}
  \label{thm:Sum->ProductIsNormalEpi}%
  For objects $A$ and $B$, the \Defn{sum/product} or \Defn{coproduct/product} comparison map
  \index[not]{c!$\SumProdComp{A}{B}$\IndSep comparison map $A+B\longrightarrow \Prdct{A}{B}$}%
  \index{coproduct/product comparison map}\index{sum/product comparison map}%
  \begin{equation*}
    \SumProdComp{A}{B}\DefEq\PrdctMapInto{\SumMapOutOf{\IdMapOn{A},\ZeroMap},\SumMapOutOf{\ZeroMap,\IdMapOn{B}}} = \SumMapOutOf{\PrdctMapInto{\IdMapOn{A},\ZeroMap}, \PrdctMapInto{\ZeroMap,\IdMapOn{B}} } \from A+B \longrightarrow \Prdct{A}{B}
  \end{equation*}
  is a normal epimorphism.
\end{proposition}
\begin{proof}
  The two inclusions $\PrdctMapInto{\IdMapOn{A},\ZeroMap}$ and $ \PrdctMapInto{\ZeroMap,\IdMapOn{B}}$ are jointly extremal-epimorphic by Corollary \ref{thm:ProtoMaltsev}, applied in the situation where $Y=\ZeroObject$. Thus the induced arrow $A+B\to A\times B$ is an extremal epimorphism. So it is a normal epimorphism by (\ref{thm:CoKer=NormalEpi=RegEpi=EffectiveEpi}) .
\end{proof}

Proposition \ref{thm:Sum->ProductIsCokernel} implies that the factor inclusions into the \emph{product} are jointly epimorphic:
\begin{equation*}
  A \XRA{(\IdMapOn{A},\ZeroMap)} \Prdct{A}{B} \XLA{(\ZeroMap,\IdMapOn{B})} B
\end{equation*}

\begin{subordinate}{}
  \begin{subsubordinate}{On the comparison map $A+B\NEpi \Prdct{A}{B}$}
    In (\ref{thm:Sum->ProductIsCokernel}), we showed that the comparison map $A+B\NEpi \Prdct{A}{B}$ is a normal epimorphism. This fact plays a key role in what follows: It provides the foundation for discussing abelian objects in a semiabelian category. It supports the development of a well behaved theory of commutators. Finally, a homological category in which $\SumProdComp{A}{A}$ is an isomorphism for every $A$ is additive.
  \end{subsubordinate}

  \begin{subsubordinate}{Split short exact sequences and semidirect products}
    As in the category $\Grps$ of groups, we will see later that the middle object $X$ in a split short exact sequence $K\NMono X\leftrightarrows Y$ may be described as a semidirect product of $K$ and $Y$.

    Functorial kernels in a homological category $\EuScript{X}$ yield a canonical equivalence between the category $\SSESCat{X}$ of split short exact sequences in $\EuScript{X}$, and the category $\SEpisIn{X}$ of sectioned epimorphisms in $\EuScript{X}$. The latter is also known as the `category of points' in $\EuScript{X}$; see (\ref{sec:AbsoluteDiagrams}) and (\ref{def:Protomodularity}).
  \end{subsubordinate}
\end{subordinate}

\begin{exercises}

\begin{exercise}[KSG in $\SetsBsd$]
  \label{exe:KSG-In-Set_*}
  In the category $\SetsBsd$ of pointed sets show the following:
  \begin{thmlist}
    \item $\SetsBsd$ fails to satisfy the \KSGInline-Axiom.
    \item If $K \XRA{k} X \XRA{p} Y$ is a short exact sequence with splitting $s$ of $p$, then $k$ and $s$ are jointly extremal-epimorphic.
  \end{thmlist}
\end{exercise}

\begin{exercise}[Weak condition for: Sectioned epimorphism is normal\ZExactTag]
  \label{exe:SplitEpi->NormalEpi-Weak}
  In a \ZExact\ category, consider a sectioned epimorphism $p\from {X\to Y}$ with section $s\from {Y\to X}$ and kernel $\kappa=\KerMap{p}$.
  \begin{equation*}
    \xymatrix@R=5ex@C=4em{
    K \ar@{{ |>}->}[r]^-{\kappa} &
    X \ar@{-{>>}}@<.5ex>[r]^-{p} &
    Y \ar@<.5ex>[l]^-s
    }
  \end{equation*}
  If $\kappa$ and $s$ are jointly epimorphic, then $p$ is a cokernel of $\kappa$.
\end{exercise}

\begin{exercise}[Biproduct in a pointed variety\ANKTag]
  \label{rem:BiProduct-Pointed}%
  Is there are non-linear variety with a non-zero object $X$ for which $\SumProdComp{X}{X}\from X+X\to \Prdct{X}{X}$ is an  isomorphism?
\end{exercise}

\begin{exercise}[Double absolute epimorphisms\HTag]
  \label{exe:Di-AbsoluteEpi->AbsolutePush-H}
  First recall (\ref{thm:PullbackOfSplitEpi-II}), then revisit Exercise~\ref{exe:DoubleAbsoluteEpi} in the context of a homological category: is a commutative square consisting of absolute epimorphisms always an absolute pushout?
\end{exercise}
\end{exercises}
\section{Exact Sequences and Base Change}%
\label{sec:ExactSequences-BaseChange}%

In the study of group or module extensions, functoriality of Yoneda's Ext-functor in its first variable relies on pulling a given short exact sequence $K\NMono X\NEpi Q$ back along a morphism $f\from R\to Q$. Thus we obtain a new short exact sequence $K\NMono X'\NEpi R$. By virtue of the \PNEInline-Axiom, we show in (\ref{thm:BaseChangePreserves-SESs}) that this construction is available in any homological category. Thus, given a short exact sequence $\varepsilon$, as in the bottom row of the diagram below, form the pullback of $q$ along an arbitrary morphism $f\from R\to Q$:
\stepcounter{theorem}%
\begin{equation}
  \label{eq:Pullback-SES}%
  \vcenter{
  \xymatrix@R=5ex@C=4em{
  f^{\ast}\varepsilon \ar[d]_{(\IdMap,\bar{f},f)}&
  K \ar@{{ |>}.>}[r]^-{k'} \ar@{:}[d] &
  \bar{X} \PullLU{rd} \ar@{.{ >>}}[r]^-{\bar{q}} \ar@{.>}[d]_{\bar{f}} &
  R \ar[d]^{f} \\
  \varepsilon &
  K \ar@{{ |>}->}[r]_-{k} &
  X \ar@{-{ >>}}[r]_-{q} &
  Q
  }
  }
\end{equation}
We find that the top row is a short exact sequence, denoted $f^{\ast}\varepsilon$. We refer to it as the short exact sequence obtained  by \Defn{pulling $\varepsilon$ back along $f$} or by \Defn{base change along $f$}. %
\index[not]{p!$f^{\ast}\varepsilon$\IndSep pullback of short exact sequence $\varepsilon$ along $f$}%

Combined with the Short 5-Lemma we will be well positioned to define Ext as a contravariant functor in the first variable; see Section~\ref{sec:(Short)5-Lemma}, which relies on technical preparation developed here.

\begin{proposition}[Base change preserves short exact sequences\HTag]
  \label{thm:BaseChangePreserves-SESs}
  If the sequence $K\NMono X\NEpi Q$ is short exact, then a morphism $f\from R\to Q$ yields a morphism of short exact sequences, as in \eqref{eq:Pullback-SES}.
\end{proposition}
\begin{proof}
  The normal epimorphism $q$ pulls back to a normal epimorphism by the \PNEInline-condition. By (\ref{thm:Pullback->IsoOfKernels}) we may represent $\KerMap{\bar{q}}$ by a map $k'$ so that the map of kernels induced by $(\bar{f},f)$ is the identity.
\end{proof}

Next, we present properties of pullbacks, including tools to identify a commutative square as a pullback. We know already that pullbacks preserve monomorphisms (\ref{exe:BaseChange-Mono-Epi-Iso}). Here, we are establish the following converse:

\begin{proposition}[Pullbacks preserve/reflect monomorphisms\HTag]
  \label{thm:PullbacksPreserve/ReflectMonos}%
  In the pullback square
  \begin{equation*}
    \xymatrix@R=5ex@C=4em{
    K \ar[r]^-{\bar{h}} \ar[d]_-{\bar{k}} \PullLU{rd} &
    L \ar[d]^-{k} \\
    X \ar[r]_-{h} & A
    }
  \end{equation*}
  the map $\bar{k}$ is a monomorphism if and only if $k$ is a monomorphism.
\end{proposition}
\begin{proof}
  We know already (\ref{exe:BaseChange-Mono-Epi-Iso}) that base change preserves monomorphisms. For the converse, suppose $\bar{k}$ is a monomorphism. Then $\Ker{k}=\Ker{\bar{k}}=0$ by combining (\ref{thm:KernelFunctor-Props}) and (\ref{thm:Pullback->IsoOfKernels}). The monomorphism recognition criterion (\ref{thm:Mono<->0-Kernel}) implies that $k$ is a monomorphism.
\end{proof}

\begin{corollary}[Pullback along normal epi preserves/reflects (normal) monos\HTag]
  \label{thm:PullbacksAlongCokerPreserve/ReflectKernels}%
  Pullback along a normal epimorphism preserves and reflects monomorphisms and normal monomorphisms. %
  \index{pullback!along normal epi preserves/reflects (normal) monos}%
  \index{base change!along normal epi preserves/reflects (normal) monos}%
\end{corollary}
\begin{proof}
  That base change preserves and reflects monomorphisms holds by (\ref{thm:PullbacksPreserve/ReflectMonos}). It preserves normal monomorphisms by (\ref{thm:PullbackPreservesNormalMonos}). Finally, pulling a normal epimorphism back along any map $\gamma$ yields a morphism of short exact sequences. If $\gamma$ pulls back to a normal monomorphism $\beta$, then we have a totally normal sequence of monomorphisms as in (\ref{fig:SESs-TotallyNormalSeqMonos}). So, $\gamma$ is a normal monomorphism because every homological category is homologically self-dual.
\end{proof}

\newpage
\begin{corollary}[Pullbacks and preservation/reflection of normal maps\HTag]
  \label{thm:PullbacksPreserve/ReflectNormalMaps}%
  \label{thm:PullbacksPreserve/ReflectPropers}
  About the pullback diagram below the following are true:
  
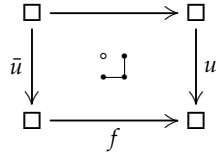
\begin{equation*}
    \xymatrix@R=5ex@C=4em{
    \DiagObj \ar[d]_-{\bar{u}} \ar[r] \PullLU{rd} &
    \DiagObj \ar[d]^-{u} \\
    \DiagObj \ar[r]_-{f} &
    \DiagObj
    }
  \end{equation*}
  \begin{thmlist}
    \item If $u$ is a normal map, then so is its pullback $\bar{u}$ along $f$.
    \item Suppose $f$ is a normal epimorphism. If $\bar{u}$ is a normal map, then so is $u$.
  \end{thmlist}
\end{corollary}
\begin{proof}
  That a pullback along an arbitrary map preserves normal maps follows from the fact that pullbacks preserve normal monomorphisms (\ref{thm:KernelFunctor-Props}) as well as image factorizations (\ref{thm:PullbacksPreserveImageFactorizations}). If $f$ is a normal epimorphism, then the pullback along $f$ also reflects normal monomorphisms (\ref{thm:PullbacksAlongCokerPreserve/ReflectKernels}) and, hence, normal maps.
\end{proof}

We may now extend (\ref{thm:PullbackRecognition-KernelSide-1}) as follows:

\begin{lemma}[Pullback recognition: kernel side II\HTag]
  \label{thm:PullbackRecognition-KernelSide}%
  \cite{DBourn2001}\quad In the commutative diagram below, suppose the sequence at the top is short exact, and $\alpha'=\KerMap{q'}$. Then the following hold:
  
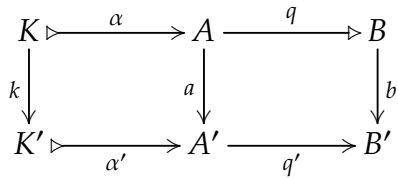
\begin{equation*}
    \xymatrix@R=5ex@C=4em{
    K \ar@{{ |>}->}[r]^-{\alpha} \ar[d]_{k} &
    A \ar@{-{ >>}}[r]^-{q} \ar[d]_{a} &
    B \ar[d]^{b} \\
    K' \ar@{{ |>}->}[r]_-{\alpha'} &
    A' \ar[r]_-{q'} &
    B'
    }
  \end{equation*}
  \begin{enumerate}[(i)]
    \item The left hand square is a pullback if and only if $b$ is a monomorphism. %
          \index{pullback!recognition: kernel side II}
    \item If $b$ is a monomorphism, then $k$ is a monomorphism if and only if $a$ is a monomorphism.
  \end{enumerate}
\end{lemma}
\begin{proof}
  (i) If $b$ is a monomorphism, the left hand square is a pullback by (\ref{thm:PullbackRecognition-KernelSide-1}). Conversely, if the left hand square is a pullback then (\ref{thm:NormalMono-Props}.\ref{thm:Kernel(gf)}) tells us that $ \alpha=\Ker{q'a}$. As $q$ is a cokernel of $\alpha$, we see that $bq=q'a$ is the \NEM image factorization (\ref{thm:OFS-(NormalEpis,Monos)}) of $ q'a$. Consequently, $b$ is a monomorphism.

  (ii) If $b$ is a monomorphism then the square on the left is a pullback by (i). The result now follows, since pullbacks preserve and reflect \eqref{thm:PullbacksPreserve/ReflectMonos} monomorphisms.
\end{proof}

\begin{corollary}[Monomorphism recognition\HTag]
  \label{thm:MonomorphismRecognition}%
  \label{thm:proper-pullback-left}
  Consider a morphism of exact sequences. %
  \index{monomorphism!recognition in map of exact sequences}
  \begin{equation*}
    \xymatrix@R=5ex@C=3em{
    A \ar[r]^-{\alpha} \ar[d]_{a} \ar@{}[rd]|-{\text{(L)}} &
    B \ar@{-{ >>}}[r]^{\beta} \ar[d]_{b} &
    C \ar[r] \ar[d]^{c} &
    0 \\
    X \ar[r]_-{\xi} &
    Y \ar@{-{ >>}}[r]_{\eta} &
    Z \ar[r] & 0
    }
  \end{equation*}
  If the square (L) is a pullback, then $c$ is a monomorphism.
\end{corollary}
\begin{proof}
  The maps $\alpha$ and $\xi$ are normal maps because the horizontal sequences are exact in $B$ and $Y$, respectively. Thus, upon inserting the normal factorizations of $\alpha$ and $\xi$, we obtain this commutative diagram:
  \begin{equation*}
    \xymatrix@R=5ex@C=3em{
    A \ar@/^2ex/[rr]^-{\alpha} \ar[d]_{a} \ar@{-{ >>}}[r] \ar@{}[rd]|-{\text{(U)}} &
    \Img{\alpha}=\Ker{\beta} \ar@{{ |>}->}[r] \ar@<2.5ex>[d] \ar@{}[rd]|-{\text{(V)}}&
    B \ar@{-{ >>}}[r]^{\beta} \ar[d]_{b} &
    C \ar[d]^{c} \\
    X \ar@/_2ex/[rr]_-{\xi} \ar@{-{ >>}}[r] &
    \Img{\xi}=\Ker{\eta} \ar@{{ |>}->}[r] &
    Y \ar@{-{ >>}}[r]_{\eta} &
    Z
    }
  \end{equation*}
  As (L) is a pullback, we infer with Proposition~\ref{thm:PullbacksPreserveImageFactorizations} that both (U) and (V) are pullbacks, and so $c$ is a monomorphism by (\ref{thm:PullbackRecognition-KernelSide}).
\end{proof}

The following proposition extends categorical pullback cancellation; see \ref{thm:Pullbacks,ConcatenatedSquares}.

\begin{proposition}[2-out-of-3 property for pullbacks\HTag]
  \label{thm:SAPullbackCancellationI}
  \label{thm:PullbackCancellation-2-OutOf-3}
  \label{thm:Pullback-2-OutOf-3}%
  \cite[p.~242]{FBorceuxDBourn2004}, \cite[p.~36f]{Borceux-Semiab}, \cite{Janelidze-Sobral-Tholen}\quad In the commutative diagram below assume that $x$ is a normal epimorphism.
  \index{pullback!cancellation $2$-out-of-$3$}%
  \begin{equation*}
    \xymatrix@R=5ex@C=3em{
    A \ar[r] \ar[d] &
    B \ar[r] \ar[d] &
    C \ar[d] \\
    X \ar@{-{ >>}}[r]_-{x} &
    Y \ar[r] &
    Z
    }
  \end{equation*}
  If two out of three of the commutative rectangles in this diagram are pullbacks, then so is the third.
\end{proposition}
\begin{proof}
  By Proposition~\ref{thm:Pullbacks,ConcatenatedSquares} we only need to consider the case where the outer rectangle and the square on the left are pullbacks. To see that the right hand square is a pullback as well, consider the commutative diagram below.
  \begin{equation*}
    \xymatrix@!0@R=7ex@C=4em{
    A \ar[dd]_{a=b^{\ast}} \ar[rr]
    \ar[rd]_{\cong}^{\bar{x}^{\ast}f} &&
    B \ar[dd]_(.7){b} \ar[rr] \ar[rd]^{f} &&
    C \ar@{=}[rd] \ar[dd] \\
    & x^{\ast}y^{\ast}C \ar[ld]_{\hat{c}} \ar[rr]^(.3){\overline{x}}|\hole &&
    y^{\ast}C \ar[ld]_{\bar{c}} \ar[rr]|\hole &&
    C \ar[ld]^{c} \\
    X \ar[rr]_-{x} &&
    Y \ar[rr]_{y} &&
    Z
    }
  \end{equation*}
  The bottom square on the right is constructed as a pullback. We want to show that the comparison map $f$ is an isomorphism. Indeed, via the \PNEInline-Axiom, the normal epimorphism $x$ pulls back along $\bar{c}$ to the normal epimorphism $\bar{x}$. By commutativity of the top left square $f$ is a normal epimorphism, using (\ref{thm:NormalEpis-Props-Normal}). Next, the upper left square is a pullback by pullback cancellation (\ref{thm:Pullbacks,ConcatenatedSquares}), applied to the left hand side of the diagram. As pullbacks reflect monomorphisms (\ref{thm:PullbacksPreserve/ReflectMonos}) we see that $f$ is a monomorphism. So $f$ is an isomorphism by (\ref{thm:IsomorphismRecognition}).
\end{proof}

Our next objective is to establish a converse of (\ref{thm:Pullback->IsoOfKernels}): if a morphism of short exact sequences presents an isomorphism of kernel objects, then the square of cokernels is a pullback. Actually, we already did part of the work in (\ref{thm:PullbackRecognition-SplitSES}) where we showed that
the desired conclusion holds for a morphism of \emph{split short exact sequences}. To make the transition, we employ the following strategy: we prove an intermediate lemma in which kernels are replaced by kernel pairs. The kernel pair construction turns an ordinary short exact sequence into a split short exact sequence to which we can apply (\ref{thm:PullbackRecognition-SplitSES}). So, here is the transitional lemma, a theorem of Barr and Kock.

\begin{theorem}[Barr--Kock\HTag]
  \label{thm:BarrKock}%
  In the diagram below, assume that the square (R) commutes, with $f$ be a normal epimorphism. Then form the left hand side by taking kernel pairs of $f$ and $f'$, respectively.
  \begin{equation*}
    \xymatrix@R=5ex@C=3em{
    \KrnlPr{f} \ar[d]_-{u} \ar@<1ex>[r]^-{p_{1}} \ar@<-1ex>[r]_-{p_{2}} &
    X \ar@{->}[d]_{x} \ar[l]|-{e} \ar@{-{ >>}}[r]^{f} \ar@{}[rd]|-{\text{(R)}} &
    Y \ar[d]^{y} \\
    \KrnlPr{f'}\ar@<1ex>[r]^-{p'_{1}} \ar@<-1ex>[r]_-{p'_{2}} &
    X' \ar[r]_-{f'} \ar[l]|-{e'} &
    Y'
    }
  \end{equation*}
  If at least one of squares given by $ xp_1=p'_{1}u$ or $ xp_2=p'_{2}u$ is a pullback, then the square on the right is a pullback as well.
\end{theorem}
\begin{proof}
  Consider the case where the square given by $xp_1=p'_{1}e$ is a pullback. The remaining case is similar.
  \begin{equation*}
    \xymatrix@R=5ex@C=3em{
    \KrnlPr{f} \PullLU{rd} \ar[r]^{u} \ar[d]_{p_1} &
    \KrnlPr{f'} \PullLU{rd} \ar[r]^{p_{2}'} \ar[d]^-{p'_{1}} &
    X' \ar[d]^{f'} \\
    X \ar[r]_-{x} &
    X' \ar[r]_{f'} &
    Y'
    }
  \end{equation*}
  Both squares in the diagram above are pullbacks. So the outer rectangle is a pullback as well \eqref{thm:Pullbacks,ConcatenatedSquares}. Now this outer rectangle reappears in the commutative diagram below.
  \begin{equation*}
    \xymatrix@R=5ex@C=3em{
    \KrnlPr{f} \PullLU{rd} \ar[r]^{p_{2}} \ar[d]_{p_1} &
    X \ar[r]^{x} \ar@{-{ >>}}[d]^-{f} &
    X' \ar[d]^{f'} \\
    X \ar@{-{ >>}}[r]_-{f} &
    Y \ar[r]_{y} &
    Y'
    }
  \end{equation*}
  Since $f$ is a normal epimorphism, the 2-out-of-3 property for pullbacks (\ref{thm:SAPullbackCancellationI}) shows that the right hand square is a pullback, as desired.
\end{proof}

In the statement of (\ref{thm:BarrKock}) it is important to notice here that the morphism $f'$ is \emph{not} assumed to be a normal epimorphism. This is used in Proposition~(\ref{thm:PullbackFromKerIso}) below.

\begin{proposition}[Pullback recognition by kernel isomorphism\HTag]
  \label{thm:PullbackFromKerIso}%
  In the commutative diagram below, suppose the top sequence is short exact, and $l=\KerMap{g}$. %
  \index{pullback!recognition by kernel isomorphism}%
  \begin{equation*}
    \xymatrix@R=5ex@C=3em{
    0 \ar[r] &
    A \ar@{{ |>}->}[r]^-{k} \ar[d]_-{a} &
    B \ar@{-{ >>}}[r]^-{f} \ar[d]_-{b} &
    C \ar[d]^-{c} \ar[r] &
    0 \\
    0\ar[r] &
    X \ar@{{ |>}->}[r]_-{l} &
    Y \ar[r]_-{g} &
    Z
    }
  \end{equation*}
\end{proposition}
Then the map $a$ is an isomorphism if and only if the square on the right is a pullback.
\begin{proof}
  We use kernel pairs to reduce the right side pullback recognition for a morphism of short exact sequences to right side pullback recognition for a morphism of split short exact sequences: take kernel pairs of $f$ and $g$ to obtain the diagram below.
  \begin{equation*}
    \xymatrix@R=5ex@C=4em{
    \KrnlPr{f} \ar[d]_-{r} \ar@<1ex>[r]^-{p_{1}} \ar@<-1ex>[r]_-{p_{2}} &
    B \ar@{->}[d]_{b} \ar[l]|-{e} \ar@{-{ >>}}[r]^{f} &
    C \ar[d]^{c} \\
    \KrnlPr{g} \ar@<1ex>[r]^-{p'_{1}} \ar@<-1ex>[r]_-{p'_{2}} &
    Y \ar[r]_-{g} \ar[l]|-{e'} &
    Z
    }
  \end{equation*}
  Then $a$ is part of the following morphism of split short exact sequences.
  \begin{equation*}
    \xymatrix@R=5ex@C=3em{
    A \ar@{{ |>}->}[r] \ar[d]_-{a}^{\cong} &
    \KrnlPr{f} \ar@<.5ex>@{-{>}}[r]^-{p_{1}} \ar[d]_-{r} &
    B \ar[d]^-{b} \ar@<.5ex>[l]^-{e} \\
    X \ar@{{ |>}->}[r] &
    \KrnlPr{g} \ar@<.5ex>@{-{>}}[r]^-{p'_{1}} &
    Y \ar@<.5ex>[l]^-{e'}
    }
  \end{equation*}
  Pullback recognition for split short exact sequences \eqref{thm:PullbackRecognition-SplitSES} shows that the square on the right is a pullback. Now, the claim follows with Theorem \ref{thm:BarrKock} of Barr--Kock.
\end{proof}

\begin{proposition}[Normal epimorphism in morphism of short exact sequences\HTag]
  \label{thm:MorphismSESs-Props}
  \label{thm:NormalEpisInSES}%
  A morphism of short exact sequences has the properties stated below. %
  \index{normal!epimorphism!in short exact sequence}%
  \begin{equation*}
    \xymatrix@R=5ex@C=4em{
    K \ar@{{ |>}->}[r]^-{k} \ar[d]_-{\kappa} &
    X \ar@{-{ >>}}[r]^-{q} \ar[d]_-{\xi} &
    Q \ar[d]^-{\rho} \\
    L \ar@{{ |>}->}[r]_-{l} &
    Y \ar@{-{ >>}}[r]_-{r} &
    R
    }
  \end{equation*}
  \begin{thmlist}
    \item \label{thm:NormalEpisInSES-aIso}%
    If $\kappa$ is an isomorphism, then $\xi$ is a normal epimorphism if and only $\rho$ is one such.
    \item \label{thm:NormalEpisInSES-cIso}%
    If $\rho$ is an isomorphism, then $\kappa$ is a normal epimorphism if and only $\xi$ is one such.
  \end{thmlist}
\end{proposition}
\begin{proof}
  (\ref{thm:NormalEpisInSES-aIso})\quad If $\xi$ is a normal epimorphism, then so is $\rho$ by (\ref{thm:NormalEpis-Props-Normal}). For the converse, note first that the square on the right is a pullback by \eqref{thm:PullbackFromKerIso}. If $\rho$ is a normal epimorphism, the pullbacks preserve it via \PNEInline, and so $\xi$ is a normal epimorphism as well.

  (\ref{thm:NormalEpisInSES-cIso})\quad If $\rho$ is an isomorphism, then \eqref{thm:PullbackRecognition-KernelSide-1} the left hand square is a pullback. If $\xi$ is a normal epimorphism, then it pulls back to the normal epimorphism $\kappa$. Now assume that $\kappa$ is a normal epimorphism. The commutative diagram of short exact sequences below results from applying functoriality of image factorizations (\ref{thm:OFS-Factorizations->Functorial}) to $\xi$ and $\rho$.
  \begin{equation*}
    \xymatrix@!0@R=7ex@C=4em{
    K \ar@{{ |>}->}[rr]^-{k} \ar@{-{ >>}}[dd]_-{\kappa} \ar[rd] &&
    X \ar@{-{ >>}}[rr]^-{q} \ar[dd]_(.3){\xi}\ar@{-{ >>}}[rd]^{u} &&
    Q \ar[dd]^(.3){\rho} \ar@{=}[rd] \\
    & \Ker{q} \ar@{-{ >>}}[ld]_{\hat{\kappa}} \ar@{{ |>}->}[rr]|\hole_(.3){\hat{l}} &&
    I \ar@{-{ >>}}[rr]|\hole_(.3){s} \ar@{{ >}->}[ld]^{m} &&
    C \ar@{{ >}->}[ld]^{c}_{\cong} \\
    L \ar@{{ |>}->}[rr]_-{l} &&
    Y \ar@{-{ >>}}[rr]_-{r} &&
    R
    }
  \end{equation*}
  The maps $\hat{\kappa}$ and $s$ are normal epimorphisms by (\ref{thm:NormalEpis-Props-Normal}). But $\hat{\kappa}$ is also a monomorphism because $l\hat{\kappa}=m\hat{l}$ is monic. So, $\kappa$ is an isomorphism by (\ref{thm:IsomorphismRecognition}). By part (i), $m$ is a normal epimorphism. As $m$ is also a monomorphism, it is an isomorphism, again by (\ref{thm:IsomorphismRecognition}). But then $\xi$ is a normal epimorphism, as was to be shown.
\end{proof}

\begin{corollary}[Normal monomorphism in morphism of short exact sequences\HTag]
  \label{thm:NormalMonoInSES-Morphism}%
  In the situation of \eqref{thm:PullbackFromKerIso} suppose $a$ is an isomorphism, and $c$ is a normal monomorphism. Then $b$ is also a normal monomorphism. %
  \index{normal!monomorphism in morphisms of SESs}%
\end{corollary}
\begin{proof}
  From \eqref{thm:PullbackFromKerIso} we know that the right hand square is a pullback. So the claim follows from the fact that pullbacks preserve normal monomorphism (\ref{thm:Pullback->IsoOfKernels}).
\end{proof}

\begin{subordinate}{}
  \begin{subsubordinate}{On the Barr--Kock Theorem}
    The original version of Theorem~\ref{thm:BarrKock}---see \cite[Lemma A.5.8]{FBorceuxDBourn2004}, \cite{Barr}, \cite{DBourn2000}---was considered in the much more general context of a regular category (in the sense of Definition~\ref{def:RegularCategory}). We will, however, not need that stronger version of the result for our purposes.
  \end{subsubordinate}
\end{subordinate}

\begin{exercises}
\begin{exercise}[Equivalent split extensions]
  \label{exe:EquivalentSplitExtensions}
  For an arbitrary abelian group $A$ show that the following three short exact sequences are equivalent:
  \begin{equation*}
    \xymatrix@R=0.5ex@C=4em{
    A \ar@{{ |>}->}[r]^-{\InclsnOf{1}} &
    A \oplus A \ar@{-{ >>}}[r]^-{\PrjctnOnto{2}} &
    A \\
    A \ar@{{ |>}->}[r]^-{\PrdctMapInto{\IdMap_A,-\IdMap_A}} &
    A \oplus A \ar@{-{ >>}}[r]^-{\SumMapOutOf{\IdMap_A,\IdMap_A}} &
    A \\
    A \ar@{{ |>}->}[r]^-{\InclsnOf{2}} &
    A \oplus A \ar@{-{ >>}}[r]^-{\PrjctnOnto{1}} &
    A
    }
  \end{equation*}
\end{exercise}

\begin{exercise}
  Give an example of two short exact sequences in a  homological category whose sum is not exact.
  \index{sum!of short exact sequences}
\end{exercise}

\begin{exercise}
  Find an example in the category of groups which shows that the converse of  (\ref{thm:Pushout->IsoOfCoKers}) does not hold: i.e.\ in the diagram in (\ref{thm:Pushout->IsoOfCoKers}) assume $y$ is an isomorphism. Then the left hand square need not be a pushout. - Use (\ref{thm:PullbackFromKerIso}) to conclude that properties of pushouts in a homological category are not just duals of properties of pullbacks.
\end{exercise}
\end{exercises}
\section[Normal Pushouts]{Normal Pushouts in Homological Categories}
\label{sec:NormalPushouts-Homological}%

For the purpose of constructing di-extensions, normal pushouts have proven to be a useful tool because they are di-extensive. In a homological category, normal pushouts are particularly easy to recognize: Via the recognition criterion of normal pushouts  \ref{thm:SemiNormalPushout-Recognition}, a pushout is di-extensive if and only if it is normal . Proving this proposition and deriving consequences is the objective of this section. - We will repeatedly refer to the diagram below for notation of a morphism of short exact sequences.
\stepcounter{theorem}
\begin{equation}\label{fig:NormalPushoutSES-H}
  \vcenter{
  \xymatrix@R=5ex@C=4em{
  K \ar@{{ |>}->}[r]^-{k} \ar[d]_{\kappa} \ar@{}[rd]|-{\text{(L)}} &
  X \ar@{-{ >>}}[r]^-{q} \ar[d]_{\xi} \ar@{}[rd]|-{\text{(R)}} &
  Q \ar[d]^{\rho} \\
  L \ar@{{ |>}->}[r]_-{l} &
  Y \ar@{-{ >>}}[r]_-{r} &
  R
  }
  }
\end{equation}

\begin{proposition}[Seminormal pushout recognition\HTag]
  \label{thm:SemiNormalPushout-Recognition}%
  In a homological category, given a morphism of short exact sequences, as in \eqref{fig:NormalPushoutSES-H}, the square (R) is a seminormal pushout if and only if $\kappa$ is a normal epimorphism. %
  \index{seminormal!pushout recognition}%
\end{proposition}
\begin{proof}
  If (R) is a seminormal pushout, then $k$ is a normal epimorphism by (\ref{thm:HDS<->SemiNormalPushoutProp}). Conversely, suppose  $k$ is a normal epimorphism. Consider this commutative diagram whose bottom right square is the pullback of $r$ along $\rho$.
  \begin{equation}\label{fig:NormalPushoutSES-H-Pull}
    \vcenter{
    \xymatrix@R=5ex@C=3em{
    K \ar@{{ |>}->}[rr]^-{k} \ar[dd]_{\kappa} \ar[dr]_{\kappa} &&
    X \ar@{-{ >>}}[rr]^-{q} \ar[dd]_(0.65){\xi} \ar[rd]^{\bar{\xi}} &&
    Q \ar[dd]_(.3){\rho} \ar@{=}[rd] \\
    & L \ar@{{ |>}->}[rr]|\hole^(.3){\hat{k}} &&
    P \ar[ld]^(0.4){\bar{\rho}} \ar@{-{ >>}}[rr]|\hole_(.25){\bar{r}} &&
    Q \ar[dl]^(0.4){\rho} \\
    L \ar@{{ |>}->}[rr]_-{l} \ar@{=}[ru] &&
    Y \ar@{-{ >>}}[rr]_-{r} &&
    R
    }
    }
  \end{equation}
  So, the normal epimorphism $r$ pulls back to a normal epimorphism $\bar{r}$. This makes the upward facing part of the diagram a morphism of short exact sequences. So, $\bar{\xi}$ is a normal epimorphism by (\ref{thm:NormalEpisInSES}), which means that (R) is a seminormal pushout.
\end{proof}

\begin{corollary}[Pullback stability of (semi)normal pushouts\HTag]
  \label{thm:NormalPushout-PullbackStabilty}
  In a homological category, a commutative cube of maps has the following properties: %
  \index{seminormal!pushout - pullback stability}%
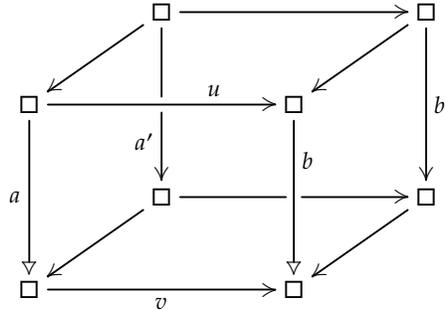
  \begin{equation*}
    \xymatrix@R=4ex@C=3em{
    & \DiagObj \ar[rr] \ar[dd]|\hole_(0.7){a'} \ar[ld] &&
    \DiagObj \ar[ld] \ar[dd]^{b'} \\
    \DiagObj \ar@{-{ >>}}[dd]_{a} \ar[rr]^(.7){u} &&
    \DiagObj \ar@{-{ >>}}[dd]^(.3){b} \\
    & \DiagObj \ar[ld] \ar[rr]|\hole &&
    \DiagObj \ar[ld] \\
    \DiagObj \ar[rr]_-{v} &&
    \DiagObj
    }
  \end{equation*}
  \begin{thmlist}
    \item \label{thm:NormalPushout-PullbackStabilty-SemiNormal}%
    Suppose the front face is a seminormal pushout with $a$ and $b$ normal epimorphisms. If the two squares on the side are pullbacks, then the back face is a seminormal pushout in which $a'$ and $b'$ are normal epimorphisms.
    \item \label{thm:NormalPushout-PullbackStabilty-Normal}%
    Suppose the front face is a normal pushout. If the side faces as well as the top and bottom faces are pullbacks, then the back face is a normal pushout.
  \end{thmlist}
\end{corollary}
\begin{proof}
  (\ref{thm:NormalPushout-PullbackStabilty-SemiNormal})\quad If the front face is a seminormal pushout, then the induced map $\Ker{a}\to \Ker{b}$ is a normal epimorphism. If the side faces are pullbacks, then both $a'$ and $b'$ are normal epimorphisms. Moreover, we have isomorphisms $\Ker{a}\cong \Ker{a'}$ and $\Ker{b}\cong \Ker{b'}$ by (\ref{thm:Pullback->IsoOfKernels}). But then the induced map $\Ker{a'}\to \Ker{b'}$ is a normal epimorphism. So the back face is a seminormal pushout by (\ref{thm:SemiNormalPushout-Recognition}). The proof (\ref{thm:NormalPushout-PullbackStabilty-Normal}) is similar.
\end{proof}

\begin{corollary}[Di-extensive pushouts are normal\HTag]
  \label{thm:DiExtensivePush=NormalPush}%
  A pushout square of normal epimorphisms is a di-extensive if and only if it is a normal pushout. %
  \index{di-extensive!pushout in \HTag}%
\end{corollary}
\begin{proof}
  We know from (\ref{thm:NormalPushout->DiExtensive}) that every normal pushout in a homologically self-dual category, hence in a homological category, is di-extensive. With Proposition \ref{thm:SemiNormalPushout-Recognition}, we see that a di-extensive pushout is a normal pushout.
\end{proof}

\begin{corollary}[Homological $\implies$ \DPNInline\HTag]
  \label{thm:Homological->DPN}
  Every homological category has the \DPNInline-property. %
  \index{homological!category satisfies \DPNInline}%
  \index[acr]{d!\DPNInline\IndSep property that dinversion preserves normal maps}%
\end{corollary}
\begin{proof}
  Given a normal map with antinormal decomposition, as in (L) of \eqref{fig:NormalPushoutSES-H}, then $\kappa$ is a normal epimorphism. So (R) is a normal pushout, hence is di-extensive. With (\ref{thm:Dinversion-DiExtensivePush/Pull}), we see that dinversion preserves the normal map $\xi k$.
\end{proof}

\begin{corollary}[Sectioned epi of normal epis is normal pushout\HTag]
  \label{thm:SplitEpiOfNormalEpis->NormalPush}%
  \label{thm:NormalEpiOfSplitEpis->NormalPush}%
  If square (R) below suppose $(q,r)\from \SctndEpi{\xi}{\alpha}\to \SctndEpi{\rho}{\beta}$ is a morphism of sectioned epimorphisms. Then (R) is a normal pushout.
  \begin{equation*}
    \xymatrix@R=6ex@C=5em{
    K \ar@{{ |>}->}[r]^-{k} \ar[d]_{\kappa} &
    X \ar@{-{ >>}}[r]^-{q} \ar@{-{ >>}}@<-0.5ex>[d]_{\xi} \ar@{}[rd]|-{\text{(R)}}&
    Q \ar@{-{ >>}}@<-0.5ex>[d]_{\rho} \\
    L \ar@{{ |>}->}[r]_-{l} &
    Y \ar@{-{ >>}}[r]_-{r} \ar@<-0.5ex>[u]_{\alpha} &
    R \ar@<-0.5ex>[u]_{\beta}
    }
  \end{equation*}
\end{corollary}
\begin{proof}
  We claim that $\kappa$ is an absolute epimorphism. Indeed, a map $\sigma\from L\to K$ with $k\sigma = \alpha l$ comes from the fact that $q\alpha l=\beta rl=\ZeroMap$. That $\kappa\sigma=\IdMapOn{L}$ comes from the monic property of $l$ via the computation
  \begin{equation*}
    l \sigma \kappa = \xi k \sigma = \xi \alpha l = l
  \end{equation*}
  Thus $\sigma$ is a section of $\kappa$. So, $\kappa$ is a normal epimorphism by (\ref{thm:SplitEpi->Normal-Homological}). Hence, (R) is a normal pushout by (\ref{thm:SemiNormalPushout-Recognition}).
\end{proof}

\begin{corollary}[Normal components in morphism of short exact sequences\HTag]
  \label{thm:MorphismOfSESs-Properness}
  \label{thm:MorphismOfSESs-NormalComponents}%
  A morphism of short exact sequences, such as \eqref{fig:NormalPushoutSES-H} has the following properties. %
  \index{normal!map in morphism of short exact sequences}%
  \begin{enumerate}[(i)]
    \item If $k$ is an isomorphism, then $u$ is normal, respectively a normal monomorphism, if and only if $v$ is normal, respectively a normal monomorphism.
    \item If $k$ is a normal epimorphism, then $u$ is normal if and only if $v$ is normal.
    \item If $v$ is a monomorphism, and $u$ is normal, then $k$ is normal.
  \end{enumerate}
\end{corollary}
\begin{proof}
  (i)\quad If $k$ is an isomorphism, then (R) is a pullback (\ref{thm:PullbackFromKerIso}). Since $r$ is a normal epimorphism, pulling back along it preserves and reflects normal monomorphisms (\ref{thm:PullbacksAlongCokerPreserve/ReflectKernels}) and normal maps by (\ref{thm:PullbacksPreserve/ReflectNormalMaps}).

  (ii)\quad If $k$ is a normal epimorphism, then (R) is a seminormal pushout by (\ref{thm:SemiNormalPushout-Recognition}). Thus, if we factor the given morphism of short exact sequences via the pullback of $r$ along $v$, then the map $\bar{u}$ in the diagram \eqref{fig:NormalPushoutSES-H-Pull} is a normal epimorphism. Thus $\Img{u}=\Img{\bar{v}}$. So $u$ is normal if and only if $\bar{v}$ is normal. By (i), $\bar{v}$ is normal if and only if $v$ is normal.

  (iii)\quad If $v$ is a monomorphism, then (L) is a pullback (\ref{thm:PullbackRecognition-KernelSide-1}). It preserves the normal map $u$ by (\ref{thm:MorphismOfSESs-NormalComponents}).  So $k$ is normal.
\end{proof}

Another way of formulating (\ref{thm:MorphismOfSESs-NormalComponents}.ii) is: In a homological category, seminormal pushouts preserve and reflect normal maps.

\begin{corollary}[Normal pushout of coproduct/product comparison maps\HTag]%
  \label{thm:CoKernels->NormalPushOfComparisonMaps}%
  \label{thm:CoKernels->RegularPushOfComparisonMaps}
  A pair of normal epimorphisms $u\from X\to X'$ and $v\from Y\to Y'$ induces a normal pushout diagram of coproduct/product comparison maps: %
  \begin{equation*}
    \xymatrix@R=5ex@C=4em{
    \CoPrdct{X}{Y} \ar[d]_{\CoPrdct{u}{v}} \ar@{-{ >>}}[r]^-{\SumProdComp{X}{Y}} &
    \Prdct{X}{Y} \ar[d]^{\Prdct{u}{v}} \\
    \CoPrdct{X'}{Y'} \ar@{-{ >>}}[r]_-{\SumProdComp{X'}{Y'}} &
    \Prdct{X'}{Y'}
    }
  \end{equation*}
\end{corollary}
\begin{proof}
  Via Proposition~\ref{thm:Sum->ProductIsCokernel}, we know that the horizontal maps are normal epimorphisms. The map $ \CoPrdct{u}{v}\from \CoPrdct{X}{Y}\to \CoPrdct{X'}{Y'}$ is a normal epimorphism because colimits commute with colimits. The map $\Prdct{u}{v}$ is a normal epimorphism (\ref{thm:NormalEpis-Props-Normal}). To see that the square is a normal pushout, it suffices (\ref{thm:SemiNormalPushout-Recognition}) to show that the map
  \begin{equation*}
    \SumProdComp{X}{Y}|\from \Ker{\CoPrdct{u}{v}} \longrightarrow \Ker{ \Prdct{u}{v} } \cong \Prdct{\Ker{u}}{\Ker{v}}
  \end{equation*}
  is a normal epimorphism. This is indeed the case: $\CoPrdct{\Ker{u}}{\Ker{v}}$ maps into $\Ker{\CoPrdct{u}{v}}$ and its coproduct/product comparison map is a normal epimorphism to $\Prdct{\Ker{u}}{\Ker{v}}$. By (\ref{thm:NormalEpis-Props-Normal}) the map $ {\Ker{\CoPrdct{u}{v}}\to \Prdct{\Ker{u}}{\Ker{v}}}$ is a normal epimorphism as well, and the claim follows.
\end{proof}

Any morphism of short exact sequences allows us to construct an associated seminormal pushout as follows:

\begin{lemma}[(Semi)normal pushout from morphism of short exact sequences\HTag]
  \label{thm:SemiNormalPushoutFrom-SES}
  \label{thm:SemiRegularPushoutFrom-SES}
  Given a morphism of short exact sequences as on the left, the associated square on the right is a seminormal pushout.
  \begin{equation*}
    \xymatrix@R=5ex@C=3em{
    K \ar@{{ |>}->}[r]^-{k} \ar[d]_{\kappa} &
    X \ar@{-{ >>}}[r]^-{q} \ar[d]_{\xi} &
    Q \ar[d]^{\rho} &&
    L+X \ar@{-{ >>}}[r]^-{\SumMapOutOf{0,q}} \ar[d]_{\SumMapOutOf{l,\xi}}&
    Q \ar[d]^{\rho} \\
    L \ar@{{ |>}->}[r]_-{l} &
    Y \ar@{-{ >>}}[r]_-{r} &
    R &&
    X \ar@{-{ >>}}[r]_-{r}&
    R
    }
  \end{equation*}
  The square on the right is a normal pushout if and only if $\rho$ is a normal epimorphism.
\end{lemma}
\begin{proof}
  Consider the commutative diagram below; its floor is the pullback of $r$ along $\rho$, and the horizontal sequences are short exact.
  
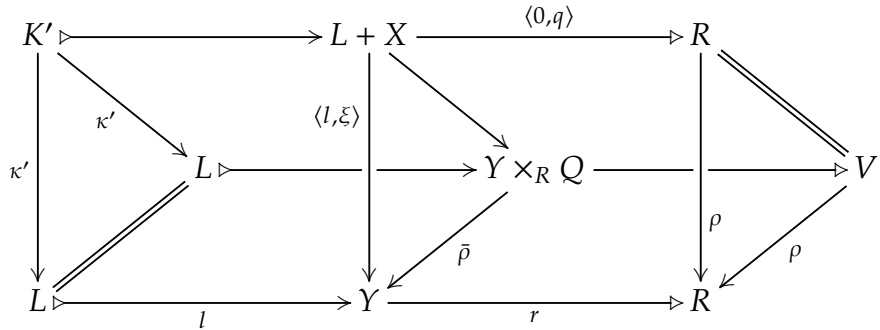
\begin{equation*}
    \xymatrix@!@C=1.5em@R=1ex{
    K' \ar@{{ |>}->}[rr] \ar[dd]_{\kappa'} \ar[rd]_{\kappa'} &&
    L+X \ar@{-{ >>}}[rr]^-{\SumMapOutOf{0,q}} \ar[dd]_(0.3){\SumMapOutOf{l,\xi}} \ar[rd] &&
    R \ar@{=}[rd] \ar[dd]^(0.7){\rho} \\
    & L \ar@{{ |>}->}[rr]|\hole \ar@{=}[ld] &&
    Y\prdct_R Q \ar@{-{ >>}}[rr]|\hole \ar[ld]^{\bar{\rho}} &&
    V \ar[ld]^{\rho} \\
    L \ar@{{ |>}->}[rr]_-{l} &&
    Y \ar@{-{ >>}}[rr]_-{r} &&
    R
    }
  \end{equation*}
  The map $\kappa'$ is a normal epimorphism because it is sectioned by the coproduct inclusion $\InclsnOf{L}\from L\to L+X$ factored through to $K'$. This renders the front face on the right a seminormal pushout by (\ref{thm:SemiNormalPushout-Recognition}). So, the comparison map $L+X\to Y\prdct_R Q$ is a normal epimorphism. If now $\rho$ is a normal epimorphism, then so are $\bar{\rho}$ and, hence, $\SumMapOutOf{l,\xi}$. So, the square on the right is a normal pushout, as was to be shown.
\end{proof}

The following corollary complements (\ref{thm:JEE-then-cokernel-NE}).

\begin{corollary}[Normal epimorphism recognition - cokernel side\HTag]%
  \label{thm:CoKerRecognition-CoKernelSide}
  \label{thm:NormalEpiRecognition-CoKernelSide}%
  In the morphism of short exact sequences in \eqref{fig:NormalPushoutSES-H}, the maps  $l$ and $\xi$ are jointly extremally epimorphic if and only if $\rho$ is a normal epimorphism. %
  \index{normal!epimorphism - recognition}%
\end{corollary}
\begin{proof}
  If the maps $l$ and $\xi$ are jointly extremally epimorphic, then $\rho$ is a normal epimorphism by (\ref{thm:JEE-then-cokernel-NE}). Conversely, if $\rho$ is a normal epimorphism, then so is $\SumMapOutOf{l,\xi}$ by (\ref{thm:SemiRegularPushoutFrom-SES}). So, $l$ and $\xi$ are jointly extremally epimorphic.
\end{proof}

\begin{proposition}[Normal pushout recognition by kernel pair\HTag]
  \label{thm:NormalPushoutRecognizeKernelPair}
  \label{thm:RegularPushoutRecognizeKernelPair}
  Assume that the diagram below consists of a commutative square (R) of normal epimorphisms, together with the kernel pairs of $f$ and $f'$.
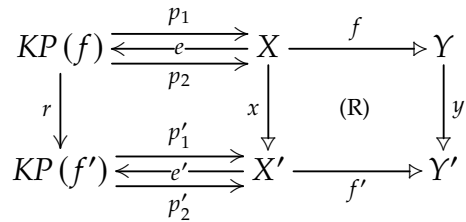
  \begin{equation*}
    \xymatrix@R=5ex@C=4em{
    \KrnlPr{f} \ar[d]_-{r} \ar@<1ex>[r]^-{p_{1}} \ar@<-1ex>[r]_-{p_{2}} &
    X \ar@{-{ >>}}[d]_{x} \ar[l]|-{e} \ar@{-{ >>}}[r]^{f} \ar@{}[rd]|-{\text{(R)}} &
    Y \ar@{-{ >>}}[d]^{y} \\
    \KrnlPr{f'} \ar@<1ex>[r]^-{p'_{1}} \ar@<-1ex>[r]_-{p'_{2}} & X' \ar@{-{ >>}}[r]_-{f'} \ar[l]|-{e'} &
    Y'
    }
  \end{equation*}
  Then the square (R) is a normal pushout if and only if $r$ is a normal epimorphism.
\end{proposition}
\begin{proof}
  As in the proof of \eqref{thm:PullbackFromKerIso}, the kernels of $p_1$ and $p'_{1}$ yield the following morphism of split short exact sequences:
  \begin{equation*}
    \vcenter{\xymatrix@R=5ex@C=3em{
    0 \ar[r] &
    \Ker{f} \ar@{{ |>}->}[r] \ar[d]_-{k} &
    \KrnlPr{f} \ar@<.5ex>@{-{ >>}}[r]^-{p_{1}} \ar[d]_-{r} \ar@{}[rd]|-{\text{(S)}} &
    X \ar@{-{ >>}}[d]^-{x} \ar@{{ >}->}@<.5ex>[l]^-{e} \ar[r] &
    0 \\
    0 \ar[r] & \Ker{f'} \ar@{{ |>}->}[r] &
    \KrnlPr{f'} \ar@<.5ex>@{-{ >>}}[r]^-{p'_{1}} &
    X' \ar@{{ >}->}@<.5ex>[l]^-{e'} \ar[r] & 0
    }}
  \end{equation*}
  By (\ref{thm:MorphismOfSESs-NormalComponents}.ii), square (S) is a normal pushout if and only if $k$ is a normal epimorphism. Now if $k$ is a  normal epimorphism, then (\ref{thm:MorphismOfSESs-NormalComponents}.ii) tells us that $r$ is a normal epimorphism as well. On the other hand, if $r$ is a normal epimorphism, then (\ref{thm:NormalEpiOfSplitEpis->NormalPush}) says that (S) is a normal pushout, and so $k$ is a normal epimorphism.
\end{proof}

\begin{proposition}[Pullback / seminormal pushout recognition\HTag]
  \label{thm:Pullback/SemiNormalPushoutRecognition}%
  If in the morphism of short exact sequence \eqref{fig:NormalPushoutSES-H} the map $v$ is an isomorphism, and $k$ is a normal epimorphism, then (L) is a pullback and a seminormal pushout. %
  \index{seminormal pushout!recognition II}\index{pushout!recognition II}%
\end{proposition}
\begin{proof}
  From (\ref{thm:NormalEpisInSES}) we see that $u$ is a normal epimorphism. Next, square (L) is a pullback follows from (\ref{thm:PullbackRecognition-KernelSide-1}). But then $\Ker{\alpha}\to \Ker{\beta}$ is an isomorphism, and so (L) is a seminormal pushout by (\ref{thm:SemiNormalPushout-Recognition}).
\end{proof}

\begin{corollary}[Normal pushout recognition\HTag]
  \label{thm:NormalPushOut-Recognize-H}%
  In the morphism of short exact sequence presented in \eqref{fig:NormalPushoutSES-H}, square (R) is a normal pushout if and only if $k$ and $v$ are normal epimorphisms. \NoProof%
  \index{normal!pushout recognition}%
\end{corollary}

\bigskip

\begin{subordinate}{}
  \begin{subsubordinate}{Normal epimorphisms in $\SEpisIn{X}$}
    \index{category!of split epimorphisms}%
    In a homological category $\Ctgry{X}$, any diagram as in (\ref{thm:NormalEpiOfSplitEpis->NormalPush}) is actually a normal epimorphism in the functor category $\SEpisIn{X}$ of sectioned epimorphisms in $\Ctgry{X}$. This is so because (co)limits in that functor category are computed object-wise. Thus we may paraphrase (\ref{thm:NormalEpiOfSplitEpis->NormalPush}) by saying that a normal epimorphism in the category of sectioned epimorphisms over $\Ctgry{X}$ has an underlying normal pushout in the $\Ctgry{X}$.

    In particular, any double sectioned epimorphism, when considered as a square in $\Ctgry{X}$, is a normal pushout---compare with (\ref{thm:AbsolutePush/Pull-SectionedMorInSEpi(X)}) and~(\ref{thm:ProtoMaltsev}).
  \end{subsubordinate}

  \begin{subsubordinate}{Alternate view of the axiom \KSGInline}
    Corollary (\ref{thm:NormalEpiRecognition-CoKernelSide}) specializes to axiom \KSGInline\ by choosing $q\DefEq \rho \DefEq \IdMapOn{R}$. In this case, Lemma~\ref{thm:SemiRegularPushoutFrom-SES} is especially interesting because we find this morphism of split short exact sequences
    \begin{equation*}
      \xymatrix@!@C=3em@R=1ex{
      K' \ar@{{ |>}->}[r] \ar[d]_{\varepsilon} &
      L+R \ar@<-.5ex>@{-{ >>}}[r]_-{\SumMapOutOf{0,\IdMapOn{R}}} \ar@{-{ >>}}[d]_{\SumMapOutOf{l,\xi}} &
      R \ar@{=}[d] \ar@<-.5ex>[l]_-{\InclsnOf{R}}\\
      L \ar@{{ |>}->}[r]_-{l} &
      Y \ar@<-.5ex>@{-{ >>}}[r]_-{q} &
      R \ar@<-.5ex>[l]_-{\xi}
      }
    \end{equation*}
    Following \cite{Bourn-Janelidze:Semidirect},  we will interpret $Y$ as the `semidirect product'  $Y=L\rtimes_{\varepsilon} R$ of $L$ and $R$ with respect to the `internal action' given by $\varepsilon\from K'\to L$.
  \end{subsubordinate}
\end{subordinate}

\begin{exercises}

\begin{exercise}[Short exact sequence / regular pushout-pullback\HTag]
  \label{exe:SES-RegPush/Pull}
  In a homological category, show that a diagram $X \XRA{x} Y \XRA{y} Z$ is a short exact sequence if and only if the square below is a seminormal pushout and a pullback.
  \begin{equation*}
    \xymatrix@R=5ex@C=4em{
    X \ar[r] \ar[d]_{x} &
    0 \ar[d] \\
    Y \ar[r]_-{y} &
    Z
    }
  \end{equation*}
\end{exercise}

\begin{exercise}[Morphism of short exact sequences - special case\HTag]
  \label{exe:SES-Morphism-I}
  In the morphism of short exact sequences below, assume that $a$ is a normal epimorphism.
  \begin{equation*}
    \xymatrix@R=4ex@C=2em{
    K \ar@{{ |>}->}[rr]^-{k} \ar@{-{ >>}}[d]_{\kappa}&&
    X \ar@{-{ >>}}[rr]^-{q} \ar[d]_{\xi} &&
    Q \ar[d]^{\rho} \\
    L\ar@{{ |>}->}[rr]_-{l} &&
    Y \ar@{-{ >>}}[rr]_-{r} &&
    R
    }
  \end{equation*}
  Show that the induced map $\kappa\from \Ker{\xi}\to \Ker{\rho}$ is a normal epimorphism. - Hint: Apply Exercise \ref{exe:FactorMorphisZeroMaps} to the morphism $(\kappa,\beta,\eta)$ of composite zero maps.
\end{exercise}

\begin{exercise}[Detection of normal pullbacks fails\HTag]
  \label{thm:NormalPullback-Recognize}
  In a homological category consider a commutative square of normal monomorphisms. If the map of cokernels of two opposing arrows is a normal monomorphism show that the square need not be be normal pullback.
\end{exercise}

\begin{exercise}[Linear category from normal pullbacks\HTag]
  \label{exe:PullbacksNormal->Linear}
  In a homological category $\EuScript{X}$, assume that for all objects $X$ and $Y$ the pullback of the canonical inclusions $\PrdctMapInto{\IdMapOn{X},0}\from  X\to \Prdct{X}{Y}$ and $\PrdctMapInto{0,\IdMapOn{Y}}\from Y\to \Prdct{X}{Y}$ is a normal pullback. Show that $\EuScript{X}$ is a linear category. %
  \index{linear!category}%
\end{exercise}
\end{exercises}
\newpage
\section[The (Short) 5-Lemma]{The (Short) 5-Lemma}
\label{sec:(Short)5-Lemma}

In Section \ref{sec:SESMaps} we proved a version of the Short $5$-Lemma, namely (\ref{thm:Short-5-Primordial}), which we termed  `primordial' because we needed to assume the map $u$ in a morphism $(\kappa,\xi,\rho)$ of short exact sequences is a normal map. In the structurally richer environment of a homological category, we can eliminate this assumption and prove the Short $5$-Lemma (\ref{thm:Short5}) in the form familiar from abelian categories. It then provides a stepping stone for proving (variants of) the $5$-Lemma; see (\ref{thm:5-Lemma}) and (\ref{thm:5-Lemma,Easy}).

The Short $5$-Lemma provides a useful tool for identifying when a given morphism is an isomorphism. We close this section by explaining how it plays a key role in the definition of the functor $\Ext^{1}$.

\begin{theorem}[Short 5-Lemma\HTag]
  \label{thm:Short5}%
  \cite[p.~275f]{FBorceuxDBourn2004}\quad Consider a morphism of short exact sequences. %
  \index{Short 5-Lemma}%
  \begin{equation*}
    \xymatrix@R=5ex@C=4em{
    K \ar@{{ |>}->}[r]^{k} \ar[d]_{\kappa} &
    X \ar@{-{ >>}}[r]^{q} \ar[d]_{\xi} \ar@{}[rd]|-{(R)}&
    Q \ar[d]^{\rho} \\
    L \ar@{{ |>}->}[r]_{\lambda} &
    Y \ar@{-{ >>}}[r]_{r} &
    R
    }
  \end{equation*}
  If $\kappa$ and $\rho$ are isomorphisms, then $\xi$ is an isomorphism as well.
\end{theorem}
\begin{proof}
  As $\kappa$ is an isomorphism, the square (R) is a pullback by (\ref{thm:PullbackFromKerIso}). Thus $\xi$ is the pullback of the isomorphism $\rho$. So $\xi$ is an isomorphism as well.
\end{proof}

Since the proof of the Short $5$-Lemma is very concise, let us emphasize that, ultimately, it relies on the \KSGInline-property of homological categories: For a split short exact sequence it says that the middle object is generated by the end objects. If then those generating objects are communicated isomorphically by a morphism of split short exact sequences, then so is the middle object. This is the content of the Split Short 5-Lemma (\ref{thm:SplitShort5}). In the previous section, we explained how the kernel pair construction may be used to lift split short exact sequence information to short exact sequence information. This was critical in proving the pullback recognition criterion (\ref{thm:PullbackFromKerIso}) which led to the above short proof of (\ref{thm:Short5}).

Some classical diagram lemmas follow immediately from the Short $5$-Lemma:

\begin{corollary}[5-Lemma\HTag]
  \label{thm:5-Lemma}%
  Consider a morphism of exact sequences. %
  \index{$5$-Lemma}
  \begin{equation*}
    \xymatrix@R=5ex@C=3em{
    \DiagObj \ar@{->>}[d]_{a} \ar[r] &
    \DiagObj \ar[d]_{b}^{\cong} \ar[r] &
    \DiagObj \ar[d]_{c} \ar[r] &
    \DiagObj \ar[d]^{\cong}_{d} \ar[r] &
    \DiagObj \ar@{{ >}->}[d]_{e} \\
    \DiagObj \ar[r] &
    \DiagObj \ar[r] &
    \DiagObj \ar[r] &
    \DiagObj \ar[r] &
    \DiagObj
    }
  \end{equation*}
  If $a$ is an epimorphism, $e$ a monomorphism, $b$, $d$ isomorphisms, then $c$ is an isomorphism.
\end{corollary}
\begin{proof}
  Via pushout and pullback recognition, we reduce the proof of the $5$-Lemma to an application of the Short $5$-Lemma (\ref{thm:Short5}). As the rows are exact we obtain image factorizations at selected points as depicted in the commutative diagram below.
  \begin{equation*}
    \xymatrix@R=3ex@C=2.5em{
    &&& \DiagObj \ar@{{ |>}->}[dr] \ar[dd]|\hole_(.7){u} &&
    \DiagObj \ar@{{ |>}->}[dr] \ar[dd]|\hole_(.7){v} \\
    \DiagObj \ar@{ ->>}[dd]_{a} \ar[rr] &&
    \DiagObj \ar[dd]_{b}^{\cong} \ar[rr] \ar@{-{ >>}}[ru] &&
    \DiagObj \ar[dd]_{c} \ar[rr] \ar@{-{ >>}}[ru] &&
    \DiagObj \ar[dd]^{\cong}_{d} \ar[rr] &&
    \DiagObj \ar@{{ >}->}[dd]_{e} \\
    &&& \DiagObj \ar@{{ |>}->}[dr] &&
    \DiagObj \ar@{{ |>}->}[dr] \\
    \DiagObj \ar[rr] &&
    \DiagObj \ar[rr] \ar@{-{ >>}}[ru]&&
    \DiagObj \ar[rr] \ar@{-{ >>}}[ru] &&
    \DiagObj \ar[rr] &&
    \DiagObj
    }
  \end{equation*}
  The square with vertical maps $b$ and $u$ is a pushout because $a$ is an epimorphism \eqref{thm:PushoutRecognize-Categorical}, implying that $u$ is an isomorphism. The square with vertical maps $v$ and $d$ is a pullback because $e$ is a monomorphism \eqref{thm:PullbackRecognition-KernelSide-1}, implying that $v$ is an isomorphism. Now the Short 5-Lemma \eqref{thm:Short5} tells us that $c$ is an isomorphism.
\end{proof}

Frequently, the following special case of the 5-Lemma is sufficient:

\begin{corollary}[`Easy' 5-Lemma\HTag]
  \label{thm:5-Lemma,Easy}%
  In the morphism of exact sequences below assume that $a$, $b$, $d$, $e$ are isomorphisms. %
  \index{$5$-Lemma!easy}%
  \begin{equation*}
    \xymatrix@R=5ex@C=3em{
    \DiagObj \ar[d]_{a}^{\cong} \ar[r] &
    \DiagObj \ar[d]_{b}^{\cong} \ar[r] &
    \DiagObj \ar[d]_{c} \ar[r] &
    \DiagObj \ar[d]^{\cong}_{d} \ar[r] &
    \DiagObj \ar[d]_{e}^{\cong} \\
    \DiagObj \ar[r] &
    \DiagObj \ar[r] &
    \DiagObj \ar[r] &
    \DiagObj \ar[r] &
    \DiagObj
    }
  \end{equation*}
  Then $c$ is an isomorphism as well. \NoProof
\end{corollary}

In the remainder of this section, we develop consequences of the Short $5$-Lemma. We start by improving substantially upon our basic product recognition tool (\ref{thm:ProductRecognition}): In a homological category, if a normal monomorphism admits a section, then this section is a product projection.

\begin{proposition}[Product recognition via kernel splitting\HTag]
  \label{thm:ProductRecognition-SplitKer->Product}%
  Consider a short exact sequence such as %
  \index{product!recognition via kernel splitting}%
  \begin{equation*}
    \xymatrix@R=5ex@C=3em{
    0 \ar[r] &
    A \ar@<.5ex>@{{ |>}->}[r]^{k}  &
    B \ar@{-{ >>}}[r]^{f}  \ar@<.5ex>@{-{ >>}}[l]^{q} &
    C \ar[r] &
    0 }
  \end{equation*}
  If $qk=1_A$, then $\PrdctMapInto{q,f}\from B \to A\times C$ is an isomorphism, so that $(B,q,f)$ forms a product of $A$ and $C$.
\end{proposition}
\begin{proof}
  Together, the maps $q$ and $f$ form this morphism of short exact sequences:
  \begin{equation*}
    \xymatrix@R=5ex@C=3em{
    0 \ar[r] &
    A \ar@<.5ex>@{{ |>}->}[r]^{k} \ar@{=}[d] &
    B \ar@{-{ >>}}[r]^{f} \ar[d]_{\PrdctMapInto{q,f}} \ar@<.5ex>@{-{ >>}}[l]^{q} &
    C \ar[r] \ar@{=}[d] &
    0 \\
    0 \ar[r] &
    A \ar@{{ |>}->}[r]_-{\PrdctMapInto{1_A,0}} &
    A\times C \ar@{-{ >>}}[r]_-{\pi_C} &
    C \ar[r] & 0
    }
  \end{equation*}
  The Short 5-Lemma \ref{thm:Short5} shows that $\PrdctMapInto{q,f}$ is an isomorphism.
\end{proof}

\begin{proposition}[Pullbacks reflect isomorphisms\HTag]
  \label{thm:BaseChangeReflectsIsos}
  \label{thm:PullBackReflectsIsos}%
  Consider the morphism of normal epimorphisms below. %
  \index{base change!reflects isomorphisms}
  \begin{equation*}
    \xymatrix@R=5ex@C=3em{
    X \ar[r]^-{u} \ar@{-{ >>}}[d]_{p} &
    Y \ar@{-{ >>}}[d]^{q} \\
    R \ar@{=}[r] &
    R
    }
  \end{equation*}
  If pulling this diagram back along $f\from S\to R$ induces an isomorphism $\bar{u} \from \bar{X}\to \bar{Y}$, then $u$ is an isomorphism.
\end{proposition}
\begin{proof}
  The situation at hand leads to the commutative diagram below.
  \begin{equation*}
    \xymatrix@!0@R=6ex@C=3.5em{
    & \bar{L} \ar[rr]^-{\cong} \ar@{{ |>}->}[dd]|\hole &&
    L \ar@{{ |>}->}[dd] \\
    \bar{K} \ar[ru]^{\bar{\kappa}}_{\cong} \ar[rr]^(.7){\cong} \ar@{{ |>}->}[dd] &&
    K \ar[ru]^{\kappa} \ar@{{ |>}->}[dd] \\
    & \bar{Y} \ar[rr]|\hole \ar@{-{ >>}}[dd]|\hole_(.7){\bar{q}} &&
    Y \ar@{-{ >>}}[dd]_{q} \\
    \bar{X} \ar[ru]^{\bar{u}}_{\cong} \ar[rr] \ar@{-{ >>}}[dd]_{\bar{p}} &&
    X \ar[ru]_{u} \ar@{-{ >>}}[dd]_(.7){p} \\
    & S \ar[rr]|\hole^(.3){f} &&
    R \\
    S \ar[rr]_-{f} \ar@{=}[ru] &&
    R \ar@{=}[ru]
    }
  \end{equation*}
  The front and back faces of the bottom cube are pullbacks, and the top cube is formed from the kernels of $p$, $q$, $\bar{p}$, $\bar{q}$. Therefore (\ref{thm:Pullback->IsoOfKernels}) the maps $\bar{K}\to K$ and $\bar{L}\to L$ are isomorphisms. By the 5-Lemma $\bar{\kappa}$ is an isomorphism. It follows that $\kappa$ is an isomorphism. Now, apply the Short 5-Lemma (\ref{thm:Short5}) to the right hand side of the diagram to see that $u$ is an isomorphism. This was to be shown.
\end{proof}

This is closely related to another view on homological categories, as explained in Section~\ref{sec:Protomodular-SEpi(X)->X}; see, in particular, Definition~\ref{def:Protomodularity}. As a sibling to the pullback cancellation results (\ref{thm:Pullbacks,ConcatenatedSquares}) and (\ref{thm:Pullback-2-OutOf-3}), we have:

\begin{theorem}[2-out-of-3 property for pullbacks II\HTag]
  \label{thm:SAPullbackCancellation-II}%
  \cite[p.~242]{FBorceuxDBourn2004}, \cite[p.~36f]{Borceux-Semiab}, \cite{Janelidze-Sobral-Tholen}\quad In the commutative diagram below let $e$ be a normal epimorphism. %
  \index{pullback!cancellation $2$-out-of-$3$}%
  \begin{equation*}
    \xymatrix@R=5ex@C=3em{
    A \ar[r] \ar[d] &
    B \ar[r] \ar@{-{ >>}}[d]^-{e} &
    C \ar[d] \\
    X \ar[r] &
    Y \ar[r] &
    Z
    }
  \end{equation*}
  If two out of three of the commutative rectangles in this diagram are pullbacks, then so is the third.
\end{theorem}
\begin{proof}
  By Proposition~\ref{thm:Pullbacks,ConcatenatedSquares} it only remains to consider the case where the outer rectangle and the square on the left are pullbacks. 	In the commutative diagram below, note that $\bar{c}$ is a normal epimorphism because $e$ is one such; see (\ref{thm:NormalEpis-Props-Normal}.\ref{thm:Cokernels-Props.gfNEpi->gNEpi}).
  \begin{equation*}
    \xymatrix@!0@R=7ex@C=4em{
    A \ar[dd]_{a=\bar{b}} \ar[rr]
    \ar[rd]_{\cong}^{x^{\ast}f} &&
    B \ar@{-{ >>}}[dd]_(.7){e} \ar[rr] \ar[rd]^{f} &&
    C \ar@{=}[rd] \ar[dd] \\
    & x^{\ast}y^{\ast}C \ar[ld]_{\hat{c}} \ar[rr]^(.3){\bar{x}}|\hole &&
    y^{\ast}C \ar[ld]_{\bar{c}} \ar[rr]|\hole &&
    C \ar[ld]^{c} \\
    X \ar[rr]_-{x} &&
    Y \ar[rr]_-{y} &&
    Z
    }
  \end{equation*}
  Now we have two normal epimorphisms to pull back along $x$:
  \begin{enumerate}[(1)]
    \item $e$ pulls back to $a$ because the left hand square is a pullback.
    \item $\bar{c}$ pulls back to $\hat{c}$.
  \end{enumerate}
  The map $x^{\ast}f$ is an isomorphism because both, the front outer diagram and the bottom outer diagram are pullbacks along one and the same map $yx$. But then $f$ is an isomorphism because (\ref{thm:BaseChangeReflectsIsos}) pulling the normal epimorphisms $e$ and $\bar{c}$ back along $x$ reflects isomorphisms.
\end{proof}

We close this section by taking a first step toward the construction of a fundamental invariant used in homological algebra: the extension functor. We begin by introducing relevant terminology.

\begin{terminology}[Extension / morphism of extensions]
  \label{def:Extension}
  Recall (\ref{def:ShortExactSequence-Basic}) that given a short exact sequence
  \begin{equation*}
    \xymatrix@R=5ex@C=2em{
    \varepsilon &&
    K \ar@{{ |>}->}[rr] &&
    U \ar@{-{ >>}}[rr] &&
    V
    }
  \end{equation*}
  we refer to $\varepsilon$ as an \Defn{extension of $K$ over $Q$}, or a \Defn{cover of $Q$ by $K$}.
  \index{extension}\index{extension!morphism}\index{morphism!of extensions}

  A morphism from an extension $\varepsilon$ to an extensions $\mu$ is given by a morphism of the underlying short exact sequences:
  \begin{equation*}
    \xymatrix@R=5ex@C=2em{
    \varepsilon \ar[d]_{(k,u,v)} &&
    K \ar@{{ |>}->}[rr] \ar[d]_{k} &&
    U \ar@{-{ >>}}[rr] \ar[d]_{u} &&
    V \ar[d]^{v} \\
    \mu &&
    L \ar@{{ |>}->}[rr] &&
    X \ar@{-{ >>}}[rr] &&
    V
    }
  \end{equation*}
\end{terminology}

\begin{definition}[Equivalent extensions]
  \label{def:Extension-Equivalence}
  Two extensions $\varepsilon$ and $\mu$ of $K$ over $Q$ are \Defn{equivalent} if they admit a morphism of the form $(\IdMapOn{K},\beta,\IdMapOn{Q})$: %
  \index{equivalent!extensions}%
  \begin{equation*}
    \xymatrix@R=5ex@C=2em{
    \varepsilon \ar[d]_{(\IdMapOn{K},u,\IdMapOn{V})} &&
    K \ar@{{ |>}->}[rr] \ar[d]_{\IdMapOn{K}} &&
    X \ar@{-{ >>}}[rr] \ar[d]_{u} &&
    V \ar[d]^{\IdMapOn{V}} \\
    \mu &&
    K \ar@{{ |>}->}[rr] &&
    Y \ar@{-{ >>}}[rr] &&
    V
    }
  \end{equation*}
\end{definition}

In (\ref{def:Extension-Equivalence}), the map $u$ is an isomorphism by the Short $5$-Lemma (\ref{thm:Short5}). Thus an equivalence of extensions provides a very selected form of isomorphism between their middle objects, resulting in an equivalence relation on the class of all extensions of $K$ and over $Q$. This leads to:

\begin{definition}[Extension equivalence class\HTag]
  \label{def:ExtensionObject}
  Given objects $K$ and $Q$, we define %
  \index[not]{e!$\ExtMini{V}{K}$\IndSep equivalence classes of extensions}
  \begin{equation*}
    \ExtMini{V}{K} \DefEq \Set{\text{extensions of $K$ over $V$} } / \text{`equivalence'}
  \end{equation*}
\end{definition}

Let's analyze what happens in the definition of $\ExtMini{V}{K}$: Often times, the collection of extensions of $K$ over $V$ forms a proper class on which `equivalence of extensions' is an equivalence relation. Even in the category $\AbGrps$ of abelian groups, equivalence classes are proper classes themselves. So, strictly speaking, $\ExtMini{V}{K}$ is a conglomerate. In many varieties, it is possible to choose representatives of the equivalence classes in $\ExtMini{V}{K}$ which form a proper set. In those cases, $\ExtMini{V}{K}$ may be treated as a proper set, and
\begin{equation*}
  \ExtMini{-}{K}\from \Ctgry{X} \to \SetsBsd\from V \mapsto \ExtMini{V}{K}.
\end{equation*}
is a contravariant functor. - For a summary of the early history of Ext we refer the reader to the notes in Mac Lane \cite[p. 103]{SMacLane1995}.

\begin{subordinate}{}
  \begin{subsubordinate}{Retractionable normal monomorphism vs.\ sectionable normal epimorphism}
    Proposition~\ref{thm:ProductRecognition-SplitKer->Product} says that whenever a normal monomorphism admits a left inverse, then this left inverse is a product projection. Recalling that the notion of short exact sequence is self-dual, let's examine this!

    The dual of (\ref{thm:ProductRecognition-SplitKer->Product}) is: `If a normal epimorphism admits a right inverse, then this right inverse is a product inclusion'. While true in abelian categories, in general, it is false in the category $\Grps$ of groups. Why? - A right inverse of $f$ is a monomorphism which, in a semiabelian environment, \emph{is not necessarily normal}. Consequence: if $f$ admits a right inverse $i$ which happens to be a kernel, then $i$ is a product inclusion\dots\ as may be seen by dualizing the proof of (\ref{thm:ProductRecognition-SplitKer->Product}).
  \end{subsubordinate}
\end{subordinate}

\begin{exercises}

\begin{exercise}[Biproduct recognition]
  \label{exe:Sum=Product-If-SummandInclusionNormal} 
  \label{exe:BiProduct-Recognition}%
  Suppose in a semiabelian category $\SACtgry{X}$, the inclusion of an object $X$ into the sum $X+Y$ is normal. Show that the canonical map $X+Y\to \Prdct{X}{Y}$ is an isomorphism. %
  \index{biproduct!recognition}
\end{exercise}

\begin{exercise}
  Show that these two extensions in the category $\AbGrps$ of abelian groups are equivalent by filling in the missing arrow in the middle.
  \begin{equation*}
    \xymatrix@R=5ex@C=4em{
    \ZMod{2} \ar@{=}[d] \ar@{{ |>}->}[r]^-{\Dgnl} &
    \ZMod{2}\oplus \ZMod{2} \ar@{-{ >>}}[r]^-{+} \ar@{.>}[d]|-{?} &
    \ZMod{2} \ar@{=}[d] \\
    \ZMod{2} \ar@{{ |>}->}[r]_-{\PrdctMapInto{\IdMap,0}} &
    \ZMod{2}\oplus \ZMod{2} \ar@{-{ >>}}[r]_-{\PrjctnOnto{2}} &
    \ZMod{2}
    }
  \end{equation*}
\end{exercise}

\begin{exercise}[$\textit{Ext}(\ZMod{p},\ZMod{p})$]
  \label{exe:Ext(Z/p,Z/p)}
  For a prime number $p\in \ZNr$, show that there are $p$ distinct equivalence classes of abelian group extensions of $\ZMod{p}$ over $\ZMod{p}$.
\end{exercise}

\begin{exercise}[Inequivalent extensions with isomorphic middle objects]
  \label{exe:Extensions-SameMiddle-Inequivalent}
  In each of the categories familiar to you, e.g., $\AbGrps$, $\Grps$, \dots find inequivalent extensions with isomorphic middle objects; see the remark after (\ref{def:Extension-Equivalence}).
\end{exercise}

\begin{exercise}[Short $5$-lemma in $\SetsBsd$]
  \label{exe:Short5-In-Set_*}
  Compare with Exercise \ref{exe:Set_*CoSemiabelian}. Show that the short $5$-lemma holds in the category $\SetsBsd$ of pointed sets, even though $\SetsBsd$ does \emph{not} satisfy the the \KSGInline-condition, hence is not a homological category.
\end{exercise}
\end{exercises}
\section[The \texorpdfstring{$(\Prdct{3}{3})$}{(3x3)}-Lemma]{The \texorpdfstring{$(\Prdct{3}{3})$}{(3x3)}-Lemma}
\label{sec:3x3-Lemma-Homological}

In Section~\ref{sec:DinversionPreservesNormalMaps}, we already found necessary and sufficient conditions under which the four border cases of the $\Prdct{3}{3}$-Lemma hold. In the axiomatically richer environment of a homological category, we show that the `middle cases' of the $(\Prdct{3}{3})$-Lemma are valid as well; see (\ref{thm:(3x3)-LemmaSemiAb}).

Given a di-extension, we then explain in (\ref{thm:Di-Extension-InternalStructure})  how its initial, central, and terminal objects are related.

\begin{theorem}[$(\Prdct{3}{3})$-Lemma\HTag]
  \label{thm:(3x3)-Lemma-H}
  \label{thm:(3x3)-LemmaSemiAb}
  \cite[p.~45f]{Borceux-Semiab}, \cite[p.~279]{FBorceuxDBourn2004}\quad In the commutative diagram below, assume that each column is a short exact sequence, and that $eb=\ZeroMap$. %
  \index{$(3\prdct 3)$-Lemma}%
  \begin{equation}
    \label{eq:(3x3)-Lemma}%
    \vcenter{\xymatrix@R=5ex@C=3em{
    M \ar@{}[rd]|-{\text{(NM)}} \ar[r]^-{a} \ar@{{ |>}->}[d]_-{u} &
    L \ar[r]^-{d} \ar@{{ |>}->}[d]^-{v} &
    I \ar@{{ |>}->}[d]^-{w} \\
    K \ar[r]_-{b} \ar@{-{ >>}}[d]_-{x} &
    X \ar@{}[rd]|-{\text{(NE)}} \ar[r]^-{e} \ar@{-{ >>}}[d]_-{y} &
    Q \ar@{-{ >>}}[d]^-{z} \\
    J \ar[r]_-{c} &
    R \ar[r]_-{f} &
    S }}
  \end{equation}
  If any two of the rows are exact, then so is the third.
\end{theorem}
\begin{proof}
  \emph{Border cases: the middle row is known to be short exact}\quad From (\ref{thm:Homological->DPN}), we know that homological categories enjoy the \DPNInline-property. So, if the middle row is exact, then the claim follows from (\ref{thm:DPN-PreservationNormalMaps-Border(3x3)}).

  \emph{Middle case: the top and bottom rows are known to be short exact}\quad Since $d$ and $f$ are normal epimorphisms, (\ref{thm:NormalPushOut-Recognize-H}) shows that (NE) is a normal pushout. So, $e$ is a normal epimorphism. Thus, we find ourselves in this situation:
  
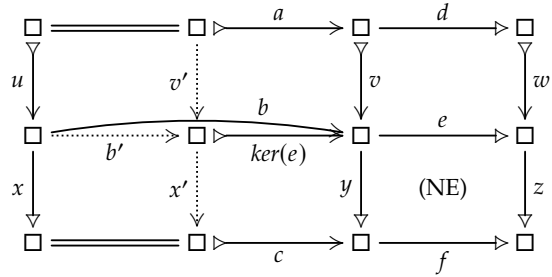
\begin{equation*}
    \xymatrix@R=5ex@C=4em{
    \DiagObj \ar@{=}[r] \ar@{{ |>}->}[d]_{u} &
    \DiagObj \ar@{{ |>}->}[r]^-{a} \ar@{.>}[d]_-{v'} &
    \DiagObj \ar@{-{ >>}}[r]^-{d} \ar@{{ |>}->}[d]^-{v} &
    \DiagObj \ar@{{ |>}->}[d]^{w} \\
    \DiagObj \ar@{.>}[r]_-{b'} \ar@{-{ >>}}[d]_{x} \ar@/^1ex/[rr]^(.7){b} &
    \DiagObj \ar@{{ |>}->}[r]_-{\KerMap{e}} \ar@{.>}[d]_-{x'} &
    \DiagObj \ar@{}[rd]|-{\text{(NE)}} \ar@{-{ >>}}[r]^-{e} \ar@{-{ >>}}[d]_-{y} &
    \DiagObj \ar@{-{ >>}}[d]^-{z} \\
    \DiagObj \ar@{=}[r] &
    \DiagObj \ar@{{ |>}->}[r]_-{c} &
    \DiagObj \ar@{-{ >>}}[r]_-{f} &
    \DiagObj
    }
  \end{equation*}
  Via the kernel property of $\KerMap{e}$ we obtain maps $v'$ and $b'$, unique with the properties
  \begin{equation*}
    va = \KerMap{e} v' \quad\text{and}\quad \KerMap{e}b' =b.
  \end{equation*}
  Here, we need the hypothesis $eb=\ZeroMap$. Finally, the kernel property of $c$ yields $x'$, unique with $cx' = y\KerMap{e}$. A border case of the $(\Prdct{3}{3})$-Lemma shows that the sequence $x'v'$ is short exact.

  We claim that the two squares on the left commute. Indeed,
  \begin{equation*}
    \KerMap{e}b'u =bu = va = \KerMap{e}v'.
  \end{equation*}
  So, the monomorphic property of $\Ker{e}$ shows that $v' = b'u$. Moreover,
  \begin{equation*}
    cx = yb = y \KerMap{e} b' = c x' b'
  \end{equation*}
  So, the monomorphic property of $c$ shows that $x = x' b'$. By the Short 5-Lemma (\ref{thm:Short5}), $b'$ is an isomorphism, and so we see that the middle row of the original diagram is short exact. - This completes the proof.
\end{proof}

Next, we explain how, in a di-extension, the initial, central, and terminal objects are related. This involves the meet and join of the normal subobjects of the central object.

\begin{proposition}[Internal structure of a di-extension\HTag]
  \label{thm:9Diagram-JntEpiOnKernel,JntMonoOnCoKernel}
  \label{thm:Di-Extension-InternalStructure}
  If $K$, $L\normal X$ are normal subobjects that fit a di-extension, such as \eqref{eq:(3x3)-Lemma}, then the sequences below are short exact.
  \begin{equation*}
    \xymatrix@R=1ex@C=4em{
    K\meet L \ar@{{ |>}->}[r]^-{va=bu} &
    X \ar@{-{ >>}}[r]^-{(y,e)} &
    R\times_{S}Q \\
    K\join L \ar@{{ |>}->}[r] &
    X \ar@{-{ >>}}[r]_-{fy=ze} &
    S
    }
  \end{equation*}
  In particular, $K\meet L = M$ and $K\join L=\Ker{fy}$ are normal subobjects of $X$.
\end{proposition}
\begin{proof}
  To establish the first of the short exact sequences, note that $M=K\meet L$ since the top left square in~\eqref{eq:(3x3)-Lemma} is a pullback. Next, construct the commutative diagram below by pulling the vertical sequence on the right of \eqref{eq:(3x3)-Lemma} back along $f$.
  \begin{equation*}
    \xymatrix@R=5ex@C=4em{
    M \ar@{{ |>}->}[r]^-{a} \ar@{=}[d] &
    L \ar@{{ |>}->}[d]_-{v} \ar@{-{ >>}}[r]^-{d} \PullLU{rd} &
    I \ar@{{ |>}->}[d] \ar@{=}[r] &
    I \ar@{{ |>}->}[d]^{w} \\
    M \ar@{{ |>}->}[r] &
    X \ar@{-{ >>}}[r]_-{\gamma} \ar@{-{ >>}}[d]_{y} &
    R\times_S Q \ar@{-{ >>}}[d] \ar@{-{ >>}}[r] \PullLU{rd}&
    Q \ar@{-{ >>}}[d]^{z} \\
    &
    R \ar@{=}[r] &
    R \ar@{-{ >>}}[r]_-{f} &
    S
    }
  \end{equation*}
  Then the comparison map $\gamma$ induced by $e$ and $y$ is a normal epimorphism as the square (NE) in (\ref{eq:(3x3)-Lemma}) is a regular pushout. The top middle square is then a pullback by (\ref{thm:PullbackRecognition-KernelSide}). Consequently,  $\Ker{\gamma}= \Ker{d}= M= K\meet L$ by (\ref{thm:Pullback->IsoOfKernels}).

  To prove that the second sequence is exact, first consider the morphism of short exact sequences
  \begin{equation*}
    \xymatrix@R=5ex@C=3em{
    K\meet L \ar@{{ |>}->}[r]^-{a} \ar@{{ |>}->}[d]_-{u} \PullLU{rd} &
    L \ar@{-{ >>}}[r]^-{d} \ar@{{ |>}->}[d]^-{v'} &
    I \ar@{.>}[d]^i \\
    K \ar@{{ |>}->}[r]_-{b'} &
    K\join L \ar@{-{ >>}}[r]_-{q} &
    (K\join L)/K
    }
  \end{equation*}
  where $v'$ and $b'$ are the induced factorizations of $v$ and $b$ through the join. Further, $q$ is the cokernel of the normal monomorphism $b'$. The induced map $i$ is a monomorphism by (\ref{thm:PullbackRecognition-KernelSide}), and is a normal epimorphism by (\ref{thm:NormalEpiRecognition-CoKernelSide}). So, $i$ is an isomorphism.  Now consider the $(3\prdct 3)$-diagram
  \begin{equation*}
    \xymatrix@R=5ex@C=3em{
    K \ar@{=}[d] \ar@{{ |>}->}[r]^-{b'} &
    K\join L \ar@{.>}[d] \ar@{-{ >>}}[r]^-{q} &
    (K\join L)/K \ar[d]^{wi} \\
    K \ar[d] \ar@{{ |>}->}[r]^b &
    X \ar@{-{ >>}}[r]^{e} \ar@{.>}[d]^{ze} &
    Q \ar@{-{ >>}}[d]^{z} \\
    0 \ar[r] & S \ar@{=}[r] & S
    }
  \end{equation*}
  and use the $(\Prdct{3}{3})$-Lemma (\ref{thm:(3x3)-LemmaSemiAb}) to conclude that its middle column is a short exact sequence.
\end{proof}

\begin{corollary}[The join in a di-extension\HTag]
  \label{thm:9Diagram-JntEpiOnKernel,JntMonoOnCoKernel-ExactSequence}%
  If $K$ and $L$ are normal subobjects of $X$,  that fit a di-extension such as~\eqref{eq:(3x3)-Lemma}, then the two exact sequences of Proposition~\ref{thm:9Diagram-JntEpiOnKernel,JntMonoOnCoKernel} may be spliced to form the exact sequence
  \begin{equation*}
    \xymatrix@R=5ex@C=3em{
    0 \ar[r] &
    K\meet L \ar@{{ |>}->}[r] & K\join L \ar[r] &
    R\times_{S}Q \ar@{-{ >>}}[r] & S \ar[r] & 0.
    }
  \end{equation*}
\end{corollary}
\begin{proof}
  This is an application of the Pure Snake Lemma \ref{thm:SnakeLemma-Pure}.
\end{proof}

\begin{corollary}[$2$-step extensions from di-extension\HTag]
  \label{thm:2-StepExtensionsFromDoubleExtension}
  \label{thm:2-StepExtensionsFromDiExtension}
  Associated to a di-extension such as \eqref{eq:(3x3)-Lemma} is this commutative diagram of $2$-step extensions:
  \begin{equation*}
    \xymatrix@R=5ex@C=3em{
    0 \ar[r] &
    M \ar@{{ |>}->}[r] \ar@{=}[d] &
    L \ar[r]^-{wd} \ar@{{ |>}->}[d]_{i_{L}} \ar@{}[rd]|-{\text{(I)}} &
    Q \ar@{-{ >>}}[r] &
    S \ar@{=}[d]  \ar[r] &
    0 \\
    0 \ar[r] &
    M \ar@{{ |>}->}[r] \ar@{=}[d]  &
    K \join L \ar[r]^-{\xi} \ar@{}[rd]|-{\text{(II)}} &
    Q\times_S R \ar@{-{ >>}}[r] \ar@{-{ >>}}[u]_{\PrjctnOnto{Q}} \ar@{-{ >>}}[d]^{\PrjctnOnto{R}} &
    S \ar@{=}[d] \ar[r] &
    0\\
    0 \ar[r] &
    M \ar@{{ |>}->}[r] &
    K \ar[r]_-{cx} \ar@{{ |>}->}[u]^{b'} &
    R \ar@{-{ >>}}[r] &
    S \ar[r] &
    0
    }
  \end{equation*}
\end{corollary}
\begin{proof}
  The $2$-step extension in the middle comes from applying the Pure Snake Lemma \ref{thm:SnakeLemma-Pure} to this morphism of short exact sequences:
  \begin{equation*}
    \xymatrix@R=5ex@C=3em{
    0 \ar[r] &
    M=K\meet L \ar@{{ |>}->}[r]^-{va=bu} \ar@{{ |>}->}[d] &
    X \ar@{-{ >>}}[r] \ar@{=}[d] &
    Q\times_{S}R \ar@{-{ >>}}[d] \ar[r] & 0
    \\
    0 \ar[r] &
    K\join L \ar@{{ |>}->}[r] &
    X \ar@{-{ >>}}[r]_-{ze=fy} &
    S \ar[r] & 0
    }
  \end{equation*}
  We know that $\xi= (e,y)\Comp i$. This implies the commutativity of square (I) and (II).
\end{proof}

\begin{corollary}[Chinese Remainder Theorem\HTag]
  \label{thm:ChineseRemainderThm}
  If $K$\*, $L\normal X$ are such that $K\join L=X$, then there is a functorial isomorphism
  \[
    \frac{X}{K\meet L}\cong \frac{X}{L}\times \frac{X}{K}.
  \]
\end{corollary}
\begin{proof}
  As in the proof of (\ref{thm:Noether}), we obtain the di-extension
  \begin{equation*}
    \xymatrix@R=5ex@C=3em{
    K\meet L \ar@{{ |>}->}[r] \ar@{{ |>}->}[d] &
    L \ar[r] \ar@{{ |>}->}[d] &
    {L}/(K\meet L) \ar@{=}[d] \\
    K \ar@{{ |>}->}[r] \ar@{-{ >>}}[d] &
    X \ar@{-{ >>}}[r] \ar@{-{ >>}}[d] &
    X/K \ar[d] \\
    K/(K\meet L) \ar@{=}[r] &
    X/L \ar[r] &
    0
    }
  \end{equation*}
  The result now follows from the first sequence in Proposition~\ref{thm:9Diagram-JntEpiOnKernel,JntMonoOnCoKernel}.
\end{proof}

\begin{corollary}[Quotient of join is product\HTag]
  \label{thm:QuotientOfJoin}
  If $K$\*, $L\normal K\join L\leq X$, then there is a functorial isomorphism
  \[
    \frac{K\join L}{K\meet L}\cong \frac{K}{K\meet L}\times \frac{L}{K\meet L}.
  \]
\end{corollary}
\begin{proof}
  This is a combination of (\ref{thm:ChineseRemainderThm}) and (\ref{thm:Noether}).
\end{proof}

\begin{subordinate}{On the join of two normal subobjects}
  In a homological category, if $K$ and $L$ are normal subobjects of an object $X$ which fit into a di-extension, then we found in  (\ref{thm:9Diagram-JntEpiOnKernel,JntMonoOnCoKernel}) that their join is again normal in $X$.

  This prompts the question: can we say anything about the join $K\join L$ of two arbitrarily given normal subobjects $K$ and $L$ of $X$? The answer is `yes', and we will come back to this topic in Section \ref{sec:LatticeOfSubjects}.
\end{subordinate}

\begin{exercises}

\begin{exercise}[Kernel pair of a normal epimorphism\HTag]
  \label{exe:KernelPairOfCoKernel}%
  Given a short exact sequence
  \begin{equation*}
    \xymatrix@R=1ex@C=3em{
    N \ar@{{ |>}->}[r]^-{\alpha} &
    X \ar@{-{ >>}}[r]^-{f} &
    Q
    }
  \end{equation*}
  show that $\Prdct{N}{N} \to \KrnlPr{f} \to Q$ is a short exact sequence as well.

\end{exercise}

\begin{exercise}[Product of short exact sequences\HTag]
  \label{exe:ProductSESs}%
  In a homological category $\Ctgry{X}$ do the following:
  \begin{enumerate}[(i)]
    \item Use the $(\Prdct{3}{3})$-Lemma (\ref{thm:(3x3)-Lemma-H}) to show that the product of two short exact sequences in a semiabelian category is short exact.
          \index{product!of short exact sequences}%
          \index{exactness!of product of short exact sequence }%
    \item For a fixed object $X$ in $\SACtgry{X}$, determine if the functor  $ \Prdct{X}{-}\from \SACtgry{X}\to \SACtgry{X}$  is exact.
    \item For a fixed object $X$, and a short exact sequence  $ A\XRA{\alpha} B \XRA{\beta} C$  in $\SACtgry{X}$, show that the sequences below are short exact.
          $$
            \Prdct{X}{A} \XRA{\IdMapOn{X}\times \alpha} \Prdct{X}{B} \XRA{\beta\PrjctnOnto{B}} C \quad\text{and}\quad A \XRA{(0,\alpha)} \Prdct{X}{B} \XRA{\IdMapOn{X}\times \beta} \Prdct{X}{C}
          $$
  \end{enumerate}
\end{exercise}

\begin{exercise}[Third Isomorphism Theorem from $(\Prdct{3}{3})$-Lemma\HTag]
  \label{exe:ThirdIsoThm}%
  If $K$ and $N$ are normal subobjects of $X$ such that $K\leq N$, then $N/K$ is normal in $X/K$, and the sequence below is short exact. %
  \index{Isomorphism Theorem!Third}\index{Third Isomorphism Theorem}%
  \begin{equation*}
    \xymatrix@R=5ex@C=4em{
    N/K \ar@{{ |>}->}[r] &
    X/K \ar@{-{ >>}}[r] &
    X/N
    }
  \end{equation*}
\end{exercise}

\begin{exercise}[Image / kernel in $\SESCat{X}$\HTag]
  \label{thm:Image/Kernel-In-SES(X)}
  In the category $\SESCat{X}$ associated to a homological category $\Ctgry{X}$ the following conditions are equivalent for a morphism $(\kappa,\xi,\rho)\from \varepsilon\to \varphi$:
  \begin{enumerate}[(i)]
    \item $(\kappa,\xi,\rho)$ admits a pointwise image factorization in $\SESCat{X}$.
    \item $(\kappa,\xi,\rho)$ has a pointwise kernel in $\SESCat{X}$.
    \item The sequence $\Ker{\kappa}\to \Ker{\xi}\to \Ker{\rho}$ in $\Ctgry{X}$ is short exact.
  \end{enumerate}
\end{exercise}

\begin{exercise}[Normal closure of image / Cokernel in $\SESCat{X}$\HTag]
  \label{exe:NormalClosureImage/CoKernel-In-SES(X)}
  In the category $\SESCat{X}$ associated to a homological category $\SACtgry{X}$ the following conditions are equivalent for a morphism $(\kappa,\xi,\rho)\from \varepsilon\to \varphi$:
  \begin{thmlist}
    \item The sequence of normal closures $\overline{\Img{\kappa}}\to \overline{\Img{\xi}}\to \overline{\Img{\rho}}$ in $\Ctgry{X}$ is short exact.
    \item $(\kappa,\xi,\rho)$ has a pointwise cokernel in $\SESCat{X}$.
  \end{thmlist}
\end{exercise}

\begin{exercise}[Normal image / kernel / cokernel in $\SESCat{X}$\HTag]
  \label{thm:NormalImage/Kernel/CoKernel-In-SES(X)}
  In the category $\SESCat{X}$ associated to a homological category $\Ctgry{X}$ the following conditions are equivalent for a morphism $(\kappa,\xi,\rho)\from \varepsilon\to \varphi$ for which $\kappa,\xi,\rho$ have normal images in $\Ctgry{X}$.
  \begin{thmlist}
    \item $(\kappa,\xi,\rho)$ admits a pointwise image factorization in $\SESCat{X}$.
    \item $(\kappa,\xi,\rho)$ admits a pointwise kernel in $\SESCat{X}$.
    \item $(\kappa,\xi,\rho)$ has a pointwise cokernel in $\SESCat{X}$.
    \item The sequence $\Ker{\kappa}\to \Ker{\xi}\to \Ker{\rho}$ in $\Ctgry{X}$ is short exact.
  \end{thmlist}
\end{exercise}
\end{exercises}
\section[Commuting Morphisms]{Commuting  Morphisms}
\label{sec:CommutingMorphisms}%

Two subgroups $K$ and $L$ of a group $G$ \emph{commute} if the identity $kl=lk$ holds for every $k\in K$ and $l\in L$. In this section make the notion of commuting subobjects and, more generally, of commuting maps, categorical. This provides a foundation for the study of abelian group objects in a homological category, see Section \ref{sec:CommutativeObjects-H}, and enables a study of central subobjects later on. The approach taken here is based on the work of Bourn~\cite{DBourn1996,DBourn2002}, Huq~\cite{Huq} and Borceux--Bourn~\cite{FBorceuxDBourn2004}.

Recall (\ref{sec:Surjectivity}) that a pair of maps $ (r\from A\to X, s\from B\to X)$ in a category $\Ctgry{C}$ is \Defn{jointly extremal-epimorphic} provided in any commutative diagram
\index{jointly!extremal-epimorphic family of maps}%
\begin{equation*}
  \xymatrix@!0@=4em{
  & M \ar@{ >->}[d]^- m \\
  A \ar[r]_-{r} \ar[ur] &
  X & B \ar[l]^-{s} \ar[ul]
  }
\end{equation*}
in which $m$ is a monomorphism, it is necessarily an isomorphism. %
The intuition is that `$r$~and $s$ together generate all of $X$'. If $\Ctgry{C}$ has binary sums, then this is equivalent to saying that the arrow $\SumMapOutOf{r,s}\from A+B\to X$ is an extremal epimorphism.


A pointed category with finite limits is said to be \Defn{unital} when for any pair of objects $A$, $B$, the product inclusions
\begin{equation*}
  \xymatrix@R=6ex@C=4em{
  A \ar[r]^-{(\IdMapOn{A},\ZeroMap)} &
  \Prdct{A}{B} &
  B \ar[l]_-{(\ZeroMap,\IdMapOn{B})}
  }
\end{equation*}
are jointly extremal-epimorphic. This is the case, for example, in the category of groups, the category of rings, and categories of modules over a ring, as is easily checked by hand. In fact, all homological categories are unital, as follows immediately from the \KSGInline-axiom (Definition~\ref{def:HomologicalCategory}). Here, the condition is equivalent to asking that the comparison map $\SumProdComp{A}{B}\from A+B\to \Prdct{A}{B}$ is an extremal epimorphism, for all $A$ and $B$.

In unital categories, the notion of commuting maps is well defined:

\begin{definition}[Commuting morphisms]
  \label{def:CommutingMaps}
  Two maps $f\from A\to X$ and $g\from B\to X$ in a unital category \Defn{commute} if there exists a morphism $\mu_{f,g}$ which renders the diagram
  \index{commutating!morphisms}
  \begin{equation*}
    \xymatrix@R=6ex@C=4em{
    A \ar[r]^-{(\IdMapOn{A},\ZeroMap)} \ar@/_2ex/[dr]_{f} &
    \Prdct{A}{B} \ar[d]^-{\mu_{f,g}} &
    B \ar[l]_-{(\ZeroMap,\IdMapOn{B})}  \ar@/^2ex/[dl]^{g}\\
    & X
    }
  \end{equation*}
  commutative. The map $\mu_{f,g}$ is called the \Defn{multiplication\footnote{In \cite{FBorceuxDBourn2004} this map is called the \Defn{collaborator} of $f$ and $g$.} of $f$ by $g$}. %
  \index{multiplication!of $f$ by $g$}%
\end{definition}

Using the comparison map $\SumProdComp{A}{B}\from A+B\to \Prdct{A}{B}$ in the context of a homological category, the maps $f$ and $g$ in (\ref{def:CommutingMaps}) commute if and only if the universal map $\SumMapOutOf{f,g}\from A+B\to X$ factors through $\SumProdComp{A}{B}$ as shown in the square below, which is a seminormal pushout by (\ref{thm:SemiNormalPushout-Recognition}).
\begin{equation*}
  \xymatrix@R=5ex@C=4em{
  A + B \ar@{-{ >>}}[r]^-{\SumProdComp{A}{B}} \ar[d]_{\SumMapOutOf{f,g}} &
  \Prdct{A}{B} \ar[d]^{\mu_{f,g}} \\
  X \ar@{=}[r] &
  X
  }
\end{equation*}

In Proposition \ref{thm:CommutingMaps-Props} and its corollaries (\ref{thm:CommutingMaps-CommutingImages}) and (\ref{thm:ImageOfCommutingSubobj}) we present basic properties of commuting maps.

\begin{proposition}[Properties of commuting maps]
  \label{thm:CommutingMaps-Props}
  Given commuting maps $f\from A\to X$ and $g\from B\to X$ the following hold:
  \begin{enumerate}[(i)]
    \item The multiplication map $\mu_{f,g}$ of $f$ by $g$ is unique.
    \item If $\varphi\from A'\to A$ and $\psi\from B'\to B$ are arbitrary maps, then $f\varphi$ commutes with $g\psi$, and $\mu_{f\varphi,g\psi} = \mu_{f,g}\Comp (\Prdct{\varphi}{\psi})$.
    \item If $\xi\from X\to Y$ is arbitrary, then $\xi f$ commutes with $\xi g$, and $\mu_{\xi f,\xi g} = \xi\mu_{f,g}$.
    \item If $u\from R\to Y$ commutes with $v\from S\to Y$, then $\Prdct{f}{u}\from \Prdct{A}{R}\to \Prdct{X}{Y}$ commutes with $\Prdct{g}{v}\from \Prdct{B}{S}\to \Prdct{X}{Y}$. If $\tau_{R,B}\from \Prdct{R}{B}\to \Prdct{B}{R}$ denotes the twist map, then $\mu_{\Prdct{f}{u},\Prdct{g}{v}} = (\Prdct{\mu_{f,g}}{\mu_{u,v}})\Comp (\IdMapOn{A}\prdct \tau_{R,B}\prdct\IdMapOn{S})$.
  \end{enumerate}
\end{proposition}

\begin{proof}
  (i)\quad By assumption, the inclusions $(\IdMapOn{A},0)\from A\to \Prdct{A}{B}$ and $(0,\IdMapOn{B})\from B\to \Prdct{A}{B}$ are jointly epimorphic. So, if $\mu_{f,g}$ exists, then it is unique.

  (ii), (iii), (iv) are left as exercises.
\end{proof}

\begin{definition}[Central map]
  \label{def:CentralMap}
  A morphism $f\from A\to X$ in a unital category is called \Defn{central} if it commutes with the identity map on $X$. %
  \index{central!morphism}\index{morphism!central}%
\end{definition}

\begin{proposition}[Properties of central maps]
  \label{thm:CentralMaps-Props}
  For central maps $f\from A\to X$ and $g\from B\to X$ the following hold:
  \begin{enumerate}[(i)]
    \item $f$ commutes with $g$, and $\mu_{f,\IdMapOn{X}}\Comp (\Prdct{\IdMapOn{X}}{g} )= \mu_{f,g}=\mu_{\IdMapOn{X},g}\Comp (\Prdct{f}{\IdMapOn{X}})$.
    \item $\SumMapOutOf{f,g}\from A+B\to X$ is central and, if $H=(\Prdct{\IdMapOn{A}}{\tau_{B,X}})\Comp (\Prdct{\SumProdComp{A}{B}}{\IdMapOn{X}})$, then
          \begin{equation*}
            \mu_{f,\IdMapOn{X}}\Comp (\Prdct{\IdMapOn{A}}{\mu_{\IdMapOn{X},g}})\Comp H = \mu_{\SumMapOutOf{f,g},\IdMapOn{X}} = \mu_{\IdMapOn{X},g}\Comp (\Prdct{\mu_{f,\IdMapOn{X}}}{\IdMapOn{B}})
          \end{equation*}
  \end{enumerate}
\end{proposition}
\begin{proof}
  We give a compact proof involving coproducts; the general case is left as an exercise. A proof of (i) may be read off of this commutative diagram:
  \begin{equation*}
    \xymatrix@R=5ex@C=3em{
    \Prdct{A}{B} \ar[d]_{\IdMapOn{A}\prdct g}&
    A+B \ar@{-{ >>}}[l]_{\SumProdComp{A}{B}} \ar[d]^{\IdMapOn{A}+g}&
    A+B \ar@{=}[l] \ar@{=}[r] \ar[dd]|-{\SumMapOutOf{f,g}}&
    A+B \ar@{-{ >>}}[r]^-{\SumProdComp{A}{B}} \ar[d]_{f+\IdMapOn{B}} &
    \Prdct{A}{B} \ar[d]^{f\prdct \IdMapOn{B}} \\
    \Prdct{A}{X} \ar[d]_{\mu_{f,\IdMapOn{X}}} &
    A+X \ar@{-{ >>}}[l]_-{\SumProdComp{A}{X}} \ar[d]^{\SumMapOutOf{f,\IdMapOn{X}}} &&
    X+B \ar@{-{ >>}}[r]^-{\SumProdComp{X}{B}} \ar[d]_{\SumMapOutOf{\IdMapOn{X},g}} &
    \Prdct{X}{B} \ar[d]^{\mu_{\IdMapOn{X},g}} \\
    X \ar@{=}[r] &
    X \ar@{=}[r] &
    X \ar@{=}[r] &
    X \ar@{=}[r] &
    X
    }
  \end{equation*}
  To check the claim of (ii), consider this commutative diagram:
  \begin{equation*}
    \xymatrix@R=6ex@C=4em{
    (A+B)\prdct X \ar[r]^-{\SumProdComp{A}{B}\times \IdMapOn{X}} \ar@/_3ex/[rr]_-{H} &
    \PPrdct{A}{B}{X} \ar[r]^-{\IdMapOn{A}\prdct \tau_{B,X}} &
    \PPrdct{A}{X}{B} \ar[r]^-{\IdMapOn{A}\times \mu_{\IdMapOn{X},g}} \ar[d]_{\mu_{f,\IdMapOn{X}}\times \IdMapOn{B}} &
    \Prdct{A}{X} \ar[d]^{\mu_{f,\IdMapOn{X}}} \\
    && \Prdct{X}{B} \ar[r]_-{\mu_{\IdMapOn{X},g}} &
    X
    }
  \end{equation*}
  Using (i), we find that
  \begin{equation*}
    \mu_{\SumMapOutOf{f,g},\IdMapOn{X}}\Comp \PrdctMapInto{\IdMapOn{A+B},0} = \mu_{f,g}\Comp \SumProdComp{f}{g} = \SumMapOutOf{f,g}
  \end{equation*}
  Further, $\mu_{\SumMapOutOf{f,g},\IdMapOn{X}}\Comp (0,\IdMapOn{X})\from X\to (A+B)\prdct X$ equals $\IdMapOn{X}$. So, (ii) holds.
\end{proof}

Finally, we single out the special case where the identity map on an object commutes with itself:

\begin{definition}[Commutative object]
  \label{def:CommutativeObject}
  An object in a unital category is called \Defn{commutative} if $\IdMapOn{X}$ commutes with itself. %
  \index{commutative!object}\index{object!commutative}%
\end{definition}

A commutative object always carries the structure of a commutative unital magma which is unique:

\begin{proposition}[Commutative object is commutative monoid]
  \label{thm:CommutativeObject->CommutativeMonoid-H}
  In a unital category $\Ctgry{C}$, any commutative object carries the structure of an internal commutative monoid $(X,\mu)$.
\end{proposition}
\begin{proof}
  We treat the case where $\Ctgry{C}$ has coproducts. The general case is left as an exercise.

  If $\mu\from \Prdct{X}{X}\to X$ multiplies the identity map on $X$ by itself, then $\mu$ automatically satisfies the requirements for a unitary magma structure on $X$. To see that the multiplication $\mu$ is commutative, consider the twist map on $\Prdct{X}{X}$:
  \begin{equation*}
    \tau=(\PrjctnOnto{2},\PrjctnOnto{1})\from \Prdct{X}{X}\longrightarrow \Prdct{X}{X}
  \end{equation*}
  Then $\mu\tau$ is also a unitary magma structure on $X$. Such a structure is unique, and so $\mu=\mu\tau$; i.e., $\mu$ is commutative.

  To see that $\mu$ is associative, we want to know that the rectangle below commutes.
  \begin{equation*}
    \xymatrix@R=5ex@C=4em{
    X+(\Prdct{X}{X}) \ar@{-{ >>}}[r]^-{\gamma\DefEq\gamma_{X,\Prdct{X}{X}}} &
    X\prdct (\Prdct{X}{X}) \ar[r]^-{\alpha}_-{\cong} &
    X\prdct X \prdct X \ar[r]^-{\IdMap\prdct \mu} \ar[d]_{\mu\prdct\IdMap} &
    \Prdct{X}{X} \ar[d]^{\mu} \\
    && \Prdct{X}{X} \ar[r]_-{\mu} &
    X
    }
  \end{equation*}
  Since $\alpha\gamma$ is an epimorphism, it suffices to show that the composites
  \begin{equation*}
    \mu(\IdMap\prdct\mu)\alpha\gamma \qquad \text{and}\qquad \mu(\mu\prdct \IdMap)\alpha\gamma
  \end{equation*}
  are equal. This is, indeed, the case. On the summand $\Prdct{X}{X}$ we find that both composites equal~$\mu$, while on the summand $X$ we find that both composites equal $\IdMapOn{X}$. Thus the square commutes, and this proves that $\mu$ is associative; i.e., $(X,\mu)$ is an internal commutative monoid.
\end{proof}

\subsection{Commuting subobjects}
In the special case where commuting maps $f$ and $g$ represent $A$ and $B$ as subobjects, we say that the subobjects commute. %
\index{commuting!subobjects}%
This is legitimate,  since any two monomorphisms representing the same subobject are related by a unique isomorphism. The terminology of commuting maps / subobjects is motivated by examples such as (\ref{exa:CommutingSubgroups}), (\ref{exa:LieCommutingSubalgebras}), and (\ref{exa:AssociativeAlgs-CommutingSubobjects}). In the context of a homological category, commuting maps and commuting subobjects are related via the following results.

\begin{proposition}[Normal quotient of commuting maps\HTag]\label{thm:NormalQuotientOfCommutingMaps}
  Given commuting maps $f\from A\to X$ and $g\from B\to X$  in a homological category, if $f = f'p$ and $g=g'q$ with normal epimorphisms $p$ and $q$, then $f'$ commutes with $g'$, and $\mu_{f,g} = \mu_{f',g'}\Comp (\Prdct{p}{q}) $.\NoProof
\end{proposition}

\begin{corollary}[Commuting maps vs.\ commuting images\HTag]
  \label{thm:CommutingMaps-CommutingImages}
  Two maps $f\from A\to X$ and $g\from B\to X$ commute if and only if the image of $f$ commutes with the image of $g$.
\end{corollary}
\begin{proof}
  Using the image factorizations $A \XRA{p} I \XRA{i} X$ and $B\XRA{q} J \XRA{j} X$ of $f$ and $g$ respectively, the claim follows from (\ref{thm:NormalQuotientOfCommutingMaps}).
\end{proof}

\begin{corollary}[Direct images preserve commuting subobjects\HTag]
  \label{thm:ImageOfCommutingSubobj}
  Let $K$ and $L$ be commuting subobjects of $X$. If $f\from X\to Y$ is a normal epimorphism, then $f(K)$ and $f(L)$ commute as well. \NoProof
\end{corollary}

Let us now turn to situations in which two normal subobjects commute.

\begin{theorem}[Normal subobjects with zero intersection commute\HTag]
  \label{thm:couniv}
  \label{thm:NormalSubobjects-0Intersection-Commute}
  If two subobjects $K$ and $L$ of $X$ have intersection $0$ and are normal in their join, that is $K, L\normal K\join L\leq X$, then $K$ and $L$ commute.
\end{theorem}
\begin{proof}
  With $K\meet L=0$, (\ref{thm:QuotientOfJoin}) tells us that
  \begin{equation*}
    K\join L \cong \dfrac{K\join L}{K\meet L} \cong \dfrac{K}{K\meet L}\times \dfrac{L}{K\meet L} \cong \Prdct{K}{L}
  \end{equation*}
  So, if $k$ and $l$ represent $K$ and $L$ as subobjects of $X$, then $\SumMapOutOf{k,l}$ factors uniquely through $\Prdct{K}{L}$.
\end{proof}

\begin{corollary}[Quotient to make normal subobjects commute\HTag]
  \label{thm:quotient of intersection}
  Given normal subobjects $K$ and $L$ of $X$, let $q\from  X \to Y\DefEq X/(K\meet L)$ be the quotient map. Then $q(K)$ and $q(L)$ are commuting normal subobjects of $Y$.
\end{corollary}
\begin{proof}
  We show that $q(K)$ and $q(L)$ are normal subobjects of $Y$ whose intersection is $0$. So, the claim follows with (\ref{thm:NormalSubobjects-0Intersection-Commute}). Indeed, both of $K$ and $L$ contain $K\meet L$ as a subobject. So, (\ref{thm:Hofmann}) implies that $q(K)$ and $q(L)$ are normal subobjects of $Y$. Now construct the cube below from the pullback faces on the left and on the right, together with $q$ and its restrictions to $K$ and~$L$.
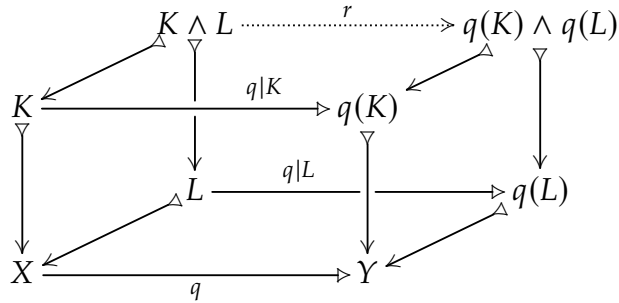
  \begin{equation*}
    \xymatrix@!@R=-6ex@C=0em{
    & K\meet L \ar@{.>}[rr]^-{r} 	\ar@{{ |>}->}[dd]|\hole \ar@{{ |>}->}[ld] &&
    q(K)\meet q(L) \ar@{{ |>}->}[ld] \ar@{{ |>}->}[dd] \\
    K \ar@{{ |>}->}[dd] \ar@{-{ >>}}[rr]^(.7){q|K} &&
    q(K) \ar@{{ |>}->}[dd] \\
    & L \ar@{{ |>}->}[ld] \ar@{-{ >>}}[rr]|\hole^(.3){q|L} &&
    q(L) \ar@{{ |>}->}[ld] \\
    X \ar@{-{ >>}}[rr]_-{q} &&
    Y
    }
  \end{equation*}
  Then the horizontal maps of the front and bottom faces have $K\meet L$ as a common kernel. So these faces are pullbacks by \ref{thm:PullbackFromKerIso}. The pullback property of the right hand face yields a map $r\from (K\meet L)\to (q(K)\meet q(L))$, which renders the entire diagram commutative. The back face is seen to be a pullback by combining pullback concatenation (\ref{thm:Pullbacks,ConcatenatedSquares}) and the $2$-of-$3$-property of pullbacks in a homological category (\ref{thm:Pullback-2-OutOf-3}). But then $r$ is a normal epimorphism by the \PNEInline-condition. Its image is $q(K)\meet q(L)=q(K\meet L)=0$. This implies the claim.
\end{proof}

\begin{example}[Commuting subgroups]
  \label{exa:CommutingSubgroups}
  If $K$ and $L$ are subgroups of a (multiplicative) group $X$, such that every element of $K$ commutes with every element of $L$, then the inclusion functions of $K$ and $L$ in $X$ conspire to define a group morphism
  \begin{equation*}
    \mu=\mu_{K,L}\from \Prdct{K}{L} \to X \from (k,l)\mapsto k\cdot l.
  \end{equation*}
  Conversely, if the inclusions of $K$ and $L$ into $X$ admit a multiplication map $\mu\from \Prdct{K}{L}\to X$, then
  \begin{equation*}
    k\cdot l = \mu(k,1)\cdot \varphi(1,l) = \mu\left( (k,1)\cdot (1,l) \right) =\mu\left( (1,l)\cdot (k,1) \right) = \mu(1,l)\cdot \mu(k,1) = l\cdot k.
  \end{equation*}
  So $K$ and $L$ commute in $X$. --- In short, `different factors in a product of groups commute with one another' is what we use for a categorical definition of commuting subgroups and commuting morphisms.

  Applied to the special case where $K=L=X$, we see that $X$ commutes with itself as a subobject if and only if the group operation on $X$ is commutative.
\end{example}

\begin{example}[Commuting Lie subalgebras]
  \label{exa:LieCommutingSubalgebras}
  Two sub-Lie-algebras $K$ and $L$ of a Lie algebra $X$ commute as subobjects of $X$ if and only if the bracket $[k,l]\in X$ vanishes whenever $k\in K$ and $l\in L$. To see this, consider first the case where the inclusions of $K$ and $L$ in $X$ admit a multiplication $\mu=\mu_{K,L}$. We then find
  \[
    \mu(k,l)=\mu((k,0)+(0,l))=\mu(k,0)+\mu(0,l)=k+l.
  \]
  Consequently,
  \begin{equation*}
    [k,l] = [\mu(k,0),\mu(0,l)]=\mu([k,0],[0,l])=\mu(0,0)=0.
  \end{equation*}
  On the other hand, if $[k,l]=0$ whenever $k\in K$ and $l\in L$, then the linear map defined by $\psi(k,l)=k+l$ is a morphism of Lie algebras, because
  \begin{align*}
    [\psi(k,l),\psi(k',l')] & =[k+l,k'+l']=[k,k']+[k,l']+[l,k']+[l,l'] \\
                            & =[k,k']+[l,l']=\psi([k,k'],[l,l'])       \\
                            & =\psi([(k,l),(k',l')]).
  \end{align*}
  So $\psi$ is the multiplication induced by $K$ and $L$.
\end{example}

\begin{example}[Commuting subobjects of associative algebras]
  \label{exa:AssociativeAlgs-CommutingSubobjects}
  If $K$ and $L$ are subalgebras of an associative algebra $X$, then subobjects $K$ and $L$ commute if and only if for each $k\in K$ and $l\in L$ we have that $kl=0$.		The argument is similar to the one given for Lie algebras.
\end{example}

\begin{subordinate}{}
  \begin{subsubordinate}{On the origin of `commuting morphisms'}
    In \cite{DBourn1996,DBourn2002}, D.~Bourn builds on earlier work of S.~A.~Huq~\cite{Huq}, extending definitions and results from a context which is essentially equivalent to that of semiabelian categories, to the new and far more general setting of unital categories.

    Indeed, the property of being a unital category is much weaker than being homological: a variety of universal algebras is unital if and only if it contains a unique constant $0$, and a binary operation $+$ satisfying the identities $x+0=x=0+x$; see~\cite{FBorceuxDBourn2004}.
  \end{subsubordinate}

  \begin{subsubordinate}{Weaker conditions}
    First of all, let us notice that the notions of commuting maps and commutative object are fully viable whenever the maps in (\ref{def:CommutingMaps}) are jointly epimorphic. In particular, Proposition \ref{thm:CommutativeObject->CommutativeMonoid-H} still holds. Such categories are called \emph{weakly unital} in \cite{NMFPhD}.
  \end{subsubordinate}
\end{subordinate}

\begin{exercises}
\begin{exercise}[Properties of central maps\HTag]
  \label{exe:CentralMaps-Props}
  For arbitrary objects $A$ and $X$, show the following:
  \begin{enumerate}[(i)]
    \item The zero map $\ZeroMap\from A\to X$ is central.
    \item $f\from A\to X$ is central if and only if $f$ commutes with every map $g\from B\to X$; compare (\ref{thm:CentralMaps-Props}).
    \item If $f\from A\to X$ is central, and $g\from X\to Y$ is a normal epimorphism, then $gf\from A\to Y$ is central.
  \end{enumerate}
\end{exercise}

\begin{exercise}[Normal closure of diagonal\HTag]
  \label{exe:Diagonal-NormalClosure}
  Given an object $X$, show that the inclusion of $X$ into the normal closure $\bar{X}$ of the diagonal in $\Prdct{X}{X}$ is a split monomorphism.
\end{exercise}

\begin{exercise}[Complete the proof of (\ref{thm:CentralMaps-Props}) and (\ref{thm:CommutativeObject->CommutativeMonoid-H}).]
  Write versions of these proofs which don't involve coproducts.
\end{exercise}
\end{exercises}
\newpage
\section[Commutative Objects in Homological Categories]{Commutative Objects in Homological Categories}
\label{sec:CommutativeObjects-H}%

In Section \ref{sec:CommutingMorphisms}, we introduced the concept of commutative object. Then we showed in (\ref{thm:CommutativeObject->CommutativeMonoid-H}) that it is always a commutative monoid. Here, we prove the remarkable result (\cite[3.D]{Bourn-Gran-CategoricalFoundations}, \cite{FBorceuxDBourn2004}) which says that, in a homological category, a commutative object is actually an internal abelian group; see (\ref{thm:AbelianGroupObject-Recognize-PSA}). For this reason, we shall refer to a commutative object in a homological category as an \emph{abelian object}.

To this end, we rely on the following lemma.

\begin{lemma}[Unitary magma structure is product projection\HTag]
  \label{thm:MultiplicationInUnitaryMagma-ProductProjection}
  If $(X,\mu)$ is a unitary magma structure in a homological category, then $(X\times X, \PrjctnOnto{1},\mu)$ is a product.
\end{lemma}
\begin{proof}
  We use the 2-out-of-3 property for pullbacks (\ref{thm:SAPullbackCancellation-II}) on this diagram:
  \begin{equation*}
    \xymatrix@R=5ex@C=3em{
    X \ar[r]_-{\PrdctMapInto{\ZeroMap,\IdMapOn{X}}} \ar[d] \ar@/^3ex/[rr]^-{\IdMapOn{X}} &
    X\times X \ar[r]_-{\mu} \ar@{-{ >>}}[d]^-{\PrjctnOnto{1}} &
    X \ar[d] \\
    0 \ar[r] &
    X \ar[r] &
    0
    }
  \end{equation*}
  The outer rectangle as well as the square on the left are pullbacks. Since $\PrjctnOnto{1}$ is a split epimorphism, it is a normal epimorphism, and so the theorem tells us that the square on the right is a pullback. So, $(\Prdct{X}{X}, \PrjctnOnto{1},\mu)$ is a product.
\end{proof}

Summarizing what we found about commutative objects in Section \ref{sec:CommutingMorphisms} combined with present developments, we obtain:

\begin{theorem}[Internal abelian group in a homological category\HTag]
  \label{thm:AbelianGroupObject-Recognize-PSA}%
  For an object $X$ in a homological category the following hold:
  \begin{thmlist}
    \item $X$ admits an internal unitary magma structure $(X,\mu)$ if and only if $\IdMapOn{X}$ commutes with itself.
    \item If $(X,\mu)$ is an internal unitary magma structure on $X$, then $\mu$ is commutative and associative, and there exists a unique map $i\from X\to X$ such that $(X,\mu,i)$ is an internal abelian group.
  \end{thmlist}%
  \index{internal!unitary magma structure}\index{internal!abelian group structure}%
\end{theorem}
\begin{proof}
  In view of Proposition \ref{thm:CommutativeObject->CommutativeMonoid-H}, it only remains to find an inversion $i\from X\to X$ so that $(X,\mu,i)$ is an internal abelian group. By Lemma~\ref{thm:MultiplicationInUnitaryMagma-ProductProjection}, $(X\times X,\PrjctnOnto{1},\mu)$ is a product. With (\ref{thm:ProductRecognition}) we obtain a unique map $\PrdctMapInto{\IdMapOn{X},i}\from X\to \Prdct{X}{X}$ so that $\PrjctnOnto{1}\Comp \PrdctMapInto{\IdMapOn{X},i} = \IdMapOn{X}$ and $\mu\Comp \PrdctMapInto{\IdMapOn{X},i} =\ZeroMap$. So, $i$ is a multiplicative right inverse for $\mu$. Since $\mu$ is commutative, $i$ is also a multiplicative left inverse. With (\ref{thm:InternalMonoid-InverseUnique}) we conclude that $i\from X\to X$ is the one and only two-sided inverse operation on the commutative monoid $X$ which turns it into an internal abelian group.
\end{proof}

\begin{definition}
  An \Defn{abelian object} in a homological category is an object $X$ which admits a (necessarily unique) internal abelian group structure $(X,\mu,i)$.
\end{definition}

By the above, any commutative object in a homological category is an abelian object.

\begin{corollary}[Intersection of commuting subobjects is abelian\HTag]
  \label{Example Abelian Intersection}
  \label{thm:IntersectionCommutingSubobjects-Abelian}
  If subobjects $K$ and $L$ of $X$ commute, then $K\meet L$ is an abelian object.
\end{corollary}
\begin{proof}
  The cube below commutes and explains why the multiplication induced by the inclusions $k$ and $l$ of $K$ and $L$ into $X$ give rise to a multiplication $\mu_{k\meet l,k\meet l}$ of the inclusion $k\meet l$ of $K\meet L$ into $X$ with itself.
  \begin{equation*}
    \xymatrix@!0@R=7ex@C=9ex{
    & (K\meet L)\times(K\meet L) \ar@{{ >}->}[rr] \ar@{.>}[dddl]|(0.335)\hole|-{\mu_{k\meet l,k\meet l}}
    &&
    \Prdct{K}{L} \ar[dddl]|-{\mu_{k,l}}\\
    K\meet L \ar@{{ >}->}[ru] \ar@{{ >}->}[dd] \ar@{{ >}->}[rr] &&
    K \ar@{{ >}->}[ru] \ar@{{ >}->}[dd]_(.3){k} \\
    & K\meet L \ar@{{ >}->}[uu]|\hole \ar@{{ >}->}[ld] \ar@{{ >}->}[rr]|\hole|(0.67)\hole &&
    L \ar@{{ >}->}[ld]^-{l} \ar@{{ >}->}[uu] \\
    X \ar@{=}[rr] &&
    X
    }
  \end{equation*}
  Proposition~\ref{thm:Sum->ProductIsCokernel} tells us that the vertical arrow on the left in the diagram
  
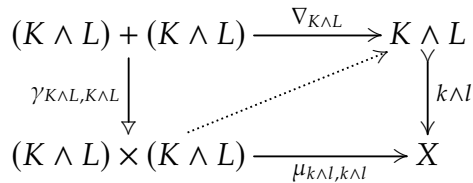
\begin{equation*}
    \xymatrix@R=5ex@C=4em{
    (K\meet L)+(K\meet L) \ar@{-{ >>}}[d]_-{\gamma_{K\meet L,K\meet L}} \ar[r]^-{\nabla_{K\meet L}} &
    K\meet L \ar@{{ >}->}[d]^-{k\meet l}\\
    (K\meet L)\times(K\meet L) \ar@{.>}[ru] \ar[r]_-{\mu_{k\meet l,k\meet l}} &
    X
    }
  \end{equation*}
  is a normal and, hence, a strong epimorphism. The induced lifting is a multiplication of $1_{K\meet L}$ with itself. Theorem \ref{thm:AbelianGroupObject-Recognize-PSA} implies that $K\meet L$ is an internal abelian group.
\end{proof}

Combining (\ref{thm:IntersectionCommutingSubobjects-Abelian}) with Corollary~\ref{thm:quotient of intersection} yields:

\begin{corollary}[Direct image of intersection\HTag]\label{thm:direct image of intersection}
  Given normal subobjects $K$ and $L$ of $X$, let $q\from  X \to Y\DefEq X/(K\meet L)$ be the quotient map. Then $q(K)\meet q(L)\leq Y$ is an abelian object.\NoProof
\end{corollary}

\begin{corollary}[Subdiagonal of subobject of an abelian object\HTag]
  \label{thm:SubDiagonal-AbelianObject}
  If $f\from C\to A$ represents $C$ as a subobject of an abelian object $(A,\mu,i)$, then the sequence below is short exact.
  \begin{equation*}
    \xymatrix@R=5ex@C=5em{
    C \ar@{{ |>}->}[r]^-{\PrdctMapInto{\IdMapOn{C},f}} &
    \Prdct{C}{A} \ar@{-{ >>}}[r]^-{\mu\circ ((i\circ f)\times \IdMapOn{A})} &
    A
    }
  \end{equation*}
\end{corollary}
\begin{proof}
  With the information  available, we construct this commutative diagram:
  \begin{equation*}
    \xymatrix@R=5ex@C=5em{
    0 \ar[r]  \ar[d] &
    A \ar@{=}[r] \ar@{{ |>}->}[d]_{\PrdctMapInto{0,\IdMapOn{A}}} &
    A \ar@{=}[d] \\
    C \ar[r]^-{\PrdctMapInto{\IdMapOn{C},f}} \ar@{=}[d] &
    \Prdct{C}{A} \ar@{-{ >>}}[r]^-{\mu\Comp ((i\Comp f)\times \IdMapOn{A})} \ar@{-{ >>}}[d]_{\PrjctnOnto{C}}&
    A \ar[d] \\
    C \ar@{=}[r] &
    C \ar[r] &
    0
    }
  \end{equation*}
  The columns, as well as the top and bottom rows are short exact. The composite of the maps in the middle row vanishes. With the $(3\times 3)$-lemma (\ref{thm:(3x3)-LemmaSemiAb}) we conclude that the middle row is short exact as well.
\end{proof}

As a companion to Theorem \ref{thm:AbelianGroupObject-Recognize-PSA} we have the Proposition \ref{thm:AbGroupViaNormalDiagonal} below.

\begin{proposition}[Recognizing abelian object via normal diagonal\HTag]
  \label{thm:MagmaViaNormalDiagonal} 
  \label{thm:AbGroupViaNormalDiagonal}
  An object $X$ admits an internal abelian group structure $(X,\mu,i)$ if and only if the diagonal $\DgnlOn{X}$ is a normal monomorphism. If so, then the sequence below is short exact.
  \begin{equation*}
    \xymatrix@R=5ex@C=2em{
    X \ar@{{ |>}->}[rr]^-{\Delta_X} &&
    X\times X \ar@{-{ >>}}[rr]^-{\mu(i\times 1_X)} &&
    X
    }
  \end{equation*}
\end{proposition}
\begin{proof}
  If $X$ admits an internal abelian group structure then the stated sequence is short exact by Corollary (\ref{thm:SubDiagonal-AbelianObject}). Conversely, if the monomorphism $\DgnlOn{X}$ is normal, then we construct a multiplication on $X$ from a suitable choice of $D\from \Prdct{X}{X}\to Q$ representing $\CoKerMap{\DgnlOn{X}}$. This construction proceeds through the following steps.
  \begin{enumerate}
    \item First we show that we may choose $Q=X$ in a manner so that $D$ may be interpreted as `addition of first coordinate minus second coordinate'.
    \item From $D$ we extract a candidate $i\from X\to X$ for an internal inverse operation on $X$.
    \item Then $\mu \DefEq D(\Prdct{\IdMapOn{X}}{i})$ is the addition on $X$.
  \end{enumerate}
  \emph{Step 1}\quad With and arbitrary representation $q\from \Prdct{X}{X}\to Q$ of $\CoKerMap{\DgnlOn{X}}$, we obtain the short exact sequence in the commutative diagram below.
  \begin{equation*}
    \xymatrix@R=5ex@C=4em{
    X \ar@{{ |>}->}[r]^-{\DgnlOn{X}} \ar@{=}[d] &
    \Prdct{X}{X} \ar@{-{ >>}}[r]^-{q} \ar[d]_{\PrjctnOnto{2}} \PullLU{rd} &
    Q \ar[d] \\
    X \ar@{=}[r] &
    X \ar[r] &
    0
    }
  \end{equation*}
  The right hand square is a pullback, implying that $q\Comp (\IdMapOn{X},0)\from X\to Q=\Ker{Q\to 0}$ is an isomorphism. Therefore $D\DefEq [q(\IdMapOn{X},0)]^{-1}q$ also represents a cokernel for $\DgnlOn{X}$ with the property that $D(\IdMapOn{X},0)=\IdMapOn{X}$.

  \emph{Step 2}\quad Since $D(\IdMapOn{X},0)=\IdMapOn{X}$ while $D\DgnlOn{X}=0$, we take $i\DefEq D(0,\IdMapOn{X})$ as a candidate for the inverse operation $i$. It is seen to be an isomorphism on $X$ via the following diagram:
  \begin{equation*}
    \xymatrix@R=5ex@C=4em{
    X \ar@{{ |>}->}[r]^-{\DgnlOn{X}} \ar@{=}[d] &
    \Prdct{X}{X} \ar@{-{ >>}}[r]^-{D} \ar[d]_{\PrjctnOnto{1}} \PullLU{rd} &
    X \ar[d] \\
    X \ar@{=}[r] &
    X \ar[r] &
    0
    }
  \end{equation*}
  Again, the right hand square is a pullback, implying that $i\DefEq D\Comp (0,\IdMapOn{X})$ is the induced isomorphism $\Ker{\PrjctnOnto{1}}\to \Ker{X\to 0}$. To show that the proposed Step 3 is applicable, we need to confirm that $i^2 = \IdMapOn{X}$. To this end, we check that the diagram below commutes by verifying that this is the case for each of the maps in the jointly extremal collection $(\IdMapOn{X},0,0)$, $(0,\IdMapOn{X},0)$, $(0,0,\IdMapOn{X})$.
  \begin{equation*}
    \xymatrix@C=6em{
    \PPrdct{X}{X}{X} \ar[d]_{(\PrjctnOnto{1},\PrjctnOnto{3})} \ar[r]^-{\IdMapOn{X}\times \DgnlOn{X}\times \IdMapOn{X}} &
    X\prdct X\prdct X\prdct X \ar[r]^-{\IdMapOn{X}\times \IdMapOn{X}\times (\PrjctnOnto{2},\PrjctnOnto{1})} & \PPrdct{X}{X}{X}\prdct X \ar[d]^{D\times D}\\
    \Prdct{X}{X} \ar[r]_-{D} &
    X &
    \Prdct{X}{X} \ar[l]^-{D}
    }
  \end{equation*}
  Now, we find:
  \begin{equation*}
    \begin{aligned}
      i^2 & = [D\Comp (0,\IdMapOn{X})]\Comp [D\Comp (0,\IdMapOn{X})]                                                                                                                                 \\
          & = D\Comp (0,D\Comp (0,\IdMapOn{X}) )                                                                                                                                                     \\
          & = D\Comp (D\Comp (\IdMapOn{X},\IdMapOn{X})\ ,\ D\Comp (0,\IdMapOn{X}) )                                                                                                                  \\
          & = D \Comp (\Prdct{D}{D})\Comp (\PPrdct{\IdMapOn{X}}{\IdMapOn{X}}{(\PrjctnOnto{2},\PrjctnOnto{1})})\Comp (\PPrdct{\IdMapOn{X}}{\DgnlOn{X}}{\IdMapOn{X}})\Comp (\IdMapOn{X},\IdMapOn{X},0) \\
          & = D\Comp (\PrjctnOnto{1},\PrjctnOnto{3})\Comp (\IdMapOn{X},\IdMapOn{X},0)                                                                                                                \\
          & = D\Comp (\IdMapOn{X},0) = \IdMapOn{X}
    \end{aligned}
  \end{equation*}
  \emph{Step 3}\quad With $\mu\DefEq D\Comp (\Prdct{\IdMapOn{X}}{i})$ we now find $\mu\Comp (\IdMapOn{X},0) = D\Comp (\IdMapOn{X},0) = \IdMapOn{X}$ by design, and
  \begin{equation*}
    \begin{aligned}
      \mu\Comp (0,\IdMapOn{X}) & = D\Comp (\IdMapOn{X},i)\Comp (0,\IdMapOn{X})             \\
                               & = D\Comp (0, i)                                           \\
                               & = D\Comp (0 , D\Comp (0,\IdMapOn{X}))                     \\
                               & = [D\Comp (0,\IdMapOn{X})] \Comp [D\Comp (0,\IdMapOn{X})] \\
                               & = i^2 = \IdMapOn{X}
    \end{aligned}
  \end{equation*}
  Thus $\mu$ is an internal magma structure on $X$ and, hence, determines uniquely the structure of an internal abelian group object on $X$.
\end{proof}

If $(A,\mu,i)$ is an internal abelian group object in a homological category then we use Proposition \ref{thm:AbGroupViaNormalDiagonal} to identify the kernel of $\mu$:

\begin{definition}[Antidiagonal in abelian group object\HTag]
  \label{def:AbGroupObject-Antidiagonal}
  If $(A,\mu,i)$ is an abelian group object, then the \Defn{antidiagonal} in $\Prdct{A}{A}$ is the subobject represented by $\PrdctMapInto{i,\IdMapOn{A}}\from A\to \Prdct{A}{A}$. %
  \index{antidiagonal}\index{abelian group object!antidiagonal}%
\end{definition}

\begin{proposition}[Kernel of multiplication in an abelian group object\HTag]
  \label{thm:AbGrpObject-Ker(Multiplication)}
  If $(A,\mu,i)$ is an abelian group object, then its antidiagonal is the kernel of $\mu$.
\end{proposition}
\begin{proof}
  From (\ref{thm:AbGroupViaNormalDiagonal}) we know that $\DgnlOn{A}=\Ker{\mu\Comp (\Prdct{i}{\IdMapOn{A}})}$. In the composite
  \begin{equation*}
    \left[ \mu \Comp (\Prdct{i}{\IdMapOn{A}})\right]\Comp \DgnlOn{A} =  \mu \Comp \left[ (\Prdct{i}{\IdMapOn{A}})\Comp \DgnlOn{A}\right] = \mu\Comp \PrdctMapInto{i,\IdMapOn{A}}
  \end{equation*}
  the map $(\Prdct{i}{\IdMapOn{A}})$ is an isomorphism, and this implies the claim.
\end{proof}

\begin{subordinate}{}
  \begin{subsubordinate}{On recognizing abelian objects via a normal diagonal map}
    Proposition \ref{thm:AbGroupViaNormalDiagonal}	is a variation on Proposition~3.2.15 in~\cite{FBorceuxDBourn2004}, where the concept of a normal subobject is strictly weaker from what we consider here. The two approaches coincide in the context of a semiabelian category.
  \end{subsubordinate}

  \begin{subsubordinate}{On the collection of all abelian objects in a homological, respectively semiabelian, category}
    In Section~\ref{sec:CatsOfAbelianObjects}, we study the global structure of the collection of these abelian objects. In a homological category $\Ctgry{X}$, they form a left almost abelian category which is full and replete in $\Ctgry{X}$; see Theorem \ref{thm:AbGrp-Objects-H}. If $\Ctgry{X}$ is semiabelian, then the abelian objects form an abelian category.
  \end{subsubordinate}
\end{subordinate}
\newpage
\section[Alternate Characterizations of Homological Categories]{Alternate Characterizations of Homological Categories}
\label{sec:HomologicalCats-AlternateAxioms}

We offer several alternate characterizations of what separates normal categories from homological ones. For example (\ref{thm:N-implies-SFL-then-H}), in the context of a normal category, turning the conclusion of Short $5$-Lemma into a structural axiom, implies the validity of the \KSGInline-condition. We emphasize that this is far from true in an arbitrary \ZExact\ category which does not satisfy the \AENInline-condition. This is explained in detail in Section~\ref{sec:Protomodular-SEpi(X)->X}.

Building on (\ref{thm:N-implies-SFL-then-H}), we present in (\ref{thm:Hofmann})  further criteria which characterize the separation between normal and homological categories. Among these, we single out criterion (IV) which already appeared in the work of Hofmann \cite[Ax. XII]{FHofmann1960-GKat}. For additional information on early works generalizing the framework of abelian categories to non-commutative settings, see \cite{Janelidze-Marki-Tholen}.

\begin{proposition}[Short $5$-Lemma implies the \KSGInline-Axiom\NTag]
  \label{thm:N-implies-SFL-then-H}%
  Any normal category in which the Split Short $5$-Lemma holds is homological.
\end{proposition}
\begin{proof}
  We verify the \KSGInline-property. Let $m\from M\to X$ be a monomorphism through which $k$ and $x$ factor via morphisms $\tilde{k}\from K\to M$ and $\tilde{x}\from Q\to M$. Then we obtain the diagram below in which $k=m\tilde{k}$ and $x=m\tilde{x}$:
  \begin{equation*}
    \xymatrix@R=5ex@C=4em{
    K \ar@{{ |>}->}[r]^-{\tilde{k}} \ar@{=}[d] &
    M \ar@{-{ >>}}@<.5ex>[r]^-{q\Comp m} \ar@{{ >}->}[d]_{m} &
    Q \ar@{=}[d] \ar@<.5ex>[l]^-{\tilde{x}}\\
    K \ar@{{ |>}->}[r]_-{k} &
    X \ar@{-{ >>}}@<+.5ex>[r]^-{q} &
    Q \ar@<.5ex>[l]^-{x}
    }
  \end{equation*}
  Then $\tilde{x}$ is a section for $qm$. Further, $\tilde{k}=\KerMap{qm}$ follows by direct verification. Thus, the diagram above is a morphism of split short exact sequence to which we apply the Split Short $5$-Lemma, and conclude that $m$ is an isomorphism. So, the \KSGInline-condition holds.
\end{proof}

\begin{corollary}[Conditions which are equivalent to \KSGInline\NTag]
  \label{thm:Normal->Homological-Conditions}%
  \label{thm:Hofmann}
  In a normal category $\Ctgry{X}$, the following conditions are equivalent:
  \begin{tfae}
    \item $\Ctgry{X}$ is homological.
    \item The Short $5$-Lemma holds in $\Ctgry{X}$.
    \item The Split Short $5$-Lemma holds in $\Ctgry{X}$.
    \item $\Ctgry{X}$ satisfies \Defn{Hofmann's Axiom} \cite[Ax. XII]{FHofmann1960-GKat}: whenever in a morphism of short exact sequences, as in \eqref{eq:Hofmann} below, the map $\xi$ is a monomorphism and $\rho$ is a normal monomorphism, then $\xi$ is a normal monomorphism.
    \item For a morphism of short exact sequences in $\Ctgry{X}$,
    \begin{equation}\label{eq:Hofmann}
      \vcenter{
      \xymatrix@R=5ex@C=4em{
      K \ar@{{ |>}->}[r]^-{k} \ar@{=}[d] &
      X \ar@{-{ >>}}[r]^-{q} \ar[d]_{\xi} &
      Q \ar@{{ |>}->}[d]^{\rho} \\
      K  \ar@{{ |>}->}[r]_{\xi k} &
      Y  \ar@{-{ >>}}[r]_-{r} &
      R
      }}
    \end{equation}
    whenever $\rho$ is a normal monomorphism, then $\xi$ is a normal monomorphism.
    \item For a morphism of short exact sequences in $\Ctgry{X}$,
    \begin{equation}\label{eq:CoHof}
      \vcenter{
      \xymatrix@R=5ex@C=4em{
      K \ar@{{ |>}->}[r]^-{k} \ar@{-{ >>}}[d]_{\kappa}  &
      X \ar@{-{ >>}}[r]^-{r\xi} \ar[d]^{\xi}&
      Q \ar@{=}[d] \\
      L \ar@{{ |>}->}[r]_-{l} &
      Y \ar@{-{ >>}}[r]_-{r} &
      Q
      }
      }
    \end{equation}
    whenever $\kappa$ is a normal epimorphism, then $\xi$ is a normal epimorphism.
  \end{tfae}
\end{corollary}
\begin{proof}
  The implication (I) $\Rightarrow$ (II) is Theorem \ref{thm:Short5}; (II) $\Rightarrow$ (III) is clear. (III) $\Rightarrow$ (I) is  (\ref{thm:N-implies-SFL-then-H}); (I) $\Rightarrow$ (V) by (\ref{thm:MorphismOfSESs-NormalComponents}.i); (I) $\Rightarrow$ (VI) by (\ref{thm:NormalPushOut-Recognize-H}); (V) $\Rightarrow$ (IV) is clear.

  To see that (IV) $\Rightarrow$ (II), consider the commutative diagram below which is constructed from the morphism of short exact sequences in the front. Assuming that $\kappa$ and $\rho$ are isomorphisms, we need to show that $\xi$ is an isomorphism.
  
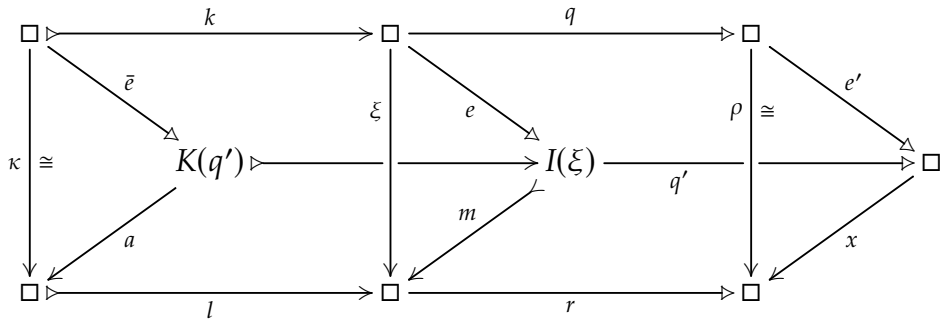
\begin{equation*}
    \xymatrix@!@R=3ex@C=3em{
    \DiagObj \ar@{{ |>}->}[rr]^-{k} \ar[dd]_{\kappa}^{\cong} \ar@{-{ >>}}[rd]^{\bar{e}} &&
    \DiagObj \ar@{-{ >>}}[rr]^-{q} \ar[dd]_(.3){\xi} \ar@{-{ >>}}[rd]_{e} &&
    \DiagObj \ar[dd]_(.3){\rho}^(.3){\cong} \ar@{-{ >>}}[rd]^{e'}\\
    & \Ker{q'} \ar@{{ |>}->}[rr]|\hole \ar[dl]^{a} &&
    \Img{\xi} \ar@{-{ >>}}[rr]|\hole_(0.3){q'} \ar@{{ >}->}[dl]_{m} &&
    \DiagObj \ar[dl]^{x} \\
    \DiagObj \ar@{{ |>}->}[rr]_-{l} &&
    \DiagObj \ar@{-{ >>}}[rr]_{r} &&
    \DiagObj
    }
  \end{equation*}
  In the above diagram, the composite $\xi=me$ is the image factorization of $\xi$. The upward square on the right is the pushout of $e$ along $q$. So, $q'$ and $e'$ are normal epimorphisms. The pushout property of this square yields the map $x$, uniquely rendering the right hand side of the diagram commutative. But then $e'$ is a monomorphism as well, hence is an isomorphism.

  The universal property of the kernel $\Ker{q'}$ yields $\bar{e}$ uniquely rendering the upward square on the left commutative. This face is a pullback since $e'$ is monic. So $\bar{e}$ is a normal epimorphism via the \PNEInline-property. The universal property of the kernel $l$ yields $a$ uniquely rendering the downward square on the left commutative. The monic property of $l$ shows that the triangle on the left commutes as well. But then the normal epimorphism $\bar{e}$ is also monic, hence is an isomorphism. The Primordial Short $5$-Lemma (\ref{thm:Short-5-Primordial})  shows that $e$ is an isomorphism.

  Thus $a$ and $x$ are isomorphisms as well. By hypothesis (IV), $m$ is a normal monomorphism. Another application of the Primordial Short $5$-Lemma shows that $m$ is an isomorphism. So, $\xi$ is an isomorphism, as was to be shown.

  To see that (VI) $\Rightarrow$ (II), consider the situation where $\kappa$ in \eqref{eq:CoHof} is an isomorphism. Then $\xi$ is a normal epimorphism by (VI), hence an isomorphism by the Primordial Short $5$-Lemma \ref{thm:Short-5-Primordial}, which holds in homologically self-dual categories and, hence, in normal categories.
\end{proof}

Note how (V) and (VI) supplement the \HSDInline-conditions in Section~\ref{sec:HomologicalSelfDuality}.

\begin{subordinate}{On Hofmann's approach to generalizing abelian categories}
  In \cite{FHofmann1960-GKat}, F.~Hofmann presents an axiomatic framework which, while inspired by that of abelian categories, includes non-abelian categories such as $\Grps$ as well. He works with (normal) monomorphisms and (normal) epimorphisms under the assumption that every morphism has an essentially unique epi/mono factorization; his axioms (V) and (VI). Later he also assumes that every epimorphism is normal; Ax (XI). Combined with his Ax (XII), this implies that the \ANNInline-condition holds. So, as we would predict, he is then able to prove that the border cases of the $(\Prdct{3}{3})$-Lemma hold. - The axiom we single out in (\ref{thm:Normal->Homological-Conditions}) is his Ax (XII.2).
\end{subordinate}
\newpage
\section[Protomodularity]{Protomodularity and the Fibration $\SEpisIn{X}\to \Ctgry{X}$}
\label{sec:Protomodular-SEpi(X)->X}

In this section, we explain the relationship between homological categories, as presented here, and Bourn's pioneering work~\cite{DBourn1991} which led to the crucial concept of \emph{protomodularity}. Given a category $\Ctgry{X}$, Bourn analyzed the associated category whose objects  are sectioned epimorphisms\footnote{As we know from Example \ref{exa:NxN-add->N}, such a sectioned epimorphism need not a normal map.}; see (\ref{subsec:AbsoluteEpis}) for background.

The category  $\SEpisIn{X}$  is a functor category. To see this,  let $\RtrctCat$ denote the category whose graph is \dots
\begin{equation*}
  \RtrctCat \DefEq\quad     \xymatrix@R=5ex@C=5em{
  0 \ar@<+0.5ex>[r]^-{\pi} &
  1 \ar@<+0.5ex>[l]^-{\sigma}
  }
  \qquad\qquad \pi\sigma=\IdMapOn{1}
\end{equation*}
Then $\SEpisIn{X}=\Ctgry{X}^{\RtrctCat}$. Consequently:

\begin{proposition}[(Co)Limits in $\SEpisIn{X}$]
  \label{thm:SEpi(X)-(Co)Limits}%
  Given a category $\Ctgry{X}$, assume that $\Ctgry{X}$ admits (co)limits over a small category $\SmallCtgry{D}$. Then $\SEpisIn{X}$ admits (co)limits over $D$, and these are computed object-wise.\NoProof
\end{proposition}

\begin{corollary}[Inheritance of properties of $\SEpisIn{X}$ from $\Ctgry{X}$]
  \label{thm:SEpi(X)InheritsFrom-X}%
  If $\Ctgry{X}$ is \ZExact, linear, di-exact, normal, homological, semiabelian, or abelian, then so is $\SEpisIn{X}$.\NoProof
\end{corollary}

The inclusion of $\Ord{0}\to \RtrctCat$ which sends $0$ to $1$ induces the forgetful functor ${\pi\from \SEpisIn{X}\to \Ctgry{X}}$. It sends a sectioned epimorphism $\SctndEpi{p}{x}$ where $p\from X\to Y$ to its codomain $Y$. Proposition~\ref{thm:PullbackOfSplitEpi} and \ref{thm:PushOutOfSplitMono} together show that the functor $\pi$ is a Grothendieck bifibration\footnote{For a finitely complete category $\Ctgry{X}$, the functor $\pi\from \SEpisIn{X}\to \Ctgry{X}$ is also known as the \Defn{fibration of points.}} as soon as pullbacks of split epimorphisms and pushouts of split monomorphisms exist in~$\Ctgry{X}$. For a given object $R$ in $\Ctgry{X}$, we consider the fiber of $\pi$ over~$R$:

\begin{definition}[Sectioned category over a fixed object]%
  \label{def:SEpi_R(X)}%
  For a fixed object $R$ in a category $\Ctgry{X}$, the category of sectioned epimorphisms over $R$ is $\SEpisInOver{X}{R}$ has %
  \index{category!of sectioned epis over $R$}%
  \index[not]{s!$\SEpisInOver{X}{R}$\IndSep category of sectioned epis in $\Ctgry{X}$ over $R$}%
  \begin{enumerate}[(i)]
    \item \emph{objects}\quad all objects of $\SEpisIn{X}$  of the form  $(p\from E\to R$ with section $x\from R\to E)$;
    \item \emph{morphisms}\quad all morphisms of sectioned epimorphisms of the form
          \begin{equation}\label{diag:MapIn-Pt_R(X)}%
            \vcenter{\xymatrix@R=5ex@C=5em{
            X \ar[r]^-{h} \ar@<-.5ex>[d]_{q} &
            Y \ar@<-.5ex>[d]_{r} \\
            R \ar@<-.5ex>[u]_{x} \ar@{=}[r] &
            R \ar@<-.5ex>[u]_{y}
            } }
          \end{equation}
          \index{morphism!in $\SEpisInOver{X}{R}$}%
  \end{enumerate}
\end{definition}

\begin{proposition}[$\SEpisInOver{X}{R}$ is pointed]
  \label{thm:SEpi(X)Pointed}%
  For any object $R$ in any category $\Ctgry{X}$, the category $\SEpisInOver{X}{R}$ has a zero object given by $\SctndEpi{\IdMapOn{R}}{\IdMapOn{R}}$. \NoProof
\end{proposition}

If the category $\Ctgry{X}$ be finitely bicomplete, then a morphism $f\from Q\to R$, yields the \Defn{change-of-base functors} $f^*\from \SEpisInOver{X}{R}\to \SEpisInOver{X}{Q}$ and $f_*\from \SEpisInOver{X}{R}\to \SEpisInOver{X}{Q}$ defined, respectively, by pulling back and pushing out, as in (\ref{thm:PullbackOfSplitEpi}) and (\ref{thm:PushOutOfSplitMono}).

\begin{definition}[Conservative functor]
  \label{def:ConservativeFunctor}%
  A functor $F\from \Ctgry{X}\to \Ctgry{Y}$ is \Defn{conservative} if it \Defn{reflects isomorphisms}.
  \index{functor!conservative}\index{conservative!functor}%
\end{definition}

Thus $F$ is conservative if, whenever it sends a morphism $h$ in $\Ctgry{X}$ to an isomorphism $F(h)$ in $\Ctgry{Y}$, then $h$ is an isomorphism. Here is Bourn's view of protomodularity, cf.~(\ref{thm:PullBackReflectsIsos}):

\begin{definition}[Protomodularity]
  \label{def:Protomodularity}%
  \cite{DBourn1991,DBourn2000}\quad A finitely complete category $\Ctgry{X}$ is \Defn{(Bourn) protomodular} if, for every morphism $f\from Q\to R$ in $\Ctgry{X}$, its change-of-base functor $f^{\ast}\from \SEpisInOver{X}{R} \to \SEpisInOver{X}{Q}$ is conservative.
\end{definition}

This condition is particularly interesting when $\Ctgry{X}$ is pointed, because the change-of-base functor along $\ZeroMap\from \ZeroObject\to R$
\begin{equation*}
  \ZeroMap^{\ast}\from \SEpisInOver{X}{R}\longrightarrow \SEpisInOver{X}{\ZeroObject} \cong \Ctgry{X}
\end{equation*}
returns the kernel of a given sectioned epimorphism as in this morphism in $\SEpisIn{X}$:
\begin{equation*}
  \xymatrix@R=5ex@C=4em{
  \Ker{q} \ar@{{ |>}->}[r] \ar@<-0.5ex>[d]_{\ZeroMap=\bar{q}} &
  X \ar@<-0.5ex>[d]_{q} \\
  \ZeroObject \ar[r]_-{\ZeroMap} \ar@<-0.5ex>[u]_{\bar{x}=\ZeroMap}&
  R \ar@<-0.5ex>[u]_{x}
  }
\end{equation*}
Thus $\ZeroMap^{\ast}$ takes a sectioned epimorphism $\SctndEpi{q}{r}$ and sends it to $\Ker{q}$, viewed as a sectioned epimorphism over $\ZeroObject$. We may thus think of it as the \Defn{kernel functor} and denote it $K_R\from \SEpisInOver{X}{R} \to \Ctgry{X}$. Accordingly, a morphism $(h,\IdMapOn{R})\from \SctndEpi{q}{x}\to \SctndEpi{r}{y}$ yields the commutative diagram below.
\begin{equation*}
  \xymatrix@R=5ex@C=4em{
  \Ker{q} \ar[r]^-{\KerMap{q}} \ar[d]_{K_{R}(h)} &
  X \ar[d]^-{h} \ar@<0.5ex>[r]^-{q} &
  R \ar@<0.5ex>[l]^-{x} \ar@{=}[d] \\
  \Ker{r} \ar[r] \ar[r]_-{\KerMap{r}} &
  Y \ar@<0.5ex>[r]^-{r} &
  R \ar@<0.5ex>[l]^-{y}
  }
\end{equation*}

\begin{proposition}[Protomodularity for pointed categories]\label{thm:ProtomodularityPointed}
  A finitely complete pointed category $\Ctgry{X}$ is protomodular if and only if, for each object $R$ of~$\Ctgry{X}$, the kernel functor ${K_R\from \SEpisInOver{X}{R} \to \Ctgry{X}}$ reflects isomorphisms.
\end{proposition}
\begin{proof}
  To see why $\Leftarrow$ is true, consider change-of-base along $f\from Q\to R$, as in the diagram below.
  
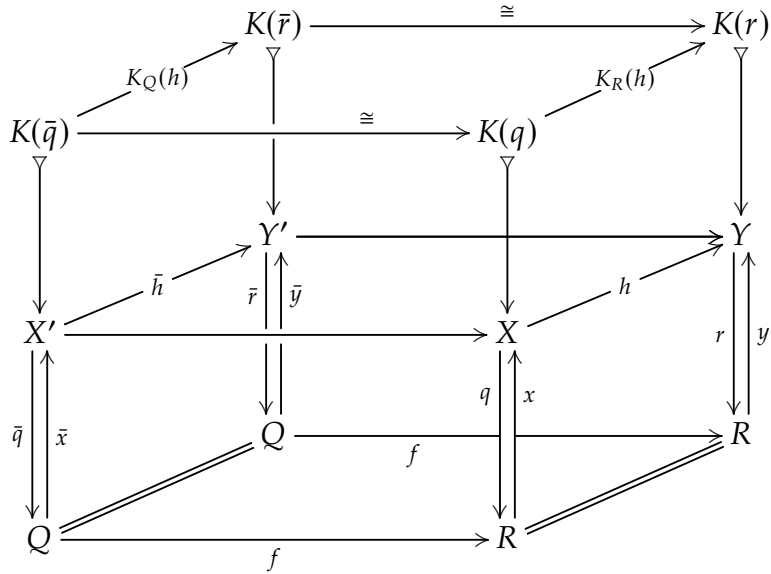
\begin{equation*}
    \xymatrix@R=4ex@C=5em{
    & \Ker{\bar{r}} \ar[rr]^-{\cong} \ar@{{ |>}->}[dd]|\hole &&
    \Ker{r} \ar@{{ |>}->}[dd] \\
    \Ker{\bar{q}} \ar@{{ |>}->}[dd] \ar[rr]^(0.7){\cong} \ar[ru]|-{\ K_{Q}(h)\ } &&
    \Ker{q} \ar@{{ |>}->}[dd] \ar[ru]|-{\ K_{R}(h)\ }\\
    & Y' \ar[rr] \ar@<-0.5ex>[dd]|\hole_(0.3){\bar{r}} \ar[rr] &&
    Y \ar@<-0.5ex>[dd]_{r} \\
    X' \ar[rr] \ar@<-0.5ex>[dd]_{\bar{q}} \ar[ru]|-{\ \bar{h}\ } &&
    X \ar@<-0.5ex>[dd]_(0.3){q} \ar[ru]|-{\ h\ } \\
    & Q \ar[rr]|\hole_(0.3){f} \ar@<-0.5ex>[uu]|\hole_(0.7){\bar{y}} &&
    R \ar@<-0.5ex>[uu]_{y} \\
    Q \ar[rr]_-{f}\ar@{=}[ru] \ar@<-0.5ex>[uu]_{\bar{x}} &&
    R \ar@{=}[ru] \ar@<-0.5ex>[uu]_(0.7){x}
    }
  \end{equation*}
  The front and back faces of the bottom block are constructed as pushouts. Therefore, the top two horizontal maps are isomorphisms; which yields the identity $K_{Q}\Comp f^{\ast} \cong K_{R}$.

  If $\bar{h}$ is an isomorphism, then we must show that $h$ is an isomorphism. But if $\bar{h}$ is an isomorphism, then so is $K_{Q}(\bar{h})$ and, hence, $K_{R}(h)$. Since $K_{R}$ is an isomorphism, so is $h$. This means that $f^{\ast}$ is conservative, as was to be shown.
\end{proof}

In any homological category, the \KSGInline-condition is satisfied by definition, while (\ref{thm:SplitShort5}) implies that $K_{R}$ reflects isomorphisms for every object $R$. So by (\ref{thm:ProtomodularityPointed}), any homological category is protomodular. Actually, the equivalence between protomodularity and \KSGInline\ is true in a far more general context:

\begin{corollary}[Protomodularity vs.\ the \KSGInline-condition]
  \label{thm:Protomodular<->KSG}%
  A finitely complete pointed category $\Ctgry{X}$ is Bourn protomodular if and only if the
  \KSGInline-condition holds.
\end{corollary}
\begin{proof}
  First suppose $\Ctgry{X}$ is protomodular. To see that the \KSGInline-condition is satisfied, consider a sectioned epimorphism $\SctndEpi{q}{x}$, whose section and kernel lift through a monomorphism $m$, as in this diagram:
  \begin{equation*}
    \xymatrix@R=5ex@C=4em{
    K \ar@{{ |>}->}[r]^-{\kappa} \ar@{=}[d] &
    M \ar@<.5ex>@{-{>}}[r]^-{q\Comp m} \ar@{{ >}->}[d]^{m} &
    Q \ar@{=}[d] \ar@<.5ex>[l]^-{\sigma} \\
    K \ar@{{ |>}->}[r]_-{k} &
    X \ar@<.5ex>@{-{>}}[r]^-{q} &
    Q \ar@<.5ex>[l]^-{x}
    }
  \end{equation*}
  As in the proof of (\ref{thm:N-implies-SFL-then-H}), we see that $k=\KerMap{qm}$. By protomodularity, $m$ is an isomorphism. So $k$ and $x$ are extremally epimorphic; i.e., the \KSGInline-condition holds.

  Now assume that \KSGInline\ holds and consider the diagram
  \begin{equation*}
    \xymatrix@R=5ex@C=4em{
    K \ar@{{ |>}->}[r]^-{k} \ar[d]^-{\overline{h}}_-{\cong} &
    X \ar[d]^-{h} \ar@<0.5ex>[r]^-{q} &
    R \ar@<0.5ex>[l]^-{x} \ar@{=}[d] \\
    L \ar@{{ |>}->}[r]_-{l} &
    Y \ar@<0.5ex>[r]^-{r} &
    R \ar@<0.5ex>[l]^-{y}
    }
  \end{equation*}
  where $k$, $l$ are the respective kernels of $q$ and $r$ and the restriction of $h$ to a map $K\to L$ is an isomorphism. Taking kernel pairs vertically, by (\ref{thm:KernelPair-Monos}) we find a following split epimorphism of effective equivalence relations, together with its kernel.
  \begin{equation*}
    \xymatrix@R=5ex@C=4em{
    K \ar[r]^-{\overline{k}} \ar@{=}@<-1ex>[d] \ar@{=}@<1ex>[d] &
    \KrnlPr{h} \ar@<-1ex>[d] \ar@<1ex>[d] \ar@<0.5ex>[r]^-{\overline{q}} &
    R \ar@<0.5ex>[l]^-{\overline{x}} \ar@{=}@<-1ex>[d] \ar@{=}@<1ex>[d]   \\
    K \ar[r] \ar[r]_-{k}  \ar@{=}[u] &
    X \ar@<0.5ex>[r]^-{q} \ar[u]|(.45){\hole}|-{m}|(.55){\hole} &
    R \ar@<0.5ex>[l]^-{x} \ar@{=}[u]
    }
  \end{equation*}
  By the \KSGInline-property, the upward pointing arrow $m$ is a monomorphism. Again from (\ref{thm:KernelPair-Monos}), we deduce that $h$ is a monomorphism. Now we may use the \KSGInline-property once more to see that $h$ is an isomorphism.
\end{proof}

We may ask if a pointed finitely complete category in which the conclusion of the Split Short $5$-Lemma holds satisfies the \KSGInline-condition. The answer is `no', as is demonstrated by the following example.

\begin{example}[Split Short 5-Lemma but not \KSGInline]
  \label{exa:Hopf-Algs-SplitShort-5}
  Theorem 4.1 in \cite{Molnar} implies that a version of the Split Short 5-Lemma (where the split epimorphisms in the diagram are supposed to be normal epis) holds in the category consisting of all Hopf algebras over a given field. However, by the results of~\cite{GM-VdL1}, this category does not satisfy the \KSGInline-condition.
\end{example}

Conclusion: The protomodularity condition is subtly stronger than asking that the Split Short $5$-Lemma holds. On the other hand, recall (\ref{thm:N-implies-SFL-then-H}) which says, in a normal category, the Split Short $5$-Lemma is equivalent to the \KSGInline-condition.

\subsection{Borceux--Bourn homological categories}\label{sec:BBHom}
In the notes to Section~\ref{sec:HomologicalCats-Axioms} we called a category is \Defn{Borceux--Bourn homological} when it is pointed, finitely complete, admits coequalizers of kernel pairs, preserves normal epimorphisms under base change, and has the \KSGInline-property. We explain now that this agrees with the original definition of Borceux and Bourn: a pointed, regular (\ref{def:RegularCategory}) and protomodular category.

\begin{proposition}[BB-homological categories]
  \label{thm:BBHom}
  A category is Borceux--Bourn homological if and only if it is pointed, regular and protomodular.
\end{proposition}
\begin{proof}
  Corollary~\ref{thm:Protomodular<->KSG} takes care of the equivalence between protomodularity and the \KSGInline-condition. Then, in both situations---pointed plus \KSGInline\ plus the given (co)limits on the one hand, pointed plus protomodular plus the given (co)limits on the other---normal and regular epimorphisms coincide, so that regularity is equivalent to pullback-stability of normal epimorphisms.
\end{proof}

\begin{subordinate}{}
  \begin{subsubordinate}{On the classifying properties of the fibration of points}
    A red thread in the foundational work of D.~Bourn is the interpretation of categorical-algebraic properties in terms of the fibration of points. We saw that protomodularity amounts to the condition that the change-of-base functors are conservative. Another example is the fact that a pointed finitely bicomplete category is abelian if and only if the change-of-base functors are equivalences (\ref{exe:AbelianIffEquiv}). See \cite{BB} and the references there for a proof of this fact and for further examples.
  \end{subsubordinate}

  \begin{subsubordinate}{Relationship to internal actions}
    The category of points over an object $R$ may be interpreted as a category of non-abelian modules over $R$. We will develop this point of view within the context of \emph{internal $R$-actions}.
  \end{subsubordinate}
\end{subordinate}

\begin{exercises}
\begin{exercise}[Isomorphisms in $\SEpisInOver{X}{R}$]
  \label{exe:IsosIn-SEpi_R(X)}%
  Given an arbitrary category $\Ctgry{X}$, identify the isomorphisms in $\SEpisInOver{X}{R}$.
\end{exercise}

\begin{exercise}[(Co)limits in $\SEpisInOver{X}{R}$]
  \label{exe:(Co)LimitsIn-SEpi_R(X)}%
  If the category $\Ctgry{X}$ admits (co)limits modeled on a small category $\SmallCtgry{D}$, show that so does $\SEpisInOver{X}{R}$.
\end{exercise}

\begin{exercise}[Abelianness and the functors $K_R$]\label{exe:AbelianIffEquiv}
  Show that a pointed finitely bicomplete category is abelian if and only if each of the change-of-base functors $K_R$ is an  equivalence.
\end{exercise}
\end{exercises}
\section[The Mal'tsev Property]{The Mal'tsev Property}
\label{sec:Maltsev}

We explain how being a homological category relates to a property of internal relations~(\ref{sec:InternalGraphs/Relations}): each internal reflexive relation is automatically an equivalence relation. This connects to work of Mal'tsev in the context of universal algebra~\cite{Maltsev-Sbornik}.

\begin{definition}[Mal'tsev category]
  \label{def:MaltsevCategory}%
  A finitely complete category is a \Defn{Mal'tsev category} if every reflexive relation is an equivalence relation. %
  \index{Mal'tsev category}\index{category!Mal'tsev}
\end{definition}

\begin{example}[Mal'tsev varieties]
  \label{exa:MaltsevVariety}%
  A variety of algebras is a Mal'tsev category if and only if it is a Mal'tsev variety; 
  that is, there exists a ternary operation $p(x,y,z)$ such that
  \begin{equation*}
    p(x,x,z)=z \qquad \text{and}\qquad p(x,z,z)=x.
  \end{equation*}
  This operation defines a non-commutative affine structure via the interpretation $p(x,y,z)$ gives the `vector sum' of $x$ and $z$ relative to $y$. The defining equations may thus be seen as a version of the parallelogram rule. For example, for a group $G$, setting
  \begin{equation*}
    p(x,y,z) \DefEq xy^{-1}z
  \end{equation*}
  turns $G$ into a non-commutative affine space: any $y\in G$ may be turned into a neutral element by placing it into the second coordinate.
\end{example}

\begin{example}[Monoids]
  \label{exa:MonoidsNotMaltsev}%
  The variety of monoids is not a Mal'tsev variety. Indeed, the order relation $\leq$ may be considered as an internal relation on the monoid $\NNr$ of natural numbers. It is a reflexive relation which is not symmetric.
\end{example}

The following result explains the terminology used in~(\ref{thm:ProtoMaltsev}):

\begin{proposition}[Characterization of Mal'tsev categories]
  \label{thm:CharMaltsev}%
  For a finitely complete category $\Ctgry{X}$, the following conditions are equivalent:
  \begin{tfae}
    \item in a double sectioned epimorphism
    \[
      \vcenter{\xymatrix@R=5ex@C=4em{
      U \ar@<-.5ex>@{-{ >>}}[r]_-g \ar@{-{ >>}}@<-.5ex>[d]_-{u} &
      V \ar@{-{ >>}}@<-.5ex>[d]_-{v} \ar@<-.5ex>@{{ >}->}[l]_-n \\
      X \ar@<-.5ex>@{-{ >>}}[r]_-f \ar@<-.5ex>@{{ >}->}[u]_-s &
      Y \ar@<-.5ex>@{{ >}->}[u]_-t \ar@<-.5ex>@{{ >}->}[l]_-m
      }}
    \]
    which forms a pullback square $fu=vg$, the sections $s$ and $n$ are jointly extremal-epimorphic;
    \item every reflexive relation is transitive;
    \item every reflexive relation is symmetric;
    \item $\SACtgry{X}$ is a Mal'tsev category; see (\ref{def:MaltsevCategory}).
  \end{tfae}
\end{proposition}
\begin{proof}
  We start with the proof that (I) implies (II). Let $(R,d_1,d_2,e)$ be a reflexive relation on an object $X$ and consider the pullback
  \[
    \vcenter{\xymatrix@R=5ex@C=4em{
    R\times_X R \ar@<-.5ex>@{-{ >>}}[r]_-{\pi_2} \ar@{-{ >>}}@<-.5ex>[d]_-{\pi_1} &
    R \ar@{-{ >>}}@<-.5ex>[d]_-{d_1} \ar@<-.5ex>@{{ >}->}[l]_-{e_2} \\
    R \ar@<-.5ex>@{-{ >>}}[r]_-{d_2} \ar@<-.5ex>@{{ >}->}[u]_-{e_1} &
    X. \ar@<-.5ex>@{{ >}->}[u]_-e \ar@<-.5ex>@{{ >}->}[l]_-e
    }}
  \]
  By assumption, $e_1$ and $e_2$ are jointly extremal-epimorphic. Hence in the square
  \begin{equation*}
    \xymatrix@C=4em{
    R+R \ar[d]_-{\SumMapOutOf{e_1,e_2}} \ar[r]^-{\nabla_R} &
    R \ar@{{ >}->}[d]^-{(d_1,d_2)} \\
    R\times_X R \ar[r]_-{(d_1\pi_1,d_2\pi_2)} \ar@{.>}[ru]|-{\ \rho\ } &
    X\times X}
  \end{equation*}
  there is a lifting $\rho$ showing that $R$ is a transitive relation.

  For the other implications, given any relation $(R,d_1,d_2)$ from $X$ to $Y$, we shall write $S$ for the reflexive relation on $R$ defined by $(a,b)S(c,d)$ if and only if $aRd$. In other words, $S$ is the pullback of $(d_1,d_2)\from R\to X\times Y$ along the map $d_1\times d_2\from R\times R\to X\times Y$ as in the diagram
  \begin{equation*}
    \xymatrix@R=5ex@C=4em{
    S \PullLU{rd} \ar@{{ >}.>}[d]_-{(d'_1,d'_2)} \ar@{.>}[r] & R \ar@{{ >}->}[d]^-{(d_1,d_2)}\\
    R\times R \ar[r]_-{d_1\times d_2} & X\times Y
    }
  \end{equation*}
  For the proof that (II) implies (III), let $R$ be a reflexive relation on an object $X$ and $S$ the induced reflexive relation. Assume that $aRb$. Then $(b,b)S(a,b)$ and $(a,b)S(a,a)$. Since $S$ is transitive, we have $(b,b)S(a,a)$, so that $bRa$. It follows that $R$ is symmetric.

  To see that (III) implies (II), let $R$ be a reflexive relation on an object $X$ and $S$ the induced reflexive relation. Assume that $aRb$ and $bRc$. Then $cRb$ by symmetry of $R$. Hence $(a,a)S(c,b)$, so that $(c,b)S(a,a)$ by symmetry of $S$. It follows that $cRa$, so that $aRc$ and $R$ is a transitive relation.

  We see that either (II) or (III) implies (IV). To show that (IV) implies (I), consider a double sectioned epimorphism as in the statement of condition (I). Consider a monomorphism $m\from M\to U$ through which $s$ and $n$ factor as $s'\from {X\to M}$ and $n'\from V\to M$ such that $ms'=s$ and $mn'=n$. We have to prove that $m$ is an isomorphism. This amounts to showing that $aMb$ for all $(a,b)\from A\to X\times V$. Note that $M$ is indeed a relation from $X$ to $V$. It induces a reflexive relation $S$ on $M$ as above. This $S$ is an equivalence relation on $M$ by assumption. Now $(a,tfa)S(mvb,tfa)$ and $(mvb,tfa)S(mvb,b)$ so that $(a,tfa)S(mvb,b)$ which implies $aRb$.
\end{proof}

We know from Exercise \ref{exe:MorphismOfSplitEpis-Properties} that a diagram as in (\ref{thm:CharMaltsev}) is always pushout of sectioned epis and a pullback of sectioned monos. In the following corollary we characterize when it is a pullback of sectioned epimorphisms.

\begin{corollary}[When is a split epi of split epis a pullback?]
  \label{thm:PullbackRecognition-SplitEpiOfSplitEpis}
  In a Mal'tsev category a diagram as in (\ref{thm:CharMaltsev}) is a pullback of sectioned epimorphisms if and only if $u$ and $g$ are jointly monomorphic.
\end{corollary}
\begin{proof}
  If the diagram is a pullback, then $u$ and $g$ are seen as jointly monomorphic. Conversely, they are jointly monomorphic, then the comparison to the pullback of $v$ along $f$. This comparison is an extremal epimorphism and a monomorphism, hence an isomorphism.
\end{proof}

\begin{corollary}[Homological implies Mal'tsev\HTag]
  \label{thm:HomologicalThenMaltsev}
  Any homological category is a Mal'tsev category. %
  \index{homological category!is Mal'tsev category}%
\end{corollary}
\begin{proof}
  This combines Proposition~\ref{thm:CharMaltsev} with Corollary~\ref{thm:ProtoMaltsev}.
\end{proof}

\begin{subordinate}[Comment]{on the origins of the concept of a Mal'tsev category}
  Mal'tsev discovered~\cite{Maltsev-Sbornik} that a variety satisfies what we now call the Mal'tsev identities (\ref{exa:MaltsevVariety}) exactly when all of its reflexive relations are equivalence relations. Since this latter condition can be formulated in purely categorical terms, it has been studied in settings which are far wider than varieties of algebras: see, for instance, \cite{Carboni-Lambek-Pedicchio,Carboni-Kelly-Pedicchio,FBorceuxDBourn2004,EGJVdL} and the references there.
\end{subordinate}
\chapter[Semiabelian Categories]{Semiabelian Categories}
\label{chap:SACats}

We introduce the environment of semiabelian category to combine features of di-exact categories and homological categories: A semiabelian category is a homological category which is also di-exact. Thus, in a semiabelian category all of the diagram lemmas established in Chapter \ref{chap:HomologicalCats} are valid. In addition , useful computational tools become available, such as Theorem \ref{thm:(Co)Ker(ProperMapLESs)}, plus results that will be developed here.

From a logician's perspective, it is of interest that, through the validity of the di-exactness axiom, semi-abelian categories allow self-dual proofs of results which are still valid is a homological context, but do not permit a self-dual proof there. An example of this phenomenon is the Snake Lemma.

\emph{Historical perspective}\quad Let us review the key steps leading to the conception of semiabelian categories: In the 1960's, Gerstenhaber, Huq, Orzech \cite{Huq, Gerstenhaber, Orzech} looked for a generalization of Grothendieck's abelian categories which would also include the categories of groups, Lie-algebras, loops, non-unitary rings, etc.\ They presented several systems of structural axioms. However, at the time, these competing systems of structural axioms appeared dissimilar, and no favorite emerged over the next couple of decades.

Some thirty to forty years later, Janelidze, Márki, and Tholen discovered that these dissimilar seeming structural axiom systems were all equivalent. Their discovery was based on more recent developments in categorical algebra and, in particular, on the pivotal concepts of `Barr exactness'~\cite{Barr}, and `Bourn protomodularity'~\cite{DBourn1991, DBourn2001}.
From these unifying insights emerged the much younger concepts of \emph{semiabelian category}, introduced in~\cite{Janelidze-Marki-Tholen}, and \emph{homological category}, introduced in~\cite{FBorceuxDBourn2004}. %

We now have a fairly fine grained homological view of categorical algebra. A relationship of many of the well investigated algebraic categories to this view is presented in the introduction to the current part of the text, starting from page \pageref{sec:Examples-Overview}.

\bigskip\bigskip\bigskip

\begin{center}
  \textbf{Leitfaden for Chapter \ref{chap:SACats}}
\end{center}

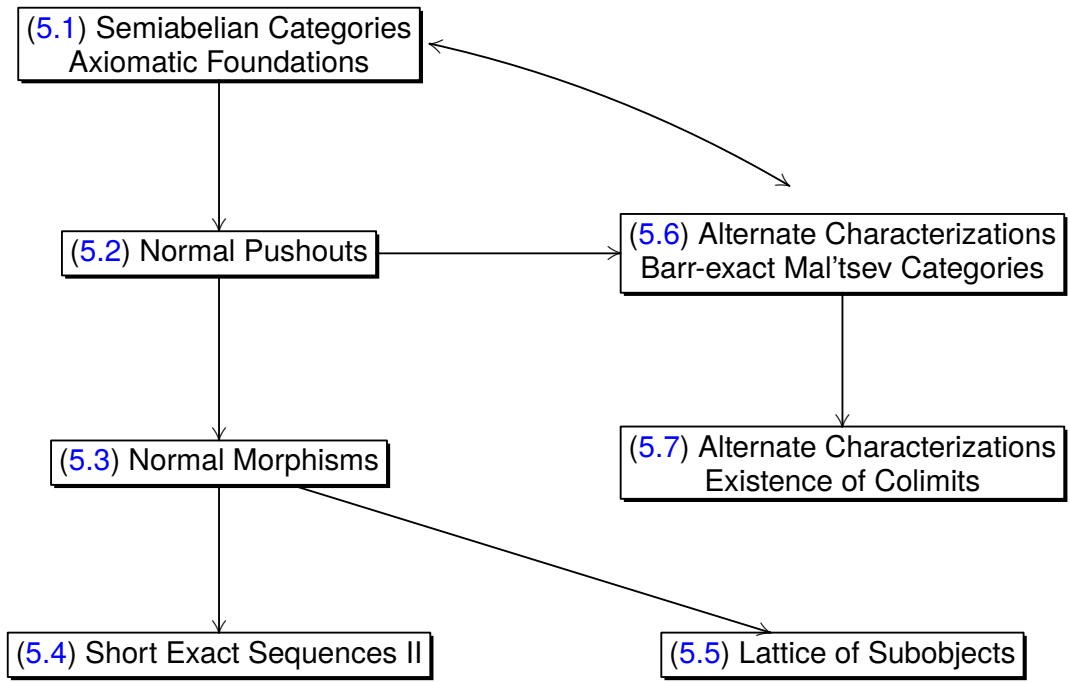
\bigskip

\begin{equation*}
  \xymatrix@R=9ex@C=6em{
  *+[F-,]{\txt{\sffamily (\ref{sec:SACats}) Semiabelian Categories \\ \sffamily Axiomatic Foundations}}\ar[d] \ar@{<->}@/^3ex/@<2ex>[rd] \\
  *+[F-,]{\txt{\sffamily (\ref{sec:NormalPushouts-SACats}) Normal Pushouts}} \ar[d] \ar[r] &
  *+[F-,]{\txt{\sffamily (\ref{sec:ExactMaltsev}) Alternate Characterizations \\ \sffamily Barr-exact Mal'tsev Categories}} \ar[d] \\
  *+[F-,]{\txt{\sffamily (\ref{sec:NormalMaps-SA}) Normal Morphisms}} \ar[d] \ar[rd] &
  *+[F-,]{\txt{\sffamily (\ref{sec:ExistenceColimits}) Alternate Characterizations \\ \sffamily Existence of Colimits}} \\
  *+[F-,]{\txt{\sffamily (\ref{sec:Cat-SESs-SACat}) Short Exact Sequences II}} &
  *+[F-,]{\txt{\sffamily (\ref{sec:LatticeOfSubjects}) Lattice of Subobjects}}
  }
\end{equation*}
\newpage
\section[Axiomatic Foundations]{Semiabelian Categories}
\label{sec:SACats}

With the concept of a semiabelian category, we combine features of di-exact categories and homological categories. For convenience, here is a complete set of structural axioms which characterize a semiabelian category.

\begin{definition}[Semiabelian category]
  \label{def:SACategory}%
  A category $\Ctgry{X}$ is semiabelian if the following structural axioms are satisfied:
  \begin{ulist}
    \item	\emph{Zero object:} $\Ctgry{X}$ has a zero object; see (\ref{def:0-Object}).
    \item \emph{Functorially  finitely complete:} For every functor $F\from J\to \Ctgry{X}$, with $J$ finite, $\LimOf{F}$ exists; see (\ref{sec:Limits-CoLimits}).
    \item \emph{Functorially finitely cocomplete:} For every functor $F\from J\to \Ctgry{X}$, with $J$ finite, $\CoLimOf{F}$ exists; see (\ref{sec:Limits-CoLimits}).
    \item \emph{\PNEInline\ Pullbacks preserve normal epimorphisms:}\quad  The pullback $\bar{g}$ of a normal epimorphism $g$ along any morphism $f$ in $\Ctgry{X}$ is a normal epimorphism.
    \begin{equation*}
      \xymatrix@R=5ex@C=4em{
      P \ar@{-{ >>}}[d]_-{\bar{g}} \ar[r] \PullLU{rd} &
      Z \ar@{-{ >>}}[d]^-{g} \\
      X \ar[r]_-{f} &
      Y
      }
    \end{equation*}
    \item \emph{\KSGInline\ Kernel and section generate:} For a morphism $q\from X\to Q$ with any given section $x\from {Q\to X}$ (so that $qx=\OneMapOn{Q}$), %
    \index[acr]{k!\KSGInline\IndSep for $q\from X\to Q$ sectioned by $x$, images of $x$ and $\KerMap{q}$ generate $X$}%
    \begin{equation*}
      \xymatrix@C=4em{
      K \ar@{{ |>}->}@[blue][r]^-{\color{blue} k} &
      X \ar@<0.5ex>[r]^-{q} &
      Q \ar@[blue]@<0.5ex>[l]^-{\color{blue} x}
      }
    \end{equation*}
    the morphisms $x$ and $k=\KerMap{q}$ are jointly extremal-epimorphic; see (\ref{def:ExtremalEpimorphism}).
    \item \emph{\ANNInline}\quad Every antinormal map is normal.
  \end{ulist}
  \index{homological category}\index{semiabelian category}\index{category!semiabelian}\index{category!homological}%
  \index[acr]{a!\ANNInline\IndSep property that antinormal maps are normal}%
  \index[acr]{s!{\color{MidnightBlue} $\EuRoman{S{\kern-0.15ex}A}$}\IndSep semiabelian category}%
\end{definition}

Thus, among the types of categories we have introduced, the most prominent ones are:
\begin{ulist}
  \item Every finitely bicomplete category with a zero object is \emph{\ZExact}.
  \item \emph{Di-exact categories}\quad are \ZExact, and satisfy the condition \ANNInline.
  \item \emph{Normal categories}\quad are finitely bicomplete, have a zero object, and satisfy the conditions \PNEInline\ and \AENInline---absolute epimorphisms are normal.
  \item \emph{Homological categories}\quad are finitely bicomplete, have a zero object, and satisfy the \PNEInline\ and \KSGInline-conditions.
  \item \emph{Semiabelian categories}\quad are homological and satisfy the \ANNInline-condition. Thus a semiabelian category is at the same time homological and di-exact.
\end{ulist}

We already discussed the nature of these structural axioms individually. We point out that any variety of algebras is bicomplete. So, it is automatically \ZExact\ when it has a zero object. Further, if a pointed variety is homological, then it is automatically semiabelian. This is remarkable, because outside a homological context, a \ZExact\ variety need \emph{not} be di-exact.

\begin{subordinate}{}
  In Definition \ref{def:SACategory}, we characterized a semiabelian category by the requirements that it has a zero object, is finitely bicomplete, pullbacks preserve normal epimorphisms \PNEInline, it satisfies the \KSGInline-condition, and antinormal maps are normal \ANNInline. Sections~\ref{sec:ExactMaltsev} and \ref{sec:ExistenceColimits} lead up to Theorem~\ref{thm:JMTvsUS} which tells us that this collection of structural axioms is equivalent to the ones proposed by Janelidze--Márki--Tholen in \cite{Janelidze-Marki-Tholen}.

  An early version of the concept of a semiabelian category goes back to Huq \cite{Huq} in the late 1960's. Around the same time other competing axiom systems were proposed. All of these developments remained largely dormant until Janelidze--Márki--Tholen in \cite{Janelidze-Marki-Tholen} managed to connect the pivotal concept of Barr exactness with the work of Huq and others. We draw the reader's attention to the introduction of~\cite{Janelidze-Marki-Tholen} for a more detailed historical perspective.

  Next to being valid in any variety of algebras, Barr exactness had long been understood to form part of the structure of an abelian category, and also that of an elementary topos. Rooted in the study of categories of sectioned epimorphisms, called `points' by Bourn, he studied what is now called \emph{Bourn protomodularity}, that is \KSGInline\ and its equivalent versions---see (\ref{sec:HomologicalCats-AlternateAxioms}) and (\ref{sec:Protomodular-SEpi(X)->X}). From the beginning, he recognized the relevance of protomodularity toward the goals of Huq and others and provided further insight into this topic in \cite{DBourn1996,DBourn2001}.

  Full recognition of Bourn's accomplishment only arrived with the works of Bourn--Janelidze \cite{Bourn-Janelidze:Semidirect} and, in particular, \cite{Janelidze-Marki-Tholen} where it is explained how protomodularity acts as the link between the development surrounding Barr exactness and the works of Huq and others. --- In the words of Janelidze--Márki--Tholen, their main result is that `old $=$ new'. %
  \index{`old $=$ new'}

  As a consequence of this discovery, the need to compare seemingly competing notions of semiabelian categories vanished: via `old $=$ new', there is a single notion left, which, as we shall see, is equivalent to the one defined above.

\end{subordinate}

\begin{exercises}
\begin{exercise}[Cokernel and retraction co-generate in $\SetsBsd$]
  \label{exe:(CRG)Set_*}%
  In the category $\SetsBsd$ of pointed sets, show that any monomorphism is a normal map which is also retractable. Then consider a short exact sequence
  \begin{equation*}
    \xymatrix@R=5ex@C=4em{
    K \ar@{{ |>}->}[r]^-{k} &
    X \ar@{-{ >>}}[r]^-{q} &
    Q
    }
  \end{equation*}
  in $\SetsBsd$ and show that $X\cong K+Q$.
\end{exercise}

\begin{exercise}[$\SetsBsd$ is co-semiabelian]
  \label{exe:Set_*CoSemiabelian}%
  We know already that the category $\SetsBsd$ of pointed sets is bicomplete. Show that it is also co-semiabelian and, hence, that $\SetsBsdOp$ is semiabelian: %
  \index{$\SetsBsd$!co-semiabelian}\index{co-semiabelian category!$\SetsBsd$}%
  \begin{thmlist}
    \item \emph{(PNM)}\quad Pushouts preserve normal monomorphisms.
    \item \emph{(CRC) Cokernel and retraction co-generate}\quad Given a retracted monomorphism $\RtrctdMono{m}{r}$, let $q\DefEq \CoKerMap{m}$. Then $r$ and $q$ are jointly extremally monomorphic.
    \item \emph{\ANNInline}\quad Every antinormal map is normal.
  \end{thmlist}
\end{exercise}

\begin{exercise}[$\SetsBsd$ is not p-exact]
  \label{exe:Set_*Not-p-Exact}%
  In $\SetsBsd$ show that there are morphisms which are not normal; i.e., show that $\SetsBsd$ is not a p-exact category. %
  \index{$\SetsBsd$!not p-exact}%
\end{exercise}
\end{exercises}
\section[Normal Pushouts]{Normal Pushouts in a Semiabelian Category}
\label{sec:NormalPushouts-SACats}

In a homologically self-dual category, we know from Section \ref{sec:NormalPushouts/Pullbacks} that every (semi-)normal pushout square actually enjoys the universal property expected of a pushout. A special feature of a semiabelian category is that every pushout of normal epimorphisms is a normal pushout.

\begin{proposition}[Pushout of normal epimorphisms in a semiabelian category\SATag]
  \label{thm:PushoutNormalEpis-SA}%
  In a semiabelian category, every pushout of a normal epimorphism along another is a normal pushout.
\end{proposition}
\begin{proof}
  Via the \ANNInline-condition, we know from (\ref{thm:DPN-DiExtensivePull/Push}) that every pushout of normal epimorphism $q$ along another $\xi$ is di-extensive; that is we obtain the morphism of short exact sequences below.
  \begin{equation*}
    \xymatrix@R=5ex@C=4em{
    \Ker{q} \ar@{-{ >>}}[d]_{\kappa} \ar@{{ |>}->}[r] &
    X \ar@{-{ >>}}[r]^-{q} \ar@{-{ >>}}[d]_-{\xi} \PushRD{rd} &
    Q \ar@{-{ >>}}[d]^{\underline{\xi}} \\
    \Ker{\underline{q}} \ar@{{ |>}->}[r] &
    Y \ar@{-{ >>}}[r]_-{\underline{q}} &
    R
    }
  \end{equation*}
  In particular, the map $\kappa$ is a normal epimorphism.	The square on the right is a normal pushout square by (\ref{thm:NormalPushOut-Recognize-H}).
\end{proof}

Via normal pushouts we obtain the following view of the separation between homological and semiabelian categories:

\begin{proposition}[\ANNInline-condition in terms of normal pushouts\HTag]
  \label{thm:Pushout,AllCokernels->NormalPushout}
  For a homological category $\Ctgry{X}$, the following conditions are equivalent.
  \begin{tfae}
    \item $\Ctgry{X}$ is semiabelian.
    \item $\Ctgry{X}$ is di-exact.
    \item The category $\Ctgry{X}$ has the \ANNInline-property: any antinormal composite is a normal map.
    \item In $\Ctgry{X}$, the pushout of two normal epimorphisms with a common domain is a normal pushout.
  \end{tfae}
\end{proposition}
\begin{proof}
  (I), (II) and (III) are equivalent by definition. (III) $\implies$ (IV) was proved in (\ref{thm:PushoutNormalEpis-SA}). For the proof that (IV) implies (III), given an antinormal map composed of a normal monomorphism $k$ and a normal epimorphism $r$, construct the morphism of short exact sequences below from the pushout of $\CoKerMap{k}$ along $f$.
  \begin{equation*}
    \xymatrix@!0@R=8ex@C=7em{
    K \ar@{{ |>}->}[r]^-{k} \ar@{.>}[d]_-{\kappa} &
    X \ar@{-{ >>}}[r]^-{\CoKerMap{k}} \ar@{-{ >>}}[d]_-{f} \PushRD{rd} &
    Q \ar@{-{ >>}}[d] \\
    L \ar@{{ |>}->}[r]_-{l} &
    Y \ar@{-{ >>}}[r] &
    R
    }
  \end{equation*}
  By hypothesis, the square on the right is a normal pushout. So, $\kappa$ is a normal epimorphism. This means that $l \varphi$ is an image factorization of $fk$. Thus the \ANNInline-condition is satisfied.
\end{proof}

\begin{subordinate}{Relation with Barr-exactness}
  The equivalence in (\ref{thm:Pushout,AllCokernels->NormalPushout}) will allow us in (\ref{sec:ExactMaltsev}) to prove the following: %
  \begin{ulist}
    \item a semiabelian category is always a Barr-exact Mal'tsev category, %
    %
    \item In a finitely bicomplete pointed category satisfying the \PNEInline-condition, the \ANNInline-condition is equivalent to Barr-exactness.
  \end{ulist}

  This amounts to saying that \emph{a homological category is di-exact if and only it it is Barr exact}.

  (\ref{thm:Pushout,AllCokernels->NormalPushout}) is closely related to the question whether a join of two normal subobjects is normal and, in turn, to the $(3\times 3)$-lemma (\ref{thm:(3x3)-LemmaSemiAb}), one of the key results of Section~\ref{sec:3x3-Lemma-Homological}. In particular, we prove in (\ref{thm:JoinPreservesNormalMonos<->Ax6}) that in a homological category joins of pairs of normal subobjects are normal if and only if the \ANNInline-condition is satisfied.
\end{subordinate}

\begin{exercises}
\begin{exercise}[Normal pushout vs. square of normal epimorphisms\SATag]
  \label{exe:SquareNormalEpis-NotPushout}
  In a semiabelian category give an example of a commutative square of normal epimorphisms  which is not a pushout square. %
  \index{regular!pushout}
\end{exercise}
\end{exercises}
\section[Normal Morphisms]{Normal Morphisms in a Semiabelian Category}
\label{sec:NormalMaps-SA}

Let us recall (\ref{def:NormalMap}) that a morphism $f\from X\to Y$ in a \ZExact\ category $\Ctgry{X}$ is normal if it admits a factorization $f=me$ in which $e$ is a normal epimorphism and $m$ is a normal monomorphism. Here we establish properties of normal morphisms which hold in semiabelian categories but may not hold in a \ZExact\ category, and we derive consequences thereof.

The next result is a companion to (\ref{thm:PullbackRecognition-KernelSide}); the difference being that we do not require $\alpha$ or $\xi$ to be kernels here. However, we are requiring them to appear in exact sequences. So, each of them is a normal map.

\begin{proposition}[Induced map of cokernel of square\SATag]
  \label{thm:CoKernelsOfSquares-Normal}%
  If in the morphism of exact sequences below the map $b$ is normal, then so is $\CoKer{\alpha,\beta}$.
  \begin{equation*}
    \xymatrix@R=5ex@C=3em{
    0 \ar[r] & \Ker{\alpha} \ar[d]_{\Ker{\alpha,\beta}} \ar@{{ |>}->}[r] &
    A \ar[r]^-{\alpha} \ar[d]_{a} &
    B \ar@{-{ >>}}[r] \ar[d]_{b} &
    \CoKer{\alpha} \ar[d]^{\CoKer{\alpha,\beta}} \ar[r] &
    0 \\
    0 \ar[r] &
    \Ker{\beta} \ar@{{ |>}->}[r] &
    X \ar[r]_-{\beta} &
    Y \ar@{-{ >>}}[r] &
    \CoKer{\beta} \ar[r] &
    0
    }
  \end{equation*}
\end{proposition}
\begin{proof}
  Inserting (normal) image factorizations into the square on the right, we find the following situation:
  \begin{equation*}
    \xymatrix@R=5ex@C=5em{
    B \ar@{-{ >>}}[r]^-{\pi} \ar@{-{ >>}}[d]_{e} &
    \CoKer{\alpha} \ar@{-{ >>}}[d]^{\varepsilon} \\
    \DiagObj \ar[r]^-{\lambda} \ar@{{ |>}->}[d]_{m} &
    \DiagObj \ar@{{ >}->}[d]^{\mu} \\
    Y \ar@{-{ >>}}[r]_{\rho} &
    \CoKer{\beta}
    }
  \end{equation*}
  The composite $\varepsilon\pi$ is a normal epimorphism, and so $\lambda$ is a normal epimorphism. Consequently, the $\mu\lambda$ is the image factorization of the antinormal composite $\rho m$. Via the \ANNInline-condition we see that $\mu$ is a normal monomorphism. So $\CoKer{\alpha,\beta}$ is a normal map.
\end{proof}

\begin{example}[Conclusion of (\ref{thm:CoKernelsOfSquares-Normal}) not for kernels]
  \label{exa:KernelsOfSquaresNormal-No}%
  In the setting of (\ref{thm:CoKernelsOfSquares-Normal}) the reader may wonder whether $\Ker{\alpha,\beta}$ is normal whenever $a$ is. In general, the answer is `no', since this would imply that a composite of two normal monomorphisms is again normal. Here is an explicit example in $\Grps$: Let $\CyclcGrp{2}=\Set{-1,+1}$ be the cyclic group of order $2$, and let $\ZNr[\CyclcGrp{2}]$ be its integral group ring. When regarded as a module over itself, there is an underlying action of $\CyclcGrp{2}$ on $\ZNr[\CyclcGrp{2}]$. With the associated semidirect product we may form the square of proper maps on the right, along with the indicated kernels.
  \begin{equation*}
    \xymatrix@R=5ex@C=4em{
    \ZNr.(+1) \ar@{{ |>}->}[r] \ar[d]_{\Ker{\alpha,\beta}} &
    \ZNr[\CyclcGrp{2}] \ar@{-{ >>}}[r] \ar@{{ |>}->}[d] &
    \ZNr \ar[d] \\
    \ZNr[\CyclcGrp{2}]\rtimes \CyclcGrp{2} \ar@{=}[r] &
    \ZNr[\CyclcGrp{2}]\rtimes \CyclcGrp{2} \ar[r] &
    0
    }
  \end{equation*}
  Evidently, the map $\Ker{\alpha,\beta}$ is not normal.
\end{example}

\begin{proposition}[Pushouts along normal epis preserve normal maps\SATag]
  In a semiabelian category, a pushout along a normal epimorphism preserves normal maps.
\end{proposition}
\begin{proof}
  In the commutative diagram below, we see the pushout of a normal map $f=me$ along a normal epimorphism $\varphi$ broken down as a concatenation of two steps.
  \begin{equation*}
    \xymatrix@R=5ex@C=4em{
    \DiagObj \ar@{-{ >>}}[r]_-{e} \ar@{-{ >>}}[d]_{\varphi} \ar@/^2ex/[rrr]^-{f} \PushRD{rd}&
    \DiagObj \ar@{=}[r] \ar@{-{ >>}}[d]_{\underline{\varphi}}&
    \DiagObj \ar@{{ |>}->}[r]_-{m} \ar@{-{ >>}}[d]^{\rho} &
    \DiagObj \ar@{-{ >>}}[d]^{\underline{\varphi}} \\
    \DiagObj \ar@{-{ >>}}[r]^-{\underline{e}} \ar@/_4ex/[rrr]_-{\underline{f}} \ar@{.>}@/_2ex/[rr]|-{\ \varepsilon\ }&
    \DiagObj \ar@{-->}[r]^-{t} &
    \Ker{q} \ar@{{ |>}->}[r]^-{\mu} &
    \DiagObj \ar@{-{ >>}}[r]_-{q} &
    \CoKer{\underline{f}}
    }
  \end{equation*}
  The long rectangle and the left hand square are constructed as pushouts. So, the underlined maps are normal epimorphisms. The antinormal composite $\underline{\varphi}m$ is a normal map by the \ANNInline-condition. Its normal factorization $\mu\rho$ may be constructed as
  \begin{equation*}
    \mu = \KerMap{\CoKerMap{\underline{\varphi}m}} = \CoKerMap{\underline{\varphi}f}=\CoKerMap{\underline{f}\varphi} = \CoKerMap{\underline{f}}
  \end{equation*}
  Here, we used repeatedly that a map $\alpha$, and $\alpha$ precomposed by a normal epimorphism have the same cokernels; see (\ref{thm:NormalEpi-Props}). The map $\underline{f}$ factors through $\mu$ via $\varepsilon$. To see that $\varepsilon$ is a normal epimorphism, we argue as follows. The universal property of the pushout on the left yields the map $t$ which renders the entire diagram commutative. We conclude that $t$ is a normal epimorphism via (\ref{thm:NormalEpi-Props}). So, $t$ is a composite of two normal epimorphisms, hence is itself a normal epimorphism by (\ref{thm:NormalEpis-Props-Normal}). - Thus $f$ pushes forward to the normal map $\mu\varepsilon$.
\end{proof}

\begin{subordinate}{}
  \begin{subsubordinate}{On the \ANNInline-condition}
    In a homological category, the conclusion of (\ref{thm:CoKernelsOfSquares-Normal}) is equivalent to the \ANNInline-condition. Indeed, consider an antinormal map $qk$ as in this diagram.
    \begin{equation*}
      \xymatrix@R=5ex@C=4em{
      K \ar@{{ |>}->}[r]^-{k} \ar@{-{ >>}}[d]_{e} &
      X \ar@{-{ >>}}[d]^{q} \\
      I \ar@{{ >}->}[r]_-{m} &
      Q
      }
    \end{equation*}
    The composite $me$ is the image factorization of $qk$ as in (\ref{thm:NEM-Img-Fact-Existence}). Supplementing the kernels of $e$ and $q$, we obtain a pullback diagram, since $m$ is monic. Via the conclusion of (\ref{thm:CoKernelsOfSquares-Normal}), we see that $m$ is normal as well, as was to be shown.
  \end{subsubordinate}

  \begin{subsubordinate}{Normal diagonal maps and semiabelian vs.\ abelian}
    \label{rem:DiagonalsProper<->abelian}
    A defining property of an abelian category $\SACtgry{A}$ is that every monomorphism is a normal monomorphism. Therefore every morphism in $\SACtgry{A}$ is normal. Perhaps surprisingly, this exact same property also separates a semiabelian category from being abelian.
    An even more minimalistic criterion is: A semiabelian category $\SACtgry{X}$ is abelian if and only if every diagonal map $\DgnlOn{X}\from X\to X\times X$ is a normal monomorphism. %
  \end{subsubordinate}

\end{subordinate}

\begin{exercises}

\begin{exercise}[Product of normal monomorphisms\SATag]
  \label{exe:ProductNormalMonos}
  Given monomorphisms $f\from A\to X$ and $g\from A\to Y$ in a semiabelian category, show the following:
  \begin{enumerate}[(i)]
    \item If $\PrdctMapInto{f,g}\from A\to \Prdct{X}{Y}$ is normal, then $f$ and $g$ are normal.
    \item If $f$ and $g$ are normal, find an example in which $\PrdctMapInto{f,g}$ is not normal.
    \item In the category $\Grps$ of groups, show that $f$ and $g$ are normal if and only if $f$ is central in $X$ and $g$ is central in $Y$.
  \end{enumerate}
\end{exercise}
\end{exercises}
\newpage
\section[The Category of Short Exact Sequences - II]{The Category of Short Exact Sequences - II}
\label{sec:Cat-SESs-SACat}

In a di-exact category $\Ctgry{X}$ and, hence, in a semiabelian category a short exact sequence in $\SESCat{X}$ is a di-extension in $\Ctgry{X}$; see Proposition~\ref{thm:ANN<->PointwiseSES}. Here, we show that $\SESCat{X}$ is homological whenever $\Ctgry{X}$ is semiabelian.

\begin{theorem}[$\Ctgry{X}$ semiabelian implies $\SESCat{X}$ is homological\SATag]
  \label{thm:X-SemiAbelian->SES(X)Homological}%
  If $\Ctgry{X}$ is a semiabelian category, then the category $\SESCat{X}$ of short exact sequences in $\Ctgry{X}$ is homological. %
  \index[not]{s!$\SESCat{X}$\IndSep category of short exact sequences in $\Ctgry{X}$}
\end{theorem}
\begin{proof}
  We already know from Theorem \ref{thm:SES(C)IsP-Category} that $\SESCat{X}$ is a category with zero-object which is finitely bicomplete.

  To verify that the \PNEInline-condition holds in $\SESCat{X}$, we show that it holds in the equivalent category $\NMonoCat{X}$. Consider the commutative cube below.
  
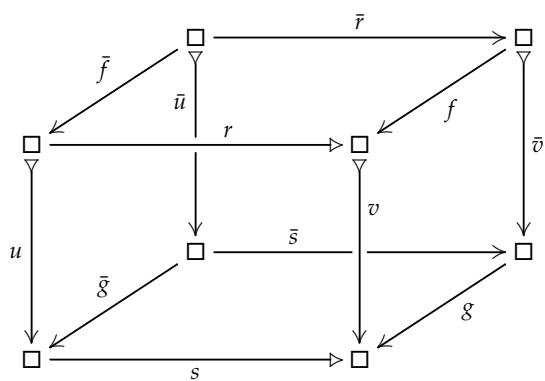
\begin{equation*}
    \xymatrix@R=5ex@C=4em{
    & \DiagObj \ar@{{ |>}->}[dd]|\hole _(0.3){\bar{u}} \ar[ld]_{\bar{f}} \ar[rr]^-{\bar{r}} &&
    \DiagObj \ar@{{ |>}->}[dd] ^{\bar{v}} \ar[ld]^{f} \\
    \DiagObj \ar@{{ |>}->}[dd]_{u} \ar@{-{ >>}}[rr]^(0.6){r} &&
    \DiagObj \ar@{{ |>}->}[dd]^(0.3){v} \\
    & \DiagObj \ar[ld]_{\bar{g}} \ar[rr]|\hole^(0.3){\bar{s}} &&
    \DiagObj \ar[ld] ^(0.4){g}\\
    \DiagObj \ar@{-{ >>}}[rr]_-{s} &&
    \DiagObj
    }
  \end{equation*}
  The front face is a normal epimorphism in $\NMonoCat{X}$, given by a pair $(r,s)$ of normal epimorphism $(r,s)$ in $\Ctgry{X}$; see (\ref{thm:CoKernelsInNEpi(X)-DiExact}). Pulling $(r,s)$ back along the morphism $(f,g)$ in $\NMonoCat{X}$ is accomplished via the bottom and top pullback squares of the cube by (\ref{thm:SESCat(X)-(Co)Limits}). The normal epimorphisms $r$ and $s$ pull back to normal epimorphisms $\bar{r}$ and $\bar{s}$ by (\ref{SA-Ax:NormalEpis-BaseChangePreserved}). Thus, the normal epimorphism $(r,s)$ in $\NMonoCat{X}$ pulls back to the normal epimorphism $(\bar{r},\bar{s})$, again (\ref{thm:CoKernelsInNEpi(X)-DiExact}).

  The \KSGInline-condition is satisfied because in $\NMonoCat{X}$, sectioned epimorphisms, kernels and monomorphisms are all given by their pointwise counterparts.
\end{proof}

\begin{proposition}[Adjunctions with the category of morphisms\SATag]
  \label{thm:Adjunctions-SES(X)<->X^[1]}%
  Associated to semiabelian category $\SACtgry{X}$ are these adjunctions
  \begin{equation*}
    \xymatrix@C=4em{\ArrowCat{\SESCat{X}} \ar@<1ex>[r] \ar@{}[r]|-{\bot} & \SESCat{\SESCat{X}} \ar@<1ex>[l]}
    \qquad\qquad
    \xymatrix@C=4em{\ArrowCat{\SESCat{X}} \ar@<1ex>[r] \ar@{}[r]|-{\top} & \SESCat{\SESCat{X}} \ar@<1ex>[l]}
  \end{equation*}
  The right adjoint of the adjunction on the left sends a short exact sequence of short exact sequences, i.e., a di-extension, to its normal monomorphic part. Its left adjoint constructs from a morphism $(\mu,\xi,\eta)$ of short exact sequences the di-extension shown on the right below.
  \begin{equation*}
    \xymatrix@R=5ex@C=4em{
    &&& \Ker{\kappa} \ar@{{ |>}->}[r]^-{\hat{k}} \ar@{{ |>}->}[d]_{\alpha} &
    \Ker{\xi} \ar@{{ |>}->}[d]_{\beta} \ar@{-{ >>}}[r] &
    \CoKer{\hat{k}} \ar@{{ |>}->}[d]^{\gamma} \\
    K \ar@{{ |>}->}[r]^-{k} \ar[d]_{\kappa} &
    X \ar@{-{ >>}}[r]^-{q} \ar[d]_{\xi} &
    Q \ar[d]^{\rho} &
    K \ar@{{ |>}->}[r]^-{k} \ar[d] &
    X \ar@{-{ >>}}[r]^-{q} \ar[d] &
    Q \ar[d] \\
    L \ar@{{ |>}->}[r]_-{l} &
    Y \ar@{-{ >>}}[r]_-{r} &
    R &
    \CoKer{\alpha} \ar@{{ |>}->}[r] &
    \CoKer{\beta} \ar@{-{ >>}}[r] &
    \CoKer{\gamma}
    }
  \end{equation*}
  The left adjoint of the adjunction on the right `forgets to normal epimorphisms of short exact sequences'. The right adjoint constructs from $(\kappa,\xi,\rho)$ a di-extension by taking the cokernel of its kernel. \NoProof
\end{proposition}

\begin{subordinate}{}
  \begin{subsubordinate}{$\SESCat{X}$ if $\EuScript{X}$ is abelian}
    For the category $\AbGrps$ of abelian groups, its associated category $\SESCat{\AbGrps}$ is preabelian but fails to be abelian. It is homological, but fails to be semiabelian.
  \end{subsubordinate}
\end{subordinate}
\section[The Lattice of Subobjects]{The Lattice of Subobjects}%
\label{sec:LatticeOfSubjects}

In any category $\SACtgry{X}$ the subobjects of a given object $X$ form a class $\SubObjcts{X}$ which is partially ordered by inclusion `$\leq$'. Here, we analyze $(\SubObjcts{X},\leq)$ further in the two settings where $\Ctgry{X}$ is homological, respectively semiabelian.

The partial order $\SubObjcts{X}$ has $X$ itself as its greatest element, and has the zero objects as its least element. We  investigate the additional structures given by `meet' and `join'.

In a homological category $\Ctgry{X}$, the meet of two subobjects $K$ and $L$ of $X$ exists, and may be constructed as the pullback of $K\to X\leftarrow L$; see (\ref{exe:IntersectionConstruct-Via-Pullback}). If $\Ctgry{X}$ has underlying sets, then this yields the set-intersection of $K$ and $L$.

In (\ref{def:Meet/Join-Subobjects}), we characterized the join $K\join L$ of $K$ and $L$ as the least subobject of $X$ containing $K$ and $L$, provided such a least subobject exists. Under appropriate conditions, it may be constructed as the `intersection' of all subjects of $X$ containing $K$ and $L$; see (\ref{thm:Meet/Join-InCompleteCat}).  In a homological category, the situation is simpler: $K\join L$ may be constructed as the image of the canonical map $K+L\to X$; see (\ref{thm:JoinSubobjects-Construction}). Thus, meet and join determine a bounded lattice structure on $\SubObjcts{X}$.

Whenever the class $\SubObjcts{X}$ is small for every $X$ in $\Ctgry{X}$, we obtain a category $\SubObjctCat{X}$ with
\begin{ulist}
  \item \emph{objects}\quad the bounded lattices $\SubObjcts{X}$, with $X$ in $\Ctgry{X}$, and
  \item \emph{morphisms}\quad the order preserving functions $\SubObjcts{X}\to \SubObjcts{Y}$.
\end{ulist}
The operation which sends an object $X$ to its lattice of subobjects is a contravariant functor $\SubObjctsPullFunc\from \Ctgry{X}\to \SubObjctCat{X}$. This is so because pullbacks preserve monomorphisms.

Similarly, we have the lattice $\NSubObjcts{X}$ of normal subobjects of $X$. The meet of two normal subobjects is again normal, by (\ref{thm:NormalSubobjects-Intersection}). Their join, however, need not be normal in general. In a homological category $\Ctgry{X}$ we show that joins of normal subobjects are always normal if and only if $\Ctgry{X}$ is semiabelian; see (\ref{thm:DiExtension-iff-JoinNormal}).

As applications of these general properties of the join and meet operations, we show:
\begin{ulist}
  \item (\ref{thm:JoinInSplitSES})\quad Given a split extension, the join of the end objects is equals the middle object. \HTag
  \item (\ref{thm:SplitExtViaSubobj})\quad A normal subobject $K\leq X$ and a subobject $Y\leq X$ determine a split short exact sequence if and only if $K\meet Y=\ZeroObject$ and $K\join Y=X$. \HTag
  \item (\ref{thm:Noether})\quad The Second Isomorphism Theorem. \HTag
\end{ulist}
Turning to details:

\begin{proposition}[Construction of join of subobjects via image factorization\HTag]
  \label{thm:JoinMorphisms-Construction} 
  \label{thm:JoinSubobjects-Construction}
  In a homological category, the join of subobjects $m\from M\Mono X$ and $n\from N\Mono X$ exists and is represented by the monomorphism $i\from I\Mono X$ in the image factorization %
  \index{join!of subobjects: construction}
  \begin{equation*}
    \xymatrix@R=5ex@C=4em{
    M+N \ar@/^3ex/[rr]^-{\SumMapOutOf{m,n}} \ar@{-{ >>}}[r]_-{q=\SumMapOutOf{m',n'}} &
    I \ar@{{ >}->}[r]_-{i} &
    X
    }
  \end{equation*}
  of $\SumMapOutOf{m,n}\from M+N\to X$.
\end{proposition}
\begin{proof}
  From (\ref{thm:CoKer=NormalEpi=RegEpi=EffectiveEpi}) we know that the normal epimorphism $q=\SumMapOutOf{m',n'}$ is an extremal epimorphism. This means that $m'$ and $n'$ are jointly extremally epimorphic. By (\ref{thm:ExtremalEpiCharacterizationOfJoin}), $i\from I\to X$ represents the join of $m$ and $n$.
\end{proof}

In following corollary is a well known property of many familiar varieties of algebras.

\begin{corollary}[Join of normal subobjects\SATag]
  \label{thm:JoinNormal}
  In a semiabelian category, the join of two normal subobjects of a given object exists and is again normal.
\end{corollary}
\begin{proof}
  The pullback square of two normal subobjects $K$ and $L$ of $X$ forms the initial square of a di-extension by (\ref{thm:DPN-DiExtensivePull/Push}). Then Proposition \ref{thm:9Diagram-JntEpiOnKernel,JntMonoOnCoKernel} explains why $K\join L$ is normal in $X$.
\end{proof}

\begin{corollary}[Universal map from sum with normal constituents\SATag]
  \label{thm:ProperMapByAddingProperMaps}
  \label{thm:SumNormalMaps->NormalMap}
  \label{thm:UMapFromSum-NormalConstituents}
  If $f\from X\to Z$ and $g\from Y\to Z$ are normal, then so is $\SumMapOutOf{f,g}\from X+Y\to Z$.
\end{corollary}
\begin{proof}
  By (\ref{thm:ImageFactorizationCommutesBinarySums}) the image of $\SumMapOutOf{f,g}$ is the join of the images of $f$ and $g$. Now the claim follows from (\ref{thm:JoinNormal}).
\end{proof}

\begin{proposition}[Split short exact sequence and join\HTag]
  \label{thm:JoinInSplitSES}%
  \label{thm:JoinMorphisms-Properties}
  In a split short exact sequence, as below, $X=K\join Y$.
  \begin{equation*}
    \xymatrix@C=4em{K \ar@{{ |>}->}[r]^-{k} &
    X \ar@{-{ >>}}@<-.5ex>[r]_-{f} &
    Y \ar@{{ >}->}@<-.5ex>[l]_-{s}}
  \end{equation*}
\end{proposition}
\begin{proof}
  The structure maps from $K$ and $Y$ into their sum yield this commutative diagram:
  
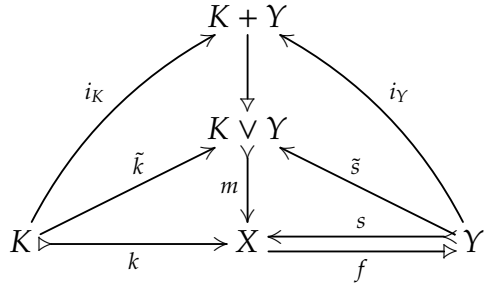
\begin{equation*}
    \xymatrix@R=5ex@C=5em{
    & K+Y \ar@{-{ >>}}[d] \\
    & K\join Y \ar@{{ >}->}[d]_{m} \\
    K \ar@{{ |>}->}[r]_-{k} \ar@/^2ex/[ruu]^{\InclsnOf{K}} \ar[ru]^(0.55){\tilde{k}} &
    X \ar@{-{ >>}}@<-.5ex>[r]_-{f} &
    Y \ar@{{ >}->}@<-.5ex>[l]_-{s} \ar@/_2ex/[luu]_{\InclsnOf{Y}} \ar[lu]_(0.55){\tilde{s}}
    }
  \end{equation*}
  By Axiom~\ref{SA-Ax:KSG}, $k$ and $s$ are jointly extremally epimorphic. So, $m$ is an isomorphism; i.e., $K\join Y=X$.
\end{proof}

\begin{corollary}[Split extension as a pair of subobjects\HTag]
  \label{thm:SplitExtViaSubobj}
  For any object $X$ in a homological category the following conditions are equivalent:
  \begin{tfae}
    \item $X$ is the middle object in a split short exact sequence:
    \begin{equation*}
      \xymatrix@C=4em{
      K \ar@{{ |>}->}[r]^-{k} &
      X \ar@{-{ >>}}@<-.5ex>[r]_-{f} &
      Y \ar@{{ >}->}@<-.5ex>[l]_-{s}
      }
    \end{equation*}
    \item $K\normal X$ and  $Y\leq X$ such that $K\meet Y=0$ and $K\join Y=X$.
  \end{tfae}
\end{corollary}
\begin{proof}
  (i) $\Rightarrow$ (ii)\quad With (\ref{thm:JoinInSplitSES}), we see that $K\join Y=X$. Further, $k=\Ker{f}$ is a normal subobject of $X$, while  $s$ represents $Y$ as a subobject of $X$. To see that 	$K\meet Y=0$, consider this pullback diagram:
  \begin{equation*}
    \xymatrix@R=5ex@C=4em{
    K\meet Y \ar@{{ |>}->}[r]^-{\bar{k}} \ar@{{ >}->}[d]_{\bar{s}} \PullLU{rd} &
    Y \ar@{{ >}->}[d]^{s} \\
    K \ar@{{ |>}->}[r]_-{k} &
    X
    }
  \end{equation*}
  The map $\bar{k}$ is a monomorphism by (\ref{thm:KernelFunctor-Props}). Also, we have $s\bar{k}=k\bar{s}$ and, hence,
  \begin{equation*}
    \bar{k} = fs\bar{k} = fk\bar{s} = \ZeroMap
  \end{equation*}
  Thus $K\meet Y$ is the zero subobject of $X$.

  (ii) $\Rightarrow$ (i)\quad From $K\normal X$ and $Y\leq X$, we construct this commutative diagram:
  \begin{equation*}
    \xymatrix@R=5ex@C=4em{
    0 \ar[r] \ar[d] \PullLU{rd} &
    Y \ar@{=}[r] \ar@{{ >}->}[d]^-{\sigma} &
    Y \ar[d]^{f\sigma} \\
    K \ar@{{ |>}->}[r]_-{k} &
    X \ar@{-{ >>}}[r]_-{f} &
    \CoKer{k}
    }
  \end{equation*}
  The pullback construction (\ref{exe:IntersectionConstruct-Via-Pullback}) of $K\meet Y$ yields the square on the left. So, $f\sigma$ is  a monomorphism by Lemma~\ref{thm:PullbackRecognition-KernelSide}.  By Corollary~\ref{thm:NormalEpiRecognition-CoKernelSide} it is a normal epimorphism as well. By (\ref{thm:IsomorphismRecognition}) it is an isomorphism. So, $s\DefEq \sigma\Comp (f\sigma)^{-1}$ is a section of $f$.
\end{proof}

\begin{theorem}[Second Isomorphism Theorem\HTag]  
  \label{thm:Noether}%
  \index{Second Isomorphism Theorem}\index{Isomorphism Theorem!Second}%
  For subobjects $K$, $L\leq X$ with $L\normal (L\join K)$, we have that
  \[
    \frac{K}{L\meet K}\cong\frac{L\join K}{L}.
  \]
\end{theorem}
\begin{proof}
  From the subobjects $K$ and $L$ of $L\join K$, we construct the commutative square on the left as a pullback. The map $\bar{\lambda}$ is a normal monomorphism by (\ref{thm:KernelFunctor-Props}). So, taking cokernels, yields this morphism of short exact sequences:
  \begin{equation*}
    \xymatrix@R=5ex@C=3em{
    L\meet K \PullLU{rd} \ar@{{ |>}->}[r]^-{\bar{\lambda}} \ar@{{ >}->}[d]_{\bar{u}} &
    K \ar@{-{ >>}}[r] \ar@{{ >}->}[d]^{u} &
    K/(L\meet K) \ar[d]^{v} \\
    L \ar@{{ |>}->}[r]_{\lambda} &
    L\join K \ar@{-{ >>}}[r] &
    (L\join K)/L,
    }
  \end{equation*}
  The map $v$ is a monomorphism by Lemma~\ref{thm:PullbackRecognition-KernelSide}. By the join recognition criterion (\ref{thm:ExtremalEpiCharacterizationOfJoin}), the inclusions of $K$ and $L$ are jointly extremal-epimorphic. With (\ref{thm:NormalEpiRecognition-CoKernelSide}) we see that $v$ is a normal epimorphism. By (\ref{thm:IsomorphismRecognition}), $v$ is an isomorphism.
\end{proof}

\begin{proposition}[Normality of join in a homological category\HTag]
  \label{thm:DiExtension-iff-JoinNormal}
  Given an object $X$ in a homological category, the following conditions are equivalent for normal subobjects $K$ and $L$ of $X$:
  \begin{tfae}
    \item $K$ and $L$ determine a di-extension with central object $X$.
    \item $K\join L$ is normal in $X$.
    \item The pushout of $X\to X/K$ and $X\to X/L$ is a normal pushout.
  \end{tfae}
\end{proposition}
\begin{proof}
  (I) $\Rightarrow$ (II) follows via  Proposition~\ref{thm:9Diagram-JntEpiOnKernel,JntMonoOnCoKernel}.

  (I) $\Rightarrow$ (III) follows via the (semi)normal pushout recognition criterion (\ref{thm:SemiNormalPushout-Recognition}).

  (III) $\Rightarrow$ (I) follows from via (\ref{thm:DPN-DiExtensivePull/Push}).

  (II) $\Rightarrow$ (I) If $K\join L$ is normal in $X$, then we obtain this morphism of short exact sequences:
  \begin{equation*}
    \xymatrix@R=5ex@C=4em{
    L \ar@{{ |>}->}[r] \ar@{{ |>}->}[d] &
    X \ar@{-{ >>}}[r] \ar@{=}[d] &
    X/L \ar@{-{ >>}}[d] \\
    K\join L \ar@{{ |>}->}[r] &
    X \ar@{-{ >>}}[r] &
    X/(K\join L)
    }
  \end{equation*}
  Combining the Pure Snake Lemma with (\ref{thm:Noether}), we obtain this exact sequence:
  \begin{equation*}
    \xymatrix@R=5ex@C=4em{
    L \ar@{{ |>}->}[r] &
    K\join L \ar[r] \ar@{-{ >>}}[d] &
    X/L \ar@{-{ >>}}[r] &
    X/(K\join L) \\
    & K\join L/L \ar@{=}[r] &
    K/(K\meet L) \ar@{{ |>}->}[u]
    }
  \end{equation*}
  Similarly, we obtain the short exact sequence $L/K\meet L \NMono X/K \NEpi X/K\meet L$. This yields the di-extension.
  \begin{equation*}
    \xymatrix@R=4ex@C=3em{
    K\meet L \ar@{{ |>}->}[r] \ar@{{ |>}->}[d] &
    K \ar@{-{ >>}}[r] \ar@{{ |>}->}[d] &
    K/L\meet K \ar@{{ |>}->}[d] \\
    L \ar@{{ |>}->}[r] \ar@{-{ >>}}[d] &
    X \ar@{-{ >>}}[r] \ar@{-{ >>}}[d] &
    X/L \ar@{-{ >>}}[d] \\
    L/L\meet K \ar@{{ |>}->}[r] &
    X/K \ar@{-{ >>}}[r] &
    X/K\join L
    }
  \end{equation*}
  The proof of (\ref{thm:DiExtension-iff-JoinNormal}) is complete.
\end{proof}

\begin{corollary}[Join of normal subobjects is normal iff di-exact\HTag]
  \label{thm:JoinPreservesNormalMonos<->Ax6}
  \label{thm:JoinPreservesNormalMonos<->DiExact}%
  In a homological category $\Ctgry{X}$, the join of a pair of normal subobjects of an object $X$ is normal if and only if $\Ctgry{X}$ is di-exact. \NoProof
\end{corollary}

\begin{corollary}[Universal map from sum normal iff di-exact\HTag]
  \label{thm:JoinProperMaps<->Ax6}
  \label{thm:SumNormalMaps<->DiExact}
  \label{thm:UMapFromSumNormal<->DiExact}%
  In a homological category $\Ctgry{X}$, the sum of two normal maps is a normal map if and only if $\Ctgry{X}$ is di-exact. \NoProof
\end{corollary}
\begin{proof}
  We know (\ref{thm:UMapFromSum-NormalConstituents}) that, in a semiabelian category the sum of two normal maps is normal. Conversely, (\ref{thm:JoinPreservesNormalMonos<->DiExact}) shows that, if the sum of normal maps is normal, then $\Ctgry{X}$ is semiabelian.
\end{proof}

Consequently, a homological category is semiabelian if and only if every universal map out of a binary sum with normal constituents is again a normal map.

Next, we infer the well-known fact that the lattice of $\NSubObjcts{X}$ of normal subobjects of an object $X$ in a semiabelian category is a \Defn{modular} lattice. This means that for all $K$, $L$, $M\normal X$ where $M\leq L$, the equality $M\join (K\meet L)=(M\join K)\meet L$ holds. Here we freely use (\ref{thm:JoinNormal}).%
\index{modular lattice}\index{lattice!modular}


\begin{proposition}[Modularity of lattice of normal subobjects\SATag]
  \label{thm:Modularity}
  \label{thm:LatticeNormalSubobjects-Modular}%
  For normal subobjects $K$, $L$, $M\normal X$ with $M\leq L$,
  \begin{equation*}
    M\join (K\meet L)=(M\join K)\meet L
  \end{equation*}
  as normal subobjects of $X$.
\end{proposition}
\begin{proof}
  Consider this commutative diagram whose rows are short exact:
  \begin{equation*}
    \xymatrix@R=5ex@C=3em{
    K\meet L \ar@{{ |>}->}[r] \ar@{{ |>}->}[d] &
    K \ar@{-{ >>}}[r] \ar@{{ |>}->}[d] &
    K/(K\meet L) \ar[d]^{v} \\
    M\join (K\meet L) \ar@{{ |>}->}[r] \ar@{{ |>}->}[d] &
    M\join K \ar@{-{ >>}}[r] \ar@{{ |>}->}[d] &
    (M\join K)/(M\join (K\meet L)) \ar[d]^{w} \\
    L \ar@{{ |>}->}[r] &
    K\join L \ar@{-{ >>}}[r] &
    (K\join L)/L
    }
  \end{equation*}
  Note that by (\ref{thm:JoinNormal}), all subobjects in this diagram are normal. The desired identity follows, once we show that the rectangle on the bottom left is a pullback. Indeed, the rectangle $(K\meet L)\rightrightarrows (M\join K)$ is a pullback. As in the proof of Corollary~\ref{thm:Noether} the composite $w\Comp v$ is an isomorphism. Thus $v$ is monic, and the top left square is a pullback by (\ref{thm:PullbackRecognition-KernelSide-1}).

  We claim that $v$ is an isomorphism as well. Here, we used that the extremal epimorphism $M+K\to M\join K$ factors through $(M\join (M\meet L))+K$, which implies that the inclusions of $M\join (M\meet L)$ and $K$ in $M\join K$ are jointly extremally epimorphic as well.

  Consequently, $w$ is an isomorphism. By Proposition~\ref{thm:PullbackRecognition-KernelSide-1} the bottom left square is a pullback, as claimed.
\end{proof}

The following example shows that working with normal subobjects is essential for the modularity result (\ref{thm:Modularity}) to hold.

\begin{example}[Non-modularity of the lattice of all subobjects]
  \label{exa:LatticeSubojects-NotModular}%
  The lattice of subgroups of a group need not be modular. For an example consider the lattice of subgroups of the dihedral group of order $8$.
\end{example}

\begin{exercises}

\begin{exercise}[Properties of meet and join]
  \label{exe:Meet-Join-NotDistributive}
  Determine whether all subobject lattices in a semiabelian category are distributive; i.e.\ are the constructions $A\meet (B\join C)$ and $(A\meet B)\join (A\meet C)$ on subobjects $A$, $B$, $C$ of a given object $X$ always isomorphic?
\end{exercise}

\end{exercises}
\section[Alternate Characterizations: Barr-exact Mal'tsev Categories]{Alternate Characterizations:\\ Barr-exact Mal'tsev Categories}
\label{sec:ExactMaltsev}

The first important insight presented in this section is Theorem \ref{thm:ReflexiveRelationInSACategory}, which says that, in a semiabelian category, any reflexive relation $(R,d_1,d_2,e)$ on an object $X$ is an effective equivalence relation, so a kernel pair of a morphism $q\from X\to Q$. For convenience, we can always let~$q$ be the coequalizer of $d_1$ and $d_2$, and thus obtain a diagram
\begin{equation*}
  \xymatrix@R=5ex@C=4em{
  R \BiCart{rd} \ar@{-{ >>}}[r]^{d_1} \ar@{-{ >>}}[d]_{d_2} &
  X \ar@{-{ >>}}[d]^{q}\\
  X \ar@{-{ >>}}[r]_q & Q
  }
\end{equation*}
which is both a pushout and a pullback. This implies that a semiabelian category is at the same time a Mal'tsev category (Definition~\ref{def:MaltsevCategory}) and a Barr-exact category (Definition~\ref{def:BarrExact}). Actually, the former conclusion holds in homological categories; see (\ref{thm:HomologicalThenMaltsev}). But the latter happens to be characteristic of the separation between homological and semiabelian categories. In other words, \emph{a homological category is di-exact if and only if it is Barr-exact} (\ref{thm:BarrExactIffDiExact-InSACat}): proving this is the main objective of this section.

\begin{example}{Barr-exactness of $\Grps$}
  \label{exa:Grps-Is-BarrExact}%
  \index{Barr exact!$\Grps$}%
  We explain why the category $\Grps$ of groups is Barr exact: Suppose $R\subseteq X\times X$ is an internal reflexive relation on a group $X$. To prove that $R$ is a symmetric relation, take $(x,y)\in R$. Then also $(x^{-1},x^{-1})$ and $(y^{-1},y^{-1})$ are in $R$, by reflexivity. Using that $R$ is a subgroup of $X\times X$, we see that $(y^{-1},x^{-1})=(x^{-1},x^{-1})(x,y)(y^{-1},y^{-1})\in R$. Then also $(y,x)\in R$, for the same reason. Hence $R$ is symmetric. Similar reasoning shows that $R$ is transitive.
\end{example}

\begin{theorem}[Semiabelian categories are Barr-exact Mal'tsev\SATag]
  \label{thm:ReflexiveRelationInSACategory}%
  In a semiabelian category $\Ctgry{X}$, every reflexive relation is an effective equivalence relation. As a consequence, $\Ctgry{X}$ is a Barr-exact Mal'tsev category. %
  \index{semiabelian category!is Barr-exact}\index{semiabelian category!is Mal'tsev category}%
\end{theorem}
\begin{proof}
  Let $(R,d_{1},d_{2},e)$ be a reflexive relation on an object $X$ of $\SACtgry{X}$. Since $d_1$ and $d_2$ are normal epimorphisms, their pushout is a normal pushout (\ref{thm:Pushout,AllCokernels->NormalPushout}):
  \begin{equation*}
    \xymatrix@R=5ex@C=4em{
    R \PushRD{rd} \ar@{-{ >>}}[r]^{d_1} \ar@{-{ >>}}[d]_{d_2} &
    X \ar@{-{ >>}}[d]^{q'} \\
    X \ar@{-{ >>}}[r]_q &
    Q
    }
  \end{equation*}
  Furthermore, $q=qd_2e=q'd_1e=q'$. So taking the pullback of $q$ along $q'$, yields the kernel pair of $q$. The comparison arrow $R\to \KrnlPr{q}$ is now both a normal epimorphism (normal pushout property of the square) and a monomorphism ($R$ is a relation). So (\ref{thm:IsomorphismRecognition}), it is an isomorphism. In particular, the reflexive relation $R$ is a kernel pair.
\end{proof}

Theorem \ref{thm:ReflexiveRelationInSACategory} allows us to obtain a recognition criterion for reflexive coequalizers\index{reflexive coequalizer}\index{coequalizer!reflexive}, valid in semi-abelian categories (but not in homological categories!). Let
\begin{equation}\label{diag:CoeqRecog}
  \vcenter{\xymatrix@!0@R=5em@C=4em{
  R \ar@<1ex>[rr]^-{d_1} \ar@<-1ex>[rr]_-{d_2} \ar[rd]_-{r=\PrdctMapInto{d_1,d_2}} &&
  B \ar[ll]|-{\ e\ } \ar[dl] \ar[r]^-{f} & A\\
  & \KrnlPr{f} \ar@<1ex>[ru]^(.4){\pi_{1}} \ar@<-1ex>[ru]_(.4){\pi_{2}}}}
\end{equation}
be a reflexive graph with its coequaliser $f$, the induced kernel pair $(\KrnlPr{f},\pi_{1},\pi_{2})$ and the comparison morphism $r$. By (\ref{thm:Pushout,AllCokernels->NormalPushout}), also $r$ is normal epimorphism. But in fact, the converse also holds:

\begin{proposition}[Coequalizer recognition criterion\SATag]\label{thm:CoequalizerRecognitionSA}
  In \eqref{diag:CoeqRecog}, suppose a normal epimorphism $f$ satisfies $fc=fd$. Then $f$ is the coequalizer of $c$ and $d$ if and only if $r$ is a normal epimorphism.
\end{proposition}
\begin{proof}
  One implication is explained above; the other is a variation on the proof of Proposition~\ref{thm:NormalEpi-Props}, item (ii).
\end{proof}

To continue, recall the direct image of a relation along a normal epimorphism (\ref{def:DirectImageOfRelation}).

\begin{proposition}[Normal pushout recognition via kernel pairs\HTag]
  \label{thm:ImageKernelPairThenPushout,AllCokernels->NormalPushout}%
  In a homological category the following conditions are equivalent:
  \begin{tfae}
    \item A direct image of an effective equivalence relation is is again an effective equivalence relation.
    \item Every pushout of a normal epimorphism along a normal epimorphism is a normal pushout.
  \end{tfae}
\end{proposition}
\begin{proof}
  The proof follows the lines of the proof of (\ref{thm:Pushout,AllCokernels->NormalPushout}): it suffices to replace kernels by kernel pairs, and use (\ref{thm:RegularPushoutRecognizeKernelPair}) instead of (\ref{thm:NormalPushOut-Recognize-H}).
\end{proof}

Recall, however, that in a homological category, equivalence relations are always preserved by direct images: Lemma~\ref{thm:directimageeqrel}.

\begin{proposition}[Barr-exactness implies the \ANNInline-condition\HTag]
  \label{thm:BarrExactimplies(Ax.6)}%
  In a homological category, if every equivalence relation is a kernel pair, then the \ANNInline-condition holds.
\end{proposition}
\begin{proof}
  We know by Lemma~(\ref{thm:directimageeqrel}) that the direct image of an equivalence relation along a normal epimorphism is still an equivalence relation. So if every equivalence relation is effective, then effective equivalence relations are preserved by direct images. Hence every pushout of a normal epimorphism along a normal epimorphism is a normal pushout by Proposition~(\ref{thm:ImageKernelPairThenPushout,AllCokernels->NormalPushout}), so that Proposition~(\ref{thm:Pushout,AllCokernels->NormalPushout}) gives the \ANNInline-condition.
\end{proof}

\begin{theorem}[Semiabelian iff Barr exact for homological categories\HTag]
  \label{thm:BarrExactIffDiExact-InSACat}%
  A homological category is semiabelian if and only if it is Barr exact.\NoProof
\end{theorem}
\newpage
\section[Alternate Characterizations: On the Existence of Colimits]{Alternate Characterizations:\\ On the Existence of Colimits}
\label{sec:ExistenceColimits}

The Mal'tsev property allows us to prove that in a finitely complete pointed category $\Ctgry{X}$ which satisfies the \PNEInline, \KSGInline, and \ANNInline-conditions, the existence of binary coproducts and coequalizers of kernel pairs implies that $\Ctgry{X}$ is finitely cocomplete. It is this condition which separates Barr-exact Borceux-Bourn--homological categories from semiabelian ones.

Here, we clarify the role played by the existence of coproducts: we show that any Barr-exact Mal'tsev category with finite coproducts is finitely cocomplete. In (\ref{thm:ColimitsBis}) we show that, if binary sums exist in a BB-homological category $\Ctgry{X}$ where every equivalence relation is effective, then $\Ctgry{X}$ is finitely cocomplete. So, it is semiabelian.

\begin{proposition}[Colimits by coproducts and coequalizers of kernel pairs]\label{thm:Colimits}
  Any Barr-exact Mal'tsev category with finite coproducts is finitely cocomplete.
\end{proposition}
\begin{proof}
  We know, see for instance~\cite{SMacLane1998}, that every finite colimits may be constructed as the coequalizer of a pair of parallel arrows between suitable finite sums. Here all finite coproducts exist by assumption. In this context, any coequalizer may be obtained as the coequalizer of a reflexive graph. Indeed, given a pair of parallel arrows $p_1$, $p_2\from G\to X$, we may obtain the coequalizer of $p_1$, $p_2$ as the coequalizer of
  \begin{equation*}
    \xymatrix@R=5ex@C=5em{
    G+X \ar@<1ex>[r]^-{\langle p_1, 1_X\rangle} \ar@<-1ex>[r]_-{\langle p_2, 1_X\rangle} &
    X \ar[l]|-{\ \iota_{2}\ }
    }
  \end{equation*}
  because for any morphism $f\from X\to Y$, we have $f\Comp p_1=f\Comp p_2$ if and only if $f\Comp \langle p_1, 1_X\rangle=f\Comp \langle p_1, 1_X\rangle$.

  So, we must construct coequalizers of reflexive graphs from coequalizers of kernel pairs. If every kernel pair admits a coequalizer, then this is done as follows. Given a reflexive graph
  \begin{equation*}
    \xymatrix@R=5ex@C=4em{
    G \ar@<1ex>[r]^-{p_1} \ar@<-1ex>[r]_-{p_2} &
    X \ar[l]|{\ i\ }
    }
  \end{equation*}
  we wish to obtain the coequalizer of $p_1$, $p_2$. First we factor $(p_1,p_2)\from {G\to X\times X}$ as a regular epimorphism $q\from {G\to R}$, followed by a monomorphism $(d_1,d_2)\from {R\to X\times X}$. This is possible because a Barr-exact category is always regular, so that it admits regular image factorizations: see~(\ref{thm:ImageInRegular}). We thus obtain a reflexive relation $(R,d_1,d_2,e)$, where $e=qi\from X\to R$ is a common splitting for $d_1$ and $d_2$. Now, a coequalizer of $d_1$, $d_2$ is also a coequalizer of $p_1$, $p_2$. So we reduced the problem to constructing coequalizers of reflexive relations. Since we are in a Barr-exact Mal'tsev category, these are all effective equivalence relations. This completes the proof.
\end{proof}

\begin{corollary}[Existence of finite colimits]
  \label{thm:ColimitsBis}%
  In a BB-homological category where every equivalence relation is effective, the existence of binary coproducts implies that $\Ctgry{X}$ is finitely cocomplete.
\end{corollary}
\begin{proof}
  In a pointed category with binary coproducts, all finite coproducts exist, since they may be obtained as repeated binary coproducts---with the exception of the nullary case, where the coproduct is given by the zero object. By (\ref{thm:HomologicalThenMaltsev}), we find ourselves in a Barr-exact Mal'tsev category. Hence we may apply Proposition~\ref{thm:Colimits}, which proves that the category is finitely cocomplete.
\end{proof}

\begin{theorem}[Janelidze--Márki--Tholen definition of semiabelian categories]
  \label{thm:JMTvsUS}%
  A category is semiabelian in the sense of (\ref{def:SACategory}) if and only if it is pointed, Barr-exact and Bourn-protomodular with finite coproducts.
\end{theorem}
\begin{proof}
  Proposition~\ref{thm:BBHom} already tells us that a category is pointed, finitely complete, with coequalizers of kernel pairs, and such that \PNEInline\ and \KSGInline\ hold, if and only if it is pointed, regular and protomodular.

  Together with (\ref{thm:ReflexiveRelationInSACategory}) this implies that a semiabelian category is always pointed, Barr-exact and Bourn-protomodular with finite coproducts.

  For the converse, we use (\ref{thm:ColimitsBis}) on those axioms to find a homological category---which is then semiabelian by Proposition~\ref{thm:BarrExactimplies(Ax.6)}.
\end{proof}
\chapter[Internal Structures]{Internal Structures}
\label{chap:InternalStructures}

\begin{center}
  \textbf{Leitfaden for Chapter \ref{chap:InternalStructures}}
\end{center}

\bigskip

\begin{equation*}
  \xymatrix@R=9ex@C=6em{
  *+[F-,]{\txt{\sffamily (\ref{sec:InternalMagmas-Monoids-Groups}) Internal Magmas, Monoids, Groups}}\ar[d] \\
  *+[F-,]{\txt{\sffamily (\ref{sec:BiProducts}) Biproducts and Linear Categories}} \\
  *+[F-,]{\txt{\sffamily (\ref{sec:InternalGraphs/Relations}) Internal Graphs and Internal Relations}}
  }
\end{equation*}
\newpage
\section{Internal Magmas, Monoids, Groups}
\label{sec:InternalMagmas-Monoids-Groups}

Whenever we equip a set with the structural data for a group $G$, then every set $X$ yields a group $\Hom{X}{G}$ under pointwise multiplication: the product of copies of $G$ indexed by the elements of $X$. Thus $\Hom{-}{G}$ is a contravariant functor from the category of sets into the category of groups. Here we explain how this phenomenon applies, more generally, to categories with a which admit finite products and, hence, have a terminal object.

\begin{definition}[Internal magma, monoid, group]
  \label{def:InternalMagma/Monoid/Group}
  Let $\Ctgry{X}$ be a category with finite products and terminal object $t$. An \Defn{internal magma} structure on an object $X$ in $\Ctgry{X}$ is given by a couple $(X,\mu)$, with $\mu\from \Prdct{X}{X}\to X$ an arbitrary morphism. A unit in $X$ is given by a morphism $e\from \tau\to X$. An \Defn{(internal) unitary magma} $(X,\mu,e)$ is a magma $(X,\mu)$ with a unit $e$ for which the diagrams (U) below commute. %
  \index{internal!unitary magma}%
  \index{internal!magma}%
  \index{internal!monoid}\index{monoid}%
  \index{internal!group}\index{group}%
  \index{commutative!internal magma/monoid/group}%
  \begin{equation*}
    \xymatrix@R=5ex@C=4em{
    X \ar[d]_{\DgnlOn{X}} \ar@{=}[r] \ar@{}[dr]|-{\ \text{(U)}\ } &
    X \ar@{}[dr]|-{\ \text{(U)}\ } &
    X \ar@{=}[l] \ar[d]^{\DgnlOn{X}} &
    X\prdct X\prdct X \ar[r]^-{\Prdct{\IdMapOn{X}}{\mu}} \ar[d]_{\Prdct{\mu}{\IdMapOn{X}}} \ar@{}[dr]|-{\ \text{(A)}\ } &
    \Prdct{X}{X} \ar[d]^{\mu} \\
    \Prdct{X}{X} \ar[r]_-{\IdMapOn{X}\times \hat{e}} &
    \Prdct{X}{X} \ar[u]_{\mu} &
    \Prdct{X}{X} \ar[l]^-{\hat{e}\times \IdMapOn{X}} &
    \Prdct{X}{X} \ar[r]_-{\mu} &
    X
    }
  \end{equation*}
  Here $\hat{e}$ is the composite $X\to \tau \XRA{e} X$. An internal unitary magma is an \Defn{internal monoid} if diagram (A) above commutes. An internal \Defn{group} is a quadruple $(X,\mu,e,i)$ in which $(X,\mu,e)$ is an internal monoid, and the map $i\from X\to X$ renders the diagram below commutative. %
  \index{inverse map!of internal group}%
  \begin{equation*}
    \xymatrix@R=5ex@C=4em{
    X \ar[r]^-{\DgnlOn{X}} \ar[d]_{\DgnlOn{X}} \ar[rrd]|-{\ \hat{e}\ } &
    \Prdct{X}{X} \ar[r]^-{\Prdct{\IdMapOn{X}}{i}} &
    \Prdct{X}{X} \ar[d]^{\mu} \\
    \Prdct{X}{X} \ar[r]_-{\Prdct{i}{\IdMapOn{X}}} &
    \Prdct{X}{X} \ar[r]_-{\mu} &
    X
    }
  \end{equation*}
  Any of these internal structures is \Defn{commutative} if the diagram below commutes.
  \begin{equation*}
    \xymatrix@!@C=-0.4em@R=0.8em{
    \Prdct{X}{X} \ar[rrrr]^-{\mu} \ar[rrd]_{\tau} &&&&
    X \\
    && \Prdct{X}{X} \ar[rru]_{\mu}
    }
  \end{equation*}
  Here $\tau\DefEq (\PrjctnOnto{2},\PrjctnOnto{1})$ is the \Defn{twist map}. %
  \index{twist map}%
\end{definition}

\begin{example}[Magma / monoid / group]
  \label{exa:Magma/Monoid/Group-From-InternalMMG}
  An internal group in the category $\Sets$ of sets is what we normally call a group; similarly for the structures of (commutative) magma or monoid. %
  \index{group}%
  \index{monoid}%
  \index{magma}%
  \NoProof
\end{example}

In (\ref{def:InternalMagma/Monoid/Group}), if the category $\Ctgry{X}$ is pointed, then $e=\ZeroObject$ is the unique unit in every object $X$. The commutativity of the diagram involving the map $i$ may be conceptualized by saying that $i$ is a two-sided inverse with respect to $\mu$. The following lemma asserts that, if an internal monoid admits a two-sided inverse operation, then this inverse operation is unique.

\begin{lemma}[Inverse operation on monoid is unique]
  \label{thm:InternalMonoid-InverseUnique}
  In the setting of (\ref{def:InternalMagma/Monoid/Group}), if $i$, $i'\from X\to X$ are two-sided inverse operations turning the monoid $(X,\mu,e)$ into an internal group, then $i=i'$.
\end{lemma}

Whenever an object $X$ in a category $\Ctgry{X}$ carries one of the internal structures defined in (\ref{def:InternalMagma/Monoid/Group}), then $\HomIn{\Ctgry{X}}{-}{X}$ is canonically a functor from $\Ctgry{X}^{\op}$ to the category of those structures, and vice versa:

\begin{theorem}[Internal magma yields magma Hom-sets]
  \label{thm:InternalMagma->Hom-Magma}
  In a category $\Ctgry{X} $, an internal unitary magma structure $(X,\mu,e)$ on an object $X$ determines a factorization $M$ of $\HomIn{\Ctgry{X}}{-}{X}$ through the category $\Magmas$ of unitary magmas, as in the diagram
  \[
    \xymatrix@R=2em@C=4em{&\Magmas \ar[d]^-U \\
    \Ctgry{X}^{\op} \ar@{-->}[ru]^-{M} \ar[r]_-{\HomIn{\Ctgry{X}}{-}{X}} &\SetsBsd}
  \]
  Furthermore, %
  \index[not]{m!$\Magmas$\IndSep category of unitary magmas}
  \begin{enumerate}[(i)]
    \item If $\mu$ is associative, respectively commutative, then so is the unitary magma $M(A)$, for every $A$ in $\Ctgry{X}$.
    \item If $(X,\mu,e,i)$ is an internal group, then $M\from \Ctgry{X} ^{\op}\to \Grps$.
    \item If $(X,\mu,e,i)$ is an internal commutative group, then $M\from \Ctgry{X} ^{\op}\to \AbGrps$.
  \end{enumerate}
\end{theorem}
\begin{proof}
  First, we note that every $A$ in $\Ctgry{X}$, the set $H_{A}\DefEq \HomIn{\Ctgry{X}}{A}{X}$ has a distinguished element $e_{A}$, namely the composite $A\to t\XRA{e} X$. Moreover, fore every $u\from A\to B$ in $\Ctgry{X}$, we have the commuting diagram
  \begin{equation*}
    \xymatrix@R=5ex@C=4em{
    A \ar[d]_{u} \ar[r] &
    t \ar[r]^-{e} &
    X \\
    B \ar[ru]
    }
  \end{equation*}
  This shows that $u^{\ast}(e_{B})=e_{A}$ so that the contravariant functor $\HomIn{\Ctgry{X}}{-}{X}$ actually takes values in the category $\SetsBsd$ of pointed sets. A unitary magma structure on $M(A)\coloneq H_{A}$ is given by $\mu_{A}\from \Prdct{H_{A}}{H_{A}}\to H_{A}$, the function which sends $f$, $g\from A\to X$ to the composite
  \begin{equation*}
    \xymatrix@R=5ex@C=4em{
    A \ar[r]^-{\mu_{A}(f,g)} \ar[d]_{\DgnlOn{A}} &
    X \\
    \Prdct{A}{A} \ar[r]_-{\Prdct{f}{g}} &
    \Prdct{X}{X} \ar[u]_{\mu}
    }
  \end{equation*}
  To see that $M$ is a functor, we observe that a map $u\from B\to A$ renders the diagram below commutative.
  \begin{equation*}
    \xymatrix@R=5ex@C=4em{
    B \ar[r]^-{u} \ar[d]_{\DgnlOn{B}}&
    A \ar[r]^-{\mu_{A}(f,g)} \ar[d]_{\DgnlOn{A}} &
    X \\
    \Prdct{B}{B} \ar[r]_-{\Prdct{u}{u}} &
    \Prdct{A}{A} \ar[r]_-{\Prdct{f}{g}} &
    \Prdct{X}{X} \ar[u]_{\mu}
    }
  \end{equation*}
  Defining $u^{\ast}\from \HomIn{\Ctgry{X}}{A}{X}\to \HomIn{\Ctgry{X}}{B}{X}$, by $u^{\ast}(f)\DefEq fu$, we conclude that $u^{\ast}\mu_{A}(f,g) = \mu_B(u^{\ast}f,u^{\ast}g)$. Via similar reasoning, associativity and commutativity of $\mu_{A}$ follow from the corresponding properties of $\mu$. If $i$ is the inverse operation of the internal group $X$, then $i_{\ast}\from H_{A}\to H_{A}$ is the inverse operation on the associative monoid $H_{A}$. (iii) follows by similar reasoning.
\end{proof}

\begin{theorem}[Magma Hom-sets yield internal magma]
  \label{thm:MagmaHom-sets->InternalMagma}
  In a category $\Ctgry{X}$ with finite products, a factorization of $H\DefEq\HomIn{\Ctgry{X} }{-}{X}$ as a composite
  \[
    \xymatrix@R=2em@C=4em{&\Magmas \ar[d]^-U \\
    \Ctgry{X}^{\op} \ar@{-->}[ru]^-{M} \ar[r]_-{\HomIn{\Ctgry{X}}{-}{X}} &\SetsBsd}
  \]
  determines uniquely the structure of an internal unitary magma $(X,\mu,e)$ on $X$ with the property that
  \begin{equation*}
    \mu_{A}(f,g) = \mu\Comp (\Prdct{f}{g}) \Comp \DgnlOn{A},\qquad \text{for all}\quad f,g\from A\to X.
  \end{equation*}
  Here $(H_{A},\mu_{A})$ is the internal magma structure on $H_{A}\coloneq H(A)=\HomIn{\Ctgry{X}}{A}{X}$.
  If the codomain of $M$ is the category of (commutative) monoids, then $\mu$ is (commutative) associative. If the codomain of $M$ is $\Grps$, then there is a unique $i\from X\to X$ such that $(X,\mu,e,i)$ is an internal group in $\Ctgry{X}$. If the codomain of $M$ is $\AbGrps$, then $(X,\mu,e,i)$ is an internal abelian group.
\end{theorem}
\begin{proof}
  To construct the internal multiplication on $X$, let $\PrjctnOnto{1}$, $\PrjctnOnto{2}\from \Prdct{X}{X}\to X$ be the projections onto the first and second factor. Then $\mu\DefEq \mu_{X\times X}(\PrjctnOnto{1},\PrjctnOnto{2})\from \Prdct{X}{X}\to X$ is an internal magma structure on $X$. Then the magma operation $\mu_{A}$ on $H_{A}$ is given by $\mu$ because any two morphisms $f$, $g\from A\to X$ determine $\PrdctMapInto{f,g}\from A\to \Prdct{X}{X}$ for which the functor property of $M$ renders this diagram commutative:
  \begin{equation*}
    \xymatrix@R=5ex@C=4em{
    \Prdct{M_{X\times X}}{M_{X\times X}} \ar[r]^-{\mu_{X\times X}} \ar[d]_{(f,g)^{\ast}\prdct (f,g)^{\ast}} &
    M_{X\times X} \ar[d]^{(f,g)^{\ast}} &
    (\PrjctnOnto{1},\PrjctnOnto{2}) \ar@{|->}[r] \ar@{|->}[d] &
    \mu \ar@{|->}[d] \\
    \Prdct{M_A}{M_A} \ar[r]_-{\mu_A} &
    M_A &
    (f,g) \ar@{|->}[r] &
    *\txt{$\mu\Comp (f,g)$ \\ $\mu_A(f,g)$}
    }
  \end{equation*}
  If $M$ takes values in the category of unitary magmas, then each $M_{A}$ contains a designated element $e_{A}\from A\to X$ such that
  \begin{equation*}
    \mu_{A}(f,e_{A}) = f = \mu_{A}(e_{A},f)
  \end{equation*}
  If $\tau$ is a terminal object in $\Ctgry{X}$, let $t_{A}\from A\to \tau$ be the unique map. Then functoriality of $M$ yields:
  \begin{equation*}
    e_{A} = t_{A}^{\ast}(e_{\tau}) = e_{\tau}\Comp t_{A},\qquad A \XRA{t_{A}} \tau \XRA{e_{\tau}} X
  \end{equation*}
  Then $(X,\mu,e)$ is an internal unitary magma structure on $X$ because neutrality of $e$ on the right follows from this computation.
  \begin{equation*}
    \mu\Comp (\Prdct{\IdMapOn{X}}{\hat{e}})\Comp \DgnlOn{X} = \mu\Comp \PrdctMapInto{\IdMapOn{X},\hat{e}} = \mu_{X}\PrdctMapInto{\IdMapOn{X},\hat{e}} = \IdMapOn{X}
  \end{equation*}
  Neutrality of $e$ on the left follows similarly.

  To see that $\mu$ is unique, suppose $\mu'\from \Prdct{X}{X}\to X$ is another map with
  \begin{equation*}
    \mu_A(f,g)= \mu'\Comp (\Prdct{f}{g})\Comp \DgnlOn{A}
  \end{equation*}
  for all $f$, $g\from A\to X$. Then, using that $(\PrjctnOnto{1}\prdct \PrjctnOnto{2})\Comp\DgnlOn{X\times X} =\IdMapOn{X\times X}$, we see:
  \begin{equation*}
    \mu = \mu\Comp (\PrjctnOnto{1}\prdct \PrjctnOnto{2})\Comp\DgnlOn{X\times X} = \mu_{X\times X}(\PrjctnOnto{1}, \PrjctnOnto{2})= \mu'\Comp (\PrjctnOnto{1}\prdct \PrjctnOnto{2})\Comp \DgnlOn{X\times X} = \mu'.
  \end{equation*}
  Next, if $M\from \Ctgry{X} ^{\op}\to \Monoids$, we need to show that $\mu$ is associative. We know that the diagram below commutes. So, the projections $\PrjctnOnto{1}$, $\PrjctnOnto{2}$, $\PrjctnOnto{3}\from X^3\to X$ are mapped as shown on the right.
    {\small
      \begin{equation*}
        \xymatrix@R=5ex@C=3.5em{
        M_{X^3}\prdct M_{X^3}\prdct M_{X^3} \ar[r]^-{\IdMap\times \mu_{X^3}} \ar[d]_{\mu_{X^3}\times \IdMap} &
        \Prdct{M_{X^3}}{M_{X^3}} \ar[d]^{\mu_{X^3}} &
        (\PrjctnOnto{1},\PrjctnOnto{2},\PrjctnOnto{3}) \ar@{|->}[r] \ar@{|->}[d] &
        (\PrjctnOnto{1},\mu\Comp (\PrjctnOnto{2},\PrjctnOnto{3})) \ar@{|->}[d] \\
        \Prdct{M_{X^3}}{M_{X^3}} \ar[r]_-{\mu_{X^3}} &
        M_{X^3} &
        (\mu\Comp (\PrjctnOnto{1},\PrjctnOnto{2}),\PrjctnOnto{3}) \ar@{|->}[r] &
        *\txt{$\mu\Comp (\PrjctnOnto{1},\mu\Comp (\PrjctnOnto{2},\PrjctnOnto{3}))$ \\ $\mu\Comp (\mu\Comp (\PrjctnOnto{1},\PrjctnOnto{2}), \PrjctnOnto{3})$ }
        }
      \end{equation*}}%
  Spelling out in detail how the two composites on the bottom right are formed, we arrive at this diagram:
  \begin{equation*}
    \xymatrix@R=5ex@C=4.5em{
    X^3 \ar[r]^-{\DgnlOn{X^3}} &
    \Prdct{X^3}{X^3}\prdct {X^3} \ar[r]^-{\PrjctnOnto{1}\times\PrjctnOnto{2}\times \PrjctnOnto{3}} &
    \Prdct{X}{X}\prdct {X} \ar[r]^-{\IdMapOn{X}\times \mu} \ar[d]_{\mu\times \IdMapOn{X}} &
    \Prdct{X}{X} \ar[d]^{\mu} \\
    && \Prdct{X}{X} \ar[r]_-{\mu} &
    X
    }
  \end{equation*}
  The clockwise composite $X^3\to X$ is $\mu\Comp (\PrjctnOnto{1},\mu\Comp (\PrjctnOnto{2},\PrjctnOnto{3}))$, and the counterclockwise composite is $\mu\Comp (\mu\Comp (\PrjctnOnto{1},\PrjctnOnto{2}), \PrjctnOnto{3})$. Now, notice that $(\PrjctnOnto{1}\times\PrjctnOnto{2}\times \PrjctnOnto{3})\Comp \DgnlOn{X^3}=\IdMapOn{X^3}$, and the associativity of $\mu$ follows.

  Finally, if $M\from \Ctgry{X} ^{\op}\to \Grps$, we just established the structure $(X,\mu,e)$ of an internal monoid. For it to become an internal group structure, we need an internal inverse operation $i$ on $X$. Let $I\from M\Rightarrow M$ be the natural transformation, unique by (\ref{thm:InternalMonoid-InverseUnique}), such that $(M_A,\mu_A,E_{A},I_A)$ is a group object structure on $M_{A}=\HomIn{\Ctgry{X} }{A}{X}$. We claim that $i\DefEq I_X(\IdMapOn{X})$ yields an internal group $(X,\mu,e,i)$. To see this, note that this composite is $\hat{e}$:
  \begin{equation*}
    \xymatrix@R=1ex@C=4em{
    M_X \ar[r]^-{\DgnlOn{X}} &
    \Prdct{M_X}{M_X} \ar[r]^-{\IdMapOn{X}\times I_{X}} &
    \Prdct{M_X}{M_X} \ar[r]^-{\mu_X} &
    M_X \\
    \IdMapOn{X} \ar@{|->}[r] &
    (\IdMapOn{X},\IdMapOn{X}) \ar@{|->}[r] &
    (\IdMapOn{X},i) \ar@{|->}[r] &
    \mu_X(\IdMapOn{X},i) = \hat{e}
    }
  \end{equation*}
  Since $\mu_X(\IdMapOn{X},i)$ is given by the composite $X \XRA{\DgnlOn{X}} \Prdct{X}{X} \XRA{\IdMapOn{X}\times i} \Prdct{X}{X} \XRA{\mu} X$, we see that $i$ is a right inverse. A similar argument shows that $i$ acts as left inverse, and the proof is complete.
\end{proof}

\begin{corollary}[$i^2=i\Comp i=\IdMap$ for inverse operation of internal group]
  \label{thm:InternalGroup:i^2=i}
  If $(X,\mu,e,i)$ is an internal group, then $i^2 =\IdMapOn{X}$. %
  \index{internal group!inverse operation $i^2=\IdMap$}
\end{corollary}
\begin{proof}
  From Theorem \ref{thm:InternalMagma->Hom-Magma} and the proof of part (ii), we know that $G\DefEq \Hom{X}{X}$ is a group. The left inverse of $f\in G$ is given by $i\Comp f$. It follows that, for all $g\in G$,
  \begin{equation*}
    \mu\Comp (f,i\Comp f) = e = \mu\Comp (i\Comp f,f),\quad \text{hence}\quad \mu\Comp (i^2\Comp f,i\Comp f) = e = \mu\Comp (f,i\Comp f)
  \end{equation*}
  This means that, in the group $G$,
  \begin{equation*}
    f\cdot (i\Comp f) = (i^2\Comp f)\cdot (i\Comp f)\quad \text{and so}\quad f = i^2\Comp f
  \end{equation*}
  This means that $i^2$ is a left identity in the monoid $(\Hom{X}{X},\circ)$, and so $i^2=\IdMapOn{X}$.
\end{proof}

\begin{definition}[Morphism of internal magmas/monoids/groups]
  \label{def:InternalMagmas/Monoids/Groups-Morphism}
  In a category $\Ctgry{X}$ with finite products, let $S$ denote any of the internal structures (commutative) magma, monoid, group. A \Defn{morphism of $S$-objects} $X$ and $Y$ is given by a morphism $f\from X\to Y$ which commutes with the $S$-structures diagrams in Definition \ref{def:InternalMagma/Monoid/Group}. %
  \index{morphism!of internal magma/monoid/group}
\end{definition}

A morphism $f\from X\to Y$ of $S$-objects induces, for every object $U$ in $\Ctgry{X} $ a morphism of $S$-structures on Hom-sets.
\begin{equation*}
  f_{\ast}\from \HomIn{\Ctgry{X} }{U}{X}\longrightarrow \HomIn{\Ctgry{X} }{U}{Y}
\end{equation*}
Remarkable is that $f$ is a morphism of any of the structures (commutative) unitary magma, monoid, group if and only if $f$ is a morphism of the underlying unitary magmas.

\begin{subordinate}[Comment]{on the approach taken in this section}
  In our discussion of internal objects, we used elementary methods only. However, we also  point out that a treatment which is both more efficient and more conceptual is possible using Yoneda embeddings---see \cite{FBorceuxDBourn2004}, for instance.

\end{subordinate}

\begin{exercises}

\begin{exercise}[Recognizing a morphism of internal magmas/monoids/groups]
  \label{exe:InternalMagmas/Monoids/Groups-MorphismRecognize}
  In the setting of Definition \ref{def:InternalMagmas/Monoids/Groups-Morphism}, show the following about an $\Ctgry{X} $-morphism $f\from X\to Y$:
  \begin{enumerate}[(i)]
    \item For every object $U$ in $\Ctgry{X} $, $f$ induces a morphism of $S$-structures $f_{\ast}\from \HomIn{\Ctgry{X} }{U}{X}\to \HomIn{\Ctgry{X}}{U}{Y}$.
    \item $f$ is a morphism of underlying $S$-structures if and only if $f$ is a morphism of underlying unitary magmas.
  \end{enumerate}
\end{exercise}
\end{exercises}
\section{Biproducts and Linear Categories}
\label{sec:BiProducts}

Here we provide categorical foundations for a discussion of internal abelian group objects %
and abelian categories. 
Thus we are working in a category $\Ctgry{X}$ which is enriched over pointed sets; see Section \ref{sec:PointedSets}.

Already the environment of a pointed category supports a discussion of biproducts familiar from abelian categories: If objects $X$ and $Y$ admit both a coproduct and a product, then there is a canonical comparison map $\SumProdComp{X}{Y}\from \CoPrdct{X}{Y}\to \Prdct{X}{Y}$. Whenever $\SumProdComp{X}{Y}$ is an isomorphism, we can `compress' it into a single object via (\ref{thm:BiProduct-Existence}), called a biproduct of $X$ and $Y$, denoted $\BiPrdct{X}{Y}$.

Quite useful is that a zero object in $\Ctgry{X}$ can be recognized by the property that it is neutral with respect to the formation of biproducts: An object $I$ is a zero object if and only if, for every $A$ in~$\Ctgry{X}$, $\BiPrdct{A}{I}=A=\BiPrdct{I}{A}$.

A bit deeper is that, whenever an object $X$ in a pointed category $\Ctgry{X}$ admits finite biproducts, then the map $\mu=\SumMapOutOf{\IdMapOn{X},\IdMapOn{X}}\from \BiPrdct{X}{X}\to X$ turns $X$ into an internal commutative monoid object in $\Ctgry{X}$ with `neutral element' given by the morphism $\ZeroMap\from \ZeroObject\to X$. Moreover, $\mu$ is the one and only internal monoid map on $X$; see (\ref{thm:CommutativeInternalMonoid-From-BiProduct}). While this kind of result is familiar from additive categories, we see here that it holds in much greater generality.

Proposition \ref{thm:CommutativeInternalMonoid-From-BiProduct} motivates the concept of `linear category' (\ref{def:LinearCategory}), that is a pointed category in which the biproduct $\BiPrdct{X}{Y}$ of any two objects exists. Consequently, every object in a linear category $\Ctgry{X} $ carries a unique structure of an internal commutative monoid. So, $\Ctgry{X}$ has Hom-sets which are commutative monoids. Moreover, this commutative monoid structure on Hom-sets responds biadditively to composition of morphisms, and so $\Ctgry{X} $ is enriched in the category $\CMon$ of commutative monoids.

Since categories which are $\CMon$\*-enriched are $\SetsBsd$\*-enriched, the concept of a biproduct makes sense. Remarkable is that the commutative monoid structure on Hom-sets enables us to identify a biproduct structure diagram without referring to either of its universal properties; see (\ref{thm:BiProduct-By-CMonEnrichment}). As an immediate application we show that, whenever a category $\SACtgry{X}$ is enriched over $\CMon$, then any product cone and any coproduct cocone admits a unique expansion to a biproduct; see (\ref{thm:BiProduct-From-(Co)Product}).

This implies (\ref{thm:Linear<->EnrichedInCommutativeMonoids}) that a category with finite sums or finite products is linear if and only if it is enriched in the category $\CMon$ of commutative monoids. Consider then a finitely complete category and form the full subcategory $\CMon(\Ctgry{X})$ of internal commutative monoids. As the product of any two commutative monoids is again one such, $\CMon(\Ctgry{X})$ has finite products, hence is linear; see (\ref{thm:InternalCommMonoidsLinear}).

To get started, let $\Ctgry{X}$ be a category, enriched over pointed sets. Assume that, for a pair $X$ and $Y$ of objects in $\Ctgry{X}$, both their coproduct $(X+Y,\InclsnOf{X},\InclsnOf{Y})$ and their product $(\Prdct{X}{Y},\PrjctnOnto{X},\PrjctnOnto{Y})$ exist. Using the zero-map of $\Ctgry{X}$ enables us to construct the map
\begin{equation*}
  \SumProdComp{X}{Y}\DefEq\PrdctMapInto{\SumMapOutOf{\IdMapOn{X},\ZeroMap},\SumMapOutOf{\ZeroMap,\IdMapOn{Y}}} = \SumMapOutOf{\PrdctMapInto{\IdMapOn{X},\ZeroMap}, \PrdctMapInto{\ZeroMap,\IdMapOn{Y}} } \from X+Y \to \Prdct{X}{Y}.
\end{equation*}
We refer to $\SumProdComp{X}{Y}$ as the \Defn{coproduct/product comparison map} of $X$ and $Y$.

\begin{proposition}[Properties of the coproduct/product comparison map]
  \label{thm:CPMap-Properties}
  In a category enriched over pointed sets, the coproduct/product comparison map of $X$ and $Y$ fits into the commutative diagram
  
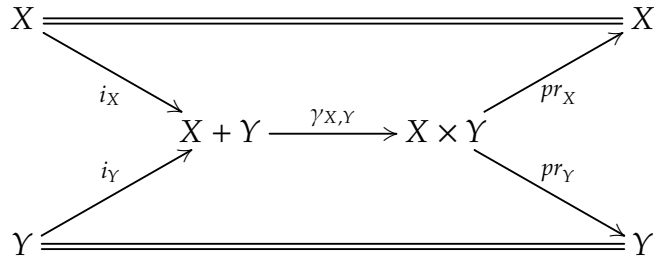
\begin{equation*}
    \xymatrix@R=5ex@C=4em{
    X \ar[rd]_{\InclsnOf{X}} \ar@{=}[rrr] &&&
    X \\
    & X+Y \ar[r]^-{\SumProdComp{X}{Y}} &
    \Prdct{X}{Y} \ar[ru]_{\PrjctnOnto{X}} \ar[rd]^{\PrjctnOnto{Y}} \\
    Y \ar[ru]^{\InclsnOf{Y}} \ar@{=}[rrr] &&&
    Y
    }
  \end{equation*}
  with the properties $\PrjctnOnto{Y}\SumProdComp{X}{Y} i_X = 0_{YX}$ and $\PrjctnOnto{X}\SumProdComp{X}{Y}i_Y=0_{XY}$.
\end{proposition}
\begin{proof}
  The see that the top part of the diagram commutes, we compute:
  \begin{equation*}
    \PrjctnOnto{X}\SumProdComp{X}{Y}\InclsnOf{X} = \PrjctnOnto{X} \SumMapOutOf{\PrdctMapInto{\IdMapOn{X},0},\PrdctMapInto{0,\IdMapOn{Y} } } \InclsnOf{X} = \PrjctnOnto{X} \PrdctMapInto{\IdMapOn{X},0} = \IdMapOn{X}
  \end{equation*}
  That the bottom part of the diagram commutes follows via a similar argument. Now
  \begin{equation*}
    \PrjctnOnto{Y}\SumProdComp{X}{Y}\InclsnOf{X} = \PrjctnOnto{Y}(\SumMapOutOf{\IdMapOn{X},0_{YX}}, \SumMapOutOf{0_{XY},\IdMapOn{Y}}) \InclsnOf{X} = \SumMapOutOf{\IdMapOn{X},0_{YX}}\InclsnOf{X} = 0_{YX}\InclsnOf{X} = 0_{YX}.
  \end{equation*}
  A similar computation shows that $\PrjctnOnto{X}\SumProdComp{X}{Y}\InclsnOf{Y}=0_{XY}$.
\end{proof}

If $\SumProdComp{X}{Y}$ is an epimorphism, the $\PrjctnOnto{Y}=\CoKerMap{\PrdctMapInto{\IdMapOn{X},\ZeroMap}}$; see Exercise \ref{exe:PoductProjectionNormalEpi}. Next, we single out the special case where the coproduct/product comparison map is the identity map.

\begin{definition}[Biproduct]
  \label{def:BiProduct}
  In a $\SetsBsd$-enriched category $\Ctgry{X}$ we say that objects $X$ and $Y$ \Defn{admit a biproduct} if and only if there exists an object $B$ in $\Ctgry{X}$, together with a commutative diagram of maps
  \begin{equation*}
    \vcenter
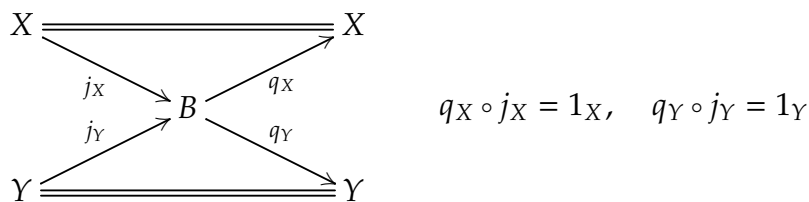
{\xymatrix@R=3ex@C=4em{
    X \ar@{=}[rr] \ar[rd]_{j_X} &&
    X \\
    & B \ar[ru]_{q_X} \ar[rd]^{q_Y} \\
    Y \ar[ru]^{j_Y} \ar@{=}[rr] &&
    Y
    }}
    \qquad q_X\Comp j_X=1_X, \quad q_Y\Comp j_Y=1_Y
  \end{equation*}
  such that the following conditions hold:
  \begin{enumerate}[(i)]
    \item $(B,j_X,j_Y)$ is a coproduct of $X$ and $Y$;
    \item $(B,q_X,q_Y)$ is a product of $X$ and $Y$;
    \item $q_Yj_X=\ZeroMap_{YX}$, and $q_Xj_Y=\ZeroMap_{XY}$.
  \end{enumerate}
  In this situation, we write $B= \BiPrdct{X}{Y}$. %
  \index[not]{b!$\BiPrdct{X}{Y}$\IndSep biproduct of $X$ and $Y$}%
  \index{biproduct!in pointed category}%
\end{definition}

\begin{remark}[Biproduct inclusions and projections are canonical]
  \label{rem:Notation Inclusion Projection Biproduct}%
  In Definition~\ref{def:BiProduct}, the interpretation of $B$ as a product of $X$ and $Y$ allows us to write $j_X=(1_X,0)$ and $j_Y=(0,1_Y)$. Dually, $q_X=\langle 1_X,0\rangle$ and $q_Y=\langle 0,1_Y\rangle$.
\end{remark}

\begin{proposition}[Existence of biproducts]
  \label{thm:BiProduct-Existence}
  In a $\SetsBsd$-enriched category $\Ctgry{X}$, objects $X$ and $Y$ admit a biproduct if and only if both the sum and the product of $X$ and $Y$ exist, and the coproduct/product comparison map for $X$ and $Y$ is an isomorphism.
\end{proposition}
\begin{proof}
  If $X$ and $Y$ admit a biproduct $\BiPrdct{X}{Y}$, then $X+Y\cong \BiPrdct{X}{Y}\cong \Prdct{X}{Y}$ by universality. This implies that the coproduct/product comparison map for $X$ and $Y$ is an isomorphism. Conversely, if $\SumProdComp{X}{Y}$ is an isomorphism, set
  \begin{equation*}
    B\DefEq X+Y,\quad j_X\DefEq i_X,\quad j_Y\DefEq i_Y,\quad q_X\DefEq \PrjctnOnto{X}\SumProdComp{X}{Y},\quad q_Y\DefEq \PrjctnOnto{Y}\SumProdComp{X}{Y}
  \end{equation*}
  Then properties (i) and (iii) are satisfied by design. It remains to show that $(B,q_X,q_Y)$ is a product of $X$ and $Y$. If $f\from A\to X$ and $g\from A\to Y$ be given, then there exists $\PrdctMapInto{f,g}\from A\to \Prdct{X}{Y}$ unique with
  \begin{equation*}
    \PrjctnOnto{X}\PrdctMapInto{f,g} = f \qquad \text{and}\qquad q_Y\PrdctMapInto{f,g} = g.
  \end{equation*}
  Put $\phi\DefEq \SumProdComp{X}{Y}^{-1}\PrdctMapInto{f,g}\from A\to B$. Then
  \begin{equation*}
    q_X\phi = \PrjctnOnto{X} \PrdctMapInto{f,g} = f \qquad \text{and}\qquad q_Y\phi = \PrjctnOnto{Y} \PrdctMapInto{f,g} = g.
  \end{equation*}
  Moreover, if $\psi\from A\to B$ is arbitrary with $q_X\psi =f$ and $q_Y\psi = g$, then $\SumProdComp{X}{Y}\psi = (f,g)$, and so $\phi=\psi$.
\end{proof}

\begin{corollary}[Zero object via biproduct, I]
  \label{thm:Char-ZeroObject-In-PointedCat}
  In a category enriched over pointed sets, an object $I$ is a zero object if and only if for each object $A$, the diagram
  \begin{equation*}
    \vcenter
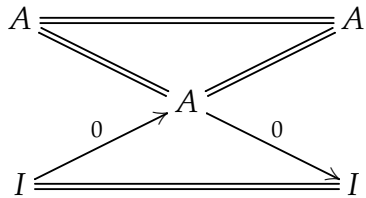
{\xymatrix@R=3ex@C=4em{
    A \ar@{=}[rr] \ar@{=}[rd] &&
    A \\
    & A \ar@{=}[ru] \ar[rd]^{\ZeroMap} \\
    I \ar[ru]^{\ZeroMap} \ar@{=}[rr] &&
    I
    }}
  \end{equation*}
  commutes and represents $A$ as a biproduct of $A$ and $I$.
\end{corollary}
\begin{proof}
  If $I$ is a zero object and $A$ is arbitrary then the diagram commutes; we encourage the reader to verify the biproduct recognition conditions of Proposition~\ref{thm:BiProduct-Existence}.

  For the converse, we prove that $I$ is terminal. Then $I$ is initial by duality. For a given object $A$, we show that $\ZeroMap\from A\to I$ is the only map. Indeed, for an arbitrary $a\from A\to I$, the product property of the right hand part of the diagram tells us that the diagram of solid arrows below is rendered commutative by the filler $c$.
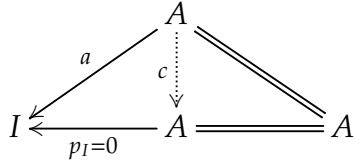
  \begin{equation*}
    \xymatrix@R=5ex@C=4em{
    & A \ar[ld]_{a} \ar@{=}[rd] \ar@{.>}[d]_{c}\\
    I &
    A \ar[l]^-{p_I=0} \ar@{=}[r] &
    A
    }
  \end{equation*}
  It follows that $c=\IdMapOn{A}$ and, hence, $a = p_I=\ZeroMap$. - This completes the proof.
\end{proof}

\begin{corollary}[Zero object via biproduct, II]
  In a category enriched over pointed sets, an object $I$ is a zero object if and only if, for each object $A$, every one of the maps
  \[
    \xymatrix@C=4em{
    A\oplus I \ar@<.5ex>[r]^-{\langle 1_A,0\rangle} &
    \ar@<.5ex>[l]^-{(1_A,0)} A \ar@<.5ex>[r]^-{(0,1_A)} &
    I\oplus A \ar@<.5ex>[l]^-{\langle 0,1_A\rangle}
    }
  \]
  is an isomorphism. \NoProof
\end{corollary}

\begin{remark}[Finite Biproducts]
  \label{rem:FiniteBiproducts}
  If a category admits binary biproducts, then it admits biproducts for any finite collection of objects. The details are a bit intricate, and we organize an approach in Exercise \ref{exe:Binary(Co)Products->Finite(Co)Products}.
\end{remark}

Let us now turn to a fundamental relationship between biproducts in a pointed category and internal monoids.

\begin{proposition}[Commutative internal monoid from biproduct]
  \label{thm:CommutativeInternalMonoid-From-BiProduct}
  In a pointed category $\Ctgry{X}$, suppose an object $X$ admits a biproduct $X\oplus X$. Then
  \[
    \mu\DefEq \SumMapOutOf{\IdMapOn{X},\IdMapOn{X}}\from \BiPrdct{X}{X}\to X
  \]
  forms an internal commutative monoid $(X,\mu)$ with unit $e \from \ZeroObject\to X$. Moreover, $\mu$ is the only internal unitary magma structure on $X$.
\end{proposition}
\begin{proof}
  Via the sum property, diagram (U) in (\ref{def:InternalMagma/Monoid/Group}) shows that $\mu$ is unitary and unique, and that diagram (C) commutes. Note that here we make implicit use of Remark~\ref{rem:Notation Inclusion Projection Biproduct}. - This was to be shown.
\end{proof}

\begin{definition}[Linear category]
  \label{def:LinearCategory}
  A pointed category $\SACtgry{X}$ is called \Defn{linear} or \Defn{half-additive} if for any pair of objects $X$ and $Y$ the sum and product exists and if the comparison map $\SumProdComp{X}{Y}\from X+Y\to \Prdct{X}{Y}$ is an isomorphism. %
  \index{linear category}\index{category!linear}
\end{definition}

In other words, a linear category is a category \emph{with finite biproducts}: it is pointed and for any finite collection of objects $X_1$, \dots, $X_n$, $n\geq 1$, their biproduct $X_1\oplus\cdots \oplus X_n$ exists; see (\ref{rem:FiniteBiproducts}).

\begin{proposition}[Linear category and commutative monoids]\label{thm:LinearThenCMon}
  If $\Ctgry{X}$ is a linear category, then we have an equivalence $\CMon(\Ctgry{X})=\Ctgry{X}$.
\end{proposition}
\begin{proof}
  If $X$ is an object of $\Ctgry{X}$, then Proposition~\ref{thm:CommutativeInternalMonoid-From-BiProduct} tells us that the fold map $\SumMapOutOf{\IdMapOn{X},\IdMapOn{X}}$ on $X$ is an internal commutative monoid on $X$, and it is the only such structure. So, it remains to check that each morphism $f\from X\to Y$ in $\Ctgry{X}$ is a morphism of internal monoid structures. To see this, verify that the identity $\langle1_Y,1_Y\rangle \Comp(f+f)=f\Comp\langle1_X,1_X\rangle$ holds on the summands of $X+X$.
\end{proof}

With Corollary~\ref{thm:InternalCommMonoidsLinear} below we provide a converse to this statement: each category of internal monoids in a category with finite products is linear.

Combined, Propositions \ref{thm:InternalMagma->Hom-Magma} and \ref{thm:CommutativeInternalMonoid-From-BiProduct} show that all objects in a linear category $\Ctgry{X} $ carry a unique internal commutative monoid structure and, hence, that each Hom-set $\HomIn{\Ctgry{X} }{A}{X}$ carries a commutative monoid structure which responds biadditively to composition. In other words, $\Ctgry{X}$ is a category \emph{enriched over commutative monoids}.

We are now going to prove the converse, and obtain Proposition~\ref{thm:Linear<->EnrichedInCommutativeMonoids} which says that a category with finite sums or finite products is linear if and only if it is enriched over the category $\CMon$ of commutative monoids.

\begin{lemma}[Recognizing a biproduct via $\CMon$-enrichment]
  \label{thm:BiProduct-By-CMonEnrichment}
  In a category enriched over commutative monoids, a commutative diagram
  \begin{equation}\label{BiproductDiagram}
    \vcenter{
    \xymatrix@R=3ex@C=4em{
    X \ar@{=}[rr] \ar[rd]_{i_X} &&
    X \\
    & B \ar[ru]_{p_X} \ar[rd]^{p_Y} \\
    Y \ar[ru]^{i_Y} \ar@{=}[rr] &&
    Y
    }}\qquad p_X\Comp i_X=1_X, \quad p_Y\Comp i_Y=1_Y
  \end{equation}
  describes $B$ as a biproduct of $X$ and $Y$ if and only if %
  \index{biproduct!in linear category} %
  \begin{equation*}
    p_Yi_X=\ZeroMap_{YX} \qquad \text{and}\qquad p_Xi_Y=\ZeroMap_{XY} \qquad \text{and}\qquad i_Xp_X + i_Yp_Y=\IdMapOn{B}.
  \end{equation*}
\end{lemma}
\begin{proof}
  Suppose the diagram describes $B$ as a biproduct. Then
  \begin{equation*}
    p_X(i_Xp_X + i_Yp_Y)=p_X  \qquad \text{and}\qquad p_Y(i_Xp_X + i_Yp_Y)=p_Y
  \end{equation*}
  so that $i_Xp_X + i_Yp_Y=\IdMapOn{B}$ by the universal property of the product.

  Conversely, if the morphisms in the diagram satisfy the stated identities, then $i_X$ and $i_Y$ are structure maps for the sum of $X$ and $Y$. To see this, suppose $f_X\from X\to Z$ and $f_Y\from Y\to Z$ are given. Then $\varphi\DefEq f_Xp_X + f_Yp_Y$ satisfies $\varphi i_X = f_X$ and $\varphi i_Y = f_Y$. To see that the map $\varphi$ is unique with this property, suppose $\psi\from B\to Z$ satisfies $\psi i_X=f_X$ and $\psi i_Y = f_Y$. Then
  \begin{equation*}
    \psi = \psi (i_Xp_X + i_Yp_Y) = f_Xp_X + f_Yp_Y = \varphi(i_Xp_X + i_Yp_Y) = \varphi
  \end{equation*}
  Dual reasoning shows that $p_X$ and $p_Y$ are structure maps for a product.
\end{proof}

\begin{proposition}[Biproduct from product or coproduct]
  \label{thm:BiProduct-From-(Co)Product}
  For any two objects $X$ and $Y$ in a category enriched over $\CMon$ the following hold:
  \begin{enumerate}[(i)]
    \item Any product cone $X \XLA{p_X} B \XRA{p_Y} Y$ admits an expansion to a biproduct diagram \eqref{BiproductDiagram}. Moreover, this expansion is unique.
    \item Any coproduct cocone $X \XRA{i_X} B \XLA{i_Y} Y$ admits an expansion to a biproduct diagram \eqref{BiproductDiagram}. Moreover, this expansion is unique.
  \end{enumerate}
\end{proposition}
\begin{proof}
  We prove (i), the verification of (ii) is dual. As candidates for coproduct structure maps, we propose $i_X=(1_X,0)$ and $i_Y=(0,1_Y)$. It then remains to check that the resulting diagram of shape \eqref{BiproductDiagram} meets the requirements of the biproduct recognition criterion (\ref{thm:BiProduct-By-CMonEnrichment}). By design, diagram \eqref{BiproductDiagram} commutes, and the diagonal composites vanish. Via the universal property of the product, the following computations imply that $i_Xp_X + i_Yp_Y = \IdMapOn{B}$:
  \begin{equation*}
    p_X(i_Xp_X + i_Yp_Y) = p_X + \ZeroMap = p_X\IdMapOn{X} \qquad \text{and}\qquad p_Y(i_Xp_X + i_Yp_Y) = \ZeroMap + p_Y = p_Y\IdMapOn{B}
  \end{equation*}
  Finally, to see that the maps $i_X$ and $i_Y$ are the only maps which complete $p_X$ and $p_Y$ to a biproduct, we compute
  \begin{equation*}
    j_X = \IdMapOn{X}j_X = (i_Xp_X + i_Yp_Y)j_X = i_X + \ZeroMap = i_X
  \end{equation*}
  A similar argument shows that $j_Y=i_Y$, and this completes the proof.
\end{proof}

\begin{corollary}[Initial / terminal object is zero-object]
  \label{thm:ZeroObject-In-LinearCat}
  If in a category enriched over commutative monoids an initial or a terminal object exists, then this category is pointed.
\end{corollary}
\begin{proof}
  We show that an initial object is also a terminal object. Dualizing will then complete the proof. So, consider the situation where $I$ is initial and $A$ is any object. Note first that the left hand side of the diagram below is the cocone of $A$ as a coproduct of $A$ and $I$.
  \begin{equation*}
    \vcenter{\xymatrix@R=3ex@C=4em{
    A \ar@{=}[rr] \ar@{=}[rd]_{i_A} &&
    A \\
    & A \ar@{.>}[ru]_{p_A} \ar@{.>}[rd]^{p_I} \\
    I \ar[ru]^{\ZeroMap=i_I} \ar@{=}[rr] &&
    I
    }}
  \end{equation*}
  By (\ref{thm:BiProduct-From-(Co)Product}) there exist $p_A\from A\to A$ and $p_I\from A\to I$ such that the right hand side of the diagram is the structure diagram of the product of $A$ by $I$. Clearly $p_A=1_A$. As $p_Ii_A=\ZeroMap$, we have $p_I=\ZeroMap$. The result now follows from Corollary~\ref{thm:Char-ZeroObject-In-PointedCat}.
\end{proof}

\begin{corollary}[Linear if and only if enriched in commutative monoids]
  \label{thm:Linear<->EnrichedInCommutativeMonoids}
  A category $\SACtgry{X}$ with either finite sums or finite products is linear if and only if it is enriched in the category of commutative monoids.
\end{corollary}
\begin{proof}
  If $\SACtgry{X}$ is linear then (\ref{thm:CommutativeInternalMonoid-From-BiProduct}) is the key to showing that $\SACtgry{X}$ is enriched over commutative monoids.

  Conversely, if $\SACtgry{X}$ is enriched over commutative monoids. Then the existence of either finite sums or finite products implies that $\SACtgry{X}$ is pointed (by Corollary~\ref{thm:ZeroObject-In-LinearCat}, since a nullary sum is an initial object, while a nullary product is a terminal object) and has biproducts (by Proposition~\ref{thm:BiProduct-From-(Co)Product}). So, $\Ctgry{X}$ is linear.
\end{proof}

\begin{corollary}[Internal commutative monoids form linear category]
  \label{thm:InternalCommMonoidsLinear}
  The category $\CMon(\Ctgry{X})$ of internal commutative monoids in any category with finite products $\Ctgry{X}$ is linear.
\end{corollary}
\begin{proof}
  Motivated by (\ref{thm:Linear<->EnrichedInCommutativeMonoids}), we show that $\CMon(\Ctgry{X})$ admits binary products via the following construction. For commutative monoids $(A,\mu)$ and $(B,\nu)$\*, their product $\Prdct{A}{B}$ in $\SACtgry{X}$ is again a commutative monoid, with the structure map
  \begin{equation*}
    \xymatrix@R=5ex@C=4em{
    (\Prdct{A}{B})\prdct (\Prdct{A}{B}) \ar[r]^-{\IdMapOn{A}\times \tau \times \IdMapOn{B}} &
    (\Prdct{A}{A})\prdct (\Prdct{B}{B}) \ar[r]^-{\mu\times \nu} &
    \Prdct{A}{B}
    }
  \end{equation*}
  Here $\tau=(\PrjctnOnto{A},\PrjctnOnto{B})\from \Prdct{B}{A} \to \Prdct{A}{B}$ denotes the twist map. - The proof is left to the reader.
\end{proof}

\begin{example}[Categories enriched in commutative monoids]
  \label{exa:CategoriesEnrichedInCommutativeMonoids}%
  Examples of linear categories include the category of commutative monoids, the category $\AbGrps$ of abelian groups, or the category of modules over a given ring. In the latter two examples Hom-sets are actually abelian groups.
\end{example}

Summarizing, we state:

\begin{theorem}[Recognizing a category of commutative monoid objects]
  \label{thm:CatCommutativeMonoidsRecognize}
  For a pointed category $\SACtgry{X}$ with finite products the following are equivalent:
  \begin{enumerate}[(i)]
    \item $\Ctgry{X}$ is a linear category;
    \item $\Ctgry{X}$ is enriched in the category $\CMon$ of commutative monoids;
    \item every object in $\Ctgry{X}$ carries the structure of an internal commutative monoid such that, for every object $A$ every morphism $f\from X\to Y$ induces a morphism of commutative monoids $f_{\ast}\from \Hom{A}{X}\to \Hom{A}{Y}$;
    \item $\Ctgry{X}\simeq \CMon(\Ctgry{X})$.
  \end{enumerate}
\end{theorem}
\begin{proof}
  The equivalence of (i) and (ii) is given by (\ref{thm:Linear<->EnrichedInCommutativeMonoids}). The implication (ii) $\implies$ (iii) is given by (\ref{thm:MagmaHom-sets->InternalMagma}). The implication (iii) $\implies$ (ii) is given by (\ref{thm:InternalMagma->Hom-Magma}). The equivalence between (iv) and (i) is given by Corollary~\ref{thm:InternalCommMonoidsLinear} and Proposition~\ref{thm:LinearThenCMon}.
\end{proof}

\begin{subordinate}{on the comparison map $\SumProdComp{X}{Y}$}

  Given objects $X$ and $Y$, we constructed the coproduct/product comparison map explicitly as $\SumProdComp{X}{Y}\DefEq\PrdctMapInto{\SumMapOutOf{\IdMapOn{X},\ZeroMap},\SumMapOutOf{\ZeroMap,\IdMapOn{Y}}} = \SumMapOutOf{\PrdctMapInto{\IdMapOn{X},\ZeroMap}, \PrdctMapInto{\ZeroMap,\IdMapOn{Y}} }$.  Remarkably, the maps $\SumProdComp{X}{Y}$ are isomorphisms if there exists any natural family of isomorphisms $\psi_{X,Y} \from X+Y\cong \Prdct{X}{Y}$; see \cite[3]{SLack2012-Isos}.
\end{subordinate}

\begin{exercises}

\begin{exercise}[Internal monoid in preadditive category]
  \label{exe:InternalMonoid-PreAdditiveCat}%
  Show that an internal monoid in a category that is enriched in $\AbGrps$ is an internal abelian group.
\end{exercise}

\begin{exercise}[Product projection normal epi\ZExactTag]
  \label{exe:PoductProjectionNormalEpi}%
  If the comparison map $\SumProdComp{X}{Y}\from X+Y\to \Prdct{X}{Y}$ is an epimorphism, then the projection $\PrjctnOnto{Y}\from \Prdct{X}{Y}\to Y$ is the cokernel of the inclusion $\PrdctMapInto{\IdMapOn{X},\ZeroMap}$. %
  \index{product!projection is normal epi}\index{normal!epi from project projection}%
\end{exercise}

\begin{exercise}[(Normal) subobjects in $\NNr$]
  \label{exe:(Normal)SubsIn-N}
  In the commutative monoid $\NNr = \Set{0,1,2,\dots}$ show the following: %
  \index{subobject!of $\NNr$}\index{normal!subobject of $\NNr$}%
  \begin{thmlist}
    \item A subset $S\subseteq \NNr$ is a subobject of $\NNr$ if and only if there are $a_{1},\dots ,a_{k}\in \NNr$ such that
    \begin{equation*}
      S=\SetSlct{t_{1}a_{1}+\cdots +t_{k}a_{k}}{t_{1},\dots ,t_{k}\in \NNr}
    \end{equation*}
    \item A subset $S\subseteq \NNr$ is a normal subobject of $\NNr$ if and only if there is $a\in \NNr$ such that $S = a\NNr$.
    \item In $\NNr$, the join (as subobjects) of two normal subobjects need not be a normal subobject. %
    \index{join!of (normal) subobjects in $\NNr$}%
    \item In $\CMon$ a pullback of normal monomorphisms need not be a normal pullback. %
    \index{normal!pullback in $\CMon$}%
  \end{thmlist}
\end{exercise}

\begin{exercise}[Biproduct induces special relation\HTag]
  \label{exe:BiproductThenRelation}
  In a homological category, use (\ref{thm:SAPullbackCancellation-II}) to show that if $X\times X\cong X+X$ is a biproduct (\ref{def:BiProduct}), then $(X\oplus X,\pi_1,\nabla_X)$ is a relation (\ref{sec:InternalGraphs/Relations}) on $X$. (Compare with~\ref{thm:MultiplicationInUnitaryMagma-ProductProjection}.)
\end{exercise}
\end{exercises}
\section{Internal Graphs and Internal Relations}%
\label{sec:InternalGraphs/Relations}%

The notion of (equivalence) relation admits categorical conceptualization as in (\ref{def:InternalRelation}). Thus, an internal relation on an object $X$ in a category with binary products $\Ctgry{X}$ is a subobject of the product $\Prdct{X}{X}$. We present the notion of internal relation, along with any of the optional properties of being reflexive, symmetric, or transitive. Thus we arrive at the concept of an internal equivalence relation.

Given a morphism $f\from X\to Y$ in a category $\Ctgry{X}$, we are most interested in the relation given by its kernel pair $(\KrnlPr{f},d_1,d_2)$. If $\Ctgry{X}$ is a variety of algebras, then this relation on $X$ is an equivalence relation, and its equivalence classes are the fibers of $f$.

More generally, a graph from $X$ to $Y$ is a map $g\from G\to \Prdct{X}{Y}$. So, a relation is a special kind of graph. These concepts actually admit a formulation which does not involve products at all:

\begin{definition}[Internal graph]
  \label{def:InternalGraph}%
  In any category an \Defn{internal (directed) graph from $X$ to $Y$} is given by a span $X\XLA{g_{1}} G\XRA{g_{2}} Y$. %
  \index{internal!graph from $X$ to $Y$}\index{graph}\index{graph!internal}%
\end{definition}

\begin{definition}[Internal relation]
  \label{def:InternalRelation}%
  In any category, an \Defn{internal relation from $X$ to $Y$} is an internal graph $X\XLA{g_{1}} G\XRA{g_{2}}$ in which the pair $(g_{1},g_{2})$ is jointly monomorphic. An \Defn{internal relation on $X$} is an internal relation from $X$ to $X$.  %
  \index{relation}\index{relation!from $X$ to $Y$}\index{internal!relation from $X$ to $Y$}%
  \index{relation!on $X$}%
\end{definition}

In a category with binary products, an internal graph is fully determined by a map $g=\PrdctMapInto{g_{1},g_{2}}\from {G\to \Prdct{X}{Y}}$. Such $g$ is an internal relation if and only if it represents a subobject of $\Prdct{X}{Y}$.

\begin{notation}[$aRb$]
  \label{not:Relation}%
  Given an internal relation $X\XLA{d_{1}} R\XRA{d_{2}} Y$, whenever a span $X \XLA{x} A \XRA{y}Y$ factors through $(d_1,d_2)$ via a map $f\from A\to R$, then the map $f$ is unique, and we write $aRb$. %
  \index[not]{r!$aRb$\IndSep $a$ is related to $b$}%
\end{notation}

\begin{example}[Relation on a set]
  When working with a relation on a set $X$, given by a subset of $\Prdct{X}{X}$, a pair $(x,y)\in \Prdct{X}{Y}$ belongs to $R$ if and only if the corresponding function $\PrdctMapInto{x,y}\from A\to \Prdct{X}{X}$ factors through~$R$. This motivates the notation $xRy$ in general.
\end{example}

Recall from Section \ref{sec:SubObjects-QuotientObjects} that subobjects of a given object are partially ordered. Thus, in a category with binary products, we obtain an induced partial ordering of the relations from $X$ to $Y$. Proceeding analogously, we also obtain a partial ordering on such relations in categories which may not admit binary products. As such adjustments are straightforward to make in general, we will henceforth consider graphs and relations in categories with binary products.

\begin{example}[Discrete relation]
  \label{exa:Relation-Discrete}%
  The \Defn{discrete relation} on an object $X$ is given by the diagonal map $\DgnlOn{X}=\PrdctMapInto{1_X,1_X}\from X\to \Prdct{X}{X}$.  %
  \index{discrete!relation}\index{relation!discrete}%
\end{example}

Discrete relations appear `naturally' as the kernel pairs of monomorphisms; see (\ref{thm:KernelPair-Monos}).

\begin{definition}[Reflexive relation]
  \label{def:InternalRelation-Reflexive}
  A relation $d=(d_{1},d_{2})\from R\to \Prdct{X}{X}$ on an object $X$ is \Defn{reflexive} if the discrete relation on $X$ factors through $d$. %
  \index{reflexive relation}\index{relation!reflexive}
\end{definition}

In (\ref{def:InternalRelation-Reflexive}), the relation $R$ is reflexive exactly when the diagram below has a (necessarily unique) filler $e$:
\begin{equation*}
  \xymatrix@R=5ex@C=2em{
  X \ar[rd]_-{\DgnlOn{X}} \ar@{.>}[rr]^-{e} &&
  R \ar@{{ >}->}[dl]^-{d} \\
  & \Prdct{X}{X}
  }
\end{equation*}
Interpreted in $\Sets$, this means that every $x\in X$ is related to itself. In the notation of (\ref{not:Relation}), a relation $d\from R\to \Prdct{X}{X}$ is reflexive if and only if $aRa$ for all $a\from A\to X$. - The discrete relation on $X$ is least among all reflexive relations on $X$.

Note also, that, for every reflexive relation $(R,d_1,d_2,e)$ on $X$, both $d_1$ and $d_2$ are absolute epimorphisms sectioned by $e$.

\begin{definition}[Opposite / symmetric relation]
  \label{def:Opposite/SymmetricRelation}%
  The \Defn{opposite\footnote{This terminology is in line with the concept of the opposite of a category: both are special cases of the opposite of a directed graph.} relation} of a relation $(R,d_{1},d_{2})$ on $X$ is the relation $(R^{\circ},d^{\circ}_{1},d^{\circ}_{2})\DefEq (R,d_{2},d_{1})$. The relation $(R,d_{1},d_{2})$ is said to be \Defn{symmetric} when $R= R^{\circ}$ as subobjects of $X\times X$. %
  \index{opposite!relation}\index{symmetric!relation}\index{relation!opposite}\index{relation!symmetric}%
\end{definition}

Thus, a relation $r\from R\to \Prdct{X}{X}$ is symmetric if the twist map $\tau\DefEq \PrdctMapInto{\PrjctnOnto{2},\PrjctnOnto{1}}\from \Prdct{X}{X}\to \Prdct{X}{X}$ lifts to a map $\tilde{\tau}\from R\to R$ as in this commutative diagram:
\begin{equation*}
  \xymatrix@R=5ex@C=5em{
  R \ar@{{ >}->}[r]^-{r} \ar@{.>}[d]_{\tilde{\tau}} &
  \Prdct{X}{X} \ar[d]^{\tau} \\
  R \ar@{{ >}->}[r]_-{r} &
  \Prdct{X}{X}
  }
\end{equation*}
Whenever such $\tilde{\tau}$ exists, it is necessarily unique. Moreover, it equals its own inverse, as $\tau\Comp \tau=\IdMapOn{\Prdct{X}{X}}$.

In the notation of \eqref{not:Relation}, we may express symmetry as the condition that $aRb$ implies $bRa$. Thus, $aRb$ iff $aR^\circ b$ iff $bRa$ when $R$ is symmetric. Conversely, $d_1Rd_2$ implies $d_2Rd_1$, so that $R^{\circ}\leq R$, from which $R=R^{\circ}$ follows.

\begin{definition}[Transitive relation]
  \label{def:InternalRelation-Transitive}
  In a category with pullbacks, a relation $R$ is called \Defn{transitive} if the diagram below can be filled by a (necessarily unique) map $\rho$ as indicated. %
  \index{transitive relation}\index{relation!transitive}%
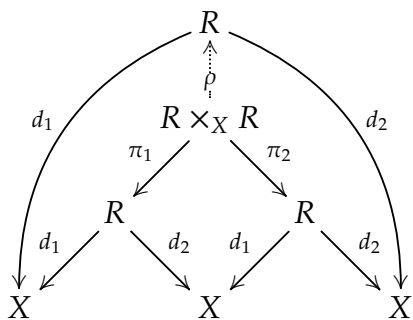
  \begin{equation*}
    \xymatrix@!0@=3em{
    && R \ar@/_2em/[dddll]_-{d_{1}} \ar@/^2em/[dddrr]^-{d_{2}}\\
    && R\times_{X}R \ar@{.>}[u]|(.4){\rho} \ar[dl]_-{\pi_{1}} \ar[dr]^-{\pi_{2}}\\
    & R \ar[dl]_-{d_{1}} \ar[dr]^-{d_{2}} && R \ar[dl]_-{d_{1}} \ar[dr]^-{d_{2}} \\
    X && X && X}
  \end{equation*}
  In this diagram, the square in the center is a pullback.
\end{definition}

In the notation of (\ref{not:Relation}), transitivity amounts to the condition that $aRb$ and $bRc$ imply $aRc$. Indeed, if $aRb$ and $bRc$ for $a$, $b$, $c\from A\to X$, then $\PrdctMapInto{a,b}$ and $\PrdctMapInto{b,c}\from A\to R$ are such that $d_2\Comp\PrdctMapInto{a,b}=d_1\Comp \PrdctMapInto{b,c}$. The universally induced factorization $\PrdctMapInto{a,b,c}\from A\to R\times_XR$ composes to a morphism $\rho\Comp \PrdctMapInto{a,b,c}\from A\to R$ that satisfies $d_i\Comp \rho\Comp \PrdctMapInto{a,b,c}=d_i\Comp \pi_i\Comp \PrdctMapInto{a,b,c}$, which is equal to $a$ when $i=1$ and to $c$ when $i$ is $2$. This shows that $aRc$.

\begin{definition}[Equivalence relation]
  \label{def:Equivalence Relation}%
  In a category with pullbacks, an internal relation is an \Defn{internal equivalence relation} it is reflexive, symmetric and transitive. %
  \index{equivalence relation}%
\end{definition}

\begin{example}[Indiscrete relation]%
  \label{exa:Relation-Indiscrete}%
  The \Defn{indiscrete relation} on $X$ is represented by the identity $\IdMapOn{\Prdct{X}{X}}\from \Prdct{X}{X}\to \Prdct{X}{X}$. It is an equivalence relation on $X$, and is greatest among all relations on $X$. %
  \index{indiscrete!relation}\index{relation!indiscrete}%
\end{example}

\begin{example}[Kernel pair is equivalence relation]
  \label{exa:KernelPair-IsEquivalenceRelation}
  In a category with pullbacks, the kernel pair of a morphism $f\from X\to Y$ is an equivalence relation. %
  \index{kernel pair!equivalence relation}\index{equivalence relation!given by kernel pair}%
\end{example}
\begin{proof}
  \emph{Reflexivity} of the kernel pair relation is assured via the map $s$ in the diagram below.
  
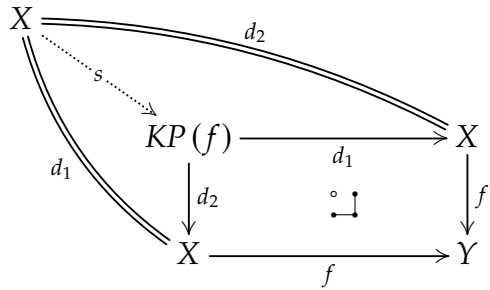
\begin{equation*}
    \xymatrix@R=5ex@C=3em{
    X \ar@{=}@/^2ex/[rrrd]^-{d_2} \ar@{=}@/_2ex/[rdd]_{d_1} \ar@{.>}[rd]|-{\ s\ } \\
    & \KrnlPr{f} \ar[d]^{d_2} \ar[rr]_-{d_1} \PullLU{rrd} &&
    X \ar[d]^{f} \\
    & X \ar[rr]_-{f} &&
    Y
    }
  \end{equation*}
  \emph{Symmetry} of the kernel pair relation is established via the map $T$ below.
  
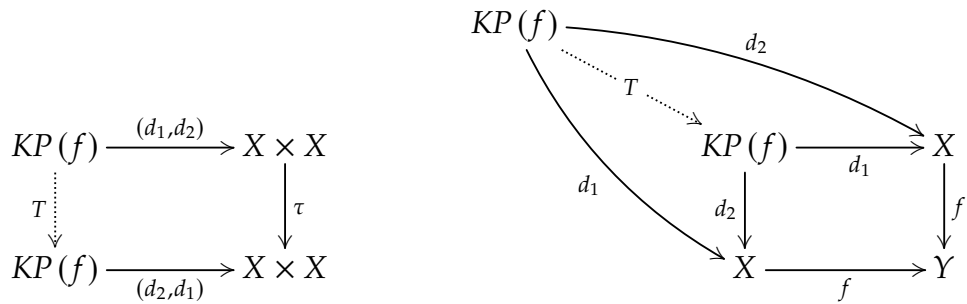
\begin{equation*}
    \xymatrix@R=5ex@C=4em{
    && \KrnlPr{f} \ar@/^2ex/[rrd]^-{d_2} \ar@/_2ex/[ddr]_{d_1} \ar@{.>}[rd]|-{\ T\ } \\
    \KrnlPr{f} \ar@{.>}[d]_{T} \ar[r]^-{\PrdctMapInto{d_1,d_2}} &
    \Prdct{X}{X} \ar[d]^{\tau} &&
    \KrnlPr{f} \ar[r]_-{d_1} \ar[d]_{d_2} &
    X \ar[d]^{f} \\
    \KrnlPr{f} \ar[r]_-{\PrdctMapInto{d_2,d_1}} &
    \Prdct{X}{X} &&
    X \ar[r]_-{f} &
    Y
    }
  \end{equation*}
  \emph{Transitivity} of the kernel pair relation holds because the solid arrows in the diagram
  
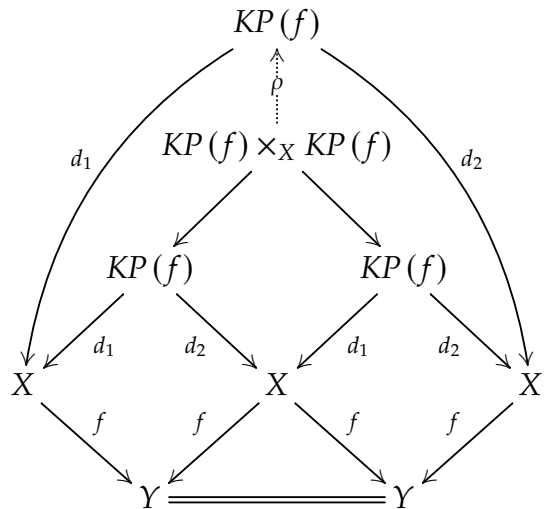
\begin{equation*}
    \xymatrix@R=5ex@!C=2em{
    && \KrnlPr{f} \ar@/^4ex/[rrddd]^{d_2} \ar@/_4ex/[llddd]_{d_1} \\
    && \KrnlPr{f}\prdct_X \KrnlPr{f} \ar[ld] \ar[rd] \ar@{.>}[u]|-{\ \rho\ } \\
    & \KrnlPr{f} \ar[ld]^{d_1} \ar[rd]_{d_2} &&
    \KrnlPr{f} \ar[ld]^{d_1} \ar[rd]_{d_2} \\
    X \ar[rd]^{f} &&
    X \ar[ld]_{f} \ar[rd]^{f} &&
    X \ar[ld]_{f}\\
    & Y \ar@{=}[rr] &&
    Y
    }
  \end{equation*}
  yield a universal map $\rho$ which renders the entire diagram commutative.
\end{proof}

We will see below that there are equivalence relations which are not realizable as a kernel pair. Thus, whenever an equivalence relation is realizable as the kernel pair of some map, it is special, hence merits a name of their own:

\begin{definition}[Effective equivalence relation]
  An equivalence relation which may be obtained as the kernel pair of some morphism is called an \Defn{effective} equivalence relation. %
  \index{effective!equivalence relation}%
\end{definition}

Not every equivalence relation is effective, as the following example shows.

\begin{example}[Example of a non-effective equivalence relation]
  \label{exa:NonEffEqRel}%
  In the category $\Grps(\Tops)$ of topological groups, a subobject is a subgroup whose inclusion is continuous. Hence an equivalence relation $(R,d_1,d_2,e)$ on a topological group $X$ is an equivalence relation (in the usual, set-theoretic sense) on the group underlying $X$, equipped with a topology such that the inclusion $(d_1,d_2)\from R\to \Prdct{X}{X}$ is a continuous group homomorphism. This allows for the topology on $R$ to be any refinement of the subspace topology of its image in $\Prdct{X}{X}$.

  On the other hand, if $R$ is an effective equivalence relation, then necessarily the group $R$ carries the subspace topology induced by $\Prdct{X}{X}$, as all limits in $\Tops$ do. So, as soon as a set $X$, equipped with the indiscrete topology, has two elements or more, the largest equivalence relation $\nabla_X={(X\times X,\pi_1,\pi_2,1_X)}$ on $X\times X$ equipped with the discrete topology is a non-effective equivalence relation on $X$ in the category $\Grps(\Tops)$.
\end{example}

\begin{definition}[Barr exactness]
  \label{def:BarrExact}%
  A regular category (\ref{def:RegularCategory}) is called \Defn{Barr exact} if every equivalence relation is effective. %
  \index{Barr exact category}\index{category!Barr exact}%
\end{definition}

It is known that all algebraic varieties of algebras are Barr exact. By the above, we know that the category of topological groups is not. Here is another example:

\begin{example}[$\AbGrpsTF$ is not Barr exact]
  \label{exa:Ab_tfNotExact}%
  We give an example of a non-effective equivalence relation in the additive category $\AbGrpsTF$ of torsion-free abelian groups. While $\AbGrpsTF$ is a homological category (\ref{sec:RelatedCategoricalStructurs-AbelianCats}), it is not abelian. We prove that it is not Barr exact by considering the relation $(R,d_1,d_2)$ on $\ZNr$ defined by $mRn$ if and only if  $m-n=2k$ for some $k\in \ZNr$. In $\AbGrpsTF$, the coequalizer $q\from \ZNr\to Q$ of the projections $d_1$, $d_2$ is zero, so that the kernel pair of $q$ is $\ZNr\prdct \ZNr$, which contains $R$ as a proper subobject.
\end{example}

\begin{proposition}[Effective equivalence relation as a special square]
  \label{thm:EffectiveEquivalenceRelation-Properties}
  In a finitely bicomplete category, the following conditions are equivalent for an equivalence relation $(R,d_1,d_2,e)$ on an object $X$: %
  \index{effective!equivalence relation: properties}%
  \begin{tfae}
    \item $(R,d_1,d_2)$ is an effective equivalence relation.
    \item $(d_1,d_2)$ is the kernel pair of some morphism $f\from X\to Y$.
    \item $(d_1,d_2)$ is its kernel pair of the coequalizer $q$ of $d_1$ and $d_2$.
    \item The pushout of $d_1$ and $d_2$ yields  a bicartesian square:
    \begin{equation*}
      \xymatrix@R=5ex@C=4em{
      R \BiCart{rd} \ar@{-{ >>}}[r]^{d_1} \ar@{-{ >>}}[d]_{d_2} & X \ar@{-{ >>}}[d]^{q}\\
      X \ar@{-{ >>}}[r]_{q} & Q
      }
    \end{equation*}
  \end{tfae}
\end{proposition}
\begin{proof}
  The equivalence between (I) and (II) holds by definition of an effective equivalence relation. (III) $\implies$ (II) holds by definition.

  (II) $\implies$ (III)\quad If $R$ is the kernel pair of $f$, we find ourselves in this situation:
  \begin{equation*}
    \xymatrix@R=5ex@C=4em{
    R \ar[r]^-{d_{1}} \ar[d]_{d_{2}} \PushRD{rd} &
    X \ar@/^2ex/[rdd]^{f} \ar[d]_{q} \\
    X \ar[r]^-{q} \ar@/_2ex/[rrd]_{f} &
    Q \ar@{.>}[rd]^(0.4){t} \\
    && Y
    }
  \end{equation*}
  The coequalizer of $d_{1}$ and $d_{2}$ is given by the pushout of $d_{1}$ along $d_{2}$; see (\ref{exe:CoEqualizer-Pushout}). So, there is a universal map $t$ which renders the entire diagram commutative. Via $t$, the pullback property of the square $R\rightrightarrows Y$ yields the pullback property of the square $R\rightrightarrows Q$. So, $(d_{1},d_{2})$ is the kernel pair of $q$.

  (IV) $\implies$ (III) holds by definition. Conversely, if $(d_1,d_2)$ is the kernel pair of the coequalizer $q$ of $d_1$ and $d_2$. Then the associated commutative square is a pullback and, reasoning as in the implication from (II) to (III), a pushout. - This completes the proof.
\end{proof}

The following proposition (see~\cite[Proposition 3.9]{EverVdL2} and \cite[Corollary 6]{DBourn2001}) is an interpretation of (\ref{thm:EffectiveEquivalenceRelation-Properties}) in homological categories. It refines some ideas used in the construction of the image factorization in regular categories (\ref{thm:ImageInRegular}).

\begin{proposition}[Reflexive graph in homological category\NTag]
  \label{thm:coeq(d,c)=coker(c circ k_d)}%
  In a normal category $\Ctgry{X}$, consider a reflexive graph on an object $X$:
  \begin{equation*}
    \xymatrix@R=5ex@C=5em{
    R \ar@<1ex>[r]^-{d_1} \ar@<-1ex>[r]_-{d_2} & X \ar[l]|-{\ e\ }
    }\qquad\qquad
    d_1e=1_X=d_2e
  \end{equation*}
  Write $k\from K\to R$ for the kernel of $d_1$. Then the following hold:
  \begin{enumerate}[(i)]
    \item The cokernel of  $d_2k\from K\to X$ is the coequalizer of $d_1$ and $d_2$.
    \item $(R,d_1,d_2)$ is a relation on $X$ if and only if $d_2k\from K\to X$ is a monomorphism.
    \item If $\Ctgry{X}$ is homological, then $(R,d_1,d_2)$ is an effective equivalence relation on $X$ if and only if $d_2k\from K\to X$ is a normal monomorphism.
  \end{enumerate}
\end{proposition}
\begin{proof}
  The maps $d_{1}$ and $d_{2}$ are normal epimorphisms by the \AENInline-condition (\ref{def:NormalCat}). From (\ref{exe:CoEqualizer-Pushout}) we know that the pushout of $d_{1}$ along $d_{2}$ yields the coequalizer  $q\from X\to Q$ of $d_1$ and $d_2$:
  \begin{equation*}
    \xymatrix@R=5ex@C=4em{
    K \ar@{=}[d] \ar@{{ |>}->}[r]^-{k} &
    R \ar@{-{ >>}}[r]^-{d_1} \ar@{-{ >>}}[d]_-{d_2} \PushRD{rd} &
    X \ar[d]^-{q} \\
    K \ar[r]_-{d_2k} &
    X \ar[r]_-{q} & Q
    }
  \end{equation*}
  The top row is a short exact sequence. So $q=\CoKerMap{d_{2}k}$ by (\ref{thm:NormalEpi-Props}.\ref{thm:CoKer(gf)}). - This proves (i).

  (ii)\quad If $R$ is a relation then $d_{1}$ and $d_{2}$ are jointly monic. So, $d_1k=0$ and $d_2k$ are jointly monic as well. So $d_2k$ is a monomorphism. To see the converse, take the factorization of $\PrdctMapInto{d_{1},d_{2}}\from R\to X\times X$ as a normal epimorphism $e$ followed by a monomorphism $(m_1,m_2)$. Then the kernel $l\from L\to X$ of $e$ is such that $d_2l=m_2el=0$. So, it factors through $k$ as a morphism $l'\from L\to K$ such that $kl'=l$. Now $d_2kl'= m_2el=0$, so if $d_2k$ is a monomorphism then $l'=0$. It follows that $e$ is an isomorphism and $(R,d_1,d_2)$ is a relation on $X$.

  (iii)\quad If $R$ is an effective equivalence relation then the square on the right is a pullback. So, $\Ker{q}$ may be represented by $d_2k$. In particular, $d_{2}k$ is a normal monomorphism. Conversely, if $d_2k$ is a normal monomorphism, then it is the kernel of its cokernel $q$, which implies that the square on the right is a pullback by Proposition~\ref{thm:PullbackFromKerIso}. So, $(d_{1},d_{2})$ is an effective equivalence relation on $X$.
\end{proof}

\begin{definition}[Direct image of a relation\HTag]
  \label{def:DirectImageOfRelation}%
  In a homological category the \Defn{direct image} of a relation $(R,d_{1},d_{2})$ on an object $X$ under a normal epimorphism $f\from{X\to Y}$ is $(f_{\ast}R,d'_{1},d'_{2})$, given by the image factorization of $(\Prdct{f}{f})\PrdctMapInto{d_1,d_2}$, as in the commutative diagram below. %
  \index{direct image!of a relation under a normal epimorphism}%
  \begin{equation*}
    \xymatrix@R=5ex@C=4em{
    R \ar@{{ >}->}[r]^-{\PrdctMapInto{d_1,d_2}} \ar@{-{ >>}}[d] &
    \Prdct{X}{X} \ar@{-{ >>}}[d]^{\Prdct{f}{f}} \\
    f_{\ast}R \ar@{{ >}->}[r]_-{\PrdctMapInto{d'_{1},d'_{2}}} &
    \Prdct{Y}{Y}
    }
  \end{equation*}
\end{definition}

\begin{lemma}[Direct image preserves equivalence relations\HTag]
  \label{thm:directimageeqrel}
  The `direct image' operation preserves equivalence relations.
\end{lemma}
\begin{proof}[Proof sketch]
  First verify directly that the direct image of a reflexive, symmetric relation along a normal epimorphism is also reflexive and symmetric. To explain why transitivity is preserved,  show that when $S$ is the image of $R$ along a normal epimorphism $f\from {X\to Y}$, then the induced map $R\times_{X}R\to S\times_{Y}S$ is a normal epimorphism. This follows from Corollary~\ref{thm:NormalPushout-PullbackStabilty} combined with Corollary~\ref{thm:SplitEpiOfNormalEpis->NormalPush}.
\end{proof}

\begin{subordinate}{}
  \begin{subsubordinate}{On the relevance of Barr exactness to categorical algebra}
    We will show that all varieties of algebras are Barr exact. They are also a key ingredient in the definition of an abelian category---see Chapter~\ref{chap:AbelianCategories}. Further, in Section~\ref{sec:ExactMaltsev}, we prove that a homological category is Barr exact if and only if it is di-exact; i.e., every antinormal map is normal. As a consequence, every homological variety of algebras is semiabelian.
  \end{subsubordinate}
\end{subordinate}
\chapter[Abelian Categories]{Abelian Categories}
\label{chap:AbelianCategories}%

Efforts toward categorically capturing selected features of abelian groups go back at least as far as Mac Lane in the late 1940's. Beginning with the observation that interpreting certain properties of abelian groups in diagram language reveals an otherwise less obvious form of duality, he conceives of this possibility \cite[p.~494]{MacLane:Duality}: \emph{`One of our chief objectives is that of providing a background in which the proofs for axiomatic homology theory become exactly dual to those for cohomology theory.'}

Then Cartan--Eilenberg published the first-ever textbook on homological algebra \cite{Cartan-Eilenberg}. While its focus is on deriving functors between categories of modules, in its appendix, Buchsbaum applies Mac Lane's axiomatization ideas to the particulars of homological algebra; see also \cite{Buchsbaum:ExactCats}. The result is an abstract categorical setting in which the methods of homological algebra, as presented by Cartan--Eilenberg, are applicable.

Grothendieck's notion of abelian category in \cite{Tohoku} extends the scope of homological algebra considerably. For example, if $\Ctgry{A}$ is an abelian category, then so is the category of functors from a small category $J$ to $\Ctgry{A}$. Thus any category of presheaves with values in an abelian category forms an abelian category, as does any category of sheaves. Thus, Grothendieck provides a setting whose scope encompasses in a unifying manner what Cartan--Eilenberg presented in module categories as well as Leray's sheaf cohomology \cite{JLeray1_1946,JLeray2_1946}. Standard references for abelian categories include \cite{Tohoku,Freyd,Borceux:Cats,Weibel,PAluffi2009}.

Here, we encounter abelian categories from at least three points of view:
\begin{enumerate}[(A)]
  \item Coming from a background of abelian categories, homological categories and semiabelian categories appear as vast generalizations.
  \item Coming from a background of \ZExact, normal, homological, or semiabelian categories, abelian categories appear as a specialization, with particularly nice features.
  \item Every homological category $\Ctgry{X}$ has a \emph{left almost abelian} core $\AbCoreOf{X}$, the full subcategory of objects in $\Ctgry{X}$ which admit an internal abelian group structure. If $\Ctgry{X}$ is semiabelian, then its abelian core is an abelian category, contained in $\Ctgry{X}$ as a full and replete subcategory. It is closed under subobjects and quotients. 
\end{enumerate}

Abelian group objects in a homological category are particularly easy to identify: We show that every internal magma in a homological category is automatically an abelian group object. `Tierney's equation'
\[
  \text{abelian = additive + exact}
\]
(Theorem~\ref{thm:AdditivePlusp-exactIsAbelian}) expresses that a Barr exact category is an abelian category if and only if it admits biproducts and a zero object. In particular, for a semiabelian category to be abelian, it suffices that for any two objects $A$ and $B$ the comparison morphism $\gamma_{A,B}\from {A+B\to A\times B}$ from Proposition~\ref{thm:Sum->ProductIsCokernel} is not just a regular epimorphism, but also a monomorphism, so that binary sums and binary products coincide (they are \Defn{biproducts}~\ref{def:BiProduct}). A category that satisfies this condition is called a \Defn{linear} category (\ref{def:LinearCategory}), so that
\[
  \text{abelian = semiabelian + linear.}
\]
We saw that a pointed category is linear precisely when it is enriched in commutative monoids (\ref{thm:Linear<->EnrichedInCommutativeMonoids}), so additivity is equivalent to the existence of inverses to the sum in the hom-monoids. The calculus of relations allows us to view this fact differently: a homological category is abelian if and only if every monomorphism is normal (\ref{thm:AbelianViaSemiabelianWithNormalMonos}). Yet another viewpoint connects right away with the approach in Chapter~\ref{chap:Di-ExactCats}: we have that
\[
  \text{abelian = additive + di-exact.}
\]
On the other hand, the fundamental non-self-dual nature of semiabelian categories is reflected in the fact that
\[
  \text{abelian = semiabelian + co-semiabelian.}
\]
All of these diverse takes on the concept of an abelian category will be will be treated in the current chapter.

We will use select combinations of following properties a category $\Ctgry{X}$ may have:
\begin{enumerate}[({P}.1)]
  \item \label{Property:HasFiniteProducts}%
        \begin{enumerate}[(a)]
          \item $\Ctgry{X}$ has finite products.
          \item $\Ctgry{X}$ has kernels.
          \item $\Ctgry{X}$ is finitely complete.
        \end{enumerate}
  \item \label{Property:HasFiniteSums}%
        \begin{enumerate}[(a)]
          \item $\Ctgry{X}$ has finite sums.
          \item $\Ctgry{X}$ has cokernels.
          \item $\Ctgry{X}$ is finitely cocomplete.
        \end{enumerate}
  \item \label{Property:NormalMonosClosedUnderCobaseChange}%
        Pushouts preserve normal monomorphisms.
  \item \label{Property:NormalEpisClosedUnderBaseChange}%
        \PNEInline: Pullbacks preserve normal epimorphisms.
  \item \label{Property:MonoIsNormal}%
        Every monomorphism in $\Ctgry{X}$ is a normal monomorphism. %
  \item \label{Property:EpiIsNormal}%
        Every epimorphism in $\Ctgry{X}$ is a normal epimorphism.
  \item \label{def:Ab-Enriched}%
        $\Ctgry{X}$ is enriched in the category $\AbGrps$ of abelian groups. This means: for all objects $X$ and~$Y$ in~$\Ctgry{X}$, the set $\HomIn{\Ctgry{X}}{X}{Y}$ carries the structure of an abelian group. Furthermore, with respect to these abelian group structures, the composition function
        \begin{equation*}
          \HomIn{\Ctgry{X}}{Y}{Z}\prdct \HomIn{\Ctgry{X}}{X}{Y} \longrightarrow \HomIn{\Ctgry{X}}{X}{Z}
        \end{equation*}
        is bilinear.
\end{enumerate}

Using combinations of the above properties, a hierarchy of types of categories appears in the literature. Here are the definitions of these types of categories. A discussion of their basic properties forms part of the body of this chapter.

\begin{definition}[Hierarchy leading to abelian categories]
  \label{def:PreAdditiveCat}\quad %
  \begin{enumerate}[({Ab}.1)]
    \item A category $\Ctgry{X}$ is \Defn{preadditive} if it is enriched in $\AbGrps$; i.e.\ (P.\ref{def:Ab-Enriched}) holds. %
          \index{preadditive category}\index{category!preadditive}%
          By forgetting structure, we have that $\Ctgry{X}$ is enriched in $\CMon$; cf.\ Section~\ref{sec:BiProducts}.
    \item A category $\Ctgry{X}$ is \Defn{additive} if it preadditive and admits finite sums; i.e.\ (P.\ref{Property:HasFiniteSums}.a) holds. By the results in (\ref{sec:BiProducts}) then, $\Ctgry{X}$ is linear, so we obtain Proposition \ref{thm:BiProduct-From-(Co)Product-Additive} saying that this is equivalent to the existence of finite products (P.\ref{Property:HasFiniteProducts}.a). Furthermore, the category is pointed, and each finite sum is actually a finite biproduct, see (\ref{thm:BiProduct-From-(Co)Product-Additive}). On the other hand, in any additive category, the \KSGInline-condition holds (\ref{thm:AddThenProto}) and by self-duality, also the dual of this condition.%
          \index{additive!category}\index{category!additive}%
    \item\label{def:PreAbelianCat} A category $\Ctgry{X}$ is \Defn{preabelian} if it is additive and every morphism has a kernel and a cokernel. So a preabelian category is, in particular, z-exact. In other words, a category is preabelian when it is enriched in $\AbGrps$, has a zero object, and is finitely complete and finitely cocomplete; see (\ref{thm:PreAbelianCategory-Recognize}).\index{preabelian category}\index{category!preabelian}
          The validity of \KSGInline\ now implies that a preabelian category is both protomodular and co-protomodular (\ref{sec:Protomodular-SEpi(X)->X}).
    \item \label{def:LAACat} A category is \Defn{left almost abelian} if it is preabelian and (P.\ref{Property:NormalEpisClosedUnderBaseChange}) holds: a pullback of a normal epimorphism is again a normal epimorphism, condition \PNEInline. Since also \KSGInline\ holds, a left almost abelian category is always homological (\ref{sec:HomologicalCats-Axioms}). Dually, a \Defn{right almost abelian} category is preabelian with pushout-stable normal monomorphisms. It is always a co-homological category. An \Defn{almost abelian category} satisfies both conditions. Left almost abelian, right almost abelian, and almost abelian categories coincide with so-called left quasi-abelian, right quasi-abelian and quasi-abelian categories, respectively. We show that a category is quasi-abelian/almost abelian when it is simultaneously homological and co-homological.
    \item \label{def:AbelianCat} A category is \Defn{abelian} if it is a preabelian category in which conditions (P.\ref{Property:MonoIsNormal}) and (P.\ref{Property:EpiIsNormal}) hold. We show (\ref{thm:AdditivePlusp-exactIsAbelian}) that in an additive context, this is equivalent to p-exactness; i.e., all morphisms are normal. As a consequence, an additive category is abelian if and only if it is di-exact (\ref{thm:AdditivePlusDiexactIsAbelian}).\index{abelian category}\index{category!abelian}\index{p-exact category}\index{category!p-exact}\index{di-exact category}\index{category!di-exact} We further explain that a category is abelian if and only if it is both semiabelian and co-semiabelian, if and only if it is both additive and Barr-exact.
  \end{enumerate}
\end{definition}

\bigskip\bigskip\bigskip

\begin{center}
  {\bf Organization of the sections in this chapter}
\end{center}

\bigskip

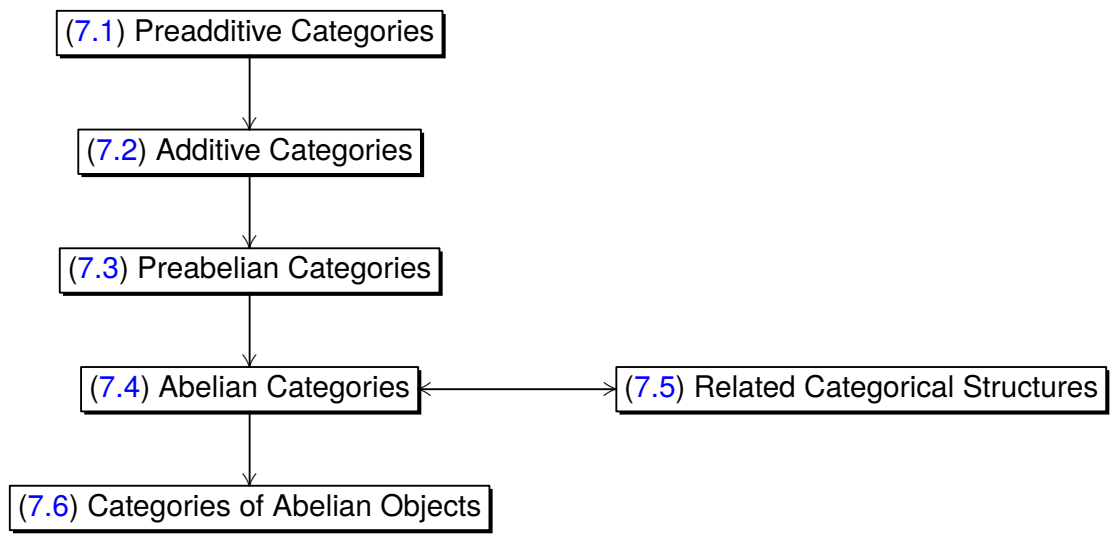
\begin{equation*}
  \xymatrix@R=5ex@C=4em{
  *+[F-,]{\txt{\sffamily (\ref{sec:PreAdditiveCats}) Preadditive Categories}} \ar[d] \\
  *+[F-,]{\txt{\sffamily (\ref{sec:AdditiveCategory}) Additive Categories}} \ar[d] \\
  *+[F-,]{\txt{\sffamily (\ref{sec:PreAbelianCategories}) Preabelian Categories}} \ar[d] \\
  *+[F-,]{\txt{\sffamily (\ref{sec:AbelianCategories}) Abelian Categories}} \ar[d] \ar@{<->}[r] &
  *+[F-,]{\txt{\sffamily (\ref{sec:RelatedCategoricalStructurs-AbelianCats}) Related Categorical Structures}} \\
  *+[F-,]{\txt{\sffamily (\ref{sec:CatsOfAbelianObjects}) Categories of Abelian Objects}}
  }
\end{equation*}
\newpage
\section[Preadditive Categories]{Preadditive Categories}
\label{sec:PreAdditiveCats}%

By definition, a preadditive category is enriched in abelian groups and, hence, is found among categories which are enriched in commutative monoids; see Section \ref{sec:BiProducts}. Thus, in a preadditive category $\Ctgry{X}$, any finite product is also a coproduct, and any finite coproduct is also a product; see (\ref{thm:BiProduct-From-(Co)Product}). This means that in a preadditive category, properties (P.\ref{Property:HasFiniteSums}.a) and (P.\ref{Property:HasFiniteProducts}.a) are equivalent, and each implies the existence of a zero object.

Worded differently, for objects $A$ and $B$ in a preadditive category $\Ctgry{X}$, their sum exists if and only if their product exists, and this happens if and only if the comparison map $\SumProdComp{A}{B}\from A+B\to \Prdct{A}{B}$ is an isomorphism. By (\ref{thm:BiProduct-By-CMonEnrichment}), this is equivalent to the existence of a commutative diagram
\stepcounter{theorem}%
\begin{equation}\label{eqn:Biproduct-PreAdditive}
  \vcenter{\xymatrix@R=3ex@C=4em{
  A \ar@{=}[rr] \ar[rd]_{i_A} &&
  A \\
  & X \ar[ru]_{p_A} \ar[rd]^{p_B} \\
  B \ar[ru]^{i_B} \ar@{=}[rr] &&
  B
  }}
\end{equation}
in which $(X,i_A,i_B)$ is a coproduct of $A$ and $B$, and $(X,p_A,p_B)$ is a product of $A$ and $B$, while $p_Bi_A=\ZeroMap_{BA}$, and $p_Ai_B=\ZeroMap_{AB}$.

Moreover, one can recognize a biproduct via diagram~\eqref{eqn:Biproduct-PreAdditive}, plus the condition that $i_Ap_A+i_Bp_B=\IdMapOn{X}$; i.e.\ one obtains a universal property via a criterion which makes no mention whatsoever of universality (!); see (\ref{thm:BiProduct-By-AbEnrichment}). It enables us to show that any binary (co-)product diagram admits an expansion to a biproduct diagram which, moreover, is unique; see (\ref{thm:BiProduct-From-(Co)Product}).

In particular, an initial object or a terminal object in a preadditive category is automatically a zero-object; see (\ref{thm:ZeroObject-In-LinearCat}). However, the preadditive categories are segregated into two classes, namely the ones with zero object and the ones without; see Exercise~\ref{exa:PreAdditiveCat-UnderlyingEnrichment-Set_*}.

Via the additive inverse operation in the Hom-groups of a preadditive category $\Ctgry{X}$ we obtain a new criterion for when two morphisms $f$, $g\from A\to B$ are equal: this happens if and only if $g-f=\ZeroMap$. Taking advantage of this criterion, we obtain nice recognition criteria for detecting a monomorphism, respectively an epimorphism, as in (\ref{thm:PreAdditiveCat-Mono/Epi-Recognize}).

\begin{lemma}[Recognizing a biproduct via $\AbGrps$-enrichment]
  \label{thm:BiProduct-By-AbEnrichment}%
  In a preadditive category diagram \eqref{eqn:Biproduct-PreAdditive} describes $X$ as a biproduct of $A$ and $B$ if and only if %
  \index{biproduct!in preadditive category} %
  \begin{equation*}
    p_Bi_A=\ZeroMap_{BA} \qquad \text{and}\qquad p_Ai_B=\ZeroMap_{AB} \qquad \text{and}\qquad i_Ap_A + i_Bp_B=\IdMapOn{X}.
  \end{equation*}
\end{lemma}
\begin{proof}
  This is a special case of Lemma \ref{thm:BiProduct-By-CMonEnrichment}.
\end{proof}

\begin{proposition}[Biproduct from product or coproduct]
  \label{thm:BiProduct-From-(Co)Product-Additive}
  For any two objects $A$ and $B$ in a preadditive category the following hold:
  \begin{enumerate}[(i)]
    \item Any product cone $A \XLA{p_A} X \XRA{p_B} B$ admits an expansion to a biproduct diagram \eqref{eqn:Biproduct-PreAdditive}. Moreover, this expansion is unique.
    \item Any coproduct cocone $A \XRA{i_A} X \XLA{i_B} B$ admits an expansion to a biproduct diagram \eqref{eqn:Biproduct-PreAdditive}. Moreover, this expansion is unique. \NoProof
  \end{enumerate}
\end{proposition}

\begin{proposition}[Recognizing mono/epimorphism in a preadditive category]%
  \label{thm:PreAdditiveCat-Mono/Epi-Recognize}%
  In a preadditive category $\Ctgry{X}$ the following hold:
  \begin{enumerate}[(i)]
    \item A morphism $ f\from A\to B$ in $\Ctgry{X}$ is a monomorphism if and only if for every $z\from Z\to A$, $ fz=0$ implies $z=0$.
    \item A morphism $ f\from A\to B$ in $\Ctgry{X}$ is an epimorphism if and only if for every $ \omega\from B\to Z$, $ \omega f=0$ implies $ \omega=0$.
  \end{enumerate}
\end{proposition}
\begin{proof}
  Suppose $ fz=0$ implies $z\from Z\to A$. To see that $f$ is a monomorphism, consider $u$, $v\from Z\to A$ with $fu=fv$. Then $f(v-u)=0$, implying that $v-u=0$ and, hence, $u=v$. The proof of (ii) is similar.
\end{proof}

\begin{exercises}

\begin{exercise}
  \label{exe:EmbeddingCatInPreAdditiveCat}%
  Show that every category $\Ctgry{X}$ can be embedded into a preadditive category. Hint: Consider each Hom-set as a basis for a free abelian group.
\end{exercise}

\begin{exercise}[Preadditive categories not pointed]
  \label{exa:PreAdditiveCat-UnderlyingEnrichment-Set_*}
  Show that a preadditive category need not have a zero object. Still, it has an underlying enrichment in $\SetsBsd$. %
  \index{preadditive category!need not be pointed}%
\end{exercise}
\end{exercises}
\newpage
\section{Additive Categories}
\label{sec:AdditiveCategory}%

An \Defn{additive category} is a preadditive category which admits finite sums. So, it has a zero object because the sum over the empty category is an initial object, hence also a terminal object by (\ref{thm:ZeroObject-In-LinearCat}). Further, via (\ref{thm:BiProduct-From-(Co)Product}) the existence of finite sums implies that finite products as well as finite biproducts exist. Thus a category is additive if and only if its opposite is additive.

One feature an additive category has which need not be shared by a preadditive one is that the addition in Hom-sets is given by the diagram in Proposition \ref{thm:AdditiveCategoryHom}.

An alternate perspective of additive categories follows from the development in Section~\ref{sec:BiProducts}. There we proved that a pointed category with binary sums and products is a linear category (\ref{def:LinearCategory}) if and only if it is enriched in commutative monoids; see (\ref{thm:Linear<->EnrichedInCommutativeMonoids}). Thus a category is additive if and only if it is linear and its Hom-monoids are groups. As a consequence every object carries a canonical internal abelian group structure; see (\ref{thm:MagmaHom-sets->InternalMagma}). --- We use the tag {\color{OliveGreen} $\EuRoman{A{\kern-0.15ex}dd}$} to identify additive categories.

As in Theorem \ref{thm:MagmaHom-sets->InternalMagma}, addition in Hom-sets may be computed as follows:

\begin{proposition}[Hom-addition in an additive category\AddTag]%
  \label{thm:AdditiveCategoryHom}%
  In an additive category the sum of two morphisms $f$, $g\from A\to B$ may be computed via this commutative diagram:
  \begin{equation*}
    \xymatrix@R=5ex@C=3em{
    A \ar[r]^-{f+g} \ar[d]_{\DgnlOn{A}} &
    B \\
    {A}\oplus{A} \ar[r]_-{{f}\oplus{g}} &
    B\oplus B \ar[u]_{\FoldOn{B}}
    }
  \end{equation*}
  In other words, $f+g=\FoldOn{B}\Comp({f}\oplus{g})\Comp\DgnlOn{A}$.\NoProof
\end{proposition}

Additive categories need not admit kernels or cokernels. The following result can help identify situations where they do exist. Observe how the presence of a zero object extends the scope of the criterion by which we identify monomorphisms and epimorphisms in a preadditive category (\ref{thm:PreAdditiveCat-Mono/Epi-Recognize}):

\begin{proposition}[(Co-)equalizers in a additive category\AddTag]
  \label{thm:(Co)Equalizers-AdditiveCat}
  In an additive category two maps $f$, $g\from A\to B$ admit a (co-)equalizer if and only if $g-f$ admits a (co-)kernel. When this happens, (co-)equalizers of $f$ and $g$ and (co-)kernels of $g-f$ agree. \NoProof
\end{proposition}

For the following result recall that a monomorphism is regular when it is the equalizer of a pair of maps.

\begin{corollary}[Regular monos and epis are normal\AddTag]
  \label{thm:RegularVsNormalInAdditive}
  In an additive category, every regular epimorphism is a normal epimorphism, and every regular monomorphism is a normal monomorphism.\NoProof
\end{corollary}

\begin{proposition}[Recognizing mono/epi-morphisms in an additive category\AddTag]%
  \label{thm:AdditiveCat-Mono/Epi-Recognize}%
  In an additive category the following hold:
  \begin{enumerate}[(i)]
    \item If $ f\from A\to B$ has a kernel, then $f$ is a monomorphism if and only if $\KerMap{f}= (0\to A)$.
    \item If $ f\from A\to B$ has a cokernel, then $f$ is an epimorphism if and only if $\CoKerMap{f} = (B\to 0)$.
  \end{enumerate}
\end{proposition}
\begin{proof}
  (i) If $f$ is a monomorphism, then (\ref{thm:Ker(Mono)=0}) tells us that $\Ker{f} = \ZeroMap$. Conversely, if $u\from Z\to A$ is such that $fu=\ZeroMap$, then $u$ factors uniquely through $\Ker{f}=\ZeroMap$. So $u=\ZeroMap$ and $f$ is a monomorphism by (\ref{thm:PreAdditiveCat-Mono/Epi-Recognize}). --- Dualizing, we verify (ii).
\end{proof}

\begin{corollary}[Hom-recognition of kernels\AddTag]
  \label{thm:AdditiveCat-HomRecognition-Of-Kernels}
  If $f\from A\to B$ is a morphism in an additive category, then $\kappa\from K\to A$ is a kernel for $f$ if and only if, for each object $U$ in $\Ctgry{X}$, the induced sequence below is exact.
  \begin{equation*}
    0 \longrightarrow \HomIn{\Ctgry{X}}{ U }{ K } \XRA{\kappa_{\ast}} \HomIn{\Ctgry{X}}{ U }{ A } \XRA{f_{\ast}} \HomIn{\Ctgry{X}}{ U }{ B } \NoProofDiag
  \end{equation*}
\end{corollary}

\begin{corollary}[Hom-recognition of cokernels\AddTag]
  \label{thm:AdditiveCat-HomRecognition-Of-CoKernels}
  If $f\from A\to B$ is a morphism in an additive category, then $\omega\from B\to Q$ is a cokernel for $f$ if and only if, for each object $W$ in $\Ctgry{X}$, the induced sequence below is exact.
  \begin{equation*}
    0 \longrightarrow \HomIn{\Ctgry{X}}{Q }{ W } \XRA{\omega^{\ast}} \HomIn{\Ctgry{X}}{ B }{ W } \XRA{f^{\ast}} \HomIn{\Ctgry{X}}{ A }{ W } \NoProofDiag
  \end{equation*}
\end{corollary}

\begin{proposition}[Pullback / pushout via exactness\AddTag]
  \label{thm:AdditiveCategory-Pullback/Pushout-Recognition}
  In an additive category $\Ctgry{X}$, given the diagram (C) below, construct the sequence on the right as indicated.
  \begin{equation*}
    \vcenter{\xymatrix@R=5ex@C=4em{
    A \ar[r]^-{a} \ar[d]_{b}\ar@{}[rd]|-{\ \text{(C)}\ } &
    U \ar[d]^{u} \\
    V \ar[r]_-{v} &
    Z
    }}\qquad\qquad \vcenter{\xymatrix@R=5ex@C=4em{
    A \ar[r]^-{\PrdctMapInto{a,b}} &
    U\oplus V \ar[r]^-{\SumMapOutOf{u,-v}} &
    Z
    }}
  \end{equation*}
  Then the following hold:
  \begin{enumerate}[(i)]
    \item (C) commutes if and only if the composite $\SumMapOutOf{u,-v}\Comp \PrdctMapInto{a,b} =\ZeroMap$.
    \item (C) is a pullback diagram if and only if $\PrdctMapInto{a,b}=\KerMap{\SumMapOutOf{u,-v}}$.
    \item (C) is a pushout diagram if and only if $\SumMapOutOf{u,-v}=\CoKerMap{\PrdctMapInto{a,b}}$.
    \item (C) is both a pullback and a pushout diagram if and only if the sequence is short exact; see (\ref{def:ShortExactSequence-Basic}).
  \end{enumerate}
\end{proposition}
\begin{proof}
  (i) follows from $\SumMapOutOf{u,-v}\Comp \PrdctMapInto{a,b}= (ua) + (-v)b= ua -vb$; see (\ref{exe:Addition-Composites}) and (\ref{exe:AbCat-Inverse-Composition}). To see (ii), observe that, for any $r=(r_U,r_V)\from R\to U\oplus V$,
  \begin{equation*}
    \ZeroMap=\SumMapOutOf{u,-v}\Comp (r_U,r_V) = (ur_U) - (vr_V)\quad \text{if and only if}\quad ur_U=vr_V
  \end{equation*}
  Thus the pullback property of diagram (C) corresponds to the kernel property of $\PrdctMapInto{a,b}$. Statement	(iii) follows from (ii) by duality, and (iv) combines (ii) and (iii).
\end{proof}

\begin{corollary}[Pullbacks and monos, pushouts and epis\AddTag]
  \label{thm:AbelianCategory-PullbackPreserves/ReflectsMono/PushoutPreserves/ReflectsEpis}
  For a commutative diagram (C) in additive category the following hold:
  \begin{equation*}
    \xymatrix@R=5ex@C=4em{
    A \ar[r]^-{a} \ar[d]_{b} \ar@{}[rd]|-{\ \text{(C)}\ } &
    U \ar[d]^{u} \\
    V \ar[r]_-{v} &
    Z
    }
  \end{equation*}
  \begin{enumerate}[(i)]
    \item If (C) is a pullback diagram then it preserves and reflects monomorphisms.
    \item If (C) is a pushout diagram then it preserves and reflects epimorphisms.
  \end{enumerate}
\end{corollary}
\begin{proof}
  (i) We know from (\ref{thm:KernelFunctor-Props}) that $\Ker{u}\cong \Ker{b}$. Thus, if one of these two maps is monic, then both kernels vanish, and so (\ref{thm:AdditiveCat-Mono/Epi-Recognize}) both maps are monomorphisms. --- The proof of (ii) is dual.
\end{proof}

\subsection[Related categorical properties]{Related categorical properties}
\label{subsec:AdditiveCats-RelatedStructures}%

We saw in (\ref{thm:RegularVsNormalInAdditive}) that, in an additive category, normal and regular epimorphisms coincide. This property is a prominent feature of normal categories, see (\ref{thm:CoKer=NormalEpi=RegEpi=EffectiveEpi}). So, we may ask whether additive categories normal. The answer is `no' for two reasons. Firstly, an additive category need not admit finite limits nor colimits. Second, pullbacks need not preserve normal epimorphisms.

The key difference between normal categories and homological ones is the validity of the \KSGInline-condition; see (\ref{def:HomologicalCategory}). Perhaps surprisingly, all additive categories satisfy it.

\begin{proposition}[Additive categories and the \KSGInline-condition\AddTag]
  \label{thm:AddThenProto}
  \label{thm:Additive->KSG}%
  In every additive category the \KSGInline-condition holds. %
  \index{KSG!in additive categories}\index{additive category!satisfies \KSGInline-condition}%
\end{proposition}
\begin{proof}
  Consider a map $p$, sectioned by $s$, and a kernel $\kappa=\KerMap{p}$.
  \begin{equation*}
    \xymatrix@C=4em{
    K \ar@{{ |>}->}[r]^-{\kappa} &
    X \ar@<-.5ex>[r]_-{p} &
    Y \ar@<-.5ex>[l]_-{s}
    }
  \end{equation*}
  Assume $m\from M\to X$ is a monomorphism through which $\kappa$ and $s$ factor as $\lambda\from K\to M$ and $t\from Y\to M$, respectively; that is $m\Comp \lambda=\kappa$ and $m\Comp t=s$. We notice that $p\Comp (\IdMapOn{X}-s\Comp p)=0$. Thus there exists $q\from X\to K$ unique with $\kappa\Comp q=\IdMapOn{X}-s\Comp p$. Therefore:
  \begin{equation*}
    \IdMapOn{X} =\kappa\Comp q+s\Comp p=m\Comp\lambda\Comp q+m\Comp t\Comp p=m\Comp(\lambda\Comp q+ t\Comp p)
  \end{equation*}
  This means that the monomorphism $m$ is also an absolute epimorphism, hence an isomorphism by (\ref{exe:SplitEpis}).
\end{proof}

\begin{exercises}

\begin{exercise}[(Co-)equalizers in an additive category]
  \label{exe:(Co)Equalizers-Additive}
  Prove (\ref{thm:(Co)Equalizers-AdditiveCat}): In an additive category two morphisms $f$, $g\from A\to B$ admit a (co-)equalizer if and only if $g-f\from A\to B$ admits a (co-)kernel. %
  \index{coequalizer!in additive category}\index{equalizer!in additive category}%
\end{exercise}

\begin{exercise}[Additive inverse operation and composition in Hom-groups]
  \label{exe:AbCat-Inverse-Composition}
  Given composable maps $A\XRA{f} B \XRA{g}C$ in an additive category show that
  \begin{equation*}
    g\Comp (-f) = -(g\Comp f) = (-g)\Comp f
  \end{equation*}
\end{exercise}

\begin{exercise}[Hom-addition in an additive category]
  \label{exe:HomAddition-In-Additive}
  For morphisms $f$, $g\from A\to B$ in an additive category show that $f+g=\FoldOn{B}\Comp (f\oplus g)\Comp \DgnlOn{A}$; see (\ref{thm:AdditiveCategoryHom}).
\end{exercise}

\begin{exercise}[Alternate description of Hom-addition]
  \label{exe:HomAddition-In-Additive-Alternate}
  For morphisms $f$, $g\from A\to B$ in an additive category show that $\FoldOn{B}\Comp \PrdctMapInto{f,g} = f+g = \SumMapOutOf{f,g}\Comp \DgnlOn{A}$.
\end{exercise}

\begin{exercise}[Addition of composites]
  \label{exe:Addition-Composites}
  In an additive category, consider morphisms
  \begin{equation*}
    u\from A\to B,\quad v\from A\to C \qquad \text{and}\qquad f\from B\to D,\quad g\from C\to D
  \end{equation*}
  Show that $(fu)+(gv) = \SumMapOutOf{f,g}\Comp \PrdctMapInto{u,v}$.
\end{exercise}
\end{exercises}
\newpage
\section[Preabelian Categories]{Preabelian Categories}
\label{sec:PreAbelianCategories}

By definition (Ab.\ref{def:PreAbelianCat}), a category is preabelian if it is enriched in abelian groups, admits finite sums, and has kernels and cokernels. With (\ref{thm:(Co)Equalizers-AdditiveCat}) we conclude that it is finitely bicomplete.

Equivalently, a category is preabelian if and only if it is additive and z-exact. Thus every morphism admits a normal mono factorization and a normal epi factorization. So, the notion of (short) exact sequence may be adopted from Section \ref{sec:ExactSeqs}. %
\index{preabelian category}\index{category!preabelian}%

Here is how we can recognize preabelian categories among preadditive categories:

\begin{proposition}[Recognizing a preabelian category among preadditive ones]
  \label{thm:PreAbelianCategory-Recognize}%
  For a preadditive category $\Ctgry{X}$, the following conditions are equivalent: %
  \index{preabelian category!recognize}%
  \begin{tfae}
    \item $\Ctgry{X}$ is preabelian;
    \item $\Ctgry{X}$ has a zero object and is finitely bicomplete;
    \item $\Ctgry{X}$ has a zero object and admits arbitrary pullbacks and pushouts.
  \end{tfae}
\end{proposition}
\begin{proof}
  (II) implies (I), (III) follows from the definitions. To see that (I) implies (II), recall the general fact that a category with finite sums and coequalizers is finitely cocomplete. In the case at hand, a preabelian category has finite products and coproducts. Further, it admit (co)equalizers via Proposition~\ref{thm:(Co)Equalizers-AdditiveCat}. - Similar reasoning shows that (III) implies~(II).
\end{proof}

\begin{proposition}[Recognizing mono/epi-morphisms in a preabelian category]%
  \label{thm:PreabelianCat-Mono/Epi-Recognize}%
  In a preabelian category $\Ctgry{X}$, a morphism $f$ has the following properties:
  \begin{enumerate}[(i)]
    \item $f$ is a monomorphism if and only if $\KerMap{f}= (0\to A)$.
    \item $f$ is an epimorphism if and only if $\CoKerMap{f} = (B\to 0)$.
  \end{enumerate}
\end{proposition}
\begin{proof}
  In the preabelian category $\Ctgry{X}$, $\KerMap{f}$ exists. So, (i) follows from (\ref{thm:AdditiveCat-Mono/Epi-Recognize}). Dual reasoning proves (ii).
\end{proof}

\subsection[Related categorical structures]{Related categorical structures}
\label{subsec:PreAbelianCat-RelatedStructures}%

We know already that an additive category satisfies the \KSGInline-condition. Thus, a preabelian category is protomodular (\ref{thm:Protomodular<->KSG}). Noting that additive and preabelian categories are self-dual notions, we conclude that additive categories also satisfy the co\KSGInline-condition, and that preabelian categories are also coprotomodular.

The following result provides a partial converse in the context of homological (\ref{sec:HomologicalCats-Axioms}) linear categories:

\begin{theorem}[Recognizing preabelian categories among homological ones\HTag]
  \label{thm:HomologicalPlusLinearIsAdditive}
  \label{thm:Homological+Cndtn->PreAbelian}%
  A homological satisfying any of the conditions below is preabelian.
  \begin{thmlist}
    \item $\Ctgry{X}$ is enriched in the category $\Magmas$ of unitary magmas.
    \item $\Ctgry{X}$ is linear.
  \end{thmlist}
\end{theorem}
\begin{proof}
  (i) If $\Ctgry{X}$ is enriched in $\Magmas$, then every object admits a unitary magma structure by (\ref{thm:MagmaHom-sets->InternalMagma}.i). Thus every object is an internal abelian group by (\ref{thm:AbelianGroupObject-Recognize-PSA}). With (\ref{thm:MagmaHom-sets->InternalMagma}.iii) we see that $\Ctgry{X}$ is enriched in $\AbGrps$. So, $\Ctgry{X}$ is preadditive and finitely bicomplete. It is preabelian by~(\ref{thm:PreAbelianCategory-Recognize}).

  (ii) follows from (i) because every linear category is enriched in the category $\CMon$ of commutative monoids and, hence, in $\Magmas$.
\end{proof}

\begin{exercises}
\begin{exercise}[Split extension in a preabelian category - recognize]%
  \label{exe:SplitExtensionInPreAbCat-Recognize}%
  In a preabelian category, consider a short exact sequence $K\NMono X\NEpi Y$. Prove that if either $K\NMono X$ is a retractable monomorphism or if $X\NEpi Y$ is a sectionable epimorphism, then $X\cong K\oplus Y$. In either case, there is an isomorphism of short exact sequences
  \begin{equation*}
    \xymatrix@R=5ex@C=4em{
    K \ar@{{ |>}->}[r] \ar@{=}[d] &
    X \ar@{-{ >>}}[r] \ar[d]^{\cong} &
    Y \ar@{=}[d] \\
    K \ar@{{ |>}->}[r]_-{\PrdctMapInto{\IdMapOn{K},\ZeroMap}} &
    K\oplus Y \ar@{-{ >>}}[r]_-{\PrjctnOnto{Y}} &
    Y
    }
  \end{equation*}
\end{exercise}
\end{exercises}
\newpage
\section{Abelian Categories}
\label{sec:AbelianCategories}

Building on previous sections, we present a view of abelian categories. In the process, we take advantage of the relationship with homological and semiabelian categories. We refer the reader to \cite{Tohoku,Freyd,Alligators,Weibel,PAluffi2009} for a self-contained development of the foundations of abelian categories.

\begin{definition}[Abelian category]
  \label{def:AbelianCategory}%
  A preabelian category $\Ctgry{X}$ is \Defn{abelian} if every monomorphism is a kernel (condition (P.\ref{Property:MonoIsNormal})) and every epimorphism is a cokernel (condition (P.\ref{Property:EpiIsNormal})). %
  \index{abelian!category}\label{category!abelian}
\end{definition}

\begin{lemma}[Pullback of normal epi is pushout\AbTag]\label{thm:PbNormalEpiIsPo}
  In an abelian category, any pullback of a normal epimorphism is also a pushout; any pushout of a normal monomorphism is also a pullback.
\end{lemma}
\begin{proof}
  We consider a commutative diagram (C) as in (\ref{thm:AbelianCategory-PullbackPreserves/ReflectsMono/PushoutPreserves/ReflectsEpis}). If (C) is a pullback and $u$ is a normal epimorphism, then the associated sequence $0\to A\longrightarrow U\oplus V \longrightarrow Z \to 0$ is short exact, by item (ii) in (\ref{exe:Epimorphisms-Composite}) combined with the hypothesis (P.\ref{Property:EpiIsNormal}) that all epimorphisms are normal. It follows that (C) is both a pullback and a pushout diagram. The second statement follows by duality.
\end{proof}

\begin{proposition}[Pullback-stability of normal epimorphisms\AbTag]
  \label{thm:AbelianCategory-NormalEpiPullbackStable-NormalMonoPushoutStable}%
  In an abelian category, normal epimorphisms are pullback-stable, and normal monomorphisms are pushout-stable: properties (P.\ref{Property:NormalEpisClosedUnderBaseChange}) and (P.\ref{Property:NormalMonosClosedUnderCobaseChange}) hold.
\end{proposition}
\begin{proof}
  We consider a commutative diagram (C) as in (\ref{thm:AbelianCategory-PullbackPreserves/ReflectsMono/PushoutPreserves/ReflectsEpis}). If $u$ is a normal epimorphism, then (\ref{thm:Pushout->IsoOfCoKers}) tells us that $\CoKer{b} = \CoKer{u} = 0$. Thus $b$ is a epimorphism by (\ref{thm:PreabelianCat-Mono/Epi-Recognize}), and a normal epimorphism by the assumption (P.\ref{Property:EpiIsNormal}) that all epimorphisms are normal. The second statement follows by duality.
\end{proof}

\begin{theorem}[Abelian implies homological\AbTag]\label{thm:LAAthenHomological}
  Any abelian category is homological and co-homological.
\end{theorem}
\begin{proof}
  The given category is pointed and bicomplete by (\ref{thm:PreAbelianCategory-Recognize}). Condition \KSGInline\ holds by (\ref{thm:Additive->KSG}), while \PNEInline\ follows from (\ref{thm:AbelianCategory-NormalEpiPullbackStable-NormalMonoPushoutStable}).
\end{proof}

\begin{theorem}[Every map is normal\AbTag]
  \label{thm:EveryMapNormal-Ab}
  In an abelian category every morphism is a normal map.
\end{theorem}
\begin{proof}
  We know that any abelian category is homological by (\ref{thm:LAAthenHomological}). So given $f\from A\to B$, we may consider its image factorization $f=me$, where $m$ is a monomorphism and $e$ is a normal epimorphism. Now by (P.\ref{Property:MonoIsNormal}), $m$ is normal as well.
\end{proof}

\begin{theorem}[Abelian implies semiabelian\AbTag]
  \label{thm:Abelian->Semiabelian}%
  Any abelian category is semiabelian and co-semiabelian. %
  \index{abelian category!is semiabelian}%
\end{theorem}
\begin{proof}
  The given category is homological by Theorem \ref{thm:LAAthenHomological}. Since every map in $\Ctgry{X}$ is normal (\ref{thm:EveryMapNormal-Ab}), the \ANNInline-condition is satisfied. So, $\Ctgry{X}$ is semiabelian (\ref{def:SACategory}). Further, $\Ctgry{X}$ is co-semiabelian, because an abelian category is also co-abelian.
\end{proof}

\begin{exercises}
\begin{exercise}[Isomorphism recognition in abelian categories]
  \label{exe:IsomorphismRecognize-Abelian}
  In an abelian category show that a morphism $f\from A\to B$ is an isomorphism if and only if it is an epimorphism and a monomorphism. %
  \index{isomorphism!recognition in abelian category}%
\end{exercise}

\begin{exercise}[Formal difference in abelian categories]
  \label{exe:FormalDifferenceInAb}
  Recall the construction of the difference object $D(X)$ in normal categories (\ref{def:DifferenceObject}), and the associated formal difference of maps $f$, $g\from X\to Y$. In an abelian category $\Ctgry{X}$ show the following:
  \begin{thmlist}
    \item The difference object functor $D$ is naturally equivalent to the identity functor on $\Ctgry{X}$.
    \item Via the natural equivalence in (i), the formal difference $f- g$ agrees with the difference given by the abelian group enrichment of $\Ctgry{X}$.
  \end{thmlist}
\end{exercise}
\end{exercises}
\section[Related categorical structures]{Related categorical structures}
\label{sec:RelatedCategoricalStructurs-AbelianCats}

In this section, we elaborate on the relationship between abelian categories and the categorical structures we developed earlier, notably p-exact, homological, and semiabelian categories, as well as categories which are enriched in abelian groups. As a refinement, we insert the following intermediate concepts: \emph{almost abelian categories}~\cite{Rump0}, \emph{quasi-abelian categories}~\cite{Yoneda-Exact-Sequences,Rump0,MR1779315}, and the corresponding one-sided notions.

\begin{definition}[Left almost abelian category]
  \label{def:LeftAlmostAbelianCat}%
  A preabelian category is \Defn{left almost abelian} or \Defn{left quasi-abelian} if pullbacks preserve normal epimorphisms; Property (P.\ref{Property:NormalEpisClosedUnderBaseChange}) holds. %
  \index{left!almost abelian category}\index{left!quasi-abelian category}%
\end{definition}

Thus, a preabelian category is left almost abelian if it satisfies the \PNEInline-condition. --- Dually:

\begin{definition}[Right almost abelian category]
  \label{def:RightAlmostAbelianCat}%
  A preabelian category is \Defn{right almost abelian} or \Defn{right quasi-abelian} if pushouts preserve normal monomorphisms; Property (P.\ref{Property:NormalMonosClosedUnderCobaseChange}) holds. %
  \index{right!almost abelian category}\index{right!quasi-abelian category}%
\end{definition}

\begin{definition}[Almost abelian category]
  \label{def:AlmostAbelianCat}%
  A preabelian category is \Defn{almost abelian} or \Defn{quasi-abelian} if it is both left and right almost abelian. %
  \index{almost abelian category}\index{quasi-abelian category}\index{category!almost abelian}%
  \index{category!quasi abelian}%
\end{definition}

From (\ref{thm:AbelianCategory-NormalEpiPullbackStable-NormalMonoPushoutStable}), it follows immediately that every abelian category is almost abelian:

\begin{lemma}[Abelian implies almost abelian]
  \label{thm:Abelian->AlmostAbelian}%
  Every abelian category is an almost abelian category.\NoProof
\end{lemma}

\begin{example}[Examples of almost abelian categories]
  \label{exa:AlmostAbelianCats}%
  The categories of topological abelian groups, locally compact abelian groups, topological vector spaces, normed spaces, Banach spaces, locally convex spaces and Fréchet spaces are almost abelian---see~\cite{Rump0}. %
  \index{almost!abelian categories - examples}%
\end{example}

\begin{theorem}[Homological category vs.\ left almost abelian category]
  \label{thm:LeftAlmostAbelian-Homological}%
  A preabelian category is left almost abelian if and only if it is homological. Further, a homological category is left almost abelian if and only if it is additive.
\end{theorem}
\begin{proof}
  By definition, a homological category satisfies the \PNEInline-condition. Therefore, any homological preabelian category is left almost abelian.

  To see that a left almost abelian category is homological, recall first that any preabelian category has a zero object and is finitely bicomplete, see Proposition~\ref{thm:PreAbelianCategory-Recognize}. By definition of a left almost abelian category, the \PNEInline-condition is satisfied. The \KSGInline-condition holds in any additive category by Proposition~\ref{thm:Additive->KSG} and, hence, in any preabelian category.

  Now consider a homological category. If it is left almost abelian, then it is additive by definition. If it is additive, then it is preabelian because it is pointed and finitely bicomplete. Finally, pullbacks preserve normal epimorphisms by the \PNEInline-condition of homological categories. So, it is left almost abelian.
\end{proof}

By combining (\ref{thm:LeftAlmostAbelian-Homological}) and its dual, we see that an almost abelian category is both homological and co-homological. A result due to G.~Janelidze~\cite{JRosickyWTholen2007} shows that the converse also holds.

\begin{theorem}[Almost abelian = homological + co-homological]
  \label{thm:AlmostAbelianIffHomoPlusCoHomo}%
  A category is almost abelian if and only if it is both homological and co-homological.
\end{theorem}
\begin{proof}
  By (\ref{thm:LeftAlmostAbelian-Homological}), it suffices to prove that a category $\Ctgry{X}$ which is both homological and co-homological is additive. Given two objects $A$ and $B$, consider this diagram in $\Ctgry{X}$:
  \begin{equation*}
    \xymatrix@R=5ex@C=4em{
    A \ar@<+.5ex>[r]^-{\InclsnOf{A}} &
    A+ B \ar@<+.5ex>[r]^-{\SumMapOutOf{\ZeroMap , \IdMapOn{B}}} \ar@<+.5ex>[l]^-{\SumMapOutOf{\IdMapOn{A},\ZeroMap}} &
    B \ar@<+.5ex>[l]^-{\InclsnOf{B}}
    }
  \end{equation*}
  In the homological $\Ctgry{X}^{\op}$, it is the product decomposition diagram from (\ref{thm:ProductRecognition}). Thus the diagram consists of two short exact sequences. The property of being short exact is self dual. So, the diagram, as constructed in $\Ctgry{X}^{\op}$ consists of two short exact sequences. Via (\ref{thm:ProductRecognition}), the middle object is the product of the end objects. So, the diagram above is a biproduct diagram. This means that the category $\Ctgry{X}$ is linear and $\Ctgry{X}$ homological as well. So, it is preabelian by (\ref{thm:Homological+Cndtn->PreAbelian}). Since is is homological and co-homological, it is almost abelian.
\end{proof}

We now consider characterizations of abelianness in relation to additive, p-exact, homological, and semiabelian categories.

\begin{theorem}[Abelian is additive plus p-exact]
  \label{thm:AdditivePlusp-exactIsAbelian}%
  A category is abelian if and only if it is p-exact and additive. %
  \index{abelian category!if and only if p-exact and additive}%
\end{theorem}
\begin{proof}
  An abelian category is additive by definition. Maps are normal by Theorem~\ref{thm:EveryMapNormal-Ab}. So, an abelian category is also p-exact.

  Conversely, given an additive p-exact category, we know that (P.\ref{Property:MonoIsNormal}) and (P.\ref{Property:EpiIsNormal}) hold. So, it has kernels and cokernels by  or by Exercise \ref{exe:PuppeExactCat-EquivalentCharacterization}. So, it is abelian.
\end{proof}

\begin{theorem}[Abelian is additive plus di-exact]
  \label{thm:AdditivePlusDiexactIsAbelian}%
  A category is abelian if and only if it is di-exact and additive. %
  \index{abelian category!if and only if p-exact and diexact}%
\end{theorem}
\begin{proof}
  An abelian category is additive by definition. It is di-exact by (\ref{thm:EveryMapNormal-Ab}).

  Conversely, in a pointed additive category, a morphism $f\from X\to Y$ may be factored as the normal monomorphism $\PrdctMapInto{\IdMapOn{X},0}\colon X\to X\oplus Y$ followed by the normal epimorphism $\MapOutOf{f, \IdMapOn{Y}}\colon X\oplus Y\to Y$. So, $f$ is normal by the \ANNInline-condition. Hence the given category is both p-exact and additive, hence abelian by (\ref{thm:AdditivePlusp-exactIsAbelian}).
\end{proof}

\begin{corollary}[Abelian is additive plus semiabelian]
  \label{thm:AdditivePlusSemiabelianIsAbelian}%
  A category is abelian if and only if it is semiabelian and additive.
\end{corollary}
\begin{proof}
  If a category $\Ctgry{X}$ is abelian, then it is additive by definition. It is semiabelian by (\ref{thm:Abelian->Semiabelian}). Conversely, if $\Ctgry{X}$ is semiabelian and additive, then it is di-exact and additive, hence abelian by (\ref{thm:AdditivePlusDiexactIsAbelian}).
\end{proof}

\begin{corollary}[Abelian is semiabelian plus co-semiabelian]\label{thm:SemiabelianPlusCoSemiabelianIsAbelian}
  A category is abelian if and only if it is semiabelian and co-semiabelian.
\end{corollary}
\begin{proof}
  If a category $\Ctgry{X}$ is abelian, then it is semiabelian and co-semiabelian by (\ref{thm:Abelian->Semiabelian}).

  Conversely, if $\Ctgry{X}$ is semiabelian and co-semiabelian then it is almost abelian by (\ref{thm:AlmostAbelianIffHomoPlusCoHomo}). In particular, it is additive, so it is abelian by (\ref{thm:AdditivePlusSemiabelianIsAbelian}).
\end{proof}

\begin{theorem}[Abelian is additive plus Barr-exact]
  \label{thm:AdditivePlusBarrExactIsAbelian}%
  A category is abelian if and only if it is Barr-exact and additive.
\end{theorem}
\begin{proof}
  By (\ref{thm:AdditivePlusSemiabelianIsAbelian}), a category is abelian if and only if semiabelian and additive. By (\ref{thm:ReflexiveRelationInSACategory}), any semiabelian category is Barr-exact.

  Conversely, by definition, any additive category $\Ctgry{X}$ is pointed with binary products and coproducts. Furthermore (\ref{thm:Additive->KSG}), \KSGInline\ holds in $\Ctgry{X}$. Barr-exactness adds finite completeness. Hence, being pointed, regular and protomodular (\ref{thm:Protomodular<->KSG}), the category~$\Ctgry{X}$ is BB-homological by (\ref{thm:BBHom}). On the other hand, by (\ref{thm:ProtoMaltsev}) and (\ref{thm:CharMaltsev}), the Mal'tsev property holds for $\Ctgry{X}$. So $\Ctgry{X}$ has all finite colimits (\ref{thm:Colimits}). Hence $\Ctgry{X}$ is an additive homological category. Di-exactness follows from (\ref{thm:BarrExactimplies(Ax.6)}), and abelianness from (\ref{thm:AdditivePlusDiexactIsAbelian}).
\end{proof}

This result is often called the \Defn{Tierney Equation} after Myles Tierney who discovered it.

\begin{example}[Almost abelian but not abelian]
  \label{exa:AlmostAbelian}%
  In \cite[Corollary in Section 4]{Rump0}, it is shown that a category is almost abelian if and only if it is the torsion or torsion-free class in an abelian category. For instance, the categories of torsion and torsion-free abelian groups are almost abelian categories which are not Barr-exact. %
  \index{almost!abelian categories - examples}%
\end{example}

In view of (\ref{thm:AdditivePlusBarrExactIsAbelian}), a variety of algebras is abelian if and only if it is an additive category. The following more precise analysis appears in~\cite{Freyd}:

\begin{example}[Abelian varieties of algebras]%
  \label{exa:AbelianVarietiesOfAlgebras}%
  A variety of algebras is an abelian category if and only if it is equivalent to a category of modules over a ring $R$. If $F(1)$ denotes the free object on a single generator, then $R$ is the ring of endomorphisms $\EndRng{F(1)}$ of $F(1)$. %
  \index{abelian!varieties of algebras}%
  \NoProof
\end{example}

\begin{subordinate}{Palamodov and Raïkov semiabelian categories}
  In the literature, the term `semiabelian category' is used to refer to at least three non-equivalent kinds of structure:
  \begin{ulist}
    \item Semiabelian categories in the sense of Janelidze--Márki--Tholen \cite{Janelidze-Marki-Tholen} are discussed here in Chapter~\ref{chap:SACats}.
    \item Semiabelian categories in the sense of Raïkov~\cite{Raikov:Semi-Abelian} are the same as the almost abelian categories of Definition~\ref{def:AlmostAbelianCat}. %
    \index{Raïkov!semiabelian category}\index{semiabelian category!in sense of Raïkov}%
    \item Semiabelian categories in the sense of Palamodov~\cite{VPPalamodov1968}, see also~\cite{VPPalamodov1971,Gruson,Rump0}. The full definition is given below; see (\ref{def:Palamodov-Semiabelian}). Any almost abelian category is Palamodov-semiabelian~(\ref{thm:AlmostAbelianImpliesPSA}). But, as explained in~\cite{Rump}, the converse does not hold. For example, the category of bornological spaces is Palamodov-semiabelian but not almost abelian~\cite{JBonetSDierolf2005-Bornological}.
  \end{ulist}

  For context for the definition of Palamodov-semiabelian categories, recall that, in a \ZExact\ category, every morphism $f$ admits a normal epi factorization and a normal mono factorization whose factoring objects are uniquely related by a unique comparison map
  \begin{equation*}
    \NENMComp{f}\from \CoKer{\KerMap{f}} \to \Ker{\CoKerMap{f}}
  \end{equation*}

  \begin{definition}[Palamodov-semiabelian category]
    \label{def:Palamodov-Semiabelian}%
    A preabelian category is \emph{Palamodov-semiabelian} if for every morphism  $f\from X\to Y$, the comparison map $\NENMComp{f}$ is simultaneously a monomorphism and an epimorphism. %
    \index{Palamodov!semiabelian category}\index{semiabelian category!in sense of Palamodov}%
  \end{definition}

  This is clearly equivalent to asking that every morphism in the given category has both an image and a coimage factorization (see (\ref{def:ImageFactorization})). Since homological categories have images, while co-homological categories have coimages, any almost abelian category is Palamodov-semiabelian by (\ref{thm:AlmostAbelianIffHomoPlusCoHomo}).

  \begin{theorem}[(Almost) abelian implies Palamodov-semiabelian]\label{thm:AlmostAbelianImpliesPSA}
    Any (almost) abelian category is Palamodov-semiabelian.\NoProof
  \end{theorem}

  \begin{proposition}[Images vs.\ pullback of normal epi is epi]
    In a preabelian category, pullbacks send normal epimorphisms to epimorphisms if and only if every morphism admits an image factorization.
  \end{proposition}
  \begin{proof}
    First suppose that pullbacks send normal epimorphisms to epimorphisms. We show that any morphism $f$ admits a decomposition $f=v e$, with $e$ a normal epimorphism and $v$ a monomorphism. Consider its normal epi factorization as in this diagram:
    \begin{equation*}
      \xymatrix@R=5ex@C=4em{
      \Ker{f} \ar@{{ |>}->}[r]^-{k} \ar@{-{>>}}[d]_{\bar{e}} \PullLU{rd} &
      X \ar[r]^-{f} \ar@{-{ >>}}[d]^{e} &
      Y \ar@{=}[d] \\
      \Ker{v} \ar@{{ |>}->}[r]_-{l} &
      \CoKer{k} \ar[r]_-{v} &
      Y
      }
    \end{equation*}
    We know that $l\bar{e}=ek=\ZeroMap$. So, $\Ker{v}=\ZeroObject$ because $\bar{e}$ is epic, as a pullback (\ref{thm:PullbackRecognition-KernelSide-1}) of the normal epimorphism $e$. Hence $v$ is monic by (\ref{thm:PreAdditiveCat-Mono/Epi-Recognize}.i).

    For the converse, we first prove that when image factorizations exist, any pullback of a normal epimorphism is also a pushout (cf.\ (\ref{thm:PbNormalEpiIsPo})). Consider a diagram (C) such as in (\ref{thm:AdditiveCategory-Pullback/Pushout-Recognition}) and assume that $u$ is a normal epimorphism. Let $r\colon P\to Z$ be the induced comparison from the pushout of $a$ and $b$ to $Z$. Taking $f$ to be $\SumMapOutOf{u,-v}$, we have that $\NENMComp{f}$ and $r$ coincide by points (ii) and (iii) in (\ref{thm:AdditiveCategory-Pullback/Pushout-Recognition}) combined. Now $r=\NENMComp{f}$ is a monomorphism by hypothesis, though which the strong epimorphism $u$ factors. Hence $r$ is an isomorphism, so that (C) is a pushout.

    Now it follows from (\ref{thm:Pushout->IsoOfCoKers}) that a pullback of any normal epimorphism has a trivial cokernel, so that by (\ref{thm:PreabelianCat-Mono/Epi-Recognize}) it is an epimorphism.
  \end{proof}

  We see that Palamodov-semiabelian categories may be characterized in terms of a weak form of the conditions (P.\ref{Property:NormalMonosClosedUnderCobaseChange}) and (P.\ref{Property:NormalEpisClosedUnderBaseChange}).

  \begin{corollary}[Characterization of Palamodov-semiabelian categories]\label{thm:Palamodov-By-Stability}
    A preabelian category is Palamodov-semiabelian if and only if pushouts send normal monos to monos and pullbacks send normal epis to epis.\NoProof
  \end{corollary}
\end{subordinate}

\begin{exercises}

\begin{exercise}[Alternate verification of \KSGInline\ in an abelian category]
  Show that any abelian category satisfies \KSGInline\ by showing that $\CoKer{m}=0$ in the diagram in the proof of (\ref{thm:AddThenProto}); i.e.\ $m$ is a monomorphism and an epimorphism.
\end{exercise}

\begin{exercise}[Normality check in pointed additive context]
  Check the claim of normality for the two morphisms in (\ref{thm:AdditivePlusDiexactIsAbelian}).
\end{exercise}
\end{exercises}
\section[Categories of abelian objects]{Categories of abelian objects}
\label{sec:CatsOfAbelianObjects}%

When considering abelian group objects in a finitely complete category $\Ctgry{X}$, we encounter two quite different scenarios:
\begin{ulist}
  \item Category $\Ctgry{X}$ is such that an object may carry two or more distinct abelian group structures. This happens, for example in the category $\Sets$ of sets. The category of all abelian group objects in $\Sets$ is the category of abelian groups. In general, the collection of all internal abelian group objects in $\Ctgry{X}$ is an additive category; see (\ref{thm:AbelianGroupObjects->AdditiveCat}).
  \item Category $\Ctgry{X}$ is such that if an object $X$ carries an abelian group structure, then this structure is unique and is a property of $X$. In this case the collection of all abelian group objects in $\Ctgry{X}$ determines a full and replete subcategory of $\Ctgry{X}$ which we call the \Defn{abelian core} $\AbCoreOf{X}$ of $\Ctgry{X}$. This happens, for example in every homological category, where the abelian core consists of all abelian objects. Below, we develop basic properties of $\AbCoreOf{X}$.
\end{ulist}

\begin{proposition}[Abelian group objects form additive category]
  \label{thm:AbelianGroupObjects->AdditiveCat}%
  Given a category $\Ctgry{X}$ with finite products, the collection of abelian group objects $\AbGrps(\Ctgry{X})$ in $\Ctgry{X}$ forms an additive category.
\end{proposition}
\begin{proof}
  With Theorem~\ref{thm:MagmaHom-sets->InternalMagma} we see that $\AbGrps(\Ctgry{X})$ is enriched in $\AbGrps$. With (\ref{thm:BiProduct-From-(Co)Product}), we see that $\Ctgry{X}$ has finite biproducts, hence is an additive category.
\end{proof}

In the remainder of this section we consider abelian group objects in a homological or in a semiabelian category $\Ctgry{X}$. From Section \ref{sec:CommutativeObjects-H} we already know that, (a) whenever an object $X$ in $\Ctgry{X}$ admits the structure of an internal unitary magma, then this structure is unique, and (b) this magma structure on $X$ is commutative, associative, and associated with a unique inverse operation, so that it is actually an internal abelian group object. In (\ref{sec:CommutativeObjects-H}) we called such $X$ an \emph{abelian object}.

\begin{theorem}[Abelian objects in a homological category\HTag]
  \label{thm:AbGrp-Objects-H}
  In a homological category $\Ctgry{X}$ the abelian objects form a left almost abelian category $\AbCoreOf{X}$. It is full and replete in $\Ctgry{X}$, and is closed under subobjects and quotients.
\end{theorem}
\begin{proof}
  We begin by showing that $\AbCoreOf{X}$ is full in $\Ctgry{X}$. So, let $(X,\mu)$ and $(Y,\nu)$ be abelian group objects, and let $f\from X\to Y$ be a morphism in $\Ctgry{X}$. To see that $f$ is also a  morphism $(X,\mu)\to (Y,\nu)$ consider the diagram below.
  
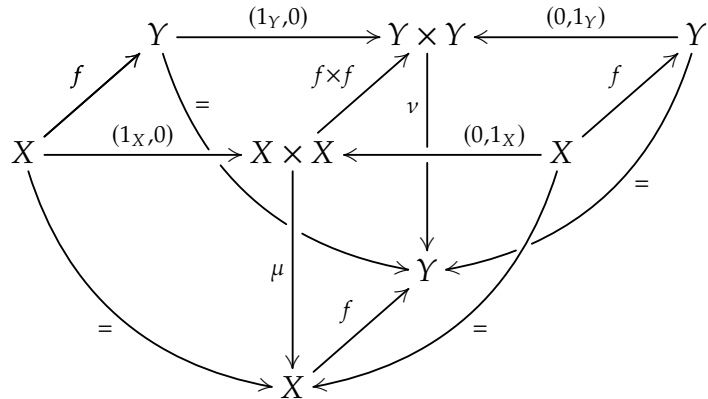
\begin{equation*}
    \xymatrix@!@R=1ex@C=1em{
    & Y \ar[rr]^-{(\IdMapOn{Y},0)} \ar@/_4ex/[rrdd]|(0.34)\hole|(0.64)\hole^(0.2){=} &&
    \Prdct{Y}{Y} \ar[dd]|\hole_(0.3){\nu}&&
    Y \ar[ll]_-{(0,\IdMapOn{Y})} \ar@/^4ex/[ddll]|(0.75)\hole^(0.4){=}\\
    X \ar[ru]^{f} \ar[rr]^-{(\IdMapOn{X},0)} \ar@/_4ex/[ddrr]_{=} \ar[ru]^{f} &&
    \Prdct{X}{X} \ar[ru]^{f\times f} \ar[dd]_{\mu} &&
    X \ar[ru]^{f} \ar[ll]_(0.25){(0,\IdMapOn{X})} \ar@/^4ex/[ddll]^{=}\\
    &&& Y \\
    && X \ar[ru]^{f}
    }
  \end{equation*}
  The top squares and the front/back facing triangles commute. The identity $\nu\Comp (\Prdct{f}{f}) = f\mu$ follow via checking on the factors of $\Prdct{X}{X}$. This means that $f$ is a morphism of internal abelian group objects; i.e.\ $\AbCoreOf{X}$ is full in $\Ctgry{X}$.

  That $\AbCoreOf{X}$ is closed under subobjects follows from (\ref{thm:CommutingMaps-Props}.ii). It remains to show that $\AbCoreOf{X}$ is closed under quotients. Given $(A,\mu)$ and a normal epimorphism $q\from A\NEpi Q$ in $\Ctgry{X}$, let $\kappa\from K\to A$ be the kernel of $q$. By (\ref{thm:Product-SESs}), we obtain this diagram of short exact sequences
  \begin{equation*}
    \xymatrix@R=5ex@C=4em{
    \Prdct{K}{K} \ar@{{ |>}->}[r]^-{\kappa\times \kappa} \ar[d]_{\sigma} &
    \Prdct{A}{A} \ar@{-{ >>}}[r]^-{q\times q} \ar[d]^-{\mu} &
    \Prdct{Q}{Q} \ar@{.>}[d]^{\nu} \\
    K \ar@{{ |>}->}[r]_-{\kappa} &
    A \ar@{-{ >>}}[r]_-{q} &
    Q
    }
  \end{equation*}
  The left hand square commutes because the subobject $\kappa$ of $A$ belongs to the full subcategory $\AbCoreOf{X}$.  Then $q\mu(\Prdct{\kappa}{\kappa})=\ZeroMap$. Thus there is a map $\nu\from \Prdct{q}{q}\to Q$ unique with $q\mu = \nu(\Prdct{q}{q})$. By (\ref{thm:CommutingMaps-Props}.iii) $q\mu$ is the multiplication of $q$ by $q$. Using $q=\IdMapOn{Q}\Comp q$ in (\ref{thm:NormalQuotientOfCommutingMaps}) shows that $\nu$ is the multiplication of $\IdMapOn{Q}$ by itself.

  It remains to show that $\AbCoreOf{X}$ is finitely bicomplete. Finite products exist by (\ref{thm:CommutingMaps-Props}.iv), and these are biproducts by (\ref{thm:BiProduct-From-(Co)Product-Additive}). We conclude that $\AbCoreOf{X}$ has coequalizers because it is closed under quotients. So $\AbCoreOf{X}$ is finitely cocomplete. Further $\AbCoreOf{X}$ has equalizers: If $f\from A\to B$ is a morphism in $\AbCoreOf{X}$, then its kernel $\kappa\from K\to A$ in $\Ctgry{X}$ is a subobject of $A$, hence belongs to $\AbCoreOf{X}$. The defining property of `kernel' tells us that $\kappa$ is the kernel of~$f$ in $\AbCoreOf{X}$ as well. So, $\AbCoreOf{X}$ is finitely complete.

  Pullbacks and normal epimorphisms in $\AbCoreOf{\Ctgry{X}}$ are as in $\Ctgry{X}$. We know that, in $\Ctgry{X}$, pullbacks preserve normal epimorphisms. So, this property also holds in $\AbCoreOf{\Ctgry{X}}$. Hence $\AbCoreOf{\Ctgry{X}}$ satisfies condition (P.\ref{def:LAACat}). So, it is left almost abelian.
\end{proof}

\begin{corollary}[Recognizing (left almost) abelian categories\HTag]
  \label{thm:PSA-Cat-Recognize}
  \label{thm:AbelianViaSemiabelianWithNormalMonos}
  \label{thm:AbelianCat-Recognize}
  For a homological category $\Ctgry{X}$ the following hold:
  \begin{enumerate}[(i)]
    \item $\Ctgry{X}$ is left almost abelian if and only if the diagonal map $\DgnlOn{X}\from X\to \Prdct{X}{X}$ is a normal monomorphism for every object $X$;
    \item $\Ctgry{X}$ is abelian if and only if every monomorphism is normal.
  \end{enumerate}
\end{corollary}
\begin{proof}
  (i) By Proposition~\ref{thm:AbGroupViaNormalDiagonal}, the diagonal $\DgnlOn{X}$ is normal precisely when $X$ is an internal abelian group object of $\Ctgry{X}$. So $\Ctgry{X}=\AbCoreOf{X}$ when this is true for all objects of $\Ctgry{X}$. By Theorem~\ref{thm:AbGrp-Objects-H}, this is equivalent to $\Ctgry{X}$ being left almost abelian.

  (ii) If every monomorphism is normal, then $\Ctgry{X}$ is left almost abelian by (i), and we only need to show that every epimorphism is normal. Given an epimorphism $f\from A\to B$ in $\Ctgry{X}$, consider its image factorization:
  \begin{equation*}
    \xymatrix@R=5ex@C=4em{
    A \ar@{-{ >>}}[r]_-{q} \ar@{-{>>}}@/^2ex/@<0.5ex>[rr]^-{f} &
    I \ar@{{ >}->}[r]_-{m} &
    B
    }
  \end{equation*}
  We just saw that the monomorphism $m$ is normal. It is an epimorphism by (\ref{exe:Epimorphisms-Composite}.ii), hence an isomorphism by (\ref{thm:IsomorphismRecognition}.iii). So, $f$ is normal by (\ref{thm:NormalEpi-Props}.ii). - This completes the proof.
\end{proof}

\begin{theorem}[Internal abelian group objects in a semiabelian category\SATag]
  \label{thm:AbGrp-Objects-SA}
  In a semiabelian category $\Ctgry{X}$, the abelian objects form an abelian category $\AbCoreOf{X}$. It is full and replete in $\Ctgry{X}$, and is closed under subobjects and quotients.
\end{theorem}
\begin{proof}
  The category $\Ctgry{X}$ is homological, so that by Corollary~\ref{thm:PSA-Cat-Recognize}, we only need to show that all monomorphisms in the left almost abelian category $\AbCoreOf{X}$ are normal. So, consider a monomorphism $m\from M\to X$. By (\ref{thm:SubDiagonal-AbelianObject}), the subdiagonal map $\PrdctMapInto{\IdMapOn{M},m}\from M\to \Prdct{M}{X}$ is a normal monomorphism in $\Ctgry{X}$. Then via the \ANNInline-condition, the commutative diagram
  \begin{equation*}
    \xymatrix@R=5ex@C=4em{
    M \ar@{{ |>}->}[r]^-{\PrdctMapInto{\IdMapOn{M},m}} \ar@{=}[d] &
    \Prdct{M}{X} \ar@{-{ >>}}[d]^{\PrjctnOnto{X}} \\
    M \ar@{{ >}->}[r]_-{m} &
    X
    }
  \end{equation*}
  in $\Ctgry{X}$ shows that $m$ is a normal monomorphism in $\Ctgry{X}$, so that it is a normal monomorphism in $\AbCoreOf{X}$ as well.
\end{proof}

\begin{example}[Examples of abelian cores]
  From (\ref{exa:CommutingSubgroups}) we learn that a group is an abelian object in $\Grps$ if and only if it is an abelian group. Hence $\AbCoreOf{\Grps}=\AbGrps$.

  Let $R$ be a ring and let $\Lie_R$ be the semiabelian variety of Lie algebras over $R$. Example~\ref{exa:LieCommutingSubalgebras} explains that a Lie algebra is abelian if and only if its bracket is trivial. Thus we find the category $\ModulesOver{R}$ of $R$-modules (equipped with a zero multiplication) as~$\Ab(\Lie_R)$.

  The case of associative algebras over a ring $R$ is very similar to the case of Lie algebras. An object is abelian  if and only if the algebra multiplication is zero, so that we find the category of $R$-modules.

  A loop which carries an internal unitary magma structure is the same thing as a unitary magma that carries an abelian group structure; in other words, a loop admits such structure if and only if it is an abelian group. Hence a loop is an abelian object in $\Loop$ if and only if its multiplication is at the same time associative and commutative. In other words, $\AbCoreOf{\Loop}=\AbGrps$.
\end{example}

Since any subobject of an abelian object in a semiabelian category is again an abelian object, abelianness of the full subcategory of all abelian objects tells us that such a subobject is automatically normal, which gives us the next result. We chose to add a direct proof to this argument, because it is interesting in its own right.

\begin{corollary}[Subobject of abelian object is normal\SATag]
  \label{thm:SubobjectAbelianObjectIsNormal}
  In a semiabelian category, any subobject of an abelian object is a normal subobject.
\end{corollary}
\begin{proof}
  Let $m\from M\to X$ represent a subobject of an abelian object $X$. We denote its abelian group structure $(X,\mu,i)$ and write $t=\mu\Comp (\IdMapOn{X} \prdct i)\from X\times X\to X$. Pulling back $t$ along~$m$, we find a commutative square, which we can complete to the diagram below.
  \begin{equation*}
    \xymatrix@R=5ex@C=4em{
    M \ar@{.>}[r]_-k \ar[d]_{m} \ar@/^1em/@{=}[rr] &
    R \ar[r]_-{\bar{t}} \ar[d]_r \PullLU{rd} &
    M \ar[d]^-{m} \\
    X \ar@/_1em/@{=}[rr] \ar[r]^-{(1_X,0)} &
    X\times X \ar[r]^-{t} &
    X
    }
  \end{equation*}
  The map $k$ comes via the universal property of the pullback. But then the left hand square is a pullback as well by (\ref{thm:Pullbacks,ConcatenatedSquares}). Then $k$ is a normal monomorphism because so is $(\IdMapOn{X},0)$.

  By definition, $t\Comp \DgnlOn{X} = 0$. So, the universal property of the pullback yields the morphism $(\DgnlOn{X},0)\from X\to R$, unique with $r(\DgnlOn{X},0)=\DgnlOn{X}$ and $\bar{t}(\DgnlOn{X},0)=0$.  Setting $r_1=\PrjctnOnto{1}\Comp r$ we see that
  \begin{equation*}
    r_1 (\DgnlOn{X},0) = \PrjctnOnto{1}r(\DgnlOn{X},0) = \PrjctnOnto{1}\DgnlOn{X}= \IdMapOn{X}
  \end{equation*}
  This means that $r_1$ is a split epimorphism, hence a normal epimorphism. Further,
  \begin{equation*}
    r_1k = \PrjctnOnto{1}r k = \PrjctnOnto{1}(\IdMapOn{X} m) = m = m\IdMapOn{M}
  \end{equation*}
  With the \ANNInline-condition we conclude that $m$ is a normal monomorphism. This was to be shown.
\end{proof}
\appendix
\part[Appendix]{Appendix}
\chapter[Background from Category Theory]{Background from Category Theory}
\label{chap:Categorical-Preliminaries}%

We assume that the reader is familiar with the concepts of category and functor, and has had contact with examples of categories such as the category $\Sets$ of sets, the category $\Grps$ of groups, etc.\ Further we assume some basic familiarity with the concepts of limit and colimit. For references on categories, we recommend \cite{SMacLane1998}, \cite{AHS:Cats}, and \cite{Borceux:Cats1,Borceux:Cats2}.

Here we present more specialized categorical concepts which will be used extensively. Also, we introduce notational conventions.
\section[Limits and Colimits]{Limits and Colimits}
\label{sec:Limits-CoLimits}%

A category $J$ is called \Defn{small} if its class of morphisms if a proper set. This implies that the object set of $J$ is proper. Similarly, $J$ is \Defn{finite} if its class of morphisms and, hence, the class of its objects is a finite set. Thus, every finite category is small. %
\index{category!small}\index{category!finite}%
\index{finite!category}\index{small!category}%

\subsection{Limits and colimits}
\label{subsec:Limits/Colimits}

Let $\EuScript{X}$ be an arbitrary category, and let $J$ be a small category. The limit of a functor $F\from J\to \EuScript{X}$ is given by a \Defn{limit cone} $\lambda\from \LimOfOver{F}{J}\Rightarrow F$ which we may visualize as the blue part of this commutative diagram in $\EuScript{X}$: %
\index[not]{l!$\LimOfOver{F}{J}$\IndSep limit of functor $F\from J\to \EuScript{X}$}%
\begin{equation*}
  \xymatrix@R=6ex@C=2em{
  & A \ar[d]_{\alpha_j} \ar[rrd]|\hole_(.3){\alpha_k} \ar@{.>}[rr]^{a} &&
  {\color{blue} \LimOfOver{F}{J}} \ar@[blue][lld]^(.3){\color{blue} \lambda_{j}} \ar@[blue][d]^{\color{blue} \lambda_k} \\
  {\color{blue} \cdots} & {\color{blue} F(j)} \ar@[blue][rr]_{\color{blue} F(u)} &&
  {\color{blue} F(k)} &
  {\color{blue} \cdots}
  }
\end{equation*}
The limit cone $\lambda$ is terminal among all cones $\alpha\from A\Rightarrow F$: Every such cone factors uniquely through $\LimOfOver{F}{J}$. This implies that the morphisms $\lambda_j$ of the limit cone are jointly monomorphic.

Dually, the colimit of $F$ is given by a cocone $\gamma\from F\Rightarrow \CoLimOfOver{F}{J}$ which is initial among all cocones $F\Rightarrow Z$. Consequently, the family of morphisms $\gamma_j$ of the colimit cocone is jointly epimorphic. %
\index[not]{c!$\CoLimOfOver{F}{J}$\IndSep colimit of functor $F\from J\to \EuScript{X}$}

We say that $\Ctgry{X}$ is (finitely) complete if all (finite) limits in $\Ctgry{X}$ exist. Similarly, $\Ctgry{X}$ is (finitely) cocomplete if all (finite) colimits in $\Ctgry{X}$ exist. A (finitely) bicomplete category is (finitely) complete and (finitely) cocomplete. %
\index{finitely!complete category}\index{finitely!cocomplete category}\index{finitely!bicomplete category}%
\index{complete category}\index{cocomplete category}\index{bicomplete category}%
\index{category!(finitely) complete}\index{category!(finitely) cocomplete}\index{category!(finitely) bicomplete}%
\label{InText:CompleteCategory}\label{InText:CoCompleteCategory}\label{InText:BiCompleteCategory}%

\subsection{Functoriality of limits and colimits}
While (co)limits are unique up to unique isomorphism, it is convenient to assume that functorial (co)limits are available in the following sense.

\begin{definition}[Functorial (co)limits$^{\ast}$]
  \label{def:Functorial(Co)Limits}
  Fix a small category $J$, and consider the embedding $\Phi\from \Ctgry{X}\to \Ctgry{X}^J$, which sends an object $X$ in $\Ctgry{X}$ to the constant functor $ J\to \Ctgry{X}$ with value $X$. Functorial limits over $J$ in $\Ctgry{X}$ are given by a right adjoint $ \Lim_J\from \Ctgry{X}^J\to \Ctgry{X}$ of $\Phi$, and functorial colimits over $J$ in $\EuScript{X}$ are given by a left adjoint $ \CoLim^J\from \Ctgry{X}^J\to \Ctgry{X}$ of $\Phi$. A (co)limit is called \Defn{finite} when $J$ is a finite category. %
  \index{finite!(co-)limit}%
\end{definition}

\subsection{Sum and product$^{\ast}$}\label{subsec:Sum-Product}%
If $J$ is a discrete category (every morphism is an identity) %
\index{discrete!category}\index{category!discrete}%
then a functor $F\from J\to \Ctgry{X} $ has
\begin{ulist}
  \item as its limit:  $\LimOfOver{F}{J}$ the product $\FamPrdct{j\in J}{F(j)}$; %
  \index{product!as limit of discrete diagram}\index[not]{p!$\FamPrdct{j\in J}{F(j)}$\IndSep product over discrete category $J$}
  \item as its colimit: $\CoLimOfOver{F}{J}$ the coproduct $\FamCoPrdct{j\in J}{F(j)}$. %
  \index{coproduct!as colimit of discrete diagram}\index[not]{c!$\FamCoPrdct{j\in J}{F(j)}$\IndSep coproduct over discrete category $J$}
\end{ulist}
If $J=\Set{0,1}$ contains two objects only, it is standard to write $F(0)\prdct F(1)$ for the limit of $F$; i.e.\ the product of $F(0)$ by $F(1)$. A morphism from an object $A$ to this product is then uniquely determined by two maps $f\from A\to F(0)$ and $g\from A\to F(1)$. The associated map into the product is denoted
\begin{equation*}
  \PrdctMapInto{f,g}\from A\to F(0)\prdct F(1)
\end{equation*}
\index[not]{p!$\PrdctMapInto{f,g}\from A\to \Prdct{X}{Y}$\IndSep universal map into a product}%

Dually, for $F\from \Set{0,1}\to \Ctgry{X}$, we write $F(0)+F(1)$ for the colimit of $F$; i.e.\ the sum of $F(0)$ and $F(1)$. A morphism from $F(0)+F(1)$ to an object $Z$ is uniquely determined by two maps $u\from F(0)\to Z$ and $v\from F(1)\to Z$. The associated map out of the sum is denoted
\begin{equation*}
  \SumMapOutOf{u,v}\from {F(0)+F(1)\to Z}
\end{equation*}
\index[not]{s!$\SumMapOutOf{u,v}\from X+Y\to Z$\IndSep universal map out of a coproduct / sum}

\subsection{Pullbacks$^{\ast}$}
\label{subsec:Pullbacks}

Choosing $J$ to be the category $a\to t\leftarrow b$, yields the category of diagrams in $\Ctgry{X} $ of the shape $A\XRA{f} T \XLA{g} B$. `Pullback' is the name for the limit of such a diagram. We state its universal property in terms of this diagram:
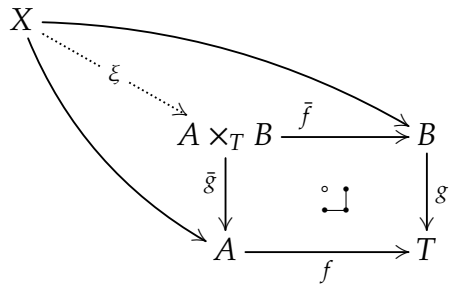
\begin{equation*}
  \xymatrix@R=5ex@C=4em{
  X \ar@/^2ex/[rrd] \ar@/_2ex/[rdd] \ar@{.>}[rd]|-{\ \xi\ } && \\
  & A\times_TB \ar[r]^(.4){\bar{f}} \ar[d]_(.4){\bar{g}} \PullLU{rd} &
  B \ar[d]^{g} \\
  & A \ar[r]_-{f} &
  T
  }
\end{equation*}
The universal property the pullback square asserts that any commutative diagram of solid arrows, as shown above, has a unique filler $\xi$ which renders the entire diagram commutative. %
\index{pullback}\index{universal property!pullback}

The morphism $\bar{g}$ depends functorially on the diagram $A\XRA{f} T \XLA{g} B$ and is said be the \Defn{base change} of $g$ along $f$. Similarly, $\bar{f}$ is said to be the base change of $f$ along $g$. %
\index{base change}%
Base change preserves monomorphisms and isomorphisms, (\ref{exe:BaseChange-Mono-Epi-Iso}), but fails to preserve epimorphisms in general. On the other hand, in many categories, such as $\Sets$, $\Grps$, module categories etc., base change does preserve epimorphisms.

\begin{proposition}[Pullbacks and concatenated squares]
  \label{thm:Pullbacks,ConcatenatedSquares}%
  In an arbitrary category consider the commutative diagram below. %
  \index{pullback!concatenated squares}%
  \begin{equation*}
    \xymatrix@R=5ex@C=3em{
    \DiagObj \ar[r] \ar[d] &
    \DiagObj \ar[r] \ar[d] &
    \DiagObj \ar[d] \\
    \DiagObj \ar[r] &
    \DiagObj \ar[r] &
    \DiagObj
    }
  \end{equation*}
  \begin{enumerate}[(i)]
    \item \emph{Pullback composition}\quad If both squares are pullbacks, so is the outer rectangle. %
          \index{pullback!composition}
    \item \emph{Pullback cancellation}\quad If the right hand square and the outer rectangle are pullbacks, then the left hand square is a pullback as well. \NoProof %
          \index{pullback!cancellation}%
  \end{enumerate}
\end{proposition}

A frequently used pullback construction is that of a kernel pair of a morphism:

\begin{definition}[Kernel pair]
  \label{def:KernelPair}%
  In a category with pullbacks, the \Defn{kernel pair} $\KrnlPr{f}=(X\times_YX,\PrjctnOnto{1},\PrjctnOnto{2})$ of a morphism $ f\from X\to Y$ is given by the pullback diagram below. %
  \index{kernel!pair of a morphism}\index[not]{k!$\KrnlPr{f}$\IndSep kernel pair of $f$}
  \begin{equation*}
    \xymatrix@R=5ex@C=4em{
    X\times_YX \ar[r]^-{\PrjctnOnto{1}} \ar[d]_-{\PrjctnOnto{2}} \PullLU{rd} &
    X \ar[d]^-{f} \\
    X \ar[r]_-{f} &
    Y
    }
  \end{equation*}
  We will sometimes abuse notation by writing $\KrnlPr{f}$ for the object $X\times_YX$.
\end{definition}

The maps $\PrjctnOnto{1}$ and $\PrjctnOnto{2}$ are jointly monomorphic. Moreover these maps have a common left inverse, namely the map $\DgnlOn{X}$ in the diagram below.
$$
  \xymatrix@R=5ex@C=4em{
  X \ar@{=}@/^2ex/[rrd] \ar@{=}@/_2ex/[rdd] \ar@{.>}[rd]|-{\DgnlOn{X} } \\
  & X\times_YX \ar[r]_-{\PrjctnOnto{1}} \ar[d]^-{\PrjctnOnto{2}} \PullLU{rd} &
  X \ar[d]^-{f} \\
  & X \ar[r]_-{f} &
  Y
  }
$$
Via the kernel pair construction we characterize monomorphisms by a limit property.

\begin{proposition}[Kernel pair and monomorphism]
  \label{thm:KernelPair-Monos}%
  In the kernel pair of a map $ f\from X\to Y$\*, the structure maps $\PrjctnOnto{1}$ and $\PrjctnOnto{2}$ have a unique common section $\DgnlOn{X}\from X\to X\times _YX$\*. Moreover, the following conditions are equivalent.
  \begin{enumerate}[(i)]
    \item $f$ is a monomorphism.
    \item $\PrjctnOnto{1}$ is a monomorphism.
    \item $\PrjctnOnto{1}$ is an isomorphism.
    \item $\DgnlOn{X}$ is an isomorphism.
    \item $\DgnlOn{X}$ is an epimorphism.
    \item $\PrjctnOnto{1}=\PrjctnOnto{2}$ \NoProof
  \end{enumerate}
\end{proposition}

\subsection{Pushouts}%
\label{subsec:Pushouts}%

The notion of pushout is dual to that of a pullback, and the notion of cobase change is dual to that of base change. %
\index{cobase change}%
Thus a pushout diagram in a category $\Ctgry{X}$ is a commutative square which has this universal property:
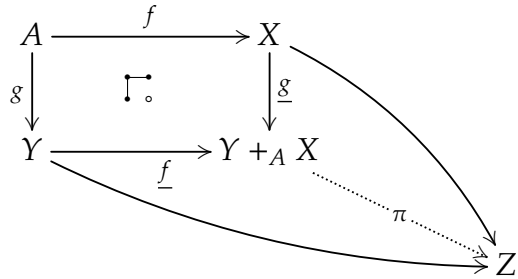
\begin{equation*}
  \xymatrix@R=5ex@C=5em{
  A \ar[r]^-{f} \ar[d]_{g} \PushRD{rd} &
  X \ar[d]^-{\underline{g}}\ar@/^2ex/[rdd] \\
  Y \ar[r]_(0.56){\underline{f}} \ar@/_2ex/[rrd] &
  Y+_AX \ar@{.>}[rd]|-{\ \ \pi\ \ } \\
  && Z
  }
\end{equation*}
For every commutative diagram of solid arrows, there exists a unique filler $\pi$ which renders the entire diagram commutative. - We refer to the morphism $\underline{g}$ as the \Defn{cobase change} of $g$ along $f$. %
\index{cobase change}

\begin{proposition}[Pushouts and concatenated squares]
  \label{thm:PushOuts,ConcatenatedSquares}%
  In an arbitrary category consider the commutative diagram below. %
  \index{pushout!concatenated squares}%
  \begin{equation*}
    \xymatrix@R=5ex@C=3em{
    \DiagObj \ar[r] \ar[d] &
    \DiagObj \ar[r] \ar[d] &
    \DiagObj \ar[d] \\
    \DiagObj \ar[r] &
    \DiagObj \ar[r] &
    \DiagObj
    }
  \end{equation*}
  \begin{enumerate}[(i)]
    \item \emph{Pushout composition}\quad If both squares are pushouts, so is the outer rectangle.
    \item \emph{Pushout cancellation}\quad If the left hand square and the outer rectangle are pushouts, then the right hand square is a pushout as well. %
          \index{pushout!cancellation}%
          \NoProof
  \end{enumerate}
\end{proposition}

\begin{subordinate}{}

  \begin{subsubordinate}{Existence, preservation, reflection of finite (co)limits}
    Sufficient conditions for the existence of finite (co-)limits are given in~\cite[Section~V.2]{SMacLane1998}. Further, any functor $F\from {\Ctgry{X}\to \Sets}$ with a left adjoint preserves all limits in $\Ctgry{X}$. If the adjunction is monadic, then $F$ reflects limits. In fact, it actually \emph{creates} them in the sense of~\cite[Theorem~23.11]{AHS:Cats}.
  \end{subsubordinate}

  \begin{subsubordinate}{Functorial (co)limits}
    It may be advantageous to work with (co)limits which respond functorially to a change of the template category $J$. An appropriate context for this discussion is: The operation which sends a small category $J$ to the functor category $\Ctgry{X}^J$ is a functor $\Phi\from \SmallCats^{\op}\to \LocSmallCats$, with $\SmallCats$ the category of small categories, and $\LocSmallCats$ the huge category of locally small categories. The associated Grothendieck category $\GrothendieckCContravOf{\Phi}$ contains $\Ctgry{X} \cong \Ctgry{X}^{\Ord{0}}$ as a full subcategory. Functorial colimits in $\Ctgry{X}$ are given by a left adjoint to the inclusion $\Ctgry{X} \to \GrothendieckCContravOf{\Phi}$.
  \end{subsubordinate}
\end{subordinate}

\begin{exercises}

\begin{exercise}[Jointly monomorphic maps]
  \label{exe:JointlyMonomorphicMaps}
  In a complete category $\SACtgry{X}$ show the following:
  \begin{thmlist}
    \item A family of morphisms $\{f_i\from A\to X_i\}_{i\in I}$ is jointly monomorphic if and only if the universal map $(f_i)\from A\to \FamPrdct{i\in I}{X_i}$ is a monomorphism. %
    \index{jointly!monomorphic family - product description}%
    \item If $F\from J\to \SACtgry{X}$ is a functor from a small category $J$ into $\SACtgry{X}$, then the structure maps $p_i\from \LimOf{F}\to F_i$ of a limit cone are jointly monomorphic. %
    \label{exe:LimitCone-StructureMapsJointlyMonomorphic}%
    \index{jointly!monomorphic family - from limit cone}%
  \end{thmlist}
\end{exercise}

\begin{exercise}[Base change and monomorphisms / epimorphisms / isomorphisms]
  \label{exe:BaseChange-Mono-Epi-Iso}%
  In any category, consider a pullback square and show the following:
  $$
    \xymatrix@R=6ex@C=4em{
    P \ar[r]^{f} \ar[d]_{m} \PullLU{rd}&
    M \ar[d]^{\mu} \\
    X \ar[r]_-{\varphi} &
    Y
    }
  $$
  \begin{enumerate}[(i)]
    \item If $\mu$ is a monomorphism, respectively an isomorphism, show that $m$ is a monomorphism, respectively an isomorphism. %
          \index{pullback!preserves monomorphism}%
          \index{monomorphism!is preserved by pullback}%
          \index{base change!preserves monomorphisms}%
          \index{base change!preserves isomorphism}%
    \item In the category $\Sets$ of sets, show that if $\mu$ is a surjective function, then so is $m$. %
          \index{base change!preserves surjective function}
  \end{enumerate}
\end{exercise}

\begin{exercise}[Pullbacks in commutative cube]
  \label{exe:PullbacksCube}
  In an arbitrary category, suppose the solid arrowed portion faces of a cube the cube below are pullback squares.
  \begin{equation*}
    \xymatrix@R=3ex@C=3em{
    & \DiagObj \ar@{.>}[rr] \ar@{.>}[dd]|\hole \ar@{.>}[ld] &&
    \DiagObj \ar[ld] \ar[dd] \\
    \DiagObj \ar[rr] \ar[dd] &&
    \DiagObj \ar[dd] \\
    & \DiagObj \ar[rr]|\hole \ar[ld] &&
    \DiagObj \ar[ld] \\
    \DiagObj \ar[rr] &&
    \DiagObj
    }
  \end{equation*}
  Then construct the top face as a pullback, and show the following:
  \begin{thmlist}
    \item There is a unique arrow rendering the entire cube commutative.
    \item All faces of the so obtained cube are pullback squares.
    \item Show that the so obtained dotted arrows form a limit cone for the solid arrowed part of the diagram.
  \end{thmlist}
\end{exercise}
\begin{exercise}[Cobase change and monomorphisms / epimorphisms / isomorphisms]
  \label{exe:CoBaseChange-Mono-Epi-Iso}%
  In any category, consider a pushout square and show the following:
  $$
    \xymatrix@R=6ex@C=4em{
    P \ar[r]^{f} \ar[d]_{e} \PushRD{rd}&
    M \ar[d]^{\varepsilon} \\
    X \ar[r]_-{\varphi} &
    Y
    }
  $$
  \begin{enumerate}[(i)]
    \item If $e$ is an epimorphism, respectively an isomorphism, show that $\varepsilon$ is an epimorphism, respectively an isomorphism. %
          \index{pushout!preserves epimorphism}%
          \index{epimorphism!is preserved by pushout}%
          \index{cobase change!preserves epimorphisms}%
          \index{cobase change!preserves isomorphism}%
    \item In the category $\Sets$ of sets, show that if $e$ is injective, then so is $\varepsilon$.
  \end{enumerate}
\end{exercise}

\begin{exercise}[Cobase change resulting in isomorphism]
  \label{exe:CoBaseChange-IsoResult}
  In any category suppose $a\from A\to B$ is the composite $A\XRA{f} X \XRA{g} B$. If $f$ is an epimorphism, and the square below is a pushout, show that $\hat{f}$ is an isomorphism.
  \begin{equation*}
    \xymatrix@R=5ex@C=4em{
    A \ar[r]^-{f} \ar[d]_{a} \PushRD{rd} &
    X \ar[d]^{\hat{a}} \\
    B \ar[r]_-{\hat{f}} &
    P
    }
  \end{equation*}
\end{exercise}

\begin{exercise}[Pushouts and concatenated squares]
  Prove (\ref{thm:PushOuts,ConcatenatedSquares}): A composite of pushout squares is a pushout square, and pushout cancellation.
\end{exercise}

\begin{exercise}[Constructing the intersection of subobjects using a pullback]
  \label{exe:IntersectionConstruct-Via-Pullback}
  In a category which admits pullbacks, let $m\from M\to X$ and $n\from N\to X$ represent subobjects of $X$. In the pullback diagram below, show that $m\bar{n}=n\bar{m}\from P\to X$ represents $M\meet N$ as a subobject of $X$. %
  \index{intersection!of subobjects: existence}\index{subobjects!intersection: existence}%
  \index{meet!of subobjects: existence}\index{subobjects!meet: existence}%
  \begin{equation*}
    \xymatrix@R=5ex@C=4em{
    P \PullLU{rd} \ar[r]^-{\bar{n}} \ar[d]_-{\overline{m}} &
    M \ar[d]^{m} \\
    N \ar[r]_-{n} &
    X
    }
  \end{equation*}
\end{exercise}

\begin{exercise}[Isomorphic kernel pairs]
  \label{exe:KernelPairs-Isomorphism}
  \cite[2.16]{Bourn-Gran-CategoricalFoundations}\quad In a category with pullbacks, suppose a morphism $g$ fits into the diagram below:
  \begin{equation*}
    \xymatrix@R=6ex@C=3em{
    X \ar@<0.5ex>[rr]^-{f} \ar[rd]_{g} &&
    Y \ar@<0.5ex>[ll]^-{s} \ar[ld]^{h} \\
    & Z
    }
  \end{equation*}
  If $fs=\IdMapOn{Y}$, show that the kernel pair of $g$ is isomorphic to the kernel pair of $f$ if and only if $h$ is monic.
\end{exercise}

\begin{exercise}[Finite (co-)products from binary (co-)products]
  \label{exe:Binary(Co)Products->Finite(Co)Products}%
  In an arbitrary category $\Ctgry{X}$ show the following:
  \begin{thmlist}
    \item If $\Ctgry{X}$ admits binary products, show that objects $X_1,X_2,\dots X_n$ admit a product.
    \item If $\Ctgry{X}$ admits binary coproducts, show that objects $X_1,X_2,\dots X_n$ admit a coproduct.
    \item If $\Ctgry{X}$ admits binary biproducts, show that objects $X_1,X_2,\dots X_n$ admit a biproduct.
    \item If the comparison map $X_1+\cdots +X_n\to X_1\prdct \cdots \prdct X_n$ is an isomorphism, show that the $n$-ary biproduct of $X_1,\dots, X_n$ exists.
  \end{thmlist}
\end{exercise}

\begin{exercise}[Pullback of maps with common codomain]
  \label{exe:Pullback-MapsCommonCodomain}
  In any category with binary products and pullbacks show the following. If $a\from A\to X$ and $b\from B\to X$ are arbitrary morphisms, then the universal objects in the pullback diagrams below are naturally isomorphic.
  \begin{equation*}
    \xymatrix@R=5ex@C=4em{
    P \ar[r]^-{\PrdctMapInto{\alpha,\beta}} \ar[d]_{p} \PullLU{rd} &
    \Prdct{A}{B} \ar[d]^{\Prdct{a}{b}} &
    Q \ar[r]^-{\varphi} \ar[d]_{\psi} \PullLU{rd} &
    B \ar[d]^{b} \\
    X \ar[r]_-{\DgnlOn{X}} &
    \Prdct{X}{X} &
    A \ar[r]_-{a} &
    X
    }
  \end{equation*}
\end{exercise}
\end{exercises}
\section[Absolute Diagrams]{Diagrams with Absolute Properties}
\label{sec:AbsoluteDiagrams}

In \cite[I.1]{SMacLane1998}, Mac Lane defines a category to be a graph with certain properties. It then makes sense to say a \Defn{diagram} $D$ in a category $\Ctgry{X}$ is a subgraph of $\Ctgry{X}$. We say a property $P$ of a diagram $D$ in $\Ctgry{X}$ is \Defn{absolute} if, for every functor $F\from \Ctgry{X}\to \Ctgry{Y}$, the diagram $F(D)$ has property $P$ as well. Elementary examples of absolute diagram properties include: %
\index{diagram}\index{absolute!diagram property}\index{diagram!with absolute property}%
\begin{ulist}
  \item An equality between (composites of) morphisms is an absolute property of a diagram. In particular, being commutative is an absolute property of a diagram $D$.
  \item Being an isomorphism is an absolute property of a diagram which consists of a single morphism.
  \item Being an idempotent is an absolute property of an endomorphism $e$ of object $X$.
\end{ulist}

Below, we single out certain diagrams with absolute properties, often (co)limit (co)cones, which occur frequently.

\subsection{Absolute epimorphisms}
\label{subsec:AbsoluteEpis}%
An epimorphism (\ref{def:Epimorphism}) in a category $\Ctgry{X}$ is absolute if $F(f)$ is an epimorphism for every functor $F\from \Ctgry{X}\to \Ctgry{Y}$. We show that this happens if and only if $f$ admits a section. %
\index{absolute!epimorphism}\index{epimorphism!absolute}%

\begin{definition}[Sectionable / sectioned morphism]
  \label{def:SectionedEpimorphism}
  A morphism $q\from X\to Q$ is called \Defn{sectionable}, also \Defn{splittable}, if there exists $x\from Q\to X$ such that $q\Comp x=\IdMapOn{Q}$. A choice of such $x$ is called a \Defn{section} or a \Defn{splitting} of $q$, and we say that $q$ is sectioned by $x$. We write $\SctndEpi{q}{x}$ for a morphism $q$ sectioned by a morphism $x$.%
  \index{epimorphism!split}\index{sectionable!epimorphism}\index{morphism!split}%
  \index{sectionable!morphism}\index{section}\index{splitting}%
  \index[not]{s!$\SctndEpi{q}{x}$\IndSep map $q$ sectioned by $x$}
\end{definition}

If $x$ is a section of $q$, then $q$ is a left inverse of $x$, and $x$ is a right inverse of $q$. Being a sectionable epimorphism is an absolute property. Direct verification shows that every sectionable morphism is an absolute epimorphism. The converse also holds:

\begin{proposition}[Sectionable epi if and only if absolute epi]
  \label{thm:Sectionable<->AbsoluteEpi}%
  A morphism $q\from {X} \to {Q}$ is sectionable if and only if it is an absolute epimorphism. In this case, every section of $q$ is an absolute monomorphism.
\end{proposition}
\begin{proof}
  If $q$ admits a section $x$, then $F(q)$ also admits a section, namely $F(x)$. So, admitting a section is an absolute property. Conversely, suppose $q\from X\to Q$ is an absolute epimorphism. Then the functor $\Hom{Q}{-}\from \Ctgry{X}\to \Sets$ sends $q$ to the epimorphism, i.e., surjective function, $q_{\ast}\from \Hom{Q}{X}\to \Hom{Q}{Q}$ in $\Sets$. So, there is a morphism $x\from Q\to X$ in $\Ctgry{X}$ such that $q\Comp x=\IdMapOn{Q}$. This makes $x$ is a section of $p$.
\end{proof}

Dually, a morphism $m\from M\to X$ is \Defn{retractable} if there exists $r\from X\to M$ such that $rm=\IdMapOn{M}$. In this situation we refer to $r$ as the retraction of $X$ to $M$ and write $\RtrctdMono{m}{r}$. Being a retractable morphism is an absolute property. Such a map is always monic. Dualizing (\ref{thm:Sectionable<->AbsoluteEpi}), we find: %
\index{retractable monomorphism}\index{morphism!retractable}%
\index[not]{r!$\RtrctdMono{m}{r}$\IndSep mono $m$ retracted by $r$}%

\begin{proposition}[Retractable if and only if absolute mono]
  \label{thm:Retractable<->AbsoluteMono}%
  A map $m\from M\to X$ in retractable if and only if $m$ is an absolute monomorphism. \NoProof%
  \index{retractable!monomorphism}\index{monomorphism!retractable}
\end{proposition}

To clarify, a sectioned epimorphism $\SctndEpi{q}{x}$ has an associated retracted monomorphism $\RtrctdMono{x}{q}$ and vice versa. We must distinguish between the two structures because, in general, neither does $q$ determine $x$ uniquely, nor does $x$ determine $q$ uniquely. So, it matters which of the two maps is given first and then coupled with a \emph{choice} of the other.

The collection of sectioned epimorphisms in $\Ctgry{X}$ forms a category $\SEpisIn{X}$. A morphism $(f,g)\from \SctndEpi{q}{x}\to \SctndEpi{r}{y}$ in this category%
\footnote{In the literature, a split epimorphism is sometimes referred to as a \emph{point}. Accordingly, the category $\SEpisIn{X}$ is known as \emph{the category of points in $\Ctgry{X}$}. }%
is given by a diagram with commuting properties stated below.
\begin{equation}
  \label{fig:SectionedEpi-Morphism}%
  \vcenter{
  \xymatrix@R=5ex@C=5em{
  X \ar[r]^-{f} \ar@<-0.5ex>[d]_{q} &
  Y \ar@<-0.5ex>[d]_{r} \\
  Q \ar[r]_-{g} \ar@<-0.5ex>[u]_{x} &
  R \ar@<-0.5ex>[u]_{y}
  }
  }
\end{equation}
With $q\Comp x=\IdMapOn{Q}$ and $r\Comp y=\IdMapOn{R}$, the vertical structures are objects in $\SEpisIn{X}$. For $(f,g)$ to form a morphism in $\SEpisIn{X}$, we require $g\Comp q = r\Comp f$, and $f\Comp x=y\Comp g$. %
\index[not]{s!$\SEpisIn{X}$\IndSep category of sectioned morphisms in $\Ctgry{X}$}

\begin{proposition}[Pullback of a sectioned epimorphism]
  \label{thm:PullbackOfSplitEpi}%
  In any category consider an epimorphism $r\from Y\to R$, sectioned by $y\from R\to Y$ and pulled back along $f\from Q\to R$: %
  \index{pullback!sectioned epimorphism}%
  \begin{equation*}
    \xymatrix@R=5ex@C=5em{
    P \ar[r]^-{\bar{f}} \ar@<-.5ex>[d]_{\bar{r}} \PullLU{rd} &
    Y  \ar@<-.5ex>[d]_{r} \\
    Q \ar[r]_-{f} \ar@<-.5ex>[u]_{\bar{y}}&
    R \ar@<-.5ex>[u]_{y}
    }
  \end{equation*}
  Then $\bar{r}$ is sectionable, and there is exactly one section $\bar{y}$ with the property $\bar{f}\bar{y} = yf$. Thus we obtain the morphism $(\bar{f},f)\from\SctndEpi{\bar{r}}{\bar{y}}\to \SctndEpi{r}{y}$ in $\SEpisIn{X}$. Moreover, $Q$ is a pullback of $y$ along $\bar{f}$.%
  \label{exe:SplitEpiStable}
  \label{exe:SplitEpi-PreservedByBaseChange}
  \NoProof
\end{proposition}

\begin{proposition}[Pushout of a sectioned monomorphism]
  \label{thm:PushOutOfSplitMono}
  In any category consider a monomorphism $x\from Q\to X$, with retraction $q\from X\to Q$ and pushed out along $f\from Q\to R$: %
  \index{pushout!of a sectioned monomorphism}%
  \begin{equation*}
    \xymatrix@R=5ex@C=5em{
    X \ar[r]^-{\underline{f}} \ar@<-.5ex>[d]_{q} \PushRU{rd} &
    P  \ar@<-.5ex>[d]_{\underline{q}} \\
    Q \ar[r]_-{f} \ar@<-.5ex>[u]_{x}&
    R \ar@<-.5ex>[u]_{\underline{x}}
    }
  \end{equation*}
  Then $\underline{x}$ is retractable, and there is exactly one retraction with the property $\underline{q}\underline{f} = fq$. Thus we obtain the morphism $(\underline{f},f)\from\SctndEpi{q}{x}\to \SctndEpi{\underline{q}}{\underline{x}}$ in $\SEpisIn{X}$. Moreover, $\underline{q}$ is a pushout of $q$ along $\underline{f}$. %
  \NoProof
\end{proposition}

Certain morphisms in $\SEpisIn{X}$ yield absolute diagrams:
\begin{proposition}[Pushout / pullback from morphism in $\SEpisIn{X}$]
  \label{thm:AbsolutePush/Pull-MorInSEpi(X)}
  For a morphism $(f,g)\from \SctndEpi{q}{x}\to \SctndEpi{r}{y}$, as in \eqref{fig:SectionedEpi-Morphism}, the following hold:
  \begin{thmlist}
    \item If $f$ is an epimorphism, then the square $X\rightrightarrows R$ is a pushout.
    \item If $f$ is a monomorphism, then the square $Q\rightrightarrows Y$ is a pullback.\NoProof
  \end{thmlist}
\end{proposition}

\begin{corollary}[Absolute pushout / pullback from sectioned morphism in $\SEpisIn{X}$]
  \label{thm:AbsolutePush/Pull-SectionedMorInSEpi(X)}%
  Consider a sectioned morphism in $\SEpisIn{X}$; i.e.,
  \[
    (f,g)\from \SctndEpi{q}{x}\to \SctndEpi{r}{y} \qquad\text{and}\qquad (\varphi,\gamma)\from \SctndEpi{r}{y}\to \SctndEpi{q}{x}
  \]
  are morphisms of sectioned epimorphisms, and $(\varphi,\gamma)\Comp (f,g)= \IdMapOn{\SctndEpi{q}{x}}$.
  \begin{equation}
    \label{fig:SectionedMapInSEpi(X)}%
    \vcenter{
    \xymatrix@R=5ex@C=5em{
    X \ar@<0.5ex>[r]^-{f} \ar@<-0.5ex>[d]_{q} &
    Y \ar@<-0.5ex>[d]_{r} \ar@<0.5ex>[l]^-{\varphi} \\
    Q \ar@<0.5ex>[r]^-{g} \ar@<-0.5ex>[u]_{x} &
    R \ar@<-0.5ex>[u]_{y}\ar@<0.5ex>[l]^-{\gamma}
    }
    }
  \end{equation}
  Then the square $X\rightrightarrows R$ is an absolute pushout, and the square $R\rightrightarrows X$ is an absolute pullback. \NoProof
\end{corollary}

\begin{proposition}[Absolute coequalizer from kernel pair of absolute epimorphism]
  \label{thm:Coeq(KP(AbsoluteEpi))IsAbsolute}%
  In any category $\Ctgry{X}$, if the kernel pair of an absolute epimorphism $f\from X\to Y$ exists, then the diagram
  \begin{equation*}
    \xymatrix@R=5ex@C=4em{
    \KrnlPr{f} \ar@<+0.5ex>[r]^-{\PrjctnOnto{1}} \ar@<-0.5ex>[r]_-{\PrjctnOnto{2}} &
    X \ar@<0.1ex>[r]^-{\ f\ } &
    Y
    }
  \end{equation*}
  is an absolute coequalizer.
\end{proposition}

\begin{exercises}
\begin{exercise}[On sectionable epimorphisms]
  \label{exe:SplitEpis}
  \label{exe:SectionableEpis}%
  In an arbitrary category $\Ctgry{X}$ show the following:
  \begin{thmlist}
    \item Every sectionable epimorphism (\ref{def:SectionedEpimorphism}) is an epimorphism.%
    \index{sectionable!epimorphism is epimorphism}%
    \item If a morphism $f$ is both an absolute epimorphism and a monomorphism, then it is an isomorphism.%
    \index{sectioned epimorphism!$+$ monomorphism is isomorphism}%
    \index{absolute epimorphism!$+$ monomorphism is isomorphism}
  \end{thmlist}
\end{exercise}

\begin{exercise}[Absolute epimorphism need not be an isomorphism]
  \label{exe:AbsoluteEpi-NeedNotBe-Iso}
  Give examples of categories in which an absolute epimorphism need not be an isomorphism.
\end{exercise}

\begin{exercise}[Morphism of sectioned epimorphisms: Properties]
  \label{exe:MorphismOfSplitEpis-Properties}
  Given a morphism $(f,g)\from \SctndEpi{q}{x} \to \SctndEpi{r}{y}$ of sectioned epimorphisms, as in \eqref{fig:SectionedEpi-Morphism}, do the following:
  \begin{enumerate}[(i)]
    \item Prove (\ref{thm:AbsolutePush/Pull-MorInSEpi(X)}.i) and (\ref{thm:AbsolutePush/Pull-MorInSEpi(X)}.ii).
    \item Prove (\ref{thm:AbsolutePush/Pull-SectionedMorInSEpi(X)}).
    \item Give an example which shows that, for (i), the assumption that $f$ is an epimorphism is essential.
    \item Give an example which shows that, for (ii), the assumption that $f$ is a monomorphism is essential.
  \end{enumerate}
\end{exercise}

\begin{exercise}[Double absolute epimorphisms]\label{exe:DoubleAbsoluteEpi}
  Is a commutative square consisting of absolute epimorphisms always an absolute pushout?
\end{exercise}

\begin{exercise}[Coequalizer by pushout]
  \label{exe:CoEqualizer-Pushout}
  Suppose two epimorphisms $f$, $g\from X\to Y$ admit a common section $e\from Y\to X$. Show that the coequalizer of $f$ and $g$ is given by the pushout
  \begin{equation*}
    \xymatrix@R=5ex@C=4em{
    X \ar@{-{ >>}}[r]^-{f} \ar@{-{ >>}}[d]_{g} \PushRD{rd} &
    Y \ar[d]^{u} \\
    Y \ar[r]_-{v} &
    Z
    }
  \end{equation*}
  If at least one of $f$ or $g$ is a normal epimorphism, conclude that $u=v$ is a normal epimorphism.
\end{exercise}
\end{exercises}
\section{Around Surjectivity}
\label{sec:Surjectivity}

A function of sets $f\from X\to Y$ is surjective if, for every $y\in Y$, its fiber $f^{-1}(y)$ under $f$ is not empty. If the objects in a category $\Ctgry{X}$ do not have underlying sets, an approach to defining surjectivity of a morphism in $\Ctgry{X}$ via the fiber of its morphisms is not available. Instead, we look for properties of functions which are independent of underlying sets and give an alternate characterization of surjectivity.

We already encountered one such property: a function $f\from X\to Y$ is surjective if and only if it has the right cancellation property in composites: $uf=vf$ implies $u=v$. This condition can be formulated in any category, and gives the definition of an epimorphism (\ref{def:Epimorphism}); see Section \ref{sec:SubObjects-QuotientObjects}. We will work with several more alternate characterizations of surjective functions. When interpreted in a category with pullbacks they all imply the epimorphism property, yet they differ subtly in other respects. Thus they are known as types of epimorphisms.

The graphic below presents the types of epimorphisms we encounter. Solid arrows indicate implication of property in any category with the appropriate properties, as indicated by their labels.

\begin{figure}[ht]
  \begin{equation*}
    \resizebox{.98\textwidth}{!}{$
      \xymatrix@R=8ex@C=3.8em{
      *+[F-,]{\txt{absolute\\epi}}\ar@{=>}[rrd]|-*+[o][F-]{\ \txt{L}\ } \ar@{=>}[rr]|-*+[o][F-]{\ \txt{L}\ } \ar@/_2ex/@{.>}[rd]|-*+[o][F-]{\ \txt{N}\ } &&
      *+[blue][F-,]{\txt{effective\\epi}} \ar@<-2ex>@{=>}[d] \\
      &*+[blue][F-,]{\txt{normal\\epi}} \ar@<-1.2ex>@{=>}[r] &
      *+[blue][F-,]{\txt{regular\\epi}} \ar@<-1.2ex>@{=>}[l]|-*+[o][F-]{\ \txt{N}\ } \ar@<-1.7ex>@{=>}[r]|-*+[o][F-]{\ \txt{L}\ } \ar@{=>}@<-2ex>[u]|-*+[o][F-]{\ \txt{L}\ }
      & *+[blue][F-,]{\txt{strong\\epi}} \ar@<-1.2ex>@{=>}[r] \ar@<-1.4ex>@{=>}[l]|-*+[o][F-]{\ \txt{N}\ }
      & *+[blue][F-,]{\txt{extremal\\epi}} \ar@<-1.2ex>@{=>}[l]|-*+[o][F-]{\ \txt{L}\ } \ar@{=>}[r]|-*+[o][F-]{\ \txt{L}\ } \ar@/^6ex/@{.>}[lll]|-*+[o][F-]{\ \txt{N}\ }
      & *+[F-,]{\txt{epi}} \ar@/^9ex/@{.>}[llll]|-*+[o][F-]{\ \txt{A}\ }
      }$}
  \end{equation*}
  \caption{Several types of epimorphisms}\label{fig:Epis}
\end{figure}
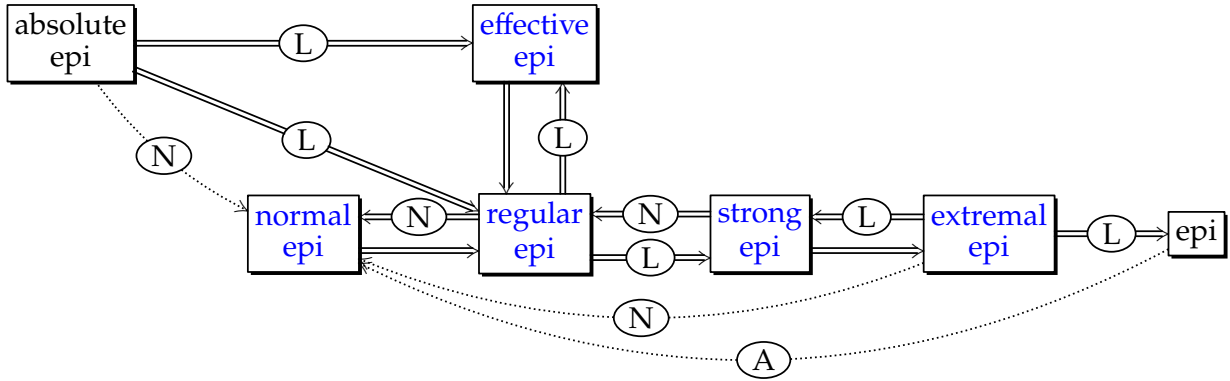
\begin{enumerate}
  \item An undecorated arrow indicates an implication which holds in every category.
  \item In a category with finite limits, the implications labeled `L' hold.
        %
        %
  \item In a normal category (\ref{def:NormalCat}) and, hence, in homological category, the implications labeled `N', or `L' hold. ; see Section~\ref{sec:ImageFactorizations}.
  \item In an abelian category all implications hold, also the one labeled `A'.
\end{enumerate}
In a normal category then, the properties of normal, regular, effective, strong, and extremal epimorphisms all coincide; see (\ref{thm:CoKer=NormalEpi=RegEpi=EffectiveEpi}). Therefore, we mostly work with three types of epimorphisms: split epimorphisms, normal epimorphisms (equivalently cokernels), and epimorphisms. Using the defining property of an equivalent type is only of occasional advantage. In the more specialized context of abelian categories, we only need to distinguish between epimorphisms and split epimorphisms.

In a variety of algebras $\Ctgry{V}$, surjective morphisms are the same as effective/regular/strong epimorphisms. On other hand, in general, a regular epimorphisms need not be a surjective algebra morphism. - Let us now turn to details.

\begin{definition}[Extremal epimorphism]%
  \label{def:ExtremalEpimorphism}%
  A morphism $f\from X\to Y$ is called an \Defn{extremal epimorphism} if it \emph{only} factors through the subobject of $Y$ represented by $\IdMapOn{Y}$. %
  \index{extremal epimorphism}\index{epimorphism!extremal}
\end{definition}

Thus an extremal epimorphism does not factor through any proper subobject of its codomain. If $f$ admits an image factorization, then $f$ is extremal if and only if its image is the codomain of $f$. It is useful to express definition (\ref{def:ExtremalEpimorphism}) in diagram language: $f$~is a extremal epimorphism if and only if, for every factorization of $f$ through a monomorphism $m$, the map $m$ is an isomorphism.
\begin{equation*}
  \xymatrix@R=5ex@C=3em{
  & M \ar@{ >->}[d]^-{m} \\
  X \ar[r]_-{f} \ar[ru]^-{x} &
  Y
  }
\end{equation*}

More generally, a pair of maps $ (r\from A\to X, s\from B\to X)$ in $\Ctgry{C}$ is \Defn{jointly extremal-epimorphic} provided in any commutative diagram
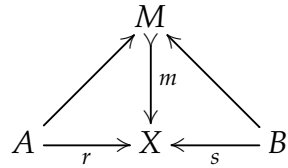
\begin{equation*}
  \xymatrix@!0@=4em{
  & M \ar@{ >->}[d]^- m \\
  A \ar[r]_-{r} \ar[ur] &
  X & B \ar[l]^-{s} \ar[ul]
  }
\end{equation*}
in which $m$ is a monomorphism, it is necessarily an isomorphism. %
\index{jointly!extremal-epimorphic family of maps}%
The intuition is that `$r$ and $s$ together generate all of $X$'. In a category with binary sums, this is equivalent to saying that the arrow $\SumMapOutOf{r,s}\from A+B\to X$ is an extremal epimorphism.

From its definition, it is not at all apparent that an extremal epimorphism is actually an epimorphism. That this is, indeed, the case is the content of the following lemma.

\begin{proposition}[Extremal epimorphism is epimorphism]
  \label{thm:ExtremalEpi->Epi}%
  In a category with equalizers, an extremal epimorphism is an epimorphism.
\end{proposition}
\begin{proof}
  Suppose $ f\from X\to Y$ is an extremal epimorphism, and consider maps $ u$, $v\from Y\to Z$ with $uf=vf$. We then find ourselves in the situation depicted below.
  
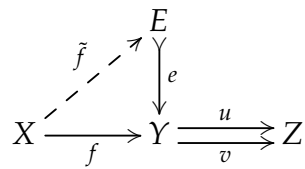
\[
    \xymatrix@R=5ex@C=3em{
    & E \ar@{ >->}[d]^-{e} \\
    X \ar[r]_-{f} \ar@{-->}[ru]^-{\tilde{f}} &
    Y \ar@<+.5ex>[r]^-{u} \ar@<-.5ex>[r]_-{v} &
    Z
    }
  \]
  Here $e$ is an equalizer of $u$ and $v$, hence is monic. Via its universal property, $f$ admits the indicated factorization $\tilde{f}$. As $f$ is an extremal epimorphism, $e$ is an isomorphism, implying that $u=v$. So $f$ is an epimorphism.
\end{proof}

The following alternate characterization of `extremal epimorphism' achieves compatibility with factorization systems (Section~\ref{sec:FactorizationSystems}) and strong epimorphisms (\ref{def:StrongEpimorphism}).

\begin{proposition}[Extremal epimorphism via lifting]
  \label{thm:ExtremalEpi-Via-Lifting}
  A map $f\from X\to Y$ in a category $\EuScript{X}$ is an extremal epimorphism if and only if every commutative square as shown below admits a filler which renders the entire diagram commutative. %
  $$
    
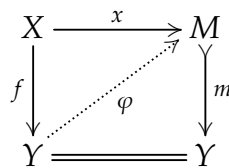
\xymatrix@R=6ex@C=4em{
    X \ar@{->}[d]_-{f} \ar[r]^-{x} & M \ar@{ >->}[d]^-{m} \\
    Y \ar@{.>}@<+.5ex>[ru]_-{\varphi} \ar@{=}[r] &
    Y
    }
  $$
  \index{epimorphism!extremal}\index{extremal epimorphism}
\end{proposition}
\begin{proof}
  The lifting condition means that $m$ is a split epimorphism and a monomorphism, hence an isomorphism; see (\ref{exe:SplitEpis}). Conversely, if $m$ is an isomorphism then $\varphi=m^{-1}$ is the one and only lifting $\varphi$.
\end{proof}

Once a category has pullbacks, the lifting property in (\ref{thm:ExtremalEpi-Via-Lifting}) turns out to be equivalent to an \emph{a priori} stronger condition, namely:

\begin{definition}[Strong epimorphism]
  \label{def:StrongEpimorphism}%
  In any category, a morphism $ f\from X\to Y$ is a \Defn{strong epimorphism} if every commutative square, with $m$ a monomorphism admits a unique filler $\varphi$ which renders the entire diagram commutative. %
  \index{epimorphism!strong}%
  \index{strong epimorphism}%
  \[
    \xymatrix@R=5ex@C=3em{
    X \ar[d]_-{f} \ar[r]^-{x} &
    M \ar@{ >->}[d]^-{m} \\
    Y \ar[r]_-{y} \ar@{.>}[ru]^-{\varphi} &
    Z}
  \]
\end{definition}

\begin{proposition}[Extremal vs. strong epimorphism]
  \label{thm:ExtremalEpi-StrongEpi}
  Every strong epimorphism is extremal.	In a category with pullbacks, every extremal epimorphism is strong. \NoProof
\end{proposition}

\begin{definition}[Regular/effective epimorphism]
  \label{def:Regular/EffectiveEpi}%
  A morphism $ f\from X\to Y$ is called a \Defn{regular epimorphism} it is the coequalizer of some parallel pair of arrows. It is called an \Defn{effective epimorphism} if it is the coequalizer of its kernel pair. %
  \index{regular!epimorphism}%
  \index{epimorphism!regular}%
  \index{effective!epimorphism}%
  \index{epimorphism!effective}%
\end{definition}

\begin{proposition}[Regular vs. effective epimorphism]
  \label{thm:EffectiveEpi=RegularEpi}
  In a category with pullbacks, a morphism is a regular epimorphism if and only if it is an effective epimorphism.
\end{proposition}
\begin{proof}
  Comparing definitions shows that an effective epimorphism is regular. To see that the converse holds, consider a regular epimorphism $ f\from{X\to Y}$\*, expressed as the coequalizer of maps $u$, $v\from{W\to X}$\*. Let $\KrnlPr{f}=(X\times_YX,\pi_{1},\pi_{2})$ be the kernel pair of $f$ and $ g\from{X\to Z}$ a morphism which satisfies $gu=gv$.

  The universal property of the kernel pair yields a map $\gamma\from W\to X\times_YX$\*, unique with the properties $ \pi_1\gamma=u$ and $ \pi_2\gamma=v$\*. We conclude
  $$
    gu = g\pi_1\gamma = g\pi_2\gamma = gv.
  $$
  So the universal property of the coequalizer $f$ yields a map $ t\from Y\to Z$\*, unique with the property $ tf =g$\*. This proves that $f$ is a coequalizer of its kernel pair.
\end{proof}

A normal epimorphism (\ref{def:NormalEpi}) is always a regular epimorphism. Whenever a category has pullbacks, the following hold.

\begin{lemma}[Absolute epi implies regular epi implies strong epi]
  \label{thm:CoKer=Strong=ExtremalEpi}%
  In an arbitrary category $\Ctgry{X}$, consider the following properties a morphism $f\from X\to Y$:
  \begin{enumerate}[(i)]
    \item $f$ is an absolute epimorphism;
    \item $f$ is an effective epimorphism;
    \item $f$ is a regular epimorphism;
    \item $f$ is an extremal epimorphism;
    \item $f$ is a strong epimorphism.
  \end{enumerate}
  If $\Ctgry{X}$ has pullbacks, then (i) $\Rightarrow$ (ii) $\Leftrightarrow$ (iii) $\Rightarrow$ (iv) $\Leftrightarrow$ (v).
\end{lemma}
\begin{proof}
  (i) $\implies$ (ii) Suppose $x\from {Y\to X}$ is a section of $f$. We show directly that $f$ satisfies the universal property required of a coequalizer of its kernel pair. If $(\KrnlPr{f},\pi_1,\pi_2)$ is the kernel pair of $f$, consider $g\from X\to Z$ with $g\pi_1=g\pi_2$. We claim that $gx\from Y\to Z$ is the unique factorization of $f$ through $f$. Uniqueness follows from the epimorphic property of $f$. To see that $g=(gx)f$, let $\eta\from X\to \KrnlPr{f}$ be the universal map induced by $\PrdctMapInto{\IdMapOn{X},sf}$.
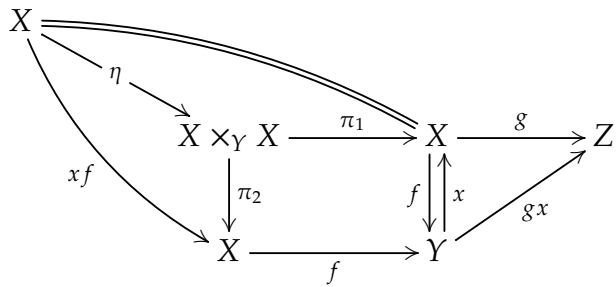
  \begin{equation*}
    \xymatrix@R=5ex@C=4em{
    X \ar@{=}@/^2ex/[rrd] \ar@/_2ex/[rdd]_{xf} \ar[rd]|-{\ \eta\ } &&\\
    & X\times_YX \ar[r]^-{\pi_1} \ar[d]^{\pi_2} &
    X \ar@<-.5ex>[d]_{f} \ar[r]^-{g} &
    Z \\
    & X \ar[r]_-{f} &
    Y \ar[ru]_{gx} \ar@<-.5ex>[u]_{x}
    }
  \end{equation*}
  Then we find, as required:
  \begin{equation*}
    g = g\pi_1\eta = g\pi_2\eta = gxf
  \end{equation*}
  (ii) $\Leftrightarrow$ (iii) is Proposition~\ref{thm:EffectiveEpi=RegularEpi}.

  (iii) $\implies$ (iv) consider a regular epimorphism $ f\from{X\to Y}$, expressed as the coequalizer of maps $u$, $v\from{W\to X}$\*. Assume $f=m\tilde{f}$ where $m$ is a monomorphism. Then by the universal property of $m$ we have $\tilde{f}u=\tilde{f}v$. Since $f$ is the coequalizer of $u$ and $v$, we have a unique $n\from Y\to M$ such that $nm=1_Y$, from which it follows that $m$ is an isomorphism.

  (iv) $\Leftrightarrow$ (v) is Proposition~\ref{thm:ExtremalEpi-StrongEpi}.
\end{proof}

\begin{proposition}[Strong epi implies regular epi in regular category]\label{thm:RegCat:StrongImpliesRegular}
  In a regular category (\ref{def:RegularCategory}), any strong epimorphism is a regular epimorphism.
\end{proposition}
\begin{proof}
  \index{epimorphism!strong}%
  \index{strong epimorphism}%
  Consider a strong epimorphism $f$ and factor it as a regular epimorphism $\varepsilon$ followed by a monomorphism $\mu$. This yields the following commutative square of solid arrows.
  \[
    \xymatrix@R=5ex@C=3em{
    X \ar[d]_-{f} \ar[r]^-{\varepsilon} &
    M \ar@{ >->}[d]^-{\mu} \\
    Y \ar@{=}[r] \ar@{.>}[ru]^-{\varphi} &
    Y}
  \]
  Since $f$ is a strong epimorphism, a dotted filler $\varphi$ exists. Now $\mu$ is both a monomorphism and a split epimorphism, which makes it an isomorphism. It follows that $f$ is a regular epimorphism.
\end{proof}

\begin{example}[An epimorphism which is not regular]
  \label{exa:Z>->Q}%
  In the category $\URngs$ of unital rings the inclusion $i\from \ZNr\to \QNr$ is an epimorphism, even though it is not a surjective ring homomorphism. To see this, consider morphisms of rings $f,g\from {\QNr\to R}$ satisfying $f\Comp i=g\Comp i$. Writing $q\in \QNr$ as a quotient $\tfrac{a}{b}$ with $a$, $b\in \ZNr$, we find:
  \begin{align*}
    f(q) & =f(\tfrac{a}{b})=f(\tfrac{1}{b}a)=f(\tfrac{1}{b})f(a)=f(\tfrac{1}{b})g(a)=f(\tfrac{1}{b})g(b\tfrac{a}{b})=f(\tfrac{1}{b})g(b)g(\tfrac{a}{b}) \\
         & =f(\tfrac{1}{b})f(b)g(\tfrac{a}{b})=f(\tfrac{b}{b})g(\tfrac{a}{b})=f(1)g(\tfrac{a}{b})=g(1)g(\tfrac{a}{b})=g(1\tfrac{a}{b})=g(q).
  \end{align*}
  Since a regular epimorphism in $\URngs$ is a quotient function of underlying sets, $i$ is an epimorphism which is not regular.
\end{example}

Example \ref{exa:Z>->Q} also shows that a morphism which is both epic and monic need not be an isomorphism.

\begin{example}[Epimorphisms in $\Grps$]
  \label{exa:EpiInGrps}
  A surjective group homomorphism is an epimorphism of underlying sets, hence is an epimorphism of groups. The converse is also true. But the proof is a bit of a challenge; see Exercise \ref{exe:Epimorphisms-In-Grps}.
\end{example}

\begin{exercises}

\begin{exercise}[`Being an epimorphism' is a colimit property]
  \label{exe:EpimorphismViaPullback}
  In any category show that $f\from X\to Y$ is an epimorphism if and only if the diagram below has the pushout property. %
  \index{epimorphism!via pushout property}\index{monomorphism!via limit property}%
  \begin{equation*}
    \xymatrix@R=5ex@C=4em{
    X \ar[r]^-{f} \ar[d]_{f} &
    Y \ar@{=}[d] \\
    Y \ar@{=}[r] &
    Y
    }
  \end{equation*}
  Then formulate and prove the dual of this statement.
\end{exercise}

\begin{exercise}[Extremal epimorphism - properties]
  \label{exe:ExtremalEpi-Props}
  Show that following holds in an arbitrary category:
  \begin{thmlist}
    \item If $g\Comp f$ is an extremal epimorphism, then $g$ is an extremal epimorphism.
    \item A morphism $f$ is an isomorphism if and only if it is an extremal epimorphism and a monomorphism.
  \end{thmlist}
\end{exercise}

\begin{exercise}[Strong epimorphism - properties]
  \label{exe:StrongEpi-Props}
  Show that following holds in an arbitrary category:
  \begin{thmlist}
    \item If $g\Comp f$ is a strong epimorphism, and $f$ is an epimorphism, then $g$ is a strong epimorphism.
    \item Any composite of two strong epimorphisms is a strong epimorphism.
  \end{thmlist}
\end{exercise}

\begin{exercise}[Strong epi / mono factorization unique]
  \label{exe:StrongEpi/MonoFactorization->Unique}
  In the commutative square below, assume that $e$ and $\varepsilon$ are strong epimorphisms, and that $m$ and $\mu$ are monomorphisms.
  
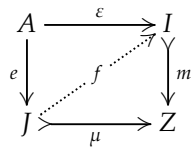
\begin{equation*}
    \xymatrix@R=5ex@C=4em{
    A \ar[r]^-{\varepsilon} \ar[d]_{e} &
    I \ar@{{ >}->}[d]^{m} \\
    J \ar@{{ >}->}[r]_-{\mu} \ar@{.>}[ru]|-{\ f\ }&
    Z
    }
  \end{equation*}
  Show that there is an isomorphism $f\from J\to I$ which is unique with the property that it renders the entire diagram commutative. %
  \index{image!factorization: strong epi/mono}%
\end{exercise}

\begin{exercise}[Extremal epimorphism characterization of join]
  \label{exe:ExtremalEpiCharacterizationOfJoin}%
  \label{thm:ExtremalEpiCharacterizationOfJoin}
  In any category, given two subobjects $m\from {M\to X}$ and $n\from {N\to X}$ of an object $X$, a monomorphism $i\from {I\to X}$ represents the join of $m$ and $n$ if and only if the following two conditions hold:
  \begin{enumerate}[(i)]
    \item $m$ and $n$ factor through $i$ via $m'\from M\to I$ and $n'\from N\to I$, and
    \item the morphisms $m'$ and $n'$ are jointly extremal-epimorphic. \NoProof
  \end{enumerate}
\end{exercise}

\begin{exercise}[Extremal epi vs.\ strong epi]
  \label{exe:ExtremalEpi+Pullbacks->StrongEpi}
  Prove (\ref{thm:ExtremalEpi-StrongEpi}): Every strong epimorphism is extremal. In a category with pullbacks, every extremal epimorphism is a strong epimorphism.
\end{exercise}

\begin{exercise}[Exact diagram of sets]
  \label{exe:ExactDiagramSets}
  A diagram of sets %
  \index{exact!diagram of sets}
  \begin{equation*}
    \xymatrix@R=5ex@C=4em{
    A \ar[r]^-{k} &
    X \ar@<+0.5ex>[r]^-{u} \ar@<-0.5ex>[r]_-{v} &
    Y
    }
  \end{equation*}
  is said to be \Defn{exact} if, for every $x\in X$ with $u(x)=v(x)$ there exists a unique $a\in A$ such that $k(a)=x$. Show the following
  \begin{enumerate}[(i)]
    \item The diagram of sets above is exact if and only if $k$ is an equalizer of $u$ and $v$ in the category $\Sets$ of sets.
    \item In the category $\Grps$ of groups, show that $\kappa\from K\to X$ is a kernel of $f\from X\to Y$ if and only if the diagram of underlying sets and functions below is exact.
          \begin{equation*}
            \xymatrix@R=5ex@C=4em{
            K \ar[r]^-{\kappa} &
            X \ar@<+0.5ex>[r]^-{f} \ar@<-0.5ex>[r]_-{\ZeroMap} &
            Y
            }
          \end{equation*}
  \end{enumerate}
\end{exercise}

\begin{exercise}[Hom-view of effective epimorphism]
  \label{exe:EffectiveEpi-HomView}%
  In a category which admits kernel pairs show that the following conditions of a morphism $f\from X\to Y$ are equivalent. %
  \index{effective!epimorphism - Hom-view}
  \begin{enumerate}[(i)]
    \item $f$ is an effective epimorphism; see (\ref{def:Regular/EffectiveEpi}).
    \item For every object $T$ in $\Ctgry{X}$ this diagram of Hom-sets is an exact diagram of sets:
          \begin{equation*}
            \xymatrix@C=5em{
            \Hom{Y}{T}\ar[r]^-{\Hom{f}{T}} &
            \Hom{X}{T} \ar@<.5ex>[r]^-{\Hom{\PrjctnOnto{1}}{T}} \ar@<-.5ex>[r]_-{\Hom{\PrjctnOnto{2}}{T}} &
            \Hom{\KrnlPr{f}}{T}
            }
          \end{equation*}
    \item \cite[p.~4.3]{DGQuillen1967}\quad For every $\alpha\from X\to T$ satisfying
          \begin{equation*}
            \left[ \text{for all}\quad u,v\from A\to X\quad \text{with}\ fu=fv\right] \implies \left[\alpha u=\alpha v\right]
          \end{equation*}
          there exists a unique $\beta\from Y\to T$ with $\beta f=\alpha$.
  \end{enumerate}
\end{exercise}

\begin{exercise}
  In an arbitrary  category with pullbacks,  assume that the front and back faces in the diagram below are pullbacks.
  
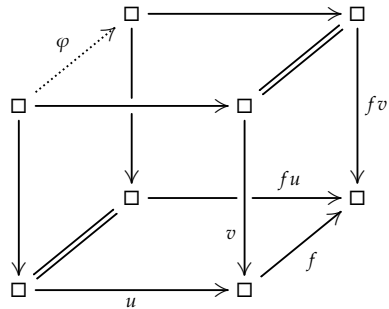
\begin{equation*}
    \xymatrix@R=5ex@C=3em{
    & \DiagObj \ar[rr] \ar[dd]|\hole &&
    \DiagObj \ar[dd]^{fv} \\
    \DiagObj \ar[rr] \ar[dd] \ar@{.>}[ru]^{\varphi} &&
    \DiagObj \ar@{=}[ru] \ar[dd]_(.7){v} \\
    & \DiagObj \ar[rr]|\hole^(.7){fu} &&
    \DiagObj \\
    \DiagObj \ar@{=}[ru] \ar[rr]_-{u} &&
    \DiagObj \ar[ru]_{f}
    }
  \end{equation*}
  \begin{enumerate}[(i)]
    \item Show that there exists a unique map $\varphi$ which renders the entire cube commutative.
    \item If $f$ is a monomorphism, show that $\varphi$ is an isomorphism.
    \item \cite[p.~173]{Bourn-Gran-CategoricalFoundations}\quad If  $u=v$  is a split epimorphism, and $\varphi$ is an isomorphism, show that $f$ is a monomorphism.
    \item \cite[p.~176]{Bourn-Gran-CategoricalFoundations}\quad In a semiabelian category, assume  $ u=v$  is a normal epimorphism. Then $\varphi$ is an isomorphism, if and only if $f$ is a monomorphism.
  \end{enumerate}
\end{exercise}
\end{exercises}
\section{Factorization Systems}
\label{sec:FactorizationSystems}

Here we collect background relevant to the theory of factorizations in general, and image factorizations in particular. Thus, in a category $\Ctgry{C}$ we consider two classes of morphisms $\EuScript{E}$ and $\EuScript{M}$, assuming formally the role played by surjections / injections in image factorizations of functions of sets. Then, an arbitrary morphism $f\from X\to Y$ in $\Ctgry{C}$ should be a composite $f=me$ with $e\in \EuScript{E}$ and $m\in \EuScript{M}$ in a functorial manner; that is, any commutative diagram of solid arrows in $\Ctgry{C}$ of the kind displayed below should admit a unique filler as indicated.
\begin{equation*}
  \xymatrix@R=5ex@C=4em{
  A \ar[r]_-{e} \ar@/^3ex/[rr]^-{f} \ar[d]_{a} &
  I \ar[r]_-{m} \ar@{.>}[d]_{\varphi} &
  B \ar[d]^{b} \\
  X \ar[r]^-{e'} \ar@/_3ex/[rr]_-{g} &
  J \ar[r]^-{n} &
  Y
  }
\end{equation*}
A composite $f=me$ is called an $(\EuScript{E},\EuScript{M})$-factorization of $f$. Functoriality of such decompositions is elegantly achieved (\ref{thm:OFS-Factorizations->Functorial}) by requiring the pair $(\EuScript{E},\EuScript{M})$ to form an orthogonal prefactorization system, see (\ref{def:OrthogonalPrefactorizationSystem}).

\begin{definition}[Left / right orthogonal]
  \label{def:Left/Right-Orthogonal}
  In a category $\Ctgry{C}$, a morphism $e\from A\to B$ is said to be \Defn{left orthogonal} to $m\from U\to V$, and $m$ is said to be \Defn{right orthogonal} to $e$ if every one of these commutative squares in $\Ctgry{C}$ admits a unique filler $\varphi$ as indicated. %
  \index{left!orthogonal}\index{right!orthogonal}
  \begin{equation*}
    \xymatrix@R=5ex@C=4em{
    A \ar[d]_{e} \ar[r]^-{u} &
    U \ar[d]^{m} \\
    B \ar[r]_-{v} \ar@{.>}[ru]|-{\ \varphi\ } &
    V
    }
  \end{equation*}
\end{definition}

\begin{definition}[Orthogonal prefactorization system]
  \label{def:OrthogonalPrefactorizationSystem}
  A pair of classes of morphisms $(\EuScript{E},\EuScript{M})$ in a category $\Ctgry{C}$ is said to form an \Defn{orthogonal prefactorization system} (OPFS) if these two conditions are satisfied: %
  \index{orthogonal!prefactorization system}\index{OPFS - orthogonal prefactorization system}%
  \begin{enumerate}[(i)]
    \item $\EuScript{E}$ is the \Defn{left orthogonal complement} of $\EuScript{M}$, meaning: $e\in \EuScript{E}$ if and only if $e$ is left orthogonal to every $m\in \EuScript{M}$. %
          \index{left!orthogonal complement}\index{orthogonal!complement}
    \item $\EuScript{M}$ is the \Defn{right orthogonal complement} of $\EuScript{E}$, meaning: $m\in \EuScript{M}$ if and only if $m$ is right orthogonal to every $e\in \EuScript{E}$. %
          \index{right!orthogonal complement}
  \end{enumerate}
\end{definition}

\begin{definition}[Orthogonal factorization system]
  \label{def:OrthogonalFactorizationSystem}
  An \Defn{orthogonal factorization system} (OFS) in a category $\Ctgry{C}$ is given by an orthogonal prefactorization system $(\EuScript{E},\EuScript{M})$ such that every morphism $\from X\to Y$ in $\Ctgry{C}$ admits are factorization $f=me$, with $m\in \EuScript{M}$ and $e\in\EuScript{E}$. %
  \index{orthogonal!factorization system}\index{OFS}%
\end{definition}

\begin{proposition}[OFS-factorizations are functorial]
  \label{thm:OFS-Factorizations->Functorial}
  If $(\EuScript{E},\EuScript{M})$ is an OFS in a category $\Ctgry{C}$, then $(\EuScript{E},\EuScript{M})$-factorizations of morphisms are functorial.
\end{proposition}
\begin{proof}
  The need for a functorial filler on the left below may be rearranged into the lifting problem presented on the right.
  
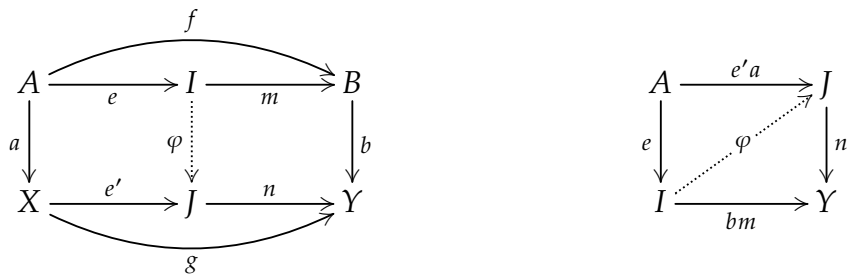
\begin{equation*}
    \xymatrix@R=5ex@C=4em{
    A \ar[r]_-{e} \ar@/^3ex/[rr]^-{f} \ar[d]_{a} &
    I \ar[r]_-{m} \ar@{.>}[d]_{\varphi} &
    B \ar[d]^{b} &&
    A \ar[d]_{e} \ar[r]^-{e'a} &
    J \ar[d]^{n} \\
    X \ar[r]^-{e'} \ar@/_3ex/[rr]_-{g} &
    J \ar[r]^-{n} &
    Y &&
    I \ar[r]_-{bm} \ar@{.>}[ru]|-{\ \varphi\ } &
    Y
    }
  \end{equation*}
  Thus $\varphi$ exists and is unique as required.
\end{proof}

In particular, (\ref{thm:OFS-Factorizations->Functorial}) tells us that $(\EuScript{E},\EuScript{M})$-factorizations of a given map are unique up to unique isomorphism.

Orthogonal prefactorization systems are easy to construct:

\begin{proposition}[Construction of OPFS]
  \label{thm:OPFS-Construct}%
  Given any class $\mathcal{W}$ of morphisms in a category $\Ctgry{C}$, each of the following two constructions yields an OPFS:
  \begin{enumerate}[(i)]
    \item $(\EuScript{E},\EuScript{M})$, with $\EuScript{E}$ the left orthogonal complement of $\mathcal{W}$, and $\EuScript{M}$ be the right orthogonal complement of $\EuScript{E}$.
    \item $(\EuScript{E},\EuScript{M})$, with $\EuScript{M}$ the right orthogonal complement of $\mathcal{W}$, and $\EuScript{E}$ be the left orthogonal complement of $\EuScript{M}$. \NoProof
  \end{enumerate}
\end{proposition}

While it is easy to check the truth of (\ref{thm:OPFS-Construct}) directly, we point out that this is an example of a Galois correspondence. Motivated by image factorizations encountered for morphisms between groups, Lie algebras, objects in an abelian category, etc., we are particularly interested in an OPFS $(\EuScript{E},\EuScript{M})$ in which $\EuScript{M}$ is the class of monomorphisms. This suggests (see also~\eqref{def:StrongEpimorphism}):

\begin{definition}[Strong epimorphism]
  \label{def:Strong-Epimorphism}
  A morphism $f$ in a category $\Ctgry{C}$ is called a \Defn{strong epimorphism} if it is left orthogonal to every monomorphism in $\Ctgry{C}$.
\end{definition}

According to (\ref{thm:OPFS-Construct}), we obtain in an arbitrary category $\Ctgry{C}$ an OPFS $(\EuScript{E},\EuScript{M})$ in which $\EuScript{E}$ is the left orthogonal complement of the class of monomorphisms in $\Ctgry{C}$, and $\EuScript{M}$ is the right orthogonal complement of $\EuScript{E}$. In general, there is then nothing that tells us that every strong epimorphism is actually an epimorphism, see Section \ref{sec:Surjectivity} nor is there any guarantee that every map in $\EuScript{M}$ is a monomorphism. However, all of this goes well under mild additional assumptions, which is what we are going to explain now.

\begin{proposition}[OPFS: (strong epis,monos) in cat with coequalizers]
  \label{thm:OPFS-(StrongEpis,Monos)}%
  In a category $\Ctgry{C}$ which admits coequalizers, an OPFS is given by $\EuScript{M}$ the class of monomorphisms of $\Ctgry{C}$, and $\EuScript{E}$ the strong epimorphisms.
\end{proposition}

\begin{corollary}[OPFS: (strong epis,monos) and strong epi is epi]
  \label{thm:OPFS-(StrongEpis,Monos)-StrongEpi->Epi}
  In a category which admits equalizers and coequalizers an OPFS $(\EuScript{E},\EuScript{M})$ is given by strong epimorphisms and monomorphisms. Moreover, every strong epimorphism is epic. \NoProof
\end{corollary}

\begin{corollary}[OPFS: in reflective subcategory]
  \label{thm:OPFS-ReflectiveSubCat}
  In a finitely complete and cocomplete category $\Ctgry{C}$, let $\Ctgry{R}$ be a full and replete reflective subcategory. Then the following hold:
  \begin{enumerate}[(i)]
    \item The OPFS $(\EuScript{E},\EuScript{M})$ on $\Ctgry{C}$ given by (strong epimorphisms, monomorphisms) restricts to an OPFS on $\Ctgry{R}$.
    \item If $(\EuScript{E},\EuScript{M})$ is a factorization system on $\Ctgry{C}$, and $\Ctgry{R}$ is closed under subobjects, then the inclusion functor $G\from \Ctgry{R}\to \Ctgry{C}$ creates image factorizations in $\Ctgry{R}$.
  \end{enumerate}
\end{corollary}
\begin{proof}
  (i) Both categories are finitely complete and cocomplete. %
  So, in each, an OPFS is given by (strong epimorphisms,monomorphism). Further, a map $m$ is a monomorphism in $\SACtgry{R}$ if and only if it is a monomorphism in $\SACtgry{X}$.
  So, for $X,Y$ are in $\Ctgry{R}$, the strong epimorphisms in $\HomIn{\Ctgry{R}}{X}{Y}$ and in $\HomIn{\Ctgry{C}}{X}{Y}$ are the same, because $\Ctgry{R}$ is full in $\SACtgry{X}$. - This proves (i)

  (ii) If $f\from X\to Y$ is a morphism in $\Ctgry{R}$, let $X \XRA{q} I \XRA{i} Y$ be is image factorization in $\SACtgry{X}$. Then $i$ and, hence $q$ are in $\Ctgry{R}$ because $\SACtgry{R}$ is closed under subobjects. By (i), $q$ is a strong epimorphism in $\Ctgry{R}$, and $i$ is a monomorphism in $\Ctgry{R}$. So, the image factorization of $f$ in $\SACtgry{X}$ is also an image factorization of $f$ in $\Ctgry{R}$.
\end{proof}

We turn to the question under which conditions the kind of OPFS encountered in (\ref{thm:OPFS-(StrongEpis,Monos)-StrongEpi->Epi}) is actually an OFS. Thus, we require of every morphism $f\from X\to Y$ that it admits a factorization
\begin{equation*}
  \xymatrix@R=5ex@C=4em{
  X \ar@{-{>>}}[r]^-{q} &
  \Img{f} \ar@{{ >}->}[r]^-{m} &
  Y
  }
\end{equation*}
with $q$ in $\EuScript{E}$ and $m$ in $\EuScript{M}$. Such a factorization is called an \Defn{image factorization} of $f$, and the \Defn{image} of $f$ is the subobject of $Y$ represented by $m$. %
\index{image!of a morphism}\index{image!factorization}

There are essentially two mutually complementary approaches:
\begin{enumerate}[(1)]
  \item Given a morphism $f\from X\to Y$, assume that the intersection $m\from I\to Y$ of all subobjects of $Y$ through which $f$ factors exists. Then $f=mp$, and $p$ can be shown to be an extremal epimorphism, hence an epimorphism which is strong. So $f=mp$ is the desired image factorization of $f$, which is then automatically functorial by (\ref{thm:OFS-Factorizations->Functorial}).
  \item We perform a construction on $f$ which achieves $f=\mu p$, with $p$ a strong epimorphism. If we can then show that $\mu$ is a monomorphism, we have the desired image factorization.
\end{enumerate}
The following two variations of approach (2) are particularly relevant to us:
\begin{enumerate}[({2}.a)]
  \item\label{item:ImagesViaRegularity} Take the kernel pair of $f$ and $p$ its coequalizer. This approach is classical~\cite{Barr-Grillet-vanOsdol}, and probably the best one can do in a regular category if the KSG-property is not guaranteed. For this, existence of coequalizers of kernel pairs suffices. See~\ref{thm:ImageInRegular}.
  \item Take the kernel of $f$ and $p$ its cokernel. This works in varieties such as module categories, groups, Lie algebras. It is the approach taken here in section \ref{sec:ImageFactorizations} in the context of a normal category, and relies heavily on the KSG-axiom which asserts that, in a split short exact sequence, the middle object is generated by the subobjects at its ends. We prefer this approach over~(2.\ref{item:ImagesViaRegularity}) because it is technically less involved.
\end{enumerate}

\subsection[Image Factorization in Regular Categories]{Image Factorization in Regular Categories}
\label{subsec:ImageFactorization-RegCats}

We explain how regular categories provide a kind of image factorization which is particularly well adapted to varieties of universal algebras; see \cite[1.14]{MGran2021-RegCats}.

\begin{definition}[Regular category]
  \label{def:RegularCategory}%
  A finitely complete category with coequalizers of kernel pairs is \Defn{regular} if pullbacks preserve regular epimorphisms. %
  \index{regular!category}\index{category!regular}%
\end{definition}

We show that this is equivalent to the existence of pullback-stable image factorizations\footnote{Here things are really different from what happens in Proposition~\ref{thm:ExtremalEpi->Cokernel}. So we need a new proof.}:

Given a morphism $f\from X\to Y$, we rely on the following construction:
\begin{equation*}
  \xymatrix@R=5ex@C=5em{
  \KrnlPr{f} \ar@<1ex>[r]^-{\pi_{1}} \ar@<-1ex>[r]_-{\pi_{2}}  &
  X \ar[l]|-{\ e\ } \ar[r]^-{q\DefEq \CoEq{\pi_{1}}{\pi_{2}}} \ar[rd]_{f} &
  I \ar@{.>}[d]^{m} \\
  && Y
  }
\end{equation*}

\begin{proposition}[Images in regular categories]
  \label{thm:ImageInRegular}%
  In a regular category, the map $m$ in the diagram above is a monomorphism.
\end{proposition}
\begin{proof}
  Consider the situation where $r$ in the diagram diagram below is an epimorphism.
  \begin{equation*}
    \xymatrix@R=5ex@C=4em{
    \KrnlPr{f} \ar@{.>}[d]_-{r} \ar@<1ex>[r]^-{\pi_{1}} \ar@<-1ex>[r]_-{\pi_{2}} &
    X \ar@{-{>>}}[d]_{q} \ar[l]|-{\ e\ } \ar@{->}[r]^-{f} &
    Y \ar@{=}[d] \\
    \KrnlPr{m} \ar@<1ex>[r]^-{\pi'_{1}} \ar@<-1ex>[r]_-{\pi'_{2}} &
    I \ar[r]_-{m} \ar[l]|-{e'} & Y
    }
  \end{equation*}
  Then $\pi_1'=\pi_2'$, because $\pi'_{1}r=q\pi_{1}=q\pi_{2}=\pi'_{2}r$. So, $m$ is a monomorphism by (\ref{thm:KernelPair-Monos}). To see that $r$ is epic consider the diagram below in which every square is a pullback:
  \begin{equation*}
    \xymatrix@!0@R=3em@C=6em{ & \KrnlPr{f} \ar[d]_-{(\pi_1,\pi_2)} \ar[r]^-{r_1} & P \ar[r]^-{r_2} \ar[d] & \KrnlPr{m} \ar[d]^-{(\pi_1',\pi_2')}\\
    & X\times X \ar[ld]_-{\pi_2} \ar[r]^-{1_X\times q} & X\times I \ar[ld]^-{\pi_2} \ar[rd]_-{\pi_1} \ar[r]^-{q\times 1_X} & I\times I \ar[rd]^-{\pi_1}\\
    X \ar[r]_-q & I && X \ar[r]_-q & I}
  \end{equation*}
  All horizontal arrows are pulled back from the regular epimorphism $q$; so they are regular epimorphisms. Thus $r=r_2r_1$ is an epimorphism as a composite of two regular epimorphisms.
\end{proof}

\begin{corollary}[Orthogonal factorization system in regular category]
  \label{thm:RegularEpi/MonoFactorizationSystemInRegCat}%
  In a regular category, an orthogonal factorization system $(\EuScript{E},\EuScript{M})$ is given by $\EuScript{E}$ the class of regular, hence strong, epimorphisms, and $\EuScript{M}$ the class of monomorphisms. \NoProof
\end{corollary}

\begin{corollary}[Recognizing a regular category]
  \label{thm:RegularCat-Recognize}
  A finitely complete category with coequalizers of kernel pairs is regular if and only if it satisfies the following two properties: %
  \begin{enumerate}
    \item Every morphism $f$ admits a factorization $f=\mu \Comp \varepsilon$ in which $\varepsilon$ is a regular epimorphism and $\mu$ is a monomorphism.
    \item Every such factorization is preserved under pullbacks. \NoProof
  \end{enumerate}
\end{corollary}

Note that in this setting, regular epimorphisms and normal epimorphisms need not coincide. On the other hand, extremal, strong, effective and regular epimorphisms do; see Section~\ref{sec:Surjectivity}.
\section[Reflective Subcategories]{Reflective Subcategories}
\label{sec:ReflectiveSubCats}

\begin{definition}[Reflective subcategory]
  \label{def:ReflectiveSubCat}%
  Given a category $\Ctgry{C}$, a subcategory $\Ctgry{R}$ of $\Ctgry{C}$ is called \Defn{reflective} if the inclusion functor $G\from \Ctgry{R}\to \Ctgry{C}$ has a left adjoint, say $T\from \Ctgry{C}\to \Ctgry{R}$.  The functor $T$ is called a \Defn{reflector} of $\Ctgry{C}$ in $\Ctgry{R}$. %
  \index{reflective subcategory}\index{reflector}%
\end{definition}

When working with a reflection of $\Ctgry{C}$ in $\Ctgry{R}$, we will always assume that $\Ctgry{R}$ is full and replete in $\Ctgry{C}$.  As is explained in \cite[IV.3]{SMacLane1998}, such a reflection may equivalently be given by a pair $(T,\AdjUnit)$, where
\begin{enumerate}
  \item $T\from \Ctgry{C}\to \Ctgry{R}$ is a functor, and
  \item $\AdjUnit\from \IdMapOn{\Ctgry{C}}\Rightarrow GT$ is a natural transformation with this universal property: For every pair of objects $X$ in $\Ctgry{C}$, and $A$ in $\Ctgry{R}$, every morphism $f\from X\to GA$ in $\Ctgry{C}$ may be uniquely as shown in the  diagram on the left:
        \begin{equation*}
          \xymatrix@R=5ex@C=4em{
          X \ar[r]^-{\AdjUnitOn{X}} \ar[d]_{f} &
          GTX \ar@{.>}[dl]^{\bar{f}} \\
          GA
          }\qquad\qquad
          \xymatrix@R=5ex@C=4em{
          X \ar[r]^-{\AdjUnitOn{X}} \ar[d]_{f} &
          TX \ar@{.>}[dl]^{\bar{f}} \\
          A
          }
        \end{equation*}
\end{enumerate}
The natural transformation $\eta$ is the \Defn{adjunction unit} of the reflection. Since $G$ is an inclusion functor, we will often simplify the exposition by omitting it. Thus, the triangle on the left above turns into the less precise but more economical triangle on the right. For $A$ in $\Ctgry{R}$, the map $\AdjUnitOn{A}\from A\to TA$ is an isomorphism; see (\ref{exe:ReflectorUnit-Iso}). Consequently, the commutative square below explains why $T$ is idempotent in the strong sense that, for every object $X$ in~$\Ctgry{C}$, the adjunction unit is an isomorphism $T\AdjUnitOn{X}=\AdjUnitOn{TX}\from TX\to T^2$:
\begin{equation*}
  \xymatrix@R=5ex@C=4em{
  X \ar[r]^-{\AdjUnitOn{X}} \ar[d]_-{\AdjUnitOn{X}} &
  TX \ar[d]^-{\AdjUnitOn{TX}}_-{\cong} \\
  TX \ar[r]_-{T\AdjUnitOn{X}}^-{\cong} &
  T^2X
  }
\end{equation*}

Limits and colimits of a category $\Ctgry{C}$ are strongly related to limits and colimits in any reflective subcategory. Here are the details.

\begin{theorem}[Colimits in a reflective subcategory]
  \label{thm:ReflectiveSubCat-Colimits}
  Let $T\from \Ctgry{C}\to \Ctgry{R}$ be a reflector, and let $G\from \Ctgry{R}\to \Ctgry{C}$ be the inclusion. For a functor $\phi\from D\to \Ctgry{C}$ from a small category $D$ into $\Ctgry{C}$ the following hold:
  \begin{enumerate}[(i)]
    \item $T$ \emph{constructs colimits in $\Ctgry{R}$}: If $\phi\from D\to \Ctgry{R}$, then $\phi\Rightarrow \CoLimOf{G\phi} \XRA{\AdjUnit}  T\CoLimOf{G\phi}$ is a colimit cocone in $\Ctgry{R}$, provided $\CoLimOf{G\phi}$ exists in $\Ctgry{C}$.
    \item $T$ \emph{preserves colimits}: If $\gamma\from \phi\Rightarrow C$ is a colimit cocone for $\phi$, then $T\Comp \gamma \from T\Comp \phi \Rightarrow TC$ is a colimit cocone in $\Ctgry{R}$.
  \end{enumerate}
\end{theorem}
\begin{proof}
  (i) follows via the universal property of $\AdjUnit$. (ii) is a property of all functors which admit a right adjoint.
\end{proof}

\begin{theorem}[Limits in a reflective subcategory]
  \label{thm:ReflectiveSubCat-Limits}
  Let $T\from \Ctgry{C}\to \Ctgry{R}$ be a reflector, and let $G\from \Ctgry{R}\to \Ctgry{C}$ be the inclusion. For a functor $\phi\from D\to \Ctgry{R}$ from a small category $D$ into $\Ctgry{R}$ the following hold:
  \begin{enumerate}[(i)]
    \item \emph{$G$ preserves limits}: If $\lambda\from L\Rightarrow \phi$ is a limit cone for $\phi$ in $\Ctgry{R}$, then its is also a limit cone for $G\phi$ in $\Ctgry{C}$.
    \item $G$ \emph{creates limits in $\Ctgry{R}$}: If $\phi\from D\to \Ctgry{R}$, and $\lambda\from L\to G\phi$ is a limit cone in $\Ctgry{C}$, then $L$ is in $\Ctgry{R}$, and $\lambda\from L\to G\phi=\phi$ is a limit cone in $\Ctgry{R}$.
  \end{enumerate}
\end{theorem}
\begin{proof}
  (i) holds because the inclusion of $\Ctgry{R}$ in $\SACtgry{X}$ is full and faithful.
  (ii) is a property of all functors which admit a left adjoint.
\end{proof}

We turn to the question: How do we recognize whether a given subcategory of $\Ctgry{C}$ is reflective? Discussing this question benefits from orthogonality between morphism and object in the following sense.

\begin{definition}[Categorical orthogonality: morphism / object]
  \label{def:Orthogonality-Morphism/Object}
  In a category $\Ctgry{C}$ an object $Z$ is \Defn{orthogonal} to a morphism $u\from A\to B$ if the induced function of morphism sets below is a bijection\index{orthogonal}
  \begin{equation*}
    u^*\from\HomIn{\Ctgry{C}}{B}{Z} \to\HomIn{\Ctgry{C}}{A}{Z},\qquad f \mapsto f\Comp u.
    \qquad\qquad \vcenter{\xymatrix{B \ar[r]^-f & Z  \\ A \ar[u]^-{u} \ar@{.>}[ru]_-{f\circ u} }}
  \end{equation*}
\end{definition}

\begin{definition}[Orthogonal complement]
  \label{def:OrthogonalComplement}
  The \Defn{orthogonal complement of a class of morphisms} $M$ in $\Ctgry{C}$ is $M^{\bot}$, the collection of all objects $Z$ which are orthogonal to every morphism in $M$. The \Defn{orthogonal complement of a class of objects} $O$ in $\Ctgry{C}$ is $O^{\bot}$, the collection of all morphisms which are orthogonal to every object in $O$. %
  \index{orthogonal complement!of a class of morphisms}%
  \index{orthogonal complement!of a class of objects}%
\end{definition}

\begin{definition}[Orthogonal pair]
  \label{def:OrthogonalPair}
  An \Defn{orthogonal pair} $(M,O)$ in a category $\EuScript{C}$ is given by a class of morphisms $M$ and a class of objects $O$ in $\Ctgry{C}$ such that %
  \index{orthogonal!pair}%
  \begin{equation*}
    M^{\bot}=O \quad\text{and}\quad M=O^{\bot}
  \end{equation*}
\end{definition}

To explain why orthogonal pairs are relevant to the discussion of reflectors on a category, we introduce the following concepts:

\begin{definition}[Equivalence / local object of reflections]
  \label{def:Equivalence/LocalObject}
  Given a reflection $(T,\AdjUnit)$ of $\Ctgry{C}$ in a subcategory $\Ctgry{R}$ of $\Ctgry{C}$, we say:
  \begin{enumerate}[(i)]
    \item A morphism $f\from X\to Y$ in $\Ctgry{C}$ is a \Defn{$T$-equivalence} if $Tf$ is an isomorphism. %
          \index{equivalence!of a reflector}%
    \item An object $Z$ in $\Ctgry{C}$ is \Defn{$T$-local} if $\AdjUnitOn{Z}\from Z\to TZ$ is an isomorphism. %
          \index{local!object w.r. to a reflector}%
  \end{enumerate}
\end{definition}

\begin{proposition}[Reflector yields orthogonal pair]
  \label{thm:Reflector->OrthogonalPair}
  If $T\from \Ctgry{C}\to \Ctgry{R}$ is a reflector from a category $\Ctgry{C}$ to a full and replete subcategory $\Ctgry{R}$, let $\mathcal{W}$ denote the class of $T$-equivalences, and let $\mathcal{L}$ be the class of $T$-local objects. Then $\mathcal{L}=\Objcts{\Ctgry{R}}$ is the object class of $\Ctgry{R}$, and $(\mathcal{W},\mathcal{L})$ is an orthogonal pair.
\end{proposition}
\begin{proof}
  We show first that every object $Z$ in $\Ctgry{R}$ is $T$-local. Indeed, the universal property of $\AdjUnitOn{Z}$ yields a left inverse $\gamma$ of $\AdjUnitOn{Z}$ via the triangle on the left,
  $$
    
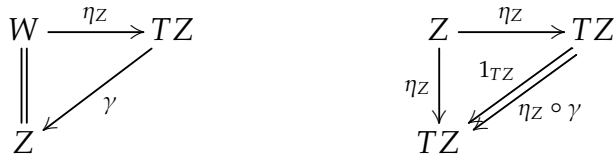
\xymatrix@R=5ex@C=3em{
    W \ar@{=}[d] \ar[r]^-{\AdjUnitOn{Z}} &
    TZ \ar[dl]^{\gamma} &&
    Z \ar[r]^-{\AdjUnitOn{Z}} \ar[d]_{\AdjUnitOn{Z}} &
    TZ \ar@<-0.1pt>[dl]_{\IdMapOn{TZ}} \ar@<3pt>[dl]^{\AdjUnitOn{Z}\Comp \gamma} \\
    Z &&&
    TZ
    }
  $$
  Via the identity $(\AdjUnitOn{Z}\Comp \gamma)\Comp \AdjUnitOn{Z}= \AdjUnitOn{Z}$,  the diagram on the right explains why $\gamma$ is also a right inverse of $\AdjUnitOn{Z}$. So $\Objcts{\Ctgry{R}}\subseteq \mathcal{L}$. Conversely, if $X$ in $\Ctgry{C}$ is $T$-local, then $\AdjUnitOn{X}\from X\to TX$ is an isomorphism to an object in $\EuScript{R}$. So, $X$ itself is in $\Ctgry{R}$ because $\Ctgry{R}$ is replete; i.e.\ $\mathcal{L}\subseteq \Objcts{\Ctgry{R}}$.

  Next, we show that every $T$-local object $Z$ is orthogonal to every  $T$-equivalence $u\from A\to B$. Given $v\from A\to Z$, the commutative diagram below shows that $v=\AdjUnitOn{Z}(Tv)(Tu)^{-1}\AdjUnitOn{B}u$.
  \begin{equation*}
    \xymatrix@R=5ex@C=4em{
    Z \ar[d]_{\AdjUnitOn{Z}}^{\cong} &
    A \ar[d]_{\AdjUnitOn{A}} \ar[r]_-{u} \ar[l]^-{v} &
    B \ar[d]^{\AdjUnitOn{B}} \ar@{.>}@/_3ex/[ll]_-{w} \\
    TZ &
    TA \ar[l]_-{Tv} \ar[r]_-{\cong}^-{Tu} &
    TB \ar@{.>}@/^3ex/[ll]^-{Tw}
    }
  \end{equation*}
  To see that the factorization of $v$ through $u$ is unique, suppose $w\from B\to Z$ satisfies $wu=v$. This gives
  $$
    w\ =\ \AdjUnitOn{Z}^{-1}T(w)\AdjUnitOn{B}\ =\ \AdjUnitOn{Z}^{-1}T(v)T(u)^{-1}\AdjUnitOn{B},
  $$
  as required. This establishes $\mathcal{L}\subseteq \mathcal{W}^{\bot}$ and $\mathcal{W}\subseteq \mathcal{L}^{\bot}$.

  Now consider the case where $u\from A\to B$ is orthogonal to every $T$-local object $Z$. We claim that $u$ is a $T$-equivalence. Consider this reflector diagram:
  \begin{equation*}
    \xymatrix@R=6ex@C=4em{
    A \ar[r]^-{u} \ar[d]_{\AdjUnitOn{A}} &
    B \ar[d]^{\AdjUnitOn{B}} \ar@{->}[ld]|-{\ \bar{\gamma}\ }\\
    TA \ar@<-2pt>[r]_-{Tu} &
    TB \ar@{.>}@<-2pt>[l]_-{\gamma}
    }
  \end{equation*}
  First $\AdjUnitOn{A}$ factors uniquely through $u$ via $\bar{\gamma}$. Then $\bar{\gamma}$ factors uniquely through $\AdjUnitOn{B}$ via $\gamma\from TB\to TA$. We claim that $\gamma$ is the inverse of $Tu$. First, we find:
  \begin{equation*}
    (Tu)\gamma \AdjUnitOn{B} u = (Tu)\bar{\gamma} u = \AdjUnitOn{A} (Tu) = \AdjUnitOn{B}u
  \end{equation*}
  So, $(Tu)\gamma \AdjUnitOn{B}=\AdjUnitOn{B}$ because $u$ is orthogonal to $TB$. Then $(Tu)\gamma=\IdMapOn{TB}$ via the universal property of $\AdjUnitOn{B}$. Similarly,
  \begin{equation*}
    \gamma(Tu)\AdjUnitOn{A} = \gamma \AdjUnitOn{B} u = \bar{\gamma} u = \AdjUnitOn{A}
  \end{equation*}
  So, $\gamma (Tu)=\IdMapOn{TB}$ follows. - So $u$ is a $T$-equivalence; i.e.\ $\mathcal{L}^{\bot}\subseteq \mathcal{W}$.

  Finally, consider the case where $Y$ in $\EuScript{C}$ is orthogonal to every $T$-equivalence. Then $Y$ is orthogonal to the $T$-equivalence $\AdjUnitOn{Y}$. So, we obtain a left inverse for $\AdjUnitOn{Y}\from Y\to TY$, which is seen to be a right inverse as well. So, $Y$ is $T$-local, and $\mathcal{W}^{\bot}\subseteq \mathcal{L}$ follows. This completes the proof.
\end{proof}

We may interpret (\ref{thm:Reflector->OrthogonalPair}) as saying that $T\from \Ctgry{C}\to \Ctgry{R}$ is, up to equivalence of categories, a model for the Gabriel--Zisman localization $\Ctgry{C}\to \Ctgry{C}[\mathcal{W}^{-1}]$; see \cite{PGabrielMZisman1967}. Thus, while in general, there is no guarantee that the Gabriel--Zisman localization of a locally small category is again locally small, in the special case where $\mathcal{W}$ is the class of equivalences of a reflector, we do obtain a locally small model for  $\Ctgry{C}\to \Ctgry{C}[\mathcal{W}^{-1}]$.

Here are the basic properties of orthogonality, and of orthogonal pairs.

\begin{proposition}[Properties of orthogonality: morphism/object]
  \label{thm:Orthogonality:Morphism/Object-Props}
  In an arbitrary category $\Ctgry{C}$ show that the following hold:
  \begin{enumerate}[(i)]
    \item For every class $M$ of morphisms in a category $\Ctgry{C}$, $M\subseteq M^{\bot\bot}$.
    \item For every class $O$ of objects in $\Ctgry{C}$, $O\subseteq O^{\bot\bot}$.
    \item If $M\subseteq N$ are classes of morphisms in $\Ctgry{C}$, then $N^{\bot}\subseteq M^{\bot}$.
    \item If $O\subseteq P$ are classes of morphisms in $\Ctgry{C}$, then $P^{\bot}\subseteq O^{\bot}$.
    \item $\bot^3=\bot$, when applied to a class of objects, as well as a class of morphisms.
    \item For any class $M$ of morphisms in $\Ctgry{C}$, $(M^{\bot\bot},O^{\bot})$ is an orthogonal pair.
    \item For any class $O$ of objects in $\Ctgry{C}$, $(O^{\bot},O^{\bot\bot})$ is an orthogonal pair.
  \end{enumerate}
\end{proposition}

\begin{proposition}[Properties of orthogonal pairs]
  \label{thm:OrthogonalPairs-Props}
  For an orthogonal pair $(M,O)$ in a category $\Ctgry{C}$ the following hold:
  \begin{enumerate}[(i)]
    \item \label{thm:OrthogonalPairProps,2of3M}%
          `$2$ out of $3$'-rule in $M$: given a composite \ $A \overset{u}{\longrightarrow} B \overset{v}{\longrightarrow} C$, if $2$ out of $u$, $v$, and $vu$ are in $M$, so is the third.
    \item \label{thm:OrthogonalPairProps,IsoTest}%
          If $u\from A\to B$ is in $M$ and $A$, $B$ are in $O$, then $u$ is an isomorphism.
    \item \label{thm:OrthogonalPairProps,LclsClsdLmts}%
          The class $O$ is closed under limits: If $\phi\from D\to \Ctgry{C}$, $D$ small, has $\phi(j)\in O$, for all $j\in \Objcts{J}$, then $\LimOf{\phi}$ is in $O$.
    \item \label{thm:OrthgnlPrPrprts,EqvlncsClsdCLmts}%
          The class $M$ is closed under arbitrary colimits: i.e.\ if $\tau\from \phi_0\Rightarrow \phi_1$ is a natural transformation of functors $\phi_0,\phi_1\from D\to \Ctgry{C}$, $D$ small, with \ $\tau_d\from \phi_0(d)\to \phi_1(d)$ \ in $M$, for all $d$ in $D$, then $\CoLimOf{\tau}\from \CoLimOf{\phi_0}\to \CoLimOf{\phi_1}$ is in $M$.
  \end{enumerate}
\end{proposition}

\begin{corollary}[When $G$ preserves/reflects monomorphisms / isomorphisms]
  \label{thm:ReflectiveSubCatInclusion-Monos/Isos}
  Let $G\from \Ctgry{R}\to \Ctgry{C}$ be the inclusion of a reflective subcategory. Then the following hold:
  \begin{enumerate}
    \item $G$ preserves and reflects isomorphisms.
    \item $G$ reflects monomorphisms.
    \item If $\Ctgry{R}$ admits pullbacks, then the inclusion $G$ preserves monomorphisms.
  \end{enumerate}
\end{corollary}
\begin{proof}
  (i) and (ii) follow from the definitions. (iii) Given a map $f\from X\to Y$ in $\Ctgry{R}$, we know that $f$ is a monomorphism if and only if  $\bar{f}$ in the kernel pair diagram below is an isomorphism.
  \begin{equation*}
    \xymatrix@R=5ex@C=4em{
    K \ar[r]^-{\bar{f}} \ar[d] \PullLU{rd}&
    X \ar[d]^{f} \\
    X \ar[r]_-{f} &
    Y
    }
  \end{equation*}
  This pullback diagram in $\Ctgry{R}$ is also a pullback diagram in $\SACtgry{X}$ by (\ref{thm:ReflectiveSubCat-Limits}.i). If $f$ is a monomorphism in $\Ctgry{R}$, then $\bar{f}$ is an isomorphism. Applying $G$ to this diagram yields a pullback diagram in $\SACtgry{X}$, and $G\bar{f}$ is an isomorphism. So $Gf$ is a monomorphism.
\end{proof}

We leave it to the reader to dualize what we said about reflective subcategories and obtain corresponding properties of coreflective subcategories.

\begin{exercises}

\begin{exercise}[Reflector unit is an isomorphism]
  \label{exe:ReflectorUnit-Iso}
  If $(T,\AdjUnit)$ is a reflection of $\Ctgry{C}$ in $\Ctgry{R}$, show that $\AdjUnitOn{A}\from A\to TA$ is an isomorphism for every $A$ in $\Ctgry{R}$.
\end{exercise}

\begin{exercise}
  Prove Propositions \ref{thm:Orthogonality:Morphism/Object-Props} and \ref{thm:OrthogonalPairs-Props}.
\end{exercise}

\begin{exercise}[Sum of equivalences]
  \label{exe:SumEquivalences->Equivalence}
  Given an orthogonal pair $(M,O)$ in a category $\Ctgry{C}$, show that any coproduct of maps in $M$ belongs also to $M$.
\end{exercise}

\begin{exercise}[Cobase change preserves equivalences]
  \label{exe:Equivalences-PreservedByCoBaseChange}
  Given an orthogonal pair $(M,O)$ in a category $\Ctgry{C}$, and $f$ in $M$, show the following:
  \begin{enumerate}[(i)]
    \item If $g\from A\to X$ is arbitrary, then the map $\bar{f}$ below is also in $M$.
          \begin{equation*}
            \xymatrix@R=5ex@C=4em{
            A \ar[r]^-{f} \ar[d]_{g} \PushRD{rd} &
            B \ar[d]^{\bar{g}} \\
            X \ar[r]_-{\bar{f}} &
            Y
            }
          \end{equation*}
    \item The common section $s$ of $t_1$ and $t_2$ in the cokernel pair construction below is in $M$.
          \begin{equation*}
            \xymatrix@R=5ex@C=4em{
            A \ar[r]^-{f} \ar[d]_-{f} \PushRD{rd} &
            B \ar@{=}@/^2ex/[rdd] \ar[d]_{t_1}  \\
            B \ar[r]^-{t_2} \ar@{=}@/_2ex/[rrd] &
            D \ar@{.>}[rd]|-{\ s\ } \\
            && B
            }
          \end{equation*}
  \end{enumerate}
  Further, if $g_1,g_2\from B\to Z$ satisfy $g_1f=g_2$, let $g\from D\to Z$ we the universal map induced by the pushout. Show that
\end{exercise}
\end{exercises}

\lhead{\bfseries\footnotesize
  References}
\chead{\fontsize{9}{10}\selectfont{\sffamily
  }}
\rhead{\fontsize{9}{10}\selectfont{\sffamily
  }}

\providecommand{\noopsort}[1]{}
\providecommand{\bysame}{\leavevmode\hbox to3em{\hrulefill}\thinspace}
\providecommand{\MR}{\relax\ifhmode\unskip\space\fi MR }
\providecommand{\MRhref}[2]{%
  \href{http://www.ams.org/mathscinet-getitem?mr=#1}{#2}
}
\providecommand{\href}[2]{#2}

%
%
\cleardoublepage
\phantomsection
\lhead{\bfseries\footnotesize
  Index}

\addcontentsline{toc}{part}{Index}
\printindex
\begin{theindex}

  \item $(3\prdct 3)$-Lemma, \hyperpage{191}
  \item $(\Prdct {3}{3})$-Lemma
  \subitem border cases, \hyperpage{122}
  \item $(\Prdct {3}{3})$-diagram, \hyperpage{71}
  \item $1$-point union of pointed sets, \hyperpage{25}
  \item $5$-Lemma, \hyperpage{64}, \hyperpage{118}, \hyperpage{186}
  \subitem easy, \hyperpage{186}
  \item $\SetsBsd $
  \subitem co-semiabelian, \hyperpage{225}
  \subitem not p-exact, \hyperpage{225}
  \item $n$-fold extension, \hyperpage{135}
  \item \ZExact \ category, \hyperpage{34}
  \item `old $=$ new', \hyperpage{224}

  \indexspace

  \item abelian
  \subitem category, \hyperpage{279}
  \subitem varieties of algebras, \hyperpage{283}
  \item abelian category, \hyperpage{270}
  \subitem if and only if p-exact and additive, \hyperpage{282}
  \subitem if and only if p-exact and diexact, \hyperpage{282}
  \subitem is semiabelian, \hyperpage{280}
  \item abelian group object
  \subitem antidiagonal, \hyperpage{208}
  \item absolute
  \subitem diagram property, \hyperpage{300}
  \subitem epimorphism, \hyperpage{140}, \hyperpage{300}
  \item absolute epimorphism
  \subitem $+$ monomorphism is isomorphism, \hyperpage{303}
  \item additive
  \subitem category, \hyperpage{269}
  \item additive category
  \subitem satisfies \KSGInline -condition, \hyperpage{275}
  \item almost
  \subitem abelian categories - examples, \hyperpage{281},
  \hyperpage{283}
  \item almost abelian category, \hyperpage{281}
  \item antidiagonal, \hyperpage{208}
  \item antinormal
  \subitem (de)composition, \hyperpage{72}
  \subitem decomposition in abelian category, \hyperpage{82}
  \subitem morphism, \hyperpage{72}
  \subitem pair, \hyperpage{72}
  \item arrow
  \subitem category, \hyperpage{44}

  \indexspace

  \item Barr exact
  \subitem $\Grps $, \hyperpage{237}
  \item Barr exact category, \hyperpage{262}
  \item base change, \hyperpage{295}
  \subitem along normal epi preserves/reflects (normal) monos,
  \hyperpage{171}
  \subitem preserves isomorphism, \hyperpage{298}
  \subitem preserves monomorphisms, \hyperpage{298}
  \subitem preserves surjective function, \hyperpage{298}
  \subitem reflects isomorphisms, \hyperpage{187}
  \item bicomplete category, \hyperpage{294}
  \item biproduct
  \subitem in linear category, \hyperpage{254}
  \subitem in pointed category, \hyperpage{251}
  \subitem in preadditive category, \hyperpage{271}
  \subitem recognition, \hyperpage{190}
  \item boundary operator of a chain complex, \hyperpage{92}

  \indexspace

  \item category
  \subitem (finitely) bicomplete, \hyperpage{294}
  \subitem (finitely) cocomplete, \hyperpage{294}
  \subitem (finitely) complete, \hyperpage{294}
  \subitem abelian, \hyperpage{270}
  \subitem additive, \hyperpage{269}
  \subitem almost abelian, \hyperpage{281}
  \subitem Barr exact, \hyperpage{262}
  \subitem di-exact, \hyperpage{74, 75}, \hyperpage{270}
  \subitem discrete, \hyperpage{294}
  \subitem finite, \hyperpage{293}
  \subitem homological, \hyperpage{161}, \hyperpage{223}
  \subitem linear, \hyperpage{253}
  \subitem Mal'tsev, \hyperpage{217}
  \subitem normal, \hyperpage{140}, \hyperpage{152}
  \subitem of arrows, \hyperpage{44}
  \subitem of morphisms, \hyperpage{44}
  \subitem of sectioned epis over $R$, \hyperpage{212}
  \subitem of split epimorphisms, \hyperpage{183}
  \subitem p-exact, \hyperpage{270}
  \subitem pointed, \hyperpage{21}
  \subitem preabelian, \hyperpage{269}, \hyperpage{277}
  \subitem preadditive, \hyperpage{269}
  \subitem quasi abelian, \hyperpage{281}
  \subitem regular, \hyperpage{315}
  \subitem semiabelian, \hyperpage{223}
  \subitem small, \hyperpage{293}
  \subitem Z.~Janelidze-normal, \hyperpage{152}
  \item central
  \subitem morphism, \hyperpage{198}
  \item chain complex, \hyperpage{92}
  \subitem boundary operator, \hyperpage{92}
  \subitem cosubnormal, \hyperpage{93}
  \subitem differential, \hyperpage{92}
  \subitem morphism, \hyperpage{98}
  \subitem normal, \hyperpage{93}
  \subitem subnormal, \hyperpage{93}
  \subitem weakly normal, \hyperpage{93}
  \item closed symmetric monoidal
  \subitem structure in $\SetsBsd $, \hyperpage{26}
  \item co-semiabelian category
  \subitem $\SetsBsd $, \hyperpage{225}
  \item cobase change, \hyperpage{296, 297}
  \subitem preserves epimorphisms, \hyperpage{299}
  \subitem preserves extremal epimorphisms, \hyperpage{147}
  \subitem preserves isomorphism, \hyperpage{299}
  \subitem preserves regular epimorphisms, \hyperpage{147}
  \subitem preserves strong epimorphisms, \hyperpage{147}
  \item cocomplete category, \hyperpage{294}
  \item coequalizer
  \subitem in additive category, \hyperpage{276}
  \subitem reflexive, \hyperpage{238}
  \item coimage
  \subitem in a normal mono factorization, \hyperpage{47}
  \item cokernel, \hyperpage{29}
  \subitem of $\ZeroObject \to X$, \hyperpage{33}
  \subitem of a composite of maps, \hyperpage{37}
  \subitem of an epimorphism, \hyperpage{33}
  \subitem of sum of maps, \hyperpage{38}
  \subitem same as extremal epimorphism, \hyperpage{145}
  \subitem same as strong epimorphism, \hyperpage{145}
  \item colimit
  \subitem in $\NEpiCat {X}$, \hyperpage{51}
  \subitem in $\NMonoCat {X}$, \hyperpage{51}
  \subitem in $\SetsBsd $, \hyperpage{26}
  \subitem in $\Tops $, \hyperpage{28}
  \subitem in $\TopsBsd $, \hyperpage{28}
  \item commutating
  \subitem morphisms, \hyperpage{197}
  \item commutative
  \subitem internal magma/monoid/group, \hyperpage{244}
  \subitem object, \hyperpage{199}
  \item commuting
  \subitem subobjects, \hyperpage{200}
  \item complete category, \hyperpage{294}
  \item composite
  \subitem of normal epimorphisms, \hyperpage{146}
  \item connecting map
  \subitem in Snake Lemma, \hyperpage{86}
  \item conservative
  \subitem functor, \hyperpage{213}
  \item coproduct
  \subitem as colimit of discrete diagram, \hyperpage{294}
  \item coproduct/product comparison map, \hyperpage{169}
  \item cosubnormal
  \subitem chain complex, \hyperpage{93}
  \subitem morphism, \hyperpage{50}
  \item cover
  \subitem of $Q$ by $K$, \hyperpage{54}

  \indexspace

  \item di-exact category, \hyperpage{74, 75}, \hyperpage{270}
  \subitem equivalent siblings, \hyperpage{128}
  \item di-extension, \hyperpage{71}
  \subitem from normal antinormal pairs, \hyperpage{73}
  \item di-extensive
  \subitem normal pullback, \hyperpage{121}
  \subitem normal pushout, \hyperpage{121}
  \subitem pullback, \hyperpage{127}, \hyperpage{133}
  \subitem pushout, \hyperpage{127}
  \subitem pushout in \HTag , \hyperpage{179}
  \item diagram, \hyperpage{300}
  \subitem with absolute property, \hyperpage{300}
  \item difference
  \subitem object, \hyperpage{153}
  \subitem object in $\AbGrps $, \hyperpage{153}
  \item differential of a chain complex, \hyperpage{92}
  \item dinversion, \hyperpage{73}
  \item direct image
  \subitem of a relation under a normal epimorphism, \hyperpage{265}
  \item discrete
  \subitem category, \hyperpage{294}
  \subitem relation, \hyperpage{259}

  \indexspace

  \item effective
  \subitem epimorphism, \hyperpage{307}
  \subitem epimorphism - Hom-view, \hyperpage{311}
  \subitem equivalence relation, \hyperpage{262}
  \subitem equivalence relation: properties, \hyperpage{263}
  \item epic morphism, \hyperpage{19}
  \item epimorphism, \hyperpage{19}
  \subitem absolute, \hyperpage{140}, \hyperpage{300}
  \subitem effective, \hyperpage{307}
  \subitem extremal, \hyperpage{305, 306}
  \subitem is preserved by pushout, \hyperpage{299}
  \subitem normal, \hyperpage{30}
  \subitem properties, \hyperpage{20}
  \subitem regular, \hyperpage{307}
  \subitem split, \hyperpage{300}
  \subitem strong, \hyperpage{307}, \hyperpage{309}
  \subitem via pushout property, \hyperpage{309}
  \item equalizer
  \subitem in additive category, \hyperpage{276}
  \item equivalence
  \subitem of a reflector, \hyperpage{319}
  \item equivalence relation, \hyperpage{260}
  \subitem given by kernel pair, \hyperpage{260}
  \item equivalent
  \subitem extensions, \hyperpage{189}
  \subitem monomorphisms, \hyperpage{20}
  \item exact
  \subitem diagram of sets, \hyperpage{310}
  \subitem node recognize, \hyperpage{55}
  \subitem sequence, \hyperpage{55}
  \subitem short exact, \hyperpage{55}
  \item exactness
  \subitem of product of short exact sequence , \hyperpage{195}
  \item extension, \hyperpage{189}
  \subitem di-, \hyperpage{71}
  \subitem morphism, \hyperpage{189}
  \subitem of $K$ by $Q$, \hyperpage{54}
  \item extremal epimorphism, \hyperpage{305, 306}
  \subitem same as cokernel, \hyperpage{145}
  \subitem same as strong epimorphism, \hyperpage{145}

  \indexspace

  \item finite
  \subitem (co-)limit, \hyperpage{294}
  \subitem category, \hyperpage{293}
  \item finitely
  \subitem bicomplete category, \hyperpage{294}
  \subitem cocomplete category, \hyperpage{294}
  \subitem complete category, \hyperpage{294}
  \item formal difference
  \subitem of two morphisms, \hyperpage{154}
  \item functor
  \subitem conservative, \hyperpage{213}

  \indexspace

  \item graph, \hyperpage{258}
  \subitem internal, \hyperpage{258}
  \item group, \hyperpage{244, 245}

  \indexspace

  \item higher extension, \hyperpage{135}
  \item homological
  \subitem category, \hyperpage{161}
  \subitem category satisfies \DPNInline , \hyperpage{179}
  \subitem self-duality, \hyperpage{111}
  \item homological category, \hyperpage{223}
  \subitem is Mal'tsev category, \hyperpage{219}
  \subitem torsion-free abelian groups, \hyperpage{163}
  \subitem torsion-free groups, \hyperpage{163}
  \item homology
  \subitem functorial, \hyperpage{99}
  \item HSD
  \subitem property of normal category, \hyperpage{142}

  \indexspace

  \item image
  \subitem factorization, \hyperpage{315}
  \subitem factorization in a normal category, \hyperpage{143}
  \subitem factorization: strong epi/mono, \hyperpage{310}
  \subitem in a normal epi factorization, \hyperpage{47}
  \subitem of a morphism, \hyperpage{315}
  \subitem of a morphism in a normal category, \hyperpage{143}
  \item image factorization, \hyperpage{143}
  \item indiscrete
  \subitem relation, \hyperpage{260}
  \item initial
  \subitem normal mono decomposition, \hyperpage{46}
  \subitem object, \hyperpage{21}
  \item internal
  \subitem abelian group structure, \hyperpage{204}
  \subitem graph from $X$ to $Y$, \hyperpage{258}
  \subitem group, \hyperpage{244}
  \subitem magma, \hyperpage{244}
  \subitem monoid, \hyperpage{244}
  \subitem relation from $X$ to $Y$, \hyperpage{258}
  \subitem unitary magma, \hyperpage{244}
  \subitem unitary magma structure, \hyperpage{204}
  \item internal group
  \subitem inverse operation $i^2=\IdMap $, \hyperpage{248}
  \item intersection
  \subitem of normal subobjects, \hyperpage{42}
  \subitem of subobjects, \hyperpage{18}
  \subitem of subobjects: existence, \hyperpage{299}
  \item inverse map
  \subitem of internal group, \hyperpage{244}
  \item inversion
  \subitem normal subobjects / quotient objects, \hyperpage{34}
  \item isomorphism
  \subitem recognition, \hyperpage{36}
  \subitem recognition in abelian category, \hyperpage{280}
  \item Isomorphism Theorem
  \subitem condition for Third, \hyperpage{111}, \hyperpage{113}
  \subitem Second, \hyperpage{234}
  \subitem Third, \hyperpage{38}, \hyperpage{195}

  \indexspace

  \item join
  \subitem of (normal) subobjects in $\NNr $, \hyperpage{257}
  \subitem of subobjects: construction, \hyperpage{19},
  \hyperpage{232}
  \subitem of subobjects: definition, \hyperpage{19}
  \item jointly
  \subitem epimorphic family of maps, \hyperpage{19}
  \subitem extremal-epimorphic family of maps, \hyperpage{197},
  \hyperpage{306}
  \subitem monomorphic family - from limit cone, \hyperpage{298}
  \subitem monomorphic family - product description, \hyperpage{298}
  \subitem monomorphic family of maps, \hyperpage{17}

  \indexspace

  \item kernel, \hyperpage{29}
  \subitem composition, \hyperpage{32}
  \subitem is monomorphism, \hyperpage{32}
  \subitem of a composite of maps, \hyperpage{36}
  \subitem of monomorphism: $\ZeroObject $., \hyperpage{35}
  \subitem of product of maps, \hyperpage{38}
  \subitem pair of a morphism, \hyperpage{296}
  \item kernel pair
  \subitem equivalence relation, \hyperpage{260}
  \item KSG
  \subitem in additive categories, \hyperpage{275}

  \indexspace

  \item lattice
  \subitem modular, \hyperpage{236}
  \item left
  \subitem almost abelian category, \hyperpage{280}
  \subitem orthogonal, \hyperpage{148}, \hyperpage{312}
  \subitem orthogonal complement, \hyperpage{148}, \hyperpage{312}
  \subitem quasi-abelian category, \hyperpage{280}
  \item limit
  \subitem in $\NEpiCat {X}$, \hyperpage{51}
  \subitem in $\NMonoCat {X}$, \hyperpage{51}
  \subitem in $\Tops $, \hyperpage{28}
  \subitem in $\TopsBsd $, \hyperpage{28}
  \item linear
  \subitem category, \hyperpage{184}
  \item linear category, \hyperpage{253}
  \item local
  \subitem object w.r. to a reflector, \hyperpage{319}
  \item long exact homology sequence, \hyperpage{102}
  \item long exact sequence, \hyperpage{56}

  \indexspace

  \item magma, \hyperpage{245}
  \item Mal'tsev category, \hyperpage{217}
  \item meet
  \subitem of normal subobjects, \hyperpage{42}
  \subitem of subobjects: construction, \hyperpage{19}
  \subitem of subobjects: definition, \hyperpage{18}
  \subitem of subobjects: existence, \hyperpage{299}
  \item modular lattice, \hyperpage{236}
  \item monic morphism, \hyperpage{17}
  \item monoid, \hyperpage{244, 245}
  \item monomorphism, \hyperpage{17}
  \subitem is preserved by pullback, \hyperpage{298}
  \subitem normal, \hyperpage{30}
  \subitem properties, \hyperpage{20}
  \subitem recognition by $0$-kernel, \hyperpage{142}
  \subitem recognition by $\ZeroObject $-kernel, \hyperpage{141}
  \subitem recognition in map of exact sequences, \hyperpage{173}
  \subitem retractable, \hyperpage{301}
  \subitem via limit property, \hyperpage{309}
  \item morphism
  \subitem central, \hyperpage{198}
  \subitem cosubnormal, \hyperpage{50}
  \subitem in $\SEpisInOver {X}{R}$, \hyperpage{213}
  \subitem normal, \hyperpage{48}
  \subitem of chain complexes, \hyperpage{98}
  \subitem of extensions, \hyperpage{189}
  \subitem of internal magma/monoid/group, \hyperpage{249}
  \subitem of zero maps - factorization, \hyperpage{57}
  \subitem retractable, \hyperpage{301}
  \subitem split, \hyperpage{300}
  \subitem split short exact sequences, \hyperpage{164}
  \subitem subnormal, \hyperpage{50}
  \subitem zero, \hyperpage{21}
  \item multiplication
  \subitem of $f$ by $g$, \hyperpage{197}

  \indexspace

  \item nilpotent
  \subitem group via cosubnormal factorization, \hyperpage{53}
  \item normal
  \subitem  map preserved by normal images, \hyperpage{78}
  \subitem category, \hyperpage{140}, \hyperpage{152}
  \subitem category - homologically self-dual, \hyperpage{142}
  \subitem chain complex, \hyperpage{93}
  \subitem closure - existence, \hyperpage{43}
  \subitem closure of a subobject, \hyperpage{42}
  \subitem epi decomposition, \hyperpage{45}
  \subitem epi decomposition of order $r$, \hyperpage{52}
  \subitem epi decomposition, terminal, \hyperpage{46}
  \subitem epi factorization, \hyperpage{46}
  \subitem epi from project projection, \hyperpage{257}
  \subitem epimorphism, \hyperpage{30}
  \subsubitem in short exact sequence, \hyperpage{176}
  \subitem epimorphism - composite, \hyperpage{146}
  \subitem epimorphism - product, \hyperpage{146}
  \subitem epimorphism - properties, \hyperpage{146}
  \subitem epimorphism - recognition, \hyperpage{182}
  \subitem epimorphism in $\NEpiCat {X}$, \hyperpage{51}
  \subitem epimorphism in $\SetsBsd $, \hyperpage{34}
  \subitem factorization, \hyperpage{48}
  \subitem factorization - recognize, \hyperpage{48}
  \subitem factorization, recognize, \hyperpage{46}
  \subitem images preserve normal monos / epis, \hyperpage{124}
  \subitem map, \hyperpage{48}
  \subitem map factored through epic map, \hyperpage{49}
  \subitem map factored through monic map, \hyperpage{48}
  \subitem map in morphism of short exact sequences, \hyperpage{180}
  \subitem mono decomposition, \hyperpage{45}
  \subitem mono decomposition of order $r$, \hyperpage{52}
  \subitem mono decomposition, initial, \hyperpage{46}
  \subitem mono factorization, \hyperpage{46}
  \subitem monomorphism, \hyperpage{30}
  \subitem monomorphism in $\NMonoCat {X}$, \hyperpage{51}
  \subitem monomorphism in $\SetsBsd $, \hyperpage{34}
  \subitem monomorphism in morphisms of SESs, \hyperpage{177}
  \subitem monomorphism: properties, \hyperpage{36}
  \subitem pullback, \hyperpage{134}
  \subitem pullback in $\CMon $, \hyperpage{257}
  \subitem pullback in $\Grps $, \hyperpage{134}
  \subitem pushout, \hyperpage{131}
  \subitem pushout recognition, \hyperpage{183}
  \subitem subobject, \hyperpage{30}
  \subitem subobject of $\NNr $, \hyperpage{257}
  \subitem subtractive category, \hyperpage{155}
  \item normal epimorphism, \hyperpage{15}
  \subitem given by $\SumMapOutOf {\IdMap ,\ZeroMap }\from X+Y\NEpi X$,
  \hyperpage{38}
  \item normal monomorphism, \hyperpage{15}
  \subitem given by $\PrdctMapInto {\IdMap ,\ZeroMap }\from X\to \Prdct {X}{Y}$,
  \hyperpage{38}
  \item normal pullback, \hyperpage{131}
  \subitem di-extensive, \hyperpage{121}
  \item normal pushout, \hyperpage{131}
  \subitem di-extensive, \hyperpage{121}
  \item normal subobject / quotient object inversion, \hyperpage{34}

  \indexspace

  \item object
  \subitem commutative, \hyperpage{199}
  \subitem initial, \hyperpage{21}
  \subitem terminal, \hyperpage{21}
  \subitem zero, \hyperpage{21}
  \item OFS, \hyperpage{148}, \hyperpage{313}
  \item OPFS - orthogonal prefactorization system, \hyperpage{312}
  \item opposite
  \subitem relation, \hyperpage{259}
  \item orthogonal, \hyperpage{318}
  \subitem complement, \hyperpage{148}, \hyperpage{312}
  \subitem factorization system, \hyperpage{148}, \hyperpage{313}
  \subitem pair, \hyperpage{319}
  \subitem prefactorization system, \hyperpage{312}
  \item orthogonal complement
  \subitem of a class of morphisms, \hyperpage{319}
  \subitem of a class of objects, \hyperpage{319}

  \indexspace

  \item p-exact category, \hyperpage{126}, \hyperpage{270}
  \item Palamodov
  \subitem semiabelian category, \hyperpage{284}
  \item pointed
  \subitem category, \hyperpage{21}
  \item preabelian category, \hyperpage{269}, \hyperpage{277}
  \subitem recognize, \hyperpage{277}
  \item preadditive category, \hyperpage{269}
  \subitem need not be pointed, \hyperpage{272}
  \item product
  \subitem as limit of discrete diagram, \hyperpage{294}
  \subitem as pullback, \hyperpage{23}
  \subitem of normal epimorphisms, \hyperpage{146}
  \subitem of short exact sequences, \hyperpage{146}, \hyperpage{195}
  \subitem projection is normal epi, \hyperpage{257}
  \subitem recognition, \hyperpage{167}
  \subitem recognition via kernel splitting, \hyperpage{187}
  \item proper
  \subitem morphism, \hyperpage{48}
  \item pullback, \hyperpage{133}, \hyperpage{295}
  \subitem along normal epi preserves/reflects (normal) monos,
  \hyperpage{171}
  \subitem cancellation, \hyperpage{295}
  \subitem cancellation $2$-out-of-$3$, \hyperpage{173},
  \hyperpage{188}
  \subitem composition, \hyperpage{295}
  \subitem concatenated squares, \hyperpage{295}
  \subitem construction of product, \hyperpage{23}
  \subitem di-extensive, \hyperpage{127}
  \subitem normal, \hyperpage{131}
  \subitem preserves monomorphism, \hyperpage{298}
  \subitem recognition by kernel isomorphism, \hyperpage{175}
  \subitem recognition, kernel side I, \hyperpage{40}
  \subitem recognition: kernel side II, \hyperpage{172}
  \subitem sectioned epimorphism, \hyperpage{302}
  \subitem seminormal, \hyperpage{131}
  \item Puppe
  \subitem exact category, \hyperpage{126}
  \item pure snake
  \subitem condition, \hyperpage{111}
  \item pushout
  \subitem cancellation, \hyperpage{297}
  \subitem concatenated squares, \hyperpage{297}
  \subitem di-extensive, \hyperpage{127}
  \subitem isomophic cokernels, \hyperpage{41}
  \subitem normal, \hyperpage{131}
  \subitem of a sectioned monomorphism, \hyperpage{302}
  \subitem of normal epimorphisms, \hyperpage{43}
  \subitem preserves epimorphism, \hyperpage{299}
  \subitem recognition - cokernel side, \hyperpage{41}
  \subitem recognition II, \hyperpage{183}
  \subitem seminormal, \hyperpage{131, 132}

  \indexspace

  \item quasi-abelian category, \hyperpage{281}
  \item quotient
  \subitem object, \hyperpage{20}

  \indexspace

  \item Raïkov
  \subitem semiabelian category, \hyperpage{284}
  \item reflective subcategory, \hyperpage{317}
  \item reflector, \hyperpage{317}
  \item reflexive coequalizer, \hyperpage{238}
  \item reflexive relation, \hyperpage{259}
  \item regular
  \subitem category, \hyperpage{315}
  \subitem epimorphism, \hyperpage{307}
  \subitem pushout, \hyperpage{227}
  \item relation, \hyperpage{258}
  \subitem discrete, \hyperpage{259}
  \subitem from $X$ to $Y$, \hyperpage{258}
  \subitem indiscrete, \hyperpage{260}
  \subitem on $X$, \hyperpage{258}
  \subitem opposite, \hyperpage{259}
  \subitem reflexive, \hyperpage{259}
  \subitem symmetric, \hyperpage{259}
  \subitem transitive, \hyperpage{260}
  \item retractable
  \subitem monomorphism, \hyperpage{301}
  \item retractable monomorphism, \hyperpage{301}
  \item right
  \subitem almost abelian category, \hyperpage{280}
  \subitem orthogonal, \hyperpage{148}, \hyperpage{312}
  \subitem orthogonal complement, \hyperpage{148}, \hyperpage{312}
  \subitem quasi-abelian category, \hyperpage{280}

  \indexspace

  \item Second Isomorphism Theorem, \hyperpage{234}
  \item section, \hyperpage{300}
  \item sectionable
  \subitem epimorphism, \hyperpage{300}
  \subitem epimorphism is epimorphism, \hyperpage{303}
  \subitem morphism, \hyperpage{300}
  \item sectioned epimorphism
  \subitem $+$ monomorphism is isomorphism, \hyperpage{303}
  \item semiabelian category, \hyperpage{223}
  \subitem in sense of Palamodov, \hyperpage{284}
  \subitem in sense of Raïkov, \hyperpage{284}
  \subitem is Barr-exact, \hyperpage{238}
  \subitem is Mal'tsev category, \hyperpage{238}
  \item seminormal
  \subitem pullback, \hyperpage{131}, \hyperpage{133}
  \subitem pushout, \hyperpage{131, 132}
  \subitem pushout - pullback stability, \hyperpage{179}
  \subitem pushout and HSD, \hyperpage{132, 133}
  \subitem pushout recognition, \hyperpage{178}
  \item seminormal pushout
  \subitem recognition II, \hyperpage{183}
  \item sequence
  \subitem short exact, \hyperpage{54}
  \item Short 5-Lemma, \hyperpage{64}, \hyperpage{117}, \hyperpage{185}
  \item short exact sequence, \hyperpage{54, 55}
  \subitem in $\SetsBsd $, \hyperpage{57}
  \subitem of chain complexes, \hyperpage{100}
  \subitem product, \hyperpage{146}
  \item small
  \subitem category, \hyperpage{293}
  \item smash product
  \subitem in $\SetsBsd $, \hyperpage{26}
  \item snake
  \subitem map/morphism, \hyperpage{86}
  \item Snake Lemma
  \subitem classical - conditionally extensive cats, \hyperpage{84}
  \subitem classical in \DPNInline -normal cat, \hyperpage{149}
  \subitem pure version, \hyperpage{83}
  \subitem pure version, relaxed, \hyperpage{90}
  \subitem relaxed version, \hyperpage{87}
  \subitem strong relaxed version, \hyperpage{106}
  \item split
  \subitem monomorphism in $\SetsBsd $, \hyperpage{27}
  \subitem short exact sequence, \hyperpage{140}
  \subitem short exact sequence - morphism, \hyperpage{164}
  \subitem short exact sequence from split epi\NTag , \hyperpage{141}
  \item splitting, \hyperpage{300}
  \item strong epimorphism, \hyperpage{307}, \hyperpage{309}
  \subitem same as cokernel, \hyperpage{145}
  \subitem same as extremal epimorphism, \hyperpage{145}
  \item subnormal
  \subitem chain complex, \hyperpage{93}
  \subitem epi decomposition, \hyperpage{50}
  \subitem mono decomposition, \hyperpage{50}
  \subitem morphism, \hyperpage{50}
  \item subobject, \hyperpage{17}
  \subitem in $\Sets $, \hyperpage{18}
  \subitem in $\Tops $, \hyperpage{18}, \hyperpage{28}
  \subitem in $\TopsBsd $, \hyperpage{28}
  \subitem normal closure, \hyperpage{42}
  \subitem of $\NNr $, \hyperpage{257}
  \item subobjects
  \subitem intersection: existence, \hyperpage{299}
  \subitem intersections, \hyperpage{18}
  \subitem join, \hyperpage{19}
  \subitem meet, \hyperpage{18}
  \subitem meet: existence, \hyperpage{299}
  \subitem union, \hyperpage{19}
  \item subtractive
  \subitem normal category, \hyperpage{155}
  \subitem variety, \hyperpage{126}, \hyperpage{142}
  \item sum
  \subitem of short exact sequences, \hyperpage{177}
  \item sum/product comparison map, \hyperpage{169}
  \item symmetric
  \subitem relation, \hyperpage{259}

  \indexspace

  \item terminal
  \subitem normal epi decomposition, \hyperpage{46}
  \subitem object, \hyperpage{21}
  \item Third Isomorphism Theorem, \hyperpage{195}
  \item topological
  \subitem monoids, \hyperpage{126}
  \subitem variety, \hyperpage{126}
  \item torsion-free
  \subitem abelian groups form homological cat, \hyperpage{163}
  \subitem groups form homological cat, \hyperpage{163}
  \item totally
  \subitem normal sequence of epimorphisms, \hyperpage{108}
  \subitem normal sequence of monomorphisms, \hyperpage{108}
  \item transitive relation, \hyperpage{260}
  \item twist map, \hyperpage{244}

  \indexspace

  \item union
  \subitem of subobjects, \hyperpage{19}
  \item universal property
  \subitem pullback, \hyperpage{295}

  \indexspace

  \item variety
  \subitem subtractive, \hyperpage{142}

  \indexspace

  \item weakly normal chain complex, \hyperpage{93}
  \item wedge in $\SetsBsd $, \hyperpage{25}

  \indexspace

  \item Z.~Janelidze-normal category, \hyperpage{152}
  \item zero
  \subitem morphism, \hyperpage{21}
  \subitem object, \hyperpage{21}

\end{theindex}

\lhead{\bfseries\footnotesize
  Notation Index}
\renewcommand{\twocolumn}[1][]{#1}
\addcontentsline{toc}{part}{Notation Index}

\begin{theindex}

  \item a
  \subitem $\ANPCat {X}$\IndSep category of antinormal pairs in $\Ctgry {X}$,
  \hyperpage{76}
  \subitem $\ArrowCat {\Ctgry {X}}$\IndSep category of morphisms in $\Ctgry {X}$,
  \hyperpage{44}

  \indexspace

  \item b
  \subitem $\BiPrdct {X}{Y}$\IndSep biproduct of $X$ and $Y$,
  \hyperpage{251}

  \indexspace

  \item c
  \subitem $(C,d^C)$\IndSep chain complex $C$ with differential $d^C$,
  \hyperpage{92}
  \subitem $\CHTops $\IndSep category of compact Hausdorff spaces,
  \hyperpage{18}
  \subitem $\CoImg {f}$\IndSep codomain of coimage map of $f$,
  \hyperpage{47}
  \subitem $\CoImgMap {f}$\IndSep coimage map of $f$, \hyperpage{47}
  \subitem $\CoKer {f}$\IndSep cokernel object of $f$, \hyperpage{29}
  \subitem $\CoKerMap {f}$\IndSep cokernel of $f$, \hyperpage{29}
  \subitem $\CoLimOfOver {F}{J}$\IndSep colimit of functor $F\from J\to \EuScript {X}$,
  \hyperpage{294}
  \subitem $\FamCoPrdct {j\in J}{F(j)}$\IndSep coproduct over discrete category $J$,
  \hyperpage{294}
  \subitem $\SumProdComp {A}{B}$\IndSep comparison map $A+B\longrightarrow \Prdct {A}{B}$,
  \hyperpage{169}
  \subitem $\bigvee _{k\in K} (X_k,x_k)$\IndSep wedge/coproduct in $\SetsBsd $,
  \hyperpage{25}

  \indexspace

  \item d
  \subitem $D(X)$\IndSep difference object of $X$, \hyperpage{153}
  \subitem $\DExCat {X}$\IndSep category of di-extensions in $\Ctgry {X}$,
  \hyperpage{76}
  \subitem $f-g$\IndSep formal difference of morphisms,
  \hyperpage{154}

  \indexspace

  \item e
  \subitem $Y\Epi Z$\IndSep epimorphism, \hyperpage{19}
  \subitem $\ExtMini {V}{K}$\IndSep equivalence classes of extensions,
  \hyperpage{189}
  \subitem $\HExCat {n}{X}$\IndSep category of $n$-fold extensions in $\Ctgry {X}$,
  \hyperpage{135}
  \subitem $\NENMComp {f}$\IndSep normal epi / normal mono factorization comparison map,
  \hyperpage{46}

  \indexspace

  \item g
  \subitem $\Grps $\IndSep category of groups, \hyperpage{22}

  \indexspace

  \item h
  \subitem $\HmlgyCoKer {n}{C}$\IndSep cokernel construction of homology,
  \hyperpage{94}
  \subitem $\HmlgyKer {n}{C}$\IndSep kernel construction of homology,
  \hyperpage{95}

  \indexspace

  \item i
  \subitem $M\meet N$\IndSep intersection/meet of $M$ and $N$,
  \hyperpage{18}
  \subitem $\Img {f}$\IndSep domain of image map of $f$,
  \hyperpage{47}
  \subitem $\ImgMap {f}$\IndSep image map of $f$, \hyperpage{47}

  \indexspace

  \item j
  \subitem $M\join N$\IndSep join/union of subobjects $M$ and $N$,
  \hyperpage{19}
  \subitem $M\meet N$\IndSep intersection/meet of $M$ and $N$,
  \hyperpage{42}

  \indexspace

  \item k
  \subitem $\Ker {f}$\IndSep kernel object of $f$, \hyperpage{29}
  \subitem $\KerMap {f}$\IndSep kernel of $f$, \hyperpage{29}
  \subitem $\KrnlPr {f}$\IndSep kernel pair of $f$, \hyperpage{296}

  \indexspace

  \item l
  \subitem $\LimOfOver {F}{J}$\IndSep limit of functor $F\from J\to \EuScript {X}$,
  \hyperpage{293}

  \indexspace

  \item m
  \subitem $M\meet N$\IndSep intersection/meet of $M$ and $N$,
  \hyperpage{18}, \hyperpage{42}
  \subitem $X\Mono Y$\IndSep monomorphism, \hyperpage{17}
  \subitem $\Magmas $\IndSep category of unitary magmas,
  \hyperpage{245}

  \indexspace

  \item n
  \subitem $K\normal X$\IndSep $K$ is normal subobject of $X$,
  \hyperpage{30}
  \subitem $\NEpiCat {X}$\IndSep category of normal epimorphisms in $\Ctgry {X}$,
  \hyperpage{44}, \hyperpage{58}, \hyperpage{114}
  \subitem $\NMonoCat {X}$\IndSep category of normal monomorphisms in $\Ctgry {X}$,
  \hyperpage{44}, \hyperpage{58}, \hyperpage{114}
  \subitem $\NQuoObjcts {X}$\IndSep category of normal quotient objects of object $X$,
  \hyperpage{34}
  \subitem $\NSubObjcts {X}$\IndSep category of normal subobjects of object $X$,
  \hyperpage{34}
  \subitem $\kappa \from K\NMono X$\IndSep $\kappa $ is normal monomorphism,
  \hyperpage{30}
  \subitem $\pi \from Y\NEpi Q$\IndSep $\pi $ is normal epimorphism,
  \hyperpage{30}

  \indexspace

  \item p
  \subitem $\FamPrdct {j\in J}{F(j)}$\IndSep product over discrete category $J$,
  \hyperpage{294}
  \subitem $\PrdctMapInto {f,g}\from A\to \Prdct {X}{Y}$\IndSep universal map into a product,
  \hyperpage{294}
  \subitem $\PwrSt {X}$\IndSep power set of set $X$, \hyperpage{18}
  \subitem $f^{\ast }\varepsilon $\IndSep pullback of short exact sequence $\varepsilon $ along $f$,
  \hyperpage{170}

  \indexspace

  \item r
  \subitem $\Rngs $\IndSep category of rings with or without unit,
  \hyperpage{22}
  \subitem $\RtrctdMono {m}{r}$\IndSep mono $m$ retracted by $r$,
  \hyperpage{301}
  \subitem $\URngs $\IndSep category of unital rings, \hyperpage{22}
  \subitem $aRb$\IndSep $a$ is related to $b$, \hyperpage{258}

  \indexspace

  \item s
  \subitem $(X,x_0)\wedge (Y,y_0)$\IndSep smash of pointed sets,
  \hyperpage{26}
  \subitem $M\leq X$\IndSep $M$ is subobject of $X$, \hyperpage{17}
  \subitem $\SESCat {X}$\IndSep category of short exact sequences in $\Ctgry {X}$,
  \hyperpage{58}, \hyperpage{114}, \hyperpage{230}
  \subitem $\SEpisIn {X}$\IndSep category of sectioned morphisms in $\Ctgry {X}$,
  \hyperpage{301}
  \subitem $\SEpisInOver {X}{R}$\IndSep category of sectioned epis in $\Ctgry {X}$ over $R$,
  \hyperpage{212}
  \subitem $\SSESCat {X}$\IndSep category of split short exact sequences,
  \hyperpage{164}
  \subitem $\SctndEpi {q}{x}$\IndSep map $q$ sectioned by $x$,
  \hyperpage{300}
  \subitem $\Sets $\IndSep category of sets, \hyperpage{23}
  \subitem $\SetsBsd $\IndSep category of sets with base point,
  \hyperpage{22, 23}
  \subitem $\SubObjcts {X}$\IndSep class of subobjects of $X$,
  \hyperpage{17}
  \subitem $\SumMapOutOf {u,v}\from X+Y\to Z$\IndSep universal map out of a coproduct / sum,
  \hyperpage{295}

  \indexspace

  \item t
  \subitem $\Tops $\IndSep category of topological spaces,
  \hyperpage{18}

  \indexspace

  \item u
  \subitem $M\join N$\IndSep join/union of subobjects $M$ and $N$,
  \hyperpage{19}

  \indexspace

  \item w
  \subitem $\bigvee _{k\in K} (X_k,x_k)$\IndSep wedge/coproduct in $\SetsBsd $,
  \hyperpage{25}

  \indexspace

  \item z
  \subitem $\ZNr $\IndSep number system of integers, \hyperpage{22},
  \hyperpage{93}
  \subitem $\ZeroMap _{YX}$\IndSep zero morphism from $X$ into $Y$,
  \hyperpage{21}
  \subitem $\ZeroObject $\IndSep zero object, \hyperpage{21}

\end{theindex}

\lhead{\bfseries\footnotesize
  Acronyms}
\addcontentsline{toc}{part}{Acronyms}

\begin{theindex}

  \item a
  \subitem $\EuRoman {A{\kern -0.2ex}N{\kern -0.15ex}K}$\IndSep exercise asking a question whose answer we do not know,
  \hyperpage{63}
  \subitem \AENInline \IndSep condition that absolute epimorphisms are normal maps,
  \hyperpage{140}
  \subitem \ANNInline \IndSep property that antinormal maps are normal,
  \hyperpage{75}, \hyperpage{223}

  \indexspace

  \item d
  \subitem \DPNInline \IndSep property that dinversion preserves normal maps,
  \hyperpage{75}, \hyperpage{179}
  \subitem {\color {Cerulean} $\EuRoman {DEx}$}\IndSep category in which antinormal composites are normal,
  \hyperpage{74}
  \subitem {\color {Cerulean} $\EuRoman {DPN}$}\IndSep category in which dinversion preserves normal maps,
  \hyperpage{75}

  \indexspace

  \item h
  \subitem \HSDInline \IndSep homological self-duality property,
  \hyperpage{75}
  \subitem {\color {Brown} $\EuRoman {H}$}\IndSep homological category,
  \hyperpage{161}

  \indexspace

  \item k
  \subitem \KSGInline \IndSep for $q\from X\to Q$ sectioned by $x$, images of $x$ and $\KerMap {q}$ generate $X$,
  \hyperpage{161}, \hyperpage{223}

  \indexspace

  \item n
  \subitem {\color {Brown} $\EuRoman {N}$}\IndSep normal category,
  \hyperpage{140}

  \indexspace

  \item p
  \subitem \PNEInline \IndSep condition that pullbacks preserve normal epimorphisms,
  \hyperpage{140}, \hyperpage{161}

  \indexspace

  \item s
  \subitem {\color {MidnightBlue} $\EuRoman {S{\kern -0.15ex}A}$}\IndSep semiabelian category,
  \hyperpage{223}

  \indexspace

  \item z
  \subitem {\color {Cerulean} $\EuRoman {z\hy Ex}$}\IndSep \ZExact \ category,
  \hyperpage{34}

\end{theindex}

\end{document}